\documentclass[11pt, reqno]{amsart}

\usepackage{amssymb}
\usepackage{mathrsfs}
\usepackage{amsmath}
\usepackage{amsthm}
\usepackage{amsfonts}
\usepackage{latexsym}

\usepackage{color}
\usepackage{txfonts}
\usepackage{enumerate}

\usepackage[colorlinks=true,
  linkcolor=blue,
  citecolor=red,
  urlcolor=magenta, backref=page]{hyperref}

\usepackage{txfonts}

\usepackage{anysize}

\allowdisplaybreaks

\newtheorem{theorem}{Theorem}[section]
\newtheorem{proposition}[theorem]{Proposition}
\newtheorem{corollary}[theorem]{Corollary}
\newtheorem{lemma}[theorem]{Lemma}

\theoremstyle{definition}
\newtheorem{definition}[theorem]{Definition}
\newtheorem{example}[theorem]{Example}

\newtheorem{algorithm}[theorem]{Algorithm}

\theoremstyle{remark}
\newtheorem{remark}[theorem]{Remark}

\numberwithin{equation}{section}

\newcommand{\A}{\mathscr A}
\newcommand{\As}{\mathcal A}
\newcommand{\C}{\mathbb C}

\newcommand{\D}{\mathscr D}
\newcommand{\vD}{\vec D}
\newcommand{\E}{\mathbb E}
\newcommand{\vE}{\vec E}
\newcommand{\Es}{\mathscr E}
\newcommand{\F}{\mathbb F}
\newcommand{\vF}{\vec F}
\newcommand{\Fs}{\mathscr F}

\newcommand{\M}{\mathfrak M}

\newcommand{\N}{\mathbb N}
\renewcommand{\P}{\mathbb P}
\newcommand{\Pc}{\mathscr P}
\newcommand{\Q}{\mathbb Q}
\newcommand{\Qs}{\mathscr Q}
\newcommand{\R}{\mathbb R}
\newcommand{\Rs}{\mathscr R}
\newcommand{\Sc}{\mathscr S}

\newcommand{\V}{\mathbb V}
\newcommand{\Vs}{\mathscr V}
\newcommand{\Xs}{\mathscr X}
\newcommand{\Ys}{\mathscr Y}
\newcommand{\Z}{\mathbb Z}

\newcommand{\vnull}{\vec 0}
\newcommand{\vone}{\vec 1}
\newcommand{\vinfty}{\vec{\infty}}
\newcommand{\va}{\vec a}
\newcommand{\val}{\vec\alpha}
\newcommand{\vb}{\vec b}
\newcommand{\vc}{\vec c}
\newcommand{\vd}{\vec d}
\newcommand{\vDe}{\vec{\Delta}}
\newcommand{\ve}{\vec e}

\newcommand{\vh}{\vec h}
\newcommand{\vi}{\vec{\textit{\i}}}
\newcommand{\vj}{\vec{\textit{\j}}}
\newcommand{\vjprime}{\vj\,{}'}
\newcommand{\vk}{\vec{k}}
\newcommand{\vell}{\vec\ell}
\newcommand{\vmu}{\vec\mu}

\newcommand{\vn}{\vec n}
\newcommand{\vnu}{\vec\nu}
\newcommand{\vN}{\vec N}
\newcommand{\vl}{\vec l}
\newcommand{\vp}{\vec p}

\newcommand{\vq}{\vec q}
\newcommand{\vr}{\vec r}
\newcommand{\vs}{\vec s}
\newcommand{\vt}{\vec t}

\newcommand{\vu}{\vec u}
\newcommand{\vv}{\vec v}

\newcommand{\eps}{\varepsilon}

\newcommand{\Atau}{\dot A^{\vs,\tau}_{\vp,\vq,\pi}}
\newcommand{\AtauCube}{\dot A^{s,\tau,\mathrm{cube}}_{p,q,\pi}}

\newcommand{\atau}{\dot a^{\vs,\tau}_{\vp,\vq,\pi}}
\newcommand{\aatau}{\dot{\alpha}^{\vs,\tau}_{\vp,\vq,\pi}}

\newcommand{\antau}{\dot a^{\vnull,\tau}_{\vp,\vq,\pi}}
\newcommand{\ann}{\dot a^{\vnull,0}_{\vp,\vq,\pi}}
\newcommand{\asn}{\dot a^{\vs,0}_{\vp,\vq,\pi}}

\newcommand{\lift}[1]{\mathfrak{I}^{#1}}

\newcommand{\BMO}{\operatorname{BMO}}
\newcommand{\co}{\operatorname{co}}

\newcommand{\esssup}{\operatornamewithlimits{ess\,sup}}

\newcommand{\id}{\mathrm{id}}
\newcommand{\one}{\mathbf{1}}
\newcommand{\supp}{\operatorname{supp}}
\newcommand{\loc}{\mathrm{loc}}

\newcommand{\Pdm}{\C_+^{m\times m}}

\newcommand{\ave}[1]{\langle #1\rangle}
\newcommand{\aveL}{\textit{\L}}
\newcommand{\Avee}[2]{\langle#1\rangle_{#2}}
\newcommand{\Ave}[3]{\langle#1\rangle_{#2,#3}}
\newcommand{\Red}[3]{[#1]_{#2,#3}}

\newcommand{\floor}[1]{\lfloor #1 \rfloor}
\newcommand{\Floor}[1]{\lfloor\!\lfloor #1 \rfloor\!\rfloor}
\newcommand{\ceil}[1]{\lceil #1 \rceil}

\newcommand{\pair}[2]{\langle #1, #2 \rangle}

\newcommand{\Bpair}[2]{\Big\langle #1, #2 \Big\rangle}

\newcommand{\abs}[1]{|#1|}
\newcommand{\Babs}[1]{\Big|#1\Big|}
\newcommand{\Norm}[2]{\|#1\|_{#2}}
\newcommand{\bNorm}[2]{\big\|#1\big\|_{#2}}
\newcommand{\BNorm}[2]{\Big\|#1\Big\|_{#2}}
\newcommand{\bbNorm}[2]{\bigg\|#1\bigg\|_{#2}}

\newcommand{\perm}[2]{(#1\ #2)}

\newcommand{\Carl}{\mathrm{Carl}}
\newcommand{\Cubes}{\mathrm{Cubes}}

\newcommand{\Open}{\mathrm{Open}}

\newcommand{\Rect}{\mathrm{Rect}}
\newcommand{\RHI}{\mathrm{RHI}}

\newcommand{\AC}{\A^{\Carl}}

\newcommand{\pdef}[1]{\F^{#1\times #1}_+}

\begin{document}

\title[Matrix-weighted multi-parameter function spaces]{Real-Variable Theory
of Matrix-Weighted Multi-Parameter Besov--Triebel--Lizorkin-Type Spaces}

\author[F. Bu]{Fan Bu}
\address{Laboratory of Mathematics and Complex Systems (Ministry of Education of China),
School of Mathematical Sciences, Beijing Normal University, Beijing 100875,
The People's Republic of China}
\curraddr{}
\email{fanbu@mail.bnu.edu.cn}

\author[Y. Chen]{Yiqun Chen}
\address{Laboratory of Mathematics and Complex Systems (Ministry of Education of China),
School of Mathematical Sciences, Beijing Normal University, Beijing 100875, The People's Republic of China}
\curraddr{}
\email{yiqunchen@mail.bnu.edu.cn}

\author[T. P. Hyt\"onen]{Tuomas P.\ Hyt\"onen (*}
\address{Department of Mathematics and Systems Analysis, Aalto University, P.O. Box 11100 (Otakaari~1), \mbox{FI-00076} Aalto, Finland}
\curraddr{}
\email{tuomas.hytonen@aalto.fi}

\thanks{*) Corresponding author}

\author[D. Yang]{Dachun Yang}
\address{Laboratory of Mathematics and Complex Systems (Ministry of Education of China),
School of Mathematical Sciences, Institute for Advanced Study, Beijing Normal University, Beijing 100875, The People's Republic of China}
\curraddr{}
\email{dcyang@bnu.edu.cn}

\author[W. Yuan]{Wen Yuan}
\address{Laboratory of Mathematics and Complex Systems (Ministry of Education of China),
School of Mathematical Sciences, Beijing Normal University, Beijing 100875, The People's Republic of China}
\curraddr{}
\email{wenyuan@bnu.edu.cn}

\thanks{This project is partially supported by the National
Natural Science Foundation of China (Grant Nos. 12431006 and 12371093),
the Beijing Natural Science Foundation (Grant No. 1262011),
the Fundamental Research Funds for the Central Universities (Grant No. 2253200028),
the Research Council of Finland (Grant Nos. 346314, 364208, and 371637), and
the China Postdoctoral Science Foundation (No. 2025M773047).}

\date{\today.}


% 46E35(1973–now)Sobolev spaces and other spaces of ``smooth'' functions, embedding theorems, trace theorems
% 42B25(1980–now)Maximal functions, Littlewood-Paley theory
% 42B20(1980–now)Singular and oscillatory integrals (Calderón-Zygmund, etc.)
% 42C40(2000–now)Nontrigonometric harmonic analysis involving wavelets and other special systems
% 47A56(1980–now)Functions whose values are linear operators (operator- and matrix-valued functions, etc., including analytic and meromorphic ones)
% 35S05(1973–now)Pseudodifferential operators as generalizations of partial differential operators
% 42B35(2000–now)Function spaces arising in harmonic analysis


\begin{abstract}
We develop a comprehensive theory for a general class of multi-parameter function spaces of Besov--Triebel--Lizorkin type, equipped with a matrix-valued weight function. We prove the equivalence of different quasi-norms, the identification of function and sequence spaces via the $\varphi$-transform of Frazier and Jawerth [J. Funct. Anal. 1990], the boundedness of almost diagonal operators and multi-parameter singular integrals under minimal assumptions, molecular and wavelet characterisations, and Sobolev-type embedding theorems. We identify matrix-weighted $L^p$ spaces, Sobolev spaces of mixed smoothness, and multi-parameter BMO spaces of Chang--Fefferman type as examples of our general scale of spaces. Thus, our result on the boundedness of multi-parameter singular integrals on these spaces is seen as a far-reaching extension, with a different method, of a recent theorem of Domelevo et al. [J. Math. Anal. Appl. 2024] on matrix-weighted $L^p$ spaces.

For the needs of this theory, we develop several tools of independent interest in the broader contexts of matrix weights and multi-parameter harmonic analysis. In several cases, previous results were restricted to integrability exponents $p\in(1,\infty)$, while the theory of Besov--Triebel--Lizorkin-type spaces naturally involves the full range $p\in(0,\infty)$. The inherent non-convexity in the range $p\in(0,1)$ requires new methods. We extend the definition of multi-parameter $A_p$ matrix weights to $p\in(0,1]$ and establish their basic properties, culminating in the $L^p$-boundedness of a matrix-weighted strong maximal operator of Christ--Goldberg type (suitably rescaled when $p\in(0,1]$) for all $p\in(0,\infty)$. For $p\in(1,\infty)$, this is a recent result of Vuorinen [Adv. Math. 2024] based on convex-set-valued techniques of Bownik and Cruz-Uribe [arXiv 2022 or Math. Ann. (to appear)]; the lack of convexity requires us to develop a new approach that works for all $p\in(0,\infty)$.

We also need and prove a multi-parameter extension of Carleson embedding-type inequalities from Frazier and Roudenko [Math. Ann. 2021] but attributed by them to F.~Nazarov. As is typical for the multi-parameter theory, the condition for this embedding involves a supremum over all open sets. We show that this is necessary, but this requires a curious nontrivial elaboration of Carleson's classical counterexample [Mittag-Leffler Rep. 1974], which we present in an appendix.
\end{abstract}

\maketitle

%\vspace{0.2cm}

\tableofcontents

\section{Introduction}

The topic of this work lies at the intersection of several broad themes of classical and modern harmonic analysis on Euclidean spaces:
\begin{enumerate}[\rm(i)]
\item\label{Ap} Muckenhoupt's theory of $A_p$ weights and its
matrix-valued extension;

\item\label{funspa} theory of function spaces based on Littlewood--Paley theory;

\item\label{multi} multi-parameter harmonic analysis.
\end{enumerate}
Each of these topics has a long and distinguished history in its own right, and we will concentrate, in this introduction, on works that deal with at least two of the said themes simultaneously.

\subsubsection*{Theory of function spaces with matrix weights: \eqref{Ap} and \eqref{funspa}}
A more detailed account of the history of these two topics and their interaction has been given in the introduction of \cite{BHYY:1a}, and we shall hence be relatively brief here.
The correct matrix-valued analogue of Muckenhoupt's $A_p$ class was identified by Nazarov, Treil, and Volberg \cite{NT96,TV97,Vol97}, and related matrix-weighted maximal function estimates were found by Christ and Goldberg \cite{CG01,Gold03}. The matrix $A_p$ condition acquired its equivalent modern formulation in the work of Roudenko \cite{Roud03}, and was extended from $p\in(1,\infty)$ to $p\in(0,1]$ by Frazier and Roudenko \cite{FR04}.

The theory of $L^p$ spaces with matrix-weights is by now quite advanced, including diverse aspects like
weighted extrapolation \cite{BCU,cs25},
sparse domination \cite{CDO18,dhl20,Lau25,NPTV17},
two-weight \cite{cim18}
and sharp one-weight inequalities \cite{DPTV,HPV19,llor23},
variational \cite{DLY21} and
multilinear estimates \cite{KN26},
compactness criteria \cite{LYZ23},
transference \cite{nie12},
and applications to partial differential equations \cite{CUMR16,ID24,IM19}.
While $p\in(1,\infty)$ in most cases, there are some results also for $0<p\leq 1$ \cite{nie25b}.

As for other functions spaces, the foundations for Besov spaces with matrix weights were already laid by Frazier and Roudenko \cite{FR04,FR08,Roud03,Roud04}, with more recent contributions including \cite{BYY}.
Restricted aspects of the theory of matrix-weighted Triebel--Lizorkin spaces already appeared in the early works of Nazarov, Treil, and Volberg \cite{NT96,Vol97} and were revisited by Isralowitz \cite{Is21} from the point-of-view of characterising the matrix-weighted $L^p$, but a systematic investigation of these spaces was only recently started by Frazier and Roudenko \cite{FR21}. Since then, the theory of matrix-weighted Besov and Triebel--Lizorkin spaces has been extended by many further works, including \cite{BX:nonreg,BX:Iran,BHYY:1a,BHYY:1b,BHYY:1c,BHYY:2b,BYYZ1,BYYZ2,lyy24,mx25,WYZ24,YYZ2}.
Matrix-weighted BMO \cite{IKP},
Hardy \cite{BCYY},
Campanato \cite{CYY},
Bourgain--Morrey \cite{bgx25,bx25},
and modulation spaces \cite{nie25a,wgx25}
have also been considered.
The recent framework of {\em directional Banach function spaces} \cite{Nie24} provides an abstraction of matrix-weighted $L^p$ spaces including their variable-exponent versions.

\subsubsection*{Multi-parameter theory with weights: \eqref{Ap} and \eqref{multi}}

For scalar weights, a starting point of their interaction with the multi-parameter theory would seem to be the work of R.~Fefferman and Stein \cite{FS82}, where the product $A_p$ weights were defined \cite[bottom of page 127]{FS82} and the boundedness of some multi-parameter singular integrals of convolution type on the corresponding weighted $L^p$ spaces established. This was then extended to general multi-parameter Calder\'on--Zygmund (or Journ\'e) operators by R.~Fefferman \cite{Fef87,Fef88}, with estimates in the smaller weight class $A_{p/2}$ \cite{Fef87} as a stepping stone for the full result \cite{Fef88}.

A more direct approach to this boundedness, via a representation of Journ\'e operators due to Martikainen \cite{Mar12}, was only recently found by Holmes, Petermichl, and Wick \cite{HPW18}, paving the way for its matrix-weighted extension by Domelevo et al. \cite{DKPS}, which laid the foundations for a systematic development of the multi-parameter theory of matrix weights.
(Some special cases of this theory has been treated earlier by Nielsen and Rasmussen \cite{nie13,nr18}. In particular, \cite{nr18} dealt with convolution type operators \`a la \cite{FS82} under more restrictive assumptions on the matrix weights and conjectured the extension of their result to all matrix $A_p$ weights.)

The boundedness of the matrix-weighted maximal operator of Christ and Goldberg \cite{CG01,Gold03} and the matrix-weighted Rubio de Francia extrapolation theorem of Bownik and Cruz-Uribe \cite{BCU} have been extended to the multi-parameter setting by Vuorinen \cite{Vuo24}, using the properties of a set-valued maximal operator established in \cite{BCU}. A technical issue that arises in the multi-parameter theory of matrix weights is the measurable dependence on the remaining variable(s) of so-called reducing operators taken with respect to the other variable(s); this is studied in detail in \cite{BCU,DKPS}.
Different versions of the matrix-weighted multi-parameter BMO space have been studied by Kakaroumpas and Soler i Gibert \cite{KS1,KS2}.

\subsubsection*{Multiparameter theory of function spaces: \eqref{funspa} and \eqref{multi}}

Here, we need to make a distinction between pure product-space theory and more complicated, ``entangled'' multi-parameter settings.
\begin{enumerate}[\rm(a)]
\item\label{pure} A theory of Besov and Triebel--Lizorkin spaces on product
domains was developed, and the boundedness of multi-parameter
singular integrals of convolution type on these spaces established
in \cite{DLM10,LuZhu13} in the unweighted and scalar $A_\infty$-weighted
settings, respectively. In \cite{CLL16,Xu22}, the boundedness of Fourier
multipliers (scalar $A_\infty$-weighted) and pseudo-differential operators
(unweighted) were obtained, respectively. The dual spaces of these spaces
(scalar $A_\infty$-weighted) have been identified in \cite{DCN19,DingZhu17}
for complementary ranges of the integrability parameter $p$.
The $\varphi$-transform characterization (unweighted)
in terms of related sequence spaces, the boundedness of
almost diagonal operators on the latter, and applications
to nonlinear wavelet approximation were obtained in \cite{GKP21}.

\item\label{mixed} Beyond pure product theory, Besov and Triebel--Lizorkin spaces adapted to different multi-parameter structures have been studied in the literature: a flag setting in \cite{DLM10,LiaoLiu14,Yang09}, a Zygmund dilation setting in \cite{LiaoLiu13}, a shearlet setting in \cite{Vera13}, and a mixture of Euclidean and parabolic homogeneities in \cite{Ding14,DingLu16,DLZ16,Tan18,Tan20}.

\item\label{vector} It is a well-known point-of-view to
interpret functions of variables $(x,y)$ on a product domain
as vector-valued functions of, say, $x$ only, taking values
in a space of functions of the other variable $y$.
In particular, this leads to the notion of
``Besov or Triebel--Lizorkin space-valued Besov
or Triebel--Lizorkin spaces'',
which have been studied for a long time
\cite{Berko85,DHP07,DenkKaip13,Linde21,Weid05}.
Such spaces arise as boundary traces of Sobolev spaces with
mixed integrability in space and time variables,
which are relevant for mixed order systems of partial differential equations.
A recurring theme in these studies is identifying (intersections of)
such vector-valued spaces with mixed-order one-parameter spaces,
where integrability is described by a vector of exponents $\vp=(p_1,\ldots,p_k)$.

\item A theory of mixed-norm Besov and Triebel--Lizorkin spaces has been developed in \cite{cgn17,cgn19a,cgn19b,gjn17,gn16,nie23}.
Building upon these developments, Cleanthous and Georgiadis \cite{cg20} introduced  mixed-norm $\alpha$-modulation spaces, which include mixed-norm Besov spaces as a special case.
More recently, Nielsen \cite{nie25a} introduced the matrix-weighted $\alpha$-modulation spaces.
\end{enumerate}

\subsection{Scope of this work}

Against the background just described, we set as the goal of this work the development of a matrix-weighted multi-parameter theory of function spaces addressing all themes \eqref{Ap}, \eqref{funspa}, and \eqref{multi}
at the same time. Given the diversity of different types of multi-parameter
theories as discussed above, attempting a comprehensive treatment in this
respect would be too much, and we concentrate here on the pure product theory on
\begin{equation}\label{Rvn}
\R^{\vn}:=\R^{n_1}\times\cdots\times\R^{n_k}
\end{equation}
with invariance with respect to the full family of independent dilations $(x_\nu)_{\nu=1}^k\mapsto(\delta_\nu x_\nu)_{\nu=1}^k$ in each variable.
Thus, the entangled dilation systems, like the flag, Zygmund, shearlet,
and mixed homogeneities, will be left outside the scope of this work.
On the other hand, while keeping the primary focus on the ``pure'' Besov
and the ``pure'' Triebel--Lizorkin spaces, we set up a framework that also includes
spaces of ``mixed'' type, as in the vector-valued point-of-view mentioned above.

To be a little more precise, recall that the classical Triebel--Lizorkin and Besov quasi-norms are essentially based on iterated $L^p\ell^q$ and $\ell^qL^p$ quasi-norms, respectively;
and these are the only permutations of the position space $L^p(\R^n)$
and the scale space $\ell^q(\Z)$ in the one-parameter situation.
However, on the product domain \eqref{Rvn}, we may in principle
have a different position space $L^{p_\nu}(\R^{n_\nu})$ and scale
space $\ell^{q_\nu}(\Z)$ for each component variable $x_{\nu}$.
Then the $2k$-sequence of spaces $L^{p_1},\ldots,L^{p_k},\ell^{q_1},
\ldots,\ell^{q_k}$ has $(2k)!$ different permutations. Among them, pure Besov
(or Besov-type) spaces
correspond to the $(k!)^2$ sequences with all $L^{p_\mu}$ on the right-hand side of
all $\ell^{q_\nu}$ (all these permutations
are called \emph{Besov-type permutations}), and pure Triebel--Lizorkin
(or Triebel--Lizorkin-type) spaces to the $(k!)^2$
sequences with all $L^{p_\mu}$ on the left of all $\ell^{q_\nu}$
(all these permutations are called \emph{Triebel--Lizorkin-type permutations}).
It is evident that these $2(k!)^2$ pure spaces are but a fraction of all $(2k)!$
possible mixed spaces when $k>1$. Already in the bi-parameter case, $2(2!)^2=8$
and $(2\cdot 2)!=24$.

The reader will agree that in order to handle such a variety of spaces, a reasonably efficient notation will be necessary. As already hinted by the prose above, a key characteristic of a function space that we consider will be the associated permutation $\pi$ of $[2k]:=\{1,\ldots,2k\}$.
With $X_\nu:=L^{p_\nu}$ and $X_{\nu+k}:=\ell^{q_\nu}$ for $\nu\in[k]$, we denote
\begin{equation}\label{Lvpvq}
(L^{\vp}\ell^{\vq})_\pi:=\prod_{\nu=1}^{2k} X_{\pi(\nu)}:=X_{\pi(1)}\Bigg(\prod_{\nu=2}^{2k} X_{\pi(\nu)}\Bigg),
\end{equation}
where the ordered ``product'' of $L^{p_\mu}$ and $\ell^{q_\nu}$ spaces is defined recursively by
\begin{equation*}
Y(Z):=\Big\{f=f(y,z)\text{ measurable}: \bNorm{y\mapsto \Norm{z\mapsto f(y,z)}{Z}}{Y} <\infty\Big\}.
\end{equation*}

A reader who feels, at this point, that this level of generality is an unnecessary complication, should consider the following: A large part of the existing theory of function spaces is written with separate statements and proofs for the Besov case and the Triebel--Lizorkin case, and so it takes roughly double the space to treat both situations. In our notation, these both pure types of spaces are conveniently covered at once, and the fact that we simultaneously cover a spectrum of other spaces is a nice bonus that we obtain almost for free.

Certainly, when it comes to the fine details of proving some of the more delicate results, we will occasionally need to pursue a case study, but even then we claim that our treatment, overall, will be more efficient than the traditional division into Besov and Triebel--Lizorkin cases, even if we only discussed these two.

\subsubsection*{Admissible permutations}

Recall that, already in the classical theory, while Besov spaces $\dot B^s_{p,q}$ have a unified definition for all $p,q\in(0,\infty]$, Triebel--Lizorkin spaces $\dot F^s_{p,q}$ are defined differently for finite and infinite $p$, essentially reflecting the fact that the ``right'' substitute for $L^\infty$ is the space $\BMO$ of functions of bounded mean oscillation, whose definition is qualitatively different from the $L^p$ spaces. Thus, the basic one-parameter theory (excluding the separate treatment of $\dot F^s_{\infty,q}$) allows all other $L^p\ell^q$ and $\ell^qL^p$ except $L^\infty\ell^q$ with $q\in(0,\infty)$.

This exclusion admits a rather natural generalisation in our setting of $(L^{\vp}\ell^{\vq})_\pi$ spaces. We say that a permutation $\pi$ is {\em admissible} for $(\vp,\vq)$ if the following holds: after the leftmost $L^{p_\mu}$ space in the ``product'' \eqref{Lvpvq}, if this space itself or any space on its right-hand side (either $L^{p_\nu}$ or $\ell^{q_\nu}$) has exponent $\infty$, then all other spaces further to the right must also have exponent $\infty$. One easily verifies that $L^\infty\ell^q$ with $q\in(0,\infty)$ is the unique non-admissible case in the one-parameter theory.
A strong maximal operator is bounded on $(L^{\vp/a}\ell^{\vq/a})_\pi$
for some $a\in(0,\infty)$ if and only if the permutation $\pi$ is admissible.
Admissibility will be a standing assumption for our theory. This still allows
all combinations of finite exponents as well as some cases of infinite exponents,
and it is completely in line with the classical theory, as just discussed.

Depending on the particular results that we are proving, we will sometimes need
to impose other conditions on $(\pi,\vp,\vq)$ besides admissibility. Nevertheless, in all main results, we always cover at least the following cases (we alway allow a general $\vq\in(0,\infty]^k$):
\begin{enumerate}[\rm(i)]
  \item matrix-weighted spaces of pure Besov-type or pure Triebel--Lizorkin-type, where $\vp=p\cdot\vone$,
  \item\label{it:unweighted} unweighted spaces of any admissible type and any mixed exponent $\vp\in(0,\infty)^k$.
\end{enumerate}
While our principal focus is in the matrix-weighted theory, some of the results seem to be new even in the unweighted case, it is therefore interesting to note that the range of indices can be completely general in this case, as stated in \eqref{it:unweighted}.

\subsubsection*{The Morrey parameter $\tau$}
As we already mentioned, the classical theory requires a separate treatment of Triebel--Lizorkin spaces $\dot F^s_{\infty,q}$, with $q<\infty$, or in other words a replacement of the inadmissible $L^\infty\ell^q$, and this is necessary to cover the important space $\BMO=\dot F^0_{\infty,2}$ as part of the theory. However, instead of addressing this special case separately, a more recent approach is to equip our space with a fourth ``Morrey'' index $\tau\in[0,\infty)$, arriving at so-called Besov{\em-type} and Triebel--Lizorkin{\em-type} spaces $\dot B^{s,\tau}_{p,q}$ and $\dot F^{s,\tau}_{p,q}$, where three-index spaces (other than $\dot F^s_{\infty,q}$) correspond to $\tau=0$. Introduced in \cite{YY08,YY10,YSY10}, these spaces cover not only the classical $\dot F^s_{\infty,q}=\dot F^{1/p}_{p,q}$ (for any $p\in(0,\infty)$, by a version of the John--Nirenberg inequality), but also new examples outside the classical Besov--Triebel--Lizorkin scales, like so-called $Q$-spaces. The recent articles \cite{BHYY:1a,BHYY:1b,BHYY:1c} already incorporate this Morrey index into the matrix-weighted one-parameter
theory, and we adopt this framework also here.

In the one-parameter theory, the Morrey index brings a new supremum over (dyadic) cubes as part of the related quasi-norms, in natural agreement with the classical
definition of the norm in $\BMO$. However, it is also well known that the norm of the most natural
\cite{CF80,CF85} (although not the only relevant; see \cite[p.\,1708]{HPW18} or \cite[Definitions 1.3--1.4]{DKPS:bmo}) version of the multi-parameter $\BMO$ space involves a supremum not just over rectangles (as the multi-parameter extension of some other one-parameter concepts based on cubes) but over general open sets.
The $\BMO$ space defined in this way is the dual of the Hardy space.
Moreover, since the class of rectangles is not nested,
R. Fefferman \cite{Fef79}, building on Carleson's construction \cite{Car},
showed that in the bi-parameter setting the dual of the Hardy space
does not coincide with the so-called rectangular BMO space.

Viewing the Chang--Fefferman BMO as the primary example that a reasonable multi-parameter theory of function spaces should cover, we define our multi-parameter spaces using a similar supremum over open sets as a replacement of the cubes in the one-parameter theory.

\subsubsection*{The weight classes and other parameters}

Recall that a key part of our goals is developing a matrix-weighted theory
of the spaces that we study. While it is known that several related results can be obtained for the largest Muckenhoupt class of $A_\infty$ weights, both in the scalar-weighted multi-parameter theory \cite{CLL16,DCN19,DingZhu17,LuZhu13} and in the matrix-weighted one-parameter theory \cite{BHYY:2b,Vol97} (some results on Besov spaces \cite{Roud03} being also available under
different assumptions from ones of Muckenhoupt type), we here restrict our
considerations to the more well-studied matrix $A_p$ class. This is in line with the majority of recent literature on matrix-weighted spaces. Since results on $A_p$ weights tend to be tied with estimates in $L^u$ with $u$ equal (or at least close, as in \cite[Theorem 5.1]{Gold03}) to $p$, we are also led to assume isotropic integrability $\vp=p\cdot\vone$ for the majority of our results.
In the unweighted setting, this restriction can be removed.
Even after making this restriction, we find it convenient to use the notation $\vp$ and $(L^{\vp}\ell^{\vq})_\pi$ to emphasise the fact that we have $k$ independent $L^p$ quasi-norms that may assume any positions before, after, and between the $\ell^{q_\nu}$ quasi-norms in \eqref{Lvpvq}. On the other hand, we rarely need to restrict the nature of the vector $\vq$ (other than the positions of possible infinities for admissibility), which is an increased generality compared to most of the
related multi-parameter articles quoted. Allowing vectorial smoothness $\vs$ seems to be standard in the existing multi-parameter theory, and we follow this practice without any trouble.

\subsection{An overview of content and methods}

\subsubsection*{The first part of the paper,} consisting of Sections \ref{sec:weights} through \ref{sec:fourier},
is devoted to a general multi-parameter theory of matrix weights, maximal functions, and other inequalities; it is needed in the subsequent analysis of Besov--Triebel--Lizorkin-type spaces, but may have independent interest.

In Section \ref{sec:weights}, we set up a theory of multi-parameter matrix $\A_p$ weights for the full range of exponents $p\in(0,\infty)$, which is relevant in the theory of Besov--Triebel--Lizorkin-type spaces. For $p\in(1,\infty)$, such a multi-parameter theory was already developed in detail by Domelevo et al.\ \cite{DKPS}, but the case of $p\in(0,1]$ seems to be new in the multi-parameter setting, extending the one-parameter theory of Frazier and Roudenko \cite{FR04}. The exponent $p=1$ presents a genuine turning point in the nature of some results, like the useful characterisation of the $\A_p$ condition in terms of the boundedness of averaging operators.

Section \ref{sec:meas} addresses an issue of measurable selection of reducing operators that arises in the multi-parameter theory. Namely, it is a standard consequence of the John ellipsoid theorem that any norm $e\mapsto r(e)$ on $\F^m\in\{\R^m,\C^m\}$ is equivalent to $e\mapsto\abs{Ae}$ for some matrix $A$ known as the reducing operator. If $r(x,\cdot)$ is a family of such norms, measurable with respect to $x\in\Omega$, then one also wants to choose the reducing operators $A(x)$ in a measurable way. For proper norms, this is known from the works of Bownik and Cruz-Uribe \cite{BCU} and Domelevo et al.\ \cite{DKPS}. However, to deal with the full range of exponents $p\in(0,\infty)$ again, we need an extension to $p$-seminorms, which requires modifications of the arguments in the lack of convexity, which plays a role in the existing results of \cite{BCU,DKPS}.

In Section \ref{sec:max and cal}, we first identify the precise conditions for a vector-valued strong maximal inequality in the general iterated $(L^{\vp}\ell^{\vq})_\pi$ spaces, leading to the notion of admissible permutations $\pi$ already mentioned. More substantially, we obtain a multi-parameter extension of certain Carleson embedding-type inequalities of the form
\begin{equation}\label{ce}
\bNorm{\{\gamma_{\vj}f_{j}\}_{\vj\in\Z^k}}{(L^{\vp}\ell^{\vq})_\pi}
\leq C \bNorm{\{ f_{\vj}\}_{\vj\in\Z^k}}{(L^{\vp}\ell^{\vq})_\pi}
\end{equation}
under suitable conditions on the multipliers $\gamma_j$,
when each function $f_{\vj}$ is piecewise constant on dyadic rectangles $R\in\mathscr D_{\vj}$.
Such estimates were introduced by Frazier and Roudenko \cite{FR21} (and attributed by them to F.~Nazarov) and played an instrumental role in the one-paramater theory of matrix-weighted Triebel--Lizorkin spaces in \cite{FR21} and their extensions in \cite{BHYY:1a,BHYY:2b}. Our multi-parameter extension presents several nontrivial novelties of independent interest:
\begin{enumerate}[(i)]
  \item The multi-parameter Carleson condition has the usual complication of involving arbitrary open sets instead of cubes as in \cite{FR21}. The proof of its sufficiency for the embedding that we need depends on a new interpolation argument with respect to $p$ instead of $q$ as in \cite{FR21}.
  \item We also prove the necessity of this condition. As expected, this necessity is related to Carleson's famous counterexample \cite{Car} but, curiously, seems to require a nontrivial elaboration of this example, which we provide in Appendix \ref{app:Carleson}. This counterexample as such is independent of matrix weights, and may have other relevance for multi-parameter harmonic analysis.
  \item Finally, we need to check that multipliers arising from matrix weights satisfy the relevant Carleson condition. We prove that this is equivalent to a matrix-weighted strong maximal inequality. This latter is recently due to Vuorinen \cite{Vuo24} for $p\in(1,\infty)$, but it is a new result requiring a  new formulation and a new approach for $p\in(0,1]$. We present a unified method that works for all $p\in(0,\infty)$, thus also providing a new proof for $p\in(1,\infty)$. The boundedness of the strong maximal operator is a fundamental result in real analysis, and likewise this matrix-weighted extension should have broader relevance for matrix-weighted harmonic analysis beyond the theory of Besov--Triebel--Lizorkin spaces.
\end{enumerate}

The short Section \ref{sec:fourier} is independent of the previous ones and collects some basic lemmas of multi-parameter Fourier analysis, like Calder\'on reproducing formulas and decay estimates for pairings of test functions, that will be repeatedly used in the remainder of the paper.

\subsubsection*{The second part of the paper,} consisting of the remaining Sections \ref{BF type spaces} through \ref{embedding}, undertakes a thorough study of the matrix-weighted multi-parameter Besov--Triebel--Lizorkin-type spaces, making essential use of the tools developed in the first part.

In Section \ref{BF type spaces}, we introduce the main object of the paper, the matrix-weighted multi-parameter Besov--Triebel--Lizorkin-type spaces $\Atau(V)$, and establish their basic properties:
\begin{enumerate}[\rm(i)]
  \item Theorem \ref{3 norms thm} on the equivalence of different quasi-norms (in particular, ones where the weight function $V(x)$ acts pointwise at each $x$, and ``averaging'' versions, where the weight acts via its associated sequence of reducing operators on each dyadic rectangle), and
  \item the $\varphi$-transform characterisation (Theorem \ref{phi}), which provides a correspondence between the function spaces $\Atau(V)$ and their analogous discrete sequence spaces $\atau(V)$.
\end{enumerate}
With the new multi-parameter tools from Section \ref{sec:max and cal} in place of their one-parameter counterparts, we could in principle follow the one-parameter theory of \cite{BHYY:1a} step-by-step to achieve this extension; however, we take the opportunity to streamline the presentation, not only by covering Besov and Triebel--Lizorkin spaces simultaneously within our framework of permuted mixed-norm spaces $(L^{\vp}\ell^{\vq})_\pi$, but also by treating function and sequence spaces in a more unified way. In fact, our Section \ref{BF type spaces} on pp.\ \pageref{BF type spaces}--\pageref{sec AD}, while essentially self-contained, is much shorter than the presentation of the corresponding one-parameter theory in \cite[Section 3, pp.\ 6140--6167]{BHYY:1a}.

Section \ref{sec AD} studies the boundedness of almost diagonal infinite matrices on the sequence spaces $\atau(V)$; thanks to the $\varphi$-transform characterisation, this is immediately relevant also for the analysis of operators on the function spaces $\Atau(V)$. The main result of this section is Theorem \ref{AD st}. For the usual Besov and Triebel--Lizorkin spaces with $\tau=0$, it provides a complete multi-parameter extension of the best available one-parameter results from \cite{BHYY:1b}. For general Besov--Triebel--Lizorkin-type spaces, our Theorem \ref{AD st} extends the one-parameter results of \cite{BHYY:1b} in the worst-case scenario among arbitrary $\A_p$ weights. While the one-parameter theory allows for a sharper statement by taking into account additional information about the $\A_p$-dimension of the weight, which we are not yet able to fully incorporate into the multi-parameter theory, we show that our results are essentially optimal, in many cases, among uniform estimates for all $\A_p$ weights.

In Section \ref{sec:molecules}, we characterise the space $\Atau(V)$ in terms of molecular decompositions (Theorem \ref{89}), and provide their first applications to the boundedness of pseudo-differential operators on $\Atau(V)$ and wavelet expansions of these spaces. With the multi-parameter tools of the previous sections at hand, these are obtained as natural extensions of their one-parameter analogues from \cite{BHYY:1c}, but additional care is needed to address the new types of mixed cancellation that arise when pairing two multi-parameter molecules: in the one-parameter case, one typically uses the smoothness of one function and vanishing moments of the other one; in the multi-parameter case, we need to couple the smoothness of $f$ with the vanishing moments of $g$ with respect to some variables, and the vanishing moments of $f$ with the smoothness of $g$ in the remaining ones.

In Section \ref{other}, we compare our general scale of $\Atau(V)$ spaces with other matrix-weighted function spaces. In particular, we show that matrix-weighted $L^p(V)$ spaces (Theorem \ref{Lp=F}), Sobolev spaces (Theorem \ref{thm:Sobolev}), and Chang--Fefferman-type multi-parameter BMO spaces (Theorem \ref{coin BMO}) all fit into the general framework of $\Atau(V)$ spaces. To cover the last important example, the more complicated definition in terms of general open sets is essential. We also compare this definition with a simpler version involving dyadic cubes or rectangles only, both in the one-parameter and in the multi-parameter situations.

In Section \ref{sec:CZO}, we prove a multi-parameter $T(1)$ theorem (Theorem \ref{T1 BF}), a general criterion for the boundedness of multi-parameter singular integral operators, on the $\Atau(V)$ spaces. Thanks to the identification of other spaces as special cases in Section \ref{other}, this recovers (Corollary \ref{T1Co}) the paraproduct-free case of the multi-parameter $T(1)$ theorem on $L^p(V)$ spaces, which was only recently obtained by Domelevo et al. \cite{DKPS} by different methods, but we simultaneously cover cases like Sobolev spaces, which are outside the scope of \cite{DKPS}. As in the first multi-parameter $T(1)$ theorem on $L^2$ by Journ\'e \cite{Journe}, our approach proceeds by induction on the number of parameters. To make this induction run smoothly, we also revisit and streamline the one-parameter results of \cite{BHYY:1c}. In particular, we formulate a condition where the usual off-diagonal kernel bounds and the weak boundedness property are essentially merged. The proof also crucially depends on the molecular characterisation from Section \ref{sec:molecules}.

The final Section \ref{embedding} presents a Sobolev-type embedding (Theorem \ref{thm-sobolev-B}) for the $\Atau(V)$ spaces, with applications to the matrix-weighted $L^p(V)$ and Sobolev spaces. These results are not difficult, but we note the failure of a na\"ive multi-parameter extension of an interpolation inequality that has been used in some related one-parameter results; this necessitates another approach.

\subsection{Notation}\label{Symbols}
We list some notation here which are used in this article.

\begin{enumerate}[]

\item $\F\in\{\R,\C\}$ (Several results will be true in both cases);

\item $|x|:=(\sum_{i=1}^{n}|x_i|^2)^{\frac12}$, where $x:=(x_1,\ldots,x_n)\in\F^n$;

\item $|\vj|:=\sum_{i=1}^k|j_i|$, where $\vj:=(j_1,\ldots,j_k)\in\Z^k$;

\item $|E|:=$ the Lebesgue measure, where $E$ is a measurable set;

\item $\N:=\{0,1,2,\ldots\}$ (set of all nonnegative integers);

\item $\mathbb Z_+:=\{1,2,\ldots\}$ (set of all positive integers);

\item $[k]:=\{1,\ldots,k\}$, where $k\in\Z_+$;

\item $\va:=(a_1,\ldots,a_k)\in\F^k$ (We mostly reserve the arrow for vectors of length $k$);

\item $\vnull:=(0,\ldots,0)$, $\vone:=(1,\ldots,1)$,
$\vinfty:=(\infty,\ldots,\infty)$ (each of length $k$);

\item $2^{\vec a}:=(2^{a_1},\ldots,2^{a_k})$, where $\vec a:=(a_1,\ldots,a_k)\in \R^k$;

\item $\vec n:=(n_1,\ldots,n_k)\in \mathbb Z_+^k$, $\abs{\vn}:=n_1+\cdots+n_k$;

\item $\R^{\vn}:=\R^{n_1}\times\cdots\times\R^{n_k}$;

\item $x=(x_1,\ldots,x_k)\in\R^{\vn}$, where $x_i\in \R^{n_i}$ for each $i\in[k]$;

\item $\mathbf 0:=(0,\ldots,0)$, $\one :=(1,\ldots,1)$ in $\R^n$ or $\R^{n_i}$;

\item $\dot\R^{\vn}:=(\R^{n_1}\setminus\{\mathbf 0\})\times\cdots
\times(\R^{n_k}\setminus\{\mathbf 0\})$;

\item $2^{\vec a}x:=(2^{a_1}x_1,\ldots,2^{a_k}x_k)\in\R^{\vn}$, where $\vec a:=(a_1,\ldots,a_k)\in \R^k$
and $x:=(x_1,\ldots,x_k)\in\R^{\vn}$;

\item $\Open(\R^{n}):=$ the collection of all open subsets of $\R^n$;

\item $\Cubes(\R^n):=\{Q=x+[0,\ell(Q))^n: x\in\R^n,\ell(Q)\in(0,\infty)\}$
(set of all cubes in $\R^n$);

\item $\ell(Q):=$ edge-length of $Q\in\Cubes(\R^n)$;

\item $a_Q:= x$ is the ``lower left corner'' of $Q:=x+[0,\ell(Q))^n\in\Cubes(\R^n)$;

\item $b_Q:= x+\ell(Q)\one$ is the ``upper right corner'' of $Q:=x+[0,\ell(Q))^n\in\Cubes(\R^n)$;

\item $c_Q:= x+\frac12\ell(Q)\one$ is the centre points of $Q:=x+[0,\ell(Q))^n\in\Cubes(\R^n)$;

\item $\D(\R^n):=\{Q=2^{-j}(m+[0,1)^n):j\in\Z,m\in\Z^n\}$ (set of all dyadic cubes in $\R^n$);

\item $\D_j(\R^n):=\{Q=2^{-j}(m+[0,1)^n):m\in\Z^n\}$
(dyadic cubes of fixed edge-length $2^{-j}$);

\item $\Rect(\R^{\vn}):=\{R=Q_1\times\cdots\times Q_k:Q_i\in\Cubes(\R^{n_i})\text{ for }i\in[k]\}$ (rectangles in $\R^{\vn}$);

\item $\vell(R):=(\ell(R_1),\ldots,\ell(R_k))$,
where $R:=R_1\times\cdots\times R_k\in\Rect(\R^{\vn})$;

\item $a_R:=(a_{Q_1},\ldots,a_{Q_k})$ is the
``lower left corner'' of $R:=Q_1\times\cdots\times Q_k\in\Rect(\R^{\vn})$;

\item $b_R:=(b_{Q_1},\ldots,b_{Q_k})$ is the ``upper right corner''
of $R:=Q_1\times\cdots\times Q_k\in\Rect(\R^{\vn})$;

\item $c_R:=(c_{Q_1},\ldots,c_{Q_k})$ is the centre points of
$R:=Q_1\times\cdots\times Q_k\in\Rect(\R^{\vn})$;

\item $\D(\R^{\vn}):=\D(\R^{n_1})\times\cdots\times
\D(\R^{n_k})$ (set of all dyadic rectangles in $\R^{\vn}$);

\item $\D_{\vj}(\R^{\vn}):=\D_{j_1}(\R^{n_1})\times\cdots
\times\D_{j_k}(\R^{n_k})$, where $\vj:=(j_1,\ldots,j_k)\in\Z^k$
(those of fixed size $2^{-\vj}$);

\item $\D_{\vj}(\Omega):=\{R\in\D_{\vj}(\R^{\vn}):R\subset\Omega\}$
(dyadic rectangles of fixed size $2^{-\vj}$ inside a set $\Omega\subset\R^{\vn}$);

\item $\Omega_{\vj}:=\bigcup_{R\in\D_{\vj}(\Omega)}^{\phantom{R}}R\subset\Omega$,
where $\vj\in\Z^k$ and $\Omega\subset\R^{\vn}$;

\item $\aveL^p(E):=$ the space $L^p(E)$ with the normalised quasi-norm $\Norm{f}{\aveL^p(E)}:=\abs{E}^{-1/p}\Norm{\one_E f}{L^p}$;

\item $L^p(V):=$ the (matrix-)weighted $L^p$ space with quasi-norm
$\Norm{f}{L^p(V)}:=\Norm{Vf}{L^p}$;

\item $\aveL^p(E,V):=$ the (matrix-)weighted $L^p$ space with
quasi-norm $\Norm{f}{\aveL^p(E,V)}:=\Norm{Vf}{\aveL^p(E)}$;

\item $\mathbf{1}_E:=$ the characteristic function of $E$;

\item $\widetilde{\mathbf{1}}_E:=|E|^{-\frac12}\mathbf{1}_E$ (the $L^2$-normalised characteristic function);

\item $\ave{f}_E:=\fint_E f(x) \,dx$ (the average of $f$ on $E$);

\item $\Ave{f}{\vj}{p}:=\sum_{R\in\D_{\vj}(\R^{\vn})}^{\phantom{R}}\one_R\Norm{f}{\aveL^p(R)}$,
where $\vj\in\Z^k$;

\item $a\wedge b:=\min(a,b)$ and $a\vee b:=\max(a,b)$, where $a,b\in\R$;

\item $\min\va:=\min\{a_i:i\in[k]\}$, where $\va:=(a_1,\ldots,a_k)\in\R^k$;

\item $\min(\va,\vb):=\min\{\min\va,\min\vb\}$, where $\va,\vb\in\R^k$;

\item $p':=\begin{cases}\infty, & \text{if }p\in(0,1], \\ p/(p-1),
& \text{if }p\in(1,\infty), \\ 1, &\text{if }p=\infty; \end{cases}$
(an extended conjugate exponent);

\item $\vh\leq\vj$ means that $h_i\leq j_i$ for every $i\in[k]$,
where $\vh:=(h_1,\ldots,h_k),\vj:=(j_1,\ldots,j_k)\in\R^k$;

\item $\perm{a}{b}:=$ a permutation that swaps the elements $a$ and $b$ (the cycle notation);

\item $\F^{l\times m}:=$ $l\times m$-matrices with entries from $\F$;

\item $\F_+^{m\times m}:=$ positive definite matrices in $\F^{m\times m}$;

\item $|A|:= \sup_{z\in\F^m} |Az|$, where $A\in \mathbb F^{l\times m}$;

\item $A^*:=$ the conjugate transpose of $A\in \mathbb F^{m\times m}$ (where the conjugation is redundant if $\F=\R$);

\item $x^\gamma := x_1^{\gamma_1} \cdots x_n^{\gamma_n}$, where
$x := (x_1, \ldots, x_n) \in \mathbb{R}^n$ and
$\gamma := (\gamma_1, \ldots, \gamma_n) \in [0,\infty)^n$;

\item $\partial^\gamma_x := (\frac{\partial}{\partial x_1})^{\gamma_1}
\cdots (\frac{\partial}{\partial x_n})^{\gamma_n}$,
where
$x := (x_1, \ldots, x_n) \in \mathbb{R}^n$ and
$\gamma := (\gamma_1, \ldots, \gamma_n) \in \mathbb{Z}_+^n$;

\item $\mathcal M:=$ the Hardy--Littlewood maximal operator
$$\mathcal Mf(x):=\sup_{Q\in\Cubes(\R^{n}),\, x\in Q} \fint_Q |f(y)| \, dy;$$

\item $\mathcal M_{\vn}:=$ the strong Hardy--Littlewood maximal operator
$$\mathcal M_{\vn}f(x):=\sup_{R\in\Rect(\R^{\vn}),\, x\in R} \fint_R |f(y)| \, dy.$$
\end{enumerate}

Finally, the symbol $C$ denotes a positive constant
that is independent of the main parameters involved,
but may vary from line to line. The symbol $A\lesssim B$
means that $A\le CB$ for some positive constant
$C$, while $A\sim B$ means $A\lesssim B\lesssim A$.
Unless explicitly stated otherwise,
any dependence of such constants $C$ on weights is only through the weight constants.
In all proofs we consistently retain the notation
introduced in the statement being proved.

\section{Matrix weights}\label{sec:weights}

In this section, we set up a theory of multi-parameter matrix $\A_p(\mathbb R^{\vec n})$ weights to serve the needs of the rest of this work. The foundations for such a theory were laid down by Domelevo et al. \cite{DKPS}, who explicitly discussed the case of bi-parameter matrix weights. While the extension from two to several parameters is relatively routine, a key novelty of our treatment is incorporating the full scale of $\A_p$ weights with $p\in(0,\infty)$, where \cite{DKPS} dealt with $p\in(1,\infty)$ only. We recall that one-parameter $\A_p$ weights with $p\in(0,1]$ were already introduced by Frazier and Roudenko \cite{FR04}, but it seems that their multi-parameter version is here studied for the first time. This is highly relevant for building a weighted theory of Besov--Triebel--Lizorkin type space in the full range of the integrability parameter $p\in(0,\infty)$, which is desirable for a complete theory of these spaces.

Especially when it comes to $p\in(0,1)$, a new difficulty is the following: When $p\in[1,\infty)$, the seemingly non-linear $\A_p$ condition is conveniently characterised in terms of the $L^p(V)$ boundedness of the linear averaging operators $f\mapsto\one_Q\fint_Q f$ when $p\in[1,\infty)$ (see \cite[Lemma 3.1]{DKPS} for $p\in(1,\infty)$ or Lemma \ref{ave ops} below for $p\in[1,\infty))$, and this linearity is key to a useful characterisation of the multi-parameter $\A_p$ condition by means of iterating the one-parameter $\A_p$ conditions in each coordinate (see \cite[Lemmas 3.6 and 3.7]{DKPS}, or Lemma \ref{prod vs coord} below).

For $p<1$, these averages are not even well defined on $L^p(V)$, and we need to look for a substitute. The surprisingly simple solution turns out to be to replace the averages $f\mapsto\one_Q\fint_Q f$ by the {\em extension operators} $f\mapsto\one_Q\otimes f$, where the right-hand side is the function of two variables $(\one_Q\otimes f)(x,y)=\one_Q(x)f(y)$, and the characterising condition is the boundedness of these operators from $L^p(V)$ to $L^p(V\otimes\one)$. Thus, we are led to consider ``decoupling'' estimates of the form
\begin{equation*}
   \fint_Q\fint_Q\abs{V(x)f(y)}^p dx\,dy
   \lesssim\fint_Q\abs{V(x)f(x)}^p dx,
\end{equation*}
comparing the action of the weight $V$ on the function $f$ with independent and equal variables, respectively; see Lemma \ref{prod vs coord}. Reformulated in terms of the boundedness of suitable linear operators as in the case $p\in(1,\infty)$, this allows us to extend the iterative characterisation of multi-parameter matrix $\A_p$ weights to the full range $p\in(0,\infty)$. Conveniently, the same approach also gives an iterative characterisation of the classical reverse H\"older inequality (Lemma \ref{RHI prod vs coord}) and a multi-parameter extension of the recent notion of $\A_p$-dimension of matrix weights from \cite{BHYY:1a} and its key application to the comparison of reducing operators over different rectangles (Lemma \ref{sharp coef}).

The fact that these multi-parameter extensions can be derived as relatively quick corollaries of their one-parameter counterparts, without the need to rebuild the entire theory from scratch, illustrates the power of the iterative description of multi-parameter weights and, at the bottom, the characterisation of the $\A_p$ condition in terms of the boundedness of suitable linear operators as sketched above. The fact that this characterisation splits into two different forms at the index $p=1$ has a curious consequence: a certain interpolation property of two weights $V_{p_i}\in \A_{p_i}$ holds whenever both $p_0,p_1\in(0,1]$ or both $p_0,p_1\in[1,\infty)$, and {\em only} in these cases (see Proposition \ref{Ap interp}). In particular, we demonstrate by example that the interpolation property fails for $p_0,p_1$ on different sides of $1$, showing that the splitting of the theory into two ranges, which require slightly different considerations, is not an artefact of our choices in setting up the framework, but a necessity of the mathematical reality.

The structure of this section is organized as follows.
In Subsection \ref{reducing}, we study
the reducing operator of multi-parameter matrix weights.
In Subsection \ref{multi Ap}, we introduce
multi-parameter matrix $\mathscr A_p(\mathbb R^{\vec n})$ weights,
and Subsection \ref{Ap dim} concerns their associated $\mathscr A_p(\mathbb R^{\vec n})$-dimensions.
In Subsection \ref{reverse}, we study the reverse H\"older inequality of matrix weights.
Finally, in Subsection \ref{sec interp}, we establish an interpolation theorem for matrix weights.

\subsection{Reducing operators} \label{reducing}
Recall that we denote by $\F$ a generic scalar field, which may be $\R$ or $\C$.
The basic idea of {\em reducing operators} is that a general norm (or even a semi-quasi-norm) on $\F^m$ can be ``reduced'' to a norm of the form $x\mapsto\abs{Ax}$ for some matrix $A\in\F^{m\times m}$, up to inessential multiplicative constants, and this is useful in reducing considerations with various norms to matrix manipulations. 

To set the stage for the variety of norm-like expressions that we will encounter, we begin with:

\begin{definition}\label{seminorm}
Let $X$ be a vector space over the scalar field $\F$.

Let $r:X\to[0,\infty)$ satisfy $r(\lambda u)=|\lambda|r( u)$ for every  $\lambda\in\F$ and $u\in X$.

For $0<p\leq 1\leq K<\infty$, we consider the following inequalities:
\begin{equation}\label{eq:quasi-tri}
  r(u+v)\leq\begin{cases}
  r(u)+r(v), & \text{triangle inequality}, \\
  [r(u)^p+r(v)^p]^{\frac1p}, & p\text{-triangle inequality}, \\
  K[r(u)+r(v)], & \text{quasi-triangle inequality (with constant $K$)}.\end{cases}
\end{equation}
If one of these holds for all $u,v\in X$, we say that $r$ is
a \emph{seminorm}, a \emph{$p$-seminorm}, or a \emph{quasi-seminorm}
(with constant $K$), respectively.
If, in addition, $r(u)=0$ only if $u=0$, then we drop the word ``semi'' in each of these notions.

If $r:X\to[0,\infty]$, then we add the word ``extended'' in front of each of these notions.
\end{definition}

A matrix $A\in\F^{m\times m}$ is said to be {\em positive definite} if $(Ax,x)>0$
for every $x\in\F^m\setminus\{\mathbf0\}$, and
{\em positive semidefinite} if $(Ax,x)\geq0$ for every $x\in\mathbb F^m$.
For every $x:=(x_1,\ldots,x_m)\in\mathbb F^m$, we denote $$|x|:=(\sum_{i=1}^m |x_i|^2)^{\frac12}.$$

\begin{proposition}\label{qn red op}
If $\|\ \|$ is a quasi-norm (resp. quasi-seminorm) on $\F^m$,
then there exists a positive definite (resp. semidefinite) matrix
$A\in\F^{m\times m}$ such that, for every $x\in \F^m$,
$
\|x\|\sim |Ax|,
$
where the positive equivalence constants depend only on $m$ and the quasi-triangle constant of $\|\ \|$.
\end{proposition}

\begin{proof}
We first consider the case where $\|\ \|$ is a quasi-norm.
By the Aoki--Rolewicz theorem, as stated e.g. in \cite[Theorem 1.2]{Malig04},
every quasi-norm $\|\ \|$ is equivalent to a $p$-norm $\|\ \|^*$,
where both $p\in(0,1]$ and the equivalence constants depend only on
the quasi-triangle constant of $\|\ \|$.
By an argument given in \cite[page 1237]{FR04},
every $p$-norm $\|\ \|^*$ on $\F^n$ is equivalent to
a norm of the form $x\mapsto |Ax|$,
where the positive equivalence constants depend only on $p$ and $m$.
A combination of these equivalences proves the present proposition in this case.

If $\|\ \|$ is only a quasi-seminorm,
then $\mathcal N:=\{x\in\F^m:\ \|x\|=0\}$ is a subspace of $\F^m$
and $\|\ \|$ is a quasi-norm on its orthogonal complement $\mathcal N^{\perp}$.
The previous argument provides us with a positive definite
linear mapping $\widetilde A$ on $\mathcal N^{\perp}$ such that
$\|x\|\sim|\widetilde Ax|$ for every $x\in \mathcal N^{\perp}$.
If $P$ is the orthogonal projection of $\F^m$ onto $\mathcal N^{\perp}$,
then $A:=\widetilde A\circ P$ provides the desired semidefinite matrix.
This finishes the proof of Proposition \ref{qn red op}.
\end{proof}

For every measurable set $\Omega\subset\mathbb R^{\vec n}$,
let $\mathscr{M}(\Omega)$ denote the set of all
measurable functions on $\Omega$.

\begin{definition}\label{Debqfs}
Let $\Omega$ be a measurable set of $\mathbb R^{\vec n}$.
We say that $X(\Omega)\subset\mathscr{M}(\Omega)$ is a \emph{quasi-Banach lattice} if
\begin{enumerate}[\rm(i)]
  \item there is an extended quasi-norm $\Norm{\ }{X(\Omega)}$ on $\mathscr{M}(\Omega)$, whose restriction to $X(\Omega)$ is a quasi-norm;
  \item\label{it:lat-ineq} whenever $f,g\in \mathscr M(\Omega)$ satisfy $|g|\le |f|$ a.e.\ on $\Omega$, then $\|g\|_{X(\Omega)}\le\|f\|_{X(\Omega)}$;
  \item\label{it:lat-incl} in the situation of \eqref{it:lat-ineq}, if $f\in X(\Omega)$, then $g\in X(\Omega)$;
  \item $X(\Omega)$ equipped with $\Norm{\ }{X(\Omega)}$ is a quasi-Banach space.
\end{enumerate}
If $\|\ \|_{X(\Omega)}$ satisfies \eqref{eq:quasi-tri} with $K=1$, then $X(\Omega)$ is called a
\emph{Banach lattice}.
\end{definition}

If $X(\Omega)$ is a quasi-Banach lattice of $\F$-valued functions, we define
\begin{equation*}
X(\Omega;\F^{l\times m})
:=\Big\{F:\ \Omega\to\F^{l\times m}\text{ measurable}:\
\|F\|_{X(\Omega)}:=\big\||F(\cdot)|\big\|_{X(\Omega)}<\infty\Big\}.
\end{equation*}
For $F\in X(\Omega;\F^{l\times m})$ and $x\in\F^m$,
we have $F(\cdot)x:\ \Omega\to\F^l=\F^{l\times 1}$, and we can set
\begin{equation*}
\|x\|:= \big\| |F(\cdot)x| \big\|_{X(\Omega)}.
\end{equation*}
This $\|\ \|$ defines a quasi-seminorm on $\F^m$
and depends on $F$.
Proposition \ref{qn red op} provides us with a positive semidefinite matrix $A\in\F^{m\times m}$ such that,
for every $x\in \F^m$,
\begin{equation}\label{reducing op}
\|x\|\sim|Ax|,
\end{equation}
where the positive equivalence constants are as in Proposition \ref{qn red op}.
Although $A$ is not unique, we denote one such choice by $[F]_{X(\Omega)}$.

Let $p\in(0,\infty)$ and $Q$ be a cube in $\R^n$.
In the case where $X(\Omega)=\aveL^{p}(Q)$
and $F\in\aveL^{p}(Q;\pdef{m})$, \eqref{reducing op} reduces to:
for every $x\in \F^m$,
\begin{align}\label{1445}
\bigg[\fint_Q |F(y)x|^p dy\bigg]^{\frac1p}
\sim |[F]_{\aveL^{p}(Q)} x|,
\end{align}
where the positive equivalence constants depend only on $m$ and $p$.
In this case, $[F]_{\aveL^{p}(Q)}$ is precisely the reducing operator
of order $p$ for $F$ introduced by Volberg \cite[p.\,450]{Vol97}.
In the sense expressed by \eqref{1445}, the reducing operator is a linearisation of the norm induced by $F$, which enables us to deal with the matrix-valued function $F$ and the vector $x$ separately.

\begin{lemma}\label{red op M}
Let $X(\Omega)$ be a quasi-Banach lattice
(of scalar-valued functions).
Then, for every $F\in X(\Omega;\F^{l\times m})$ and $M\in\F^{m\times k}$,
\begin{equation*}
\| F(\cdot)M \|_{X(\Omega)}\sim |[F]_{X(\Omega)} M|,
\end{equation*}
where the positive equivalence constants depend only on $k$, $m$,
and the quasi-triangle constant of $\|\ \|_{X(\Omega)}$.
\end{lemma}

\begin{proof}
For all functions $\{\phi_i\}_{i=1}^k$ in $X$,
it follows from the properties of a
quasi-Banach lattice quasi-norm that, for each $j\in[k]$,
$$
\|\phi_j\|_{X(\Omega)}=\big\| |\phi_j| \big\|_{X(\Omega)} \leq\Bigg\|\sum_{i=1}^k|\phi_i|\Bigg\|_{X(\Omega)}
\lesssim \sum_{i=1}^k\|\phi_i\|_{X(\Omega)}
$$
and hence
$$
\Bigg\|\sum_{i=1}^k|\phi_i|\Bigg\|_{X(\Omega)}
\sim \sum_{i=1}^k\|\phi_i\|_{X(\Omega)}.
$$
If $\{e_i\}_{i=1}^k$ are the standard unit vectors in $\F^k$, then
\begin{equation*}
|M|\sim\sum_{i=1}^k |Me_i|.
\end{equation*}
Thus, using this bound with either $F(\cdot)M$ or $[F]_{X(\Omega)}M$ in place of $M$, and combining with the other bound above, we obtain
\begin{align*}
\|F(\cdot)M\|_{X(\Omega)}
&\sim\Bigg\|\sum_{i=1}^k|F(\cdot)Me_i| \Bigg\|_{X(\Omega)}
\sim\sum_{i=1}^k\big\| F(\cdot)Me_i \big\|_{X(\Omega)}\\
&\sim\sum_{i=1}^k| [F]_{X(\Omega)}Me_i|
\sim|[F]_{X(\Omega)}M|.
\end{align*}
This finishes the proof of Lemma \ref{red op M}.
\end{proof}

Before proving the following ``Fubini theorem'' for matrix functions, we first introduce some notation.
For every measurable
\begin{equation*}
  F:\Omega_1\times\Omega_2\to\F^{k\times m},
\end{equation*}
let
\begin{equation*}
  \|(x,y)\mapsto F(x,y)\|_{X(dx;Y(dy))}
  := \big\| \|y\mapsto F(\cdot,y)\|_{Y(dy)} \big\|_{X(\Omega_1)},
\end{equation*}
where, for every $x\in \Omega_1$,
\begin{equation*}
  \|y\mapsto F(x,y)\|_{Y(dy)}
  := \|F(x,\cdot)\|_{Y(\Omega_2)}.
\end{equation*}
Below, we will be concerned with $F(x,y)=A(x)B(y)$, where $A\in X(\Omega_1;\F^{k\times l})$ and $B\in Y(\Omega_2;\F^{l\times m})$.

In what follows, for all quasi-normed linear spaces $\mathscr X,\mathscr Y$
and for every linear or sublinear operator
$T: \mathscr X \to \mathscr Y$, let
$$
\|T\|_{\mathscr X\to \mathscr Y}
:= \sup_{\|f\|_{\mathscr X} = 1} \|Tf\|_{\mathscr Y} .
$$

\begin{lemma}\label{red Fub}
Let $X(\Omega_1)$ and $Y(\Omega_2)$ be quasi-Banach lattices.
Then the following two assertions hold.
\begin{enumerate}[\rm(i)]
\item\label{it:Fub-basic} For every $A\in X(\Omega_1;\F^{k\times l})$ and $B\in Y(\Omega_2;\F^{l\times m})$,
\begin{align*}
\Norm{A,B}{X(\Omega_1),Y(\Omega_2)}
&:=\|(x,y)\mapsto A(x)B(y)\|_{X(dx;Y(dy))}
\sim \Norm{[A]_{X(\Omega_1)} B}{Y(\Omega_2)} \\
&\phantom{:}\sim \Norm{A[B^*]_{Y(\Omega_2)}}{X(\Omega_1)}
\sim\Norm{B^*,A^*}{Y(\Omega_2),X(\Omega_1)}\\
&\phantom{:}=\|(x,y)\mapsto A(x)B(y)\|_{Y(dy;X(dx))}
\sim|[A]_{X(\Omega_1)}[B^*]_{Y(\Omega_2)}|,
\end{align*}
where the positive equivalence constants are independent of $A$ and $B$.

\item\label{it:Fub-dual} Suppose further that
$(Y(\Omega_2),Y(\Omega_2)')$ is a
\emph{dual pair} of Banach lattices in the following sense:
\begin{enumerate}[\rm(a)]
\item For all $f\in Y(\Omega_2)$ and $g\in Y(\Omega_2)'$, their pointwise product satisfies
\begin{equation}\label{bhi}
\|fg\|_{L^1(\Omega_2)}\leq\Norm{f}{Y(\Omega_2)}\Norm{g}{Y(\Omega_2)'}.
\end{equation}

\item For all $f\in Y(\Omega_2)$, its quasi-norm can be recovered by
\begin{equation}\label{129}
\Norm{f}{Y(\Omega_2)}=\sup\bigg\{\Babs{\int_{\Omega_2}fg}:\
g\in Y(\Omega_2)',\ \Norm{g}{Y(\Omega_2)'}\leq 1\bigg\}.
\end{equation}
\end{enumerate}
Then, for every $A\in X(\Omega_1;\F^{k\times l})$ and $B\in Y(\Omega_2;\F^{l\times m})$,
\begin{equation*}
\Norm{A,B}{X(\Omega_1),Y(\Omega_2)}
\sim\Big\|f\mapsto A\int_{\Omega_2}B f\Big\|_{Y(\Omega_2;\mathbb F^{m\times 1})'\to X(\Omega_1;\mathbb F^{k\times 1})},
\end{equation*}
where the positive equivalence constants are independent of $A$ and $B$.
\end{enumerate}
\end{lemma}

Note that, if $Y^*$ denotes the Banach space dual of $Y$,
our assumption that $(Y,Y')$ is a dual pair implies that $Y'$
can be identified with a subspace of $Y^*$, but {\em does not}
imply that $Y'$ should consist of all of $Y^*$.
This is relevant to cover the case of $Y=L^\infty(Q)$ and $Y'=\aveL^1(Q)$.
The following remark relates the notion of a dual pair to some other concepts used in the literature that go beyond our present needs.

\begin{remark}\label{it:ii} Let $Y(\Omega)$ be a quasi-Banach lattice.
\begin{enumerate}[\rm(i)]
\item\label{satisa} Recall that the \emph{associate space (K\"othe dual)}
$Y(\Omega)_{\textrm{KD}}'$ is defined
by setting
\begin{align*}
Y(\Omega)_{\textrm{KD}}':=\Big\{f\in\mathscr{M}(\Omega):\ \|f\|_{Y(\Omega)_{\textrm{KD}}'}
:=\sup_{\|g\|_{Y(\Omega)}=1}\|fg\|_{L^1(\Omega)}<\infty\Big\}.
\end{align*}
Using this definition, we find that inequality \eqref{bhi} holds in particular with
$Y(\Omega_2)'=Y(\Omega_2)_{\textrm{KD}}' $.

\item Recall that $Y(\Omega)$ is said to have the \emph{Fatou property}
if for every sequence of positive-valued functions
$\{f_j\}_{j\in\Z_+}\subset Y(\Omega)$ with $\liminf_{j\to\infty}\|f_j\|_{Y(\Omega)}<\infty$,
one has $\liminf_{j\to\infty}f_j\in Y(\Omega)$ and
\begin{equation*}
  \Big\|\liminf_{j\to\infty}f_j\Big\|_{Y(\Omega)}
\le\liminf_{j\to\infty}\|f_j\|_{Y(\Omega)}.
\end{equation*}
If $Y(\Omega)$ is a Banach lattice with the {\em saturation property}
(see \cite[p.\,249]{LN24} for its definition), then $Y(\Omega)$
has the Fatou property if and only if $(Y(\Omega)_{\textrm{KD}}')_{\textrm{KD}}'=Y(\Omega)$
(see \cite[Theorem 3.6]{LN24}).
When $\Omega=\Omega_2$, the last identity implies that \eqref{129} holds with $Y(\Omega_2)'=Y(\Omega_2)_{\textrm{KD}}'$.
Consequently, if $Y(\Omega_2)$ is a Banach lattice with
the Fatou and saturation properties, then $(Y(\Omega_2),Y(\Omega_2)_{\textrm{KD}}')$ is a
dual pair of Banach lattices in the sense of Lemma \ref{red Fub}.
The point of the assumptions of Lemma \ref{red Fub} is simply to state the minimal {\em ad hoc} condition that we need; as we just argued, this {\em ad hoc} condition can be verified via usual conditions appearing elsewhere in the literature.

\item We point out that the Fatou property of $Y(\Omega)$ and the associate space
enable one to escape the difficulties caused by
the lack of an explicit expression of the quasi-norm of $Y(\Omega)$ under consideration,
which have proved to play an important role in \cite{CWYZ,Ho24,HCY22,SHYY}.
In the same way, Lemma \ref{red Fub} plays a key role in the proofs of Lemmas
\ref{ave ops} and \ref{sharp coef}, which are the fundamental tools of this article.
\end{enumerate}
\end{remark}

\begin{proof}[Proof of Lemma \ref{red Fub}]
For all matrices of compatible size, we have
\begin{equation}\label{AB adj}
|AB|=|(AB)^*|=|B^*A^*|.
\end{equation}
From this identity and Lemma \ref{red op M}, for each point $x$, it follows that
\begin{equation*}
\Norm{A(x)B(\cdot)}{Y(\Omega_2)}
=\Norm{ B(\cdot)^*A(x)^*}{Y(\Omega_2)}
\sim\abs{[B^*]_{Y(\Omega_2)}A(x)^*}
=\abs{A(x)[B^*]_{Y(\Omega_2)}^*}
=\abs{A(x)[B^*]_{Y(\Omega_2)}},
\end{equation*}
using the self-adjointness of $[B^*]_{Y(\Omega_2)}$ in the last step.
Taking quasi-norms in $X(\Omega_1)$ on both sides,
\begin{equation*}
\Norm{(x,y)\mapsto A(x)B(y)}{X(dx,Y(dy))}
\sim \Norm{A(\cdot)[B^*]_{Y(\Omega_2)} }{X(\Omega_1)}
\sim \abs{[A]_{X(\Omega_1)}[B^*]_{Y(\Omega_2)}}.
\end{equation*}
By \eqref{AB adj} pointwise, we have
\begin{equation*}
\Norm{B^*,A^*}{Y(\Omega_2),X(\Omega_1)}=\Norm{(x,y)\mapsto B(y)^*A(x)^*}{Y(dy;X(dx))}
=\Norm{(x,y)\mapsto A(x)B(y)}{Y(dy;X(dx))}.
\end{equation*}
Here
\begin{equation*}
\Norm{ A(\cdot)B(y)}{X(\Omega_1)}
\sim\abs{[A]_{X(\Omega_1)} B(y)}
=\abs{B(y)^*[A]_{X(\Omega_1)}}.
\end{equation*}
Taking $Y(\Omega_2)$ quasi-norms of both sides gives
\begin{align*}
\Norm{(x,y)\mapsto A(x)B(y)}{Y(dy,X(dx))}
&\sim\Norm{[A]_{X(\Omega_1)} B}{Y(\Omega_2)}
=\Norm{B^*[A]_{X(\Omega_1)}}{Y(\Omega_2)} \\
&\sim\abs{[B^*]_{Y(\Omega_2)}[A]_{X(\Omega_1)}}
=\abs{[A]_{X(\Omega_1)}[B^*]_{Y(\Omega_2)}}.
\end{align*}
Combining the several identities above, we complete the proof of \eqref{it:Fub-basic}.

Next, we show \eqref{it:Fub-dual}.
For every $f\in Y(\Omega_2;\mathbb C^{m\times 1})'$,
\begin{equation}\label{eq:AintB}
\BNorm{ A\int_{\Omega_2} B f}{X(\Omega_1)}
\sim\Babs{[A]_{X(\Omega_1)}\int_{\Omega_2}B f}
=\Babs{\int_{\Omega_2}[A]_{X(\Omega_1)} B f}.
\end{equation}
Using the assumption on the duality of $(Y(\Omega_2),Y(\Omega_2)')$,
\begin{equation*}
\begin{split}
 \sup\{\eqref{eq:AintB}:\Norm{f}{Y(\Omega_2)'}\leq 1\}
&\sim\Norm{[A]_{X(\Omega_1)} B}{Y(\Omega_2)}=\Norm{B^*[A]_{X(\Omega_1)}}{Y(\Omega_2)} \\
&\sim\abs{[B^*]_{Y(\Omega_2)}[A]_{X(\Omega_1)}}
=\abs{[A]_{X(\Omega_1)}[B^*]_{Y(\Omega_2)}},
\end{split}
\end{equation*}
which is one of the equivalent quantities in part \eqref{it:Fub-basic}.
This finishes the proof of \eqref{it:Fub-dual} and hence Lemma~\ref{red Fub}.
\end{proof}

Of particular interest to us will be the following situation:
\begin{enumerate}[\rm(i)]
\item $X=\aveL^{p}(P)$ is the $L^{p}$ space on some $P\in\Rect(\R^{\vn})$, with $p\in(0,\infty)$ and the normalized quasi-norm
\begin{equation*}
\Norm{f}{\aveL^{p}(P)}:=\frac{\Norm{\one_P f}{L^{p}}}{\Norm{\one_P}{L^{p}}},
\end{equation*}

\item $Y=\aveL^{p'}(R)$ is a similar space related to the conjugate exponent
\begin{equation*}
p':=\begin{cases}\infty, & \text{if }p\in(0,1], \\ p/(p-1), & \text{if }p\in(1,\infty),
\end{cases}
\end{equation*}

\item $B=A^{-1}$ is the pointwise inverse of $A\in L_{\loc}^{p}(\R^n;\pdef{m})$, where $\pdef{m}$ is the set of all positive definite $m\times m$ matrices
and $A\in L_{\loc}^{p}(\R^{\vn};\pdef{m})$ means $|A|^p\in L_{\loc}^1(\R^{\vn})$
(the set of all locally integrable functions on $\R^{\vn}$).
\end{enumerate}

\begin{lemma}\label{ave ops}
Let $p\in(0,\infty)$ and $V\in L_{\loc}^p(\R^{\vec n};\F_+^{m\times m})$.
Then, for every $P,R\in\Rect(\R^{\vn})$,
the quantity $\Norm{V,V^{-1}}{\aveL^p(P),\aveL^{p'}(R)}$
has the following representations as norms of averaging operators:
\begin{equation} \label{1.29}
\Norm{V,V^{-1}}{\aveL^p(P),\aveL^{p'}(R)}
\sim\begin{cases}
\Norm{f\mapsto\one_P\ave{f}_R}{\aveL^{p}(R,V)\to\aveL^p(P,V)},
& \text{if } p\in[1,\infty), \\
\Norm{f\mapsto \one_P\otimes f}{\aveL^p(R,V)\to\aveL^p(P,V;\aveL^{\min\{1,p\}}(R))},
&  \text{if } p\in(0,\infty),
\end{cases}
\end{equation}
where $\otimes$ is the tensor product and
the positive equivalence constants are independent of $V$.
\end{lemma}

\begin{proof}
If $p\in[1,\infty)$, then $Y=\aveL^{p'}(R)$ is a Banach lattice with a dual pair $Y'=\aveL^p(R)$, with respect to the normalised integration $\ave{\ }_R:=\fint_R$ over $R$.
Hence it follows from Lemma \ref{red Fub} and a change of variables that
\begin{align*}
\Norm{V,V^{-1}}{\aveL^p(P),\aveL^{p'}(R)}
\sim \Norm{f\mapsto V\ave{V^{-1}f}_R}{\aveL^{p}(R)\to\aveL^p(P)}
=\Norm{f\mapsto \ave{f}_R}{\aveL^{p}(R,V)\to\aveL^p(P,V)}.
\end{align*}
For the other bound, let $s:=\min(1,p)$,
and note that $\frac 1s=\frac 1p+\frac 1{p'}$.
Then using Lemma \ref{red Fub} again, we obtain
\begin{align*}
\Norm{V,V^{-1}}{\aveL^p(P),\aveL^{p'}(R)}
&\sim\Norm{[V]_{\aveL^p(P)}V^{-1}}{\aveL^{p'}(R)} \\
&=\sup\Big\{\Norm{[V]_{\aveL^p(P)}V^{-1} f}{\aveL^{s}(R)}:\
\Norm{f}{\aveL^p(R)}\leq 1\Big\} \\
&\sim\sup\Big\{\Norm{V,V^{-1} f}{\aveL^p(P),\aveL^{s}(R)}:\
\Norm{f}{\aveL^p(R)}\leq 1\Big\} \\
&=\sup\Big\{\Norm{\one_P\otimes f}{\aveL^p(P,V;\aveL^s(R))}:\
\Norm{f}{\aveL^p(R,V)}\leq 1\Big\},
\end{align*}
and this is the desired estimate,
which completes the proof of Lemma \ref{ave ops}.
\end{proof}

The flexibility of the operator norm in Lemma \ref{ave ops}
plays a key role in the remainder of this article.
It enables us to establish a deep connection between
one-parameter and multi-parameter weights (see Lemma \ref{prod vs coord}),
thereby reducing questions on multi-parameter weights to the
corresponding one-parameter ones.
Moreover, using Lemma \ref{ave ops}, we show some conclusions about
the weighted geometric mean of two multi-parameter weights
(see Propositions \ref{Ap interp} and \ref{BCU2}),
which play an important role in the proof of
the extrapolation theorem (see \cite[p.\,49]{BCU})
and in part of the Jones factorization theorem (see Corollary \ref{Jones} below).

\begin{lemma}\label{inverses}
Let $p,u\in(0,\infty]$ and $P\in\Rect(\R^{\vn})$.
If $V\in\aveL^{p}(P;\F^{m\times m})$ and $V^{-1}\in\aveL^{u}(P;\F^{m\times m})$, then $[V]_{\aveL^{p}(P)}$ is invertible, and for all $M\in\F^{m\times k}$,
\begin{equation*}
\abs{[V]_{\aveL^{p}(P)}^{-1}M}\lesssim\abs{[V^{-*}]_{\aveL^{u}(P)}M},
\end{equation*}
where $V^{-*}:=(V^{-1})^*$.
\end{lemma}

\begin{proof}
By the definition of the matrix norm, it suffices to prove this with $e\in\F^m$ in place of $M$.
Let $s\in(0,\infty]$ be defined by setting
\begin{equation*}
\frac{1}{s}=\frac{1}{p}+\frac{1}{u},
\end{equation*}
so that we have access to H\"older's inequality:
\begin{equation*}
\abs{e}
=\Norm{V^{-1}Ve}{\aveL^{s}(P)}
\leq\Norm{V^{-1}}{\aveL^{u}(P)}\Norm{Ve}{\aveL^{p}(P)}
\lesssim\Norm{V^{-1}}{\aveL^{u}(P)}\abs{[V]_{\aveL^{p}(P)}e}.
\end{equation*}
Hence $[V]_{\aveL^{p}(P)}e =\mathbf 0$ implies $e=\mathbf 0$,
thus $[V]_{\aveL^{p}(P)}$ is injective and therefore invertible.
Then we can compute
\begin{align*}
|[V]_{\aveL^{p}(P)}^{-1}e|^2
&=\| V^{-*}e \cdot V[V]_{\aveL^{p}(P)}^{-2}e \|_{\aveL^{s}(P)}
\leq\Norm{V^{-*}e}{\aveL^{u}(P)}\Norm{V[V]_{\aveL^{p}(P)}^{-2}e}{\aveL^{p}(P)} \\
&\lesssim |[V^{-*}]_{\aveL^{u}(P)}e|\,
|[V]_{\aveL^{p}(P)}[V]_{\aveL^{p}(P)}^{-2}e|
=|[V^{-*}]_{\aveL^{u}(P)}e|\,
|[V]_{\aveL^{p}(P)}^{-1}e|.
\end{align*}
The claimed bound is trivial for $e=\mathbf 0$, and otherwise it follows after dividing the previous computation by $\abs{[V]_{\aveL^{p}(P)}^{-1}e}$.
\end{proof}

\subsection{Multi-parameter matrix weights} \label{multi Ap}

To define an $L^p$ space with a weight, two natural definitions suggest themselves: one might think of the weight as modifying either the underlying measure, or each function in the space, as in
\begin{equation*}
\bigg[\int |f(x)|^p w(x)\,dx\bigg]^{\frac1p}
\ \ \text{or}\ \
\bigg[\int |v(x)f(x)|^p\,dx\bigg]^{\frac1p}.
\end{equation*}
If the weights are scalar-valued, the two versions are obviously equivalent, via the correspondence $w=v^p$. When $f$ is vector-valued and the weight is matrix-valued, only the second version makes sense. Even then, it has been customary to use notation derived from the first version, viewing $w$ as the primary object and writing $v$ as $v=w^{\frac1p}$.
However, in an increasing number of studies, good reasons have been noted to reverse this philosophy, and we follow this trend here.
Expressed in terms of this formalism, the matrix Muckenhoupt condition for $V$ takes the following form.

\begin{definition}\label{def Ap}
Let $\Pc$ be a collection of pairs $(P,R)$ of rectangles $P,R\in\Rect(\R^{\vn})$, let $p\in(0,\infty)$, and let $V:\R^{\vn}\to\F^{m\times m}$ be (a.e.) positive-definite-valued. Then
\begin{equation*}
[V]_{\A_{p}(\Pc)}:=\sup_{(P,R)\in\Pc}\|V,V^{-1}\|_{\aveL^{p}(P),\aveL^{p'}(R)}.
\end{equation*}
We say that $V$ belongs to $\A_{p}(\Pc)$ and write $V\in\A_{p}(\Pc)$ if $[V]_{\A_{p}(\Pc)}<\infty$.

In the specific cases
\begin{equation*}
  \Pc\in\left.\begin{cases} \Cubes^2(\R^n) :=\{(Q,Q):Q\in\Cubes(\R^n)\}, \\
    \Rect^2(\R^{\vn}) :=\{(R,R):R\in\Rect(\R^{\vn})\}\end{cases}\right\},
\end{equation*}
we denote $\A_p(\R^n):=\A_p(\Cubes^2(\R^n))$ and $\A_p(\R^{\vn}):=\A_p(\Rect^2(\R^{\vn}))$.
\end{definition}

\begin{remark}\label{rem Ap}
We use the script-font $\A_p(\R^{\vn})$ to make a distinction from the usual $A_p(\R^{\vn})$; these two classes are closely related but not identical. The relation is very simply:
\begin{equation*}
  V\in\A_p(\R^{\vn})\quad\Longleftrightarrow\quad W:=V^p\in A_p(\R^{\vn}),
\end{equation*}
i.e., $W\in A_p(\R^{\vn})$ $\Longleftrightarrow$ $V:=W^{\frac1p}\in \A_p(\R^{\vn})$. Thus, the class $\A_p(\R^{\vn})$ is the natural one in our formalism, where $V$, instead of $W$, is the primary object of study. Note that this distinction of $\A_p(\R^{\vn})$ and $A_p(\R^{\vn})$ remains even when specialised to scalar-valued weights.
\end{remark}

As in the two key examples above, most of the time, the choice of $\Pc$ will consist of pairs $(P,P)$ of two equal sets, but other cases will also be relevant. In \cite{BHYY:1a}, the notion of $A_p$-dimensions was defined, effectively involving a condition like above for
\begin{equation*}
\Cubes^2_j(\R^n):=\{(Q,2^j Q):Q\in\Cubes(\R^n)\}.
\end{equation*}
Since rectangles have $k$ independent edge-lengths, we define
\begin{equation*}
\Rect^2_{\vj}(\R^{\vn}):=\{(R,2^{\vj}R):R\in\Rect(\R^{\vn})\},\qquad
2^{\vj}R:=2^{j_1}R_1\times\cdots\times 2^{j_k}R_k.
\end{equation*}
Thus, $\Cubes^2(\R^n)=\Cubes^2_0(\R^n)$ and $\Rect^2(\R^{\vn})=\Rect^2_{\vnull}(\R^{\vn})$.

The following lemma is a generalization of \cite[Lemmas 3.6 and 3.7]{DKPS}, where the case $p\in(1,\infty)$ and $k=2$ is contained.
This can be regarded as a coordinate-wise (Fubini-type)
characterisation of multi-parameter matrix weights. Here and throughout, we denote
\begin{equation*}
  \R^{\vn}\ominus\R^{n_i}
  :=\R^{n_1}\times\cdots\times\R^{n_{i-1}}\times
  \R^{n_{i+1}}\times\cdots\times\R^{n_{k}},\qquad i\in[k].
\end{equation*}

\begin{lemma}\label{prod vs coord}
Let $p\in(0,\infty)$, $\vj\in\mathbb Z^n$,
and $k\in\mathbb N\cap[2,\infty)$. Then
\begin{equation*}
\begin{split}
\max_{i\in[k]} &\esssup_{x'\in\R^{\vn}\ominus\R^{n_i}}  [x_i\mapsto V(x'\oplus x_i)]_{\A_{p}(\Cubes^2_{j_i}(\R^{n_i}))} \\
&\lesssim [V]_{\A_{p}(\Rect^2_{\vj}(\R^{\vn}))}
\lesssim\prod_{i=1}^k \esssup_{x'\in\R^{\vn}\ominus\R^{n_i}} [x_i\mapsto V(x'\oplus x_i)]_{\A_{p}(\Cubes^2_{j_i}(\R^{n_i}))},
\end{split}
\end{equation*}
where the implicit positive constants are independent of $V$.
\end{lemma}

\begin{proof}
For both estimates, we make use of Lemma \ref{ave ops}. Let $s:=\min(1,p)$.

\textit{The lower bound:} Fix an $i\in[k]$. In the definition of $\A_{p}(\Cubes^2_{j_i}(\R^{n_i}))$, it is enough to consider all pairs $(Q,2^{j_i}Q)$, where $Q$ belongs to a sufficiently rich countable subcollection $\Qs(\R^{n_i})\subset\Cubes(\R^{n_i})$, say all cubes with rational corners. By Lemma \ref{ave ops}, we want to estimate the norms of the corresponding operators
\begin{equation}\label{pvc1}
f\mapsto\one_Q\otimes f,\quad \aveL^p(2^{j_i}Q,V(z'\oplus\cdot))\to \aveL^p(Q,V(z'\oplus\cdot);\aveL^s(2^{j_i}Q))
\end{equation}
for all $Q\in\Qs(\R^{n_i})$ and almost all $z'\in\R^{\vn}\ominus\R^{n_i}$.
In order to estimate these operator norms, it is enough to consider
all $f$ in a suitable countable collection $\Fs\subset L^\infty(\R^{n_i};\F^m)$,
say of all rational linear combinations of $\one_Q$ with $Q\in\Qs(\R^{n_i})$.
Let us then consider the countably many functions
\begin{align*}
z'\mapsto &\Norm{(x_i,y_i)\mapsto V(z'\oplus x_i)f(y_i)}{\aveL^p(Q,dx_i;\aveL^s(2^{j_i}Q,dy_i))}^p, \\
z'\mapsto &\Norm{y_i\mapsto V(z'\oplus y_i)f(y_i)}{\aveL^p(2^{j_i}Q,dy_i)}^p,
\end{align*}
where $Q\in\Qs(\R^{n_i})$ and $f\in\Fs$. They are all in $L^1_{\loc}(\R^{\vn}\ominus\R^{n_i})$, and hence almost every $z'\in \R^{\vn}\ominus\R^{n_i}$ is a Lebesgue point of each of these functions. For any such $z'$ and any sequence of cubes $P\in\Cubes(\R^{\vn}\ominus\R^{n_i})$ shrinking to $z'$, it follows that
\begin{align*}
\Norm{\one_Q\otimes f}{\aveL^p(Q,V(z'\oplus\cdot);\aveL^s(2^{j_i}Q))}
&=\Norm{V(z'\oplus x_i)f(y_i)}{\aveL^p(Q,dx_i;\aveL^s(2^{j_i}Q,dy_i))} \\
&=\lim_{P\to z'} \Norm{V(x'\oplus x_i)f(y_i)}{\aveL^p(P\times Q,dx'dx_i;\aveL^s(2^{j_i}Q,dy_i))},
\end{align*}
moreover, with $\vjprime:=(j_h)_{h\in[k]\setminus\{i\}}$,
we can freely change $\aveL^s(2^{j_i}Q,dy_i)$ to
\begin{equation*}
  \aveL^s(2^{\vjprime}P\times 2^{j_i}Q,dy'dy_i)=\aveL^s(2^{\vj}(P\times Q)),
\end{equation*}
since the integrand is independent of the additional variable $y'\in 2^{\vjprime}P$.
Clearly
\begin{equation*}
  (P\times Q,2^{\vj}(P\times Q))\in\Rect^2_{\vj}(\R^{\vn}),
\end{equation*}
and hence
\begin{align*}
&\Norm{V(x'\oplus x_i)f(y_i)}{\aveL^p(P\times Q,dx'dx_i;\aveL^s(2^{j_i}Q,dy_i))}\\
&\quad= \Norm{V(x'\oplus x_i)f(y_i)}{\aveL^p(P\times Q,dx'dx_i;\aveL^s(2^{\vj}(P\times Q),dy'dy_i))} \\
&\quad\leq [V]_{\A_p(\Rect_{\vj}^2(\R^{\vn}))}
\Norm{V(y'\oplus y_i)f(y_i)}{\aveL^p(2^{\vj}(P\times Q),dy'dy_i)}.
\end{align*}
Noting that $2^{\vjprime}P$ also shrinks regularly to $z$ (although not necessarily cubes, these rectangles still have bounded eccentricity, since the dilation parameter $\vj$ is fixed), the other half of the Lebesgue point assumption implies that
\begin{align*}
\lim_{P\to z'}\Norm{V(y'\oplus y_i)f(y_i)}{\aveL^p(2^{\vj}(P\times Q),dy'dy_i)}
&=\lim_{P\to z'}\Norm{V(y'\oplus y_i)f(y_i)}{\aveL^p(2^{\vjprime}P,dy';\aveL^p(2^{j_i}Q,dy_i))} \\
&=\Norm{V(z'\oplus y_i)f(y_i)}{\aveL^p(2^{j_i}Q,dy_i)}.
\end{align*}
Thus, we obtained the boundedness of \eqref{pvc1}, and hence the lower bound claimed in the Lemma.

\textit{The upper bound:}
Let first $p\in(0,1]$. For all $P\in\Rect(\R^{\vn})$, we should now establish the boundedness of
\begin{equation}\label{pvc2}
f\mapsto\one_P\otimes f,\quad\aveL^p(2^{\vj}P,V)\to\aveL^p(P,V;\aveL^p(2^{\vj}P))
=\aveL^p(P\times 2^{\vj}P,V\otimes\one).
\end{equation}
Let $P=P'\times P_k$ and $\vj=(\vjprime,j_k)$. From Fubini's theorem, the uniform $\A_p$-property of weights $x_k\mapsto V(x'\oplus x_k)$, and the resulting uniform boundedness of the related operators in the $x_k$ direction, it follows that
\begin{align*}
&\Norm{\one_P\otimes f}{\aveL^p(P\times 2^{\vj}P,V\otimes\one)}
=\Norm{\one_P(x)f(y)}{\aveL^p(P\times 2^{\vj}P,V(x),dx dy)} \\
&\quad=\BNorm{ \one_{P'}(x') \Norm{\one_P(x_k)f(y',y_k)}{\aveL^p(P_k\times 2^{j_k}P_k,V(x'\oplus x_k)\otimes\one,dx_kdy_k)} }{\aveL^p(P'\times 2^{\vjprime}P',dx'dy')} \\
&\quad\lesssim
\BNorm{ \one_{P'}(x') [V(x'\oplus\cdot)]_{\A_p(\Cubes_{j_k}(\R^{n_k}))}\Norm{f(y',y_k)}{\aveL^p( 2^{j_k}P_k,V(x'\oplus y_k) dy_k)} }{\aveL^p(P'\times 2^{\vjprime}P',dx'dy')} \\
&\quad =\esssup_{x'\in\R^{\vn}\ominus\R^{n_k}}[V(x'\oplus\cdot)]_{\A_p(\Cubes_{j_k}(\R^{n_k}))}
\Norm{ \one_{P'}(x')f(y)}{\aveL^p(P'\times 2^{\vj} P,V(x',y_k),dx' dy)}.
\end{align*}
Since all exponents are equal to $p$, we can permute the order of the quasi-norms as we like, bringing the next $\aveL^p(P_i\times 2^{j_i}P_i)$ into the innermost position,
and using the $\A_p$-property of the $x_i$ direction to replace the quasi-norm over this space by a quasi-norm over just $\aveL^p(2^{j_i}P_i)$. Repeating this for each $i\in[k]$, we end up with
\begin{equation*}
\Norm{\one_P\otimes f}{\aveL^p(P\times 2^{\vj}P,V\otimes\one)}
\lesssim \prod_{i=1}^k \esssup_{x'\in\R^{\vn}\ominus\R^{n_i}}[V(x'\oplus\cdot)]_{\A_p(\Cubes_{j_i}(\R^{n_i}))}
\Norm{f}{\aveL^p(2^{\vj} P,V)},
\end{equation*}
which gives the required estimate for $[V]_{\A_p(\Rect_{\vj}^2(\R^{\vn}))}$ by the equivalence of the $\A_p$ constant and the norm of the operators \eqref{pvc2}.

For $p\in(1,\infty)$, the argument is similar, except that in place of \eqref{pvc2}, we consider the operators
\begin{equation*}
A_{P,\vj}:\ f\mapsto\one_P\ave{f}_{2^{\vj}P},
\quad\aveL^p(2^{\vj}P,V)\to\aveL^p(P,V).
\end{equation*}
Note that
\begin{equation*}
A_{P,\vj}=\prod_{i=1}^k A_{P_i,j_i}^{(i)}=A_{P',\vj}'\,A_{P_k,j_k}^{(k)}
\end{equation*}
where the operator $A_{P_i,j_i}^{(i)}$ acts in the $x_i$ direction, and thus these operators commute. Hence
\begin{align*}
\Norm{A_{P,\vj}f}{\aveL^p(P,V)}
&=\BNorm{\Norm{A_{P_k,j_k}^{(k)}(A_{P',\vjprime}'f)}{\aveL^p(P_k,V(x',\cdot))}}{\aveL^p(P',dx')} \\
&\lesssim \BNorm{ [V(x'\oplus\cdot)]_{\A_p(\Cubes_{j_k}(\R^{n_k}))}\Norm{A_{P',\vjprime}'f}
{\aveL^p(2^{j_k}P_k,V(x',\cdot))}}{\aveL^p(P',dx')} \\
&\lesssim\esssup_{x'\in\R^{\vn}\ominus\R^{n_k}}[V(x'\oplus\cdot)]_{\A_p(\Cubes_{j_k}
(\R^{n_k}))}\Norm{A_{P',\vjprime}'f}{\aveL^p(P'\times 2^{j_k}P_k,V)}.
\end{align*}
Again, we can permute the order of the quasi-norms. For each $i\in[k]$ in its turn, we bring the quasi-norm over $\aveL^p(P_i)$ as the innermost object, and use the $\A_p$ condition in the $x_i$ direction to remove the
operator $A_{P_i,j_i}^{(i)}$ while changing the said innermost $\aveL^p(P_i)$ quasi-norm to the quasi-norm of $\aveL^p(2^{j_i}P_i)$. In the end, we are left with
\begin{equation*}
\Norm{A_{P,\vj}f}{\aveL^p(P,V)}
\lesssim\prod_{i=1}^k \esssup_{x'\in\R^{\vn}\ominus\R^{n_i}}[V(x'\oplus\cdot)]_{\A_p(\Cubes_{j_i}(\R^{n_i}))}
\Norm{f}{\aveL^p(2^{\vj}P,V)},
\end{equation*}
which is equivalent to the claimed bound for $[V]_{\A_p(\Rect^2_{\vj}(\R^{\vn}))}$.
\end{proof}

\subsection{The $\A_p$-dimension of multi-parameter weights} \label{Ap dim}

In \cite{BHYY:1a}, the notion of $\A_p$-dimensions of weights was introduced,
which played an important role in the theory of matrix-weighted Besov--Triebel--Lizorkin-type spaces
in the one-parameter context \cite{BHYY:1a,BHYY:1b,BHYY:1c}.
We now wish to generalise this notion to the multi-parameter weights.
Rather than repeating the theory from scratch,
we aim at using the coordinate-wise (Fubini-type) characterization of
product $\A_p$-weights from Lemma \ref{prod vs coord},
so as to be able to build on the existing theory as much as possible.
Thus, we start by collecting some useful facts about matrix weights in the isotropic case.

\begin{lemma}\label{matrix 2 scalar}
Let $p\in(0,\infty)$ and suppose that $U\in\A_p(\R^n;\F^{m\times m})$. If
$M\in\F^{m\times l}\setminus\{\mathbf O\}$, then
\begin{equation*}
    [\abs{UM}^p]_{A_\infty(\R^n)}
   \lesssim[\abs{UM}^p]_{A_{\max(1,p)}(\R^n)}
   \leq (c[U]_{\A_p(\R^n)})^p,
\end{equation*}
where the implicit positive constant in the first bound depends only on $n$, and the positive constant $c$ in the second bound only on $m$ and $l$.
\end{lemma}

\begin{proof}
The first bound is a purely scalar result for the weight $w:=\abs{UM}^p$. It follows from
\begin{equation*}
  [w]_{A_\infty(\R^n)}\lesssim[w]_{A_\infty(\R^n)}^{\exp}:=
  \sup_{Q\in\Cubes(\R^n)}\ave{w}_Q\exp\ave{\log w^{-1}}_Q\leq[w]_{A_p(\R^n)},
\end{equation*}
where the second estimate is an application of Jensen's inequality, and the first one is \cite[Proposition 2.2]{HP13} (where the notation for the two $A_\infty$ constants is different).

We then consider the second bound. Let first $p\in(0,1]$. Then, for almost all $y\in Q$,
\begin{align*}
\fint_Q\abs{U(x)M}^p dx
&=\fint_Q\abs{U(x)U(y)^{-1}U(y)M}^p dx \\
&\leq\fint_Q\abs{U(x)U(y)^{-1}}^p dx\, \abs{U(y) M}^p
\leq [U]_{\A_p(\R^n)}^p \abs{U(y)M}^p,
\end{align*}
which is precisely the classical $A_1$ condition on $\abs{UM}^p$ with constant $[U]_{\A_p(\R^n)}^p$.

We turn to the second bound with $p\in(1,\infty)$.
Let first $l=1$, i.e., $M=e\in\F^m\setminus\{\mathbf 0\}$.
The result in this case is stated in \cite[Lemma 3.2(2)]{DKPS}, referring to a ``slightly different'' formulation in \cite[Lemma 2.2]{IKP}, whose proof in turn refers to \cite{NT96,Vol97}. For completeness, we give a short argument here, based on Lemma \ref{ave ops}. By that lemma, we have
\begin{equation}\label{from ave ops}
\begin{split}
    [\abs{Ue}^p]_{A_p(\R^n)}^{\frac1p}
    = [\abs{Ue}]_{\A_p(\R^n)}
    &\sim \sup_{Q\in\Cubes(\R^n)}
    \Norm{f\mapsto\one_Q\ave{f}_Q}{\aveL^p(Q,\abs{Ue})\to \aveL^p(Q,\abs{Ue})}, \\
    [U]_{\A_p(\R^n)}
    &\sim \sup_{Q\in\Cubes(\R^n)}\Norm{F\mapsto\one_Q\ave{F}_Q}{\aveL^p(Q,U)\to \aveL^p(Q,U)}.
\end{split}
\end{equation}

Given $e\in\F^m\setminus\{\mathbf 0\}$ and $f\in\aveL^p(Q,\abs{Ue})$, we consider the vector-valued function $F:=f(\cdot)e$. Then
\begin{equation*}
\begin{split}
   \Norm{\one_Q\ave{f}_Q}{\aveL^p(Q,\abs{Ue})}
   &=\bNorm{\one_Q\abs{U\ave{f(\cdot)e}_Q}}{\aveL^p(Q)}
   =\Norm{\one_Q\ave{f(\cdot)e}_Q}{\aveL^p(Q,U)} \\
   &\leq\Norm{F\mapsto \one_Q\ave{F}_Q}{\aveL^p(Q,U)\to\aveL^p(Q,U)}
   \Norm{f(\cdot)e}{\aveL^p(Q,U)},
\end{split}
\end{equation*}
where
\begin{equation*}
  \Norm{f(\cdot)e}{\aveL^p(Q,U)}
  =\Norm{f(\cdot)U(\cdot)e}{\aveL^p(Q)}
  =\bNorm{f(\cdot)\abs{U(\cdot)e}}{\aveL^p(Q)}
  =\Norm{f}{\aveL^p(Q,\abs{Ue})}.
\end{equation*}
This shows that
\begin{equation*}
   \Norm{f\mapsto\one_Q\ave{f}_Q}{\aveL^p(Q,\abs{Ue})\to \aveL^p(Q,\abs{Ue})}
   \leq\Norm{F\mapsto\one_Q\ave{F}_Q}{\aveL^p(Q,U)\to \aveL^p(Q,U)},
\end{equation*}
which completes the proof by \eqref{from ave ops}.

Let finally $p\in(1,\infty)$ and $l\in\Z_+$ be general. Let $(e_i)_{i=1}^l$ be an orthonormal basis of $\F^l$. Then
\begin{equation*}
  \abs{UM}\leq\Big(\sum_{i=1}^l\abs{UMe_i}^2\Big)^{\frac12}
  \leq\Big(\sum_{i=1}^l\abs{UMe_i}^p\Big)^{\frac1p} l^{(\frac12-\frac1p)_+}
  \leq(l\abs{UM}^p)^{\frac1p} l^{(\frac12-\frac1p)_+}
  =l^{\frac12\vee\frac1p}\abs{UM}.
\end{equation*}
With $w:=\abs{UM}^p$ and $w_i:=\abs{UMe_i}^p$, it follows that $w^{-1}\leq l(\sum_{j=1}^l w_j)^{-1}\leq l w_i^{-1}$. Hence
\begin{equation*}
  \fint_Q w\Big(\fint_Q w^{-\frac{1}{p-1}}\Big)^{p-1}
   \leq l^{ (\frac p2-1)_+} \sum_{i=1}^l \fint_Q w_i
      \Big( \fint_Q l^{\frac{1}{p-1}}w_i^{-\frac{1}{p-1}}\Big)^{p-1}
    \leq l^{\frac p2\vee 1}\sum_{i=1}^l [w_i]_{A_p},
\end{equation*}
and thus $[w]_{A_p}\leq l^{\frac p2\vee 1+1}\max_i [w_i]_{A_p}\leq l^{2p}\max_i [w_i]_{A_p}$, recalling that $p\in(1,\infty)$.
By the previous part of the proof with $e=Me_i$, we have $[w_i]_{A_p}\leq (c[U]_{\A_p})^p$, and hence $[w]_{A_p}\leq l^{2p}(c[U]_{\A_p})^p =(c'[U]_{\A_p})^p$, where $c':=l^2 c$ depends at most on $m$ (since $c$ does) and $l$. This completes the proof.
\end{proof}

The next reverse H\"older inequality for matrix weights is essentially \cite[Lemma 2]{MR22}, but we revisit the argument to cover the case $p\in(0,1)$, which is not considered in \cite{MR22}.

\begin{lemma}\label{RHI}
Let $p\in(0,\infty)$ and $U\in\A_p(\R^n;\F^{m\times m})$. Then for all $M\in\F^{m\times k}$ and all $r\in(0,1+c/[U]_{\A_p(\R^n)}^p]$, where $c$ depends at most on the dimensions and $p$, we have,
for every $Q\in\Cubes(\R^n)$,
\begin{equation*}
\bigg[\fint_Q \abs{U(x)M}^{pr}\,dx\bigg]^{\frac{1}{pr}}
\lesssim \bigg[\fint_Q\abs{U(x)M}^{p}\,dx\bigg]^{\frac 1p},
\end{equation*}
where the implicit positive constant is universal.
\end{lemma}

\begin{proof}
According to Lemma \ref{matrix 2 scalar}, each $\abs{U(x)M}^p$ is an $A_\infty$ weight with $A_\infty$ constant at most $C[U]_{\A_p}^p$, absorbing the exponent $p$ by allowing $p$-dependence of the constant $c$. By \cite[Theorem 2.3]{HPR}, such a weight satisfies the reverse H\"older inequality for all
\begin{equation*}
r\leq 1+\frac{1}{2^{n+1}C [U]_{\A_p}^p-1},
\end{equation*}
and then in particular for all $r$ as in the assumptions of the lemma. This completes the proof.
\end{proof}

After these preliminaries, we turn to discuss the $\A_p$-dimension of multi-parameter weights, extending \cite[Definition 2.22]{BHYY:1a} in the one-parameter case.

\begin{definition}
We say that $V\in\A_{p}(\R^{\vn})$ has $\A_{p}$-dimension $\vd\in\R^k$ if
\begin{equation*}
[V]_{\A_{p}(\Rect^2_{\vj}(\R^{\vn}))}^p
\lesssim 2^{\vj\cdot\vd}\quad\text{for all}\quad\vj\in\N^k.
\end{equation*}
\end{definition}

\begin{remark}
Note that, if $k=1$ and $V=W^{\frac{1}{p}}\in\A_p(\R^n)$, then
$V\in\A_p(\R^n)$ has $\A_p$-dimension $d$ if and only if
$W\in A_p(\R^n)$ has $A_p$-dimension $d$ in the sense of \cite[Definition 2.22]{BHYY:1a}.
When $m=1$, this characterises the reverse doubling property of the weight $W$
(see \cite[Proposition 2.36]{BHYY:1a}).
\end{remark}

\begin{proposition}\label{d<n}
Every $V\in\A_{p}(\R^{\vn})$ has some $\A_{p}$-dimension $\vd$ with
\begin{equation*}
\vnull\leq \vd<\vn,
\end{equation*}
where both inequalities are understood pointwise for each component.
\end{proposition}

\begin{proof}
By \cite[Proposition 2.27]{BHYY:1a}, every $U\in\A_p(\R^n)$ has some $\A_p$-dimension $d\in[0,n)$. Moreover, the proof of \cite[Proposition 2.27]{BHYY:1a} shows that one can take $d(U)=\frac{n}{r(U)}$, where $r(U)>1$ is any number for which the reverse H\"older inequality of Lemma \ref{RHI} is satisfied. The said lemma shows in particular that we can choose the same $r>1$, and hence the same $d<n$ for any family of weights $U\in\A_p(\R^n)$ with uniform bound on $[U]_{\A_p(\R^n)}$.

We apply this observation to each of the families $\Vs_i:=\{x_i\mapsto V(x'\oplus x_i):x'\in\R^{\vn}\ominus\R^{n_i}\}$ with $i\in[k]$,
where $[U]_{\A_{p_i}(\R^{n_i})}$, for $U\in\Vs_i$, is uniformly bounded by the lower bound in Lemma \ref{prod vs coord}. Thus, all $U\in\Vs_i$ have some $\A_{p}$-dimension $d\in[0,n)$, or in other words
\begin{equation*}
\esssup_{x'\in\R^{\vn}\ominus\R^{n_i}} [x_i\mapsto V(x'\oplus x_i)]_{A_{p}(\Cubes^2_{j_i}(\R^{n_i}))}\lesssim 2^{d_i j_i}\quad
\text{for all}\quad j_i\in\N,
\end{equation*}
where $d_i\in[0,n_i)$.
Substituting these bounds into the upper bound of Lemma \ref{prod vs coord}, we obtain
\begin{equation*}
[V]_{A_{p}(\Rect^2_{\vj}(\R^{\vn}))}\lesssim\prod_{i=1}^k 2^{d_i j_i}=2^{\vd\cdot\vj},
\end{equation*}
and hence $V$ has $\A_{p}$-dimension $\vd=(d_1,\ldots,d_k)$.
\end{proof}

\begin{remark}
In \cite[Section 2.2]{BHYY:1a}, several examples are given of
the possible behaviour of the $\A_p(\R^n)$-dimension of $V\in\A_p(\R^n)$.
All these examples are scalar-valued, and can be readily lifted to
provide analogous examples concerning the behaviour of
the $\A_p(\R^{\vn})$-dimension of $V\in\A_p(\R^{\vn})$.
Indeed, for scalar-valued  weights $V_i$ on $\R^{n_i}$,
it is easy to see that $V(x):=\prod_{i=1}^k V_i(x_i)$ is in $\A_p(\R^{\vn})$
if and only if $V_i\in\A_p(\R^{n_i})$ for every $i\in[k]$,
and that $V$ has $\A_p(\R^{\vn})$-dimension $\vd$ if and only if
$V_i$ has $\A_p(\R^{n_i})$-dimension $d_i$ for every $i\in[k]$.
From the results of \cite[Section 2.2]{BHYY:1a}
it follows in particular that:
For every $\vd$ with $\vnull\leq\vd<\vn$ and for every $I\subset[k]$,
there exists a weight $V\in\A_p(\R^{\vn})$ such that
$V$ has $\A_p(\R^{\vn})$-dimension $\ve$ if and only if
\begin{enumerate}[\rm(i)]
\item $e_i\geq d_i$ for all $i\in I$, and
\item $e_j>d_j$ for all $j\in[k]\setminus I$.
\end{enumerate}
\end{remark}

\begin{lemma}\label{lem 2.18}
Let $p\in(0,\infty)$ and $V\in\A_p(\R^{\vn})$. Then for $\vj\in\N^k$,
\begin{equation*}
[V]_{\A_p(\Rect^2_{-\vj}(\R^{\vn}))}\lesssim\begin{cases} 1, & \text{if }p\in(0,1], \\ 2^{\vj\cdot\ve/p'}, &
\text{if }p\in(1,\infty)\text{ and }V^{-1}\text{ has }\A_{p'}\text{-dimension }\ve.\end{cases}
\end{equation*}
When $p\in(1,\infty)$, the condition on the $\A_{p'}$-dimension is also necessary.
\end{lemma}

\begin{proof}
The quantity on the left is
\begin{equation*}
[V]_{\A_p(\Rect^2_{-\vj}(\R^{\vn}))}=\sup_{R\in\Rect(\R^n)}\Norm{V,V^{-1}}{\aveL^p(R),
\aveL^{p'}(2^{-\vj}R)}.
\end{equation*}
If $p\in(0,1]$, then $p'=\infty$, and $\Norm{\ }{\aveL^\infty(2^{-\vj}R)}\leq\Norm{\ }{\aveL^\infty(R)}$, so that
\begin{equation*}
[V]_{\A_p(\Rect^2_{-\vj}(\R^{\vn}))}
\leq[V]_{\A_p(\Rect^2_{0}(\R^{\vn}))}\lesssim 1.
\end{equation*}
If $p\in(1,\infty)$, noting that $R\mapsto 2^{\vj}R$ is a bijection of $\Rect(\R^{\vn})$, we have
\begin{align*}
[V]_{\A_p(\Rect^2_{-\vj}(\R^{\vn}))}
&=\sup_{R\in\Rect(\R^n)}\Norm{V,V^{-1}}{\aveL^p(2^{\vj}R),\aveL^{p'}(R)} \\
&\sim\sup_{R\in\Rect(\R^n)}\Norm{V^{-1},V}{\aveL^{p'}(R),\aveL^p(2^{\vj}R)}
=[V^{-1}]_{\A_{p'}(\Rect^2_{\vj}(\R^{\vn}))}.
\end{align*}
By the definition of dimensions, the right-hand side is bounded by $2^{\vj\cdot\ve/p'}$ if and only if $V^{-1}$ has $\A_p$-dimension $\ve$,
which completes the proof of Lemma \ref{lem 2.18}.
\end{proof}

\begin{definition}\label{def Ap dim}
We say that $V\in\A_p(\R^{\vn})$ has \emph{$\A_p$-dimensions $(\vd,\ve,\vDe)$} if
\begin{enumerate}[\rm(1)]
\item $V$ has $\A_p$-dimension $\vd$,
\item $p\in(1,\infty)$ and $V^{-1}$ has $\A_{p'}$-dimension $\ve$, or $p\in(0,1]$ and $\ve=\vnull$, and
\item $\vDe=\vd/p+\ve/p'$.
\end{enumerate}
\end{definition}

The following result is a multi-parameter extension of \cite[Lemma 2.28]{BHYY:1a}, where the case $k=1$ was considered:

\begin{lemma}\label{sharp coef}
Let $p\in(0,\infty)$ and $V\in\A_{p}(\R^{\vn})$.
Then $V$ have $\A_p$-dimensions $(\vd,\ve,\vDe)$
if and only if, for every $P,R\in \Rect(\R^{\vn})$,
\begin{equation*}
\abs{[V]_{\aveL^{p}(P)}[V]_{\aveL^{p}(R)}^{-1}}
\lesssim
\prod_{i=1}^k\max \bigg\{\bigg[\frac{\ell(R_i)}{\ell(P_i)}\bigg]^{\frac{d_i}{p}},
\bigg[\frac{\ell(P_i)}{\ell(R_i)}\bigg]^{\frac{e_i}{p'}}\bigg\}
\bigg[ 1+\frac{\abs{c_{P_i}-c_{R_i}}}{\max\{\ell(P_i),\ell(R_i)\}}\bigg]^{\Delta_i},
\end{equation*}
where the implicit positive constant is independent of $P$ and $R$.
\end{lemma}

\begin{proof}
We first prove ``$\Longleftarrow$''.
From Lemma \ref{red Fub} and $V\in\A_{p}(\R^{\vn})$, it follows that
\begin{align*}
\Norm{V,V^{-1}}{\aveL^p(P),\aveL^{p'}(S)}
&=\abs{[V]_{\aveL^{p}(P)}[V^{-1}]_{\aveL^{p'}(R)}}
\leq \abs{[V]_{\aveL^{p}(P)}[V]_{\aveL^{p}(R)}^{-1}} \,
\abs{[V]_{\aveL^{p}(R)} [V^{-1}]_{\aveL^{p'}(R)}} \\
&\sim \abs{[V]_{\aveL^{p}(P)}[V]_{\aveL^{p}(R)}^{-1}} \,
\Norm{V,V^{-1}}{\aveL^p(P),\aveL^{p'}(R)}
\lesssim \abs{[V]_{\aveL^{p}(P)}[V]_{\aveL^{p}(R)}^{-1}}.
\end{align*}
Thus, for every $\vj\in\N^k$,
\begin{align*}
[V]_{\A_p(\Rect^2_{\vj}(\R^{\vn}))}
=\sup_{P\in \Rect(\R^{\vn})} \Norm{V,V^{-1}}{\aveL^p(P),\aveL^{p'}(2^{\vj}P)}
\lesssim 2^{\vj\cdot\vd/p}
\end{align*}
and
\begin{align*}
[V]_{\A_p(\Rect^2_{-\vj}(\R^{\vn}))}
=\sup_{P\in \Rect(\R^{\vn})} \Norm{V,V^{-1}}{\aveL^p(P),\aveL^{p'}(2^{-\vj}P)}
\lesssim 2^{\vj\cdot\ve/p'}.
\end{align*}
These, together with Lemma \ref{lem 2.18},
further imply that $V$ have $\A_p$-dimensions $(\vd,\ve,\vDe)$.

Next, we show ``$\Longrightarrow$''.
Let $S\in\Rect(\R^{\vn})$ be the minimal rectangle that contains $P\cup R$; thus, $\ell(S_i)\sim\ell(P_i)+\ell(R_i)+\abs{c_{P_i}-c_{R_i}}$ for every $i\in[k]$. We then write
\begin{equation*}
[V]_{\aveL^{p}(P)}[V]_{\aveL^{p}(R)}^{-1}
=   [V]_{\aveL^{p}(P)}[V]_{\aveL^{p}(S)}^{-1}\times [V]_{\aveL^{p}(S)}[V]_{\aveL^{p}(R)}^{-1}.
\end{equation*}
By Lemma \ref{inverses}, we obtain
\begin{equation*}
\begin{split}
\abs{ [V]_{\aveL^{p}(P)}[V]_{\aveL^{p}(S)}^{-1}  }
&=\abs{ [V]_{\aveL^{p}(S)}^{-1} [V]_{\aveL^{p}(P)}  } \\
&\leq \abs{ [V^{-1}]_{\aveL^{p'}(S)} [V]_{\aveL^{p}(P)}  }
=\abs{  [V]_{\aveL^{p}(P)}  [V^{-1}]_{\aveL^{p'}(S)} }.
\end{split}
\end{equation*}
From Lemma \ref{red Fub}, it follows that
\begin{equation*}
\abs{  [V]_{\aveL^{p}(P)}  [V^{-1}]_{\aveL^{p'}(S)} }
\sim \Norm{V,V^{-1}}{\aveL^p(P),\aveL^{p'}(S)}
\lesssim \Norm{V,V^{-1}}{\aveL^p(P),\aveL^{p'}(2^{\vj}P)},
\end{equation*}
where $2^{j_i}\sim\ell(S_i)/\ell(P_i)$. By definition,
\begin{equation*}
\Norm{V,V^{-1}}{\aveL^p(P),\aveL^{p'}(2^{\vj}P)}\lesssim 2^{\vj\cdot\vd/p}.
\end{equation*}

Similarly,
\begin{equation*}
\begin{split}
\abs{ [V]_{\aveL^{p}(S)}[V]_{\aveL^{p}(R)}^{-1}}
&\leq\abs{ [V]_{\aveL^{p}(S)}[V^{-1}]_{\aveL^{p}(R)}} \\
&\sim \Norm{V,V^{-1}}{\aveL^{p}(S),\aveL^{p}(R)}
\lesssim \Norm{V,V^{-1}}{\aveL^{p}(2^{\vh}R),\aveL^{p}(2^{-\vh}(2^{\vh}R))},
\end{split}
\end{equation*}
where $2^{h_i}\sim\ell(S_i)/\ell(R_i)$, and by Lemma \ref{lem 2.18},
\begin{equation*}
\Norm{V,V^{-1}}{\aveL^{p}(2^{\vh}R),\aveL^{p}(2^{-\vh}(2^{\vh}R))}\lesssim 2^{\vh\cdot\ve/p'}.
\end{equation*}

Combining the estimates, we have
\begin{equation*}
\begin{split}
|[V]_{\aveL^{p}(P)}  [V]_{\aveL^{p}(R)}^{-1}|
&\lesssim 2^{\vj\cdot\vd/p}2^{\vh\cdot\ve/p'}
= \prod_{i=1}^k 2^{j_i d_i/p+h_i e_i/p'}
\sim\prod_{i=1}^k \bigg[\frac{\ell(S_i)}{\ell(P_i)}\bigg]^{\frac{d_i}{p}}
\bigg[\frac{\ell(S_i)}{\ell(R_i)}\bigg]^{\frac{e_i}{p'}} \\
&=\prod_{i=1}^k\bigg[\frac{\ell(P_i)\vee\ell(R_i)}{\ell(P_i)}\bigg]^{\frac{d_i}{p}}
\bigg[\frac{\ell(P_i)\vee\ell(R_i)}{\ell(R_i)}\bigg]^{\frac{e_i}{p'}}
\bigg[\frac{\ell(S_i)}{\ell(P_i)\vee\ell(R_i)}\bigg]^{\frac{d_i}{p}+\frac{e_i}{p'}} \\
&\sim\prod_{i=1}^k\max\bigg\{
\bigg[\frac{\ell(R_i)}{\ell(P_i)}\bigg]^{\frac{d_i}{p}},
\bigg[\frac{\ell(P_i)}{\ell(R_i)}\bigg]^{\frac{e_i}{p'}}\bigg\}
\bigg[ 1+\frac{\abs{c_{P_i}-c_{R_i}}}{\ell(P_i)\vee\ell(R_i)}\bigg] ^{\Delta_i},
\end{split}
\end{equation*}
which is exactly the claimed bound.
\end{proof}

\begin{remark}
It was shown in \cite[Lemma 2.47]{BHYY:1a} that Lemma \ref{sharp coef} is sharp in the case $k=1$. The example there is scalar-valued, and hence allows readily the construction of analogous examples with arbitrary $k$ by considering scalar-valued weights of the form $V(x)=\prod_{i=1}^k V_i(x_i)$, which satisfy
\begin{equation*}
\abs{[V]_{\aveL^p(P)}[V]_{\aveL^p(R)}^{-1}}
=\prod_{i=1}^k  \abs{[V_i]_{\aveL^p(P_i)}[V_i]_{\aveL^p(R_i)}^{-1}}.
\end{equation*}
Since the estimate of each factor is sharp by \cite[Lemma 2.47]{BHYY:1a}, it follows that the full estimate of Lemma \ref{sharp coef} is sharp as well.
\end{remark}

It is convenient to record the following special case of Lemma \ref{sharp coef}. Note in particular that the bound on the right is independent of the $\A_p$ dimensions and holds for every $V\in\A_p(\R^{\vn})$. This simple bound is already enough for many purposes.

\begin{lemma}\label{sharp coef 2}
Let $V\in\A_{p}(\R^{\vn})$ have $\A_p$-dimensions $(\vd,\ve,\vDe)$.
Then, for every $P,R\in\D_{\vj}(\R^{\vn})$,
\begin{equation*}
\abs{[V]_{\aveL^{p}(P)}[V]_{\aveL^{p}(R)}^{-1}}
\lesssim\Big[1+\abs{2^{\vj}(c_P-c_R)}\Big]^{\abs{\vDe}}
\leq\Big[1+\abs{2^{\vj}(c_P-c_R)}\Big]^{|\vn|},
\end{equation*}
where the implicit positive constant is independent of $P$ and $R$.
\end{lemma}

\begin{proof}
Now $\ell(P_i)=\ell(R_i)=2^{-j_i}$. From Lemma \ref{sharp coef} and the arithmetic--geometric mean inequality with weights $\Delta_i/\abs{\vDe}$, we obtain
\begin{equation*}
\begin{split}
\abs{[V]_{\aveL^{p}(P)}[V]_{\aveL^{p}(R)}^{-1}}
&\lesssim
\prod_{i=1}^k \Big[1+\abs{2^{j_i}(c_{P_i}-c_{R_i})}\Big]^{\Delta_i} \\
&\leq\Bigg\{\sum_{i=1}^k \Big[1+\abs{2^{j_i}(c_{P_i}-c_{R_i})}\Big]\Bigg\}^{\abs{\vDe}}
\sim\Big[1+\abs{2^{\vj}(c_P-c_R)}\Big]^{\abs{\vDe}}.
\end{split}
\end{equation*}
Finally, since $\vnull\leq\vd,\ve<\vn$, we have $\vnull\leq\vDe=\vd/p+\ve/p'<\vn$, and hence $\abs{\vDe}< |\vn|$.
\end{proof}

For some applications, it will be enough to have a family of matrices $\{A_P\}_{P\in\D(\R^{\vn})}$ with bounds like those satisfied by $[V]_{\aveL^p(P)}$ in Lemma \ref{sharp coef} , but not necessarily arising from reducing operators of some $V\in\A_p(\R^{\vn})$. Adapting \cite[Definition 1.3]{Roud04} and \cite[Definition 2.1]{FR21}, we give this property a name:

\begin{definition}\label{doubling}
A sequence of non-negative matrices $\{A_P\}_{P\in\D(\R^{\vn})}$ is said to be
\begin{enumerate}[(i)]
\item \emph{strongly doubling of order} $(\va,\vb,\vc)\in[0,\infty)^k\times[0,\infty)^k\times[0,\infty)^k$ if
\begin{equation*}
\abs{A_P A_R^{-1}}
\lesssim
\prod_{i=1}^k\max \bigg\{\bigg[\frac{\ell(R_i)}{\ell(P_i)}\bigg]^{a_i},
\bigg[\frac{\ell(P_i)}{\ell(R_i)}\bigg]^{b_i}\bigg\}
\bigg[ 1+\frac{\abs{c_{P_i}-c_{R_i}}}{\max\{\ell(P_i),\ell(R_i)\}}\bigg]^{c_i}
\end{equation*}
for all $P,R\in\D(\R^{\vn})$;

\item \emph{weakly doubling of order} $d\in[0,\infty)$ if
\begin{equation*}
\abs{A_P A_R^{-1}}
\lesssim\Big[1+\abs{2^{\vj}(c_P-c_R)}\Big]^{d}
\end{equation*}
for all $P,R\in\D_{\vj}(\R^{\vn})$, for all $\vj\in\Z^k$.
\end{enumerate}
Here, the implicit positive constants are independent of $P$ and $R$.
\end{definition}

\begin{remark}\label{Ap doubling}
In the language of Definition \ref{doubling}, the last two Lemmas \ref{sharp coef} and \ref{sharp coef 2} say that, for each $V\in\A_p(\R^{\vn})$ of $\A_p$-dimensions $(\vd,\ve,\vDe)$, the sequence of reducing operators $\{[V]_{\aveL^p(P)}\}_{P\in\D(\R^{\vn})}$ is strongly doubling of order $(\vd/p,\ve/p',\vDe)$ and weakly doubling of order $\abs{\vDe}$. Moreover, the proof of Lemma \ref{sharp coef 2} shows that every strongly doubling sequence of order $(\va,\vb,\vc)$ is weakly doubling of order $\abs{\vc}$.
\end{remark}

\begin{remark}\label{str doubling}
Let us observe some alternative ways of writing the upper bound in the condition of strong doubling when $P\in\D_{\vi}(\R^{\vn})$ and $R\in\D_{\vj}(\R^{\vn})$. Then
\begin{equation*}
\frac{1}{\max\{\ell(P_h),\ell(R_h)\}}=\frac{1}{\max\{2^{-i_h},2^{-j_h}\}}=2^{\min(i_h,j_h)},
\end{equation*}
and thus
\begin{equation*}
\begin{split}
\prod_{h=1}^k\bigg(1+\frac{\abs{c_{P_h}-c_{R_h}}}{\max\{\ell(P_h),\ell(R_h)\}}\bigg)^{c_h}
&=\prod_{h=1}^k[1+\abs{2^{\min(i_h,j_h)}(c_{P_h}-c_{R_h})}]^{c_h} \\
&=\Omega^{\vc}(2^{\vi\wedge\vj}[c_P-c_R]),
\end{split}
\end{equation*}
where
\begin{equation}\label{omega}
\Omega^{\vc}(x):=\prod_{h=1}^k(1+\abs{x})^{c_h}.
\end{equation}

We also have
\begin{equation*}
\frac{\ell(R_h)}{\ell(P_h)}=\frac{2^{-j_h}}{2^{-i_h}}=2^{i_h-j_h},
\end{equation*}
hence
\begin{equation*}
\max\bigg\{\bigg[\frac{\ell(R_h)}{\ell(P_h)}\bigg]^{a_h},\bigg[\frac{\ell(P_h)}{\ell(R_h)}\bigg]^{b_h}\bigg\}
=\max\{2^{(i_h-j_h)a_h},2^{(j_h-i_h)b_h}\}
=2^{(i_h-j_h)_+a_h+(i_h-j_h)_- b_h},
\end{equation*}
and thus
\begin{equation*}
\prod_{h=1}^k\max\bigg\{\bigg[\frac{\ell(R_h)}{\ell(P_h)}\bigg]^{a_h},
\bigg[\frac{\ell(P_h)}{\ell(R_h)}\bigg]^{b_h}\bigg\}
=2^{(\vi-\vj)_+\cdot\va+(\vi-\vj)_-\cdot\vb}
\leq 2^{\abs{\vi-\vj}\max(\va,\vb)}.
\end{equation*}
\end{remark}

\subsection{Reverse H\"older inequalities} \label{reverse}

We have already encountered a version of the reverse H\"older inequality in Lemma \ref{RHI}. We now perform a more systematic study of this concept. For exponents $0<p\leq s<\infty$, a set $E$ of positive finite measure, and a function $V:E\to\F^{m\times m}$, we define
\begin{equation*}
[V]_{\RHI_{p,s}(E)}:=\sup_{e\in\C^m}\frac{\Norm{V(\cdot)e}{\aveL^s(E)}}{\Norm{V(\cdot)e}{\aveL^p(E)}},
\end{equation*}
taking $0/0:=0$ and $\infty/\infty:=\infty$ if necessary. For any collection $\Es$ of such sets, we can then define
\begin{equation*}
[V]_{\RHI_{p,s}(\Es)}:=\sup_{E\in\Es}[V]_{\RHI_{p,s}(E)}.
\end{equation*}
For these quantities, the following lemma is an analogue
of Lemma \ref{prod vs coord}, which can be regarded as a coordinate-wise (Fubini-type)
characterisation of multi-parameter matrix weights.

\begin{lemma}\label{RHI prod vs coord}
Let $p\in(0,\infty)$ and $k\in\mathbb N\cap[2,\infty)$. Then
\begin{equation*}
\begin{split}
\sup_{i\in[k]} &\esssup_{x'\in\R^{\vn}\ominus\R^{n_i}} [x_i\mapsto V(x'\oplus x_i)]_{\RHI_{p,s}(\Cubes(\R^{n_i}))} \\
&\leq[V]_{\RHI_{p,s}(\Rect(\R^{\vn}))}
\leq\prod_{i=1}^k \esssup_{x'\in\R^{\vn}\ominus\R^{n_i}} [x_i\mapsto V(x'\oplus x_i)]_{\RHI_{p,s}(\Cubes(\R^{n_i}))}.
\end{split}
\end{equation*}
\end{lemma}

\begin{proof}
\textit{The lower bound.} We may assume that $[V]_{\RHI_{p,s}(\Rect(\R^{\vn}))}<\infty$, for otherwise there is nothing to prove. Thus, $V\in L_{\loc}^s(\R^{\vn};\F^{m\times m})\subset L_{\loc}^p(\R^{\vn};\F^{m\times m})$. In the defining conditions, it is enough to consider suitable countable families of $e\in\F^m$ and $Q_i\in\Cubes(\R^{n_i})$. Fix some $i\in[k]$ and $Q_i\in\Cubes(\R^{n_i})$.
Let $x'\in \R^{\vn}\ominus\R^{n_i}$ be a Lebesgue point of each
\begin{equation*}
x'\mapsto\Norm{V(x'\oplus\cdot)e}{L^s(Q_i)}^s,\qquad x'\mapsto \Norm{V(x'\oplus\cdot)e}{L^p(Q_i)}^p.
\end{equation*}
Thus,
\begin{align*}
\Norm{V(x'\oplus\cdot)e}{\aveL^s(Q_i)}
&=\lim_{Q'\to x'}\Norm{Ve}{\aveL^s(Q'\oplus Q_i)}
\leq \lim_{Q'\to x'}[V]_{\RHI_{p,s}(\Rect(\R^{\vn}))}\Norm{Ve}{\aveL^p(Q'\oplus Q_i)} \\
&=[V]_{\RHI_{p,s}(\Rect(\R^{\vn}))}\Norm{V(x'\oplus\cdot)e}{\aveL^p(Q_i)},
\end{align*}
which is the claimed lower estimate.

\textit{The upper bound.} Let $R=R'\times Q_k$, where $R'\in\Rect(\R^{\vn}\ominus\R^{n_k})$ and $Q_k\in\Cubes(\R^{n_k})$. Then
\begin{align*}
\Norm{Ve}{\aveL^s(R)}
&=\bigg[\fint_{R'}\Norm{V(x'\oplus \cdot)e}{\aveL^s(Q_k)}^s\, dx'\bigg]^{1/s} \\
&\leq\bigg[\fint_{R'}[V(x'\oplus\cdot)]_{\RHI_{p,s}(\Cubes(\R^{n_k}))}^s\Norm{V(x'\oplus \cdot)e}{\aveL^p(Q_k)}^s \, dx'\bigg]^{1/s} \\
&\leq\esssup_{x'\in\R^{\vn}\ominus\R^{n_k}}[V(x'\oplus\cdot)]_{\RHI_{p,s}(\Cubes(\R^{n_k}))}
\bNorm{V(x'\oplus x_k)}{\aveL^s(R',dx';\aveL^p(Q_k,dx_k))}.
\end{align*}
Since $p\leq s$, it follows from Minkowski's inequality that
\begin{equation*}
\bNorm{V(x'\oplus x_k)}{\aveL^s(R',dx';\aveL^p(Q_k,dx_k))}
\leq\bNorm{V(x'\oplus x_k)}{\aveL^p(Q_k,dx_k;\aveL^s(R',dx'))}.
\end{equation*}
Here, for each fixed $x_k\in Q_k$, the inner quasi-norm $\Norm{V(x'\oplus x_k)}{\aveL^s(R',dx')}$ has the same form as what we started with, but in a lower dimension. Thus, we can iterate the same argument, arriving at the claimed conclusion.
\end{proof}

Combining the product characterizations of both $\A_p$ (Lemma \ref{prod vs coord}) and reverse H\"older properties (Lemma \ref{RHI prod vs coord}), we can deduce a reverse H\"older property of multi-parameter $\A_p$ weights as a quick corollary of the corresponding one-parameter result (Lemma \ref{RHI}):

\begin{corollary}\label{Ap to RHI}
If $V\in\A_p(\Rect(\R^{\vn}))$, then $V\in\RHI_{p,s}(\Rect(\R^{\vn}))$ for some $s\in(p,\infty)$, uniformly over bounded subsets of $\A_p(\Rect(\R^{\vn}))$.
\end{corollary}

\begin{proof}
Suppose that $V\in\A_p(\Rect(\R^{\vn}))$. By Lemma \ref{prod vs coord}, for each $i\in[k]$ and almost all $x'\in\R^{\vn}\ominus\R^{n_i}$, we have $V(x'\oplus\cdot)\in\A_p(\Cubes(\R^{n_i}))$ and, moreover, there is a uniform bound on the $\A_p$ constant. By Lemma \ref{RHI}, there is some $s\in(p,\infty)$, uniform over bounded subsets of $\A_p(\Rect(\R^{\vn}))$, such that these same weights satisfy $V(x'\oplus\cdot)\in\RHI_{p,s}(\Cubes(\R^{n_i}))$ with a universal bound on the $\RHI_{p,s}$ constants. By Lemma \ref{RHI prod vs coord}, it then follows that $V\in\RHI_{p,s}(\Rect(\R^{\vn}))$.
\end{proof}

\subsection{Interpolation of matrix-weighted Lebesgue spaces}\label{sec interp}

We first recall the concept of the Peetre $K$-functional.

\begin{definition}
Let $A_0$ and $A_1$ be quasi-Banach spaces.
Assume that both $A_0$ and $A_1$ are subspaces of the same Hausdorff topological vector spaces.
For all $t\in(0,\infty)$ and
$f\in A_0+A_1:=\{f_0+f_1:\ f_0\in A_0,\ f_1\in A_1\}$,
the \emph{Peetre $K$-functional} is defined by setting
$$
K(t,f):=\inf\{\|f_0\|_{A_0}+t\|f_1\|_{A_1}:\
f=f_0+f_1,\ f_0\in A_0,\ f_1\in A_1\}.
$$
For all $\theta\in(0,1)$ and $p\in(0,\infty)$,
the \emph{interpolation space} $(A_0, A_1)_{\theta,p}$ is defined by setting
$$
(A_0, A_1)_{\theta,p}:=\{f\in A_0+A_1:\
\|f\|_{(A_0, A_1)_{\theta,p}}<\infty\},
$$
where, for all $f\in A_0+A_1$,
$$
\|f\|_{(A_0, A_1)_{\theta,p}}:=\bigg\{\int_0^\infty
\Big[t^{-\theta} K(t,f)\Big]^p\,\frac{dt}{t}\bigg\}^{\frac1p}.
$$
\end{definition}

\begin{definition}
Let $A$ and $B$ be two positive definite matrices.
For all $\theta\in(0,1)$, let
$$
A\#_\theta B
:=A^{\frac12}\Big(A^{-\frac12}BA^{-\frac12}\Big)^\theta A^{\frac12}
$$
be the \emph{weighted geometric mean} of $A$ and $B$.
\end{definition}

Let $p\in(0,\infty]$ and $V:\ \R^{\vn}\to\C^{m\times m}$ be measurable.
The \emph{Lebesgue space} $L^p(V)$ is defined to be the set of
all vector-valued measurable functions $f$ on $\R^{\vn}$
such that $\|f\|_{L^p(V)}:=\|\,|Vf|\,\|_{L^p(\R^{\vn})}<\infty$.

The following theorem extends a classical interpolation theorem of
$L^p$ spaces with change of measures \cite[Theorem 5.5.1]{BL76}
to the matrix-weighted setting, which is precisely
a part of \cite[Theorem 4.4]{CC25}.

\begin{theorem}\label{interpol}
Let $p_0,p_1\in(0,\infty)$ and $V_0,V_1$ be two matrix weights
on a $\sigma$-finite measure space $(\Omega,\mu)$.
Assume that $\theta\in(0,1)$ and
\begin{equation*}
\frac1p=\frac{1-\theta}{p_0}+\frac{\theta}{p_1}.
\end{equation*}
Then the real interpolation method $(\ ,\ )_{\theta,p}$ satisfies
\begin{equation} \label{our interpol}
(L^{p_0}(V_0),L^{p_1}(V_1))_{\theta,p}
=L^p\Big((V_0^2 \#_\theta V_1^2)^{\frac12}\Big).
\end{equation}
\end{theorem}

\begin{remark}
Let all the notation be as in Theorem \ref{interpol}.
\begin{enumerate}[\rm(i)]
\item Recall that \cite[Theorem 4.4]{CC25} showed that,
if $p_0,p_1\in[1,\infty)$, then
\begin{equation} \label{cc interpol}
(L^{p_0}(V_0),L^{p_1}(V_1))_{\theta,p}
=L^p\Big(|V_1V_0^{-1}|^\theta V_0\Big),
\end{equation}
where $|V_1V_0^{-1}|:=[(V_1V_0^{-1})^*V_1V_0^{-1}]^{\frac12}
=(V_0^{-1}V_1^2V_0^{-1})^{\frac12}$.
A careful checking indicates that their proof still works
for all $p_0,p_1\in(0,\infty)$. Note that
\begin{align*}
V_0^2 \#_\theta V_1^2
=V_0\Big(V_0^{-1}V_1^2V_0^{-1}\Big)^\theta V_0
=V_0\Big|V_1V_0^{-1}\Big|^{2\theta} V_0
=\Big(|V_1V_0^{-1}|^\theta V_0\Big)^2
\end{align*}
and hence \eqref{cc interpol} coincides with \eqref{our interpol}.

\item If $V_0V_1=V_1V_0$ almost everywhere, then in this case one can show that
$(V_0^2 \#_\theta V_1^2)^{\frac12}=V_0^{1-\theta} V_1^{\theta}$
(see also \cite[Section 4.1.1]{CC25}) and hence
\eqref{our interpol} in this case has the following simpler version
$$
(L^{p_0}(V_0),L^{p_1}(V_1))_{\theta,p}
=L^p\Big(V_0^{1-\theta} V_1^{\theta}\Big).
$$
\end{enumerate}
\end{remark}

Theorem \ref{interpol} has the following application to $\A_p$ weights:

\begin{proposition}\label{Ap interp}
Let $p_0,p_1\in(0,\infty)$, $\theta\in(0,1)$, and $\frac1p=\frac{1-\theta}{p_0}+\frac{\theta}{p_1}$.
Assume that $V_i\in\A_{p_i}(\R^{\vn})$ for both $i\in\{0,1\}$,
and $V:=(V_0^2\#_\theta V_1^2)^{\frac12}$.
Then the following assertions hold.
\begin{enumerate}[\rm(i)]
\item\label{Ap interp ok}
If either both $p_0,p_1\in[1,\infty)$ or
both $p_0,p_1\in(0,1]$, then $V\in\A_p(\R^{\vn})$ and
$$
[V]_{\A_p(\R^{\vn})}\lesssim[V_0]_{\A_{p_0}(\R^{\vn})}^{1-\theta}[V_1]_{\A_{p_1}(\R^{\vn})}^\theta,
$$
where the implicit positive constant is independent of $\theta$.

\item\label{Ap interp fails}
The assumption of \eqref{Ap interp ok} is necessary for the conclusion of \eqref{Ap interp ok},
that is, if $0<p_0<1<p_1<\infty$, then there exist
$V_0\in\A_{p_0}(\R^{\vn})$ and $V_1\in\A_{p_1}(\R^{\vn})$
such that $V\notin\A_p(\R^{\vn})$.
\end{enumerate}
\end{proposition}

\begin{proof}
We first consider the case where $p_0,p_1\in[1,\infty)$.
In this case, $p\in[1,\infty)$. By Lemma \ref{ave ops}, we find that
\begin{equation*}
[V]_{\A_p(\R^{\vn})}\sim\sup_{R\in\Rect(\R^{\vn})}
\|f\mapsto\one_R\ave{f}_R\|_{\aveL^p(R;V)\to \aveL^p(R;V)}
\end{equation*}
and likewise with each $(V_i,p_i)$ in place of $(V,p)$.
From Theorem \ref{interpol}, it follows that
\begin{equation}\label{LpVinterp}
L^p(R;V)=(L^{p_0}(R;V_0),L^{p_1}(R;V_1))_{p,\theta}.
\end{equation}
This, together with \cite[Theorem 3.11.2]{BL76},
further implies that, for each $R\in\Rect(\R^{\vn})$,
\begin{equation*}
\|f\mapsto\one_R\ave{f}_R\|_{\aveL^p(R;V)\to \aveL^p(R;V)}
\leq\prod_{i=0}^1 \|f\mapsto\one_R\ave{f}_R\|_{\aveL^{p_i}(R;V_i)\to \aveL^{p_i}(R;V_i)}^{\theta_i},
\end{equation*}
where $\theta_0:=1-\theta$ and $\theta_1:=\theta$. Taking the supremum over $R\in\Rect(\R^{\vn})$ shows that
\begin{equation}\label{ApVinterp}
[V]_{\A_p(\R^{\vn})}\lesssim[V_0]_{\A_{p_0}(\R^{\vn})}^{1-\theta}[V_1]_{\A_{p_1}(\R^{\vn})}^\theta.
\end{equation}

Next, we consider the case where $p_0,p_1\in(0,1]$.
In this case, $p\in(0,1]$. By Lemma \ref{ave ops}, we conclude that
\begin{equation*}
[V]_{\A_p(\R^{\vn})}\sim\sup_{R\in\Rect(\R^{\vn})}
\|f\mapsto\one_R\otimes f\|_{\aveL^p(R;V)\to \aveL^p(R\times R;V\otimes\one)},
\end{equation*}
where $(V\otimes\one)(x,y):=V(x)$, and likewise with each $(V_i,p_i)$ in place of $(V,p)$. As before, Theorem \ref{interpol} gives interpolation property \eqref{LpVinterp}, and also the analogue with $R$ replaced by $R\times R$, each $V_i$ by $V_i\otimes\one$, and $V$ by
\begin{equation*}
\Big[(V_0\otimes\one)^2\#_\theta(V_1\otimes\one)^2\Big]^{\frac12}
=\Big(V_0^2\#_\theta V_1^2\Big)^{\frac12}\otimes\one,
\end{equation*}
where the last identity is immediate from the pointwise definition of the geometric mean. Hence \cite[Theorem 3.11.2]{BL76} implies that
\begin{equation*}
\|f\mapsto\one_R\otimes f\|_{\aveL^p(R;V)\to \aveL^p(R\times R;V\otimes\one)}
\leq\prod_{i=0}^1 \|f\mapsto\one_R\otimes f\|_{\aveL^{p_i}(R;V_i)\to \aveL^{p_i}(R\times R;V_i\otimes\one)}^{\theta_i},
\end{equation*}
and taking the supremum over $R\in\Rect(\R^{\vn})$ gives \eqref{ApVinterp} as before.

Finally, without loss of generality,
we may only consider the case where $0<p_0<1<p_1<\infty$.
The reason that we cannot repeat a reasoning as above in this case
is the qualitative turning point of the characterizing condition of
Lemma \ref{ave ops} at the exponent~$1$. Let us now see that
there is a genuine obstacle, not just an issue with the chosen method of proof.

Our counterexample will be given in terms of scalar-valued weights $v_i\in \A_{p_i}$.
(Note that, even for scalar-valued weights, there is a distinction between the classical notion of $A_p$ and our $\A_p$, which is a rescaling of the classical definition.) Then the geometric mean takes the simple form $v=v_0^{1-\theta}v_1^\theta$.
Let us also consider the one-dimensional domain $\R$ for simplicity. Let $w_i:=v_i^{p_i}\in A_{p_i}(\R)$ and
\begin{equation*}
w:=v^p=w_0^{(1-\theta)\frac{p}{p_0}}w_1^{\theta\frac{p}{p_1}}.
\end{equation*}
For scalar weights, the qualitative turning point of the $A_q$ conditions
at $q=1$ is reflected in the fact that all scalar-valued $A_q$ classes
for all $q\in(0,1]$ collapse to just the classical $A_1$.
Hence $w_0\in A_1(\R)$, $w_0\in A_{p_1}(\R)$,
and we need to show that it can happen that $w\notin A_{\max\{1,p\}}(\R)$.
We choose $w_0\equiv 1\in A_1(\R)$ and $w_1(x):=\abs{x}^\alpha$
for some $\alpha\in(0,p_1-1)$ so that $w_1\in A_{p_1}(\R)$.
Then $w(x)=\abs{x}^{\alpha\theta p/p_1}$.

If $p\in(0,1]$, it suffices to observe that $w$ is $\abs{x}$ raised to some positive power, and hence $w\notin A_1(\R)$.
If $p\in(1,\infty)$, we need to investigate this power a little more carefully. From the defining equality of $p$ and the assumption that $p_0\in(0,1)$, hence $1-\frac{1}{p_0}<0$, it follows that
\begin{equation*}
1-\frac1p
=(1-\theta)\bigg(1-\frac{1}{p_0}\bigg)+\theta\bigg(1-\frac{1}{p_1}\bigg)
<\theta\bigg(1-\frac{1}{p_1}\bigg),
\end{equation*}
and therefore $p-1<(p_1-1)\theta p/p_1$.
Due to this strict inequality, we can choose $\alpha\in(0,p_1-1)$
so that $p-1<\alpha\theta p/p_1$
and hence $w\notin A_p(\R)$.
This finishes the proof of Proposition \ref{Ap interp}.
\end{proof}

\begin{remark}\label{re2.35}
Case $p_0,p_1\in[1,\infty)$ of Proposition \ref{Ap interp}\eqref{Ap interp ok} is due to
\cite[Proposition 8.8]{BCU} by a different method,
although they also indicated in \cite[Remark 8.10]{BCU}
an alternative proof based on interpolation.
However, both their first proof of \cite[Proposition 8.8]{BCU}
and the alternative indicated in \cite[Remark 8.10]{BCU}
rely heavily on duality considerations, which are not available
for exponents below $1$. Our new approach in the proof of
Proposition \ref{Ap interp}\eqref{Ap interp ok} has the advantage of covering both cases by simple variants of the same method. Another advantage is that it allows us to stay with matrix weights throughout the argument, in contrast to changing the
point of view to their induced norms, as in the proof of \cite[Proposition 8.8]{BCU}.

The negative result of \eqref{Ap interp fails} suggests some trouble to prospective
extrapolation results in the full scale $p\in(0,\infty)$,
while extrapolation in the range $p\geq 1$
was the main application of \cite[Proposition 8.8]{BCU}.
\end{remark}

Indeed, \cite[Proposition 8.8]{BCU} includes more conclusions.
Proposition \ref{Ap interp} together with the following Proposition \ref{BCU2}
contains \cite[Proposition 8.8]{BCU}.

\begin{proposition} \label{BCU2}
Let $V_0,V_1$ be two matrix weights on $\R^{\vn}$.
Assume that $\theta\in(0,1)$ and $V:=(V_0^2\#_\theta V_1^2)^{\frac12}$.
Then the following assertions hold.
\begin{enumerate}[\rm(i)]
\item\label{BCU22}
If $p\in(0,\infty)$ and $V_0^{-1},V_1^{-1}\in \A_p(\R^{\vn})$,
then $V^{-1}\in\A_p(\R^{\vn})$.

\item\label{BCU21}
If $p\in[1,\infty)$, $V_0\in \A_p(\R^{\vn})$,
and $V_1^{-1}\in \A_1(\R^{\vn})$,
then $V\in\A_{\frac p{1-\theta}}(\R^{\vn})$.
\end{enumerate}
\end{proposition}

\begin{proof}
\eqref{BCU22}: Note that
\begin{equation}\label{ga inv}
\begin{split}
V^{-1}
&=\Big[V_0\Big(V_0^{-1}V_1^2V_0^{-1}\Big)^\theta V_0\Big]^{-\frac12}
=\Big[V_0^{-1}\Big(V_0^{-1}V_1^2V_0^{-1}\Big)^{-\theta} V_0^{-1}\Big]^{\frac12}\\
&=\Big[V_0^{-1}\Big(V_0V_1^{-2}V_0\Big)^{\theta} V_0^{-1}\Big]^{\frac12}
=\Big(V_0^{-2}\#_\theta V_1^{-2}\Big)^{\frac12}.
\end{split}
\end{equation}
This, together with $V_0,V_1\in \A_p(\R^{\vn})$ and
Proposition \ref{Ap interp}(i), further implies that $V^{-1}\in\A_p(\R^{\vn})$.
This finishes the proof of \eqref{BCU22}.

\eqref{BCU21}, case $p\in(1,\infty)$: In this case, the assumption $V_0\in\A_p(\R^{\vn})$ is equivalent to $V_0^{-1}\in \A_{p'}(\R^{\vn})$. Thus, Proposition \ref{Ap interp} applies to $V_0^{-1}$ and $V_1^{-1}$ in place of $V_0$ and $V_1$ to show that
\begin{equation*}
  (V_0^{-2}\#_\theta V_1^{-2})^{\frac12}\in\A_r(\R^{\vn}),\qquad
  \frac{1}{r}=\frac{1-\theta}{p'}+\frac{\theta}{1}=\frac{1}{p'}+\frac{\theta}{p}.
\end{equation*}
By \eqref{ga inv}, the weight above is $V^{-1}$, and hence
\begin{equation*}
  V=(V_0^{-2}\#_\theta V_1^{-2})^{-\frac12}\in \A_{r'}(\R^{\vn}),\qquad
  \frac{1}{r'}=1-\frac{1}{r}=\frac{1-\theta}{p}.
\end{equation*}
Thus, $V\in A_{\frac{p}{1-\theta}}$, as claimed.

\eqref{BCU21}, general case $p\in[1,\infty)$: It remains to consider $p=1$, but the argument that this requires will also provide another proof in the case $p\in(1,\infty)$. Feeding the assumption $V_1^{-1}\in\A_1(\R^{\vn})$ into Corollary \ref{Ap to RHI}, we find that $V_1^{-1}\in\operatorname{RHI}_{1,s}(\R^{\vn})$ for some $s\in(1,\infty)$. Fixing one such $s_0\in(1,\infty)$, it is clear from the definition of the reverse H\"older property that $[V_1^{-1}]_{\operatorname{RHI}_{1,s}(\R^{\vn})}\leq[V_1^{-1}]_{\operatorname{RHI}_{1,s_0}(\R^{\vn})}\lesssim 1$ uniformly in $s\in(1,s_0)$. Then, for all $R\in\Rect(\R^{\vn})$,
\begin{equation*}
\begin{split}
  \abs{[V_1^{-1}]_{\aveL^s(R)}[V_1]_{\aveL^{s'}(R)}}
  &\sim \bigg[ \fint_R\abs{V_1^{-1}(x)[V_1]_{\aveL^{s'}(R)}}^s dx \bigg]^{\frac 1s} \\
  &\lesssim\fint_R\abs{V_1^{-1}(x)[V_1]_{\aveL^{s'}(R)}} dx\qquad\text{since }[V_1^{-1}]_{\operatorname{RHI}_{1,s}(\R^{\vn})}\lesssim 1 \\
  &\lesssim\fint_R\abs{V_1^{-1}(x)[V_1]_{\aveL^{\infty}(R)}} dx\qquad\text{by H\"older's inequality} \\
  &\lesssim [V_1^{-1}]_{\A_1(\R^{\vn})}\lesssim 1\qquad\text{by definition and assumption}.
\end{split}
\end{equation*}
Taking the supremum over $R\in\Rect(\R^{\vn})$ shows that $[V_1^{-1}]_{\A_s(\R^{\vn})}\lesssim 1$, uniformly in $s\in(1,s_0)$. (For this uniformity, it is important to observe that the implicit positive constants above depend only on the dimensions of the matrices and the quasi-triangle constants of the spaces involved; since all involved spaces $\aveL^p(R)$, with $p\in[1,\infty]$, satisfy the usual triangle inequality with constant $1$, the latter dependence is not an issue, and we only have dependence on the dimensions of the matrices.)
By duality of the $\A_p$ conditions with $p\in(1,\infty)$, the uniform estimate $[V_1^{-1}]_{\A_s(\R^{\vn})}\lesssim 1$ is equivalent to $[V_1]_{\A_{s'}(\R^{\vn})}\lesssim 1$ uniformly in $s'\in(s_0',\infty)$.

We can hence apply Proposition \ref{Ap interp} to $V_0\in\A_p(\R^{\vn})$ (assumption) and $V_1\in\A_{s'}(\R^{\vn})$ (as we just proved) to conclude that
\begin{equation}\label{VinAr}
  V\in \A_r(\R^{\vn}),\qquad\frac{1}{r}=\frac{1-\theta}{p}+\frac{\theta}{s'},
\end{equation}
with $[V]_{\A_r(\R^{\vn})}\lesssim 1$ uniformly in $s'\in(s_0',\infty)$, hence $r\in(r_0,\frac{p}{1-\theta})$ for some $r_0\in(0,\frac{p}{1-\theta})$. To conclude, we would like to be able to pass to the limit $r\to\frac{p}{1-\theta}$. The rest of the proof consists of implementing this in a precise way.

Denote $q:=\frac{p}{1-\theta}$.
For $R\in\operatorname{Rect}(\R^{\vn})$ and $N\in\N$, dominated convergence shows that
\begin{equation*}
\begin{split}
  \fint_R\min(\abs{V(x)[V]_{\aveL^{q'}(R)}}^q,N)dx
  &=\lim_{r\to q}\fint_R\min(\abs{V(x)[V]_{\aveL^{q'}(R)}}^{r},N)dx \\
  &\lesssim\limsup_{r\to q}\fint_R\abs{V(x)[V]_{\aveL^{r'}(R)}}^{r}dx \\
  &\lesssim\limsup_{r\to q}\ [V]_{\A_{r}(\R^{\vn})}^r\lesssim 1
\end{split}
\end{equation*}
by H\"older's inequality $(q'<r'$) and dropping the minimum in the first estimate, followed by the definition of $\A_{r}(\R^{\vn})$ and the uniform bound on $[V]_{\A_r(\R^{\vn})}$ observed after \eqref{VinAr}. Letting $N\to\infty$ and using monotone convergence, we can also drop the minimum from the left-hand side to see that $[V]_{\A_q(\R^{\vn})}\lesssim 1$. Since $q=\frac{p}{1-\theta}$, this finishes the proof of \eqref{BCU21}.
\end{proof}

Applying Proposition \ref{BCU2} with $p=1$ and $\theta=\frac1{q'}$,
we obtain a part of Jones factorisation theorem for multi-parameter matrix weights; we omit the details.

\begin{corollary}\label{Jones}
Let $q\in(1,\infty)$, $V_0\in \A_1(\R^{\vn})$,
and $V_1^{-1}\in \A_1(\R^{\vn})$.
If $V_0V_1=V_1V_0$ almost everywhere,
then $V_0^{\frac1q} V_1^{\frac1{q'}}\in\A_q(\R^{\vn})$.
\end{corollary}

In the case where $k=1$, Bownik and Cruz-Uribe \cite[Theorem 1.3]{BCU}
established the full Jones factorisation theorem for matrix weights. Since we do not have any need for it here, we will not consider the multi-parameter extension.

\section{Measurable selection of reducing operators}\label{sec:meas}

The goal of this section is to obtain a measurability result
to justify various analytic manipulations that
we will perform with the reducing operators.
Recall the concept of $p$-seminorm from Definition \ref{seminorm}.
The following theorem is the main result of this section.

\begin{theorem}\label{selectA}
Let $\Omega$ be a measure space and $p\in(0,1]$.
Suppose that $$r:\Omega\times\F^m\to[0,\infty)$$ is a mapping such that
\begin{enumerate}
\item[\rm(i)] for every $x\in \Omega$,
the mapping $e\mapsto r(x,e)$ is a $p$-seminorm;

\item[\rm(ii)] for every $e\in\F^m$,
the mapping $x\mapsto r(x,e)$ is measurable.
\end{enumerate}
Then there exists a measurable mapping $A:\ \Omega\to\F_+^{m\times m}$
such that $r(x,e)\sim\abs{A(x)e}$ for every $x\in\Omega$ and $e\in\F^m$,
where the positive equivalence constants depend at most on $p$ and $m$.
\end{theorem}

Closely related results are due to
Bownik and Cruz-Uribe \cite[Sections 3--4]{BCU}
and Domelevo et al. \cite[Appendix A.4]{DKPS}, but they both deal with slightly
more restrictive assumptions and hence cannot be simply quoted
for our needs. The treatment of \cite{DKPS} discusses in detail
both $\F\in\{\R,\C\}$, but only deals with proper norms,
and the fact that their closed unit balls are compact subsets of
$\F^m$ seems to be somewhat crucial for the argument.
The approach of \cite{BCU} allows seminorms, but only deals with $\F=\R$.
Assuming only the $p$-triangle inequality in Theorem \ref{selectA}
is a generalisation compared to both \cite{BCU,DKPS}.

Our plan is to adapt the techniques of Bownik and Cruz-Uribe \cite{BCU} which,
as said, are presented in the case $\F=\R$. Rather than modify
each intermediate result to the complex case, we begin by showing that
the real version of Theorem \ref{selectA} already implies its complex version,
which then permits concentrating on $\F=\R$ for the majority of the argument.
Some lemmas can then be directly quoted from \cite{BCU}, but we
provide the necessary details for the remaining part of the argument.

\begin{lemma}\label{selectA RtoC}
If Theorem \ref{selectA} holds for $\F=\R$, then it also holds for $\F=\C$.
\end{lemma}

\begin{proof}
As usual, we may identify $\C^m\simeq\R^{2m}$.
Among the properties of a $p$-seminorm, only the homogeneity
$r(x,\lambda e)=\abs{\lambda}r(x,e)$, for all $\lambda\in\F$,
refers to the scalar field $\F$, and this requirement is clearly stronger for $\F=\C$ than for $\F=\R$. Thus, under the assumptions in the $\C^m$ case
of Theorem \ref{selectA}, all the assumptions of the $\R^{2m}$ case are satisfied, and hence also the conclusions of the $\R^{2m}$ case
of Theorem \ref{selectA} hold, under the assumption that Theorem \ref{selectA} holds in this case. On the other hand, these conclusions are {\em a priori} weaker than the claims in the $\C^m$ case: identifying $A_{\R}(x)\in\R^{2m\times 2m}$ with an operator on $\C^m$, we only obtain an $\R$-linear operator in general. Expressed in the complex notation, such an operator has the form
\begin{equation*}
e\mapsto A_0(x)e+\overline{A_1(x)e}
\end{equation*}
for some $A_0(x),A_1(x)\in\C^{m\times m}$, which are easily seen to depend measurably on $x\in\Omega$ when $A_{\R}(x)$ does. Hence the conclusion of Theorem \ref{selectA} in the real case says that, for every $x\in\Omega$ and $e\in\C^m$,
\begin{equation*}
r(x,e)\sim\big|A_0(x)e+\overline{A_1(x)e}\big|.
\end{equation*}
In particular, we may also apply this to $ie$ in place of $e$, where $i=\sqrt{-1}$.
Combined with the complex homogeneity of both $r(x,\cdot)$ and $\abs{\ }$, it follows that
\begin{align*}
r(x,e)
&=r(x,ie)
\sim \big|A_0(x)ie+\overline{A_1(x)ie}\big| \\
&=\big|i\,A_0(x)e-i\,\overline{A_1(x)e}\big|
=\big|A_0(x)e-\overline{A_1(x)e}\big|.
\end{align*}
Adding the squares of the two approximate identities for $r(x,e)$ and using the parallelogram law, we further deduce that
\begin{align*}
2r(x,e)^2
&\sim\big|A_0(x)e+\overline{A_1(x)e}\big|^2
+\big|A_0(x)e-\overline{A_1(x)e}\big|^2 \\
&=2\Big( \abs{A_0(x)e}^2+\big|\overline{A_1(x)e}\big|^2\Big)
=2\Big( \abs{A_0(x)e}^2+|A_1(x)e|^2\Big)  \\
&=2\Big( \big(A_0(x)^*A_0(x)+A_1(x)^*A_1(x)\big)e,e\Big)
=2\abs{A(x)e}^2,
\end{align*}
where $A(x):=\sqrt{A_0(x)^*A_0(x)+A_1(x)^*A_1(x)}$.
This is the claimed conclusion of Theorem \ref{selectA} in the $\C^m$ case, with the new matrix $A(x)\in\C^{m\times m}$ as defined above; it is clearly positive semidefinite and depends measurably on $x\in\Omega$, when $A_0(x)$ and $A_1(x)$ do.
This finishes the proof of Lemma \ref{selectA RtoC}.
\end{proof}

We then turn to the proof of Theorem \ref{selectA} in the real case. However, the first couple of observations are elementary,
and we still state them with generic $\F$, as this makes no difference in the proofs.
Related to the $p$-seminorms $r(x,\cdot)$, we denote their \emph{closed unit balls} $K(x)$ and their so-called \emph{polar sets}
$K^\circ(x)$, respectively, by
\begin{equation}\label{def Kx}
K(x):=\{e\in\F^m:\ r(x,e)\leq 1\}
\end{equation}
and
\begin{equation}\label{def polar}
K^\circ(x):=\{u\in\F^m:\ \abs{(u,e)}\leq 1\text{ for all }e\in K(x)\}.
\end{equation}
The \emph{interior of $K(x)$} is denoted by $\operatorname{Int}K(x)$.

\begin{lemma}\label{Kx}
Let $p\in(0,1]$.
For $p$-seminorm $e\in\F^m\mapsto r(x,e)$, its closed unit ball
$K(x)$ in \eqref{def Kx} satisfies
\begin{enumerate}[\rm(i)]
\item\label{sym} $\lambda K(x)=K(x)$ for every $\lambda\in\F$ with $\abs{\lambda}=1$;

\item\label{cloInt} $K(x)$ is closed
and $K(x)=\overline{\operatorname{Int}K(x)}$;

\item\label{contBall} $\overline{B(\mathbf{0},R(x)^{-1})}\subset K(x)$
for some $R(x)\in[0,\infty)$.
\end{enumerate}
\end{lemma}

\begin{proof}
Assertion \eqref{sym} is immediate from the homogeneity
$r(x,\lambda e)=\abs{\lambda}r(x,e)$ of the $p$-seminorm.

We now prove \eqref{contBall}.
Let $\{e_i\}_{i=1}^m$ be an orthonormal basis of $\F^m$.
By the $p$-triangle inequality of $r(x,\cdot)$
and H\"older's inequality, we conclude that, for every $e\in\F^m$,
\begin{equation}\label{R(x)}
\begin{split}
r(x,e)^p
&=r\Bigg( x,\sum_{i=1}^m(e,e_i)e_i\Bigg) ^p
\leq \sum_{i=1}^m \abs{(e,e_i)}^p r(x,e_i)^p \\
&\leq\Bigg[ \sum_{i=1}^m r(x,e_i)^{\frac{2p}{2-p}}\Bigg]^{1-\frac{p}{2}}\abs{e}^p
=:R(x)^p\abs{e}^p,
\end{split}\end{equation}
where each $r(x,e_i)$, and hence $R(x)$, is some finite number.
Thus, every $e$ with $\abs{e}\leq R(x)^{-1}$
satisfies $e\in K(x)$, i.e., the closed ball $\overline{B(\mathbf 0,R(x)^{-1})}$
is contained in $K(x)$.
(The conclusion is also valid if $R(x)=0$ with $R(x)^{-1}=\infty$.)
This proves \eqref{contBall}.

Next, we prove that $K(x)$ is closed.
From the $p$-triangle inequality of $r(x,\cdot)$ and \eqref{R(x)},
it follows that, for any $e,u\in\F^m$,
\begin{equation*}
r(x,e+u)^p
\leq r(x,e)^p+r(x,u)^p
\leq r(x,e)^p+R(x)^p\abs{u}^p.
\end{equation*}
With $e+u$ in place of $e$ and $-u$ in place of $u$, we also obtain
\begin{equation*}
r(x,e)^p\leq r(x,e+u)^p+R(x)^p\abs{u}^p
\end{equation*}
and hence
\begin{equation*}
\abs{r(x,e)^p-r(x,e+u)^p}
\leq R(x)^p\abs{u}^p.
\end{equation*}
This shows that $e\mapsto r(x,e)^p$, and then also $e\mapsto r(x,e)$,
is continuous, and thus $K(x)$, as the preimage of $[0,1]$, is closed.
This is the first half of \eqref{cloInt}.

We finally consider the second claim in \eqref{cloInt}.
If $e\in K(x)$ and $t\in(0,1)$, then $r(x,te)=t r(x,e)\leq t<1$,
so each $e\in K(x)$ can be arbitrarily well approximated by
elements $te\in\widetilde K(x):=\{e\in\F^m:\ r(x,e)<1\}$.
As a preimage of an open set, this is open, and obviously contained in $K(x)$ and thus $\operatorname{Int}K(x)$. On the other hand, if $r(x,e)=1$, then $e$ is arbitrarily close to $te\notin K(x)$ for $t>1$, and hence $e\notin\operatorname{Int}K(x)$. Thus, $\operatorname{Int} K(x)=\widetilde K(x)$, and we already observed that this is dense in $K(x)$,
which completes the proof of \eqref{cloInt} and hence Lemma \ref{Kx}.
\end{proof}

Before presenting the properties of $K^\circ(x)$,
we recall some basic concepts related to sets.
A set $E\subset\mathbb F^m$
is said to be \emph{symmetric} if $-E=E$,
and said to be \emph{convex} if
$tx+(1-t)y\in E$ for every $x,y\in E$ and $t\in(0,1)$.

\begin{lemma}\label{KxPolar}
Let $p\in(0,1]$ and $R(x)$ be the same as in Lemma \ref{Kx}\eqref{contBall}.
For $p$-seminorm $\ve\in\F^m\mapsto r(x,\ve)$,
its polar set $K^\circ(x)$ in \eqref{def polar} satisfies
\begin{enumerate}[\rm(i)]
\item\label{K0cc} $K^\circ(x)$ is closed and convex;

\item\label{K0sym} $\lambda K^\circ(x)=K^\circ(x)$ for every $\lambda\in\F$ with $\abs{\lambda}=1$;

\item\label{0inK0} $\mathbf 0\in K^\circ(x)$;

\item\label{K0subBR} $K^\circ(x)\subset\overline{B(\mathbf 0,R(x))}$.
\end{enumerate}
\end{lemma}

\begin{proof} First, from the definitions of closedness, convexity,
and $K^\circ(x)$, we immediately deduce \eqref{K0cc} through \eqref{0inK0}.
For \eqref{K0subBR}, let $u\in K^\circ(x)$.
Lemma \ref{Kx}\eqref{contBall}  shows that,
for any $e\in\mathbb F^m$ with $\abs{e}\leq R(x)^{-1}$,
we have $e\in K(x)$, which further implies $|(u,e)|\leq 1$.
Taking the supremum over all such $e$ proves that
$\abs{u}R(x)^{-1}\leq 1$ and hence
$u\in\overline{B(\mathbf 0,R(x))}$.
This finishes the proof of Lemma \ref{KxPolar}.
\end{proof}

The following useful duality statement allows us to reduce the considerations from the possibly non-convex and unbounded unit balls $K(x)$ to the more favourable polar sets $K^\circ(x)$. This differentiates our approach from that of \cite[Appendix A.4]{DKPS}, where the primary objects are the sets $K(x)$ for the larger part of their argument.

\begin{lemma}\label{dual norm}
Under the assumptions of Theorem \ref{selectA},
it follows that, for any $x\in\Omega$ and $e\in\F^m$,
\begin{equation*}
r(x,e)\sim\sup\{|(e,u)|:\ u\in K^\circ(x)\},
\end{equation*}
where the positive equivalence constants depend at most on $m$ and $p$.
\end{lemma}

\begin{proof}
Let $\co K(x)$ be the convex hull of $K(x)$.
While not equal to $K(x)$ in general,
this is still uniformly comparable to $K(x)$,
as observed in \cite[p.\,1237]{FR04}.
We recall the short argument for the reader's convenience.
By a classical theorem \cite{Car07}
(see also \cite[Theorem 17.1]{Rock}),
every $v\in\operatorname{co}K(x)$ of $K(x)\subset\F^m$,
where $\F^m=\R^m$ or $\F^m=\C^m\simeq\R^{2m}$,
is a convex combination of at most $2m+1$ vectors $v_i\in K(x)$.
Hence, for some convex coefficients $\alpha_i$,
\begin{align*}
r(x,v)^p
&=r\Bigg(x,\sum_{i=1}^{2m+1}\alpha_iv_i\Bigg)^p
\leq\sum_{i=1}^{2m+1}r(x,\alpha_iv_i)^p
=\sum_{i=1}^{2m+1}\alpha_i^p r(x,v_i)^p \\
&\leq \sum_{i=1}^{2m+1}\alpha_i^p
\leq\Bigg(\sum_{i=1}^{2m+1}\alpha_i\Bigg)^p
\Bigg(\sum_{i=1}^{2m+1}1\Bigg)^{1-p}
=1\times(2m+1)^{1-p}.
\end{align*}
Thus, $r(x,v)\leq(2m+1)^{\frac1p-1}$ for all $v\in\co K(x)$, and hence
\begin{equation*}
K(x)\subset \co K(x)\subset (2m+1)^{\frac1p-1}K(x).
\end{equation*}
Since $K(x)$ is closed by Lemma \ref{Kx}, we can also take closures above to deduce that
\begin{equation*}
K(x)\subset \overline{\co K(x)}\subset (2m+1)^{\frac1p-1}K(x).
\end{equation*}
It is clear from the definition that the polar sets of any sets $E\subset F$ satisfy the reverse containment $E^\circ\supset F^\circ$, and thus
\begin{equation*}
K^\circ(x)\supset \Big(\overline{\co K}\Big)^{\circ}(x),\ \
K^{\circ\circ}(x)\subset \Big(\overline{\co K}\Big)^{\circ\circ}(x).
\end{equation*}
It is also immediate from the definition that
$E\subset E^{\circ\circ}$ for any set,
whereas equality holds for closed, convex,
and symmetric sets by the classical bipolar theorem (see e.g.\ \cite[Theorem 14.5]{Rock}, where the definition of polar is slightly different, but is seen to coincide with the present one for symmetric sets).
The set $\overline{\co K}$ is closed and convex by definition,
and it clearly inherits the symmetry of $K(x)$,
which is guaranteed by Lemma \ref{Kx}. Thus,
\begin{equation*}
K(x)
\subset K^{\circ\circ}(x)
\subset \Big(\overline{\co K}\Big)^{\circ\circ}(x)
=\overline{\co K(x)}
\subset(2m+1)^{\frac1p-1}K(x),
\end{equation*}
and in particular
\begin{equation*}
K(x)\subset K^{\circ\circ}(x)\subset (2m+1)^{\frac1p-1}K(x).
\end{equation*}
Expressed in terms of the defining $p$-seminorm of $K(x)$
and the defining property of $K^{\circ\circ}(x)=(K^{\circ})^{\circ}(x)$,
this says that, for any $e\in \mathbb F^m$,
\begin{equation*}
(2m+1)^{-(\frac1p-1)}r(x,e)
\leq \sup\{|(e,u)|:\ u\in K^{\circ}(x)\}
\leq r(x,e).
\end{equation*}
This is a quantitative formulation of the claim that we wanted to prove, and finishes the proof of Lemma \ref{dual norm}.
\end{proof}

We then turn to the more delicate measurability issues.
To make direct use of the results from \cite{BCU},
we now concentrate on $\F=\R$, keeping in mind that
this suffices for our goal of proving Theorem \ref{selectA}
in view of Lemma \ref{selectA RtoC}. We will need
the following notion of measurability of set-valued mappings
(see, for instance, \cite[Definition 8.1.1]{AF09})
and its characterization in the subsequent lemma.

\begin{definition}\label{def F meas}
A set-valued mapping
\begin{equation}\label{K(Rm)}
F:\ \Omega\to\mathcal K(\R^m):=\{\text{closed, nonempty subsets of }\R^m\}
\end{equation}
on a measurable space $(\Omega,\As)$ is said to be {\em measurable} if
\begin{equation*}
F^{-1}(U):=\{x\in\Omega:\ U\cap F(x)\neq\varnothing\}\in\As
\end{equation*}
for every open subset $U\subset\R^m$.
\end{definition}

\begin{remark}
If $F$ is singleton-valued, i.e., $F(x)=\{f(x)\}$ for some $f:\Omega\to\R^m$ and every $x\in\Omega$, then $F^{-1}(U)=f^{-1}(U)$ coincides with the usual definition of the preimage, and $F$ is measurable in the sense of Definition \ref{def F meas} if and only if $f$ is measurable in the usual sense.

In \cite[Definition 1.1]{Graf}, ``measurable'' in the sense of Definition \ref{def F meas} is called ``weakly measurable'', which may be useful to note when comparing related results in different sources.
\end{remark}

The following lemma is a part of \cite[Theorem 8.1.4]{AF09}.

\begin{lemma}\label{equiv F meas}
A set-valued mapping $F:\ \Omega\to\mathcal K(\R^m)$
on a measurable space $(\Omega,\As)$ is measurable if and only if
\begin{equation*}
\operatorname{Graph}(F)
:=\{(x,v)\in\Omega\times\R^m:\ v\in F(x)\}\subset\Omega\times\R^m
\end{equation*}
belongs to the product $\sigma$-algebra
$\As\times\operatorname{Bor}(\R^m)$,
where $\operatorname{Bor}(\mathbb R^m)$
is the Borel $\sigma$-algebra of $\mathbb R^m$.
\end{lemma}

Now, we can show the measurability of $K^\circ$.

\begin{lemma}\label{KxPolMes}
Under the assumptions of Theorem \ref{selectA} with $\F=\R$, the mapping
\begin{equation*}
K^{\circ}:\ \Omega\to\mathcal K(\R^m),
\end{equation*}
where $K^\circ$ and $\mathcal K(\R^m)$
are as in \eqref{def polar} and \eqref{K(Rm)},
is measurable in the sense of Definition \ref{def F meas}.
\end{lemma}

\begin{proof}
By Lemma \ref{KxPolar}, we find that $K^\circ$
indeed takes values in $\mathcal K(\R^m)$.

We use the equivalent condition of Lemma \ref{equiv F meas}
with $F$ replaced by $K^\circ$ to check measurability.
Recall that $v\in K^\circ(x)$ if and only if,
for every $e\in\R^m$,
\begin{equation}\label{implic1}
r(x,e)\leq 1\Rightarrow\abs{(v,e)}\leq 1.
\end{equation}
Since both sides of \eqref{implic1} are continuous in $e$,
it is equivalent to ask this condition for all $e\in\Q^m$.
Also note from elementary logic that ``$A\Rightarrow B$''
if and only if ``$B$ or not $A$'' holds.
Hence we can write
\begin{align*}
&\{(x,v)\in\Omega\times\R^m:\ v\in K^\circ(x)\} \\
&\quad=\bigcap_{e\in\Q^m} \{ (x,v)\in\Omega\times\R^m :\
\abs{(v,e)}\leq 1 \text{ or } r(x,e)>1 \} \\
&\quad=\bigcap_{e\in\Q^m}\Big(
\Big[\Omega\times\{v\in\R^m:\ \abs{(v,e)}\leq 1\}\Big]
\cup\Big[\{ x\in\Omega:\ r(x,e)>1 \}\times\R^m\Big]\Big).
\end{align*}
This is a very explicit product-measurable subset of $\Omega\times\R^m$,
which completes the proof of Lemma~\ref{KxPolMes}.
\end{proof}

The following lemma is precisely \cite[Theorem 3.7]{BCU}.

\begin{lemma}\label{BCU37}
Suppose that $F:\ \Omega\to\mathcal K(\R^m)$ is measurable
in the sense of Definition \ref{def F meas}
and takes values in the smaller target
\begin{equation*}
\mathcal K_{\rm{bcs}}(\R^m)
:=\{E\in\mathcal K(\R^m):\
E \text{ is bounded, convex, and symmetric}\}.
\end{equation*}
Then there exists a measurable $A:\Omega\to\R^{m\times m}$ such that,
for any $x\in\Omega$, the matrix $A(x)^TA(x)$ is diagonal and
\begin{equation*}
A(x)\overline{B_{\R^m}}\subset F(x)\subset \sqrt{m}\,A(x)\overline{B_{\R^m}},
\end{equation*}
where $\overline{B_{\R^m}}:= \{x\in\mathbb R^n:\ |x|\le 1\}$.
\end{lemma}

\begin{remark}
Measurability of $F:\Omega\to\mathcal K_{\rm{bcs}}(\R^m)$ may be equivalently characterised by viewing the target as a separable metric space with the Hausdorff distance of sets, see \cite[Theorem 3.5]{BCU}, and this is the approach taken in \cite[Appendix A.4]{DKPS}. However, they largely work with the still smaller space (in the notation of \cite{BCU})
\begin{equation*}
\mathcal K_{\rm{abcs}}(\R^m)
:=\{E\in\mathcal K_{\rm{bcs}}(\R^m):\ \mathbf 0\in\operatorname{Int}E\},
\end{equation*}
and $K^\circ(x)$ may fail to be in $\mathcal K_{\rm{abcs}}(\R^m)$ under the mere assumptions of our Theorem \ref{selectA}.
\end{remark}

We are now ready to show Theorem \ref{selectA}.

\begin{proof}[Proof of Theorem \ref{selectA}]
By Lemma \ref{selectA RtoC},
it suffices to consider $\F=\R$.

Applying Lemma \ref{dual norm}, we obtain,
for all $x\in\Omega$ and $e\in\mathbb R^m$,
\begin{equation}\label{use dual norm}
r(x,e)\sim\sup\{\abs{(e,v)}: v\in K^\circ(x)\},
\end{equation}
where $K^\circ(x)\in\mathcal K_{\rm{bcs}}(\R^m)$ by Lemma \ref{KxPolar}.
Thus, Lemma \ref{BCU37} guarantees the existence of a measurable
$A: \Omega\to\R^{m\times m}$ such that,
for all $x\in\Omega$,
\begin{equation}\label{WBFx}
A(x)\overline{B_{\R^m}}\subset K^\circ(x)
\subset \sqrt{m}\,A(x)\overline{B_{\R^m}}.
\end{equation}
It easily follows from
the singular value decomposition that
$A(x)\overline{B_{\R^m}}=[A(x)A(x)^T]^{\frac12}\overline{B_{\R^m}}$,
so replacing $A(x)$ by $[A(x)A(x)^T]^{\frac12}$ if necessary,
we may assume that $A(x)$ is nonnegative definite.
Now \eqref{WBFx} implies that,
for all $x\in\Omega$ and $e\in\mathbb R^m$,
\begin{align*}
\sup\{\abs{(e,v)}:\ v\in K^\circ(x)\}
&\sim\sup\Big\{\abs{(e,A(x)u)}:\ u\in \overline{B_{\R^m}} \Big\} \\
&=\sup\Big\{\abs{(A(x)e,u)}:\ u\in \overline{B_{\R^m}} \Big\}
=\abs{A(x)e}.
\end{align*}
This, together with \eqref{use dual norm},
finishes the proof of Theorem \ref{selectA}.
\end{proof}

\section{Maximal and Carleson embedding-type inequalities} \label{sec:max and cal}

The Hardy--Littlewood maximal inequality and its vector-valued extension---the Fefferman--Stein maximal inequality \cite{FS71}---have been standard tools in the theory of Triebel--Lizorkin spaces for a long time. In the matrix-weighted theory of these spaces, a more recent addition to this toolbox consists of a group of estimates presented by Frazier and Roudenko \cite[Theorem 3.7]{FR21} and attributed by them to Fedor Nazarov; these also played a key role in the subsequent studies of the more general matrix-weighted Triebel--Lizorkin-{\em type} spaces in \cite{BHYY:1a,BHYY:2b}. These estimates may be vaguely describes as {\em Carleson embedding-type}, due their similarity with the classical Carleson embedding theorem \cite{Car62}. (While the name Carleson is not explicitly mentioned in \cite{FR21}, the suggestive notation $\Norm{\ }{\mathcal C}$ is used for a key quantity featuring in their estimates; see \cite[Lemma 3.6]{FR21}.)
The aim of this section is to develop appropriate multi-parameter extensions of these tools, which will serve a similar role as their one-parameter versions in our theory of the multi-parameter Besov--Triebel--Lizorkin-type spaces in the subsequent sections.

The structure of this section is as follows.
In Subsection \ref{sec:FS}, we study the (rescaled) strong maximal operator in the iterated mixed-norm $(L^{\vp}\ell^{\vq})_\pi$ spaces, describing its boundedness in terms of the notion of {\em admissible permutations} of the component space $L^{p_i}$ and $\ell^{q_i}$; this will place restrictions on the scales of the general Besov--Triebel--Lizorkin spaces that are covered by our theory. While the sufficiency part of the description is a straightforward application of known results and some book-keeping of indices, we also present in Subsection \ref{sec:necpart} a full set of counterexamples to show the necessity of all imposed restrictions.

In Subsection \ref{Carl}, we obtain sufficient conditions for pointwise multiplier sequences $\{\gamma_{\vj}\}_{\vj\in\Z^k}$ to satisfy the multi-parameter Carleson embedding-type inequality
\begin{equation} \label{carl equ}
\bNorm{\{\gamma_{\vj}f_{\vj}\}_{\vj\in\mathbb Z^k}}{(L^{\vp}\ell^{\vq})_\pi}
\lesssim \bNorm{\{ f_{\vj}\}_{\vj\in\mathbb Z^k}}{(L^{\vp}\ell^{\vq})_\pi},
\end{equation}
and in Subsection \ref{Carl sharp} we prove the sharpness of these conditions.

In Subsection \ref{Carl special1}, we specialise the considerations to particular multipliers $\gamma_{\vj}$ that arise from matrix-weighted considerations. For these special multipliers, our sufficient conditions for \eqref{carl equ} take a form that is similar to but formally stronger than the $\A_p$ condition, and we refer to this new condition as $\A_p^{\Carl}$. In the one-parameter case, the identity $\A_p^{\Carl}=\A_p$ follows from the boundedness of the Christ--Goldberg maximal operator \cite{CG01,Gold03}. For the general multi-parameter case, this suggests a need of a multi-parameter extension of the Christ--Goldberg maximal theorem. Such a result is recently due to Vuorinen \cite{Vuo24}, and could be directly used for our needs when $p\in(1,\infty)$. However, as usual, we cannot expect maximal function bounds for $p\in(0,1]$ unless some rescaling is applied, and this required some care in the matrix-weighted setting.

In Subsection \ref{red M}, we define and study a new reduced-operator-valued maximal function, which circumvents the above-mentioned problems, and provides an approach to prove the required inequalities in the full range $p\in(0,\infty)$. Interestingly, this also gives an alternative proof the result of \cite{Vuo24} for $p\in(1,\infty)$: while \cite{Vuo24} makes use of the convex-set-valued maximal operator of \cite{BCU} as a linearising device in the intermediate steps, our approach via the reduced-operator-valued maximal function avoids this set-valued analysis. While the preference of one over the other may be a matter of taste, an objective advantage of our new approach is covering estimates in the full range $p\in(0,\infty)$.

In the final Subsection \ref{Carl special3} of this section, we combine the results of Subsections \ref{Carl} through \ref{red M} to show that the Carleson embedding-type inequalities \eqref{carl equ} indeed hold for the particular choices of $\gamma_{\vj}$ arising from matrix weights, as needed in the subsequent applications to weighted Besov--Triebel--Lizorkin type spaces.

\subsection{Maximal inequalities}\label{sec:FS}

For a function of several variables, we denote by $\mathcal M_i$ (resp. $\mathcal M_z$) the Hardy--Littlewood maximal operator acting with respect to the $i$-th variable (resp. the variable $z$),
whichever notation seems more convenient in a particular context.
The following proposition is a generalization of
the Fefferman--Stein vector-valued maximal inequality.

For a vector of exponents $\vq=(q_1,q_2,\ldots,q_k)\in(0,\infty]^k$,
the iterated $L^{\vq}$ space is defined recursively as
\begin{equation*}
\begin{split}
  L^{(q_1)} &:=L^{q_1}\qquad\text{(the usual $L^{q_1}$ space of scalar-valued functions)},\\
  L^{(q_1,q_2,\ldots,q_k)} &:=L^{q_1}(L^{(q_2,\ldots,q_k)})
  \quad\text{(the $L^{q_1}$ space of $L^{(q_2,\ldots,q_k)}$-valued functions)},\quad k\geq 2.
\end{split}
\end{equation*}
It will be convenient to have an interpretation of the previous recursive formula also for $k=1$. For this purpose, corresponding to the empty string of exponents, we define
\begin{equation*}
  L^{\varnothing}:=\F,\qquad
 \text{hence}\qquad
  L^{(q_1)}=L^{q_1}=L^{q_1}(\F)=L^{q_1}(L^{\varnothing}),
\end{equation*}
which indeed extends the recursive formula to $k=1$.

\begin{proposition}\label{MonLpw}
Let $p\in(1,\infty)$ and $w\in A_p$
(the classical, scalar-valued Muckenhoupt class).
Let $\vq\in(1,\infty)^k\times\{\infty\}^l$ for some $k,l\in\N$,
and $L^{\vq}$ be the iteration of $L^{q_i}(\Omega_i)$ spaces
over arbitrary $\sigma$-finite measure spaces $\Omega_i$.
Then the Hardy--Little\-wood maximal operator $\mathcal M_x$ is bounded on $L^p_x(w)L^{\vq}_y$.
\end{proposition}

For our immediate needs, the point of stating this with abstract spaces $\Omega_i$ is simply to cover both $L^{q_i}(\R^{n_i})$ and $\ell^{q_i}=L^{q_i}(\Z^{n_i})$ at the same time. This result is very well known in the special case that $L^{\vq}=\ell^q$ is a usual (non-iterated) sequence space; this can be found e.g. in \cite[page 62]{CUMP:book} or \cite[Theorem 8.2.6]{HNVW2}. The version as stated is certainly known to experts, and can be proved by iterating the same argument. Since we were not able to find this explicitly stated, we provide a proof for completeness.

\begin{proof}[Proof of Proposition \ref{MonLpw}]
Let first $l=0$, i.e., there is no $L^\infty$ component in $\ell^{\vq}$.

For $k=0$, i.e., $L^{\vq}=L^{\varnothing}=\F$ according to our convention, this is the classical Muckenhoupt theorem about the boundedness of the maximal operator on $L^p(w)$ for all $w\in A_p$ (see, for instance, \cite[Theorem J.1.1]{HNVW2}). Assume for induction that we already know the said boundedness on $L^r(w;L^{\vq})$, whenever $r\in(1,\infty)$, $w\in A_r$, and $L^{\vq}$ is an iterated $L^q$ space with $\vq=(q_1,\ldots,q_k)$ of some fixed length $k\in\N$. Given a function $g$ of one more variable $z$, we apply the induction assumption to $(x,y)\mapsto g(x,z,y)$ for each fixed $z$ and take the $L^r$ norms with respect to $z$. This shows that
\begin{equation*}
\Norm{\mathcal M_x g}{L^r_z L^r_x(w) L^{\vq}_y}\lesssim  \Norm{g}{L^r_z L^r_x(w) L^{\vq}_y}
\end{equation*}
and hence, by Fubini's theorem, we obtain
\begin{equation}\label{preRdF}
\BNorm{x\mapsto  \Norm{\mathcal M_x g(x,\cdot)}{ L^r_z  L^{\vq}_y}}{L^r(w)}
\lesssim  \BNorm{x\mapsto \Norm{g(x,\cdot)}{  L^r_z L^{\vq}_y} }{L^r_x(w)}.
\end{equation}
This holds for all $w\in A_r$. Hence the pair of functions
\begin{equation*}
(f(x),h(x)):=\Big( \Norm{g(x,\cdot)}{  L^r_z L^{\vq}_y}, \Norm{\mathcal M_x g(x,\cdot)}{ L^r_z  L^{\vq}_y} \Big)
\end{equation*}
satisfy the assumptions of (the operator-free formulation of) Rubio de Francia's extrapolation theorem (see, for instance, \cite[Theorem J.2.1]{HNVW2}), which shows that the same estimate \eqref{preRdF} remains valid if we replace $L^r(w)$ by $L^p(w)$ for all $p\in(1,\infty)$ and $w\in A_p$. But this is precisely the boundedness of $\mathcal M_x$ on $L^p(w)L^r L^{\vq}=L^p(w)L^{(r,\vq)}$, where we have increased the length of the vector $(r,\vq)$ by one from the previous step.
By induction, we complete the proof of Proposition \ref{MonLpw} for $l=0$.

Let then $l\geq 1$. We write $\vq=(\vec r,\infty)$, where $\vr\in(1,\infty)^k$; hence $L^{\vq}=L^{\vr}L^\infty$. We write a function $f\in L^p(w)L^{\vq}=L^p(w)L^{\vr}L^\infty$ as $f(x,y,z)$. Now
\begin{equation*}
  \abs{f(x,y,z)}\leq\Norm{f(x,y,\cdot)}{L^\infty}=:g(x,y),\qquad
  \mathcal M_x f(x,y,z)\leq \mathcal M_x g(x,y).
\end{equation*}
Hence $\Norm{\mathcal M_x f(x,y,z)}{L^\infty_x}\leq \mathcal M_x g(x,y)$ and thus
\begin{equation*}
\begin{split}
  \Norm{\mathcal M_x f}{L^p(w)L^{\vq}}
  &=\Norm{\mathcal M_x f}{L^p(w)L^{\vr}L^\infty}
  \leq\Norm{\mathcal M_x g}{L^p(w)L^{\vr}} \\
  &\lesssim\Norm{g}{L^p(w)L^{\vr}}\quad\text{by case $l=0$ that we already proved} \\
  &=\Norm{f}{L^p(w)L^{\vr}L^\infty}
  =\Norm{f}{L^p(w)L^{\vq}}.
\end{split}
\end{equation*}
This completes the proof of Proposition \ref{MonLpw} also for $l\geq 1$.
\end{proof}

In Proposition \ref{MonLpw}, the ranges of $p$ and $\vec q$ are sharp in the following sense.

\begin{theorem}\label{MonLp}
Let $p\in(0,\infty]$ and $\vec q:=(q_1,\ldots,q_k)
\in(0,\infty]^k$ for some $k\in\mathbb Z_+$.
Let $L^{\vq}$ be
the iteration of $L^{q_i}(\Omega_i)$ spaces
over arbitrary $\sigma$-finite measure spaces $\Omega_i$.
If further we assume that,
for all $i$, there exists a countable partition of $\Omega_i$
into sets of equal measure,
then the Hardy--Little\-wood maximal operator $\mathcal M_x$
is bounded on $L^p_xL^{\vq}_y$ if and only if
\begin{equation}\label{range}
(p,\vec q)\in(1,\infty)^l\times\{\infty\}^{k+1-l}
\text{ for some }l\in\{0,\ldots,k+1\}.
\end{equation}
\end{theorem}

\begin{proof}[Proof of the sufficiency part of Theorem \ref{MonLp}]
If $p\in(1,\infty)$ and $\vec q\in(1,\infty]^k$,
then the boundedness of $\mathcal M_x$ on $L^p_x(\mathbb R^n)L^{\vq}_y$
follows from Proposition \ref{MonLpw} with $w\equiv1$.
If $p=q_1=\cdots=q_k=\infty$, note that
$$
\mathop{\operatorname{ess\,sup}}
_{y\in\Omega_1\times\cdots\times\Omega_k} \mathcal Mf(\cdot,y)
\leq \mathcal M\Bigg(\mathop{\operatorname{ess\,sup}}
_{y\in\Omega_1\times\cdots\times\Omega_k} |f(\cdot,y)|\Bigg),
$$
and the boundedness of $\mathcal M_x$ on $L^p_x(\mathbb R^n)L^{\vq}_y$ follows from the boundedness of $\mathcal M_x$ on $L^p_x(\mathbb R^n)$.
This finishes the proof of the sufficiency part of Theorem \ref{MonLp}.
\end{proof}

Since the proof of the necessity part is lengthy,
we place it in the next subsection.

\begin{remark}
Some special cases of Theorem \ref{MonLp} have already been studied.
\begin{enumerate}[\rm(i)]
\item
It is classical that $\mathcal M$ is bounded on $L^p(\mathbb R^n)$
if and only if $p\in(1,\infty]$.

\item
$\mathcal M_x$ is not bounded on $L^\infty_x(\mathbb R)\ell^2$ (see \cite[Remark 2.9]{GMT93}).

\item
$\mathcal M_x$ is not bounded on $L^p_x(\mathbb R)\ell^1$ for any $p\in(1,\infty)$
(see \cite[Proposition 2.4]{GMT93} or \cite[p.\,75, Section 5.1]{Ste93}).

\item
$\mathcal M_x$ is not bounded on $L^2_x(\mathbb R)\ell^\infty\ell^2$
(see \cite[Proposition 8.1]{NVW15}).
One of the authors (T.H.) would like to thank
Emiel Lorist for informing him about this result.
\end{enumerate}
\end{remark}

We will also need versions of Theorem \ref{MonLp},
where some of the exponents may be outside the range $(1,\infty]$.
For every $k\in\mathbb N$,
the set of all permutations of $[k]$ is denoted by $S_{[k]}$.

\begin{lemma}\label{Fefferman Stein}
Let $k\in\mathbb Z_+$, $l\in\mathbb N$, $(p_j)_{j=1}^{k+l}
=(\vp,\vq)\in(0,\infty]^k\times(0,\infty]^{l}$,
and $\pi\in S_{[k+l]}$.
Fix some $i\in[k+l]$ such that $\pi(i)\in[k]$,
and suppose that for some $h\in\{0,\ldots,k+l-i+1\}$
\begin{equation}\label{>1 infty at end}
(p_{\pi(j)})_{j=i}^{k+l}\in(1,\infty)^{(k+l-i+1)-h}\times\{\infty\}^{h},
\end{equation}
i.e., all these values are in $(1,\infty]$, and all infinite values, if any, are at the end of the sequence.
Then there exists a positive constant $C$,
which depends only on $\vp$, $\vq$, and $\vn$, such that,
for every $f\in(L^{\vp}\ell^{\vq})_\pi$,
\begin{equation*}
\Norm{ \mathcal M_{\pi(i)}f }{(L^{\vp}\ell^{\vq})_\pi}
\le C  \Norm{ f }{(L^{\vp}\ell^{\vq})_\pi}.
\end{equation*}
\end{lemma}

\begin{proof}
The claimed inequality can be spelled out as
\begin{equation*}
\Norm{  \mathcal M_{y}f(x,y,z)  }{X(dx,L^{\pi(i)}(dy,Z(dz)))} \lesssim \Norm{ f(x,y,z)  }{X(dx,L^{\pi(i)}(dy,Z(dz)))},
\end{equation*}
where $X$ and $Z$ are iterated $L^p\ell^q$ spaces consisting of the first $i-1$ and the last $k+l-i$ factors of $(L^{\vp}\ell^{\vq})_\pi$, respectively, and $\mathcal M_y$ is the maximal operator with respect to the variable $y$. [If $i=1$ or $i=k+l$, we interpret $X(dx,L^{\pi(i)}(dy,Z(dz)))$ as $L^{\pi(i)}(dy,Z(dz))$ or $X(dx,L^{\pi(i)}(dy))$, respectively, and omit the respective variable $x$ or $z$ inside the quasi-norm.]

From this we observe that it suffices to prove the boundedness of $\mathcal M_y$ on $L^{\pi(i)}(dy,Z(dz))$, as the claimed inequality then follows by applying this boundedness to $(y,z)\mapsto f(x,y,z)$ for each $x$, and taking the quasi-norm in $X(dx)$ of both sides.
Indeed, this boundedness is a special case of Theorem \ref{MonLp},
which completes the proof of Lemma \ref{Fefferman Stein}.
\end{proof}

\begin{lemma}\label{FS strong}
Let $(p_i)_{i=1}^{2k}=(\vp,\vq)\in(0,\infty]^k\times(0,\infty]^{k}$ and $\pi\in S_{[2k]}$.
Let
\begin{equation}\label{mu}
\mu(\pi):=\min\{i\in[2k]:\ \pi(i)\in[k]\},
\end{equation}
and suppose that assumption \eqref{>1 infty at end} of Lemma \ref{Fefferman Stein} is satisfied with $l=k$ and $i=\mu(\pi)$.
Then there exists a positive constant $C$,
depending only on $\vp$, $\vq$, and $\vn$, such that,
for any  $f\in (L^{\vp}\ell^{\vq})_\pi$,
\begin{equation*}
\Norm{ \mathcal M_{\vn}f }{(L^{\vp}\ell^{\vq})_\pi}
\le C  \Norm{ f }{(L^{\vp}\ell^{\vq})_\pi}.
\end{equation*}
\end{lemma}

\begin{proof}
It is well known and straightforward that
\begin{equation*}
\mathcal M_{\vn}f\leq \Bigg(\prod_{i\in[2k]:\,\pi(i)\in[k]} \mathcal M_{\pi(i)}\Bigg) f,
\end{equation*}
where the products are understood as compositions of operators. If assumption \eqref{>1 infty at end} of Lemma \ref{Fefferman Stein} is satisfied with $i=\mu(\pi)$, it follows that it is satisfied with every $i\in[2k]$ such that $\pi(i)\in[k]$. Thus, by repeated applications of Lemma \ref{Fefferman Stein}, we can extract each of the maximal operators $\mathcal M_{\pi(i)}$ from the norm, arriving at the claimed conclusion.
\end{proof}

\begin{definition}\label{adm def}
We say that $\pi\in S_{[2k]}$ is an {\em admissible permutation} for $(p_i)_{i=1}^{2k}=(\vp,\vq)\in(0,\infty]^k\times(0,\infty]^k$, if the following condition holds:
\begin{equation}\label{adm}
(p_{\pi(j)})_{j=\mu(\pi)}^{2k}\in(0,\infty)^{(2k-\mu(\pi)+1)-h}\times\{\infty\}^h
\end{equation}
for some $h\in\{0,1,\ldots,2k-\mu(\pi)+1\}$, where $\mu(\pi)$ is defined in \eqref{mu}; i.e., all infinite values of the above sequence, if any, are at the end of the sequence.
\end{definition}

\begin{theorem}\label{rescaled M}
Let $\pi\in S_{[2k]}$ be a permutation for $(\vp,\vq)\in(0,\infty]^{2k}$.
Then $\pi$ is an admissible permutation if and only if
$\mathcal M_{\vn}$ is bounded on $(L^{\frac{\vp}{a}}\ell^{\frac{\vq}{a}})_\pi$
for some $a\in(0,\infty)$.
Moreover, if $\pi$ is an admissible permutation, then,
for every $a\in(0,r(\pi,\vp,\vq))$,
$\mathcal M_{\vn}$ is bounded on $(L^{\frac{\vp}{a}}\ell^{\frac{\vq}{a}})_\pi$, where
\begin{equation}\label{rpipq-max}
r(\pi,\vp,\vq) :=\min\{p_{\pi(i)}:\ i\in\{\mu(\pi),\ldots,2k\}\}
\end{equation}
with $\mu(\pi)$ defined as in \eqref{mu}.
\end{theorem}

\begin{proof}
We first prove ``$\Longrightarrow$''.
Since admissibility depends on the position of infinities only, $\pi$ is also admissible for $(\vp/a,\vq/a)$. From the defining property of admissibility, \eqref{adm}, and the definition of $r(\pi,\vp,\vq)>a>0$, it follows that
\begin{equation*}
(p_{\pi(j)}/a)_{j=\mu(\pi)}^{2k}\in(1,\infty)^{(2k-\mu(\pi)+1)-h}\times\{\infty\}^h
\end{equation*}
for some $h\in\{0,1,\ldots,2k-\mu(\pi)+1\}$. But this is precisely condition \eqref{>1 infty at end} of Lemma \ref{Fefferman Stein} for $l=k$, $i=\pi(\mu)$, and $(\vp/a,\vq/a)$ in place of $(\vp,\vq)$. This, in turn, is the assumption of Lemma \ref{FS strong} with $(\vp/a,\vq/a)$ in place of $(\vp,\vq)$, while the claim of the proposition is also the same as the claim of the lemma with this substitution. Hence ``$\Longrightarrow$'' follows from Lemma \ref{FS strong}.

Next, we show ``$\Longleftarrow$''.
By Theorem \ref{MonLp}, we find that
if $\mathcal M_{\vn}$ is bounded on $(L^{\frac{\vp}{a}}\ell^{\frac{\vq}{a}})_\pi$
for some $a\in(0,\infty)$, then $(p_{\pi(j)}/a)_{j=\mu(\pi)}^{2k}\in(1,\infty)^{2k-\mu(\pi)+1-h}\times\{\infty\}^h$
for some $h\in\{0,\ldots,2k-\mu(\pi)+1\}$.
This further implies that $\pi$ is an admissible permutation,
which completes the proof of Theorem \ref{rescaled M}.
\end{proof}

To better understand Theorem \ref{rescaled M}, we present two simple corollaries and omit the details.
The first is the boundedness of the strong maximal operator on mixed
Lebesgue spaces.

\begin{corollary} \label{coro1}
Let $\vp\in(0,\infty]^k$.
Then $\mathcal M_{\vn}$ is bounded on $L^{\frac{\vp}{a}}$
for some $a\in(0,\infty)$ if and only if
\begin{equation}\label{equ coro1}
\vp\in (0,\infty)^{k-h}\times\{\infty\}^h
\quad\text{for some}\quad  h\in\{0,\ldots,k\}.
\end{equation}
Moreover, if \eqref{equ coro1} holds, then
$\mathcal M_{\vn}$ is bounded on $L^{\frac{\vp}{a}}$
for all $a\in(0,\min \vp)$.
\end{corollary}

\begin{remark}
In \cite[Theorem 1.4]{No19}, Nogayama gave a version of Corollary \ref{coro1} for Morrey spaces.
However, Ferenc Weisz has pointed out that the range of indices in \cite[Theorem 1.4]{No19} is not correct (see \cite[Remark 3.3]{HY21}).
\end{remark}

The second is the multi-parameter Fefferman--Stein vector-valued inequality.

\begin{corollary} \label{coro2}
Let $(\vp,\vq)\in(0,\infty]^{2k}$.
Then $\mathcal M_{\vn}$ is bounded on $L^{\frac{\vp}{a}}\ell^{\frac{\vq}{a}}$
for some $a\in(0,\infty)$, that is,
there exists a positive constant $C$ such that, for every sequence $\{f_{\vj}\}_{\vj\in\Z^k}\in L^{\frac{\vp}{a}}\ell^{\frac{\vq}{a}}$,
\begin{equation}\label{equ1 coro2}
\big\| \{\mathcal M_{\vn} f_{\vj}\}_{\vj\in\Z^k} \big\|_{L^{\frac{\vp}{a}}\ell^{\frac{\vq}{a}}}
\leq C \big\| \{f_{\vj}\}_{\vj\in\Z^k} \big\|_{L^{\frac{\vp}{a}}\ell^{\frac{\vq}{a}}},
\end{equation}
if and only if
\begin{equation}\label{equ2 coro2}
(\vp,\vq)\in (0,\infty)^{2k-h}\times\{\infty\}^h
\quad\text{for some}\quad  h\in\{0,\ldots,2k\}.
\end{equation}
Moreover, if \eqref{equ2 coro2} holds, then
\eqref{equ1 coro2} holds for all $a\in(0,\min(\vp,\vq))$.
\end{corollary}

\begin{example}\label{ex adm}
With $k=1$, consider $(p,q)\in (0,\infty]^2$. Now $S_{[2]}$ consists of two permutations, $F:=\id$ and $B:=\perm{1}{2}$, where the notation is suggestive of the fact that Triebel--Lizorkin spaces are defined in terms of an $(L^p\ell^q)_F=L^p\ell^q$ quasi-norm, and Besov spaces in terms of an $(L^p\ell^q)_B=\ell^qL^p$ quasi-norm.

For $\pi=F$, condition \eqref{adm} says: if $p=\infty$, then $q=\infty$. For $\pi=B$, condition \eqref{adm} is a tautology. In other words, the only cases that are \emph{not} admissible are when $\pi=F$, $p=\infty$, and $q<\infty$. This reflects the fact that Triebel--Lizorkin spaces with $p=\infty$ and $q<\infty$ require a different definition from all other Besov and Triebel--Lizorkin spaces, which are ``admissible'' for a simpler definition.
\end{example}

Since the main focus of this article is on weighted spaces,
we will mainly deal with the case when $\vp\in(0,\infty)^k$ and $\vq\in(0,\infty]^k$. In this case, $\pi\in S_{[2k]}$ is admissible for $(\vp,\vq)$ if and only if $(L^{\vp}\ell^{\vq})_\pi$ looks like
\begin{equation*}
\ell^{\vu}(L^{\vp}\ell^{\vv})_\sigma\ell^{\vinfty},
\end{equation*}
where
\begin{enumerate}[\rm(i)]
\item
$(\vu,\vv,\vinfty)$ is a permutation of the vector $\vq$,

\item
$\vu\in(0,\infty]^a$ may have an arbitrary mixture of finite and infinite values,

\item
$\vv\in(0,\infty)^b$ has finite values only, $\sigma\in S_{k+b}$, and

\item
$\vinfty=\{\infty\}^c$ has infinite values only.
\end{enumerate}
Here $a,b,c\in\N$ and $a+b+c=k$; if any of these values is zero, the corresponding component is understood to be missing from the space $\ell^{\vu}(L^{\vp}\ell^{\vv})_\sigma\ell^{\vinfty}$.

Thus, infinite values can appear before the first $L^p$ (where they may alternate with finite values), and as a pure sequence of infinities in the end, necessarily after the last $L^p$.

\subsection{The necessity of admissible permutations}\label{sec:necpart}

In this subsection, building on ideas from several earlier works
(see \cite{GMT93,NVW15,Ste93}), we prove
the necessity part of Theorem \ref{MonLp}.

\begin{proof}[Proof of the necessity part of Theorem \ref{MonLp}]
Assume that $\mathcal M_x$ is bounded on $L^p_x(\mathbb R^n)L^{\vq}_y$.
We first prove that $(p,\vec q)\in(1,\infty]^{k+1}$.
Let $g\in L^{\vq}$ satisfy $\|g\|_{L^{\vq}}=1$.
For every $f\in L^p(\mathbb R^n)$,
\begin{align*}
\|\mathcal Mf\|_{L^p(\mathbb R^n)}
=\|\mathcal M_x(f(x)g(y))\|_{L^p_x(\mathbb R^n)L^{\vq}_y}
\lesssim \|f(x)g(y)\|_{L^p_x(\mathbb R^n)L^{\vq}_y}
= \|f\|_{L^p(\mathbb R^n)}
\end{align*}
and hence $p\in(1,\infty]$.

Let $g\in L^{(q_2,\ldots,q_k)}$ satisfy $\|g\|_{L^{(q_2,\ldots,q_k)}}=1$.
For every $f\in L^p(\mathbb R^n)L^{q_1}(\Omega_1)$,
\begin{equation}\label{special1}
\begin{split}
\|\mathcal M_xf(x,y)\|_{L^p_x(\mathbb R^n)L^{q_1}_y(\Omega_1)}
&=\|\mathcal M_x(f(x,y)g(z))\|_{L^p_x(\mathbb R^n)L^{q_1}_y(\Omega_1)L^{(q_2,\ldots,q_k)}_z}\\
&\lesssim \|f(x,y)g(z)\|_{L^p_x(\mathbb R^n)L^{q_1}_y(\Omega_1)L^{(q_2,\ldots,q_k)}_z}
= \|f(x,y)\|_{L^p_x(\mathbb R^n)L^{q_1}_y(\Omega_1)}.
\end{split}\end{equation}
By the assumptions of the present theorem, we find that there exists a sequence
of pairwise disjoint subsets $\{F_j\}_{j\in\mathbb Z}$ in $\Omega_1$ with equal positive measure.
For every $N\in\Z_+$, $x\in\mathbb R^n$, and $y\in\Omega_1$, let
$$
f_N(x,y):=\sum_{j=1}^N
\mathbf{1}_{E_j}(x)
\mathbf{1}_{F_j}(y),
$$
where, for all $j\in[N]$,
$$
E_j:=\Big(\tfrac{j-1}{N},\tfrac{j}{N}\Big]
\times (0,1]^{n-1}.
$$
Then
$$
\|f_N(x,y)\|_{L^p_x(\mathbb R^n)L^{q_1}_y(\Omega_1)}
=\mu_1(F_1)^{\frac1{q_1}}\|\mathbf{1}_{(0,1]^n}\|_{L^p(\mathbb R^n)}
=\mu_1(F_1)^{\frac1{q_1}}.
$$
By this and \eqref{special1}, we find that
\begin{equation}\label{<infty1}
\sup_{N\in\Z_+} \|\mathcal M_xf_N(x,y)\|_{L^p_x(\mathbb R^n)L^{q_1}_y(\Omega_1)}<\infty.
\end{equation}

Next, we prove that \eqref{<infty1} fails
to hold if $q_1\in(0,1]$.
Note that
\begin{align*}
\|\mathcal M_xf_N(x,y)\|_{L^p_x(\mathbb R^n)L^{q_1}_y(\Omega_1)}
&=\Bigg\|\Bigg[\sum_{j=1}^N\mu_1(F_j)
[\mathcal M(\mathbf{1}_{E_j})]^{q_1}\Bigg]^{\frac1{q_1}}\Bigg\|_{L^p(\mathbb R^n)}\\
&=\mu_1(F_1)^{\frac1{q_1}}
\Bigg\|\Bigg[\sum_{j=1}^N
[\mathcal M(\mathbf{1}_{E_j})]^{q_1}\Bigg]^{\frac1{q_1}}\Bigg\|_{L^p(\mathbb R^n)}.
\end{align*}
Here, for any $l\in[N]$ and
$x\in(\frac{l-1}{N},\frac{l}{N}]\times(0,1]^{n-1}$,
\begin{align*}
\mathcal M(\mathbf{1}_{E_j})(x)
\gtrsim \frac{1}{|l-j|+1}
\end{align*}
and hence, for any $x\in(0,1]^n$,
\begin{align*}
\sum_{j=1}^N [\mathcal M(\mathbf{1}_{E_j})(x)]^{q_1}
\gtrsim \sum_{i=0}^{\lfloor N/2\rfloor} \frac{1}{(i+1)^{q_1}}
\sim \begin{cases}
\log (N+1), &\text{if } q_1=1,\\
N^{1-q_1}, &\text{if } q_1\in(0,1),
\end{cases}
\end{align*}
where the implicit positive constants depend only on $n$ and $q$.
Therefore, \eqref{<infty1} fails to hold if $q_1\in(0,1]$;
hence $q_1\in(1,\infty]$.
By repeating the proof of $q_1\in(1,\infty]$ with
some slight modifications, we obtain $q_l\in(1,\infty]$
for all $l\in\{2,\ldots,k\}$.
This finishes the proof of $(p,\vec q)\in(1,\infty]^{k+1}$.

To finish the proof of \eqref{range}, it remains to argue that any components of $(p,\vq)$ equal to $\infty$ must be the last entries of the vector.
To this end, we need to show the following two claims:
\begin{enumerate}[\rm(i)]
\item\label{p infty}
If $p=\infty$, then $q_1=\cdots=q_k=\infty$.

\item\label{q infty}
For every $l\in[k-1]$,
if $p\in(1,\infty)$ and $q_l=\infty$,
then $q_{l+1}=\infty$.
\end{enumerate}

We first prove \eqref{p infty}.
If $p=\infty$, we proceed by contradiction
to show $q_1=\cdots=q_k=\infty$.
Suppose, on the contrary, that $q_1\in(1,\infty)$.
For every $N\in\mathbb N$, $x\in\mathbb R^n$, and $y\in\Omega_1$, let
$$
\widetilde f_N(x,y):=\sum_{j=1}^N
\mathbf{1}_{\widetilde E_j}(x)
\mathbf{1}_{F_j}(y),
$$
where, for every $j\in[N]$,
$$
\widetilde E_j:=\Big(2^{-j},2^{-j+1}\Big]
\times (0,1]^{n-1}.
$$
Then
$$
\Big\|\widetilde f_N(x,y)\Big\|_{L^\infty_x(\mathbb R^n)L^{q_1}_y(\Omega_1)}
=[\mu_1(F_1)]^{\frac1{q_1}}\|\mathbf{1}_{(2^{-N},1]\times(0,1]^{n-1}}\|_{L^\infty(\mathbb R^n)}
=\mu_1(F_1)^{\frac1{q_1}}.
$$
By this and \eqref{special1}, we find that
\begin{equation}\label{<infty2}
\sup_{N\in\mathbb N} \Big\|\mathcal M_x\widetilde f_N(x,y)\Big\|_{L^\infty_x(\mathbb R^n)L^{q_1}_y(\Omega_1)}
<\infty.
\end{equation}
Note that, for every $N\in\mathbb N$,
$j\in[N]$, and $x\in(0,2^{-N}]\times (0,1]^{n-1}$,
\begin{align*}
\mathcal M(\mathbf{1}_{\widetilde E_j})(x)
\geq \fint_{(0,2^{-j+1}]\times I_{n-1}} \mathbf{1}_{\widetilde E_j}(z)\,dz
=\frac12,
\end{align*}
where cube $I_{n-1}\subset (0,1]^{n-1}$
satisfies $\ell(I_{n-1})=2^{-j+1}$
and $x\in (0,2^{-N}]\times I_{n-1}$.
Then, for any $x\in(0,2^{-N}]\times (0,1]^{n-1}$,
\begin{align*}
\sum_{j=1}^N [\mathcal M(\mathbf{1}_{\widetilde E_j})(x)]^{q_1}
\geq N \frac{1}{2^{q_1}},
\end{align*}
which further implies that
\begin{align*}
\Big\|\mathcal M_x\widetilde f_N(x,y)\Big\|_{L^\infty_x(\mathbb R^n)L^{q_1}_y(\Omega_1)}
&=\Bigg\|\Bigg[\sum_{j=1}^N\mu_1(F_j)
[\mathcal M(\mathbf{1}_{\widetilde E_j})]^{q_1}\Bigg]^{\frac1{q_1}}\Bigg\|_{L^\infty(\mathbb R^n)}\\
&\geq \mu_1(F_1)^{\frac1{q_1}} N^{\frac1{q_1}} \frac{1}{2}.
\end{align*}
This contradicts \eqref{<infty2}; hence $q_1=\infty$.
By repeating the proof of $q_1=\infty$ with
some slight modifications, we obtain $q_l=\infty$ for all $l\in\{2,\ldots,k\}$.
This finishes the proof of \eqref{p infty}.

Next, we show \eqref{q infty}.
Without loss of generality, we may assume $l=1$.
In this case, $p\in(1,\infty)$ and $q_1=\infty$,
we proceed by contradiction to show $q_2=\infty$.
Suppose, on the contrary, that $q_2\in(1,\infty)$.
Repeating the argument used in the proof of \eqref{special1}
with some slight modifications, we obtain,
for every $f\in L^p(\mathbb R^n)L^{q_1}(\Omega_1)L^{q_2}(\Omega_2)$,
\begin{align} \label{conj2}
\|\mathcal M_x f(x,y,z)\|_{L^p_x(\mathbb R^n)L^{q_1}_y(\Omega_1)L^{q_2}_z(\Omega_2)}
\lesssim\|f(x,y,z)\|_{L^p_x(\mathbb R^n)L^{q_1}_y(\Omega_1)L^{q_2}_z(\Omega_2)}.
\end{align}
By the assumptions of the present theorem,
we find that there exists a sequence
of pairwise disjoint subsets $\{G_i\}_{i=1}^\infty$ in $\Omega_2$
with equal positive measure.
For every $N\in\Z_+$, $x\in\mathbb R^n$,
$y\in\Omega_1$, and $z\in\Omega_2$, let
$$
g_N(x,y,z):=\mathbf{1}_{[0,1)^n}(x)
\sum_{j\in[2^N-1]}\sum_{i=1}^N
\mathbf{1}_{E_{j,i}}(x)
\mathbf{1}_{F_j}(y)\mathbf{1}_{G_i}(z),
$$
where, for every $j\in\mathbb Z$ and
$i\in[N]$, let
$$
E_{j,i}:=\Big(2^{-i}+j\cdot2^{-N},2^{-i+1}+j\cdot 2^{-N}\Big]
\times (0,1]^{n-1}.
$$
Then, for all $x\in (0,1]^n$,
\begin{align*}
\|g_N(x,y,z)\|_{L^\infty_y(\Omega_1)L^{q_2}_z(\Omega_2)}
&=\mu_2(G_1)^{\frac1{q_2}}
\Bigg\|\Bigg[\sum_{j\in[2^N-1]}\sum_{i=1}^N
\mathbf{1}_{E_{j,i}}(x)\mathbf{1}_{F_j}(y)\Bigg]^{\frac1{q_2}}
\Bigg\|_{L^\infty_y(\Omega_1)}\\
&=\mu_2(G_1)^{\frac1{q_2}}
\Bigg[\sup_{j\in[2^N-1]}
\one_{(2^{-N},1]}(x_1-j\cdot 2^{-N})
\Bigg]^{\frac1{q_2}} 
=\mu_2(G_1)^{\frac1{q_2}},
\end{align*}
and hence (noting that $g_N$ has $x$-support on $(0,1]^n$),
$$
\|g_N(x,y,z)\|_{L^p_x(\mathbb R^n)L^\infty_y(\Omega_1)L^{q_2}_z(\Omega_2)}
=\mu_2(G_1)^{\frac1{q_2}}.
$$
With \eqref{conj2}, this implies that
\begin{equation}\label{<infty3}
\sup_{N\in\mathbb N} \|\mathcal M_xg_N(x,y,z)\|_{L^p_xL^{\infty}_y(\Omega_1)L^{q_2}_z(\Omega_2)}<\infty.
\end{equation}
Let $N\in\Z_+$ and $j\in [2^N-1]$ be fixed.
For every $i\in[N]$ and $x\in(j\cdot 2^{-N},(j+1)\cdot 2^{-N}]\times (0,1]^{n-1}$, let
$$
Q_x:= (x_1-2^{-i+1}, x_1+2^{-i+1}] \times I_{n-1},
$$
where $I_{n-1}\subset(0,1]^{n-1}$ is a cube in $\mathbb R^{n-1}$
such that $Q_x\subset(0,1]^n$ is a cube in $\mathbb R^n$ and $x\in Q_x$.
Note that
\begin{align*}
\Big(2^{-i}+j\cdot2^{-N},2^{-i+1}+j\cdot 2^{-N}\Big]
&\subset ((j+1)\cdot 2^{-N}-2^{-i+1}, j\cdot 2^{-N}+2^{-i+1}]\\
&\subset (x_1-2^{-i+1}, x_1+2^{-i+1}].
\end{align*}
Then
\begin{align*}
\mathcal M(\mathbf{1}_{E_{j,i}})(x)
\geq \fint_{Q_x} \mathbf{1}_{E_{j,i}}(y)\,dy
=\frac{2^{-i}|I_{n-1}|}{2^{-i+2}|I_{n-1}|}
=\frac14
\end{align*}
and hence
\begin{equation*}
\sum_{i=1}^N [\mathcal M(\one_{E_{j,i}})(x)]^{q_2}
\geq N\frac 1{4^{q_2}}.
\end{equation*}
By the arbitrariness of $j\in [2^N-1]$,
we find that, for all $x\in(2^{-N},1]\times (0,1]^{n-1}$,
\begin{equation*}
\sup_{j\in[2^N-1]}\Bigg\{\sum_{i=1}^N [\mathcal M(\one_{E_{j,i}})(x)]^{q_2}\Bigg\}^{\frac{1}{q_2}}
\geq N^{\frac{1}{q_2}}\frac 14.
\end{equation*}
Then, for any $x\in(2^{-N},1]^n\times (0,1]^{n-1}$,
\begin{align*}
&\|\mathcal M_xg_N(x,y,z)\|_{L^{\infty}_y(\Omega_1)L^{q_2}_z(\Omega_2)}\\
&\quad=\mu_2(G_1)^{\frac1{q_2}}
\Bigg\|\Bigg[\sum_{j\in [2^N-1]}\sum_{i=1}^N
[\mathcal M(\mathbf{1}_{E_{j,i}})(x)]^{q_2}
\mathbf{1}_{F_j}(y)\Bigg]^{\frac1{q_2}}
\Bigg\|_{L^\infty_y(\Omega_1)}\\
&\quad=\mu_2(G_1)^{\frac{1}{q_2}}
\sup_{j\in [2^N-1]} \Bigg\{ \sum_{i=1}^N [\mathcal M(\one_{E_{j,i}})(x)]^{q_2}\Bigg\}^{\frac{1}{q_2}}
\geq\mu_2(G_1)^{\frac1{q_2}} N^{\frac{1}{q_2}} \frac 14,
\end{align*}
which further implies that
\begin{align*}
\|\mathcal M_xg_N(x,y,z)\|_{L^p_x(\mathbb R^n)L^{\infty}_y(\Omega_1)L^{q_2}_z(\Omega_2)}
\gtrsim \mu_2(G_1)^{\frac1{q_2}} \frac{N^{\frac{1}{q_2}}}4.
\end{align*}
This contradicts \eqref{<infty3}; hence $q_2=\infty$, which completes
the proof of (ii) and hence Theorem \ref{MonLp}.
\end{proof}

\subsection{Carleson embedding-type inequalities} \label{Carl}

In this subsection, aiming for a multi-parameter extension of the Carleson-embedding type inequalities that Frazier and Roudenko \cite[Theorem 3.7]{FR21} attribute to Fedor Nazarov, we look for conditions on a sequence of pointwise multipliers $\{\gamma_{\vj}\}_{\vj\in\Z^k}$ to obtain an inequality of the form
\begin{equation}\label{carl generic}
\bNorm{\{\gamma_{\vj}f_{\vj}\}_{\vj\in\mathbb Z^k}}{(L^{\vp}\ell^{\vq})_\pi}
\leq C \bNorm{\{ f_{\vj}\}_{\vj\in\mathbb Z^k}}{(L^{\vp}\ell^{\vq})_\pi}
\end{equation}
for all function sequences of the type that each $f_{\vj}$ is constant on each $R\in\D_{\vj}(\R^{\vn})$,
where $C$ is a positive constant independent of
$\{ f_{\vj}\}_{\vj\in\mathbb Z^k}$. Certainly, we could always estimate with
\begin{equation*}
C\leq \Norm{\{\gamma_{\vj}\}_{\vj\in\mathbb Z^k}}{L^{\infty}\ell^{\infty}},
\end{equation*}
which does not exploit the assumed piecewise constancy of
the functions $\{f_{\vj}\}_{\vj\in\Z^k}$ and is much too crude for most applications
expect the case where $p=\infty$.
The following lemma can be proved by some simple computations;
we omit the details.

\begin{lemma}\label{carl small infty}
Let $\vq\in(0,\infty]^k$ and $\pi\in S_{[2k]}$.
Then, for all measurable functions $\{\gamma_{\vj}\}_{\vj\in\mathbb Z^k}$ and
all $\{ f_{\vj}\}_{\vj\in\mathbb Z^k}\in(L^{\vinfty}\ell^{\vq})_\pi$
with each $f_{\vj}$ constant on every $R\in\D_{\vj}(\R^{\vn})$,
\begin{equation*}
\bNorm{\{\gamma_{\vj}f_{\vj}\}_{\vj\in\mathbb Z^k}}{(L^{\vec\infty}\ell^{\vq})_\pi}
\le \Norm{\{\gamma_{\vj}\}_{\vj\in\mathbb Z^k}}{L^\infty \ell^\infty}
\bNorm{\{ f_{\vj}\}_{\vj\in\mathbb Z^k}}{(L^{\vec\infty}\ell^{\vq})_\pi},
\end{equation*}
where this inequality is \emph{sharp} in the sense that
$$
\sup_{\|\{ f_{\vj}\}_{\vj\in\mathbb Z^k}\|_{(L^{\vec\infty}\ell^{\vq})_\pi}\neq 0}
\frac{\|\{\gamma_{\vj}f_{\vj}\}_{\vj\in\mathbb Z^k}\|_{(L^{\vec\infty}\ell^{\vq})_\pi}}
{\|\{ f_{\vj}\}_{\vj\in\mathbb Z^k}\|_{(L^{\vec\infty}\ell^{\vq})_\pi}}
=\Norm{\{\gamma_{\vj}\}_{\vj\in\mathbb Z^k}}{L^\infty \ell^\infty}.
$$
\end{lemma}

Now, we consider the case where $p\in(0,\infty)$.
Necessary conditions can be obtained by ``testing'' the inequality \eqref{carl generic} with special choices of functions $f_{\vj}$.
Given $\Omega\in\Open(\R^{\vn})$ and $\vi\in\Z^k$, we denote
\begin{equation*}
  \Omega_{\vi}:=\bigcup_{\genfrac{}{}{0pt}{}{R\in\mathscr D_{\vi}(\R^{\vn})}{R\subset\Omega}}R;
\end{equation*}
hence $\Omega_{\vi}\subset\Omega$, and $\one_{\Omega_{\vi}}$ is piecewise constant on the rectangles $R\in\mathscr D_{\vi}(\R^{\vn})$.

By ``testing'' \eqref{carl generic}
with $f_{\vj}=\delta_{\vj,\vi}\one_{\Omega_{\vi}}$,
we obtain the necessary condition
\begin{equation*}
\sup_{\vj\in\Z^k}\sup_{\Omega\in\Open(\R^{\vn})}
\Norm{\one_{\Omega_{\vj}}\gamma_{\vj}}{\aveL^p(\Omega)}\leq C.
\end{equation*}
In particular, noting that $P_{\vi}=P$ for $P\in\D_{\vi}(\R^{\vn})$, we have
\begin{equation}\label{sharp}
\sup_{\vj\in\Z^k}\sup_{R\in\D_{\vj}(\R^{\vn})}
\Norm{\gamma_{\vj}}{\aveL^p(R)}\leq C.
\end{equation}

In the special case of a Besov-type space $\ell^{\vq}L^p$, this simple necessary condition turns out to be sufficient as well:

\begin{proposition}\label{carl easy}
If $p\in(0,\infty)$, $\vq\in(0,\infty]^k$,
and $f_{\vj}$ is constant on each $R\in\D_{\vj}(\R^{\vn})$,
then, for every measurable function $\{\gamma_{\vj}\}_{\vj\in\mathbb Z^k}$ and
every $\{ f_{\vj}\}_{\vj\in\mathbb Z^k}\in\ell^{\vq}L^p$,
\begin{equation}\label{carl easy equ}
\bNorm{\{\gamma_{\vj}f_{\vj}\}_{\vj\in\mathbb Z^k}}{\ell^{\vq}L^p}
\leq\bigg[\sup_{\vj\in\mathbb Z^k}\sup_{R\in\D_{\vj}(\R^{\vn})}
\Norm{\gamma_{\vj}}{\aveL^p(R)}\bigg]
\bNorm{\{ f_{\vj}\}_{\vj\in\mathbb Z^k}}{\ell^{\vq}L^p},
\end{equation}
where this inequality is \emph{sharp} in the sense that
$$
\sup_{\|\{ f_{\vj}\}_{\vj\in\mathbb Z^k}\|_{\ell^{\vq}L^p}\neq 0}
\frac{\|\{\gamma_{\vj}f_{\vj}\}_{\vj\in\mathbb Z^k}\|_{\ell^{\vq}L^p}}
{\|\{ f_{\vj}\}_{\vj\in\mathbb Z^k}\|_{\ell^{\vq}L^p}}
=\sup_{\vj\in\mathbb Z^k}\sup_{R\in\D_{\vj}(\R^{\vn})}
\Norm{\gamma_{\vj}}{\aveL^p(R)}.
$$
\end{proposition}

\begin{proof}
If the supremum is zero or infinity,
then \eqref{carl easy equ} holds trivially;
otherwise we may assume by scaling that it is one.
Let $f_R$ be the constant value of $f_{\vj}$ on $R\in\D_{\vj}(\R^{\vn})$. Then
\begin{equation*}
\begin{split}
\Norm{\gamma_{\vj}f_{\vj}}{L^p(\R^{\vn})}^p
&=\sum_{R\in\D_{\vj}(\R^{\vn})}\Norm{\one_R\gamma_{\vj}f_{\vj}}{L^p(\R^{\vn})}^p
=\sum_{R\in\D_{\vj}(\R^{\vn})}\Norm{\one_R\gamma_{\vj}}{L^p(\R^{\vn})}^p\abs{f_R}^p \\
&\leq\sum_{R\in\D_{\vj}(\R^{\vn})}\abs{R}\,\abs{f_R}^p
=\Norm{f_{\vj}}{L^p(\R^{\vn})}^p.
\end{split}
\end{equation*}
Taking the $\ell^{\vq}$ quasi-norm of both sides proves \eqref{carl easy equ}.
The sharpness follows from \eqref{carl easy equ} and \eqref{sharp}.
This finishes the proof of Proposition \ref{carl easy}.
\end{proof}

For all other spaces of type $(L^{\vp}\ell^{\vq})_\pi$,
the situation becomes more complicated.
The following lemma is a multi-parameter extension of
\cite[Theorem 3.7(i)]{FR21}, which the authors of \cite{FR21}
attributed to Fedor Nazarov. The proof below
is a relatively straightforward modification,
only using the strong maximal operator in place of the classical one.

\begin{lemma}\label{carl small q}
Let $p\in(0,\infty)$, $\vnull<\vq\leq\vp:=p\cdot\vone$,
$s\in(p,\infty)$, and $\pi\in S_{[2k]}$.
Then there exists a positive constant $C$ such that,
for all measurable functions $\{\gamma_{\vj}\}_{\vj\in\mathbb Z^k}$ and
all $\{ f_{\vj}\}_{\vj\in\mathbb Z^k}\in(L^{\vp}\ell^{\vq})_\pi$
with each $f_{\vj}$ constant on every $R\in\D_{\vj}(\R^{\vn})$,
\begin{equation*}
\bNorm{\{\gamma_{\vj}f_{\vj}\}_{\vj\in\mathbb Z^k}}{(L^{\vp}\ell^{\vq})_\pi}
\le C\bigg[\sup_{\vj\in\Z^k}\sup_{R\in\D_{\vj}(\R^{\vn})}
\Norm{\gamma_{\vj}}{\aveL^s(R)}\bigg]\bNorm{\{ f_{\vj}\}_{\vj\in\mathbb Z^k}}{(L^{\vp}\ell^{\vq})_\pi}.
\end{equation*}
\end{lemma}

\begin{proof}
If the supremum is infinite or zero, there is nothing to prove. Otherwise, we may assume by scaling that it is equal to one.

Let us begin by noting that, using the notation
\begin{equation*}
\Avee{g}{\vj}:=\sum_{R\in\D_{\vj}(\R^{\vn})}\ave{g}_R\one_R,
\end{equation*}
this rescaled assumption can be stated as $\Avee{\abs{\gamma_{\vj}}^s}{\vj}\leq 1$.

Let $r\in(0,\min(\vp,\vq))$ so that
$\frac{\vp}r,\frac{\vq}r\in(1,\infty]^k$. Then we can dualise
\begin{equation}\label{carl1}
\begin{split}
  \Norm{\{\gamma_{\vj}f_{\vj}\}_{\vj\in\mathbb Z^k}}{(L^{\vp}\ell^{\vq})_\pi}^r
&=\bNorm{\{\abs{\gamma_{\vj}}^r\abs{f_{\vj}}^r\}_{\vj\in\mathbb Z^k}}{(L^{\vp/r}\ell^{\vq/r})_\pi} \\
&=\sup\Bigg\{\int_{\R^{\vn}}\sum_{\vj\in\Z^k}\abs{\gamma_{\vj}}^r\abs{f_{\vj}}^r g_{\vj}:
\Norm{\{g_{\vj}\}}{(L^{(\vp/r)'}\ell^{(\vq/r)'})_\pi}\leq 1\Bigg\}.
\end{split}
\end{equation}
Recalling that $f_{\vj}$ is constant on each $R\in\mathscr D_{\vj}$, it follows that
\begin{equation*}
\int_{\R^{\vn}}\sum_{\vj\in\Z^k}\abs{\gamma_{\vj}}^r\abs{f_{\vj}}^r g_{\vj}
=\int_{\R^{\vn}}\sum_{\vj\in\Z^k}\abs{f_{\vj}}^r \Avee{\abs{\gamma_{\vj}}^r g_{\vj}}{\vj}.
\end{equation*}
By H\"older's inequality on each $R\in\D_{\vj}(\R^{\vn})$ and our assumption of $\gamma_{\vj}$,
\begin{equation*}
\Avee{\abs{\gamma_{\vj}}^r g_{\vj}}{\vj}
\leq \Big<|\gamma_{\vj}|^{r\frac sr}\Big>_{\vj}^{\frac rs}
\Big<|g_{\vj}|^{(\frac sr)'}\Big>_{\vj}^{\frac 1{(s/r)'}}
\lesssim \Big[\mathcal M_{\vn}\Big(\abs{g_{\vj}}^{(\frac sr)'}\Big)\Big]^{\frac{1}{(s/r)'}}.
\end{equation*}
Hence
\begin{align*}
\int_{\R^{\vn}}\sum_{\vj\in\Z^k}\abs{\gamma_{\vj}}^r\abs{f_{\vj}}^r g_{\vj}
&\lesssim\int_{\R^{\vn}}\sum_{\vj\in\Z^k}\abs{f_{\vj}}^r
\Big[\mathcal M_{\vn}\Big(\abs{g_{\vj}}^{(\frac sr)'}\Big)\Big]^{\frac{1}{(s/r)'}}\\
&\leq\Big\|\{\abs{f_{\vj}}^r\}_{\vj\in\Z^k}\Big\|_{(L^{\frac{\vp}r}
\ell^{\frac{\vq}r})_\pi}
\Big\|\Big\{\Big[\mathcal M_{\vn}\Big(\abs{g_{\vj}}^{(\frac sr)'}\Big)\Big]^{\frac{1}{(s/r)'}}\Big\}_{\vj\in\Z^k}
\Big\|_{(L^{(\frac{\vp}r)'}\ell^{(\frac{\vq}r)'})_\pi} \\
&=\bNorm{\{ f_{\vj} \}_{\vj\in\Z^k}}{(L^{\vp}\ell^{\vq})_\pi}^r
\Big\|\Big\{ \mathcal M_{\vn} \Big( \abs{g_{\vj}}^{(\frac sr)'}\Big)  \Big\}_{\vj\in\Z^k}
\Big\|_{(L^{\frac{(\vp/r)'}{(s/r)'}}\ell^{\frac{(\vq/r)'}{(s/r)'}})_\pi}^{\frac 1{(s/r)'}}.
\end{align*}
Since $\infty>s>p\geq q_i>r>0$, it follows that
$\infty>\frac sr>\frac pr\geq \frac{q_i}r>1$, hence
$$1<(\frac sr)'<(\frac pr)'\leq(\frac{q_i}r)'<\infty,$$
and thus $$\vinfty>\frac{(\vq/r)'}{(s/r)'}\geq \frac{(\vp/r)'}{(s/r)'}>\vone.$$
Then, by Lemma \ref{FS strong},
we find that $\mathcal M_{\vn}$ is bounded on $(L^{\frac{(\vp/r)'}{(s/r)'}}\ell^{\frac{(\vq/r)'}{(s/r)'}})_\pi$, and we conclude that
\begin{align*}
\int_{\R^{\vn}}\sum_{\vj\in\Z^k}\abs{\gamma_{\vj}}^r\abs{f_{\vj}}^r g_{\vj}
& \lesssim \bNorm{\{ f_{\vj} \}_{\vj\in\Z^k}}{(L^{\vp}\ell^{\vq})_\pi}^r
\Big\|\Big\{ \abs{g_{\vj}}^{(\frac sr)'} \Big\}_{\vj\in\Z^k}
\Big\|_{( L^{\frac{(\vp/r)'}{(s/r)'}}\ell^{\frac{(\vq/r)'}{(s/r)'}})_\pi}^{\frac 1{(s/r)'}} \\
& = \bNorm{\{ f_{\vj} \}_{\vj\in\Z^k}}{(L^{\vp}\ell^{\vq})_\pi}^r \bNorm{\{ g_{\vj} \}_{\vj\in\Z^k} }{( L^{(\vp/r)'}\ell^{(\vq/r)'})_\pi}.
\end{align*}
Substituting back to \eqref{carl1}, we obtain
\begin{equation*}
\Norm{\{\gamma_{\vj}f_{\vj}\}_{\vj\in\Z^k}}{(L^{\vp}\ell^{\vq})_\pi}^r
\lesssim \bNorm{\{ f_{\vj} \}_{\vj\in\Z^k}}{(L^{\vp}\ell^{\vq})_\pi}^r.
\end{equation*}
This finishes the proof of Lemma \ref{carl small q}.
\end{proof}

\begin{remark}
In Lemma \ref{carl small q},
the bound
\begin{equation}\label{sharp bound}
\sup_{\vj\in\Z^k}\sup_{R\in\D_{\vj}(\R^{\vn})}
\Norm{\gamma_{\vj}}{\aveL^s(R)}
\end{equation}
can be relaxed to
$\sup_{\vj\in\Z^k}\sup_{R\in\D_{\vj}(\R^{\vn})}
\Norm{\gamma_{\vj}}{\aveL^p(R)}$
if $(L^{\vp}\ell^{\vq})_\pi$ reduces to Besov-type.
However, if $(L^{\vp}\ell^{\vq})_\pi$ is ``truly'' of Triebel--Lizorkin-type,
then the bound \eqref{sharp bound} is sharp;
see Proposition \ref{carl small q2}.
\end{remark}

In order to discuss versions of Lemma \ref{carl small q} without the restriction to $\vq\leq p\cdot\vone$, we will need some additional tools.

\begin{definition}
Let $p\in(0,\infty)$, $q\in(0,\infty]$, and $X$ be a Banach space.
Let $f:\mathbb R^{\vn}\to X$ be strongly measurable (i.e., a pointwise limit of simple functions, see \cite[Chapter 1]{HNVW1}).
Its Lorentz quasi-norm is defined as
$$
\|f\|_{L^{p,q}(X)}
:=\begin{cases}
\displaystyle
p^{\frac1q}\bigg[\int_0^\infty s^{q-1}
|\{x\in \mathbb R^{\vn}:\ \|f(x)\|_X>s\}|^{\frac qp} \, ds\bigg]^{\frac1p},
&\text{if } q\in(0,\infty),\\
\displaystyle
\sup_{s\in(0,\infty)} s |\{x\in \mathbb R^{\vn}:\ \|f(x)\|_X >s\}|^{\frac 1p},
&\text{if } q=\infty.
\end{cases}
$$
\end{definition}

Let us recall the following result concerning the real interpolation
of Banach space $X$-valued Lorentz spaces $L^{p,r}(X)$;
see \cite[Section 1.18.6, Theorem 2]{Tri78}.

\begin{lemma}\label{interp}
Let $p_0,p_1\in(1,\infty)$ satisfy $p_0\neq p_1$,
and let $r_0,r_1,r\in[1,\infty]$.
Assume that $\theta\in(0,1)$ and $\frac 1p=\frac{1-\theta}{p_0}+\frac{\theta}{p_1}$.
Then, for any Banach space $X$,
\begin{equation*}
(L^{p_0,r_0}(X),L^{p_1,r_1}(X))_{\theta,r}=L^{p,r}(X)
\end{equation*}
with equivalent norms.
\end{lemma}

\begin{remark}\label{Cwikel}
In contrast to the clean statement of Lemma \ref{interp},
interpolation of $L^{p_i,r_i}(X_i)$ with two different $X_i$
is more tricky; see \cite{Cwi74}. This issue does not appear in
Lemma \ref{interp}, where only one $X$ is involved throughout.
\end{remark}

The following lemma is a multi-parameter extension of \cite[Theorem 3.7(ii)]{FR21}, which the authors of \cite{FR21}
attributed to Fedor Nazarov. The proofs of both \cite[Theorem 3.7(ii)]{FR21} and Lemma \ref{carl >1} below rely on interpolation; however, the details are substantially different. In \cite[Theorem 3.7(ii)]{FR21}, where the sequence space has just one exponent $q$ in place of our $\vq$, the interpolation is done between the end-points $q\in\{1,\infty\}$, for a fixed $p$. In contrast, we work with a fixed $\vq$, and interpolate with respect to the parameter $p$, for two exponent $p_i$ close to the desired $p$. This is inspired by a similar interpolation strategy used in \cite[Section 8]{HMP}.

On the level of the actual result, we will now need to deal with general open sets $\Omega$, in contrast to dyadic rectangles that were enough in Lemma \ref{carl small q}.

\begin{lemma}\label{carl >1}
If $p\in(1,\infty)$, $\vq\in(1,\infty)^k$, and $s\in(p,\infty)$,
then there exists a positive constant $C$ such that,
for every sequence $\{f_{\vj}\}_{\vj\in\mathbb Z^k}\in L^p\ell^{\vq}$ and
all measurable functions $\{\gamma_{\vj}\}_{\vj\in\mathbb Z^k}$,
\begin{equation}\label{carl >1eq}
\bNorm{\{\gamma_{\vj}\Avee{f_{\vj}}{\vj}\}_{\vj\in\mathbb Z^k}}{L^p\ell^{\vq}}
\leq C
\bigg[\sup_{\Omega\in\Open(\R^{\vn})}\Big\|\sup_{\vj\in\mathbb Z^k}\one_{\Omega_{\vj}}\gamma_{\vj}\Big\|_{\aveL^s(\Omega)}\bigg]
\bNorm{\{ f_{\vj}\}_{\vj\in\mathbb Z^k}}{L^p\ell^{\vq}}.
\end{equation}
\end{lemma}

\begin{proof}
We first assume that $f_{\vj}=\Avee{f_{\vj}}{\vj}$ takes a constant value $f_R$ on each $R\in\D_{\vj}(\R^{\vn})$. Let us also write $\gamma_R:=\gamma_{\vj}|_R$ for $R\in\D_{\vj}(\R^{\vn})$, although this need not be a constant. If the supremum in the claim is zero or infinity, there is nothing to prove; otherwise we may assume by scaling that it is one. We also assume that $\Norm{\{ f_{\vj}\}}{L^p\ell^{\vq}}<\infty$, for otherwise there is nothing to prove.

For every $\kappa\in\Z$, let
\begin{equation*}
\Omega_\kappa:=\Big\{x\in\R^{\vn}:\
\Norm{\{f_{\vj}(x)\}_{\vj\in\mathbb Z^k}}{\ell^{\vq}}>2^\kappa\Big\}
\end{equation*}
and
\begin{equation*}
\Rs_\kappa:=\Big\{R\in\D(\R^{\vn}):\
\abs{R\cap\Omega_\kappa}>\frac12\abs{R}\Big\}.
\end{equation*}
We claim that, for each $R\in\mathscr D(\R^{\vn})$ with $f_R\neq 0$, there is a maximal $\kappa\in\Z$ such that $R\in\Rs_\kappa$: On the one hand, for all $x\in R$, we have $\Norm{\{f_{\vj}(x)\}_{\vj\in\mathbb Z^k}}{\ell^{\vq}}\geq\abs{f_R}$; hence $R\subset\Omega_\kappa$ and thus $R\in\Rs_\kappa$ for all $\kappa\in\mathbb Z$ with $2^\kappa<\abs{f_R}$. On the other hand,
if $R\in\Rs_\kappa$, then
\begin{equation*}
\infty>\Norm{\{f_{\vj}\}_{\vj\in\mathbb Z^k}}{L^{\vp}\ell^{\vq}}
\geq 2^\kappa\abs{\Omega_\kappa}^{\frac 1p}
\geq 2^\kappa\abs{\Omega_\kappa\cap R}^{\frac 1p}
\geq 2^\kappa(2^{-1}\abs{R})^{\frac 1p},
\end{equation*}
which cannot hold for arbitrarily large $\kappa$.
Thus, $R\in\Rs_\kappa$ for all small enough $\kappa$, and $R\notin\Rs_\kappa$ for all large enough $\kappa$. This confirms the existence of a maximal $\kappa$ with $R\in\Rs_\kappa$, and hence $R\in\Rs_\kappa\setminus\Rs_{\kappa+1}$
and $R\notin\Rs_{\widetilde\kappa}\setminus\Rs_{\widetilde\kappa+1}$
for all $\widetilde\kappa\in\mathbb Z\setminus\{\kappa\}$.
By this and Minkowski's inequality, we conclude that,
for every $x\in\mathbb R^{\vec n}$,
\begin{align*}
\Norm{\{\gamma_{\vj}(x)f_{\vj}(x)\}_{\vj\in\mathbb Z^k}}{\ell^{\vq}}
&=\Norm{\{\gamma_{R}(x)f_{R}\one_R(x)\}_{R\in\D(\R^{\vn})}}{\ell^{\vq}} \\
&=\Bigg\|\Bigg\{\sum_{\kappa\in\Z}\one_{\Rs_\kappa\setminus\Rs_{\kappa+1}}(R)
\gamma_{R}(x)f_{R}\one_R(x)\Bigg\}_{R\in\D(\R^{\vn})} \Bigg\|_{\ell^{\vq}} \\
&\leq\sum_{\kappa\in\Z}\BNorm{\Big\{\one_{\Rs_\kappa\setminus\Rs_{\kappa+1}}(R)
\gamma_{R}(x)f_{R}\one_R(x)\Big\}_{R\in\D(\R^{\vn})}}{\ell^{\vq}} \\
&\leq\sum_{\kappa\in\Z}\Big[\sup_{P\in\Rs_\kappa}\abs{\gamma_P(x)}\one_P(x)\Big]
\BNorm{\Big\{\one_{\Rs_\kappa\setminus\Rs_{\kappa+1}}(R)f_{R}\one_R(x)
\Big\}_{R\in\D(\R^{\vn})}}{\ell^{\vq}}.
\end{align*}
Taking $L^p$ norms and using Minkowski's inequality again
and then using H\"older's inequality with $\frac1p=\frac1s+\frac1r$ [where, since $s\in(p,\infty)$, also $r=\frac{s}{s-p}p\in(p,\infty)$], we obtain
\begin{align*}
\Norm{\{\gamma_{\vj}f_{\vj}\}_{\vj\in\mathbb Z^k}}{L^p\ell^{\vq}}
&\leq \sum_{\kappa\in\Z}\Big\| x\mapsto \Big[\sup_{P\in\Rs_\kappa}\abs{\gamma_P(x)}\one_P(x)\Big]
\BNorm{\Big\{\one_{\Rs_\kappa\setminus\Rs_{\kappa+1}}(R)f_{R}\one_R(x)
\Big\}_{R\in\D(\R^{\vn})}}{\ell^{\vq}}\Big\|_{L^p}\\
&\leq\sum_{\kappa\in\Z}
\BNorm{x\mapsto\sup_{P\in\Rs_\kappa}\abs{\gamma_P(x)}\one_P(x)}{L^s}
\BNorm{\Big\{x\mapsto \one_{\Rs_\kappa\setminus\Rs_{\kappa+1}}(R)f_{R}\one_R(x)\Big\}_{R\in\D(\R^{\vn})}}{L^r\ell^{\vq}} \\
&=:\sum_{\kappa\in\Z}I_\kappa\times II_\kappa.
\end{align*}

To proceed, we note that every $P\in\Rs_\kappa$ satisfies $\ave{\one_{\Omega_\kappa}}_{P}>\frac12$ by definition, and hence
\begin{equation*}
P\subset\Omega_\kappa^*
:=\Big\{x\in\R^{\vn}:\ \mathcal M_{\vn}\one_{\Omega_\kappa}(x)>\frac12\Big\}.
\end{equation*}
The lower semicontinuity of $\mathcal M_{\vn}\one_{\Omega_\kappa}$ implies that $\Omega_\kappa^*$ is open
and the strong maximal inequality (Lemma \ref{FS strong}) implies that
\begin{equation}\label{Omega*}
\abs{\Omega_\kappa^*}\lesssim\abs{\Omega_\kappa}<\infty.
\end{equation}
Hence
\begin{equation*}
\begin{split}
I_\kappa
  &\leq\bigg\|x\mapsto\sup_{P\subset\Omega_\kappa^*}\abs{\gamma_P(x)}\one_P(x)\bigg\|_{L^s}
  =\bigg\|\sup_{\vj\in\mathbb Z^k}\one_{(\Omega_\kappa^*)_{\vj}} \gamma_{\vj}\bigg\|_{L^s} \\
  &\leq\abs{\Omega_\kappa^*}^{\frac 1s}\quad\text{by \eqref{carl >1eq}, where we assumed that the supremum is $1$} \\
  &\lesssim\abs{\Omega_\kappa}^{\frac 1s}\quad\text{by \eqref{Omega*}}.
\end{split}
\end{equation*}

Concerning $II_\kappa$, for $R\notin\Rs_{\kappa+1}$,
we have $\ave{\one_{\Omega_{\kappa+1}^c}}_R\geq\frac12$ by definition,
while $R\in\Rs_{\kappa}$ implies $R\subset\Omega_{\kappa}^*$ as before. Hence
\begin{align*}
II_\kappa
&\leq\BNorm{\Big\{x\mapsto \one_{\Rs_\kappa\setminus\Rs_{\kappa+1}}(R)f_{R}\cdot 2\ave{\one_{\Omega_{\kappa+1}^c}}_R\one_R(x)\Big\}_{R\in\mathscr D(\mathbb R^{\vn})}}{L^r\ell^{\vq}}  \\
&\leq \BNorm{\Big\{\Avee{2\cdot f_{\vj}\one_{\Omega_\kappa^*\setminus\Omega_{\kappa+1}}}{\vj}\Big\}_{\vj\in\mathbb Z^k}}{L^r\ell^{\vq}}
\leq \BNorm{\Big\{\mathcal M_{\vn}(2\cdot f_{\vj}\one_{\Omega_\kappa^*\setminus\Omega_{\kappa+1}}) \Big\}_{\vj\in\mathbb Z^k}}{L^r\ell^{\vq}} \\
&\lesssim\BNorm{\Big\{ f_{\vj}\one_{\Omega_\kappa^*\setminus\Omega_{\kappa+1}} \Big\}_{\vj\in\mathbb Z^k}}{L^r\ell^{\vq}}
=\BNorm{x\mapsto \one_{\Omega_\kappa^*\setminus\Omega_{\kappa+1}}(x) \Norm{\{ f_{\vj}(x) \}_{\vj\in\mathbb Z^k} }{\ell^{\vq} } }{L^r}.
\end{align*}
For $x\notin\Omega_{\kappa+1}$, we have $ \Norm{\{ f_{\vj}(x) \} }{\ell^{\vq} }\leq 2^{\kappa+1}$ by definition, and hence
\begin{equation*}
II_{\kappa}
\lesssim \BNorm{x\mapsto \one_{\Omega_\kappa^*\setminus\Omega_{\kappa+1}}(x)  2^{\kappa+1} }{L^r}
\lesssim \abs{\Omega_\kappa^*}^{\frac 1r}2^\kappa
\lesssim \abs{\Omega_\kappa}^{\frac 1r}2^\kappa,
\end{equation*}
again by \eqref{Omega*} in the last step.

Combining the obtained estimates, we deduce that
\begin{equation}\label{carl 12}
\begin{split}
  \Norm{\{\gamma_{\vj}f_{\vj}\}_{\vj\in\mathbb Z^k}}{L^p\ell^{\vq}}
  &\lesssim\sum_{\kappa\in\Z}I_\kappa\times II_\kappa
  \lesssim\sum_{\kappa\in\Z}\abs{\Omega_\kappa}^{\frac 1s}\times\abs{\Omega_\kappa}^{\frac 1r} 2^\kappa
  =\sum_{\kappa\in\Z}\abs{\Omega_\kappa}^{\frac 1p} 2^\kappa \\
  &=\sum_{\kappa\in\Z}\abs{\{x\in\R^{\vn}:\
  \Norm{\{f_{\vj}(x)\}_{\vj\in\mathbb Z^k}}{\ell^{\vq}}>2^\kappa\}}^{\frac 1p} 2^\kappa
  \sim\Norm{\{f_{\vj}\}_{\vj\in\mathbb Z^k}}{L^{p,1}\ell^{\vq}},
\end{split}
\end{equation}
where $L^{p,1}$ is the Lorentz space.

The obtained estimate deviates from that claimed in the lemma in two ways: by having the Lorentz $L^{p,1}$ norm in place of $L^p$ on the right,
and by being only established under the additional assumption that $f_{\vj}=\Avee{f_{\vj} }{\vj}$. Both shortcomings may be lifted with the help of the interpolation result stated in Lemma \ref{interp}. Given $1<p<s<\infty$ as in the assumptions of the lemma that we are proving, we may clearly choose $p_0,p_1,\theta$ as in Lemma \ref{interp} so that in addition $1<p_0<p<p_1<s<\infty$. By the strong maximal inequality (Lemma \ref{FS strong}), we find that
\begin{equation*}
\Norm{\{\Avee{f_{\vj}}{\vj}\}_{\vj\in\mathbb Z^k}}{L^{p_i}\ell^{\vq}}
\leq \Norm{\{\mathcal M_{\vn}(f_{\vj})\}_{\vj\in\mathbb Z^k}}{L^{p_i}\ell^{\vq}}
\lesssim\Norm{\{f_{\vj}\}_{\vj\in\mathbb Z^k}}{L^{p_i}\ell^{\vq}},\quad i=0,1.
\end{equation*}
Applying Lemma \ref{interp} with $X=\ell^{\vq}$, $L^{p_i,r_i}=L^{p_i,p_i}=L^{p_i}$, and $L^{p,r}=L^{p,1}$, we deduce for the bounded linear operator $\{f_{\vj}\}_{\vj\in\mathbb Z^k}\mapsto\{\Avee{f_{\vj}}{\vj}\}_{\vj\in\mathbb Z^k}:\ L^{p_i}\ell^{\vq}\to L^{p_i}\ell^{\vq}$ the estimate
\begin{equation}\label{carl 13}
\begin{split}
\Norm{\{\Avee{f_{\vj}}{\vj}\}_{\vj\in\mathbb Z^k}}{L^{p,1}\ell^{\vq}}
&\sim\Norm{\{\Avee{f_{\vj}}{\vj}\}_{\vj\in\mathbb Z^k}}{(L^{p_0}\ell^{\vq},L^{p_1}\ell^{\vq})_{\theta,1}} \\
&\lesssim\Norm{\{f_{\vj}\}_{\vj\in\mathbb Z^k}}{(L^{p_0}\ell^{\vq},L^{p_1}\ell^{\vq})_{\theta,1}}
\sim\Norm{\{f_{\vj}\}_{\vj\in\mathbb Z^k}}{L^{p,1}\ell^{\vq}}.
\end{split}
\end{equation}

Recalling that \eqref{carl 12} was proved under the assumption that $f_{\vj}=\Avee{f_{\vj}}{\vj}$, if we now apply this bound to $\Avee{f_{\vj}}{\vj}$ in place of $f_{\vj}$, and combine this with \eqref{carl 13}, we obtain
\begin{equation}\label{carl 14}
\Norm{\{\gamma_{\vj}\Avee{f_{\vj}}{\vj}\}_{\vj\in\mathbb Z^k}}{L^p\ell^{\vq}}
\lesssim \Norm{\{\Avee{f_{\vj}}{\vj}\}_{\vj\in\mathbb Z^k}}{L^{p,1}\ell^{\vq}}
\lesssim \Norm{\{f_{\vj}\}_{\vj\in\mathbb Z^k}}{L^{p,1}\ell^{\vq}}.
\end{equation}
Since \eqref{carl 14} is valid for any exponents as in the assumptions, it is also valid for each $p_i$ in place of $p$:
\begin{equation*}
\Norm{\{\gamma_{\vj}\Avee{f_{\vj}}{\vj}\}_{\vj\in\mathbb Z^k}}{L^{p_i}\ell^{\vq}}
\lesssim \Norm{\{f_{\vj}\}_{\vj\in\mathbb Z^k}}{L^{p_i,1}\ell^{\vq}},\quad i=0,1.
\end{equation*}
Applying Lemma \ref{interp} with $X=\ell^{\vq}$, $r_i\in\{p_i,1\}$, and $r=p$, we deduce for the bounded linear operator $\{f_{\vj}\}_{\vj\in\mathbb Z^k}\mapsto\{\gamma_{\vj}\Avee{f_{\vj}}{\vj}\}_{\vj\in\mathbb Z^k}:\ L^{p_i,1}\ell^{\vq}\to L^{p_i}\ell^{\vq}$, the estimate
\begin{equation*}
\begin{split}
\Norm{\{\gamma_{\vj}\Avee{f_{\vj}}{\vj}\}_{\vj\in\mathbb Z^k}}{L^{p}\ell^{\vq}}
&\sim\Norm{\{\gamma_{\vj}\Avee{f_{\vj}}{\vj}\}_{\vj\in\mathbb Z^k}}{(L^{p_0}\ell^{\vq},L^{p_1}\ell^{\vq})_{\theta,p}} \\
&\lesssim\Norm{\{\Avee{f_{\vj}}{\vj}\}_{\vj\in\mathbb Z^k}}{(L^{p_0,1}\ell^{\vq},L^{p_1,1}\ell^{\vq})_{\theta,p}}
\sim \Norm{\{f_{\vj}\}_{\vj\in\mathbb Z^k}}{L^p\ell^{\vq}},
\end{split}
\end{equation*}
which is the assertion of Lemma \ref{carl >1}.
\end{proof}

\begin{remark} \label{10.31}
In the one-parameter case of Lemma \ref{carl >1},
the bound
$$
\mathrm{I}:=\sup_{\Omega\in\Open(\R^{n})}\bigg\|\sup_{j\in\mathbb Z}
\one_{\Omega_{j}} \gamma_j \bigg\|_{\aveL^s(\Omega)}
$$
can be replaced by
$$
\mathrm{II}:=\sup_{Q\in\mathscr D(\mathbb R^n)}
\bigg\|\sup_{j\in\mathbb Z,\,j\geq j_Q} \gamma_j\bigg\|_{\aveL^s(Q)};
$$
see \cite[Theorem 3.7(ii)]{FR21}.
Indeed, $\mathrm{I}=\mathrm{II}$ for all $s\in(0,\infty)$.
The inequality $\mathrm{II}\leq\mathrm{I}$ is obvious.
Denoting by $\mathscr D^*(\Omega)$ the collection of maximal dyadic cubes
$P\subset\Omega$, these cubes are disjoint, and we have
$$
\bigg\|\sup_{j\in\mathbb Z}
\one_{\Omega_{j}} \gamma_j \bigg\|_{\aveL^s(\Omega)}^s
= \frac{1}{|\Omega|} \sum_{P\in \mathscr D^*(\Omega)}
\int_P \sup_{j\in\mathbb Z,j\geq j_P} |\gamma_j|^s
\leq \frac{1}{|\Omega|} \sum_{P\in \mathscr D^*(\Omega)} |P| \cdot\mathrm{II}^s
=\mathrm{II}^s
$$
and hence $\mathrm{I}\leq\mathrm{II}$.
Thus, $\mathrm{I}=\mathrm{II}$ and Lemma \ref{carl >1} in the case $k=1$
coincides with \cite[Theorem 3.7(ii)]{FR21} in the case $q\in(1,\infty)$.
\end{remark}

Finally, by specialising to sequences $f_{\vj}$ that are constant on each $R\in\D_{\vj}(\R^{\vn})$, we can push the lower bound of the exponents from one to zero.

\begin{lemma}\label{carl last}
If $0<p<s<\infty$ and $\vq\in(0,\infty)^k$,
then there exists a positive constant $C$ such that,
for all measurable functions $\{\gamma_{\vj}\}_{\vj\in\mathbb Z^k}$ and
every sequence $\{f_{\vj}\}_{\vj\in\mathbb Z^k}\in L^p\ell^{\vq}$
with each $f_{\vj}$ constant on every $R\in\D_{\vj}(\R^{\vn})$,
\begin{equation*}
\bNorm{\{\gamma_{\vj}f_{\vj}\}_{\vj\in\mathbb Z^k}}{L^p\ell^{\vq}}
\lesssim \bigg[\sup_{\Omega\in\Open(\R^{\vn})}\Big\|\sup_{\vj\in\mathbb Z^k}\one_{\Omega_{\vj}}\gamma_{\vj}\Big\|_{\aveL^s(\Omega)}\bigg]
\bNorm{\{ f_{\vj}\}_{\vj\in\mathbb Z^k}}{L^p\ell^{\vq}}.
\end{equation*}
\end{lemma}

\begin{proof}
Pick some $r\in(0,\min(p,\min\vq))$ so that $1<p/r<s/r<\infty$ and $\vone<\vq/r<\vinfty$. An application of Lemma \ref{carl >1} with these rescaled exponents then shows that
\begin{align*}
\bNorm{\{\gamma_{\vj}f_{\vj}\}_{\vj\in\mathbb Z^k}}{L^p\ell^{\vq}}^r
&=\bNorm{\{\abs{\gamma_{\vj}}^r\abs{f_{\vj}}^r\}_{\vj\in\mathbb Z^k}}{L^{p/r}\ell^{\vq/r}}  \\
&=\bNorm{\{\abs{\gamma_{\vj}}^r\ave{\abs{f_{\vj}}^r}_{\vj}\}_{\vj\in\mathbb Z^k}}{L^{p/r}\ell^{\vq/r}}
\quad\text{since $f_{\vj}$
is constant on $R\in\mathscr D_{\vj}(\R^{\vn})$}, \\
&\lesssim\sup_{\Omega\in\Open(\R^{\vn})}
\BNorm{\sup_{\vj\in\mathbb Z^k}\one_{\Omega_{\vj}}\abs{\gamma_{\vj}}^r}{\aveL^{s/r}(\Omega)}
\bNorm{\{\abs{f_{\vj}}^r\}_{\vj\in\mathbb Z^k}}{L^{p/r}\ell^{\vq/r}}
\quad\text{by Lemma \ref{carl >1}} \\
&=\sup_{\Omega\in\Open(\R^{\vn})}\BNorm{\sup_{\vj\in\mathbb Z^k}\one_{\Omega_{\vj}}\gamma_{\vj}}{\aveL^{s}(\Omega)}^r
\bNorm{\{ f_{\vj} \}_{\vj\in\mathbb Z^k}}{L^{p}\ell^{\vq}}^r,
\end{align*}
which is the claim of the present lemma.
\end{proof}

Let us summarise the findings of this subsection in the following:

\begin{theorem}\label{carl thm}
Let $0<p<s<\infty$, $\vp:=p\cdot\vone$, and $\vq\in(0,\infty)^k$.
Let $\Xs$ be as in one of the cases \eqref{carl B} through \eqref{carl F} below.
Then, for all measurable functions $\{\gamma_{\vj}\}_{\vj\in\mathbb Z^k}$ and
every $\{f_{\vj}\}_{\vj\in\mathbb Z^k}\in \Xs$
with each $f_{\vj}$ constant on every $R\in\D_{\vj}(\R^{\vn})$,
\begin{equation} \label{finally}
\bNorm{\{\gamma_{\vj}f_{\vj}\}_{\vj\in\mathbb Z^k}}{\Xs}
\lesssim \bNorm{\{ f_{\vj}\}_{\vj\in\mathbb Z^k}}{\Xs}
\end{equation}
holds, with the implicit positive constants independent of $f_{\vj}$ and $\gamma_{\vj}$, in each of the following cases:
\begin{enumerate}[\rm(i)]
\item\label{carl B} $\sup_{\vj\in\mathbb Z^k}\sup_{R\in\D_{\vj}(\R^{\vn})}\Norm{\gamma_{\vj}}{\aveL^p(R)}\lesssim 1$ and $\Xs=\ell^{\vq}L^p$ is of Besov-type,

\item\label{carl arb} $\sup_{\vj\in\mathbb Z^k}\sup_{R\in\D_{\vj}(\R^{\vn})}\Norm{\gamma_{\vj}}{\aveL^s(R)}\lesssim 1$ and $\Xs=(L^{\vp}\ell^{\vq})_\pi$ is of arbitrary type with $\vq\leq\vp$,

\item\label{carl F} $\sup_{\Omega\in\Open(\R^{\vn})}\Norm{\sup_{\vj\in\mathbb Z^k}\one_{\Omega_{\vj}}\gamma_{\vj}}{\aveL^s(\Omega)}\lesssim 1$ and $\Xs=L^p\ell^{\vq}$ is of Triebel--Lizorkin-type.
\end{enumerate}
\end{theorem}

\begin{proof}
The inequality \eqref{finally} in cases \eqref{carl B} through \eqref{carl F} follows from
Proposition \ref{carl easy}, Lemma \ref{carl small q},
and Lemma \ref{carl last}, respectively.
\end{proof}

As we saw in Proposition \ref{carl easy}, the condition of case \eqref{carl B} of Theorem \ref{carl thm} is not only sufficient but also necessary for the conclusion.
Also the assumptions on $\gamma_{\vj}$ in cases \eqref{carl arb} and \eqref{carl F} of Theorem \ref{carl thm} are sharp, in a sense to be made precise in Propositions \ref{carl small q2}, \ref{carl last 2}, and \ref{carl last 2x} below. Moreover, the assumption $\vq\leq\vp$ in \eqref{carl arb} is necessary, as we will show in Propositions \ref{counterexample-gamma} and \ref{counterexamplecom}.

\subsection{Sharpness and necessity of the Carleson conditions} \label{Carl sharp}

In this subsection, we establish the sharpness of the assumptions of Theorem \ref{carl thm}.
We first prove the sharpness of the assumptions on $\gamma_{\vj}$ in cases
\eqref{carl arb} and \eqref{carl F} of Theorem \ref{carl thm},
with case \eqref{carl arb} handled in the lemma below.

\begin{proposition}\label{carl small q2}
Let $0<q<p<s<\infty$.
Then there exists a positive constant $C$ such that,
for all measurable functions $\{\gamma_j\}_{j\in\mathbb Z}$ and
all $\{f_j\}_{j\in\mathbb Z}\in L^p\ell^q$,
where each $f_{j}$ constant on every $R\in\D_{j}(\mathbb R^n)$,
\begin{equation}\label{FS-alter 2}
\bNorm{\{\gamma_jf_j\}_{j\in\mathbb Z}}{L^p\ell^q}
\le C\bigg[\sup_{j\in\Z}\sup_{R\in\D_j(\mathbb R^n)}
\Norm{\gamma_j}{\aveL^s(R)}\bigg]
\bNorm{\{ f_j\}_{j\in\mathbb Z}}{L^p\ell^q},
\end{equation}
where this inequality is sharp in the sense that
\eqref{FS-alter 2} no longer holds when $s = p$.
\end{proposition}

\begin{proof}
Inequality \eqref{FS-alter 2} is a special case of Lemma \ref{carl small q}.
Now, we prove the sharpness.
For any $j\in\mathbb Z$, let
\begin{align}\label{counterexample2}
\gamma_{-j}:=\begin{cases}
2^{\frac{jn}{p}}\mathbf 1_{[0,1)^n},&
\text{if } j\geq 0,\\
0,&\text{otherwise}.
\end{cases}
\end{align}
Then
$\sup_{j\in\Z}\sup_{R\in\D_j}
\Norm{\gamma_j}{\aveL^p(R)}=1$.
Next, we show that
\begin{align}\label{counterexample}
\sup_{\|\{ f_j\}_{j\in\mathbb Z}\|_{L^p\ell^q}\ne 0}
\frac{\|\{\gamma_jf_j\}_{j\in\mathbb Z}\|_{L^p\ell^q}}
{\|\{ f_j\}_{j\in\mathbb Z}\|_{L^p\ell^q}}
=\infty.
\end{align}
For any $N\in\mathbb N$ and $j\in\mathbb Z$, let
$$
f_{-j}^{(N)}:=\begin{cases}
2^{-\frac{jn}{p}}\one_{[0,2^{j})^n},&
\text{if } j\in[N],\\
0,&\text{otherwise}.
\end{cases}
$$
Then $\|\{\gamma_jf_j^{(N)}\}_{j\in\mathbb Z}\|_{L^p\ell^q}=N^{\frac1q}$ and
\begin{align*}
\Big\|\Big\{ f_j^{(N)}\Big\}_{j\in\mathbb Z}\Big\|_{L^p\ell^q}^p
&= \int_{[0,1)^n} \Bigg\{\sum_{j=1}^{N}\Big[f_{-j}^{(N)}\Big]^q\Bigg\}^{\frac pq}
+\sum_{l=1}^{N}\int_{[0,2^{l})^n\setminus[0,2^{l-1})^n} \cdots\\
&=\Bigg(\sum_{j=1}^{N}2^{-jn\frac qp}\Bigg)^{\frac pq}
+\sum_{l=1}^{N} \Big[2^{ln}-2^{(l-1)n}\Big]
\Bigg(\sum_{j=l}^N 2^{-jn\frac qp}\Bigg)^{\frac pq}\\
&\sim 1+\sum_{l=1}^{N} 2^{ln}2^{-ln}\sim N.
\end{align*}
Therefore,
$$
\frac{\|\{\gamma_jf_j^{(N)}\}_{j\in\mathbb Z}\|_{L^p\ell^q}}
{\|\{ f_j^{(N)}\}_{j\in\mathbb Z}\|_{L^p\ell^q}}
\gtrsim N^{\frac1q-\frac1p}\to\infty
\quad\text{as}\quad N\to\infty.
$$
This finishes the proof of \eqref{counterexample}.
Hence \eqref{FS-alter 2} no longer holds when $s = p$,
which completes the proof of Proposition \ref{carl small q2}.
\end{proof}

In Theorem \ref{carl thm}(iii), the essential range of $p$, $q$, and $s$
is $0<p<s<q<\infty$ because
\begin{itemize}
  \item the case $q\leq p$ is already covered by Theorem \ref{carl thm}(ii), and
  \item when both $p<q$ and $p<s$, we may assume that $p<s<q$, since the smaller the~$s$, the better the bound.
\end{itemize}

\begin{proposition}\label{carl last 2}
If $0<p<s<q<\infty$,
then there exists a positive constant $C$ such that,
for all measurable functions $\{\gamma_j\}_{j\in\mathbb Z}$ and
every sequence $\{f_j\}_{j\in\mathbb Z}\in L^p\ell^q$
with each $f_j$ constant on every $R\in\D_j(\mathbb R^n)$,
\begin{equation}\label{FS-alter}
\bNorm{\{\gamma_jf_j\}_{j\in\mathbb Z}}{L^p\ell^q}
\le C\bigg[\sup_{\Omega\in\Open(\R^{\vn})}\Big\|\sup_{j\in\mathbb Z}\one_{\Omega_j}\gamma_j\Big\|_{\aveL^s(\Omega)}\bigg]
\bNorm{\{ f_j\}_{j\in\mathbb Z}}{L^p\ell^q}.
\end{equation}
Moreover, this inequality is sharp in the sense that
\eqref{FS-alter} no longer holds if
the bound
$$\sup_{\Omega\in\Open(\R^{\vn})}\Big\|\sup_{j\in\mathbb Z}
\one_{\Omega_j}\gamma_j\Big\|_{\aveL^s(\Omega)}$$
is replaced by either
$$
\sup_{j\in\Z}\sup_{R\in\D_j(\R^{\vn})}
\Norm{\gamma_j}{\aveL^s(R)}
\quad\text{or}\quad
\sup_{\Omega\in\Open(\R^{\vn})}\Big\|\sup_{j\in\mathbb Z}\one_{\Omega_j}\gamma_j\Big\|_{\aveL^p(\Omega)}.
$$
\end{proposition}

\begin{proof}
Inequality \eqref{FS-alter} is a special case of Lemma \ref{carl last}.
Now, we prove the sharpness.
We first show that
\begin{equation}\label{FS-alter-2}
\bNorm{\{\gamma_jf_j\}_{j\in\mathbb Z}}{L^p\ell^q}
\lesssim \Big[\sup_{j\in\Z}\sup_{R\in\D_j(\R^{\vn})}
\Norm{\gamma_j}{\aveL^s(R)}\Big]
\bNorm{\{ f_j\}_{j\in\mathbb Z}}{L^p\ell^q}
\end{equation}
fails to hold. For any $N\in\mathbb N$ and $j\in\mathbb Z$, let
\begin{align*}
\gamma_j^{(N)}:=\begin{cases}
N^{\frac1s} \mathbf 1_{F_j},&
\text{if } j\in[N],\\
0,&\text{otherwise},
\end{cases}
\end{align*}
where $\{F_j\}_{j=1}^N$ are pairwise disjoint subsets
of $[0,1)^n$ and,
for every $j\in[N]$,
we have $|F_j|=\frac1N$ and
$\{F_j\cap Q\}_{Q\in \mathscr{D}_j([0,1)^n)}$
have equal measure.
Let $e:=(1,0,\ldots,0)\in\mathbb Z^n$ and, for each $j\in\Z$,
$$
\gamma_j(\cdot):=\sum_{N=1}^\infty \gamma_j^{(N)}(\cdot-Ne).
$$
Then $\sup_{j\in\Z}\sup_{R\in\D_j}
\Norm{\gamma_j}{\aveL^s(R)}=1$.
For each $N\in\Z_+$ and $j\in\mathbb Z$, let
$$
f_j^{(N)}:=\begin{cases}
\mathbf 1_{[0,1)^n+Ne},&
\text{if } j\in[N],\\
0,&\text{otherwise}.
\end{cases}
$$
Then
\begin{equation*}
\begin{split}
  \|\{\gamma_jf_j^{(N)}\}_{j\in\mathbb Z}\|_{L^p\ell^q}
  &=\|\{\gamma_j^{(N)}(\cdot-Ne)f_j^{(N)}\}_{j\in\mathbb Z}\|_{L^p\ell^q} \\
  &=\|\{ N^{\frac1s} \one_{F_j+Ne}\}_{j\in[N]}\|_{L^p\ell^q}
  =\| N^{\frac1s} \one_{[0,1)^n+Ne}\|_{L^p}  =N^{\frac1s}
\end{split}
\end{equation*}
and $\|\{f_j^{(N)}\}_{j\in\mathbb Z}\|_{L^p\ell^q}=\|\{ \one_{[0,1)^n+Ne}\}_{j\in[N]}\|_{L^p\ell^q}=N^{\frac1q}$.

Therefore,
$$
\frac{\|\{\gamma_jf_j^{(N)}\}_{j\in\mathbb Z}\|_{L^p\ell^q}}
{\|\{ f_j^{(N)}\}_{j\in\mathbb Z}\|_{L^p\ell^q}}
= N^{\frac1s-\frac1q}\to\infty
\quad\text{as}\quad N\to\infty
$$
and hence \eqref{FS-alter-2} fails to hold.

Next, we prove that
\begin{equation}\label{FS-alter-3}
\bNorm{\{\gamma_jf_j\}_{j\in\mathbb Z}}{L^p\ell^q}
\lesssim \Bigg[\sup_{\Omega\in\Open(\R^{\vn})}\Big\|\sup_{j\in\mathbb Z}\one_{\Omega_j}\gamma_j\Big\|_{\aveL^p(\Omega)}\Bigg]
\bNorm{\{ f_j\}_{j\in\mathbb Z}}{L^p\ell^q}
\end{equation}
fails to hold.
Let $\{\gamma_j\}_{j\in\Z}$ be as in \eqref{counterexample2}.
Thus, $\one_{\Omega_{-j}}\gamma_{-j}=\one_{\N}(j)2^{\frac{jn}{p}}\one_{\Omega_{-j}}\one_{[0,1)^n}$. For $j\geq 0$, the set $\Omega_{-j}$ is a union of the translates of $[0,2^j)^n\supseteq[0,1)^n$ contained in $\Omega$. Hence it follows that
\begin{equation*}
  \one_{\Omega_{-j}}\gamma_{-j}=
  \begin{cases} 2^{\frac{jn}{p}}\one_{[0,1)^n}, & \text{if }j\geq 0\text{ and }[0,2^j)^n\subset\Omega, \\ 0, & \text{otherwise}.\end{cases}
\end{equation*}
Note that a necessary condition for $[0,2^j)^n\subseteq\Omega$ is $2^{jn}\leq\abs{\Omega}$. Hence
\begin{equation*}
  \sup_{j\in\Z}\one_{\Omega_j}\gamma_j
  \leq\sup_{j: 1\leq 2^{jn}\leq\abs{\Omega}} 2^{\frac{jn}{p}}\one_{[0,1)^n}
  \lesssim \abs{\Omega}^{\frac1p}\one_{[0,1)^n},
\end{equation*}
and thus
\begin{equation*}
\BNorm{\sup_{j\in\Z}\one_{\Omega_j}\gamma_j}{\aveL^p(\Omega)}
\lesssim\abs{\Omega}^{\frac1p}\bNorm{\one_{[0,1)^n}}{\aveL^p(\Omega)}
=\abs{\Omega}^{\frac1p}\abs{\Omega}^{-\frac1p}\bNorm{\one_{[0,1)^n}}{L^p(\Omega)}=
1.
\end{equation*}
Hence \eqref{FS-alter-3} would imply that
$\Norm{\{\gamma_jf_j\}_{j\in\mathbb Z}}{L^p\ell^q}
\lesssim \Norm{\{ f_j\}_{j\in\mathbb Z}}{L^p\ell^q}$.
But we know from \eqref{counterexample} that this is not the case, and hence \eqref{FS-alter-3} fails to hold.
This finishes the proof of Proposition \ref{carl last 2}.
\end{proof}

In the one-parameter case of Theorem \ref{carl thm}\eqref{carl F},
the bound
$$
\sup_{\Omega\in\Open(\R^n)}\Big\|\sup_{j\in\mathbb Z}\one_{\Omega_j}\gamma_j\Big\|_{\aveL^s(\Omega)}
$$
can be equivalently replaced by
$\sup_{Q\in\mathscr D(\mathbb R^n)}
\|\sup_{j\in\mathbb Z,\,j\geq j_Q} \gamma_j\|_{\aveL^s(Q)}$;
see Remark \ref{10.31}.
However, this replacement is not valid in the multi-parameter case.

\begin{proposition}\label{carl last 2x}
If $0<p<s<q<\infty$,
then there exists a positive constant $C$ such that,
for all measurable functions $\{\gamma_{\vj}\}_{\vj\in\mathbb Z^k}$ and
every sequence $\{f_{\vj}\}_{\vj\in\mathbb Z^k}\in L^p\ell^q$
with each $f_{\vj}$ constant on every $R\in\D_{\vj}(\mathbb R^{\vec n})$,
\begin{equation}\label{FS-alterx}
\bNorm{\{\gamma_{\vj}f_{\vj}\}_{\vj\in\mathbb Z^k}}{L^p\ell^q}
\le C\sup_{\Omega\in\Open(\R^{\vn})}\Big\|\sup_{\vj\in\mathbb Z^k}\one_{\Omega_{\vj}}\gamma_{\vj}\Big\|_{\aveL^s(\Omega)}
\bNorm{\{ f_{\vj}\}_{\vj\in\mathbb Z}}{L^p\ell^q}.
\end{equation}
Moreover, this inequality is sharp in the sense that
\eqref{FS-alterx} no longer holds if $k\geq 2$ and
the bound $\sup_{\Omega\in\Open(\R^{\vn})}\|\sup_{\vj\in\mathbb Z^k}\one_{\Omega_{\vj}}\gamma_{\vj}\|_{\aveL^s(\Omega)}$
is replaced by
$$
\sup_{P\in \D(\mathbb R^{\vec n})}
\Big\|\sup_{\vj\in\mathbb Z^k,\,\vj\geq \vj_P}
\gamma_{\vj}\Big\|_{\aveL^s(P)}.
$$
\end{proposition}

\begin{proof}
Inequality \eqref{FS-alterx} is a special case of Lemma \ref{carl last}.
Next, we show the sharpness, that is,
\begin{equation}\label{FS-alter-2x}
\bNorm{\{\gamma_{\vj}f_{\vj}\}_{\vj\in\mathbb Z^k}}{L^p\ell^q}
\lesssim \sup_{P\in \D(\mathbb R^{\vec n})}
\Big\|\sup_{\vj\in\mathbb Z^k,\,\vj\geq \vj_P}
\gamma_{\vj}\Big\|_{\aveL^s(P)}
\bNorm{\{ f_{\vj}\}_{\vj\in\mathbb Z^k}}{L^p\ell^q}
\end{equation}
fails to hold.

By an elaboration of Carleson's counterexample \cite{Car} given in Lemma \ref{Carleson-p}, we find that, for every $N\in\mathbb Z_+$,
there exists a finite collection $\{R_i\}_{i\in \mathscr I_N}$ of
pairwise distinct dyadic rectangles in $[0,1)^{\vec n}$ such that
\begin{enumerate}[\rm(i)]
\item\label{it:Car-i} $\sum_{i\in \mathscr I_N} |R_i|=1$,

\item\label{it:Car-ii} $\abs{\bigcup_{i\in \mathscr I_N} R_i}<\frac1N$, and

\item\label{it:Car-iii} $\sum_{i\in \mathscr I_N,\,R_i\subset P} |R_i| \leq |P|$ for every dyadic rectangle $P\subset[0,1)^{\vec n}$,

\item\label{it:Car-iv} $\abs{\bigcup_{i\in \mathscr I_N} R_i}^{-\frac sp}\Norm{\sum_{i\in\mathscr I_N}\one_{R_i}}{\frac ps}>N$.
\end{enumerate}

We note that properties \eqref{it:Car-i} through \eqref{it:Car-iii} consist of Carleson's original counterexample (or more precisely, its simple extension from two to higher dimensions, as given in Lemma \ref{Carleson}). If $p=s$, property \eqref{it:Car-iv} would be immediate from \eqref{it:Car-i} and \eqref{it:Car-ii}. However, now that $p<s$, we see from H\"older's inequality that property \eqref{it:Car-iv} is more difficult to satisfy, which requires an elaboration of Carleson's counterexample given in Lemma \ref{Carleson-p}, based on Proposition \ref{prop:CarlNew}.

Let $g^{(N)}:=\sum_{i\in \mathscr I_N} \mathbf 1_{R_i}$.
For every $N,\kappa\in\mathbb Z_+$, define
$
\Omega^{(N,\kappa)}:=\{x\in\mathbb R^{\vec n}:\ g^{(N)}=\kappa\}.
$
Then $$\bigcup_{\kappa\in\mathbb Z_+} \Omega^{(N,\kappa)}
= \bigcup_{i\in \mathscr I_N} R_i =:\Omega^{(N)}.$$
For every $N,\kappa\in\mathbb Z_+$ and $i\in \mathscr I_N$, let
\begin{equation*}
  E^{(N)}_{i,\kappa}\subset R_i\cap\Omega^{(N,\kappa)}
\end{equation*}
be sets such that $\{E_{i,\kappa}^{(N)}\}_{i\in \mathscr I_N,\,\kappa\in\mathbb Z_+}$ forms a partition of $\Omega^{(N)}$ and
\begin{equation*}
  \abs{E^{(N)}_{i,\kappa}}=\frac{1}{\kappa}\abs{R_i\cap\Omega^{(N,\kappa)}};
\end{equation*}
such sets can be chosen, since each point of $\Omega^{(N,\kappa)}$ belongs to exactly $\kappa$ different $R_i$. (Indeed, for every non-empty intersection of the form $R_{i_1}\cap\cdots\cap R_{i_\kappa}\cap\Omega^{(N,\kappa)}$, we divide it to $\kappa$ equal parts and give one of the parts to each $E_{i_j}^{(N,\kappa)}$; each final $E_i^{(N,\kappa)}$ will be the union of all such parts when we go through all such intersections, where $R_i$ is included.)
Thus
\begin{equation}\label{eq:gN-ENk}
  g^{(N)}=\sum_{\kappa=1}^\infty\kappa\one_{\Omega^{(N,\kappa)}}
  =\sum_{\kappa\in\Z_+}\sum_{i\in\mathscr I_N}\kappa\one_{E^{(N,\kappa)}_i}.
\end{equation}

Then, let
\begin{equation*}
\gamma_{R_i}^{(N)}:= \sum_{\kappa\in\mathbb Z_+}
\kappa^{\frac1s} \mathbf 1_{E_{i,\kappa}^{(N)}},
\end{equation*}
and, for every $\vj\in\mathbb Z^k$,
$$
\gamma_{\vj}^{(N)}:= \sum_{i\in \mathscr I_N,\,R_i\in \D_{\vj}(\R^{\vn})}  \gamma_{R_i}^{(N)}.
$$
Let $e:=(1,0,\ldots,0)\in\mathbb Z^{\vec n}$ and
$$
\gamma_{\vj}(\cdot):=\sum_{N=1}^\infty \gamma_{\vj}^{(N)}(\cdot-Ne).
$$
Then, for every $P\in \D(\mathbb R^{\vec n})$,
\begin{equation*}
\begin{split}
\Bigg\|\sup_{\vj\in\mathbb Z^k,\,\vj\geq \vj_P}
\gamma_{\vj}\Bigg\|_{\aveL^s(P)}^s
&= \fint_P \sum_{N=1}^\infty \sum_{i\in \mathscr I_N,\,\vj_{R_i}\geq \vj_P}
\sum_{\kappa\in\mathbb Z_+}
\kappa \mathbf 1_{E_{i,\kappa}^{(N)}}(\cdot-Ne)\\
&= \frac{1}{|P|} \sum_{N=1}^\infty \sum_{i\in \mathscr I_N,\,\vj_{R_i}\geq \vj_P}
\sum_{\kappa\in\mathbb Z_+}
\kappa |E_{i,\kappa}^{(N)} \cap (P-Ne)| \\
&= \frac{1}{|P|} \sum_{N=1}^\infty \sum_{i\in \mathscr I_N,\,R_i \subset (P-Ne)}
\sum_{\kappa\in\mathbb Z_+} |R_i\cap\Omega^{(N,\kappa)}| \\
&= \frac{1}{|P|} \sum_{N=1}^\infty \sum_{i\in \mathscr I_N,\,R_i \subset (P-Ne)} |R_i|
\leq \frac{1}{|P|} \sum_{N=1}^\infty |(P-Ne) \cap [0,1)^{\vec n}|
\leq 1.
\end{split}
\end{equation*}
For every $N\in\mathbb Z_+$ and $\vj\in\mathbb Z^k$, let
$$
f_{\vj}^{(N)}:=\sum_{i\in \mathscr I_N,\,R_i\in \D_{\vj}(\mathbb R^{\vec n})} \mathbf 1_{R_i+Ne}.
$$
If \eqref{FS-alter-2x} holds, then, for every $N\in\mathbb Z_+$,
$$
\bNorm{\{\gamma_{\vj}f_{\vj}^{(N)}\}_{\vj\in\mathbb Z^k}}{L^p\ell^q}
\lesssim \bNorm{\{ f_{\vj}^{(N)}\}_{\vj\in\mathbb Z^k}}{L^p\ell^q}.
$$
However,
\begin{equation*}
\begin{split}
\bNorm{\{\gamma_{\vj}f_{\vj}^{(N)}\}_{\vj\in\mathbb Z^k}}{L^p\ell^q}^p
&=\|\{\gamma_{\vj}^{(N)}(\cdot-Ne)f_{\vj}^{(N)}\}_{\vj\in\mathbb Z^k}\|_{L^p\ell^q}^p \\
&=\Bigg\| \sum_{i\in \mathscr I_N} \sum_{\kappa\in\mathbb Z_+}
\kappa^{\frac 1s} \mathbf 1_{E_{i,\kappa}^{(N)}}(\cdot-Ne)
\Bigg\|_{L^p(\mathbb R^{\vec n})}^p
= \sum_{i\in \mathscr I_N} \sum_{\kappa\in\mathbb Z_+}
\kappa^{\frac ps} |E_{i,\kappa}^{(N)}|
=\Norm{g^{(N)}}{L^{\frac ps}(\R^{\vn})}^{\frac ps}
\end{split}
\end{equation*}
by \eqref{eq:gN-ENk} in the last step.
From H\"older's inequality, we infer that
\begin{equation*}
\begin{split}
\bNorm{\{ f_{\vj}^{(N)}\}_{\vj\in\mathbb Z^k}}{L^p\ell^q}^p
=\Bigg\| \Bigg(\sum_{i\in \mathscr I_N} \mathbf 1_{R_i+Ne}\Bigg)^{\frac1q} \Bigg\|_{L^p(\mathbb R^{\vec n})}^p
= \int_{\mathbb R^{\vec n}} \mathbf 1_{\Omega^{(N)}}
[g^{(N)}]^{\frac pq}
&\leq |\Omega^{(N)}|^{\frac1{(q/s)'}} \Norm{[g^{(N)}]^{\frac pq}}{L^{\frac qs}(\R^{\vn})} \\
&= |\Omega^{(N)}|^{1-\frac sq} \Norm{g^{(N)}}{L^{\frac ps}(\R^{\vn})}^{\frac pq}.
\end{split}
\end{equation*}
Therefore,
\begin{equation*}
\frac{\|\{\gamma_{\vj}f_{\vj}^{(N)}\}_{\vj\in\mathbb Z^k}\|_{L^p\ell^q}^p}
{\|\{f_{\vj}^{(N)}\}_{\vj\in\mathbb Z^k}\|_{L^p\ell^q}^p}
\geq
\abs{\Omega^{(N)}}^{\frac sq-1}\Norm{g^{(N)}}{L^{\frac ps}(\R^{\vn})}^{\frac ps-\frac pq}
=
\Big(\abs{\Omega^{(N)}}^{-\frac sp}\Norm{g^{(N)}}{L^{\frac ps}(\R^{\vn})}
\Big)^{\frac ps-\frac pq}
>N^{\frac ps-\frac pq}
\end{equation*}
using, in the last step, \eqref{it:Car-iv} and the fact that the exponent is positive, since $s<q$.
As $N\to\infty$, we see that \eqref{FS-alter-2x} fails to hold.
This finishes the proof of Proposition \ref{carl last 2x}.
\end{proof}

The following result shows that if the assumption $\vq\leq\vp$ in Theorem \ref{carl thm}\eqref{carl arb} is removed,
then even under a stronger assumption on $\gamma_j$,
the desired inequality
\begin{equation*}
\bNorm{\{\gamma_{\vj}f_{\vj}\}_{\vj\in\mathbb Z^k}}{L^p\ell^q}
\lesssim \bNorm{\{ f_{\vj}\}_{\vj\in\mathbb Z^k}}{L^p\ell^q}
\end{equation*}
fails to hold.

\begin{proposition}\label{counterexample-gamma}
Let $p\in(0,\infty)$, $q_1\in(0,\infty)$, $q_2\in(p,\infty)$, and $s\in(p,q_2)$.
Then there exist measurable functions
$\{\gamma_{\vj}\}_{\vj\in\mathbb Z^2}$ satisfying
\begin{align}\label{supomega}
\sup_{\Omega\in\Open(\R^{\vn})}\Bigg\|\sup_{\vj\in\mathbb Z^2}
\one_{\Omega_{\vj}} \gamma_{\vj} \Bigg\|_{\aveL^s(\Omega)}\le1
\end{align}
and $\{f_{\vj}\}_{\vj\in\mathbb Z^2}\in \ell^{q_1}L_x^{p}\ell^{q_2}L_y^p$
with each $f_{\vj}$ constant on every $R\in\D_{\vj}(\R^{\vn})$ such that
\begin{equation}\label{counterexample4}
\bNorm{\{\gamma_{\vj}f_{\vj}\}_{\vj\in\mathbb Z^2}}{\ell^{q_1}L_x^{p}\ell^{q_2}L_y^p}=\infty
\quad\text{and}\quad
\bNorm{\{ f_{\vj}\}_{\vj\in\mathbb Z^2}}{\ell^{q_1}L_x^{p}\ell^{q_2}L_y^p}<\infty.
\end{equation}
\end{proposition}

\begin{proof}
For any $j\in\Z_+$ and $t\in(0,\infty)$,
let $S_t(j):=\sum_{i=1}^ji^{-t}$.
Assume that $t\in(1,\frac{pq_2^{-1}-1}{ps^{-1}-1})$
and $d\in(pq_2^{-1},tps^{-1}-t+1)$.
For any $x,y\in\mathbb R$, let
\begin{align*}
f_{\vj}(x,y):=\begin{cases}
j_{2}^{-\frac dp}
\mathbf1_{(j_2,j_2+1]}(y)
\mathbf{1}_{(0,1]}(x),
&\text{if } j_2>0\text{ and }
j_1=0,\\
0,&\text{otherwise}
\end{cases}
\end{align*}
and
\begin{align*}
\gamma_{\vj}(x,y):=\begin{cases}
[j_2^{t}\zeta(t)]^{s^{-1}}\mathbf1_{(j_2,j_2+1]}(y)
\mathbf{1}_{(S_t(j_2-1),S_t(j_2)]}[\zeta(t)x],
&\text{if } j_2>0\text{ and }
j_1=0,\\
0,&\text{otherwise},
\end{cases}
\end{align*}
where $\zeta(t):=\sum_{i=1}^\infty i^{-t}$.
For $j_1=0$, any $j_2>0$, and $Q:=Q_1\times Q_2\in\mathscr D_{\vj}(\mathbb R\times\mathbb R)$,
we obviously have
\begin{align*}
\|\gamma_{\vj}\|_{\aveL^s(Q)}\le\bigg[\fint_{Q_2}\int_{0}^1j_2^{t}\zeta(t)\mathbf1_{(j_2,j_2+1]}(y)
\mathbf{1}_{(S_t(j_2-1),S_t(j_2)]}[\zeta(t)x]\,dx\,dy\bigg]^{\frac1s}=1.
\end{align*}
It is also obvious that, for any $\vj_1\neq\vj_2$,
$\operatorname{supp}\gamma_{\vj_1}\cap\operatorname{supp}\gamma_{\vj_2}=\varnothing$.
Then
$$\sup_{\Omega\in\Open(\R^{\vn})}\bigg\|\sup_{\vj\in\mathbb Z^2}
\one_{\Omega_{\vj}} \gamma_{\vj} \bigg\|_{\aveL^s(\Omega)}\le
\sup_{\vj\in\Z^2}\sup_{Q\in\mathscr D_{\vj}(\mathbb R\times\mathbb R)}
\|\gamma_{\vj}\|_{\aveL^s(Q)}\le1,
$$
which shows that the assumption \eqref{supomega} holds.
Now, we verify that \eqref{counterexample4} holds.
Indeed,
\begin{align*}
\bNorm{\{ f_{\vj}\}_{\vj\in\Z^2}}
{\ell^{q_1}L_x^{p}\ell^{q_2}L_y^p}
&=\Bigg[\int_{\mathbb R}\Bigg\{\sum_{j_2=1}^\infty\bigg[\int_{\mathbb R}j_{2}^{-d}
\mathbf1_{(j_2,j_2+1]}(y)
\mathbf{1}_{(0,1]}(x)\,dy\bigg]^{\frac{q_2}p}\Bigg\}^{\frac p{q_2}}
\,dx\Bigg]^{\frac1p}\\
&=\Bigg(\sum_{j_2=1}^\infty j_2^{-\frac{q_2d}p}\Bigg)^{\frac1{q_2}}
\bigg\{\int_{\mathbb R}\Big[\mathbf1_{(0,1]}(x)\Big]^p\,dx\bigg\}^{\frac1p}<\infty
\end{align*}
and
\begin{align*}
&\bNorm{\{\gamma_{\vj}f_{\vj}\}_{\vj\in\Z^2}}{\ell^{q_1}L_x^{p}\ell^{q_2}L_y^p}\\
&\quad=\bigg[\int_0^1\bigg\{\sum_{j_2=1}^\infty\bigg[\int_{\mathbb R} j_2^{-d}
[j_2^{t}\zeta(t)]^{ps^{-1}}\mathbf1_{(j_2,j_2+1]}(y)
\mathbf{1}_{(S_t(j_2-1),S_t(j_2)]}[\zeta(t)x]\,dy
\bigg]^{\frac{q_2}p}\bigg\}^{\frac{p}{q_2}}\,dx\bigg]^{\frac1p}\\
&\quad\sim\bigg[\sum_{j_2=1}^\infty\int_{S_t(j_2-1)}^{S_t(j_2)}
j_2^{-d+tps^{-1}}\,dx\bigg]^{\frac 1p}=\infty.
\end{align*}
This finishes the proof of Proposition \ref{counterexample-gamma}.
\end{proof}

Indeed, if $(L^{\vp}\ell^{\vq})_\pi$ looks like
$\cdots L^p\cdots\ell^q\cdots L^p\cdots$ with $q>p$,
then the estimate \eqref{finally} also fails.
This further implies that condition $\vq\le\vp$ in \eqref{key cases}
is necessary in some sense.
More precisely, we have the following proposition.

\begin{proposition}\label{counterexamplecom}
Let $0<p<s<\infty$, $\vec p:=p\cdot\vec 1$, $\vec q\in(0,\infty)^k$,
and $\pi\in S_{[2k]}$.
Assume that there exist
$0<\gamma_1<\gamma_2<\gamma_3\le 2k$
such that
$\pi(\gamma_1),\pi(\gamma_3)\in[k]$,
$\pi(\gamma_2)\in[2k]\setminus[k]$, and
$q_{\pi(\gamma_2)-k}>s.$
Then there exist measurable functions
$\{\{\gamma^{(N)}_{\vj}\}_{\vj\in\mathbb Z^k}\}_{N\in\Z_+}$ satisfying \eqref{supomega}
and
$\{\{f^{(N)}_{\vj}\}_{\vj\in\mathbb Z^k}\}_{N\in\Z_+}\subset (L^{\vec p}\ell^{\vec q})_\pi$
with each $f^{(N)}_{\vj}$ constant on every $R\in\mathscr D_{\vj}$ such that
$$
\left\|\{\gamma_{\vj}^{(N)} f^{(N)}_{\vj}\}_{\vj\in\Z^k}
\right\|_{(L^{\vp}\ell^{\vq})_\pi}
\left\|\{f^{(N)}_{\vj}\}_{\vj\in\Z^k}\right\|^{-1}_{(L^{\vp}\ell^{\vq})_\pi}
\to \infty
\quad\text{as}\quad N\to\infty.
$$
\end{proposition}

\begin{proof}
Let $a:=\pi(\gamma_1)$,
$b:=\pi(\gamma_2)-k$,
$c:=\pi(\gamma_3)$,
and $a,b,c\in[k]$ with $q_b>s$.
Let $N\in\mathbb Z_+$. For $\vj=(j_1,\dots,j_k)\in\mathbb Z^k$, define
$$
f^{(N)}_{\vj}(x)
:=
\begin{cases}
\mathbf 1_{(j_b,j_b+1]\times(0,1]^{n_c-1}}(x_c)\,
\mathbf 1_{(0,1]^{n-n_c}}(x\ominus x_c),
&\text{if } 0<j_b\le N \text{ and } \vj\ominus j_b=\vec 0,\\
0,&\text{otherwise}.
\end{cases}
$$
We partition the interval $(0,1)$ into $2^N$ disjoint intervals
$$
I_l := \Big(l2^{-N},(l+1)2^{-N}\Big), \qquad\text{for } l\in\{0,\cdots,2^N-1\}.
$$
Each $I_l$ is further subdivided into $N$ subintervals
$$
I_{l,i}
:=
\Bigg(
l2^{-N}+(i-1)\frac{2^{-N}}{N},
\;
l2^{-N}+i\frac{2^{-N}}{N}
\Bigg),
\qquad\text{for } i\in[N].
$$
Let
$$
A_i := \bigcup_{l=0}^{2^N-1} I_{l,i},
\qquad\text{for } i\in[N]
$$
and then $|A_i|=\frac1N$.
Moreover, let
$$
\gamma^{(N)}_{\vj}(x)
:=
\begin{cases}
N\,
\mathbf 1_{(j_b,j_b+1]\times(0,1]^{n_c-1}}(x_c)\,
\mathbf 1_{A_{j_b}\times(0,1]^{n_a-1}}(x_a)\,
\mathbf 1_{(0,1]^{n-n_a-n_c}}(x\ominus x_a\ominus x_c), \\
\phantom{0,}\qquad \text{if } 0<j_b\le N \text{ and } \vj\ominus j_b=\vec 0,\\
0,\qquad \text{otherwise}.
\end{cases}
$$
It is straightforward to verify that
$\{\gamma^{(N)}_{\vj}\}_{\vj\in\mathbb Z^k}$ satisfies \eqref{supomega}.
By a direct computation, we conclude that
\begin{align*}
\left\|\{f^{(N)}_{\vj}\}_{\vj\in\Z^k}\right\|_{(L^{\vp}\ell^{\vq})_\pi}=N^\frac1{q_{b}}
\quad\text{and}\quad
\left\|\{\gamma_{\vj}^{(N)} f^{(N)}_{\vj}\}_{\vj\in\Z^k}
\right\|_{(L^{\vp}\ell^{\vq})_\pi}=N^{\frac1s}.
\end{align*}
This, combined with $q_b>s$, finishes the proof of Proposition \ref{counterexamplecom}.
\end{proof}

\subsection{Carleson conditions related to matrix weights} \label{Carl special1}

As in the one-parameter theory developed by Frazier and Roudenko \cite[Section 3]{FR21}, the key example in the context of matrix weights to which we wish to apply the general Theorem \ref{carl thm} involves multipliers of the form
\begin{equation}\label{carl main ex}
\gamma_{\vj}(x):=\sum_{R\in\D_{\vj}(\R^{\vn})}\gamma_R(x)\one_R(x),\quad
\gamma_R(x):=\abs{V(x)[V]_{\aveL^p(R)}^{-1}},
\end{equation}
where $V\in\A_p(\R^{\vn})$. These are easily seen to be compatible with cases \eqref{carl B} and \eqref{carl arb} of Theorem \ref{carl thm}:

\begin{lemma}\label{Ap 2 easy carl}
For $V\in\A_p(\R^{\vn})$, there exists $s\in(p,\infty)$ such that
\begin{equation*}
\sup_{R\in\D(\R^{\vn})}\bigg[\fint_R\abs{V(x)[V]_{\aveL^p(R)}^{-1}}^s dx\bigg]^{1/s}\lesssim 1,
\end{equation*}
uniformly over bounded subsets of $\A_p(\R^{\vn})$.
\end{lemma}

\begin{proof}
By Corollary \ref{Ap to RHI}, we can pick an $s\in(p,\infty)$, uniformly over bounded subsets of $\A_p(\R^{\vn})$, such that $V\in\RHI_{p,s}(\Rect(\R^{\vn}))$, and thus
\begin{align*}
\bigg[\fint_R\abs{V(x)[V]_{\aveL^p(R)}^{-1}}^s dx\bigg]^{\frac 1s}
&\lesssim \bigg[\fint_R\abs{V(x)[V]_{\aveL^p(R)}^{-1}}^p dx\bigg]^{\frac 1p} \\
&\sim\abs{[V]_{\aveL^p(R)}[V]_{\aveL^p(R)}^{-1}}=1.
\end{align*}
This finishes the proof of Lemma \ref{Ap 2 easy carl}.
\end{proof}

Let us then discuss case \eqref{carl F} of Theorem \ref{carl thm}. In order to facilitate the discussion, we give this condition a name in the special case of functions
of the form \eqref{carl main ex}:

\begin{definition}\label{Carl prop}
For $0<p<s<\infty$, we say that $V\in\AC_{p,s}(\R^{\vn})$ if $V\in\A_p(\R^{\vn})$ and
\begin{equation}\label{carl to prove}
\sup_{\Omega\in\Open(\R^{\vn})}\fint_{\Omega}
\sup_{R\in\D(\R^{\vn}),\, x\in R\subset\Omega}
\abs{V(x)[V]_{\aveL^p(R)}^{-1}}^s \, dx\lesssim 1.
\end{equation}
We say that $V\in\AC_p(\R^{\vn})$ if $V\in\AC_{p,s}(\R^{\vn})$ for some $s\in(p,\infty)$.
\end{definition}

\begin{lemma}\label{Carl prop 1}
In the one-parameter case, $\AC_p(\R^n)=\A_p(\R^n)$ for all $p\in(0,\infty)$. Moreover, bounded subsets of $\A_p(\R^n)$ are bounded in $\AC_{p,s}(\R^n)$ for some $s\in(p,\infty)$ uniformly over such subsets.
\end{lemma}

\begin{proof}
This is essentially \cite[Lemma 3.3]{Gold03} for $p\in(1,\infty)$ and \cite[Lemma 3.3]{FR21} for $p\in(0,1]$; an inspection of both arguments shows that they give \eqref{carl to prove} with the same $s\in(p,\infty)$ as that in the reverse H\"older class $\RHI_{p,s}(\R^n)$ containing a given $V\in\A_p(\R^n)$; by Lemma \ref{RHI}, this is uniform over bounded subsets of $\A_p(\R^n)$.

More precisely, both \cite[Lemma 3.3]{Gold03} and  \cite[Lemma 3.3]{FR21} deal with \eqref{carl to prove} for $\Omega=P\in\D(\R^n)$; however, the estimate as stated is a quick consequence: Denoting by $\D^*(\Omega)$ the collection of maximal dyadic cubes $P\subset\Omega$, these cubes are disjoint, and we have
\begin{align*}
\int_{\Omega}\sup_{\genfrac{}{}{0pt}{}{R\in\D(\R^n)}{x\in R\subset\Omega}}
\abs{V(x)[V]_{\aveL^p(R)}^{-1}}^s \, dx
&=\sum_{P\in\D^*(\Omega)}\int_P \sup_{\genfrac{}{}{0pt}{}{R\in\D(\R^n)}{x\in R\subset P}}
\abs{V(x)[V]_{\aveL^p(R)}^{-1}}^s \, dx \\
&\lesssim\sum_{P\in\D^*(\Omega)}\abs{P}=\abs{\Omega},
\end{align*}
where the inequality for each $P$ is \cite[Lemma 3.3]{Gold03} or \cite[Lemma 3.3]{FR21}, depending on the value of $p$.
This finishes the proof of Lemma \ref{Carl prop 1}.
\end{proof}

The original motivation of \cite[Lemma 3.3]{Gold03} was to prove the boundedness of the matrix-weighted maximal operator
\begin{equation}\label{MV def}
\mathcal M_V f(x):=\sup_{Q\in\D(\R^n),\, x\in Q} \fint_Q\abs{V(x)V^{-1}(y)f(y)}dy
\end{equation}
on $L^u(\R^n)$ for all $u$ in some neighbourhood of $p\in(1,\infty)$, provided that $V\in\A_p(\R^n)$; this is \cite[Theorem 3.2]{Gold03}.
(The definition of $\mathcal M_Vf(x)$ as stated in \eqref{MV def} is from \cite[p.\,1996]{CG01}, where case $p=2$ of the just mentioned result was proved, but actually an equivalent definition and case $u=p\in(1,\infty)$ of the theorem already appeared in \cite[p.\,121]{NT96}.)

The following lemma shows that the reasoning of \cite{Gold03} works in both directions, and also in the multi-parameter case. Similar statements about the equivalence of maximal and Carleson embedding-type inequalities are known in different situations, see e.g. \cite{HK11}.

\begin{lemma}\label{carl equiv M}
Let $0<v<p\leq 1$ or $v=1<p<\infty$, and let $p<s<\infty$. For $V\in\A_p(\R^{\vn})$, consider the maximal operator
\begin{equation*}
\mathcal M_{\vn,V}^{(v)}f(x):=\bigg[\sup_{R\in\D(\R^{\vn}),\, x\in R}
\fint_R\abs{V(x)V^{-1}(y)f(y)}^v dy\bigg]^{\frac 1v}.
\end{equation*}
Then the following conditions are equivalent:
\begin{enumerate}[\rm(i)]
\item\label{ACps} $V\in\AC_{p,s}(\R^{\vn})$;
\item\label{M Ls} the operator $\mathcal M_{\vn,V}^{(v)}$ is bounded on $L^s(\R^{\vn})$;
\item\label{M Lu} for some $r\in(0,p)$, whose choice is uniform over bounded subsets of $\A_p(\R^{\vn})$, the operator $\mathcal M_{\vn,V}^{(v)}$ is bounded on $L^u(\R^{\vn})$ for all $u\in(r,s]$.
\end{enumerate}
\end{lemma}

\begin{proof}
\eqref{M Lu} $\Rightarrow$ \eqref{M Ls} is obvious.

\eqref{M Ls} $\Rightarrow$ \eqref{ACps}:
Let $\Omega\in\Open(\R^{\vn})$ and $R\in\mathscr D(\R^{\vn})$ be such that $R\subset\Omega$.
By Lemma \ref{inverses}, we conclude that
\begin{equation}\label{carl inv}
\begin{split}
\abs{V(x)[V]_{\aveL^p(R)}^{-1}}
&=\abs{[V]_{\aveL^p(R)}^{-1}V(x)}
\lesssim\abs{[V^{-1}]_{\aveL^v(R)}V(x)}  \\
&\sim\bigg[\fint_R\abs{V^{-1}(y)V(x)}^v \, dy\bigg]^{\frac 1v}
=\bigg[\fint_R\abs{V(x)V^{-1}(y)}^v \, dy\bigg]^{\frac 1v} \\
&\sim\sum_{i=1}^m\bigg[\fint_R\abs{V(x)V^{-1}(y)e_i}^v \, dy\bigg]^{\frac 1v},
\end{split}
\end{equation}
where $(e_i)_{i=1}^m$ is a fixed orthonormal basis of $\F^m$.
Since $R\subset\Omega$, we can insert $\one_\Omega(y)$ inside the integrand on the right-hand side of \eqref{carl inv}.
Hence
\begin{equation*}
\sup_{R\in\D(\R^{\vn}),\, x\in R\subset\Omega}\abs{V(x)[V]_{\aveL^p(R)}^{-1}}
\lesssim\sum_{i=1}^m \mathcal M_{\vn,V}^{(v)}(\one_{\Omega}e_i)(x).
\end{equation*}
An application of assumption \eqref{M Ls} to $f=\one_\Omega e_i\in L^s(\R^{\vn};\F^m)$ for each $i\in[m]$ shows that
\begin{equation*}
\bigg\{\int_{\Omega} \sup_{R\in\D(\R^{\vn}),\, x\in R\subset\Omega}
\abs{V(x)[V]_{\aveL^p(R)}^{-1}}^s dx\bigg\}^{\frac 1s}
\lesssim \sum_{i=1}^m\bigg[\int_{\R^{\vn}}\abs{\one_\Omega(y)e_i}^s dy\bigg]^{\frac 1s}
\sim\abs{\Omega}^{\frac 1s}.
\end{equation*}
Valid for every $\Omega\in\Open(\R^{\vn})$, this is the claimed property $V\in\AC_{p,s}(\R^{\vn})$.

\eqref{ACps} $\Rightarrow$ \eqref{M Lu}:
For each $R\in\D(\R^{\vn})$, we have

\begin{equation}\label{carl to M1}
\begin{split}
&\bigg[\fint_R\abs{V(x)V^{-1}(y)f(y)}^v \, dy\bigg]^{\frac 1v} \\
&\quad\leq\abs{V(x)[V]_{\aveL^p(R)}^{-1}}
\bigg\{\fint_R\abs{[V]_{\aveL^p(R)}V^{-1}(y)f(y)}^v \, dy\bigg\}^{\frac 1v}
=: I_R(x) II_R.
\end{split}
\end{equation}

Here, from H\"older's inequality, it follows that
\begin{equation*}
II_R\leq\Norm{[V]_{\aveL^p(R)}V^{-1}}{\aveL^q(R)}\Norm{f}{\aveL^r(R)},
\qquad\frac{1}{v}=\frac{1}{q}+\frac{1}{r}.
\end{equation*}
If $p\in(0,1]$, we take $q=p'=\infty$, hence $r=v\in(0,p)$, and the definition of $\A_p(\R^{\vn})$ implies that
\begin{equation*}
\Norm{[V]_{\aveL^p(R)}V^{-1}}{\aveL^\infty(R)}\lesssim 1.
\end{equation*}
If $p\in(1,\infty)$, then $V^{-1}\in\A_{p'}(\R^{\vn})$, and we choose $q\in(p',\infty)$ so that $V^{-1}\in\RHI_{p',q}(\R^{\vn})$, which is possible by Corollary \ref{Ap to RHI}. Then
\begin{equation*}
\Norm{[V]_{\aveL^p(R)}V^{-1}}{\aveL^q(R)}
=\Norm{V^{-1}[V]_{\aveL^p(R)}}{\aveL^q(R)}
\sim\Norm{V^{-1}[V]_{\aveL^p(R)}}{\aveL^{p'}(R)}\lesssim 1.
\end{equation*}
Thus, when $p\in(1,\infty)$, we have $v=1$ by assumption and $q\in(p',\infty)$ by the choice just made; hence
\begin{equation*}
\frac{1}{r}=1-\frac{1}{q}>1-\frac{1}{p'}=\frac{1}{p}
\end{equation*}
and thus $r\in(0,p)$. Therefore, in both cases $p\in(0,1]$ and $p\in(1,\infty)$, we obtain
\begin{equation*}
II_R\lesssim \Norm{f}{\aveL^r(R)}\leq\inf_{z\in R}
\Big[\mathcal M_{\vn}(\abs{f}^r)(z)\Big]^{\frac 1r}
\end{equation*}
for some $r\in(0,p)$. If $f\in L^u(\R^{\vn})$ for some $u>r$, the boundedness of $\mathcal M_{\vn}$ on $L^{u/r}(\R^{\vn})$ shows that $\Norm{f}{\aveL^r(R)}<\infty$ for every $R\in\D(\R^{\vn})$.

For $\kappa\in\Z$, let
\begin{equation*}
\Rs_\kappa:=\{R\in\D(\R^{\vn}):\
\Norm{f}{\aveL^r(R)}>2^\kappa\}.
\end{equation*}
If $R$ is not in any $\Rs_\kappa$, then $II_R\lesssim\Norm{f}{\aveL^r(R)}=0$, and \eqref{carl to M1} shows that such rectangles $R$ do not contribute to the supremum defining the maximal function $\mathcal M_{\vn,V}^{(v)}f$. If $R$ is in every $\Rs_\kappa$, then $\Norm{f}{\aveL^r(R)}=\infty$, and this we have already excluded. Thus, all relevant $R$ belong to some but not all $\Rs_\kappa$, and hence to some $\Rs_\kappa\setminus\Rs_{\kappa+1}$. Therefore
\begin{align*}
\mathcal M_{\vn,V}^{(v)}f(x)
&=\sup_{R\in\D(\R^{\vn})} \bigg[\fint_R\abs{V(x)V^{-1}(y)f(y)}^v \, dy
\bigg]^{\frac 1v} \one_R(x) \\
&\lesssim\sup_{R\in\D(\R^{\vn})} I_R(x) \Norm{f}{\aveL^r(R)} \one_R(x) \\
&\lesssim\sup_{\kappa\in\Z}\sup_{R\in\Rs_\kappa\setminus\Rs_{\kappa+1}} I_R(x) 2^\kappa \one_R(x)
\leq\sup_{\kappa\in\Z}2^{\kappa}
\sup_{R\subset\Omega^\kappa,\, x\in R} I_R(x),
\end{align*}
where
\begin{align*}
\Omega^\kappa:=\Big\{x\in\R^{\vn}:\
[\mathcal M_{\vn}(\abs{f}^r)(x)]^{\frac 1r}>2^\kappa\Big\}.
\end{align*}

Recalling the definition of $I_R(x)$, for $V\in\AC_{p,s}(\R^{\vn})$ and $u\in(r,s]$, it then follows that
\begin{align*}
\bNorm{\mathcal M_{\vn,V}^{(v)}f}{L^u}^u
&\lesssim\sum_{\kappa\in\Z} 2^{\kappa u} \int_{\R^{\vn}}\sup_{R\subset\Omega^\kappa,\, x\in R} \abs{V(x)[V]_{\aveL^p(R)}^{-1}}^u dx \\
&=\sum_{\kappa\in\Z} 2^{\kappa u}\abs{\Omega^\kappa} \fint_{\Omega_\kappa}\sup_{R\subset\Omega^\kappa,\, x\in R} \abs{V(x)[V]_{\aveL^p(R)}^{-1}}^u dx \\
&\leq\sum_{\kappa\in\Z} 2^{\kappa u}\abs{\Omega^\kappa}
\bigg[ \fint_{\Omega_\kappa}\sup_{R\subset\Omega^\kappa,\, x\in R} \abs{V(x)[V]_{\aveL^p(R)}^{-1}}^s \, dx \bigg]^{\frac us} \\
&\lesssim \sum_{\kappa\in\Z} 2^{\kappa u}\abs{\Omega^\kappa}
\sim \BNorm{(\mathcal M_{\vn}\abs{f}^r)^{\frac 1r}}{L^u(\R^{\vn})}^u
\lesssim \Norm{f}{L^u(\R^{\vn})}^u.
\end{align*}
This finishes the proof of Lemma \ref{carl equiv M}.
\end{proof}

The following corollary is a restatement of \cite[Theorem 3.2]{Gold03} when $p\in(1,\infty)$, but may be new for $p\in(0,1]$.

\begin{corollary}\label{MVbd}
Let $0<v<p\leq 1$ or $v=1<p<\infty$. If $V\in\A_p(\R^{n})$ is a one-parameter weight, then the one-parameter maximal operator $\mathcal M_{n,V}^{(v)}$ is bounded on $L^u(\R^{n})$ for all $u$ in some neighbourhood of $p$, uniformly over bounded subsets of $\A_p(\R^n)$.
\end{corollary}

\begin{proof}
By Lemma \ref{Carl prop 1}, we conclude that
$\A_p(\R^{n})=\AC_p(\R^n)$ in the one-parameter situation.
From Lemma \ref{carl equiv M}, it follows that
the latter property is equivalent to the claimed maximal inequality.

Moreover, by Lemma \ref{Carl prop 1}, bounded subsets of $V\in\A_p(\R^n)$ are bounded in $\AC_{p,s}(\R^n)$ for some $s\in(p,\infty)$ uniformly over such subsets. By Lemma \ref{carl equiv M}, this gives the boundedness of $\mathcal M_{n,V}^{(v)}$ for all such $V$ on $L^u(\R^n)$, for all $u\in(r,s]$, where $r\in(0,p)$. An inspection of the proof shows that one can take $r=v$ for $p\in(0,1]$ and $r=q'$ for $p\in(1,\infty)$, provided that $V^{-1}\in\A_{p'}(\R^n)$ belongs to $\RHI_{p',q}$. Since $V\mapsto V^{-1}$ maps bounded subsets of $\A_p(\R^n)$ into bounded subsets of $\A_{p'}(\R^n)$, and such subsets are bounded in $\RHI_{p',q}$ for some common $q\in(p',\infty)$ by Lemma \ref{RHI}, it follows that one can always choose a common $r\in(0,p)$ for all $V$ in bounded subsets of $\A_p(\R^n)$, as claimed.
This finishes the proof of Corollary \ref{MVbd}.
\end{proof}

To give an informal summary of our observations about the $\AC_p$ condition so far, we have seen that:
\begin{itemize}
\item in general, $V\in\AC_p(\R^{\vn})$ if and only if a $V$-weighted maximal operator satisfies certain inequalities (Lemma \ref{carl equiv M});
\item in the one-parameter case, both $V\in\AC_p(\R^{n})$ and the said maximal inequalities hold unconditionally (Lemma \ref{Carl prop 1} and Corollary \ref{MVbd}).
\end{itemize}
What remains, then, to complete our picture of the $\AC_p$ condition, is to understand the boundedness of the matrix-weighted maximal operator in the general multi-parameter setting. In a recent work of \cite{Vuo24}, this has been accomplished for $p\in(1,\infty)$, but his method does not seem to extend to $p\in(0,1]$. For this reason, we next develop a new framework for vector-valued maximal operators, which may have independent interest.

\subsection{Reducing operator-valued maximal functions}\label{red M}

In this subsection, we are going to establish the boundedness of
the matrix-weighted strong maximal operator.
To this end, we introduce the reducing operator-valued maximal functions.

\subsubsection{The one-parameter case}

Given $F\in L^v_{\loc}(\R^n;\F^{m\times\ell})$, we consider the following function of variables $x\in\R^n$ and $e\in\F^m$:
\begin{equation*}
r(x,e)
:=(\mathcal M^{(v)}\abs{F(\cdot)^*e})(x)
:=\sup_{Q\in\Cubes(\R^{n}), x\in Q} \bigg[\fint_Q\abs{F(y)^*e}^v \, dy\bigg]^{\frac 1v}.
\end{equation*}
Note that $F(y)^*\in\F^{\ell\times m}$, so that $F(y)^*e\in\F^\ell$ is well defined.

For each fixed $x\in\R^n$, the mapping $e\mapsto r(x,e)$ defines a $\min\{1,v\}$-seminorm on $\F^m$.
Noting that the same supremum is achieved over a fixed countable family of cubes,
say all those with rational corners, it follows that
$x\mapsto r(x,e)$ is measurable for each $e\in\F^m$.
Hence Theorem \ref{selectA} guarantees a measurable selection of
reducing operators inducing equivalent Euclidean seminorms;
we denote one such choice by $\M^{(v)}F(x)$. Hence
\begin{equation*}
\abs{[\M^{(v)}F(x)]e}\sim \mathcal M^{(v)}(|F(\cdot)^*e|).
\end{equation*}
Using the definition of the induced matrix norm $\abs{A}:=\sup_{\abs{e}\leq 1}\abs{Ae}$, this easily extends to any matrix $A\in\F^{m\times m}$ in place of the vector $e\in\F^m$:
\begin{equation}\label{MA}
\abs{[\M^{(v)}F(x)] A}\sim \mathcal M^{(v)}(|F(\cdot)^*A|).
\end{equation}

Our main application deals with vector-valued functions $f\in L^v_{\loc}(\R^n;\F^m)$. With the usual identification of $\F^m=\F^{m\times 1}$ as the space of column vectors, the expressions $f(x)^*e\in\F$, hence $\mathcal M^{(v)}(f(\cdot)^*e)$, and thus $\M^{(v)}f$ are well defined.

Nevertheless, even to make estimates with such vector-valued functions only, we will also need to deal with matrix-valued functions as intermediate steps. Extending the definition of the maximal operator $\mathcal M_V^{(v)}$ from \eqref{MV def} to such functions is immediate:
$$
\mathcal M_V^{(v)}F(x)
:=\sup_{Q\in\Cubes(\R^{n}), x\in Q}\bigg[\fint_Q\abs{V(x)V(y)^{-1}F(y)}^{v}\,dy\bigg]^{\frac 1v}.
$$

\begin{lemma}\label{MVvsMnew}
For a matrix weight $V$ and $F\in L^v_{\loc}(V;\F^{m\times\ell})$,
\begin{equation*}
\mathcal M_V^{(v)}F(x)
\sim\abs{V(x)\M^{(v)}(V(\cdot)^{-1}F(\cdot))(x)}.
\end{equation*}
\end{lemma}

\begin{proof}
Recall that $\abs{AB}=\abs{B^*A^*}$, and that both matrix weights $V(x)$ and reducing operators are self-adjoint by definition.
Applying these and \eqref{MA} with $A=V(x)$, we obtain
\begin{align*}
\mathcal M_V^{(v)}F(x)
&=\sup_{Q\ni x}\bigg[\fint_Q
\Big| [V(y)^{-1}F(y)]^* V(x) \Big|^v \, dy\bigg]^{\frac 1v} \\
&=\mathcal M^{(v)}\Big( |[V(\cdot)^{-1}F(\cdot)]^*V(x)|\Big)(x) \\
&\sim\abs{\M^{(v)}(V(\cdot)^{-1}F(\cdot))(x)V(x)}
=\abs{V(x)\M^{(v)}(V(\cdot)^{-1}F(\cdot))(x)}.
\end{align*}
This finishes the proof of Lemma \ref{MVvsMnew}.
\end{proof}

\begin{proposition}\label{newMbd}
Let $0<v<p\leq 1$ or $v=1<p<\infty$. If $V\in\A_p(\R^n)$, then for all $u$ in some neighbourhood of $p$, and all $F\in L^u(V;\F^{m\times\ell})$,
\begin{equation*}
\Norm{\M^{(v)}F}{L^u(V;\F^{m\times m})}\lesssim\Norm{F}{L^u(V;\F^{m\times\ell})}
\end{equation*}
with the implicit positive constant independent of $F$.
\end{proposition}

\begin{proof}
By Corollary \ref{MVbd}, $\mathcal M_V^{(v)}$ is bounded on $L^u(\R^n;\F^m)$ under the assumptions of the present proposition.
Each $F\in L^u(\R^n;\F^{m\times\ell})$ can be written as $F=[f_1,\ldots,f_\ell]$ with each $f_i\in L^u(\R^n;\F^m)$, and there is the identity
\begin{equation*}
V(x)V(y)^{-1}F(y)
=\Big[V(x)V(y)^{-1}f_1(y),\ldots,V(x)V(y)^{-1}f_\ell(y)\Big].
\end{equation*}
Thus, $\mathcal M_V^{(v)}F(x)\sim\sum_{i=1}^\ell \mathcal M_V^{(v)}f_i(x)$ and from this, the boundedness of $\mathcal M_V^{(v)}$ on $L^u(\R^n;\F^{m\times\ell})$ readily follows.

By Lemma \ref{MVvsMnew}, the claimed inequality is equivalent to the said boundedness, and hence also follows.
This finishes the proof of Proposition \ref{newMbd}.
\end{proof}

\subsubsection{The multi-parameter case}

For functions on $\R^{\vn}$, definitions completely analogous to the one-parameter case above can be made. The following two versions are particularly relevant. On the one hand, there is the strong maximal operator
\begin{equation*}
(\mathcal M^{(v)}_{\vn}\abs{F(\cdot)^*e})(x)
:=\sup_{R\in\Rect(\R^{\vn})}\one_R(x)
\bigg[\fint_R\abs{F(y)^*e}^v \, dy\bigg]^{\frac 1v}.
\end{equation*}
On the other hand, for each $i\in[k]$, we can consider the one-parameter maximal operator in the $i$-th variable,
\begin{equation*}
(\mathcal M^{(v)}_i\abs{F(\cdot)^*e})(x)
:=\sup_{Q\in\Cubes(\R^{n_i})} \one_Q(x_i)
\bigg[ \fint_Q\abs{F(x'\oplus y_i)^*e}^v \, dy_i \bigg]^{\frac 1v},
\end{equation*}
where $x=x'\oplus x_i\in\R^{\vn-n_i}\oplus\R^{n_i}$.

In the scalar-valued case, there is a well-known trick, going back to \cite{JMZ},
to dominate $\mathcal M^{(v)}_{\vn}$ by the composition of $\mathcal M^{(v)}_i$ for all $i\in[k]$. This also leads to a quick proof of the weighted $L^p$ boundedness of
$\mathcal M^{(v)}_{\vn}$ for product Muckenhoupt weights, as noted in \cite[(3.4)]{FS82}. In the vector-valued, matrix-weighted situation,
a simple adaptation of this idea runs into trouble, since the maximal operator of a vector-valued function is a scalar-valued one, preventing an efficient iteration. For $v=1$, this obstacle was overcome in \cite{Vuo24} with the help of the set-valued
maximal operator  introduced in \cite{BCU}. However, the framework of \cite{BCU} seems to rely somewhat heavily on convexity, and it seems unclear how to define a meaningful set-valued counterpart of $\mathcal M^{(v)}$ for $v\in(0,1)$. Our substitute is the reducing operator-valued maximal operator, which remains compatible with this quasi-Banach regime.

With considerations as in the one-parameter case, Theorem \ref{selectA} guarantees a measurable selection of reducing operators inducing Euclidean seminorms equivalent to those defined by setting
$$
e\mapsto (\mathcal M^{(v)}_{\vn}\abs{F(\cdot)^*e})(x)
\ \ \text{and}\ \
e\mapsto (\mathcal M^{(v)}_i\abs{F(\cdot)^*e})(x);
$$
we denote some such choices by $\M^{(v)}_{\vn}F(x)$ and $\M^{(v)}_i F(x)$. Hence
\begin{equation*}
\abs{(\M^{(v)}_{\vn}F(x))e}\sim \mathcal M^{(v)}_{\vn}|F(\cdot)^*e|,\qquad
\abs{(\M^{(v)}_iF(x))e}\sim \mathcal M^{(v)}_i|F(\cdot)^*e|.
\end{equation*}
As before, these easily extend to any matrix $A$ in place of the vector $e$:
\begin{equation}\label{MnA}
\abs{(\M^{(v)}_{\vn}F(x))A}\sim \mathcal M^{(v)}_{\vn}|F(\cdot)^*A|,\qquad
\abs{(\M^{(v)}_iF(x))A}\sim \mathcal M^{(v)}_i|F(\cdot)^*A|.
\end{equation}

Let us now extend the classical domination of $\mathcal M_{\vn}^{(v)}$ by $\mathcal M_1^{(v)}\circ\cdots\circ \mathcal M_k^{(v)}$ to the reducing operator-valued maximal operators.

\begin{lemma}\label{MA bd}
Let $v\in(0,1]$. For all functions $F\in L^v_{\loc}(\R^{\vn};\F^{m\times\ell})$, constant matrices $A\in\F^{m\times\ell'}$, and points $x\in\R^n$,
\begin{equation*}
\abs{(\M_{\vn}^{(v)}F)(x)A}
\lesssim\abs{((\M_1^{(v)}\circ\cdots\circ \M_k^{(v)})F)(x)A},
\end{equation*}
with the implicit positive constant independent of $F$, $A,$ and $x$.
\end{lemma}

\begin{proof}
By the first equivalence in \eqref{MnA} and the classical domination of the scalar-valued function $\abs{F(\cdot)^*A}$, we have
\begin{equation}\label{iterBd1}
\abs{(\M_{\vn}^{(v)}F)(\cdot)A}
\sim \mathcal M_{\vn}^{(v)}(F(\cdot)^*A)
\leq \mathcal M_1^{(v)}\circ\cdots\circ \mathcal M_k^{(v)}(F(\cdot)^*A).
\end{equation}
Let $F_{k+1}:=F$. Then $F_k:=\M_k^{(v)}F$ is another matrix-valued function, and thus $F_{k-1}:=\M_{k-1}^{(v)}F_k$ is well defined. In this way, we may recursively define $F_i:=\M_i^{(v)} F_{i+1}$ for all $i\in[k]$.  Applying the second equivalence in \eqref{MnA} with $F_{i+1}$ in place of $F$, we find the pointwise equivalence
\begin{equation}\label{MnInd}
\mathcal M_i^{(v)}(F_{i+1}(\cdot)^*A)\sim\abs{(\M_i^{(v)}F_{i+1}(\cdot))A}=\abs{F_i(\cdot)A}
=\abs{F_i(\cdot)^*A},
\end{equation}
noting in the last step that $F_i(\cdot):=\M_i^{(v)}F_{i+1}(\cdot)$ is self-adjoint-valued by the definition of reducing operators. Applying $\mathcal M_1^{(v)}\circ\cdots\circ \mathcal M_{i-1}^{(v)}$ (interpreted as $\abs{\ }$ if $i=1$) to both sides of \eqref{MnInd}, it follows that
\begin{equation*}
(\mathcal M_1^{(v)}\circ\cdots\circ \mathcal M_i^{(v)})(F_{i+1}(\cdot)^*A)\sim
(\mathcal M_1^{(v)}\circ\cdots\circ \mathcal M_{i-1}^{(v)})(F_i(\cdot)^*A).
\end{equation*}
Concatenating these equivalences for all $i\in[k]$ and recalling that $F_{k+1}:=F$, we obtain
\begin{align*}
(\mathcal M_1^{(v)}\circ\cdots\circ \mathcal M_k^{(v)})(F(\cdot)^*A)
&\sim\abs{F_1(\cdot)^*A}
=\abs{F_1(\cdot)A} \\
&=\abs{((\M_1^{(v)}\circ\cdots\circ \M_k^{(v)})F(\cdot))A}.
\end{align*}
In combination with \eqref{iterBd1},
this finishes the proof of Lemma \ref{MA bd}.
\end{proof}

We are now ready to prove the boundedness of the matrix-weighted strong maximal operator $\mathcal M^{(v)}_{\vn,V}$. For $p\in(1,\infty)$ and $k=2$ (the bi-parameter case), the following result was previously obtained in \cite[Theorem 3.12]{Vuo24}, and indeed this was a key inspiration for the present generalisation. For $p\in(1,\infty)$ and arbitrary $k$, the result could also be obtained by a routine extension of the
considerations in \cite{Vuo24}, but the extension to $p\in(0,1]$ seems to require the new framework that we have introduced in this section. The proof of \cite[Theorem 3.12]{Vuo24} uses the set-valued maximal function of \cite{BCU} and its natural product-space extension; our argument is modelled after this,
but using our reducing operator-valued maximal function instead.

\begin{theorem}\label{newMprod bd}
Let $0<v<p\leq 1$ or $v=1<p<\infty$, and $V\in\A_p(\R^{\vn})$.
Then $\mathcal M^{(v)}_{\vn,V}$ is bounded on $L^u(\R^{\vn})$ for all $u$ in a neighbourhood of $p$, uniformly over bounded subsets of $\A_p(\R^{\vn})$.
\end{theorem}

\begin{proof}
The proof of Lemma \ref{MVvsMnew} extends to the multi-parameter case {\em verbatim} to give the pointwise identity
\begin{equation*}
\mathcal M^{(v)}_{\vn,V}f(x)
=\abs{V(x)\M^{(v)}_{\vn,V}(V(\cdot)^{-1}f(\cdot))(x)}.
\end{equation*}
Hence
\begin{equation*}
\Norm{\mathcal M^{(v)}_{\vn,V}f}{L^u(\R^{\vn})}
=\Norm{\M^{(v)}_{\vn,V}(V(\cdot)^{-1}f(\cdot))}{L^u(V;\F^{m\times m})}.
\end{equation*}
Let $G:=V(\cdot)^{-1}f(\cdot)\in L^u(V;\F^m)=L^u(V;\F^{m\times 1})$. It suffices to prove that $\M^{(v)}_{\vn}$ is bounded from $L^u(V;\F^{m\times 1})$ to $L^u(V;\F^{m\times m})$.

Using Lemma \ref{MA bd} with $A=V(x)$, and the identity $\abs{AB}=\abs{B^*A^*}=\abs{BA}$ for self-adjoint matrices we see that
\begin{align*}
\abs{V(x)(\M^{(v)}_{\vn}G)(x)}
&=\abs{(\M^{(v)}_{\vn}G)(x)V(x)} \\
&\lesssim\abs{(\M^{(v)}_{1}\circ\cdots\circ\M^{(v)}_k G)(x)V(x)} \\
&=\abs{V(x)(\M^{(v)}_{1}\circ\cdots\circ\M^{(v)}_k G)(x)}.
\end{align*}

From Lemma \ref{prod vs coord}, we deduce that each $x_i\mapsto V(x'\oplus x_i)$ belongs to $\A_p(\R^{n_i})$, uniformly in $x'\in\R^{\vn}\ominus\R^{n_i}$. Let
\begin{equation*}
G_i:=(\M^{(v)}_{i}\circ\cdots\circ \M^{(v)}_k)G,\quad i\in[k],\qquad G_{k+1}:=G.
\end{equation*}
Then
\begin{equation*}
\Norm{\M^{(v)}_{\vn}G}{L^u(V;\F^{m\times m})}
\lesssim\Norm{(\M^{(v)}_{1}\circ\cdots\circ \M^{(v)}_k)G}{L^u(V;\F^{m\times m})}
=\Norm{G_1}{L^p(V;\F^{m\times m})}
\end{equation*}
and, for each $i\in[k]$,
\begin{align*}
\Norm{G_i}{L^p(V;\F^{m\times m})}
&=\Norm{\M^{(v)}_{i}G_{i+1}}{L^p(V;\F^{m\times m})}  \\
&=\bigg[\int_{\R^{\vn}\ominus\R^{n_i} }
\Norm{x_i\mapsto \M^{(v)}\big(G_{i+1}(x'\oplus\cdot)\big)(x_i)}{L^u(dx_i,x_i\mapsto V(x'\oplus x_i);\F^{m\times m})}^u \, dx'\bigg]^{\frac1u} \\
&\lesssim \bigg[\int_{\R^{\vn}\ominus\R^{n_i} }
\Norm{x_i\mapsto G_{i+1}(x'\oplus x_i)}{L^u(dx_i,x_i\mapsto V(x'\oplus x_i);\F^{m\times\ell})}^u \, dx'\bigg]^{\frac1u} \\
&=\Norm{G_{i+1}}{L^u(V;\F^{m\times\ell})},
\end{align*}
where $\ell=m$ for $i<k$ and $\ell=1$ for $i=k$, and we used Fubini's theorem for the identities, and Proposition \ref{newMbd}, with the $\A_p(\R^{n_i})$ weights $x_i\mapsto V(x'\oplus x_i)$ in place of $V$, for the inequality. (Note that, even if $G_{k+1}=G=V(\cdot)^{-1}f(\cdot)\in L^u(V;\F^{m\times 1})$ is vector-valued, all other $G_{i+1}$ in the argument above will be $m\times m$-matrix-valued, which is why we need Proposition \ref{newMbd} in the stated generality.)

Concatenating these estimates, we find that
\begin{equation*}
\Norm{\M^{(v)}_{\vn}G}{L^u(V;\F^{m\times m})}\lesssim\Norm{G_1}{L^u(V;\F^{m\times m})}
\lesssim\Norm{G_{k+1}}{L^u(V;\F^{m})}=\Norm{G}{L^u(V;\F^m)},
\end{equation*}
which is the required boundedness of $\M^{(v)}_{\vn}$, and completes the proof.
\end{proof}

\begin{remark}
In the framework of set-valued maximal operators in \cite{BCU,Vuo24}, a vector-valued function $f(x)$ is identified with the set-valued function $F(x):=\{\alpha f(x):\alpha\in\bar B_{\F}\}$. In our framework, where we can allow matrix-valued functions of any number of columns as input for our $\M^{(v)}$, we may simply identify $\F^m$ with $\F^{m\times 1}$ in the most obvious way. If one prefers matrices of same dimensions $m\times m$ as both input and output of $\M^{(v)}$, then one could equally well identify $f(x)$ with the rank-one matrix $F(x)=[f(x),0,\ldots,0]$, or even the positive semidefinite rank-one matrix
\begin{equation}\label{Fff}
F(x)=\frac{f(x)f(x)^*}{\abs{f(x)}}
\end{equation}
instead. The only thing that really matters is the property $\abs{F(x)^*e}=\abs{f(x)^*e}$ for all $e\in\F^m$. Note that $F(x)$ in \eqref{Fff} is an example of a measurable choice of the reducing operator for the seminorms $e\mapsto\abs{f(x)^*e}$. However, through the application of Lemma \ref{MVvsMnew}, we also need to apply $\M^{(v)}$ to functions like $V(\cdot)^{-1}F(\cdot)$, which need not be self-adjoint-valued even if both $V$ and $F$ are; hence it does not seem purposeful to restrict the treatment to self-adjoint-valued functions only.
\end{remark}

\subsection{Carleson estimates related to matrix weights} \label{Carl special3}

With the results concerning maximal operators from the previous subsection, we can now complete our discussion of Carleson embedding-type inequalities related to matrix weights.

\begin{corollary}\label{Ap always ACp}
For all $p\in(0,\infty)$, we have $\A_p(\R^{\vn})=\AC_p(\R^{\vn})$.
\end{corollary}

\begin{proof}
This is immediate from the combination of the $V$-weighted maximal inequalities for $V\in\A_p(\R^{\vn})$ (Theorem \ref{newMprod bd}) and the equivalence of these maximal inequalities with $V\in\AC_p(\R^{\vn})$ (Lemma \ref{carl equiv M}).
\end{proof}

\begin{corollary}\label{Carl Ap}
Let $\vp=p\cdot\vone$ with $p\in(0,\infty)$, let $\vq\in(0,\infty]^k$, and $\pi\in S_{[2k]}$.
Let $V\in\A_p(\R^{\vn})$ and $\gamma_{\vj}(x):=V(x)[V]_{\aveL^p(R)}^{-1}$ for all $x\in R\in\D_{\vj}$ and $\vj\in\Z^k$.
Let $f_{\vj}$ be constant on each $R\in\D_{\vj}$, for each $\vj\in\Z^k$.
Then the inequality
\begin{equation*}
\bNorm{\{\gamma_{\vj}f_{\vj}\}_{\vj\in\Z^k}}{\Xs}
\lesssim\bNorm{\{ f_{\vj}\}_{\vj\in\Z^k}}{\Xs}
\end{equation*}
holds, with the implicit positive constants independent of $f_{\vj}$, in each of the following cases:
\begin{equation*}
\Xs =\ell^{\vq}L^p,\quad\text{or}\quad
\Xs = L^p\ell^{\vq},\quad\text{or}\quad
\Xs = (L^{\vp}\ell^{\vq})_\pi\quad\text{with}\quad \vq\leq\vp.
\end{equation*}
\end{corollary}

\begin{proof}
By Corollary \ref{Ap always ACp}, we have $V\in\AC_p(\R^{\vn})$, thus $V\in\AC_{p,s}(\R^{\vn})$ for some $s\in(p,\infty)$ by Definition \ref{Carl prop}. Written out with the notation $\gamma_{\vj}$ and making elementary estimates, this means that
\begin{equation*}
1
\gtrsim\sup_{\Omega\in\Open(\R^{\vn})}
\Big\|\sup_{\vj\in\mathbb Z^k} \one_{\Omega_{\vj}} \gamma_{\vj}\Big\|_{\aveL^s(\Omega)}
\geq\sup_{\vj\in\mathbb Z^k}\sup_{R\in\D_{\vj}(\R^{\vn})}
\Norm{\gamma_{\vj}}{\aveL^s(R)}
\geq\sup_{\vj\in\mathbb Z^k}\sup_{R\in\D_{\vj}(\R^{\vn})}
\Norm{\gamma_{\vj}}{\aveL^p(R)}.
\end{equation*}
Thus, the assumptions of Theorem \ref{carl thm} are satisfied in each case under consideration, and the said theorem implies the claimed inequality.
\end{proof}

We will also need a certain dual version of the results just discussed, dealing with $[V]_{\aveL^p(R)}V(x)^{-1}$ in place of $V(x)[V]_{\aveL^p(R)}^{-1}$:

\begin{proposition}\label{dualCarl}
Let $p\in(0,\infty)$ and $V\in\A_p(\R^{\vec n})$.
\begin{enumerate}[\rm(i)]
\item\label{dualCarl1}
For some $s>p'$ if $p\in(1,\infty)$, or all $s\in(0,\infty)$ if $p\in(0,1]$, we have
\begin{equation}\label{dualCarl10}
  \sup_{\Omega\in\Open(\R^{\vn})}
  \bigg[\fint_\Omega \sup_{\genfrac{}{}{0pt}{}{R\in\D(\R^{\vn})}{x\in R\subset\Omega}}
  \abs{[V]_{\aveL^p(R)}V(x)^{-1}}^s \, dx\bigg]^{\frac1s}
  \lesssim 1.
\end{equation}

\item\label{dualCarl2} Let moreover $\vr=r\cdot\vone$ with
\begin{equation}\label{r vs p}
r\in
\begin{cases}
(0,p'], & \text{if } p\in(1,\infty), \\
(0,\infty), & \text{if } p\in(0,1],
\end{cases}
\end{equation}
let $\vq\in(0,\infty]^k$ and $\pi\in S_{[2k]}$.
Let $\gamma_{\vj}(x):=[V]_{\aveL^p(R)}V(x)^{-1}$ for all $x\in R\in\D_{\vj}(\R^{\vn})$ and $\vj\in\Z^k$.
Let $f_{\vj}$ be constant on each $R\in\D_{\vj}(\R^{\vn})$, for each $\vj\in\Z^k$.
Then the inequality
\begin{equation*}
\bNorm{\{\gamma_{\vj}f_{\vj}\}_{\vj\in\Z^k}}{\Ys}
\lesssim\bNorm{\{ f_{\vj}\}_{\vj\in\Z^k}}{\Ys}
\end{equation*}
holds, with the implicit positive constants independent of $\{f_{\vj}\}_{\vj\in\Z^k}$, in each of the following cases:
\begin{equation*}
\Ys =\ell^{\vq}L^r,\quad\text{or}\quad
\Ys = L^r\ell^{\vq},\quad\text{or}\quad
\Ys = (L^{\vr}\ell^{\vq})_\pi\quad\text{with}\quad \vq\leq\vr.
\end{equation*}
\end{enumerate}
\end{proposition}

\begin{proof}
\eqref{dualCarl1}:
If $p\in(0,1]$, then, for almost all $x\in R$,
\begin{equation*}
\abs{[V]_{\aveL^p(R)}V(x)^{-1}}
\leq \Norm{[V]_{\aveL^p(R)}V^{-1}}{L^\infty(R)}
\sim \abs{[V]_{\aveL^p(R)}[V^{-1}]_{L^\infty(R)}}\lesssim 1
\end{equation*}
by the assumption that $V\in\A_p(\R^{\vec n})$. From this, \eqref{dualCarl1} is immediate.

If $p\in(1,\infty)$, then $U:=V^{-1}\in\A_q(\R^{\vec n})$ for $q:=p'\in(1,\infty)$, and the quantity of interest may be written as
\begin{equation}\label{dualCarl11}
\begin{split}
\abs{[V]_{\aveL^p(R)}V(x)^{-1}}
 &\leq\abs{[V]_{\aveL^p(R)}[V^{-1}]_{\aveL^{p'}(R)}}\cdot
 \Babs{[V^{-1}]_{\aveL^{p'}(R)}^{-1}V(x)^{-1}} \\
 &\lesssim\abs{[V^{-1}]_{\aveL^{p'}(R)}^{-1}V(x)^{-1}}\qquad\text{since }V\in\A_p(\R^{\vn}) \\
 &=\abs{U(x)[U]_{\aveL^q(R)}^{-1}}\qquad\text{by definition of $U:=V^{-1}$ and $q:=p'$}.
\end{split}
\end{equation}
From Corollary \ref{Ap always ACp}, it follows that $U\in\AC_q(\R^{\vn})$, and hence $U\in\AC_{q,s}(\R^{\vn})$ for some $s>q=p'$ by Definition \ref{Carl prop}. This means that
\begin{equation*}
\sup_{\Omega\in\Open(\R^{\vn})}
\bigg[\fint_\Omega \sup_{\genfrac{}{}{0pt}{}{R\in\D(\R^{\vn})}{x\in R\subset\Omega}}
\abs{U(x)[U]_{\aveL^{q}(R)}^{-1} }^s dx \bigg]^{\frac1s}
\lesssim 1.
\end{equation*}
By \eqref{dualCarl11}, this implies the claim \eqref{dualCarl10}.

\eqref{dualCarl2}: In either case of \eqref{r vs p}, part \eqref{dualCarl1} provides us with $s>r$ such that \eqref{dualCarl10} holds. This allows us to apply Theorem \ref{carl thm} with $r$ in place of $p$ to obtain the claimed conclusions.
\end{proof}

\section{Some Fourier analysis} \label{sec:fourier}

In this section, we recall some basic concepts and properties of multi-parameter Fourier analysis.
Let $\Sc(\R^n)$ be the space of all Schwartz functions on $\R^n$,
equipped with the well-known topology determined by a countable family of norms,
and let $\Sc'(\R^n)$ be the set of all continuous linear functionals on $\Sc(\R^n)$,
equipped with the weak-$*$ topology.
For any $f\in L^1(\R^n)$ and $\xi\in\R^{n}$, let
$$
\widehat{f}(\xi):=\int_{\mathbb{R}^n}f(x)e^{-\mathrm i x\cdot\xi}dx
$$
denote the \emph{Fourier transform} of $f$.
This agrees with the normalisation of the Fourier transform used,
for instance, in \cite[p.\,4]{FJW91} and \cite[p.\,452]{YY10},
and allows us to quote some lemmas from these works directly,
whereas using any other normalisation (such as with $2\pi$ in the exponent)
would also necessitate slight adjustments here and there in several other formulas.

On the product space $\R^{\vn}$, we define
\begin{align*}
\Sc_0(\R^{\vn})
:=\bigg\{\phi\in\Sc(\R^{\vn}) &:\
\forall\, i\in[k],\ \forall\,\alpha_i\in\N^{n_i},\ \forall\, x'\in\R^{\vn}\ominus\R^{n_i}, \\
&\quad\int_{\R^{n_i}} x_i^{\alpha_i}\phi(x'\oplus x_i)dx_i=0\bigg\}.
\end{align*}

\begin{definition}\label{def LPf}
We say that $\varphi$ is a \emph{Littlewood--Paley function} on $\mathbb R^n$ if $\varphi\in\mathscr S(\mathbb R^n)$ satisfies
\begin{equation*}
\widehat{\varphi}(\xi)\neq 0
\quad\text{ only if }\quad \frac12\leq|\xi|\leq 2
\end{equation*}
and
\begin{equation*}
\big|\widehat\varphi(\xi)\big|\geq c>0
\quad\text{ if }\quad\frac35\leq|\xi|\leq\frac53,
\end{equation*}
where $c$ is a positive constant independent of $\xi$.
We say that $\varphi$ is a \emph{Littlewood--Paley function} on $\mathbb R^{\vn}$ if $\varphi\in\mathscr S(\R^{\vn})$ satisfies
\begin{equation*}
\varphi(x)=\varphi(x_1,\ldots,x_k)=\prod_{i=1}^k\varphi^{(i)}(x_i),
\end{equation*}
where each $\varphi^{(i)}$ is a Littlewood--Paley function on $\mathbb R^{n_i}$. These $\varphi^{(i)}$ are called the component functions of $\varphi$.

We say that $(\varphi,\psi)$ is a \emph{Littlewood--Paley pair} on $\R^n$ if both $\varphi$ and $\psi$ are Littlewood--Paley functions on $\R^n$ and, in addition, they satisfy
\begin{equation*}
\sum_{j\in\mathbb{Z}}\overline{\widehat{\varphi}\Big(2^j\xi\Big)}
\widehat{\psi}\Big(2^j\xi\Big)=1\quad
\text{ if }\quad\xi\in\mathbb{R}^n\setminus\{\vnull \}.
\end{equation*}
We say that $(\varphi,\psi)$ is a \emph{Littlewood--Paley pair} on $\R^{\vn}$ if both $\varphi$ and $\psi$ are Littlewood--Paley functions on $\R^{\vn}$ and, in addition, their component functions $\varphi^{(i)}$ and $\psi^{(i)}$ are Littlewood--Paley pairs on $\R^{n_i}$ for each $i\in[k]$.
\end{definition}

If $\vj=(j_1,\ldots,j_k)\in\Z^k$ and $\xi=(\xi_1,\ldots,\xi_k)\in\R^{\vn}$, we denote
\begin{equation*}
2^{\vj}\xi:=(2^{j_1}\xi_1,\ldots,2^{j_k}\xi_k).
\end{equation*}
If $(\varphi,\psi)$ is a Littlewood--Paley pair on $\R^{\vn}$, it follows from the definition that
\begin{equation}\label{210}
\sum_{\vj\in\mathbb{Z}^k}\overline{\widehat{\varphi}\Big(2^{\vj}\xi\Big)}
\widehat{\psi}\Big(2^{\vj}\xi\Big)=1\quad
\text{ if }\quad\xi\in \dot\R^{\vn}:=(\R^{n_1}\setminus\{\mathbf 0 \})\times\cdots\times(\R^{n_k}\setminus\{\mathbf 0 \}).
\end{equation}

Let $\varphi$ be a function on $\mathbb{R}^{\vn}$.
For each $\vj\in\mathbb{Z}^k$ and $x\in\mathbb{R}^{\vn}$, let
\begin{equation*}
\varphi_{\vj}(x):=2^{\vj\cdot \vn}\varphi(2^{\vj}x),
\end{equation*}
where $\vn:=(n_1,\ldots,n_k)$.
For each $R=2^{-\vj}([0,1)^{\vn}+h)\in\D_{\vj}(\R^{\vn})$, where $\vj\in\Z^k$ and $h\in\Z^{\vn}$,
and for each $x\in\R^{\vn}$, let
\begin{equation*}
\varphi_R(x)
:=2^{\frac{\vj\cdot\vn}2}\varphi\Big(2^{\vj}x-h\Big)
=\abs{R}^{\frac 12}\varphi_{\vj}(x-a_R),
\end{equation*}
where $a_R:=2^{-\vj}h$ is the ``lower left'' corner of $R$.
We also write $\vj_R:=\vj$ to denote the ``\emph{level}'' of $R$.

\subsection{Calder\'on reproducing formulas}

\begin{lemma}\label{cald1}
If $(\varphi,\psi)$ is a Littlewood--Paley pair on $\R^{\vn}$, then
\begin{equation*}
f=\sum_{\vj\in\Z^k}\widetilde\varphi_{\vj}*\psi_{\vj}*f,\qquad\widetilde\varphi(x):=\overline{\varphi(-x)},
\end{equation*}
for all $f\in\Fs\in\{L^2(\R^{\vn}),\Sc_0(\R^{\vn}),\Sc_0'(\R^{\vn})\}$, in each case in the sense of convergence in the usual topology of the respective space.
\end{lemma}

\begin{proof}
For both $\Fs\in\{L^2(\R^{\vn}),\Sc_0(\R^{\vn})\}$, this is most conveniently checked on the Fourier transform side, recalling that the Fourier transform is an isomorphism on both $L^2(\R^{\vn})$ and $\Sc(\R^{\vn})$. From \eqref{210} we obtain that the Fourier transforms of both sides agree pointwise on $\dot\R^{\vn}$, and the convergence in the relevant topologies follows by some routine computations, noting that
\begin{equation} \label{hat S0}
\begin{split}
\widehat{\Sc_0}(\R^{\vn})
&:=\Big\{\widehat{f}:\ f\in\Sc_0(\R^{\vn}) \Big\} \\
&\phantom{:}=\Big\{f\in\Sc(\R^{\vn}):\ \partial^{\alpha}f(\xi)=0\text{ for all }\xi\in\R^{\vn}\setminus\dot\R^{\vn}\text{ and }\alpha\in\N^{\vn}\Big\}.
\end{split}
\end{equation}
For $f\in\Sc_0'(\R^{\vn})$, the convergence then follows by dualising against each $g\in\Sc_0(\R^{\vn})$ and using the convergence in this space.
\end{proof}

\begin{lemma}\label{cald2}
If $(\varphi,\psi)$ is a Littlewood--Paley pair on $\R^{\vn}$, then
\begin{equation}\label{cald20}
f=\sum_{\vj\in\Z^k}2^{-\vj\cdot\vn}\sum_{h\in\Z^{\vn}}
\widetilde\varphi_{\vj}*f(2^{-\vj}h)\psi_{\vj}(\cdot-2^{-\vj}h)
=\sum_{R\in\D(\R^{\vn})}\pair{f}{\varphi_R}\psi_R
\end{equation}
for all $f\in\Fs\in\{L^2(\R^{\vn}),\Sc_0(\R^{\vn}),\Sc_0'(\R^{\vn})\}$, in each case in the sense of convergence of the usual topology of the respective space.
\end{lemma}

\begin{proof}
Case $k=1$ can be found in \cite[Lemma 2.1]{YY10}; the general case can be proved similarly.
\end{proof}

\begin{lemma}\label{cald3}
Let $\varphi$ be a Littlewood--Paley function on $\R^{\vn}$, and let $\gamma\in\Sc(\R^{\vn})$ satisfy
\begin{equation*}
\widehat\gamma(\xi)=1\quad\text{if}\quad\xi\in [-2,2]^{\vn},\qquad
\supp\widehat\gamma\subset (-\pi,\pi)^{\vn}.
\end{equation*}
If $f\in\Sc_0'(\R^{\vn})$, then,
for every $x,y\in\mathbb R^{\vn}$,
\begin{equation}\label{cald30}
(\varphi_{\vj}*f)(x)=\sum_{R\in\D_{\vj}(\R^{\vn})}2^{-\vj\cdot\vn}(\varphi_{\vj}*f)(a_R+y)\gamma_{\vj}(x-a_R-y).
\end{equation}
\end{lemma}

\begin{proof}
Case $k=1$ can be found in \cite[Remark 3.16]{BHYY:1a};
the general case can be proved similarly.
\end{proof}

We also record from \cite[Lemma 3.7]{BHYY:1a} the following elementary inequality,
which is often convenient for handling the full scale of $q\in(0,\infty]$ in a unified way.
We typically apply this to situations, where $\sum_{i\in\mathscr I}a_i\lesssim 1$.

\begin{lemma}\label{one fits all}
Let $q\in(0,\infty]$ and $a_i,b_i\in[0,\infty)$ for all $i\in\mathscr I$
be numbers with some countable index set $\mathscr I$. Then
\begin{equation*}
\sum_{i\in\mathscr I}a_ib_i
\leq\Bigg[\sum_{i\in\mathscr I}a_i\Bigg]^{\frac{1}{q'}}
\Bigg[\sum_{i\in\mathscr I}a_i^{\min(1,q)}b_i^q\Bigg]^{\frac1q},
\end{equation*}
where the interpretation of the second factor is $\|b\|_\infty:=\sup_{i\in\mathscr I}b_i$ for $q=\infty$.
\end{lemma}

We now obtain a useful estimate:

\begin{lemma}\label{cald4}
Let $\varphi$ be a Littlewood--Paley function on $\R^{\vn}$.
For every $M,r\in(0,\infty)$, there exists a positive constant $C$
such that, for every $f\in\Sc_0'(\R^{\vn})$, $\vj\in\Z^k$,
$R\in\D_{\vj}(\R^{\vn})$, and $x\in R$,
\begin{equation*}
\sup_{y\in R}\abs{\varphi_{\vj}*f(y)}^r\leq C(\phi_{\vj}*\abs{\varphi_{\vj}*f}^r)(x),
\end{equation*}
where, for any $y\in\R^{\vn}$,
\begin{equation}\label{cald40}
\phi(y):=(1+\abs{y})^{-M},\quad\phi_{\vj}(y):=2^{\vj\cdot\vn}\phi(2^{\vj}y).
\end{equation}
\end{lemma}

\begin{proof}
Fix $M,r\in(0,\infty)$ and let $\gamma\in\Sc(\R^{\vn})$ be as in Lemma \ref{cald3}.
Using the fast decay of $\gamma$
and the identity of Lemma \ref{cald3},
we obtain the following estimate, for all $x,y\in\mathbb R^{\vn}$,
\begin{equation}\label{cald41}
\abs{\varphi_{\vj}*f(x)}
\lesssim\sum_{R\in\D_{\vj}(\R^{\vn})}\abs{(\varphi_{\vj}*f)(a_R+y)} [1+\abs{2^{\vj}(x-a_R-y)}]^{-M}.
\end{equation}
For $M>n$, we have
\begin{equation*}
\sum_{R\in\D_{\vj}(\R^{\vn})} [1+\abs{2^{\vj}(x-a_R-y)}]^{-M}
=\sum_{m\in\Z^k}[1+\abs{2^{\vj}(x-y)-m}]^{-M}\lesssim 1
\end{equation*}
and hence, by Lemma \ref{one fits all} with $q=r$, we conclude that
\begin{equation*}
\abs{\varphi_{\vj}*f(x)}^r
\lesssim\sum_{R\in\D_{\vj}(\R^{\vn})}\abs{(\varphi_{\vj}*f)(a_R+y)}^r [1+\abs{2^{\vj}(x-a_R-y)}]^{-M\min(1,r)}.
\end{equation*}
An integration over $y\in 2^{-\vj}[0,1)^{\vn}$, in which case $a_R+y$ ranges over $R$, gives
\begin{equation*}
\abs{\varphi_{\vj}*f(x)}^r
\lesssim(\phi_{\vj}*\abs{\varphi_{\vj}*f}^r)(x),
\end{equation*}
where $\phi_{\vj}$ is as in the statement of the lemma.
Then also
\begin{equation*}
\abs{\varphi_{\vj}*f(x+z)}^r
\lesssim(\phi_{\vj}(\cdot+z)*\abs{\varphi_{\vj}*f}^r)(x).
\end{equation*}
If $\abs{2^{\vj}z}\lesssim 1$, then $\phi_{\vj}(y+z)\sim\phi_{\vj}(y)$, and hence
\begin{equation*}
\sup_{\{z\in\R^{\vn}:\ \abs{2^{\vj}z}\leq c\}}\abs{\varphi_{\vj}*f(x+z)}^r\lesssim (\phi_{\vj}*\abs{\varphi_{\vj}*f}^r)(x).
\end{equation*}
For $x\in R\in\D_{\vj}(\R^{\vn})$, all $y\in R$ can be written as $x+z$ with $\abs{2^{\vj}z}\leq c$, and hence the bound just written also gives the bound claimed in the lemma. Since $M$ can be chosen as large as we like, we can also make $M\min\{1,r\}$ as large as we like, for any given $r\in(0,\infty)$.
This finishes the proof of Lemma \ref{cald4}.
\end{proof}

\begin{remark}
On the product space $\R^{\vn}=\R^{n_1}\times\cdots\times\R^{n_k}$, it depends a bit on the context whether it is more convenient to deal with isotropic decay functions $\phi^N(x):=(1+\abs{x})^{-N}$ as in Lemma \ref{cald4}
or product-type decay functions $\omega^{\vN}(x):=\prod_{i=1}^k(1+\abs{x_i})^{-N_i}$. In situations, like above, where we can choose the decay rate as large as we like, there is no essential difference. Indeed, if we denote the reciprocals by $\Phi^N(x):=(1+\abs{x})^{N}$ and $\Omega^{\vN}(x):=\prod_{i=1}^k(1+\abs{x_i})^{N_i}$, then
\begin{equation*}
\prod_{i=1}^k(1+\abs{x_i})\geq 1+\sum_{i=1}^k\abs{x_i}\geq 1+\abs{x},
\end{equation*}
hence $\Phi^N(x)\leq\Omega^{N\cdot\vone}(x)$ and thus $\omega^{N\cdot\vone}(x)\leq\phi^{N}(x)$. On the other hand, we also have $1+\abs{x_i}\leq 1+\abs{x}$, hence $\Omega^{\vN}(x)\leq\Phi^{\abs{\vN}}(x)$, and thus $\phi^{N}(x)\leq\omega^{\vN}(x)$ for any $\vN$ with $\abs{\vN}=N$.
\end{remark}

Let us also record the following variant of Lemma \ref{cald4}.

\begin{lemma}\label{cald5}
Let $\varphi$ be a Littlewood--Paley function on $\R^{\vn}$, and $\chi\in\Sc(\R^{\vn})$ be another test function.
For every $M,r\in(0,\infty)$, there is a constant $C\in(0,\infty)$ such that, for every $f\in\Sc_0'(\R^{\vn})$, $\vj\in\Z^k$, $R\in\D_{\vj}(\R^{\vn})$, and $x\in R$,
\begin{equation*}
\sup_{y\in R}\abs{\chi_{\vj}*\varphi_{\vj}*f(y)}^r\leq C(\phi_{\vj}*\abs{\varphi_{\vj}*f}^r)(x),
\end{equation*}
where $\phi$ is as in \eqref{cald40}.
\end{lemma}

\begin{proof}
Starting as in the proof of Lemma \ref{cald4}, we fix some $\gamma\in\Sc(\R^{\vn})$ as in Lemma \ref{cald3}. But then, instead of directly estimating as in the proof of Lemma \ref{cald4}, we first take the convolution with $\chi_{\vj}$ of both sides of the identity \eqref{cald30} guaranteed by Lemma \ref{cald3}, thus obtaining
\begin{equation}\label{cald51}
(\chi_{\vj}*\varphi_{\vj}*f)(x)
=\sum_{R\in\D_{\vj}(\R^{\vn})}2^{-\vj\cdot\vn}(\varphi_{\vj}*f)(a_R+y)
(\chi_{\vj}*\gamma_{\vj})(x-a_R-y).
\end{equation}
But $\chi_{\vj}*\gamma_{\vj}=(\chi*\gamma)_{\vj}$, and the new test function $\chi*\gamma\in\Sc(\R^{\vn})$, like $\gamma$, also has polynomial decay of any desired order. Hence
\begin{equation*}
2^{-\vj\cdot\vn}\abs{(\chi_{\vj}*\gamma_{\vj})(x-a_R-y)}
\lesssim(1+\abs{2^{\vj}(x-a_R-y)})^{-M}.
\end{equation*}
Substituting this into \eqref{cald51}, we obtain
\begin{equation*}
\abs{(\chi_{\vj}*\varphi_{\vj}*f)(x)}
\lesssim\sum_{R\in\D_{\vj}(\R^{\vn})}\abs{(\varphi_{\vj}*f)(a_R+y)} (1+\abs{2^{\vj}(x-a_R-y)})^{-M},
\end{equation*}
which is the same upper bound as that obtained for $\abs{(\varphi_{\vj}*f)(x)}$ in \eqref{cald41} in the proof of Lemma \ref{cald4}. From this point on, one can repeat the same proof as for Lemma \ref{cald4}.
This finishes the proof of Lemma \ref{cald5}.
\end{proof}

\subsection{Decay of Schwartz pairings}

\begin{lemma}\label{59}
Let $\varphi,\psi\in\Sc_0(\R^{\vn})$.
For every $M,N\in\mathbb{N}$, there exists a positive constant $C$,
depending only on $M,N$ and $\varphi,\psi$, such that,
for all $\vi,\vj\in\Z^k$ and $x\in\R^{\vn}$,
\begin{equation*}
\abs{(\varphi_{\vi}*\psi_{\vj})(x)}
\leq C2^{-\abs{\vi-\vj}M}\omega_{\vi\wedge\vj}^{(N)}(x)
\leq C2^{-\abs{\vi-\vj}M}\phi_{\vi\wedge\vj}^{(N)}(x),
\end{equation*}
where
\begin{equation*}
\omega_{\vi}^{(N)}(x):=2^{\vi\cdot\vn}\omega^{(N)}(2^{\vi}x),\quad
\phi_{\vi}^{(N)}(x):=2^{\vi\cdot\vn}\phi(2^{\vi}x),
\end{equation*}
and
\begin{equation*}
\omega^{(N)}(x):=\prod_{h=1}^k(1+\abs{x_h})^{-N},
\quad\phi^{(N)}(x):=(1+\abs{x})^{-N}.
\end{equation*}
\end{lemma}

\begin{proof}
For any multi-index $\alpha=(\alpha_1,\ldots,\alpha_k)\in\N^{\vn}$, we have
\begin{equation*}
(ix)^\alpha(\varphi_{\vi}*\psi_{\vj})(x)
=(2\pi)^{-n}\int_{\R^{\vn}}\widehat\varphi(2^{-\vi}\xi)
\widehat\psi(2^{-\vj}\xi)(ix)^\alpha e^{i x\cdot\xi} \, d\xi.
\end{equation*}
Since $(ix)^\alpha e^{ix\cdot\xi}=\partial_\xi^\alpha e^{ix\cdot\xi}$, we may integrate by parts to obtain
\begin{equation*}
(ix)^\alpha(\varphi_{\vi}*\psi_{\vj})(x)
=\frac{(-1)^{\abs{\alpha}}}{(2\pi)^{n}}\sum_{\beta+\gamma=\alpha}\binom{\alpha}{\beta}\int_{\R^{\vn}}
2^{-\vi\cdot\beta}(\partial^\beta\widehat\varphi)(2^{-\vi}\xi)
2^{-\vj\cdot\gamma}(\partial^\gamma\widehat\psi)(2^{-\vj}\xi)
e^{i x\cdot\xi} \, d\xi,
\end{equation*}
where, with slight misuse of notation, the dot product  $\vi\cdot\beta$ of $\vi\in\Z^k$ and $\beta\in\N^{\vn}$ is interpreted as
\begin{equation*}
\vi\cdot\beta:=\sum_{h=1}^k i_h\abs{\beta_h},
\end{equation*}
and similarly with $\vj\cdot\gamma$.

Since $\varphi\in\Sc_0(\R^{\vn})$, its Fourier transform $\widehat\varphi$ vanishes to arbitrary order on the coordinate hyperplanes, i.e., whenever for some $h\in[k]$ we have $\xi_h=\mathbf 0\in\R^{n_h}$ for $\xi=(\xi_1,\ldots,\xi_k)\in\R^{\vn}$. And of course $\widehat\varphi$ and all its derivatives decay faster than any polynomial at infinity. Hence, for any $M\in\N$, we have
\begin{equation*}
\big|(\partial^\beta\widehat\varphi)(2^{-\vi}\xi)\big|\lesssim\prod_{h=1}^k\min
\Big\{\abs{2^{-i_h}\xi_h}^M,\abs{2^{-i_h}\xi_h}^{-M}\Big\},
\end{equation*}
and a similar bound is valid for $(\partial^\gamma\widehat\psi)(2^{-\vj}\xi)$. Writing $\vmu:=\min(\vi,\vj)$, $\vnu:=\max(\vi,\vj)$ (componentwise), it follows that
\begin{align*}
&\int_{\R^{\vn}} \big|(\partial^\beta\widehat\varphi)(2^{-\vi}\xi)\big|
\big|(\partial^\gamma\widehat\psi)(2^{-\vj}\xi)\big| d\xi \\
&\quad\lesssim\prod_{h=1}^k
\int_{\R^{n_h}}\min\Big\{\abs{2^{-i_h}\xi_h}^M,\abs{2^{-i_h}\xi_h}^{-M}\Big\}
\min\Big\{\abs{2^{-j_h}\xi_h}^M,\abs{2^{-j_h}\xi_h}^{-M}\Big\} d\xi_h \\
&\quad=\prod_{h=1}^k
\bigg( \int_{\abs{\xi_h}<2^{\mu_h}}
\abs{2^{-\mu_h}\xi_h}^M\abs{2^{-\nu_h}\xi_h}^M d\xi_h
+\int_{2^{\mu_h}<\abs{\xi_h}<2^{\nu_h}}
\abs{2^{-\mu_h}\xi_h}^{-M}\abs{2^{-\nu_h}\xi_h}^M d\xi_h \Big. \\
&\quad\quad \Big.
+\int_{\abs{\xi_h}>2^{\nu_h}}
\abs{2^{-\mu_h}\xi_h}^{-M}\abs{2^{-\nu_h}\xi_h}^{-M} d\xi_h \bigg)  \\
&\quad\lesssim\prod_{h=1}^k
\Big[2^{-(\mu_h+\nu_h)M}2^{\mu_h(2M+n_h)}+2^{(\mu_h-\nu_h)M}2^{\nu_h n_h}+2^{(\mu_h+\nu_h)M}2^{\nu_h(n_h-2M)}\Big] \\
&\quad\sim\prod_{h=1}^k 2^{(\mu_h-\nu_h)M}2^{\nu_h n_h}
=2^{-\abs{\vi-\vj}M}2^{\vnu\cdot\vn}
=2^{-\abs{\vi-\vj}M}2^{(\vnu-\vmu)\cdot\vn}2^{\vmu\cdot\vn}
\leq 2^{-\abs{\vi-\vj}(M-\abs{\vn})}2^{\vmu\cdot\vn}.
\end{align*}
Above, we can pick $M$ as large as we like, hence also $M-\abs{\vn}$ as large as we like. Thus, redefining $M$, we can simply write the result as
\begin{equation*}
\int_{\R^{\vn}} \big|(\partial^\beta\widehat\varphi)(2^{-\vi}\xi)\big|
\big|(\partial^\gamma\widehat\psi)(2^{-\vj}\xi)\big| d\xi
\lesssim 2^{-\abs{\vi-\vj}M}2^{\min\{\vi,\vj\}\cdot\vn}.
\end{equation*}
Substituting back, it follows that
\begin{equation*}
\abs{(ix)^\alpha(\varphi_{\vi}*\psi_{\vj})(x)}
\lesssim 2^{-\abs{\vi-\vj}M}2^{\min\{\vi,\vj\}\cdot\vn}\sum_{\beta+\gamma=\alpha}
\binom{\alpha}{\beta}2^{-\vi\cdot\beta}2^{-\vj\cdot\gamma},
\end{equation*}
where
\begin{equation*}
\begin{split}
\sum_{\beta+\gamma=\alpha}\binom{\alpha}{\beta}2^{-\vi\cdot\beta}2^{-\vj\cdot\gamma}
&=\prod_{h=1}^k\prod_{l=1}^{n_h}\sum_{\beta_{h,l}=0}^{\alpha_{h,l}}
\binom{\alpha_{h,l}}{\beta_{h,l}}2^{-i_h\beta_{h,l}}2^{-j_h(\alpha_{h,l}-\beta_{h,l})} \\
&=\prod_{h=1}^k\prod_{l=1}^{n_h}(2^{-i_{h}}+2^{-j_{h}})^{\alpha_{h,l}}
\sim\prod_{h=1}^k\prod_{l=1}^{n_h}(2^{-\mu_{h}})^{\alpha_{h,l}}
=2^{-\vmu\cdot\alpha}.
\end{split}
\end{equation*}
Hence, for every $\alpha\in\N^{\vn}$,
\begin{equation*}
\abs{(2^{\vmu}x)^\alpha(\varphi_{\vi}*\psi_{\vj})(x)}
\lesssim 2^{-\abs{\vi-\vj}M}2^{\vmu\cdot\vn}.
\end{equation*}

We now make a choice of $\alpha$ depending on $x$. If $\abs{2^{\mu_h}x_h}\leq 1$, we take $\alpha_h:=0$ so that $\abs{(2^{\mu_h}x_h)^{\alpha_h}}=1$. If $\abs{2^{\mu_h}x_h}>1$, let $l$ be a coordinate with $\abs{2^{\mu_h}x_{h,l}}\sim\abs{2^{\mu_h}x_h}$. Then we take $\alpha_h:=Ne_{h,l}$, where $e_{h,l}$ is the $l$-th unit vector in $\R^{n_h}$. It follows that $\abs{(2^{\mu_h}x_h)^{\alpha_h}}=\abs{2^{\mu_h}x_{h,l}}^N\sim \abs{2^{\mu_h}x_h}^N$. Hence, in either case, we obtain
\begin{equation*}
\abs{(2^{\mu_h}x_h)^{\alpha_h}}\gtrsim(1+\abs{2^{\mu_h}x_h})^N,
\end{equation*}
and thus
\begin{equation*}
\begin{split}
\abs{(2^{\vmu}x)^{\alpha}}
=\Bigg|\prod_{h=1}^k(2^{\mu_h}x_h)^{\alpha_h}\Bigg|
&\gtrsim\prod_{h=1}^k (1+\abs{2^{\mu_h}x_h})^N=\frac{1}{\omega^{(N)}(2^{\vmu}x)} \\
&\geq\Bigg(1+\sum_{h=1}^k\abs{2^{\mu_h}x_h}\Bigg)^N
\geq(1+\abs{2^{\vmu}x})^N=\frac{1}{\phi^{(N)}(2^{\vmu}x)}.
\end{split}
\end{equation*}
Hence we have proved that
\begin{equation*}
\frac{1}{\phi^{(N)}(2^{\vmu}x)}\abs{(\varphi_{\vi}*\psi_{\vj})(x)}
\leq\frac{1}{\omega^{(N)}(2^{\vmu}x)}\abs{(\varphi_{\vi}*\psi_{\vj})(x)}
\lesssim 2^{-\abs{\vi-\vj}M}2^{\vmu\cdot\vn}.
\end{equation*}
This finishes the proof of Lemma \ref{59}.
\end{proof}

\begin{lemma}\label{59x}
Let $\varphi,\psi\in\Sc_0(\R^{\vn})$.
For every $M,N\in\mathbb{N}$, there exists a positive constant $C$,
depending only on $M,N$ and $\varphi,\psi$, such that,
for all $\vi,\vj\in\Z^k$, $P\in\D_{\vi}(\R^{\vn})$, and $R\in\D_{\vj}(\R^{\vn})$,
\begin{equation*}
\abs{\pair{\varphi_P}{\psi_R}}
\lesssim\abs{P}^{\frac 12} \abs{R}^{\frac 12}
2^{-\abs{\vi-\vj}M} \omega_{\vi\wedge\vj}^{(N)}(a_P-a_R),
\end{equation*}
where $\omega^{(N)}$ is as in Lemma \ref{59}.
\end{lemma}

\begin{proof}
Recall that $\varphi_P(x)=\abs{P}^{1/2}\varphi_{\vi}(x-a_P)$, and $\psi_R$ has a similar definition. Then
\begin{equation*}
\pair{\varphi_{\vi}(\cdot-a_P)}{\psi_{\vj}(\cdot-a_R)}
=\int_{\R^{\vn}}\varphi_{\vi}(y)\psi_{\vj}(y+a_P-a_R)dy
=(\varphi_{\vi}*\widetilde\psi_{\vj})(a_R-a_P),
\end{equation*}
where $\widetilde\psi(y):=\psi(-y)$ satisfies the same assumptions as $\psi$. The claim is then immediate from Lemma \ref{59}.
This finishes the proof of Lemma \ref{59x}.
\end{proof}

\section{Multi-parameter Besov--Triebel--Lizorkin-type spaces}
\label{BF type spaces}

In this section, we introduce the matrix-weighted
multi-parameter Besov-type and Triebel--Lizorkin-type spaces $\Atau(\mathbb V)$
and establish their $\varphi$-transform characterization.
To this end, we first prove the equivalent relations among three quasi-norms
associated with $\Atau(\mathbb V)$ in Subsection \ref{W and AQ 1}.
Using these equivalent relations, we then obtain the
$\varphi$-transform characterization of $\Atau(\mathbb V)$ in Subsection \ref{phi-transform}.

Although the one-parameter case has been well studied in \cite{BHYY:1a},
the situation is different in the multi-parameter case,
since the definition of $\Atau(\mathbb V)$ involves
taking the supremum over all open sets rather than over all rectangles.
Hence, we need to find some new methods.

\subsection{Different quasi-norms and their relations}
\label{W and AQ 1}

It is customary in the theory of Besov and Triebel--Lizorkin spaces to give two parallel sets of definitions, one of function
(actually: distribution) spaces, and another one of sequence spaces. Both of them are given in terms of a quasi-norm of a
sequence of functions that are derived from the original objects. For an efficient discussion avoiding unnecessary repetition, we give a third (or first, to emphasise its fundamental
role in generating the other two) definition of quasi-norms on sequences of functions. This is particularly useful in the present discussion of matrix-weighted spaces, where these weights may be
inserted into the quasi-norms in a variety of ways. Here is the formal definition:

\begin{definition}\label{d3.1}
Let $\vs\in\R^k$, $\tau\in[0,\infty)$, $\vp:=p\cdot\vone$ with $p\in(0,\infty)$, and $\vq\in(0,\infty]^k$.
Let $\pi \in S_{[2k]}$ be an admissible permutation of $(\vp,\vq)$ in the sense of Definition \ref{adm}.

For each $\vj\in\Z^k$, let $V_{\vj}\in L_{\loc}^p(\R^{\vec n};\C^{m\times m})$, and $\V:=\{V_{\vj}\}_{\vj\in\Z^k}$.
For $\Omega\in\Open(\R^{\vn})$, let
\begin{equation*}
\D_{\vj}(\Omega):=\{R\in\D_{\vj}(\R^{\vn}):\ R\subset\Omega\}
\quad \text{and}\quad
\Omega_{\vj}:=\bigcup_{R\in\D_{\vj}(\Omega)}R\subset\Omega.
\end{equation*}
\begin{enumerate}[\rm(i)]
\item For any sequence of $\C^m$-valued functions $\{f_{\vj}\}_{\vj\in\Z^k}$, let
\begin{equation}\label{general}
\Norm{\{f_{\vj}\}_{\vj\in\Z^k}}{\aatau(\V)}:=\sup_{\Omega\in\Open(\R^{\vn})}|\Omega|^{-\tau}
\Big\|\Big\{\one_{\Omega_{\vj}}2^{\vj\cdot\vs}V_{\vj}f_{\vj}\Big\}_{\vj\in\Z^k}
\Big\|_{(L^{\vp}\ell^{\vq})_{\pi}}.
\end{equation}
Here and below, we adopt the convention that
$|\Omega|^{-\tau}=\infty$ whenever $|\Omega|=0$.

\item Let $\varphi$ be a Littlewood--Paley function on $\R^{\vn}$.
For each distribution $f\in\Sc_{0}'(\R^{\vn};\C^m)$, let
\begin{equation*}
\Norm{f}{\Atau(\mathbb V)}^\varphi :=\Norm{\{\varphi_{\vj}*f\}_{\vj\in\Z^k}}{\aatau(\V)}.
\end{equation*}

\item For each sequence $t:=\{t_Q\}_{Q\in\D(\R^{\vn})}$ of vectors $t_Q\in\C^m$, let
\begin{equation}\label{tj}
\Norm{t}{\atau(\mathbb V)} :=\Norm{\{t_{\vj}\}_{\vj\in\mathbb Z^k}}{\aatau(\V)},
\qquad t_{\vj}:=\sum_{R\in\D_{\vj}(\R^{\vec n})}{t}_R\widetilde{\one}_R,
\quad\widetilde{\one}_R:=\abs{R}^{-\frac 12}\one_R.
\end{equation}
\end{enumerate}
If $V_{\vj}\equiv V$ is the same function for all $\vj\in\Z^k$, we write just $V$ instead of $\V$ in these quasi-norms.
If $V\equiv 1$, we write $\R^{\vn}$ instead of $V$,
and in this case, the assumption $\vec p=p\cdot1$ can be relaxed to $\vec p\in(0,\infty)^k$.
\end{definition}

\begin{remark}\label{open vs D}
Even when specialised to the one-parameter case $k=1$,
this definition differs slightly from previous definitions of
Besov--Triebel--Lizorkin-type spaces (see \cite{YY08,YY10,YSY10} for the unweighted and \cite{BHYY:1a} for the matrix-weighted case).
That is, we take the supremum in \eqref{general} over all open sets, in contrast to the much smaller collection of dyadic cubes in \cite{BHYY:1a,YY08,YY10,YSY10}.
(While dyadic cubes, strictly speaking, are not open sets, estimating them by slightly larger open sets containing them, it is easy to see that the supremum over open sets in \eqref{general} always dominates the supremum over dyadic cubes. We will frequently use this fact in what follows.)
As a consequence of this difference, the two definitions differ in many
cases (see Theorem \ref{compare A}).
The reasons for taking the supremum over open sets
rather than rectangles in the definition of $\Atau(\R^{\vn})$ are two-fold.
It is well known that, when $k\geq 2$, the most natural multi-parameter BMO space
was introduced by Chang and Fefferman \cite{CF80} in terms of
taking the supremum over all open sets, which is the dual
of the multi-parameter Hardy space.
Note that when $k\geq 2$ the space $\dot F^{0,\frac12}_{2,2}(\R^{\vn})$ in this case coincides with
Chang--Fefferman BMO (see Theorem \ref{coin BMO}), which differs from the corresponding
rectangular version (see Theorem \ref{compare B}).
Another reason is that the dual space of the multi-parameter
Triebel--Lizorkin sequence space is precisely $\atau(\R^{\vn})$
instead of their rectangular version (see Corollary \ref{dual}).
The key factor behind these is that although any open set in
$\R^{\vn}$ can be decomposed into the union of maximal dyadic
rectangles, they are not nested and hence may have too much
overlaps, which leads to their futility.
Thus, to capture the whole information of functions or distributions
on $\R^{\vn}$, we need to take the supremum over all open sets
which can cover all regions of arbitrary shape,
including collections of rectangles at
different scales, avoiding the loss
of control caused by the non-nestedness of rectangles.

Nevertheless, the two definitions coincide in important special cases. First, when $\tau=0$, corresponding to usual Besov--Triebel--Lizorkin (instead of Besov--Triebel--Lizorkin-\emph{type}) spaces, the supremum is obviously achieved by taking a sequence of sets approaching the whole $\R^{\vn}$ (or, in the case of dyadic cubes, the quadrant that gives the maximal value), which amounts to formally substituting $\Omega=\R^{\vn}$.

Another important case consists of Triebel--Lizorkin-type spaces with $(L^{\vp}\ell^{\vq})_\pi=L^p\ell^{\vq}$ at the critical value $\tau=\frac1p$. When $k=1$, letting $\mathscr D^*(\Omega)$ be the collection of maximal dyadic cubes $P\subset\Omega$, it follows from their disjointness that
\begin{equation*}
\Big\|\Big\{\one_{\Omega_{j}}2^{js}V_{j}f_{j}\Big\}_{j\in\Z}\Big\|_{L^p\ell^q}^p
=\sum_{P\in\mathscr D^*(\Omega)}
\Big\|\Big\{\one_{P}2^{js}V_{j}f_{j}\Big\}_{\vj\in\Z}\Big\|_{L^p\ell^q}^p
\leq K^p\sum_{P\in\mathscr D^*(\Omega)}\abs{P}
= K^p\abs{\Omega},
\end{equation*}
where
\begin{equation*}
K:=\sup_{P\in\mathscr D(\R^n)}\abs{P}^{-\frac1p}\Big\|\Big\{\one_{P}2^{js}V_{j}f_{j}\Big\}_{j\in\Z}\Big\|_{L^p\ell^q}.
\end{equation*}
Hence the supremum over all open $\Omega$ is the same as the supremum over all $P\in\mathscr D(\R^n)$ also in this case.
\end{remark}

\begin{definition}\label{ave}
For every $p\in(0,\infty]$ and $\vj\in\Z^k$,
we define two versions of local $L^p$ averaging operators over the dyadic rectangles $\D_{\vj}(\R^{\vn})$ as follows:
\begin{equation*}
\Ave{f}{\vj}{p}:=\sum_{R\in\D_{\vj}(\R^{\vn})}\one_R\Norm{f}{\aveL^p(R)},\qquad
\Red{V}{\vj}{p}:=\sum_{R\in\D_{\vj}(\R^{\vn})}\one_R[V]_{\aveL^p(R)}.
\end{equation*}
We denote by $[V]_p:=\{\Red{V}{\vj}{p}\}_{\vj\in\Z^k}$ the sequence of functions obtained by the second averaging method. In particular, $\Atau([V]_p)$ and $\atau([V]_p)$ will denote, respectively, $\Atau(\V)$ and $\atau(\V)$ as in Definition \ref{d3.1} for $\V=\{V_{\vj}\}_{\vj\in\Z^k}$ with $V_{\vj}=\Red{V}{\vj}{p}$.

We will also identify $[V]_p$ with the sequence $\{[V]_{\aveL^p(R)}\}_{R\in\D(\R^{\vn})}$ or with the diagonal matrix $\{[V]_{\aveL^p(R)}\delta_{P,R}\}_{P,R\in\D(\R^{\vn})}$, both of which carry the same information as the sequence of piecewise constant functions $\{\Red{V}{\vj}{p}\}_{\vj\in\Z^k}$.
\end{definition}

The following theorem relating the three types of quasi-norms is the main result of this subsection.

\begin{theorem}\label{3 norms thm}
Let $\vs\in\R^k$, $\tau\in[0,\infty)$, $\vp\in(0,\infty)^k$, and $\vq\in(0,\infty]^k$.
Assume that $\pi\in S_{[2k]}$ is admissible for $(\vp,\vq)$.
For each $\vj\in\Z^k$, let either $f_{\vj}=\varphi_{\vj}*f$ for some Littlewood--Paley function $\varphi$ and $f\in\Sc_0'(\R^{\vn};\mathbb C^m)$, or $f_{\vj}=t_{\vj}$ be as in \eqref{tj}.
Then we have the following conclusions:
\begin{enumerate}[\rm(i)]
\item\label{it:f-avef}
Under the above assumptions,
\begin{equation*}
\bNorm{\big\{f_{\vj}\big\}_{\vj\in\Z^k}}{\aatau}
\sim\BNorm{\Big\{\Ave{f_{\vj}}{\vj}{\infty}\Big\}_{\vj\in\Z^k}}{\aatau}
\end{equation*}
with the positive equivalence constants independent of
$\{f_{\vj}\}_{\vj\in\Z^k}$.

\item\label{it:[V]<V}
If we further assume that $\vp =p\cdot\vone$ for some $p\in(0,\infty)$ and $V\in\A_p(\R^{\vn})$, then
\begin{equation}\label{3 norms equ}
\bNorm{\big\{\Red{V}{\vj}{p}f_{\vj}\big\}_{\vj\in\Z^k}}{\aatau}
\sim\BNorm{\Big\{\Ave{\Red{V}{\vj}{p}f_{\vj}}{\vj}{\infty}\Big\}_{\vj\in\Z^k}}{\aatau}
\lesssim\bNorm{\big\{ V f_{\vj} \big\}_{\vj\in\Z^k} }{\aatau}
\end{equation}
with the implicit positive constants independent of
$\{f_{\vj}\}_{\vj\in\Z^k}$.

\item\label{it:[V]=V}
The final $\lesssim$ in \eqref{3 norms equ} may be replaced by $\sim$ under any of the following additional assumptions:
\begin{equation}\label{key cases}
(L^{\vp}\ell^{\vq})_\pi=\ell^{\vq}L^p,\quad\text{or}\quad
(L^{\vp}\ell^{\vq})_\pi=L^p\ell^{\vq},\quad\text{or}\quad
\vq\leq\vp.
\end{equation}
\end{enumerate}
\end{theorem}

\begin{corollary}\label{3 norms cor}
Under the assumptions of Theorem \ref{3 norms thm},
including \eqref{key cases} in the weighted case,
for every $f\in\Sc_0'(\R^{\vn};\C^m)$,
the following equivalences hold:
\begin{equation}\label{3 norms A}
\Norm{f}{\Atau([V]_p)}^\varphi\sim
\BNorm{\Big\{\Ave{\Red{V}{\vj}{p}(\varphi_{\vj}*f)}{\vj}{\infty}\Big\}_{\vj\in\Z^k}}{\aatau}\sim
\Norm{f}{\Atau(V)}^\varphi\qquad
\end{equation}
where the positive equivalence constants are independent of $f$.
Moreover, for every $t\in(\C^m)^{\D(\R^{\vn})}$,
\begin{equation}\label{2 norms a}
\Norm{t}{\atau([V]_p)}\sim
\Norm{t}{\atau(V)}
\end{equation}
where the positive equivalence constants are independent of $t$.
\end{corollary}

\begin{proof}[Proof of Corollary \ref{3 norms cor} assuming Theorem \ref{3 norms thm}]
The first chain of quasi-norm equivalences corresponds to the case of Theorem \ref{3 norms thm} with $f_{\vj}=\varphi_{\vj}*f$, and the second to the case where $f_{\vj}=t_{\vj}$ is as in \eqref{tj}, thus constant on each $R\in\D_{\vj}(\R^{\vn})$. [The equivalence of the other quasi-norms to the one involving $\Ave{\ }{\vj}{\infty}$ is also true but uninteresting in the second case, since $\Ave{\Red{V}{\vj}{p}t_{\vj}}{\vj}{\infty}=\abs{\Red{V}{\vj}{p}t_{\vj}}$, when the function inside $\Ave{\ }{\vj}{\infty}$ is already constant on each $R\in\D_{\vj}(\R^{\vn})$.]
This finishes the proof of Corollary \ref{3 norms cor}.
\end{proof}

\begin{remark}\label{3 norms rem}
To the best of our knowledge, Theorem \ref{3 norms thm} and Corollary \ref{3 norms cor} appear to be new in the multi-parameter case (i.e., $k>1$), even for scalar weights. For $k=1$, \eqref{3 norms A} and \eqref{2 norms a} roughly correspond to \cite[Theorems 3.9 and 3.27]{BHYY:1a}, but recall from Remark \ref{open vs D} that the definitions here and in \cite{BHYY:1a} differ even in this case.
\end{remark}

The first ``$\lesssim$'' (actually, ``$\leq$'') in \eqref{3 norms equ} is essentially trivial, and holds for arbitrary functions $V_{\vj}$ in place of $\Red{V}{\vj}{p}$: indeed, $\abs{g}\leq\Ave{g}{\vj}{\infty}$, and hence
\begin{equation}\label{equi easy}
\BNorm{\big\{V_{\vj}f_{\vj} \big\}_{\vj\in\Z^k} }{\aatau}
\leq\BNorm{\big\{\Ave{V_{\vj}f_{\vj}}{\vj}{\infty} \big\}_{\vj\in\Z^k} }{\aatau}.
\end{equation}
However, the other bounds will require more work, and the following lemma will play an important role:

\begin{lemma}\label{convo bd A bd}
Let $\vs\in\R^k$, $\tau\in[0,\infty)$, $\vp\in(0,\infty)^k$, $\vq\in(0,\infty]^k$,
and $\pi\in S_{[2k]}$ be an admissible permutation of $(\vp,\vq)$.
Suppose that two sequences $\{f_{\vj}\}_{\vj\in\Z^k}$ and $\{g_{\vj}\}_{\vj\in\Z^k}$
of measurable functions on $\R^{\vn}$ satisfy:
there exist $r\in(0,\min(\vp,\vq))$,
$C\in(0,\infty)$, and $M\in(n(1+\tau),\infty)$ such that,
for every $\vj\in\Z^k$ and $x\in\R^{\vn}$,
\begin{equation}\label{143}
\abs{f_{\vj}(x)}^r\leq C \Big( \phi_{\vj}*\abs{g_{\vj}}^r \Big) (x),
\end{equation}
where $\phi(\cdot):=(1+\abs{\cdot})^{-M}$.
Then there exists a positive constant $\widetilde{C}$, depending only on the parameters of the space and on $C,M,r$, such that
\begin{equation*}
\Norm{\{f_{\vj}\}_{\vj\in\Z^k}}{\aatau}
\leq \widetilde{C} \Norm{\{g_{\vj}\}_{\vj\in\Z^k}}{\aatau}.
\end{equation*}
\end{lemma}

\begin{proof}
By considering $2^{\vj\cdot\vs}f_{\vj}$ and $2^{\vj\cdot\vs}g_{\vj}$ in place of $f_{\vj}$ and $g_{\vj}$, the general case follows from the special case $\vs=\vnull$; thus, we will assume that $\vs=\vnull$ in the rest of the proof.

Recall that
\begin{equation*}
\D_{\vj}(\Omega)=\{R\in\D_{\vj}(\R^{\vn}):\ R\subset\Omega\}
\ \ \text{and}\ \
\Omega_{\vj}=\bigcup_{R\in\D_{\vj}(\Omega)} R\subset\Omega.
\end{equation*}
We introduce the auxiliary sets
\begin{equation*}
A_0:=[-1,1]^{\vn},\qquad A_h:=[-2^h,2^h]^{\vn}\setminus[-2^{h-1},2^{h-1}]^{\vn},\quad h\in\Z_+.
\end{equation*}
These lead to the splitting
\begin{equation*}
\phi_{\vj}=\sum_{h=0}^\infty\one_{2^{-\vj}A_h}\phi_{\vj},\qquad
\one_{\Omega_{\vj}}\abs{f_{\vj}}^r\leq \one_{\Omega_{\vj}}\sum_{h=0}^\infty (\one_{2^{-\vj}A_h}\phi_{\vj})*\abs{g_{\vj}}^r.
\end{equation*}
For $R\in\D_{\vj}(\Omega)$, we investigate
\begin{equation*}
\one_R(x)(\one_{2^{-\vj}A_h}\phi_{\vj})*\abs{g_{\vj}}^r(x)
=\one_R(x)\int_{2^{-\vj}A_h}2^{\vj\cdot\vn}\phi(2^{\vj}y)\abs{g_{\vj}(x-y)}^r dy.
\end{equation*}
If $x\in R$ and $y\in 2^{-\vj}A_h$, then $x-y\in 3 R^{(h)}$, where $R^{(h)}\in\D_{\vj-h\vone}(\R^{\vn})$ is the $h$-th dyadic ancestor of $R$ in the dyadic system $\{\D_{\vj+i\vone}\}_{i\in\Z}$.
Since $R\subset\Omega$, every $z\in 3 R^{(h)}$ satisfies
\begin{equation*}
\mathcal M_{\vn}\one_{\Omega}(z)
\geq\fint_{3 R^{(h)}}\one_R
=\frac{\abs{R}}{\abs{3 R^{(h)}}}
=3^{-n}2^{-hn}.
\end{equation*}
Thus,
\begin{equation*}
3 R^{(h)}
\subset\Omega^h
:=\{x\in\mathbb R^{\vn}:\
(\mathcal M_{\vn}\one_{\Omega})(x)\geq 3^{-n}2^{-hn}\}.
\end{equation*}
By the strong maximal inequality (see \cite{CF75}),
we conclude that, for all $\lambda \in(0,\infty)$,
\begin{equation*}
|\{x\in\mathbb R^{\vn}:\
(\mathcal M_{\vn}\one_{\Omega})(x)\geq \lambda\}|
\lesssim\int_{\R^{\vn}}\frac{\one_\Omega}{\lambda}
\bigg( 1+\log_+^{k-1}\frac{\one_\Omega}{\lambda} \bigg)
=\frac{\abs{\Omega}}{\lambda}
\bigg( 1+\log_+^{k-1}\frac{1}{\lambda} \bigg),
\end{equation*}
and hence
\begin{equation}\label{Oh bd}
\abs{\Omega^h}\lesssim 2^{hn}(1+h^{k-1})\abs{\Omega}.
\end{equation}

Since $3 R^{(h)}$ is a union of $3^n 2^{hn}$ rectangles of $\D_{\vj}(\R^{\vn})$, we even have
\begin{equation*}
3 R^{(h)}\subset\Omega^h_{\vj}:=(\Omega^h)_{\vj},
\end{equation*}
where the latter set has the same relation to $\Omega^h$ as $\Omega_{\vj}$ has to the original $\Omega$. Altogether, we observed that
\begin{equation*}
x-y\in 3 R^{(h)}\subset\Omega^h_{\vj},
\end{equation*}
and hence
\begin{equation*}
\one_R(x)\Big( (\one_{2^{-\vj}A_h}\phi_{\vj})*\abs{g_{\vj}}^r\Big) (x)
=\one_R(x) \Big( (\one_{2^{-\vj}A_h}\phi_{\vj})*(\one_{\Omega^h_{\vj}} \abs{g_{\vj}^r}) \Big) (x).
\end{equation*}
Summing over all $R\in\D_{\vj}(\Omega)$, we deduce that
\begin{equation}\label{insert Ohj}
\one_{\Omega_{\vj}} \Big((\one_{2^{-\vj}A_h}\phi_{\vj})*\abs{g_{\vj}}^r\Big)
=\one_{\Omega_{\vj}} \Big((\one_{2^{-\vj}A_h}\phi_{\vj})*(\one_{\Omega^h_{\vj}}\abs{g_{\vj}}^r)\Big).
\end{equation}

We then investigate, for every $x\in\mathbb R^{\vn}$,
\begin{equation*}
(\one_{2^{-\vj}A_h}\phi_{\vj})(x)
=\one_{A_h}(2^{\vj}x)2^{\vj\cdot\vn}\phi(2^{\vj}x)
=(\one_{A_h}\phi)_{\vj}(x),
\end{equation*}
where
\begin{equation*}
(\one_{A_h}\phi)(x)
=\one_{A_h}(x)(1+\abs{x})^{-M}
\lesssim\one_{A_h}(x)2^{-hM}\leq 2^{-hn}\one_{[-2^h,2^h]^{\vn}}(x)2^{h(n-M)}.
\end{equation*}
Since
\begin{equation*}
( 2^{-hn}\one_{[-2^h,2^h]^{\vn}})_{\vj}*g(x)\lesssim \mathcal M_{\vn}g(x),
\end{equation*}
it follows that
\begin{equation*}
(\one_{2^{-\vj}A_h}\phi_{\vj})*\Big(\one_{\Omega^h_{\vj}}\abs{g_{\vj}}^r\Big)
\lesssim 2^{h(n-M)} \mathcal M_{\vn}\Big(\one_{\Omega^h_{\vj}}\abs{g_{\vj}}^r\Big).
\end{equation*}
Substituting this into \eqref{insert Ohj}, we obtain
\begin{equation*}
\one_{\Omega_{\vj}} \Big((\one_{2^{-\vj}A_h}\phi_{\vj})*\abs{g_{\vj}}^r\Big)
\lesssim 2^{h(n-M)} \mathcal M_{\vn}\Big(\one_{\Omega^h_{\vj}}\abs{g_{\vj}}^r\Big).
\end{equation*}
Summing over $h$ then gives
\begin{equation}\label{insert Ohj 2}
\one_{\Omega_{\vj}}(\phi_{\vj}*\abs{g_{\vj}}^r)
\lesssim\sum_{h=0}^\infty 2^{h(n-M)} \mathcal M_{\vn}\Big(\one_{\Omega^h_{\vj}}\abs{g_{\vj}}^r\Big).
\end{equation}

We will use this estimate to bound the quasi-norm
\begin{equation*}
\begin{split}
\Norm{\{f_{\vj}\}_{\vj\in\Z^k}}{\antau}^r
&=\bigg[\sup_{\Omega\in\Open(\R^{\vn})}\abs{\Omega}^{-\tau}
\Norm{\{\one_{\Omega_{\vj}}f_{\vj}\}_{\vj\in\Z^k}}{(L^{\vp}\ell^{\vq})_\pi}\bigg]^r\\
&=\sup_{\Omega\in\Open(\R^{\vn})}\abs{\Omega}^{-\tau r}
\Norm{\{\one_{\Omega_{\vj}}\abs{f_{\vj}}^r \}_{\vj\in\Z^k}}{(L^{\frac{\vp}r}\ell^{\frac{\vq}r})_\pi}.
\end{split}
\end{equation*}
Here $\frac{\vp}r,\frac{\vq}r\in(1,\infty]^k$ and $\pi$ is an admissible permutation of $( \frac{\vp}r, \frac{\vq}r )$; hence we have access to the triangle inequality and the strong maximal inequality (Lemma \ref{FS strong}) in the following computation:
\begin{align*}
\Big\|\{\one_{\Omega^h_{\vj}}\abs{f_{\vj}}^r\}_{\vj\in\Z^k}
\Big\|_{(L^{\frac{\vp}r}\ell^{\frac{\vq}r})_\pi}
&\lesssim\bNorm{\{\one_{\Omega_{\vj}}\phi_{\vj}*\abs{g_{\vj}}^r \}_{\vj\in\Z^k}}{(L^{\frac{\vp}r}\ell^{\frac{\vq}r})_\pi}
\quad\text{by assumption \eqref{143}} \\
&\lesssim\Bigg\|\Bigg\{\sum_{h=0}^\infty 2^{h(n-M)} \mathcal M_{\vn}\Big(\one_{\Omega^h_{\vj}}\abs{g_{\vj}}^r\Big)\Bigg\}_{\vj\in\Z^k} \Bigg\|_{(L^{\frac{\vp}r}\ell^{\frac{\vq}r})_\pi}
\quad\text{by \eqref{insert Ohj 2}} \\
&\leq \sum_{h=0}^\infty 2^{h(n-M)}\bNorm{\Big\{\mathcal M_{\vn}\Big(\one_{\Omega^h_{\vj}}\abs{g_{\vj}}^r\Big)\Big\}_{\vj\in\Z^k} }{(L^{\frac{\vp}r}\ell^{\frac{\vq}r})_\pi}
\quad\text{by the triangle inequality} \\
&\lesssim \sum_{h=0}^\infty 2^{h(n-M)}
\Big\|\Big\{\one_{\Omega^h_{\vj}}\abs{g_{\vj}}^r\Big\}_{\vj\in\Z^k}
\Big\|_{(L^{\frac{\vp}r}\ell^{\frac{\vq}r})_\pi}\quad\text{by Lemma \ref{FS strong}}.
\end{align*}
Here
\begin{align*}
\Big\|\Big\{\one_{\Omega^h_{\vj}}\abs{g_{\vj}}^r\Big\}_{\vj\in\Z^k}
\Big\|_{(L^{\frac{\vp}r}\ell^{\frac{\vq}r})_\pi}^{\frac 1r}
&=\Norm{\{\one_{\Omega^h_{\vj}}g_{\vj}\}_{\vj\in\Z^k} }{(L^{\vp}\ell^{\vq})_\pi}
\leq\abs{\Omega^h}^\tau\Norm{\{g_{\vj}\}_{\vj\in\Z^k}}{\antau} \\
&\lesssim \Big[ 2^{hn}(1+h^{k-1})\abs{\Omega}\Big]^{\tau}
\Norm{\{g_{\vj}\}_{\vj\in\Z^k}}{\antau}\quad\text{by \eqref{Oh bd}}.
\end{align*}
Substituting back,
\begin{align*}
\Norm{\{\one_{\Omega_{\vj}}\abs{f_{\vj}}^r \}_{\vj\in\Z^k}}{(L^{\frac{\vp}r}\ell^{\frac{\vq}r})_\pi}
\lesssim \Bigg[ \sum_{h=0}^\infty 2^{h(n-M)} 2^{hn\tau}(1+h^{k-1})^{\tau}\Bigg]
\abs{\Omega}^{\tau r}\Norm{\{g_{\vj}\}_{\vj\in\Z^k}}{\antau}^r
\sim \abs{\Omega}^{\tau r}\Norm{\{g_{\vj}\}_{\vj\in\Z^k}}{\antau}^r
\end{align*}
for $M>n(1+\tau)$. Thus,
\begin{equation*}
\Norm{\{f_{\vj}\}_{\vj\in\Z^k}}{\antau}^r
= \sup_{\Omega\in\Open(\R^{\vn})}\abs{\Omega}^{-\tau r}\Norm{\{\one_{\Omega_{\vj}}\abs{f_{\vj}}^r \}_{\vj\in\Z^k}}{(L^{\frac{\vp}r}\ell^{\frac{\vq}r})_\pi}
\lesssim \Norm{\{g_{\vj}\}_{\vj\in\Z^k}}{\antau}^r,
\end{equation*}
which completes the proof of Lemma \ref{convo bd A bd}.
\end{proof}

The following lemma allows us to apply Lemma \ref{convo bd A bd} to relevant concrete cases.

\begin{lemma}\label{cases of convo bd}
For any $r,M\in(0,\infty)$ and $f\in\mathscr S_0'(\R^{\vn})$, there is a constant $C\in(0,\infty)$ such that the following holds: If either
\begin{equation}\label{case1 of convo bd}
f_{\vj}=\Ave{\varphi_{\vj}*f}{\vj}{\infty},\qquad g_{\vj}=\abs{\varphi_{\vj}*f},
\end{equation}
or $f_{\vj}=g_{\vj}=\abs{t_{\vj}}$ where $t_{\vj}$ is constant on each $R\in\D_{\vj}(\R^{\vn})$, then $f_{\vj}$ and $g_{\vj}$ satisfy the assumptions of Lemma \ref{convo bd A bd} with parameters $r,M,C$.
\end{lemma}

\begin{proof}
Case \eqref{case1 of convo bd} is just a restatement of Lemma \ref{cald4}.

For the other case, if $\phi$ as in Lemma \ref{convo bd A bd}, we observe that $\phi_{\vj}(x-y)\sim\abs{R}^{-1}$ if $x,y\in R\in\D_{\vj}(\R^{\vn})$.
If $f_{\vj}=g_{\vj}=\abs{t_{\vj}}$ and $t_{\vj}$ takes a constant value $t_R$ on each $R\in\D_{\vj}(\R^{\vn})$, it is then immediate that, for $x\in R\in\D_{\vj}(\R^{\vn})$,
\begin{align*}
\abs{f_{\vj}}^r(x)=\abs{t_R}^r
=\fint_R\abs{t_R}^r dy
\sim\int_R\phi_{\vj}(x-y)\abs{g_{\vj}(y)}^r dy
\leq\int_{\R^{\vn}}\phi_{\vj}(x-y)\abs{g_{\vj}(y)}^r dy.
\end{align*}
This finishes the proof of the second case
and hence Lemma \ref{cases of convo bd}.
\end{proof}

We can now prove the first equivalence of Theorem \ref{3 norms thm}, even under somewhat more general assumptions:

\begin{proposition}\label{sup vs ave}
Let $\vs\in\R^k$, $\tau\in[0,\infty)$, $\vp\in(0,\infty)^k$, $\vq\in(0,\infty]^k$,
and $\pi\in S_{[2k]}$ be an admissible permutation of $(\vp,\vq)$.
Let each $V_{\vj}$ take a constant value $V_R$ on each $R\in\D_{\vj}(\R^{\vn})$, where the sequence $\{V_R\}_{R\in\D(\R^{\vn})}$ is weakly doubling (Definition \ref{doubling}).
Let $f_{\vj}$ be either constant on each $R\in\D_{\vj}(\R^{\vn})$, or $f_{\vj}=\varphi_{\vj}*f$ for some Littlewood--Paley function $\varphi$ and $f\in\Sc_0'(\R^{\vn};\C^m)$.
Then
\begin{equation*}
\bNorm{ \big\{ \Ave{V_{\vj}f_{\vj}}{\vj}{\infty} \big\}_{\vj\in\Z^k} }{\aatau}
\sim\bNorm{ \big\{ V_{\vj}f_{\vj} \big\}_{\vj\in\Z^k} }{\aatau}
\end{equation*}
with the positive equivalence constants independent of
$\{f_{\vj}\}_{\vj\in\Z^k}$.
\end{proposition}

\begin{proof}
We observed ``$\geq$'' in \eqref{equi easy}. Turning to the main case ``$\lesssim$'', Lemma \ref{cases of convo bd} says that,
for every $x\in R\in\D_{\vj}(\R^{\vn})$,
\begin{equation*}
\sup_{y\in R}\abs{f_{\vj}(y)}^r\lesssim\phi_{\vj}*\abs{f_{\vj}}^r(x)
=\sum_{P\in\D_{\vj}(\R^{\vn})}\int_P\phi_{\vj}(x-z)\abs{f_{\vj}(z)}^r dz.
\end{equation*}
With $V_R f_{\vj}$ in place of $f_{\vj}$ (noting that it satisfies the same assumptions), expanding $V_R=V_R {V_P}^{-1}V_P$ on the right, this gives
\begin{equation}\label{sup vs ave1}
\sup_{y\in R}\abs{V_R f_{\vj}(y)}^r\lesssim
\sum_{P\in\D_{\vj}(\R^{\vn})}\int_P\phi_{\vj}(x-z)\abs{V_R V_P^{-1}}^r\abs{V_P f_{\vj}(z)}^r dz.
\end{equation}
By the assumption that $\{V_P\}_{P\in\D(\R^{\vn})}$ is weakly doubling (Definition \ref{doubling}), say of order $d$, it follows that
\begin{equation*}
\abs{V_R V_P^{-1}}\lesssim(1+\abs{2^{\vj}(c_R-c_P)})^d
\sim(1+\abs{2^{\vj}(x-z)})^d,
\end{equation*}
when $x\in R$ and $z\in P$. Thus, in the integrand in \eqref{sup vs ave1},
\begin{equation*}
\phi_{\vj}(x-z)\abs{V_R V_P^{-1}}^r
\lesssim\phi_{\vj}(x-z)\Big[1+\abs{2^{\vj}(x-z)}\Big]^{dr}
=\widetilde{\phi}_{\vj}(x-z),
\end{equation*}
where $\widetilde{\phi}$ is another function like $\phi$, but with decay $M-dr$ in place of $M$.
Noting also that $V_P=V_{\vj}(z)$ for $z\in P\in\D_{\vj}(\R^{\vn})$ and $V_R=V_{\vj}(y)$ for $y\in R\in\D_{\vj}(\R^{\vn})$, we deduce that
\begin{equation*}
\sup_{y\in R}\abs{V_{\vj}(y) f_{\vj}(y)}^r\lesssim\int_{\R^{\vn}}\widetilde{\phi}_{\vj}(x-z)\abs{V_{\vj}f_{\vj}(z)}^r dz,\quad\text{for all }x\in R\in\D_{\vj}(\R^{\vn}),
\end{equation*}
or in other words,
\begin{equation*}
\Ave{V_{\vj}f_{\vj}}{\vj}{\infty}^r
\lesssim \widetilde{\phi}_{\vj}*\abs{V_{\vj}f_{\vj}}^r.
\end{equation*}
This agrees with the assumption of Lemma \ref{convo bd A bd} with $\Ave{ V_{\vj} f_{\vj} }{\vj}{\infty}$ in place of $f_{\vj}$ and $\abs{V_{\vj}f_{\vj}}$ in place of $g_{\vj}$, as well as $\widetilde{\phi}_{\vj}$ in place of $\phi_{\vj}$. Since we can choose the decay $M$ of $\phi$ as large as we like to begin with, we can also make $M-dr$ as big as we want, and thus this $\widetilde{\phi}_{\vj}$ is just as good. Hence an application of Lemma \ref{convo bd A bd} completes the proof of Proposition \ref{sup vs ave}.
\end{proof}

A variation of the same argument also yields another comparison from Theorem \ref{3 norms thm}:

\begin{proposition}\label{sup vs V}
Let $\vs\in\R^k$, $\tau\in[0,\infty)$, $p\in(0,\infty)$, $\vp:=p\cdot\vec 1$, $\vq\in(0,\infty]^k$,
and $\pi\in S_{[2k]}$ be an admissible permutation of $(\vp,\vq)$.
Assume that $V\in\A_p(\R^{\vn})$
and $f_{\vj}$ is either constant on each $R\in\D_{\vj}(\R^{\vn})$, or $f_{\vj}=\varphi_{\vj}*f$ for some Littlewood--Paley function $\varphi$ and $f\in\Sc_0'(\R^{\vn};\C^m)$.
Then
\begin{equation*}
\bNorm{\big\{\Ave{\Red{V}{\vj}{p}f_{\vj}}{\vj}{\infty}\big\}_{\vj\in\Z^k} }{\aatau}
\lesssim
\bNorm{\big\{V f_{\vj}\big\}_{\vj\in\Z^k} }{\aatau}
\end{equation*}
with the implicit positive constant independent of
$\{f_{\vj}\}_{\vj\in\Z^k}$.
\end{proposition}

\begin{proof}
Repeating the beginning of the proof of Proposition \ref{sup vs ave}, but now with the specific coefficients $V_R=[V]_{\aveL^p(R)}$, and expanding further with $V_R=V_R V_P^{-1}V_P V(z)^{-1}V(z)$, we obtain, in place of \eqref{sup vs ave1},
\begin{equation}\label{sup vs ave2}
\begin{split}
\sup_{y\in R}\abs{V_R f_{\vj}(y)}^r
&\lesssim\sum_{P\in\D_{\vj}(\R^{\vn})}\int_P\phi_{\vj}(x-z)\abs{V_R V_P^{-1}}^r\abs{V_P V(z)^{-1}}^r\abs{V(z) f_{\vj}(z)}^r dz \\
&\lesssim\sum_{P\in\D_{\vj}(\R^{\vn})}\int_P\widetilde{\phi}_{\vj}(x-c_P)
\abs{V_P V(z)^{-1}}^r\abs{V(z) f_{\vj}(z)}^r dz,
\end{split}
\end{equation}
where $\widetilde{\phi}_{\vj}$ is another function like $\phi_{\vj}$, just like in the proof of Proposition \ref{sup vs ave}, and we used the observation that this function is almost constant on $P\in\D_{\vj}(\R^{\vn})$. If $p\in(0,1]$ and $z\in P$, then $V\in\A_p(\R^{\vn})$ guarantees that
\begin{equation*}
\abs{V_P V(z)^{-1}}
\leq\Norm{V_P V^{-1}}{\aveL^\infty(P)}\lesssim 1.
\end{equation*}
Thus, we directly deduce that
\begin{equation}\label{sup vs V1}
\Ave{\Red{V}{\vj}{p}f_{\vj}}{\vj}{\infty}^r
\lesssim \widetilde{\phi}_{\vj}*\abs{Vf_{\vj}}^r,\quad
\widetilde{\phi}(y):=(1+\abs{y})^{-(M-nr)},
\end{equation}
noting that $([V]_{\aveL^p(P)})_{P\in\D(\R^{\vn})}$ is weakly doubling of order $n$ by Lemma \ref{sharp coef 2}.

If $p\in(1,\infty)$, we use H\"older's inequality with exponent $p'/r$ and its conjugate, to the result that
\begin{align*}
&\int_P \abs{V_P V(z)^{-1}}^r \abs{V(z) f_{\vj}(z)}^r dz \\
&\quad\leq \bigg[\int_P \abs{V_P V(z)^{-1}}^{p'} dz \bigg]^{\frac{r}{p'}}
\bigg[\int_P \abs{V(z)f_{\vj}(z)}^{u} dz\bigg]^{1-\frac{r}{p'}},
\quad u:=\frac{rp'}{p'-r},
\end{align*}
where $V\in\A_p(\R^{\vn})$ guarantees that
\begin{equation*}
\bigg[\int_P \abs{V_P V(z)^{-1}}^{p'} dz \bigg]^{\frac{r}{p'}}
=\abs{P}^{\frac{r}{p'}}
\bigg[\fint_P \abs{V_P V(z)^{-1}}^{p'} dz\bigg]^{\frac{r}{p'}}
\lesssim\abs{P}^{\frac{r}{p'}}.
\end{equation*}
Another H\"older's inequality with the same exponents then allows us to continue the estimate \eqref{sup vs ave2} with
\begin{align*}
\sup_{y\in R}\abs{V_R f_{\vj}(y)}^r
&\lesssim \sum_{P\in\D_{\vj}(\R^{\vn})}\widetilde\phi_{\vj}(x-c_P)\abs{P}^{\frac{r}{p'}}
\bigg[ \int_P \abs{V(z) f_{\vj}(z)}^{u} dz \bigg]^{1-\frac{r}{p'}} \\
&\lesssim \Bigg[\sum_{P\in\D_{\vj}(\R^{\vn})}\widetilde\phi_{\vj}(x-c_P)\abs{P}\Bigg]^{ \frac{r}{p'} }
\Bigg[\sum_{P\in\D_{\vj}(\R^{\vn})}\widetilde\phi_{\vj}(x-c_P)\int_P \abs{V(z)f_{\vj}(z)}^{u} dz \Bigg]^{1-\frac{r}{p'}} \\
&\sim \bigg[\int_{\R^{\vn}}\widetilde\phi_{\vj}(x-z) dz\bigg]^{\frac{r}{p'}}
\bigg[\int_{\R^{\vn}} \widetilde\phi_{\vj}(x-z)\abs{V(z) f_{\vj}(z)}^{u} dz\bigg]^{1-\frac{r}{p'}} \\
&\lesssim \bigg[\int_{\R^{\vn}} \widetilde\phi_{\vj}(x-z)\abs{V(z) f_{\vj}(z)}^{u} dz\bigg]^{1-\frac{r}{p'}}.
\end{align*}
Raising to power $\frac{p'}{p'-r}$, we obtain
\begin{equation}\label{sup vs V2}
\Ave{\Red{V}{\vj}{p}f_{\vj}}{\vj}{\infty}^u
\lesssim \widetilde\phi_{\vj}*\abs{Vf_{\vj}}^u,\quad u=\frac{rp'}{p'-r},
\end{equation}
where $\widetilde\phi$ is as in \eqref{sup vs V1}. We proved \eqref{sup vs V2} under the assumption that $p\in(1,\infty)$, but we see that it also agrees with \eqref{sup vs V1} when $p\in(0,1]$, with the usual interpretations, so we have \eqref{sup vs V2} for all $p\in(0,\infty)$. As before, we can choose $r$ as small and $M$ as large as we like; hence we can also make $u$ as small and $M-nr$ as large as we like. Thus, the assumptions of Lemma \ref{convo bd A bd} are satisfied with $f_{\vj}=\Ave{\Red{V}{\vj}{p}f_{\vj}}{\vj}{\infty}$ and $g_{\vj}=\abs{Vf_{\vj}}$, and an application of that lemma shows that
\begin{equation*}
\bNorm{\big\{\Ave{\Red{V}{\vj}{p}f_{\vj}}{\vj}{\infty}\big\}_{\vj\in\Z^k}}{\aatau}
\lesssim\bNorm{\big\{Vf_{\vj}\big\}_{\vj\in\Z^k}}{\aatau}.
\end{equation*}
This is what we wanted to prove.
\end{proof}

\begin{remark}
For the case that $f_{\vj}$ is constant on each $R\in\D_{\vj}(\R^{\vn})$, a different (perhaps easier) proof could be given, along the lines of the corresponding argument in the one-parameter case, see \cite[Theorem 3.4]{FR21} and \cite[Theorem 3.28]{BHYY:1a}. Since we want to deal with both cases anyway, we prefer the proof above that covers both of them simultaneously. This is in contrast to \cite{BHYY:1a,FR21}, where separate arguments were used for function and sequence spaces.
\end{remark}

It remains to prove the last claim in Theorem \ref{3 norms thm}, i.e., the estimate of $\Norm{\{Vf_{\vj}\}_{\vj\in\Z^k}}{\aatau}$ by the other quasi-norms under the additional assumption \eqref{key cases}.

\begin{proposition}\label{V vs sup}
Let $\vs\in\R^k$, $\tau\in[0,\infty)$, $p\in(1,\infty)$, $\vp :=p\cdot\vone$, $\vq\in(0,\infty]^k$, and $\pi\in S_{[2k]}$ be admissible for $(\vp,\vq)$.
Suppose that $V\in \A_p(\R^{\vn})$ and that we are in one of the cases \eqref{key cases}. Then
\begin{equation*}
\bNorm{\{Vf_{\vj}\big\}_{\vj\in\Z^k}}{\aatau}
\lesssim\bNorm{\big\{\Red{V}{\vj}{p}f_{\vj}\big\}_{\vj\in\Z^k}}{\aatau}
\end{equation*}
with the implicit positive constant independent of
$\{f_{\vj}\}_{\vj\in\Z^k}$.
\end{proposition}

\begin{proof}
By definition,
\begin{equation*}
\BNorm{\Big\{ Vf_{\vj} \Big\}_{\vj\in\Z^k}}{\aatau}
=\sup_{\Omega\in\Open(\R^{\vn})}\abs{\Omega}^{-\tau}\BNorm{\Big\{ \one_{\Omega_{\vj}}2^{\vj\cdot\vs}Vf_{\vj} \Big\}_{\vj\in\Z^k}}{(L^{\vp}\ell^{\vq})_\pi }.
\end{equation*}
Here
\begin{align*}
\BNorm{\Big\{ \one_{\Omega_{\vj}}2^{\vj\cdot\vs}Vf_{\vj} \Big\}_{\vj\in\Z^k}}{(L^{\vp}\ell^{\vq})_\pi }
=\BNorm{\Big\{ V\Red{V}{\vj}{p}^{-1} \big(\one_{\Omega_{\vj}}2^{\vj\cdot\vs}\Red{V}{\vj}{p}f_{\vj}\big) \Big\}_{\vj\in\Z^k}}{(L^{\vp}\ell^{\vq})_\pi }
\leq\bNorm{\big\{ \gamma_{\vj} g_{\vj} \big\}_{\vj\in\Z^k}}{(L^{\vp}\ell^{\vq})_\pi },
\end{align*}
where
\begin{equation*}
\gamma_{\vj}:= \abs{V\Red{V}{\vj}{p}^{-1}},\quad g_{\vj}:=\one_{\Omega_{\vj}}2^{\vj\cdot\vs}\Ave{\Red{V}{\vj}{p}(\varphi_{\vj}*f)\big)}{\vj}{\infty}.
\end{equation*}
These satisfy the relevant assumptions of Corollary \ref{Carl Ap} with $g_{\vj}$ in place of $f_{\vj}$ of the corollary: $\gamma_{\vj}$ is exactly as in the corollary, and $g_{\vj}$ is constant on $R\in\mathscr D_{\vj}$. Hence the said corollary applies to show that
\begin{equation*}
\begin{split}
\bNorm{\big\{ \gamma_{\vj} g_{\vj} \big\}_{\vj\in\Z^k}}{(L^{\vp}\ell^{\vq})_\pi }
&\lesssim\bNorm{\big\{ g_{\vj} \big\}_{\vj\in\Z^k}}{(L^{\vp}\ell^{\vq})_\pi }
=\BNorm{\Big\{ \one_{\Omega_{\vj}}2^{\vj\cdot\vs}\Ave{\Red{V}{\vj}{p}f_{\vj}}{\vj}{\infty} \Big\}_{\vj\in\Z^k}}{(L^{\vp}\ell^{\vq})_\pi }\\
&\leq\abs{\Omega}^{\tau}\BNorm{\Big\{ \Ave{\Red{V}{\vj}{p}f_{\vj}}{\vj}{\infty} \Big\}_{\vj\in\Z^k}}{\aatau}.
\end{split}
\end{equation*}
Substituting back, we obtain
\begin{equation*}
\BNorm{\Big\{ Vf_{\vj} \Big\}_{\vj\in\Z^k}}{\aatau}
\lesssim\BNorm{\Big\{ \Ave{\Red{V}{\vj}{p}f_{\vj}}{\vj}{\infty} \Big\}_{\vj\in\Z^k}}{\aatau}.
\end{equation*}
Finally, for $V\in\A_p(\R^{\vn})$ as in the assumptions, the related
sequence $\{[V]_{\aveL^p(R)}\}_{R\in\D(\R^{\vn})}$ of reducing operators
is weakly doubling by Remark \ref{Ap doubling}, and hence
\begin{equation*}
\BNorm{\Big\{ \Ave{\Red{V}{\vj}{p}f_{\vj}}{\vj}{\infty} \Big\}_{\vj\in\Z^k}}{\aatau}
\lesssim  \BNorm{\Big\{ \Red{V}{\vj}{p}f_{\vj}\Big\}_{\vj\in\Z^k}}{\aatau}
\end{equation*}
by Proposition \ref{sup vs ave}.
This finishes the proof Proposition \ref{V vs sup}.
\end{proof}

The proof of Theorem \ref{3 norms thm} now consists of collecting the results of the various propositions:

\begin{proof}[Proof of Theorem \ref{3 norms thm}]
  Claim \eqref{it:f-avef} is Proposition \ref{sup vs ave} with $V_{\vj}\equiv1$.

  Claim \eqref{it:[V]<V} follows by combining Propositions \ref{sup vs ave} and \ref{sup vs V}, noting that, by Remark \ref{Ap doubling}, the assumptions of Proposition \ref{sup vs ave} are implied by those of Theorem \ref{3 norms thm}.

  Claim \eqref{it:[V]=V} is Proposition \ref{V vs sup}.
\end{proof}

\subsection{The $\varphi$-transform characterization}
\label{phi-transform}

In Theorem \ref{3 norms thm}, we proved the equivalence of several different norms of sequences of functions $f_{\vj}$, either arising from the Littlewood--Paley components $f_{\vj}=\varphi_{\vj}*f$ of a distribution, or from discrete coefficient sequences indexed by the dyadic rectangles; however, we always compared norms of two object of the same type. The aim of this subsection is to complement these results by comparing the norm of a distribution with the norm of a related sequence of coefficients, thereby obtaining equivalences between objects of two seemingly different types. Such equivalences play a key role in the subsequent sections, where operators on function or distribution spaces are studied via their discrete counterparts on sequence spaces. The identification of the function and sequences spaces is obtained via the so-called $\varphi$-transform, first introduced by Frazier and Jawerth \cite{FJ90}.

Our main result concerning its mapping properties of the matrix-weighted spaces is Theorem \ref{phi} below. Before getting there, we need some preparations, and we begin with the following elementary bound for individual entries $t_R$ of $t\in\atau(\V)$:

\begin{lemma}\label{one coef}
Let $\vs\in\R^k$, $\tau\in[0,\infty)$, $\vp\in(0,\infty)^k$, and $\vq\in(0,\infty]^k$.
Let $\V=\{V_{\vj}\}_{\vj\in\Z^k}$ be a sequence of $\C^{m\times m}$-valued functions,
where each $V_{\vj}$ takes the constant value $V_R$ on $R\in\D_{\vj}(\R^{\vn})$. For $t\in(\C^m)^{\D(\R^{\vn})}$, we have
\begin{equation*}
\abs{V_R t_R}
\leq 2^{-\vj\cdot\vn(\tau+\frac12)}
2^{-\vj\cdot(\vs-\vn/\vp)}
\Norm{t}{\atau(\V)},
\end{equation*}
where $  \vn/\vp:=(n_1/p_1,\ldots,n_k/p_k).$
\end{lemma}

\begin{proof}
A lower bound for the supremum defining $\Norm{t}{\atau(\V)}$ is obtained by taking $\Omega=R\in\D_{\vj}(\R^{\vn})$ (Although $R$ is not an open set,
a simple modification shows that this choice is still valid.) and dropping all functions except the $\vj$-th for this fixed $\vj\in\Z^k$. Thus,
\begin{equation*}
\begin{split}
\Norm{t}{\atau(\V)}
&\geq\abs{R}^{-\tau}\bNorm{\one_R 2^{\vj\cdot\vs} |V_{\vj}t_{\vj}|}{L^{\vp}(\R^{\vn})}
=\abs{R}^{-\tau}\bNorm{\one_R 2^{\vj\cdot\vs} \abs{R}^{-\frac12} |V_Rt_R|}{L^{\vp}(\R^{\vn})} \\
&=\abs{R}^{-\tau-\frac12} 2^{\vj\cdot\vs} |V_Rt_R|\,
\Norm{\one_R }{L^{\vp}(\R^{\vn})}
=2^{\vj\cdot\vn(\tau+\frac12)} 2^{\vj\cdot\vs}\abs{V_R t_R}
2^{-\vj\cdot\vn/\vp}.
\end{split}
\end{equation*}
This finishes the proof of Lemma \ref{one coef}.
\end{proof}

\begin{lemma}\label{Tpsi welldef}
Let $\vs\in\R^k$, $\tau\in[0,\infty)$, $\vp\in(0,\infty)^k$, and $\vq\in(0,\infty]^k$.
Let $\V=\{V_{\vj}\}_{\vj\in\Z^k}$ be a sequence of $\C^{m\times m}$-valued functions, where each $V_{\vj}$ takes a constant value $V_R$ on each $R\in\D_{\vj}(\R^{\vn})$, and the sequence $\{V_R\}_{R\in\D(\R^{\vn})}$ is strongly doubling (Definition \ref{doubling}). If $t\in\atau(\V)$ and $\psi\in\Sc_0(\R^{\vn})$, then
\begin{equation*}
T_{\psi}t:=\sum_{R\in\D(\R^{\vn})}t_R\psi_R
\end{equation*}
converges in $\Sc_0'(\R^{\vn};\C^m)$. More precisely, for each $\varphi\in\Sc_0(\R^{\vn})$,
\begin{equation*}
\sum_{R\in\D(\R^{\vn})}\abs{t_R\pair{\psi_R}{\varphi}}
\lesssim \Norm{t}{\atau(\V)},
\end{equation*}
where the implicit positive constant may depend on the space $\atau(\V)$ and $\varphi$ but is independent of $t$.
In particular, these estimates hold
if $\V\equiv1$ or $\vp=p\cdot\vec1$ for some $p\in(0,\infty)$, $V\in\A_p(\R^{\vn})$,
and $\V=[V]_p$ (i.e., $V_R=[V]_{\aveL^p(R)}$).
\end{lemma}

\begin{proof}
Once the other claims are proved, the last claim is immediate from the strong doubling property of the reducing operators of $\A_p(\R^{\vn})$ weights (Remark \ref{Ap doubling}).

Let $P:=[0,1)^n$, so that $\varphi=\varphi_P$. Inserting the identity as $V_PV_P^{-1}V_RV_R^{-1}$, we can estimate
\begin{equation*}
\begin{split}
\sum_{R\in\D(\R^{\vn})}\abs{t_R\pair{\psi_R}{\varphi}}
&\leq\abs{V_P^{-1}}\sum_{R\in\D(\R^{\vn})}
\abs{V_P V_R^{-1}} \, \abs{V_R t_R\pair{\psi_R}{\varphi}} \\
&\sim\sum_{\vj\in\Z^k}\sum_{R\in\D_{\vj}(\R^{\vn})}
\abs{V_PV_R^{-1}} \, \abs{V_R t_R\pair{\psi_R}{\varphi_P}}
=: \mathrm{I},
\end{split}
\end{equation*}
since $\abs{V_P^{-1}}=\abs{V_{[0,1)^n}^{-1}}$ is just a positive number depending on the quantities admissible in the implicit positive constants.

Since $\{V_R\}_{R\in\D(\R^{\vn})}$ is strongly doubling, say of order $(\va,\vb,\vc)\in[0,\infty)^{3k}$, Remark \ref{str doubling} shows that
\begin{equation*}
  \abs{V_P V_R^{-1}}\lesssim 2^{\abs{\vj}A} \Omega^{\vc} (2^{\vj\wedge\vnull}a_R),
  \qquad\Omega^{\vc}(x):=\prod_{i=1}^k(1+\abs{x_i})^{c_i},
\end{equation*}
where $A:=\max(\abs{\va},\abs{\vb})$.

On the other hand, Lemma \ref{one coef} guarantees that
\begin{equation*}
  \abs{V_R t_R}
  \lesssim 2^{-\vj\cdot\vn(\tau+\frac12)}
  2^{-\vj\cdot(\vs-\vn/\vp)}
  \Norm{t}{\atau(\V)}.
\end{equation*}
Finally, Lemma \ref{59x} says that
\begin{equation*}
   \abs{\pair{\psi_R}{\varphi_P}}\lesssim 2^{-\vj\cdot\vn/2}2^{-\abs{\vj}M}\omega_{\vj\wedge\vnull}^{\vN}(a_R),
\end{equation*}
where $\omega^{\vN}=\Omega^{-\vN}$. Hence
\begin{equation*}
\begin{split}
  |V_PV_R^{-1}| \, |V_R t_R\pair{\psi_R}{\varphi_P}|
  &\lesssim 2^{\abs{\vj}(A+B-M)}\omega_{\vj\wedge\vnull}^{\vN-\vc}(a_R)\quad
  \text{with}\quad B:=n |\tau+1|+\abs{\vs-\vn/\vp} \\
  &=2^{-\abs{\vj}M'}\omega_{\vj\wedge\vnull}^{\vN'}(a_R)\quad
  \text{with}\quad \begin{cases} M':=M-A-B\\ \vN':=\vN-\vc.\end{cases}
\end{split}
\end{equation*}
Since $M$ and the components of $\vN$ can be chosen as large as we like, the same is true of $M'$ and $\vN'$. Now
\begin{equation*}
\sum_{R\in\D_{\vj}(\R^{\vn})}\omega_{\vj\wedge\vnull}^{\vN'}(a_R)
=\sum_{k\in\Z^{\vn}}\omega_{\vj\wedge\vnull}^{\vN'}(2^{-\vj}\vk)
=\prod_{h=1}^k\sum_{k_h\in\Z^{n_h}}\omega_{j_h\wedge 0}^{N_h'}(2^{-j_h}k_h).
\end{equation*}
Here, for $N_h'>n_h$,
\begin{equation*}
\begin{split}
\sum_{k_h\in\Z^{n_h}}\omega_{j_h\wedge 0}^{N_h'}(2^{-j_h}k_h)
&=2^{(j_h\wedge 0) n_h}\sum_{k_h\in\Z^{n_h}}(1+2^{-(j_h)_+}\abs{k_h})^{-N_h'} \\
&\sim 2^{(j_h\wedge 0) n_h} 2^{(j_h)_+ n_h}=2^{j_h n_h}.
\end{split}
\end{equation*}
Thus,
\begin{equation*}
\begin{split}
\mathrm{I}
\lesssim\sum_{\vj\in\Z^k}2^{-\abs{\vj}M'}\sum_{R\in\D_{\vj}(\R^{\vn})}\omega_{\vj\wedge\vnull}^{\vN'}(a_R)
\lesssim\sum_{\vj\in\Z^k}2^{-\abs{\vj}M'}2^{\vj\cdot\vn}\lesssim1,
\end{split}
\end{equation*}
if $M'>\abs{\vn}$.
This finishes the proof of Lemma \ref{Tpsi welldef}.
\end{proof}

\begin{lemma}\label{shift}
Let $\vs\in\R^k$, $\tau\in[0,\infty)$, $\vp\in(0,\infty)^k$, $\vq\in(0,\infty]^k$,
and $\pi\in S_{[2k]}$ be a permutation of $(\vp,\vq)$.
Then, for every sequence of $\C^m$-valued functions $\{f_{\vi}\}_{\vi\in\Z^k}$
and every $\vh\in\Z^k$ with $\abs{\vh}\lesssim 1$,
\begin{equation*}
\Norm{\{f_{\vi+\vh}\}_{\vi\in\Z^k}}{\aatau(\R^{\vn})}
\sim\Norm{\{f_{\vi}\}_{\vi\in\Z^k}}{\aatau(\R^{\vn})},
\end{equation*}
where the positive equivalence constants are independent of $\{f_{\vi}\}_{\vi\in\Z^k}$.
\end{lemma}

\begin{proof}
By definition,
\begin{equation*}
\Norm{\{f_{\vi+\vh}\}_{\vi\in\Z^k}}{\aatau}
=\sup_{\Omega\in\Open(\R^{\vn})}\abs{\Omega}^{-\tau}\bNorm{\big\{\one_{\Omega_{\vi}}2^{\vi\cdot\vs}
f_{\vi+\vh}\big\}_{\vi\in\Z^k}}{(L^{\vp}\ell^{\vq})_\pi}.
\end{equation*}
The quasi-norm of $(L^{\vp}\ell^{\vq})_\pi$ is obviously invariant under translations. Hence, making the change of variable $\vj:=\vi+\vh$, we obtain
\begin{equation}\label{shift1}
\begin{split}
\bNorm{\big\{\one_{\Omega_{\vi}}2^{\vi\cdot\vs}f_{\vi+\vh}\big\}_{\vi\in\Z^k}}{(L^{\vp}\ell^{\vq})_\pi}
&=\bNorm{\big\{\one_{\Omega_{\vj-\vh}}2^{(\vj-\vh)\cdot\vs}f_{\vj}\big\}_{\vj\in\Z^k}}
{(L^{\vp}\ell^{\vq})_\pi} \\
&=2^{-\vh\cdot\vs}\bNorm{\big\{\one_{\Omega_{\vj-\vh}}2^{\vj\cdot\vs}f_{\vj}
\big\}_{\vj\in\Z^k}}{(L^{\vp}\ell^{\vq})_\pi},
\end{split}
\end{equation}
where $2^{-\vh\cdot\vs}\sim 1$ for $\abs{\vh}\lesssim 1$ and a fixed $\vs\in\R^k$.

Next, we aim at showing that
\begin{equation}\label{shift2}
\Omega_{\vj-\vh}\subset(\Omega^{\vh})_{\vj}
\end{equation}
for a suitable set $\Omega^{\vh}$. Recall that
\begin{equation*}
\Omega_{\vj-\vh}:=\bigcup_{R\in\D_{\vj-\vh}(\Omega)}R,\qquad
\D_{\vj-\vh}(\Omega):=\Big\{R\in\D_{\vj-\vh}(\R^{\vn}):\ R\subset \Omega\Big\}.
\end{equation*}
Consider one of these rectangles $R=R_1\times\cdots\times R_k\in\D_{\vj-\vh}(\Omega)$. We are going to cover $R$ by a disjoint union of rectangles $P=P_1\times\cdots\times P_k\in\D_{\vj}(\R^{\vn})$ as follows: If $h_l\geq 0$, then $R_l\in\D_{j_l-h_l}(\R^{\vn})$ is partitioned by cubes of $\D_{j_l}(\R^{\vn})$, and we let $P_l$ run through each of them. If $h_l<0$, then $R_l\in\D_{j_l-h_l}(\R^{\vn})$ is contained in a unique $P_l\in\D_{j_l}(\R^{\vn})$, and this is our choice of $P_l$. It is clear that $R$ is contained in the union of all $P$ formed like this. On the other hand, if $P$ is one of these rectangles, then
\begin{equation*}
\frac{\abs{P\cap R}}{\abs{P}}
=\prod_{l=1}^k\frac{\abs{P_l\cap R_l}}{\abs{P_l}}
=\prod_{\genfrac{}{}{0pt}{}{l=1}{h_l<0}}^k\frac{\abs{R_l}}{\abs{P_l}}
=\prod_{\genfrac{}{}{0pt}{}{l=1}{h_l<0}}^k 2^{h_l n_l}
=2^{-\vh_-\cdot\vn}
>2^{-\vh_-\cdot\vn-1}.
\end{equation*}
Since $R\subset\Omega$, this shows that
\begin{equation*}
P\subset \Big\{x\in\R^{\vn}:\ (\mathcal M_{\vn}\one_{\Omega})(x)
> 2^{-\vh_-\cdot\vn-1}\Big\}=:\Omega^{\vh},
\end{equation*}
where $\Omega^{\vh}$ is open. Since $P\in\D_{\vj}(\R^{\vn})$, it also follows that $P\subset (\Omega^{\vh})_{\vj}$. Since $R$ is contained in the union of such $P$, we have $R\subset (\Omega^{\vh})_{\vj}$, and since this holds for every $R\in\D_{\vj-\vh}(\Omega)$, we obtain \eqref{shift2}, as desired. Moreover, the strong maximal inequality guarantees that
$\abs{\Omega^{\vh}}\lesssim\abs{\Omega}$
for $\abs{\vh}\lesssim 1$.

We can now return to the computation \eqref{shift1}:
\begin{equation*}
\begin{split}
\bNorm{\big\{\one_{\Omega_{\vj-\vh}}2^{\vj\cdot\vs}f_{\vj}\big\}_{\vj\in\Z^k}}{(L^{\vp}\ell^{\vq})_\pi}
&\leq\bNorm{\big\{\one_{(\Omega^{\vh})_{\vj}}2^{\vj\cdot\vs}f_{\vj}\big\}_{\vj\in\Z^k}}
{(L^{\vp}\ell^{\vq})_\pi} \\
&\leq\abs{\Omega^{\vh}}^{\tau}\bNorm{\big\{f_{\vj}\big\}_{\vj\in\Z^k}}{\aatau}
\lesssim\abs{\Omega}^{\tau}\bNorm{\big\{f_{\vj}\big\}_{\vj\in\Z^k}}{\aatau}.
\end{split}
\end{equation*}
Substituting this into \eqref{shift1}, we obtain
\begin{equation*}
\Norm{\{f_{\vi+\vh}\}_{\vi\in\Z^k}}{\aatau}
\lesssim\Norm{\{f_{\vi}\}_{\vi\in\Z^k}}{\aatau}.
\end{equation*}
If we apply this to $g_{\vi}:=f_{\vi+\vh}$ in place of $f_{\vi}$ and the shift $-\vh$ in place of $+\vh$, we obtain the reverse inequality.
This finishes the proof of Lemma \ref{shift}.
\end{proof}

The following theorem is the $\varphi$-transform characterisation of matrix-weighted multi-parameter Triebel--Lizorkin-type spaces, the main result of this subsection, which extends various related results in the earlier literature. In the unweighted case with Morrey parameter $\tau=0$, the result is due to Frazier and Jawerth \cite[Theorem 2.2]{FJ90} in the one-parameter case, and due to Georgiadis et al. \cite[Theorem 4.2]{GKP21} in the bi-parameter ($k=2$) case for pure Besov or Triebel--Lizorkin spaces.
In the matrix-weighted one-parameter case, it was proved by Frazier and Roudenko \cite[Theorem 1.1]{FR21} with $\tau=0$ and by Bu et al. \cite[Theorem 3.29]{BHYY:1a} with general Morrey parameter $\tau\in[0,\infty)$.

\begin{theorem}\label{phi}
Let $s\in\R^k$, $\tau\in[0,\infty)$, $\vp\in(0,\infty)^k$, $\vq\in(0,\infty]^k$, and $\pi\in S_{[2k]}$ be admissible for $(\vp,\vq)$.
Let $\varphi,\psi,\varrho\in\Sc (\R^{\vn})$ be Littlewood--Paley functions, and $\widetilde\varphi(x):=\overline{\varphi(-x)}$ for every $x\in\mathbb{R}^n$.
Let $\V=\{V_{\vj}\}_{\vj\in\Z^k}$, where each $V_{\vj}$ takes a constant value $V_R\in\C^{m\times m}$ on each $R\in\D_{\vj}(\R^{\vn})$, and let $\{V_R\}_{R\in\D(\R^{\vn})}$ be strongly doubling (Definition \ref{doubling}).
Then the operators
\begin{equation*}
S_\varphi:f\mapsto\big\{\pair{f}{\varphi_R}\big\}_{R\in\D(\R^{\vn})},\quad
\Atau(\V)^{\widetilde{\varphi}}\to\atau(\V)
\end{equation*}
and
\begin{equation*}
T_\psi:t\mapsto\sum_{R\in\D(\R^{\vn})}t_R\psi_R,\quad
\atau(\V)\to \Atau(\V)^{\varrho}
\end{equation*}
are bounded.
If $(\varphi,\psi)$ is a Littlewood--Paley pair, then
$T_\psi\circ S_\varphi$ is the identity on $\Atau(\V)^{\widetilde{\varphi}}$.
In particular, these results hold
if $\V\equiv1$ or $\vp=p\cdot\vec1$ for some $p\in(0,\infty)$, $V\in\A_p(\R^{\vn})$,
and $\V=[V]_p$ (i.e., $V_R=[V]_{\aveL^p(R)}$).
\end{theorem}

\begin{proof}
Once the other claims are proved, the last claim is immediate from the strong doubling property of the reducing operators of $\A_p$ weights (Remark \ref{Ap doubling}).

\textit{The boundedness of $S_\varphi$}: For $R\in\D_{\vj}(\R^{\vn})$,
\begin{equation*}
t_R:=\pair{f}{\varphi_R}
=\abs{R}^{\frac12}\pair{f}{\varphi_{\vj}(\cdot-a_R)}
=\abs{R}^{\frac12}\int_{\R^{\vn}}f(y)\widetilde\varphi_{\vj}(a_R-y)dy
=\abs{R}^{\frac12}\widetilde\varphi_{\vj}*f(a_R).
\end{equation*}
Hence
\begin{equation*}
\abs{R}^{-\frac12}\abs{V_R t_R}=\abs{V_R\widetilde\varphi_{\vj}*f(a_R)}
\leq\sup_{y\in R}\abs{V_R\widetilde\varphi_{\vj}*f(y)}
\end{equation*}
and thus, simply rewriting the previous line in a different notation,
\begin{equation*}
\abs{V_{\vj}t_{\vj}}\leq\Ave{V_{\vj}(\widetilde\varphi_{\vj}*f)}{\vj}{\infty}.
\end{equation*}
Then, using definitions, the previous pointwise bound, as well as Proposition \ref{sup vs ave}, we obtain
\begin{align*}
\Norm{t}{\atau(\V)}
&=\Norm{\{V_{\vj}t_{\vj}\}_{\vj\in\Z^k}}{\aatau}
\leq\Norm{\{\Ave{V_{\vj}(\widetilde\varphi_{\vj}*f)}{\vj}{\infty}\}_{\vj\in\Z^k}}{\aatau}  \\
&\sim\Norm{\{V_{\vj}(\widetilde\varphi_{\vj}*f)\}_{\vj\in\Z^k}}{\aatau}
=\Norm{f}{\Atau(\V)}^{\widetilde\varphi}.
\end{align*}

\textit{The boundedness of $T_\psi$}: We already verified in Lemma \ref{Tpsi welldef} that $T_\psi t$ is well defined as a distribution in $\Sc_0'(\R^{\vn};\C^m)$. Hence its $\Atau(\V)^{\varrho}$ quasi-norm is a well defined number in $[0,\infty]$, which we need to estimate. Now
\begin{align*}
\varrho_{\vi}*T_\psi t(x)
&=\sum_{R\in\D(\R^{\vn})}t_R\varrho_{\vi}*\psi_R(x)
=\sum_{\vj\in\Z^k}\sum_{R\in\D_{\vj}(\R^{\vn})} \abs{R}^{\frac 12}t_R\varrho_{\vi}*\psi_{\vj}(x-a_R) \\
&=\sum_{\Norm{\vh}{\infty}\leq 1}\sum_{R\in\D_{\vi+\vh}(\R^{\vn})} \abs{R}^{\frac 12}t_R\varrho_{\vi}*\psi_{\vi+\vh}(x-a_R),
\end{align*}
since the Fourier supports of Littlewood--Paley functions imply that $\varrho_{\vi}*\psi_{\vj}=0$ unless $\Norm{\vi-\vj}{\infty}\leq 1$.

For $x\in P\in\D_{\vi}(\R^{\vn})$, we then compute
\begin{align*}
\abs{V_P\varrho_{\vi}*T_\psi t(x)}
&\leq\sum_{\Norm{\vh}{\infty}\leq 1}\sum_{R\in\D_{\vi+\vh}(\R^{\vn})} \abs{R}^{\frac 12} \abs{V_Pt_R\varrho_{\vi}*\psi_{\vi+\vh}(x-a_R)} \\
&\leq\sum_{\Norm{\vh}{\infty}\leq 1} \sum_{R\in\D_{\vi+\vh}(\R^{\vn})} \abs{R}^{\frac 12} |V_PV_R^{-1}|\,
\abs{V_R t_R} \, |\varrho_{\vi}*\psi_{\vi+\vh}(x-a_R)|.
\end{align*}
The last factor is dominated by Lemma \ref{59} and
the factor $\abs{V_P V_R^{-1}}$ by the assumption of
strong doubling and Remark \ref{str doubling}.
Observing in both cases that there is no essential difference between $2^{\vi}$ and $2^{\vi+\vh}$, we obtain
\begin{align*}
\abs{V_P\varrho_{\vi}*T_\psi t(x)}
&\lesssim \sum_{\Norm{\vh}{\infty}\leq 1} \sum_{R\in\D_{\vi+\vh}(\R^{\vn})}
\abs{R}^{\frac 12}
\Omega(2^{\vi}(a_P-a_R)) \abs{V_R t_R} \phi_{\vi+\vh}(x-a_R) \\
&\lesssim\sum_{\Norm{\vh}{\infty}\leq 1}\sum_{R\in\D_{\vi+\vh}(\R^{\vn})}\abs{R}^{\frac 12}
\abs{V_R t_R}\phi_{\vi+\vh}'(a_P-a_R),
\end{align*}
where $\phi'$ is another function like $\phi$, perhaps with a slower decay due to the absorption of the moderately growing factor $\Omega$, but since we can choose the decay of $\phi$ as fast as we like, the same is true for $\phi'$. We also noted that $\phi_{\vi+\vh}$ is essentially constant on $P\in\D_{\vi}(\R^{\vn})$ to replace $x\in P$ by $a_P$.

Noting that $\abs{R}=2^{-(\vi+\vh)\cdot\vn}$ is the inverse of the scaling factor in front of $\phi_{\vi+\vh}'$, we can further write
\begin{align*}
\abs{V_P\varrho_{\vi}*T_\psi t(x)}
&=\sum_{\Norm{\vh}{\infty}\leq 1}\sum_{R\in\D_{\vi+\vh}(\R^{\vn})}
\abs{R}^{-\frac 12} \abs{V_R t_R} \phi'(2^{\vi+\vh}(a_P-a_R)) \\
&=\sum_{\Norm{\vh}{\infty}\leq 1}\sum_{R\in\D_{\vi+\vh}(\R^{\vn})}
\Big| (V_{\vi+\vh}t_{\vi+\vh})(a_R) \Big|
\phi'(2^{\vi+\vh}(a_P-a_R)).
\end{align*}

For $r<\min\{1,p,q_i:i\in[k]\}$ we obtain, denoting by $\phi'':=(\phi')^r$ yet another function like $\phi$ and $\phi'$,
\begin{align*}
\abs{V_P\varrho_{\vi}*T_\psi t(x)}^r
&\lesssim\sum_{\Norm{\vh}{\infty}\leq 1}\sum_{R\in\D_{\vi+\vh}(\R^{\vn})}
\big|(V_{\vi+\vh}t_{\vi+\vh})(a_R)\big|^r \phi''(2^{\vi+\vh}(a_P-a_R)) \\
&\lesssim\sum_{\Norm{\vh}{\infty}\leq 1}\sum_{R\in\D_{\vi+\vh}(\R^{\vn})}\abs{R}
\big|(V_{\vi+\vh}t_{\vi+\vh})(a_R)\big|^r \phi''_{\vi+\vh}(a_P-a_R) \\
&\sim\sum_{\Norm{\vh}{\infty}\leq 1}(\phi_{\vi+\vh}''*\abs{V_{\vi+\vh}t_{\vi+\vh}}^r)(x),
\end{align*}
where in the last step we used the facts that $V_{\vi+\vh}t_{\vi+\vh}$ is constant on each $R\in\D_{\vi+\vh}(\R^{\vn})$ and that $\phi_{\vi+\vh}''$ is essentially constant at the same length scales.

Recalling that $x\in P\in\D_{\vi}(\R^{\vn})$ were arbitrary above, we have shown the pointwise estimate
\begin{equation}\label{Tpsi x}
\abs{V_{\vi}(\varrho_{\vi}*T_\psi t)}^r
\lesssim\sum_{\Norm{\vh}{\infty}\leq 1}\phi_{\vi+\vh}''*\abs{V_{\vi+\vh}t_{\vi+\vh}}^r.
\end{equation}
Hence
\begin{align*}
\Norm{T_\psi t}{\Atau(\V)}^{\varrho}
&=\Norm{\{V_{\vi}(\varrho_{\vi}*T_\psi t)\}_{\vi\in\Z^k}}{\aatau} \\
&\lesssim \sum_{\Norm{\vh}{\infty}\leq 1} \bNorm{\big\{(\phi_{\vi+\vh}''*\abs{V_{\vi+\vh}t_{\vi+\vh}}^r)^{1/r}\big\}_{\vi\in\Z^k}}{\aatau}
\quad\text{by \eqref{Tpsi x}} \\
&\lesssim \sum_{\Norm{\vh}{\infty}\leq 1} \bNorm{\big\{V_{\vi+\vh}t_{\vi+\vh}\big\}_{\vi\in\Z^k}}{\aatau}
\quad\text{by Lemma \ref{convo bd A bd}} \\
&\lesssim \bNorm{\big\{V_{\vi}t_{\vi}\big\}_{\vi\in\Z^k}}{\aatau}=\Norm{t}{\atau(\V)}\quad\text{by Lemma \ref{shift}.}
\end{align*}

\textit{The identity of $T_\psi\circ S_\varphi$}: According to Lemma \ref{cald4}, we have
\begin{equation*}
f=\sum_{R\in\D(\mathbb R^{\vec n})}\pair{f}{\varphi_R}\psi_R
\end{equation*}
in $\Sc_0'(\R^{\vn};\C^m)$, for every $f$ in this space. When $f\in\Atau(\V)^{\tilde\varphi}\subset\Sc_0'(\R^{\vn};\C^m)$, the right-hand side coincides with $T_\psi\circ S_\varphi(f)\in\Atau(\V)^{\varrho}\subset\Sc_0'(\R^{\vn};\C^m)$. Here we can in particular choose $\varrho=\widetilde\varphi$.
This finishes the proof of Theorem \ref{phi}.
\end{proof}

\begin{corollary}\label{38}
Let $s\in\R^k$, $\tau\in[0,\infty)$, $\vp\in(0,\infty)^k$, $\vq\in(0,\infty]^k$, and $\pi\in S_{[2k]}$ be admissible for $(\vp,\vq)$.
Let $\V=\{V_{\vj}\}_{\vj\in\Z^k}$, where each $V_{\vj}$ takes a constant value $V_R\in\C^{m\times m}$ on each $R\in\D_{\vj}(\R^{\vn})$, and let $\{V_R\}_{R\in\D(\R^{\vn})}$ be strongly doubling (Definition \ref{doubling}).
If $\varphi$ and $\varrho$ are two Littlewood--Paley functions,
then the quasi-norms $\Norm{\cdot }{\Atau(\V)}^{\varphi}$
and $\Norm{\cdot}{\Atau(\V)}^{\varrho}$ are equivalent. Thus, the space
\begin{equation*}
\Atau(\V):=\Big\{f\in\Sc_0'(\R^{\vn};\C^m):\
\Norm{f}{\Atau(\V)}^{\varphi}<\infty\Big\}
\end{equation*}
is independent of the choice of the Littlewood--Paley function $\varphi$.
\end{corollary}

\begin{proof}
Given a Littlewood--Paley function $\varphi$, we can find another Littlewood--Paley function $\psi$ such that $(\varphi,\psi)$ is a Littlewood--Paley pair.
By Theorem \ref{phi}, for $f\in\Sc_0'(\R^{\vn};\C^m)$, we then have
\begin{equation*}
\Norm{f}{\Atau(\V)}^{\varrho}
=\Norm{T_\psi\circ S_\varphi f}{\Atau(\V)}^{\varrho}
\lesssim\Norm{S_\varphi f}{\atau(\V)}
\lesssim\Norm{f}{\Atau(\V)}^{\widetilde\varphi}.
\end{equation*}
Applying this to the Littlewood--Paley function $\widetilde\varphi$ in place of $\varphi$ and noting that $\widetilde{\widetilde{\varphi}}=\varphi$, we obtain $\Norm{f}{\Atau(\V)}^{\varrho} \lesssim\Norm{f}{\Atau(\V)}^{\varphi}$, and by symmetry also the reverse estimate.
\end{proof}

While the space $\Atau(\V)$ is independent of the choice of the Littlewood--Paley function $\varphi$ in its definition, the quasi-norms $\Norm{\ }{\Atau(\V)}^{\varphi}$ are only equivalent, not equal, for different choices of $\varphi$. Nevertheless, this ambiguity is relatively minor, and we mostly drop the Littlewood--Paley function from the notation of the quasi-norms $\Norm{\ }{\Atau(\V)^{\varphi}}$, writing just $\Norm{\ }{\Atau(\V)}$, following the usual practice in the area.

\begin{corollary}\label{convergence}
Let $s\in\R^k$, $\tau\in[0,\infty)$, $\vp\in(0,\infty)^k$, $\vq\in(0,\infty]^k$, and $\pi\in S_{[2k]}$ be admissible for $(\vp,\vq)$.
Let $\V=\{V_{\vj}\}_{\vj\in\Z^k}$, where each $V_{\vj}$ takes a constant value $V_R\in\C^{m\times m}$ on each $R\in\D_{\vj}(\R^{\vn})$, and let $\{V_R\}_{R\in\D(\R^{\vn})}$ be strongly doubling (Definition \ref{doubling}).
If both $\{f_N\}_{N\in\Z_+}$ and $f$ belong to $\Atau(\V)$ and $\lim_{N\to\infty}\Norm{f_N-f}{\Atau(\V)}=0$,
then $f_N\to f$ in $\Sc_0'(\R^{\vn};\C^m)$ as $N\to\infty$.
In particular, this convergence holds
if $\V\equiv1$ or $\vp=p\cdot\vec1$ for some $p\in(0,\infty)$, $V\in\A_p(\R^{\vn})$,
and $\V=[V]_p$ (i.e., $V_R=[V]_{\aveL^p(R)}$).
\end{corollary}

\begin{proof}
Once the first claim is proved, the last claim is immediate from the strong doubling property of the reducing operators of $\A_p$ weights (Remark \ref{Ap doubling}).

For $f\in\Atau(\V)\subset\Sc_0'(\R^{\vn};\C^m)$ and $\phi\in\Sc_0(\R^{\vn})$, we can apply the Calder\'on reproducing formula \eqref{cald20} of Lemma \ref{cald2} to see that, for any Littlewood--Paley pair $(\varphi,\psi)$,
\begin{equation*}
\pair{f}{\phi}
=\sum_{R\in\D(\R^{\vn})}\pair{f}{\varphi_R}\pair{\psi_R}{\phi}
=\sum_{R\in\D(\R^{\vn})}t_R\pair{\psi_R}{\phi},
\end{equation*}
where $t_R:=\pair{f}{\varphi_R}=(S_\varphi f)_R$ are the $\varphi$-transform coefficients of $f$. A combination of Lemma \ref{Tpsi welldef} and the $\varphi$-transform characterization Theorem \ref{phi} shows that
\begin{equation*}
\abs{\pair{f}{\phi}}
\lesssim\Norm{t}{\atau(\V)}
=\Norm{S_\varphi f}{\atau(\V)}
\sim\Norm{f}{\Atau(\V)},
\end{equation*}
where the implicit positive constants may depend on $\phi$ but not on $f$. Applying this to $f_N-f$ in place of $f$ shows that $\Norm{f_N-f}{\Atau(\V)}\to 0$ implies $\pair{f_N}{\phi}\to\pair{f}{\phi}$ as $N\to\infty$.
As this holds for every $\phi\in\Sc_0(\R^{\vn})$, this is the claimed convergence $f_N\to f$ in $\Sc_0'(\R^{\vn};\C^m)$ as $N\to\infty$.
\end{proof}

For all $x:=(x_1,\ldots,x_k)\in\mathbb R^{\vec n}$
and $\vec\sigma:=(\sigma_1,\ldots,\sigma_k)\in\mathbb R^k$, let
$$
|x|^{\vec\sigma}:=\prod_{i=1}^k|x_i|^{\sigma_i}.
$$
For every $\vec\sigma\in\mathbb R^k$,
the \emph{multi-parameter lifting operator} $\dot I_{\vec\sigma}$
(see, for instance, \cite[Section 5.2.3]{Tri83})
is defined by setting, for all $f\in \mathscr S_0'(\mathbb R^{\vec n})$,
\begin{equation}\label{lift}
\dot I_{\vec\sigma}(f):=\Big(|\cdot|^{-\vec\sigma}
\widehat{f}\Big)^{\vee},
\end{equation}
where the symbol $\vee$ denotes the \emph{inverse Fourier transform},
which is defined by setting, for any $f\in\mathscr S'(\R^{\vn})$ and $\xi\in\R^{\vn}$,
$\check{f}(\xi):=(2\pi)^{-n}\widehat{f}(-\xi)$.
Applying a basic calculation, we obtain that $\dot I_{\vec\sigma}$ maps
$\mathscr S_0'(\mathbb R^{\vec n})$ onto itself.

\begin{corollary}\label{257}
Let $s\in\R^k$, $\tau\in[0,\infty)$, $\vp\in(0,\infty)^k$, $\vq\in(0,\infty]^k$, and $\pi\in S_{[2k]}$ be admissible for $(\vp,\vq)$.
Let $\V=\{V_{\vj}\}_{\vj\in\Z^k}$, where each $V_{\vj}$ takes a constant value $V_R\in\C^{m\times m}$ on each $R\in\D_{\vj}(\R^{\vn})$, and let $\{V_R\}_{R\in\D(\R^{\vn})}$ be strongly doubling (Definition \ref{doubling}).
Then, for every $\sigma\in\R^k$, the lifting operator $\dot I_\sigma$ maps $\dot A_{\vec p,\vec q,\pi}^{\vec s,\tau}(\V)$ isomorphically onto $\dot A_{\vec p,\vec q,\pi}^{\vec s+\vec \sigma,\tau}(\V)$.
Moreover, for every $f\in\mathscr{S}_0'(\mathbb R^{\vec n};\C^m)$,
$$
\|\dot I_\sigma f\|_{\dot A_{\vec p,\vec q,\pi}^{\vec s+\vec \sigma,\tau}(\V)}
\sim \|f\|_{\dot A_{\vec p,\vec q,\pi}^{\vec s,\tau}(\V)},
$$
where the positive equivalence constants are independent of $f$.
In particular, these results hold
if $\V\equiv1$ or $\vp=p\cdot\vec1$ for some $p\in(0,\infty)$, $V\in\A_p(\R^{\vn})$,
and $\V=[V]_p$ (i.e., $V_R=[V]_{\aveL^p(R)}$).
\end{corollary}

\begin{proof}
Let $\varphi$ be a Littlewood--Paley function on $\mathbb R^{\vec n}$.
Observe that, by the definition of $\dot I_\sigma$
and properties of the inverse Fourier transform, we have,
for every $\vj\in\mathbb Z^k$ and $f\in \mathscr{S}_0'(\mathbb R^{\vec n};\C^m)$,
\begin{align}\label{258}
\varphi_{\vj}*\Big(\dot I_\sigma f\Big)
= \varphi_{\vj}*\Big(|\cdot|^{-\vec\sigma}\widehat{f}\Big)^\vee
= \Big(\widehat{\varphi_{\vj}}
|\cdot|^{-\vec\sigma}\widehat{f}\Big)^\vee.
\end{align}
Let $\psi:=(|\cdot|^{-\vec\sigma}\widehat{\varphi})^\vee$.
Since $\varphi$ is a Littlewood--Paley function on $\mathbb R^{\vec n}$,
it follows that $\psi$ is also a Littlewood--Paley function on $\mathbb R^{\vec n}$.
Moreover, note that
\begin{align*}
\psi_{\vj}*f
&=\Big(\widehat{\psi_{\vj}}\widehat{f}\Big)^\vee
=\Big[ 2^{-\vj\cdot\vec n}\widehat{\psi}\Big( 2^{-\vj}\cdot\Big)\widehat{f}\Big]^\vee\\
&=2^{\vj\cdot\vec\sigma}\Big[ 2^{-\vj\cdot\vec n} |\cdot|^{-\vec\sigma}
\widehat{\varphi}\Big( 2^{-\vj}\cdot\Big)\widehat{f}\Big]^\vee
=2^{\vj\cdot\vec\sigma}\Big( |\cdot|^{-\vec\sigma}
\widehat{\varphi_{\vj}}\widehat{f}\Big)^\vee.
\end{align*}
From this and \eqref{258}, we infer that $\varphi_{\vj}*(\dot I_{\vec\sigma} f)=2^{-\vj\cdot\vec\sigma}\psi_{\vj}*f$ and hence
$$
\|\dot I_\sigma f\|^\varphi_{\dot A_{\vec p,\vec q,\pi}^{\vec s+\vec \sigma,\tau}(\V)}
=\|f\|^\psi_{\dot A_{\vec p,\vec q,\pi}^{\vec s,\tau}(\V)}.
$$
Noting that the assumptions of Corollary \ref{257} include those of Corollary \ref{38},
we can apply the said corollary to conclude that both quasi-norms above are
independent of the particular Littlewood--Paley functions $\varphi,\psi$.
This finishes the proof of Corollary \ref{257}.
\end{proof}

Due to Theorem \ref{3 norms thm},
all the results in this subsection
remain valid with $\Atau(\V)$ and $\atau(\V)$ replaced by $\Atau(V)$ and $\atau(V)$.
We list the results that will be used later, but omit the proofs.

\begin{corollary}\label{phiV}
Let $s\in\R^k$, $\tau\in[0,\infty)$, $\vp\in(0,\infty)^k$, $\vq\in(0,\infty]^k$, and $\pi\in S_{[2k]}$ be admissible for $(\vp,\vq)$.
Let $\varphi,\psi,\varrho\in\Sc (\R^{\vn})$ be Littlewood--Paley functions, and $\widetilde\varphi(x):=\overline{\varphi(-x)}$ for every $x\in\mathbb{R}^n$.
If $V\equiv1$ or if $\vp=p\cdot\vec 1$, $V\in \A_p(\R^{\vn})$,
and $(\vp,\vq,\pi)$ satisfy \eqref{key cases},
then the operators
\begin{equation*}
S_\varphi:\ \Atau(V)^{\widetilde{\varphi}}\to\atau(V)
\quad\text{and}\quad
T_\psi:\ \atau(V)\to \Atau(V)^{\varrho}
\end{equation*}
are bounded.
If $(\varphi,\psi)$ is a Littlewood--Paley pair, then
$T_\psi\circ S_\varphi$ is the identity on $\Atau(V)^{\widetilde{\varphi}}$.
\end{corollary}

\begin{corollary}\label{257V}
Let $s\in\R^k$, $\tau\in[0,\infty)$, $\vp\in(0,\infty)^k$, $\vq\in(0,\infty]^k$, and $\pi\in S_{[2k]}$ be admissible for $(\vp,\vq)$.
If $V\equiv1$ or if $\vp=p\cdot\vec 1$, $V\in \A_p(\R^{\vn})$,
and $(\vp,\vq,\pi)$ satisfy \eqref{key cases},
then, for every $\sigma\in\R^k$, $\dot I_\sigma$ maps $\dot A_{\vec p,\vec q,\pi}^{\vec s,\tau}(V)$ isomorphically onto $\dot A_{\vec p,\vec q,\pi}^{\vec s+\vec \sigma,\tau}(V)$.
Moreover, for every $f\in\mathscr{S}_0'(\mathbb R^{\vec n};\C^m)$,
$$
\|\dot I_\sigma f\|_{\dot A_{\vec p,\vec q,\pi}^{\vec s+\vec \sigma,\tau}(V)}
\sim \|f\|_{\dot A_{\vec p,\vec q,\pi}^{\vec s,\tau}(V)},
$$
where the positive equivalence constants are independent of $f$.
\end{corollary}

\section{Almost diagonal operators}\label{sec AD}

Given an infinite matrix $B=\{b_{PR}\}_{P,R\in\D(\R^{\vn})}$ with $b_{PR}\in\F$ and an infinite matrix $t=\{t_R\}_{R\in\D(\R^{\vn})}$ with $t_R\in\F^m$, both indexed by the dyadic rectangles $\D(\R^{\vn})$, we say that $Bt$ is well defined, if
\begin{equation*}
(Bt)_P:=\sum_{R\in\D(\R^{\vn})}b_{PR}t_R
\end{equation*}
is absolutely convergent for every $P\in\D(\R^{\vn})$, and in this case we define the new infinite vector $Bt=\{(Bt)_P\}_{P\in\D(\R^{\vn})}$ by the formula just given. A basic question is to provide sufficient conditions on $B$ that guarantee its bounded action on all $t$ in a given sequence space.
Starting with \cite{FJ90}, this question has been extensively studied in the one-parameter Besov and Triebel--Lizorkin sequence spaces,
and useful conditions for the boundedness are known in terms of {\em almost diagonal} estimates of $B$.
The aim of this section is to extend this theory to the matrix-weighted multi-parameter Besov--Triebel--Lizorkin-type sequence spaces $\atau(\V)$. Thanks to the $\varphi$-transform characterisation (Theorem \ref{phi}), this will also provide a tool for analysing various operators on the corresponding function or distribution spaces $\Atau(\V)$, and we will make extensive use of this in the subsequent sections.

Following \cite[Definition 3.1]{BHYY:1b} (which, in turn, is based on several earlier works, notably \cite[p.~53]{FJ90}), we define the special matrices
\begin{equation}\label{bDEF}
b^{DEF}_{i;PR}:=\bigg[1+\frac{\abs{c_P-c_R}}{\ell(P)\vee\ell(R)}\bigg]^{-D}
\begin{cases}
\displaystyle \Big[\frac{\ell(P)}{\ell(R)}\Big]^E,
& \text{if }\ell(P)\leq\ell(R), \\
\displaystyle \Big[\frac{\ell(R)}{\ell(P)}\Big]^F,
& \text{if }\ell(R)<\ell(P),
\end{cases}
\end{equation}
with parameters $D,E,F\in\R$ and indexed by dyadic cubes $P,R\in\D(\R^{n_i})$. On the product space $\R^{\vn}$, we then set
\begin{equation*}
b^{\vD \vE \vF}_{\vn;PR}:=\prod_{i=1}^k b^{D_i,E_i,F_i}_{i;P_i,R_i},
\end{equation*}
with vector parameters $\vD=(D_i)_{i=1}^k,\vE=(E_i)_{i=1}^k,\vF=(F_i)_{i=1}^k\in\R^k$ and indexed by dyadic rectangles $P,R\in\D(\R^{\vn})$, where e.g.\ $P=P_1\times\cdots\times P_k$ with $P_i\in\mathscr D(\R^{n_i})$ for each $i\in[k]$.

A matrix $B=(b_{PR})_{P,R\in\D(\R^{\vn})}$ is said to be $(\vD,\vE,\vF)$-almost diagonal, if
\begin{equation*}
\abs{b_{PR}}\lesssim b^{\vD \vE \vF}_{\vn;PR},
\end{equation*}
with the implicit positive constant independent of $P$ and $R$.
Starting with \cite{FJ90}, the one-parameter version of this notion
has been studied in several articles. The rather natural extension to several parameters just given seems to have been first introduced in \cite[Definition 5.1]{GKP21} in the case of $k=2$ and a specific choice of $(\vD,\vE,\vF)$. Indeed, they proved in \cite[Theorem 5.2]{GKP21} that if $B$ is $(\vD,\vE,\vF)$-almost diagonal with
\begin{equation}\label{GKPcond}
\vD>\frac{\vn}{r},\quad
\vE>\vs+\frac{\vn}{2},\quad
\vF>\frac{\vn}{r}-\vs-\frac{\vn}{2},
\end{equation}
then $B$ is bounded on $\dot a^{\vs}_{p,q}(\R^{\vn})$,
where $r=\min(1,p)$ in the Besov case,
and $r=\min(1,p,q)$ in the Triebel--Lizorkin case.

We will now discuss the extension of this result to the matrix-weighted multi-parameter spaces.
The remainder of this section is organized as follows.
In Subsection \ref{1st approach}, we present the first approach
to establishing the boundedness of almost diagonal operators on $\atau(\V)$ through a reduction to the unweighted case.
After analysing the advantages and the disadvantages of this approach, we introduce the second approach to this boundedness:
in Subsection \ref{2nd approach} for the case $\tau=0$,
and in Subsection \ref{tau approach} for the case $\tau\in(0,\infty)$.
Finally, in Subsection \ref{sharpness},
we prove the sharpness of the boundedness result obtained via the second approach.

\subsection{The first approach}\label{1st approach}

Using the quasi-norm equivalence of Corollary \ref{3 norms cor}, followed by the identification of $[V]_p=\{\Red{V}{\vj}{p}\}_{\vj\in\Z^k}$ with the diagonal matrix $\{[V]_{\aveL^p(R)}\delta_{P,R}\}_{P,R\in\D(\R^{\vn})}$ as explained in Definition \ref{ave}, we obtain
\begin{equation*}
\Norm{t}{\atau(V)}\sim\Norm{t}{\atau([V]_p)}
=\Norm{[V]_p t}{\atau(\R^{\vn})}.
\end{equation*}
Hence, a prospective estimate
$  \Norm{Bt}{\atau(V)}\lesssim\Norm{t}{\atau(V)}$
is equivalent to the unweighted bound
$  \Norm{\widetilde B t}{\atau(\R^{\vn})}\lesssim\Norm{t}{\atau(\R^{\vn})}$
for the new matrix
\begin{equation}\label{tildeB}
\widetilde B:=[V]_p B [V]_p^{-1} =\Big\{ [V]_{\aveL^p(P)} b_{P,R} [V]_{\aveL^p(R)}^{-1}
\Big\}_{ P,R\in\D(\R^{\vn})}.
\end{equation}

\begin{lemma}\label{lem AD simple}
Let $\vD,\vE,\vF\in \mathbb R^k$.
If $B$ is $(\vD,\vE,\vF)$-almost diagonal and $V\in\A_p(\R^{\vn})$ has $\A_p$-dimensions $(\vd,\ve,\vDe)$, then $\widetilde B$ defined in \eqref{tildeB} is
$(\vD-\vDe,\vE-\frac{\vd}{p},\vF-\frac{\ve}{p'})$-almost diagonal.
In particular, $\widetilde B$ is  $(\vD-\frac{\vn}{\min(1,p)},\vE-\frac{\vn}p,\vF-\frac{\vn}{p'})$-almost diagonal for every $V\in\A_p(\R^{\vn})$.
\end{lemma}

\begin{proof}
This is immediate from Definition \ref{def Ap dim} of $\A_p$-dimensions,
the related estimate for the term $\abs{[V]_{\aveL^p(P)} [V]_{\aveL^p(R)}^{-1}}$
provided by Lemma \ref{sharp coef}, and the general upper
bound on $\A_p$-dimensions from Proposition \ref{d<n}. Note in particular that
\begin{equation*}
\vDe=\frac{\vd}{p}+\frac{\ve}{p'}
<\vn \Big( \frac{1}{p}+\frac{1}{p'} \Big)
=\left.\begin{cases} \vn, & \text{if}\quad p\in(1,\infty) \\
{\displaystyle\frac {\vn}{p}}, & \text{if}\quad p\in(0,1]\end{cases}\right\}
=\frac{\vn}{\min(1,p)}.
\end{equation*}
This finishes the proof of Lemma \ref{lem AD simple}.
\end{proof}

\begin{corollary}\label{cor AD simple}
Let $\vs\in\R^k$, $\tau\in[0,\infty)$, $\vp=p\cdot\vone$ with $p\in(0,\infty)$, and $\vq\in(0,\infty]^k$.
Suppose that, for some $(\vD,\vE,\vF)\in\R^{3k}$, every $(\vD,\vE,\vF)$-almost diagonal matrix defines a bounded operator on the unweighted $\atau(\R^{\vn})$. Then
\begin{enumerate}[\rm(i)]
\item every $(\vD+\frac{\vn}{\min(1,p)},\vE+\frac{\vn}{p},\vF+\frac{\vn}{p'})$-almost diagonal matrix acts boundedly on $\atau([V]_p)$ for all $V\in\A_p(\R^{\vn})$;

\item every $(\vD+\vDe,\vE+\frac{\vd}{p},\vF+\frac{\ve}{p'})$-almost diagonal matrix acts boundedly on $\atau([V]_p)$ for all $V\in\A_p(\R^{\vn})$ with $\A_p$-dimensions $(\vd,\ve,\vDe)$.
\end{enumerate}
\end{corollary}

\begin{proof}
This is immediate from Lemma \ref{lem AD simple} and the preceding discussion.
\end{proof}

The idea behind Corollary \ref{cor AD simple}, of reducing the weighted boundedness of a matrix $B$ to the unweighted boundedness of the conjugated matrix \eqref{tildeB} and deriving conditions for the latter from an estimate in the spirit of Lemma \ref{sharp coef}, is from \cite[Theorem 2.6]{FR21}. Several {\em advantages} of this method are the following:
\begin{enumerate}[\rm(i)]
\item The proof of the matrix-weighted estimates is very simple, provided that corresponding unweighted results are already available as input.
\item Information about unweighted estimates can be used as a black box, without having to revisit the details of their proof. In particular, the sufficient condition \eqref{GKPcond} of \cite[Theorem 5.2]{GKP21} for $\dot a^{\vs}_{p,q}$-boundedness immediately yields conditions for $\dot a^{\vs}_{p,q}(V)$-boundedness via Corollary \ref{cor AD simple}. Moreover, any future results in the unweighted case may be equally easily used to obtain new corollaries.
\item The method applies to Besov--Triebel--Lizorkin-type spaces $\atau(V)$ of quite general form, including mixed versions of not necessarily pure Besov or pure Triebel--Lizorkin-type.
\end{enumerate}

At the same time, this approach also has some {\em disadvantages}:
\begin{enumerate}[\rm(i)]
\item It requires unweighted results as input. In contrast to the one-para\-meter case, where such information is readily available in the literature (see \cite[Theorem 3.3]{FJ90} for the classical spaces with $\tau=0$, \cite[Theorem 4.1]{YY10} for $\tau<\frac1p+\varepsilon$, and \cite[Theorem 1]{YY13} for a reduction of the remaining cases to those covered in \cite{FJ90}), the available results about almost diagonal operators in the multi-parameter setting seem to be much fewer at the time of writing. (The works \cite{CLL16,DLM10,LuZhu13,Xu22}, while proving the boundedness of various classical operators on multi-parameter Besov and Triebel--Lizorkin spaces, take a somewhat different approach that does not address almost diagonal operators on the related sequence spaces.)

\item While qualitatively covering a wide range of spaces $\atau(V)$, the quantitative assumptions [i.e., the conditions on the almost diagonal parameters $(\vD,\vE,\vF)$] imposed by this method are stronger than necessary when specialised to the arguably most important case of pure Besov or pure Triebel--Lizorkin spaces. This suboptimality of the simple approach was already observed in the one-parameter case in \cite{BHYY:1b}, where a direct analysis of almost diagonal operators in the matrix-weighted spaces gave a result, \cite[Theorem 13.1]{BHYY:1b}, that not only improved on the one-parameter analogue of Corollary \ref{cor AD simple}, \cite[Theorem 4.7]{BHYY:1b}, but also gave the best possible result in several ranges of the parameters, as explained in \cite[Remark 13.2]{BHYY:1b}.
\end{enumerate}

To address these shortcomings of the simple method, we next turn to discuss a multi-parameter extension of the direct approach of \cite{BHYY:1b}.

\subsection{The second approach}\label{2nd approach}

We will now develop another approach to the boundedness of almost diagonal operators on $\atau(V)$ following \cite{BHYY:1b}, which gives the best available (and in several cases the best possible) results in the one-parameter case. We will now concentrate on the following {\em main cases}:

\begin{definition}\label{main cases}
Let $\vp\in(0,\infty)^k$, $\vq\in(0,\infty]^k$, $\pi\in S_{[2k]}$ be admissible, and $V$ be a matrix weight on $\R^{\vn}$. We say that $(\vp,\vq,\pi,V)$ belongs to {\em one of the main cases} if
either
\begin{enumerate}[\rm(i)]
\item\label{unweighted case} $V\equiv 1$ and the other parameters are arbitrary (unweighted case), or

\item\label{Besov case} $\vp=p\cdot\vone$, the weight is $V\in\A_p(\R^{\vn})$, and $(L^{\vp}\ell^{\vq})_\pi =\ell^{\vq}L^p$ (Besov case), or

\item\label{TL case} $\vp=p\cdot\vone$, the weight is $V\in\A_p(\R^{\vn})$, and $(L^{\vp}\ell^{\vq})_\pi =L^p\ell^{\vq}$ (Triebel--Lizorkin case).
\end{enumerate}
\end{definition}

The motivation for cases \eqref{Besov case} and \eqref{TL case} should be clear, as these correspond to the typical cases of interest in the theory of function spaces. The motivation for also addressing case \eqref{unweighted case} is related to our
first approach in Subsection \ref{1st approach}: as soon as we have estimates for almost diagonal operators on the unweighted $\atau(\R^{\vn})$, these immediately imply estimates on the matrix weighted $\atau(V)$, and not only those of pure Besov or pure Triebel--Lizorkin-type. On the other hand, as noted in Subsection \ref{1st approach}, the existing results even this unweighted case are not so many, and thus it is useful to provide new input that can be fed into the general machine of Corollary \ref{cor AD simple}.

\begin{remark}\label{s vs 0}
The case of general $\vs\in\R^k$ may be reduced to $\vs=\vnull$ by an argument essentially similar to the reduction of the weighted case to the unweighted one in the first approach. Indeed,
\begin{equation*}
\Norm{t}{\atau(V)}=\Norm{\lift{\vs}t}{\antau(V)},\qquad
\lift{\vs}:=\Bigg\{\delta_{P,R}\prod_{i=1}^k\ell(P_i)^{-s_i}\Bigg\}_{P,R\in\D(\R^{\vn})},
\end{equation*}
which means that
\begin{equation*}
\Norm{Bt}{\atau(V)}\lesssim\Norm{t}{\atau(V)}
\Longleftrightarrow
\Norm{B_{\vs}t}{\antau(V)}\lesssim\Norm{t}{\antau(V)},\quad
B_{\vs}:=\lift{\vs}B\lift{-\vs},
\end{equation*}
where
\begin{align*}
&B_{\vs}\text{ is } (\vD,\vE,\vF)\text{-almost diagonal }  \\
&\quad\Longleftrightarrow
B\text{ is }(\vD,\vE+\vs,\vF-\vs)\text{-almost diagonal}.
\end{align*}
\end{remark}

The following lemma is essentially a combination of algebraic manipulations and elementary inequalities; it is stated for general sequences of matrix-valued multiplying sequences to emphasise the fact that no special assumptions on the matrix weights, like $\A_p$ or similar, are used at this point.

\begin{lemma}\label{key lem}
Let $\vp\in(0,\infty)^k$, $\vq\in(0,\infty]^k$, $\pi\in S_{[2k]}$, and $B$ be $(\vD,\vE,\vF)$-almost diagonal for some $\vD,\vE,\vF\in\R^k$.
Let $a\in(0,1]$, let $H_{\vj}:\R^n\to\F^{m\times m}$ be measurable and almost everywhere invertible, and let $0<u\leq\min\{1,p_i,q_i:i\in[k]\}$.
If $t=\{t_R\}_{R\in\D(\R^{\vn})}$ is such that the right-hand side of \eqref{key eq} below is finite, then
$Bt$ is well defined and
\begin{equation}\label{key eq}
\begin{split}
\Norm{\{H_{\vj}(Bt)_{\vj}\}}{(L^{\vp}\ell^{\vq})_\pi}^u
\lesssim &\sum_{\vh\in\Z^k}\sum_{\vl\in\N^k}
\Bigg[ 2^{-(\vE-\frac{\vn}{2})\cdot\vh_-}2^{-(\vF+\frac{\vn}{2}-\frac{\vn}{a})\cdot\vh_+}
2^{-(\vD-\frac{\vn}{a})\cdot\vl}
\phantom{\fint}\Big. \\
&\times \Big.\BNorm{
\Big\{\Big[\fint_{C(\cdot,2^{\vl+\vh_+-\vi})}\abs{H_{\vi-\vh}(\cdot)t_{\vi}(y)}^a dy\Big]^{\frac1a}\Big\}_{\vi\in\Z^k}
}{(L^{\vp}\ell^{\vq})_\pi}\Bigg]^u,
\end{split}
\end{equation}
where the implicit positive constant is independent of $t$,
\begin{equation*}
C(x,2^{\vl+\vh_+-\vi})
:=\prod_{\nu=1}^k B_\nu(x_{\nu},2^{l_{\nu}+(h_{\nu})_+-i_{\nu}}),\quad
x=(x_1,\ldots,x_k),
\end{equation*}
and each $B_\nu$ is a ball in the respective $\R^{n_{\nu}}$.
\end{lemma}

\begin{proof}
Case $k=1$ is \cite[Lemma 5.4 and Remark 5.5]{BHYY:1b}. The general case follows by routine modifications of the same argument.
\end{proof}

\begin{remark}
As in \cite[Lemma 5.4]{BHYY:1b}, we have formulated Lemma \ref{key lem} for general functions $H_{\vj}$, although we will actually only apply it with $H_{\vj}\equiv V\in\A_p(\R^{\vn})$. Another natural choice that suggests itself would be to take $H_{\vj}=\Red{V}{\vj}{p}$, leading to quantities related to the quasi-norm $\Norm{\ }{\atau([V]_p)}$.

However, one should observe the mismatch of the indices $\vi-\vh$ and $\vi$ on the right-hand side of \eqref{key eq}. This is of no concern when all $H_{\vi-\vh}\equiv V$ are equal anyway. On the other hand,
with $H_{\vj}=\Red{V}{\vj}{p}$, any attempts of estimating the right-hand side of \eqref{key eq} by $\Norm{[V]_p t}{(L^{\vp}\ell^{\vq})_\pi}$ would easily ask for an application of Lemma \ref{sharp coef} which, in turn, would introduce dependence on the $\A_p$-dimensions $(\vd,\ve,\vDe)$. Working with $H_{\vj}\equiv V$ instead, we are led to applying the Carleson estimates of Corollary \ref{Carl Ap} and Proposition \ref{dualCarl} instead of Lemma \ref{sharp coef}, and this has the advantage of avoiding dependence on the $\A_p$-dimensions.
\end{remark}

Taking $H_{\vj}\equiv V\in\A_p(\R^{\vn})$, we will then estimate the quasi-norms on the right-hand side of \eqref{key eq}. Let us simply write $\vl$ in place of $\vl+\vh_+$; in any case, this is an arbitrary vector in $\N^k$. In our treatment of the one-parameter case in \cite{BHYY:1b}, we gave a separate argument in the Besov and in the Triebel--Lizorkin case, the former one being somewhat easier. Here, we treat all cases at once, and leave it to the interested reader to work out the simplifications that one can make if one is only interested in the Besov case.

\begin{lemma}\label{key lem 2}
Let $\vp=p\cdot\vone$ with $p\in(0,\infty)$, let $\vq\in(0,\infty]^k$, and let $V\in\A_p(\R^{\vn})$.
Let $\pi\in S_{[2k]}$ be an admissible permutation of $(\vp,\vq)$ (Definition \ref{adm}) and $a\in(0,r(\pi,\vp,\vq))$ (see \eqref{rpipq-max}).
Suppose, moreover, that $(\vp,\vq,\pi,V)$ belongs to one of the main cases (Definition \ref{main cases}).
Then all $t\in(\C^m)^{\D(\R^{\vn})}$ and $\vl\in\N^k$ satisfy
\begin{equation*}
\bigg\|
\bigg\{\bigg[\fint_{C(\cdot,2^{\vl-\vi})}\abs{V(\cdot)t_{\vi}(y)}^a\, dy\bigg]^{\frac1a}\bigg\}_{\vi\in\Z^k}
\bigg\|_{(L^{\vp}\ell^{\vq})_\pi}
\lesssim\Norm{Vt}{(L^{\vp}\ell^{\vq})_\pi},
\end{equation*}
where the implicit positive constant is independent of $t$ and $\vl$.
\end{lemma}

\begin{remark}
Note that $r(\pi,\vp,\vq)=p$ in the Besov case, and $r(\pi,\vp,\vq)=\min\{p,q_\nu:\nu\in[k]\}$ in the Triebel--Lizorkin case. Moreover, the condition that $\pi$ be admissible is automatic in the Besov case, while in the Triebel--Lizorkin case it requires that $\vq\in(0,\infty)^{k-h}\times\{\infty\}^{h}$ for some $h\in\{0,1,\ldots,k\}$, i.e., all infinite components of $\vq$, if any, are at the end of the sequence.
\end{remark}

\begin{proof}[Proof of Lemma \ref{key lem 2}]
We may clearly assume that $\Norm{Vt}{(L^{\vp}\ell^{\vq})_\pi}<\infty$, for otherwise there is nothing to prove.

Let us first take care of the unweighted case $V\equiv 1$. Then we can simply dominate the averages over the rectangles $C(\cdot,2^{\vl-\vi})$ by the strong maximal operator $\mathcal M_{\vn}$ and apply Theorem \ref{rescaled M} to deduce that
\begin{align*}
\bigg\|\Big\{\Big[\fint_{C(\cdot,2^{\vl-\vi})}\abs{t_{\vi}(y)}^a dy\Big]^{\frac1a}\Big\}_{\vi\in\Z^k}
\bigg\|_{(L^{\vp}\ell^{\vq})_\pi}
&\leq\Norm{\{[\mathcal M_{\vn}(\abs{t_{\vi}}^a)]^{\frac1a}\}_{\vi\in\Z^k}}{(L^{\vp}\ell^{\vq})_\pi}
=\Norm{\{\mathcal M_{\vn}(\abs{t_{\vi}}^a)\}_{\vi\in\Z^k}}{(L^{\frac{\vp}{a}}\ell^{\frac{\vq}{a}})_\pi} ^{\frac1a} \\
&\lesssim\Norm{\{\abs{t_{\vi}}^a\}_{\vi\in\Z^k}}{(L^{\frac{\vp}{a}}\ell^{\frac{\vq}{a}})_\pi} ^{\frac1a}
=\Norm{\{t_{\vi}\}_{\vi\in\Z^k}}{(L^{\vp}\ell^{\vq})_\pi},
\end{align*}
which completes the short argument in the unweighted case.

We then turn to the more difficult matrix-weighted cases, where we assume that $(L^{\vp}\ell^{\vq})_\pi$ takes one of the forms $L^p\ell^{\vq}$ or $\ell^{\vq}L^p$.

As in the proof of \cite[Lemma 8.5]{BHYY:1b}, for each $\nu\in[k]$, the ``one-third trick'' \cite[Lemma 3.2.26]{HNVW1} allows one to pick $3^{n_{\nu}}$ shifted systems of dyadic cubes $\D^{\alpha_\nu}(\R^{n_\nu})$, $\alpha_\nu\in\{0,\frac13,\frac23\}^{n_\nu}$, so that each ball $B_{\nu}(x_\nu,2^{l_{\nu}-i_\nu})$ is contained in some shifted dyadic cube $S_{\nu}\in\D^{\alpha_\nu}_{i_\nu-l_{\nu}-2}(\R^{n_{\nu}})$; thus, $C(x,2^{\vl-\vi})$ is contained in
\begin{equation*}
S:=\prod_{\nu=1}^k S_{\nu}\in\prod_{\nu=1}^k \D^{\alpha_\nu}_{i_\nu-l_{\nu}-2}(\R^{n_{\nu}})
=:\D^{\val}_{\vi-\vl-2\cdot\vone}(\R^{\vn}),\quad
\val=(\alpha_\nu)_{\nu=1}^k\in\{0,\tfrac13,\tfrac23\}^{\vn},
\end{equation*}
where each $\D^{\val}(\R^{\vn}):=\bigcup_{\vi\in\Z^k}\D^{\alpha_\nu}_{i_\nu}(\R^{n_{\nu}})$ is a shifted system of dyadic rectangles on $\R^{\vn}$. Then
\begin{align*}
\bigg[\fint_{C(x,2^{\vl+\vh_+-\vi})}\abs{V(\cdot)t_{\vi}(y)}^a\, dy\bigg]^{\frac1a}
\lesssim\sum_{\val\in\{0,\frac13,\frac23\}^n}\sum_{S\in\D^{\val}_{\vi-\vl-2\cdot\vone}
(\R^{\vn})}\one_S(\cdot)
\bigg[\fint_{S}\abs{V(\cdot)t_{\vi}(y)}^a \, dy\bigg]^{\frac1a}.
\end{align*}
Thanks to the quasi-triangle inequality, it is enough to estimate the quantity on the right for each fixed $\val$, and we will just write $\D(\R^{\vn})$ instead of $\D^{\val}(\R^{\vn})$.

The next issue is the quantity $V(x)t_{\vi}(y)$, where $V$ and $t_{\vi}$ are evaluated at different points, in contrast to how they should appear in $\Norm{Vt}{(L^{\vp}\ell^{\vq})_\pi}$. To fix this, we will make use of the decomposition
\begin{equation}\label{dec3}
V(x)t_{\vi}(y)
= \Big(V(x)[V]_{\aveL^p(S)}^{-1}\Big) \Big( [V]_{\aveL^p(S)} V(y)^{-1}\Big) \Big[ V(y)t_{\vi}(y) \Big].
\end{equation}
Extracting the first factor and applying Corollary \ref{Carl Ap}, we obtain
\begin{align*}
& \bbNorm{\Big\{ \sum_{S\in\D_{\vi-\vl-2\cdot\vone}(\R^{\vn})}\one_S
\Big(\fint_{S}\abs{V(\cdot)t_{\vi}(y)}^a dy\Big)^{\frac1a} \Big\}_{\vi\in\Z^k} }{(L^{\vp}\ell^{\vq})_\pi} \\
&\quad\leq\bbNorm{\Big\{ \sum_{S\in\D_{\vi-\vl-2\cdot\vone}(\R^{\vn})}\one_S\abs{V(\cdot)[V]_{\aveL^p(S)}^{-1}}
\Big(\fint_{S}\abs{[V]_{\aveL^p(S)}t_{\vi}(y)}^a dy\Big)^{\frac1a} \Big\}_{\vi\in\Z^k} }{(L^{\vp}\ell^{\vq})_\pi} \\
&\quad\lesssim\bbNorm{\Big\{ \sum_{S\in\D_{\vi-\vl-2\cdot\vone}(\R^{\vn})}\one_S
\Big(\fint_{S}\abs{[V]_{\aveL^p(S)}t_{\vi}(y)}^a dy\Big)^{\frac1a} \Big\}_{\vi\in\Z^k} }{(L^{\vp}\ell^{\vq})_\pi}
\quad\text{by Corollary \ref{Carl Ap}} \\
&\quad=\bbNorm{\Big\{ \sum_{S\in\D_{\vi-\vl-2\cdot\vone}(\R^{\vn})}\one_S
\fint_{S}\abs{[V]_{\aveL^p(S)}t_{\vi}(y)}^a dy\Big\}_{\vi\in\Z^k} }{(L^{\frac{\vp}{a}}\ell^{\frac{\vq}{a}})_\pi}^{\frac1a}.
\end{align*}
To remove the second factor in \eqref{dec3}, we dualise with a suitable $\Norm{g}{(L^{(\frac{\vp}{a})'}\ell^{(\frac{\vq}{a})'})_\pi}=1$:
\begin{align*}
&\bbNorm{\Big\{ \sum_{S\in\D_{\vi-\vl-2\cdot\vone}(\R^{\vn})}\one_S
\fint_{S}\abs{[V]_{\aveL^p(S)}t_{\vi}(y)}^a \,dy\Big\}_{\vi\in\Z^k} }{(L^{\frac{\vp}{a}}\ell^{\frac{\vq}{a}})_\pi} \\
&\quad\sim\sum_{\vi\in\Z^k}\int_{\R^{\vn}}g_{\vi}(x)\bigg(\sum_{S\in\D_{\vi-\vl-2
\cdot\vone}(\R^{\vn})}\one_S(x)
\fint_{S}\abs{[V]_{\aveL^p(S)}t_{\vi}(y)}^a \,dy\bigg)dx \\
&\quad\leq\sum_{\vi\in\Z^k}\int_{\R^{\vn}}\bigg(\sum_{S\in\D_{\vi-\vl-2\cdot\vone}(\R^{\vn})} \abs{[V]_{\aveL^p(S)}V^{-1}(y)}^a\fint_S g_{\vi}(x)dx\, \one_S(y)\bigg)     \abs{V(y)t_{\vi}(y)}^a \,dy \\
&\quad\leq\bbNorm{\Big\{
\sum_{S\in\D_{\vi-\vl-2\cdot\vone}(\R^{\vn})} \abs{[V]_{\aveL^p(S)}V^{-1}(\cdot)}^a\fint_S g_{\vi}(x)dx\, \one_S(\cdot)\Big\}_{\vi\in\Z^k}}{
(L^{(\frac{\vp}{a})'}\ell^{(\frac{\vq}{a})'})_\pi}
\Norm{\{\abs{V(\cdot)t_{\vi}(\cdot)}^a\}_{\vi\in\Z^k}}{(L^{\frac{\vp}{a}}\ell^{\frac{\vq}{a}})_\pi} \\
&\quad=\bbNorm{\Big\{
\sum_{S\in\D_{\vi-\vl-2\cdot\vone}(\R^{\vn})} \abs{[V]_{\aveL^p(S)}V^{-1}(\cdot)}\Big(\fint_S g_{\vi}(x)dx\Big)^{\frac1a} \one_S(\cdot)\Big\}_{\vi\in\Z^k}}{
(L^{a(\frac{\vp}{a})'}\ell^{a(\frac{\vq}{a})'})_\pi}^{a}
\Norm{Vt}{(L^{\vp}\ell^{\vq})_\pi}^a.
\end{align*}
We now want to apply the ``dual'' Carleson embedding inequality, Proposition \ref{dualCarl}\eqref{dualCarl2}, with multipliers $\abs{[V]_{\aveL^p(S)}V^{-1}(\cdot)}$ and exponents $r:=a(\frac{p}{a})'$ as well as $a(\frac{\vq}{a})'$ in place of $\vq$.

Let us check that the assumptions of Proposition \ref{dualCarl}\eqref{dualCarl2} are satisfied. We have $p\in(0,\infty)$ and $V\in\A_p(\R^{\vn})$ as in Proposition \ref{dualCarl}. Moreover, since $a<\min(1,p)$, it follows that $\frac{p}{a}>\max(1,p)$, thus $(\frac{p}{a})'<p'$, and then $r=a(\frac{p}{a})'<p'$ as well.
This confirms condition \eqref{r vs p} of Proposition \ref{dualCarl}\eqref{dualCarl2}.
Thus, Proposition \ref{dualCarl}\eqref{dualCarl2} applies to give the bound
\begin{align*}
&\bbNorm{\bigg\{
\sum_{S\in\D_{\vi-\vl-2\cdot\vone}(\R^{\vn})} \abs{[V]_{\aveL^p(S)}V^{-1}(\cdot)}\bigg[\fint_S g_{\vi}(x)dx\bigg]^{\frac1a} \one_S(\cdot) \bigg\}_{\vi\in\Z^k}}{
(L^{a(\frac{\vp}{a})'}\ell^{a(\frac{\vq}{a})'})_\pi}^a \\
&\quad\lesssim\bbNorm{\bigg\{
\sum_{S\in\D_{\vi-\vl-2\cdot\vone}(\R^{\vn})} \bigg[\fint_S g_{\vi}(x)dx\bigg]^{\frac1a} \one_S\bigg\}_{\vi\in\Z^k}}{
(L^{a(\frac{\vp}{a})'}\ell^{a(\frac{\vq}{a})'})_\pi}^a \\
&\quad=\bbNorm{\bigg\{
\sum_{S\in\D_{\vi-\vl-2\cdot\vone}(\R^{\vn})} \fint_S g_{\vi}(x)dx\, \one_S\bigg\}_{\vi\in\Z^k}}{
(L^{(\frac{\vp}{a})'}\ell^{(\frac{\vq}{a})'})_\pi}.
\end{align*}
We dualise this again with a suitable $\Norm{f}{(L^{\frac{\vp}{a}}\ell^{\frac{\vq}{a}})_\pi}=1$. Thus,
\begin{align*}
&\bbNorm{\bigg\{
\sum_{S\in\D_{\vi-\vl-2\cdot\vone}(\R^{\vn})} \fint_S g_{\vi}(x)dx\, \one_S\bigg\}_{\vi\in\Z^k}}{
(L^{(\frac{\vp}{a})'}\ell^{(\frac{\vq}{a})'})_\pi} \\
&\quad\sim\sum_{\vi\in\Z^k}\int_{\R^{\vn}}f_{\vi}(y)\bigg[\sum_{S\in\D_{\vi-\vl-2\cdot\vone}(\R^{\vn})} \fint_S g_{\vi}(x)dx\, \one_S(y)\bigg]dy \\
&\quad=\sum_{\vi\in\Z^k}\int_{\R^{\vn}}\bigg[\sum_{S\in\D_{\vi-\vl-2\cdot\vone}(\R^{\vn})} \fint_S f_{\vi}(y)dy\, \one_S(x)\bigg]g_{\vi}(x)dx \\
&\quad\leq\bbNorm{\bigg\{
\sum_{S\in\D_{\vi-\vl-2\cdot\vone}(\R^{\vn})} \fint_S f_{\vi}(y)dy\, \one_S\bigg\}_{\vi\in\Z^k} }{
(L^{\frac{\vp}{a}}\ell^{\frac{\vq}{a}})_\pi}
\Norm{g}{(L^{(\frac{\vp}{a})'}\ell^{(\frac{\vq}{a})'})_\pi}  \\
&\quad\leq\Norm{\{\mathcal M_{\vn}f_{\vi}\}_{\vi\in\Z^k}}{ (L^{\frac{\vp}{a}}\ell^{\frac{\vq}{a}})_\pi}\times 1
\lesssim \Norm{\{f_{\vi}\}_{\vi\in\Z^k}}{ (L^{\frac{\vp}{a}}\ell^{\frac{\vq}{a}})_\pi}=1,
\end{align*}
using the vector-valued maximal inequality, which is valid under the assumption that $\pi$ is an admissible permutation. (Note that this last step was essentially the same as the full argument in the unweighted case; in the matrix-weighted setting, it took some additional effort to get to this point!)
This finishes the proof of Lemma \ref{key lem 2}.
\end{proof}

\begin{remark}
The main working engine of the proof of Lemma \ref{key lem 2} in the matrix weighted cases consisted of the Carleson measure estimates of Corollary \ref{Carl Ap} and Proposition \ref{dualCarl}. Since the assumptions of these results allow an additional case besides the pure Besov and Triebel--Lizorkin spaces, one may ask why a version of this third case is not included in Lemma \ref{key lem 2}. The reason is the following: While the third case of both Corollary \ref{Carl Ap} and Proposition \ref{dualCarl} consists of the same condition $\vq\leq\vp=p\cdot\vone$, in the proof of Lemma \ref{key lem 2}, we apply these two results in two different spaces, essentially dual to each other. More precisely, we apply Corollary \ref{Carl Ap} with the original exponents $(\vp,\vq)$ and Proposition \ref{dualCarl} with exponents $(a(\frac{\vp}{a})',a(\frac{\vq}{a})')$.
Note that $a(\frac{\vq}{a})'\leq a(\frac{\vp}{a})'$ is equivalent to $\vp\leq\vq$, so if we have both $\vq\leq\vp$ and $a(\frac{\vq}{a})'\leq a(\frac{\vp}{a})'$, then indeed $\vq=\vp=p\cdot\vone$. Then $(L^{\vp}\ell^{\vq})_\pi=L^p\ell^p=\ell^pL^p$ by Fubini's theorem, and we would be in the scope of these pure cases anyway. Thus, there is nothing new that this third case of Corollary \ref{Carl Ap} and Proposition \ref{dualCarl} would give us in Lemma \ref{key lem 2}.
\end{remark}

From the two Lemmas \ref{key lem} and \ref{key lem 2}, we immediately obtain the boundedness of almost diagonal operators on the spaces $\asn(V)$ with $\tau=0$.

\begin{theorem}\label{AD s0}
Let $\vp\in(0,\infty)^k$, $\vq\in(0,\infty]^k$, $\vs\in\R^k$, $\pi\in S_{[2k]}$ be admissible for $(\vp,\vq)$, and $V$ be a matrix weight such that $(\vp,\vq,\pi,V)$ belongs to one of the main cases (Definition \ref{main cases}).
Let $B$ be $(\vD,\vE,\vF)$-almost diagonal with
\begin{equation}\label{DEF1}
\vD>\frac{\vn}{r},\quad
\vE>\frac{\vn}{2}+\vs,\quad
\vF>\frac{\vn}{r}-\frac{\vn}{2}-\vs,\qquad r:=r(\pi,\vp,\vq)\wedge 1
\end{equation}
with $r(\pi,\vp,\vq)$ as in \eqref{rpipq-max}.
If $t\in\asn(V)$, then $Bt$ is well defined and
\begin{equation*}
\Norm{Bt}{\asn(V)}
\lesssim\Norm{t}{\asn(V)}
\end{equation*}
with the implicit positive constant independent of $t$.
\end{theorem}

\begin{proof}
By Remark \ref{s vs 0}, it suffices to consider $\vs=\vnull$. If $(\vD,\vE,\vF)$ satisfies \eqref{DEF1}, then we can choose $a\in(0,r)$ close enough to $r$ so that $(\vD,\vE,\vF)$ still satisfies
\begin{equation}\label{DEF3}
\vD>\frac{\vn}{a},\quad
\vE>\frac{\vn}{2},\quad
\vF>\frac{\vn}{a}-\frac{\vn}{2}.
\end{equation}

Let us denote by $\mathfrak{R}$ the right-hand side of \eqref{key eq} with these parameters and $H_{\vj}\equiv V$. From the convergence of geometric series in \eqref{key eq} with parameters \eqref{DEF3}, it is immediate that
\begin{equation*}
\mathfrak{R}\lesssim\sup_{\vh\in\Z^k,\vl\in\N^k}
\bbNorm{\bigg\{\bigg[\fint_{C(\cdot 2^{\vl+\vh_+-\vi})}\abs{V(\cdot)t_{\vi}(y)}^a dy\bigg]^{\frac1a}\bigg\}_{\vi\in\Z^k}}{(L^{\vp}\ell^{\vq})_\pi}
=:\mathfrak{S}.
\end{equation*}
Now Lemma \ref{key lem 2} says that
\begin{equation*}
\mathfrak{S} \lesssim \Norm{Vt}{(L^{\vp}L^{\vq})_\pi}=\Norm{t}{\ann(V)}.
\end{equation*}

Recall that Lemma \ref{key lem} works under the condition that $\mathfrak{R}<\infty$. Since we have just verified this condition, Lemma \ref{key lem} implies that $Bt$ is well defined and satisfies
\begin{equation*}
\Norm{Bt}{\ann(V)}
=\Norm{V(Bt)}{(L^{\vp}\ell^{\vq})_\pi}
\lesssim\mathfrak{R}.
\end{equation*}
A combination of the last three estimates shows that
\begin{equation*}
\Norm{Bt}{\ann(V)}
\lesssim\mathfrak{R}
  \lesssim\mathfrak{S}
  \lesssim\Norm{t}{\ann(V)},
\end{equation*}
which concludes the proof of the proposition.
\end{proof}

\begin{remark}
Theorem \ref{AD s0} reproduces \cite[Theorem 5.2]{GKP21} in the unweighted case $V\equiv 1$.
Moreover, it extends \cite[Theorem 5.2]{GKP21} to the matrix-$\A_p$-weighted situation
with exactly the same assumptions on the almost diagonal parameters as in the unweighted case.
Even in the one-parameter case, the fact that matrix-$\A_p$-weighed estimates
can be obtained with the same almost diagonal assumptions as in
the unweighted case was only recently observed in \cite{BHYY:1b},
while earlier results in the matrix-weighted Besov and Triebel--Lizorkin spaces \cite{FR04,FR21,Roud03}
required almost diagonal assumptions with additional dependence on the properties of the weight.
In the multi-parameter context, the matrix-weighted estimates of Theorem \ref{AD s0} are completely new.
\end{remark}

\subsection{Case $\tau\in(0,\infty)$} \label{tau approach}

The extension of Theorem \ref{AD s0} to $\atau(V)$ with $\tau\in(0,\infty)$ requires some additional work. At the present, our understanding of this situation is not as complete as in the one-parameter theory \cite{BHYY:1b}. We state and prove the following result, noting however that there should be a better version yet to be discovered; see Remark \ref{rem AD st} below.

\begin{theorem}\label{AD st}
Let $\vp\in(0,\infty)^k$, $\vq\in(0,\infty]^k$, $\vs\in\R^k$,
$\tau\in[0,\infty)$, $\pi\in S_{[2k]}$ be admissible for $(\vp,\vq)$, and $V$ be a matrix weight such that $(\vp,\vq,\pi,V)$ belongs to one of the main cases (Definition \ref{main cases}).
Let $B$ be $(\vD,\vE,\vF)$-almost diagonal with
\begin{equation}\label{DEF2}
\vD>\frac{\vn}{r}+\vn\tau,\quad
\vE>\frac{\vn}{2}+\vn\tau+\vs,\quad
\vF>\frac{\vn}{r}-\frac{\vn}{2}-\vs,
\end{equation}
where $r:=r(\pi,\vp,\vq)\wedge 1$.
If $t\in\atau(V)$, then $Bt$ is well defined, and
\begin{equation*}
\Norm{Bt}{\atau(V)}
\lesssim\Norm{t}{\atau(V)},
\end{equation*}
where the implicit positive constant is independent of $t$.
\end{theorem}

\begin{proof}
By Remark \ref{s vs 0}, it suffices to consider $\vs=\vnull$.
As in the proof of Theorem \ref{AD s0}, the idea is to estimate
$\Norm{Bt}{\asn(V)}$ by Lemma \ref{key lem}, and then use
Lemma \ref{key lem 2} to estimate each of the terms on the right-hand side
of \eqref{key eq} of Lemma \ref{key lem}. The main difference
compared to Theorem \ref{AD s0} is that the quasi-norm of $\atau(V)$
for $\tau>0$ is no longer just the $(L^{\vp}\ell^{\vq})_\pi$ quasi-norm,
but involves an additional supremum over the choice of a truncating set $\Omega$.

To obtain the relevant quasi-norms of $Bt$ on the left-hand side of \eqref{key eq}, for a given $\Omega\in\Open(\R^{\vn})$, we now choose $H_{\vj}:=\one_{\Omega_{\vj}}V$. For a fixed $\Omega$, these $H_{\vj}$ are no longer a.e.\ invertible; however, for a.e.\ $x\in\R^{\vn}$ and every $\vj\in\Z^k$, we can choose $\Omega$ so that $H_{\vj}(x)$ is invertible, and this is enough for the qualitative conclusion that $Bt$ is well defined, as one can check by inspection of \cite[Remark 4.5]{BHYY:1b}. On the right-hand side of \eqref{key eq}, we will then have expressions of the form
\begin{equation}\label{withOmega}
\bbNorm{\bigg\{\one_{\Omega_{\vi-\vh}}(\cdot)\bigg[\fint_{C(\cdot,2^{\vl+\vh_+-\vi})}
\abs{V(\cdot)t_{\vi}(y)}^a dy\bigg]^{\frac1a}\bigg\}_{\vi\in\Z^k}}{
(L^{\vp}\ell^{\vq})_\pi}.
\end{equation}
We also wish to insert an indicator as a multiplier of $t_{\vi}(y)$. If $x=(x_\nu)_{\nu=1}^k\in\Omega_{\vi-\vh}$, it means by definition that $x\in R\subset\Omega$ for some $R=\prod_{\nu=1}^k R_\nu\in\D_{\vi-\vh}(\R^{\vn})$. Then suppose that
\begin{align*}
y\in C(x,2^{\vl+\vh_+-\vi})
&=C(x,2^{(\vl+\vh_-)+(\vh-\vi)})
=\prod_{\nu=1}^k B(x_{\nu},2^{l_\nu+(h_\nu)_-}\ell(R_\nu)) \\
&\subset\prod_{\nu=1}^k 3\cdot R_{\nu}^{(l_\nu+(h_\nu)_-)}
= 3\cdot R^{(\vl+\vh_-)},
\end{align*}
where $R_{\nu}^{(l_\nu+(h_\nu)_-)}$ is the $(l_\nu+(h_\nu)_-)$ generations older dyadic ancestor of $R_{\nu}$.

Note that both $R\subset\Omega$ and $R\subset 3\cdot R^{(\vl+\vh_-)}$, and hence
\begin{equation*}
\mathcal M_{\vn}\one_{\Omega}(y)
\gtrsim \fint_{3\cdot R^{(\vl+\vh_-)}}\one_{\Omega}(z)dz
=\frac{\abs{(3\cdot R^{(\vl+\vh_-)})\cap\Omega} }{ \abs{ 3\cdot R^{(\vl+\vh_-)} } }
\geq \frac{\abs{R}}{  \abs{ 3\cdot R^{(\vl+\vh_-)} } }
\sim 2^{-(\vl+\vh_-)\cdot\vn}.
\end{equation*}
Thus,
\begin{equation*}
y\in 3\cdot R^{(\vl+\vh_-)}
\subset
\widetilde\Omega
:=\{y\in\R^{\vn}:\
\mathcal M_{\vn}\one_{\Omega} (y)> c2^{-(\vl+\vh_-)\cdot\vn}\}
\end{equation*}
for some dimensional constant $c>0$.

Recalling that $R\in\D_{\vi-\vh}(\R^{\vn})$, it follows that $R^{(\vl+\vh_-)}\in \D_{(\vi-\vh)-(\vl+\vh_-)}(\R^{\vn})=\D_{\vi-\vl-\vh_+}(\R^{\vn})$, and the cube $3\cdot R^{(\vl+\vh_-)}$ is a union of $3^n$ dyadic rectangles $S\in\D_{\vi-\vl-\vh_+}(\R^{\vn})$, namely $S=R^{(\vl+\vh_+)}$ and all its neighbours. Since these $S\in\D_{\vi-\vl-\vh_+}(\R^{\vn})$ satisfy $S\subset\widetilde\Omega$, they also satisfy $S\subset\widetilde\Omega_{\vi-\vl-\vh_+}$; since $y$ belongs to one of these $S$,
it follows that $y\in\widetilde\Omega_{\vi-\vl-\vh_+}\subset\widetilde\Omega_{\vi}$, since $\vl,\vh_+\geq\vnull$.

Thus, we may insert the indicator $\one_{\widetilde\Omega_{\vi}}(y)$ in the integral in \eqref{withOmega} without changing its value. Once  we have done this, we estimate up by dropping the indicator $\one_{\Omega_{\vi-\vh}}(\cdot)$ outside the integral. It follows that
\begin{align*}
&\bbNorm{\bigg\{\one_{\Omega_{\vi-\vh}}(\cdot)\bigg[\fint_{C(\cdot,2^{\vl+\vh_+-\vi})}
\abs{V(\cdot)t_{\vi}(y)}^a dy\bigg]^{\frac1a}\bigg\}_{\vi\in\Z^k}}{
(L^{\vp}\ell^{\vq})_\pi} \\
&\quad\leq\bbNorm{\bigg\{\bigg[\fint_{C(\cdot,2^{\vl+\vh_+-\vi})}\abs{V(\cdot)
\one_{\widetilde\Omega_{\vi}}(y)t_{\vi}(y)}^a dy\bigg]^{\frac1a}\bigg\}_{\vi\in\Z^k}}{
(L^{\vp}\ell^{\vq})_\pi} \\
&\quad\lesssim\BNorm{\Big\{V\one_{\widetilde\Omega_{\vi}}t_{\vi}\Big\}_{\vi\in\Z^k}}
{(L^{\vp}\ell^{\vq})_\pi}
\leq\abs{\widetilde\Omega}^{\tau}\Norm{t}{\antau(V)},
\end{align*}
where the last two steps were applications of Lemma \ref{key lem 2} with the sequence $\{\one_{\widetilde\Omega_{\vi}}t_{\vi}\}_{\vi\in\Z^k}$ in place of $\{t_{\vi}\}_{\vi\in\Z^k}$ and the definition of the quasi-norm of $\antau(V)$.

The strong maximal inequality shows that
\begin{align*}
\abs{\widetilde\Omega}
&=|\{y\in\R^{\vn}:\
\mathcal M_{\vn}\one_{\Omega} (y)> c2^{-(\vl+\vh_-)\cdot\vn}\}|
\lesssim 2^{(\vl+\vh_-)\cdot\vn}[1+\log^{k-1}2^{(\vl+\vh_-)\cdot\vn}]\abs{\Omega} \\
&\sim 2^{(\vl+\vh_-)\cdot\vn}[1+(\vl+\vh_-)\cdot\vn]^{k-1}\abs{\Omega}.
\end{align*}

Substituting the above estimates into \eqref{key eq} of Lemma \ref{key lem} with $H_{\vj}=\one_{\Omega_{\vj}}V$, it follows that
\begin{align*}
& \Norm{\{\one_{\Omega_{\vj}}V(Bt)_{\vj}\}_{\vj\in\Z^k}}{(L^{\vp}\ell^{\vq})_\pi}^u \\
&\quad \lesssim \sum_{\vh\in\Z^k}\sum_{\vl\in\N^k}
\Big[ 2^{-(\vE-\frac{\vn}{2})\cdot\vh_-}2^{-(\vF+\frac{\vn}{2}-\frac{\vn}{a})\cdot\vh_+}
2^{-(\vD-\frac{\vn}{a})\cdot\vl}
2^{(\vl+\vh_-)\cdot\vn\tau}P(\vl,\vh)
\abs{\Omega}^{\tau}\Norm{t}{\antau}  \Big]^u,
\end{align*}
where $P(\vl,\vh):=(1+(\vl+\vh_-)\cdot\vn)^{(k-1)\tau}$. It is immediate that the series converges if and only if
\begin{equation}\label{DEF tau a}
\vD>\frac{\vn}{a}+\vn\tau,\quad
\vE>\frac{\vn}{2}+\vn\tau,\quad
\vF>\frac{\vn}{a}-\frac{\vn}{2};
\end{equation}
the polynomial factor $P(\vl,\vh)$ has no impact on this. Since Lemmas \ref{key lem} and \ref{key lem 2} allow us to choose any $a\in(0,r)$, where $r:=r(\pi,\vp,\vq)\wedge 1$, inequalities \eqref{DEF tau a} will hold under the assumed \eqref{DEF2} (with $\vs=\vnull$), as soon as $a$ is chosen close enough to $r$. Thus, we obtain
\begin{equation*}
\Norm{\{\one_{\Omega_{\vj}}V(Bt)_{\vj}\}_{\vj\in\Z^k}}{(L^{\vp}\ell^{\vq})_\pi}
\lesssim\abs{\Omega}^{\tau}\Norm{t}{\antau(V)}.
\end{equation*}
Dividing by $\abs{\Omega}^{\tau}$ and taking the supremum over $\Omega\in\Open(\R^{\vn})$ prove the claimed result for $\vs=\vnull$. This suffices by Remark \ref{s vs 0}, as we already pointed out in the beginning of the proof.
\end{proof}

\begin{remark}\label{rem AD st}
Theorem \ref{AD st} clearly extends Theorem \ref{AD s0}, which is recovered as the special case $\tau=0$. Theorem \ref{AD st} is the first result on the boundedness of almost diagonal operators on multi-parameter Besov or Triebel--Lizorkin-type spaces with $\tau\in(0,\infty)$, even in the unweighted bi-parameter case, where the only result prior to this work, \cite[Theorem 5.2]{GKP21}, dealt with $\tau=0$ only.

However, Theorem \ref{AD st} is not a complete extension of the one-parameter results of \cite{YY10} (unweighted) and \cite{BHYY:1b} (matrix-weighted), since the almost diagonal conditions imposed by those works are weaker than what one obtains by specialising \eqref{DEF2} to the one-parameter case. More precisely, \cite[Theorem 5.1]{BHYY:1b} assumes a condition like \eqref{DEF2}, but with $\tau$ replaced by a modified parameter
\begin{equation*}
\widehat\tau:=\Big(\tau-\frac{1}{p}\Big[1-\frac{d}{n}\Big]\Big)_+,
\end{equation*}
when $V$ has $\A_p$-dimension $d$. In particular, in the unweighed case $V\equiv 1$, the $\A_p$-dimension is $d=0$, and hence $\widehat\tau=0$ for all $\tau\in[0,\frac{1}{p}]$; this agrees with the assumptions of \cite[Theorem 4.1]{YY10} in this parameter range. At the moment, it remains open how to construct a multi-parameter argument that reproduces this improvement in the one-parameter case (and hopefully provides a natural extension to several parameters, taking advantage of the available theory of $\A_p$-dimensions in this setting as developed in Subsection \ref{Ap dim}).

On the other hand, the said improvements of \cite[Theorem 5.1]{BHYY:1b} over the present Theorem \ref{AD st} require the additional information about the $\A_p$-dimension of the weight, while the assumptions of Theorem \ref{AD st} are uniform over all $V\in\A_p(\R^{\vn})$. From this perspective, these assumptions are not so unreasonable, as we will show in Corollary \ref{AD st sharp2} below.
\end{remark}

Next, we want to study the composition of two almost diagonal operators.

\begin{lemma}\label{AcB}
Let $\vec D,\vec E,\vec F,\vec D^*,\vec E^*,\vec F^*\in\R^k$
satisfy that, for every $i\in[k]$,
$$
D_i,D^*_i>n_i,\
E_i+F^*_i>\min(D_i, D^*_i),\
F_i+E^*_i>n_i,\
E_i\neq E^*_i,\
F_i\neq F^*_i.
$$
If $A:=\{a_{Q,R}\}_{Q,R\in\mathscr{D}(\mathbb R^{\vec n})}$
is $(\vec D,\vec E,\vec F)$-almost diagonal
and $B:=\{b_{R,P}\}_{R,P\in\mathscr{D}(\mathbb R^{\vec n})}$ is
$(\vec D^*,\vec E^*,\vec F^*)$-almost diagonal, then
\begin{align*}
A\circ B:=\Bigg\{\sum_{R\in\mathscr{D}(\mathbb R^{\vec n})}
a_{Q,R}b_{R,P}\Bigg\}_{Q,P\in\mathscr{D}(\mathbb R^{\vec n})}
\end{align*}
is $(\min(\vec D,\vec D^*),\min(\vec E,\vec E^*),\min(\vec F,\vec F^*))$-almost diagonal.
\end{lemma}

\begin{proof}
By the definition of $(\vec D,\vec E,\vec F)$-almost diagonal, we obtain,
for every $Q,P\in\mathscr{D}(\mathbb R^{\vec n})$,
\begin{align*}
|(A\circ B)_{Q,P}|
\lesssim \sum_{R\in\mathscr{D}(\mathbb R^{\vec n})}
\prod_{i=1}^k b^{D_i,E_i,F_i}_{i;Q_i,R_i}
b^{D_i^*,E_i^*,F_i^*}_{i;R_i,P_i}
= \prod_{i=1}^k \Bigg( \sum_{R_i\in\mathscr{D}(\mathbb R^{n_i})}
b^{D_i,E_i,F_i}_{i;Q_i,R_i}
b^{D_i^*,E_i^*,F_i^*}_{i;R_i,P_i} \Bigg).
\end{align*}
Therefore, without loss of generality, we may assume $k=1$.
For every $Q,P\in\mathscr{D}(\mathbb R^n)$,
\begin{align*}
\sum_{R\in\mathscr{D}(\mathbb R^n)}
b^{D,E,F}_{Q,R} b^{D^*,E^*,F^*}_{R,P}
&=\sum_{j=-\infty}^{\min(j_Q,j_P)}
\sum_{R\in\mathscr{D}_j(\mathbb R^n)} \cdots
+\sum_{j=\min(j_Q,j_P)+1}^{\max(j_Q,j_P)}
\sum_{R\in\mathscr{D}_j(\mathbb R^n)}\cdots \\
&\quad +\sum_{j=\max(j_Q,j_P)+1}^{\infty}
\sum_{R\in\mathscr{D}_j(\mathbb R^n)} \cdots\\
&=:\mathrm{I}+\mathrm{II}+\mathrm{III}.
\end{align*}
From \cite[Appendix B.1]{Graf-mod}
(see also \cite[Lemma 3.6]{BHYY:2b}), we infer that
\begin{align*}
\mathrm{I}
&= \sum_{j=-\infty}^{\min(j_Q,j_P)}
\sum_{R\in\mathscr{D}_j(\mathbb R^n)} \bigg[1+\frac{\abs{c_Q-c_R}}{\ell(Q)\vee\ell(R)}\bigg]^{-D}
\bigg[\frac{\ell(Q)}{\ell(R)}\bigg]^E
\bigg[1+\frac{\abs{c_R-c_P}}{\ell(R)\vee\ell(P)}\bigg]^{-D^*}
\bigg[\frac{\ell(P)}{\ell(R)}\bigg]^{F^*}\\
&\lesssim \sum_{j=-\infty}^{\min(j_Q,j_P)} \Big(1+2^j\abs{c_Q-c_P}\Big)^{-\min(D,D^*)}
2^{(j-j_Q)E} 2^{(j-j_P)F^*} \\
&\lesssim \sum_{j=-\infty}^{\min(j_Q,j_P)} \bigg[1+\frac{\abs{c_Q-c_P}}{\ell(Q)\vee\ell(P)}\bigg]^{-\min(D,D^*)}
2^{-[j-\min(j_Q,j_P)]\min(D,D^*)} 2^{(j-j_Q)E} 2^{(j-j_P)F^*} \\
&\lesssim \bigg[1+\frac{\abs{c_Q-c_P}}{\ell(Q)\vee\ell(P)}\bigg]^{-\min(D,D^*)}
2^{\min(j_Q,j_P)(E+F^*)-j_QE-j_PF^*}
= b^{\min(D,D^*),E,F^*}_{Q,P}.
\end{align*}
Similarly, we obtain $\mathrm{II}\lesssim b^{\min(D,D^*),\min(E,E^*),\min(F,F^*)}_{Q,P}$
and $\mathrm{III}\lesssim b^{\min(D,D^*),E^*,F}_{Q,P}$.
This finishes the proof of Lemma \ref{AcB}.
\end{proof}

\begin{definition}\label{astAD}
Let $\vp\in(0,\infty)^k$, $\vq\in(0,\infty]^k$, $\vs\in\R^k$,
$\tau\in[0,\infty)$, and $\pi\in S_{[2k]}$ be admissible for $(\vp,\vq)$.
An infinite matrix $\{b_{Q,R}\}_{Q,R\in\mathscr D(\mathbb R^{\vec n})}$
is said to be $\dot a^{\vec s,\tau}_{\vec p,\vec q,\pi}$-almost diagonal if it is $(\vec D,\vec E,\vec F)$-almost diagonal with
$\vec D,\vec E,\vec F$ satisfying \eqref{DEF2}.
\end{definition}

\begin{proposition}\label{AcB2}
Let $\vp\in(0,\infty)^k$, $\vq\in(0,\infty]^k$, $\vs\in\R^k$,
$\tau\in[0,\infty)$, and $\pi\in S_{[2k]}$ be admissible for $(\vp,\vq)$.
If $A^{(i)}:=\{a^{(i)}_{Q,R}\}_{Q,R\in\mathscr{D}(\mathbb R^{\vec n})}$
is $\dot a^{\vec s,\tau}_{\vec p,\vec q,\pi}$-almost diagonal for both $i=1,2$, then
$A^{(1)}\circ A^{(2)}$ is also $\dot a^{\vec s,\tau}_{\vec p,\vec q,\pi}$-almost diagonal.
\end{proposition}

\begin{proof}
By the definition of the $\dot a^{\vec s,\tau}_{\vec p,\vec q,\pi}$-almost diagonal,
it follows that
for both $i=1,2$, the matrix $A^{(i)}$ is $(\vec D^{(i)},\vec E^{(i)},\vec F^{(i)})$-almost diagonal with some
\begin{equation}\label{lower bounds}
\vec D^{(i)}>\frac{\vn}{r}+\vn\tau,\quad
\vE^{(i)}>\frac{\vn}{2}+\vn\tau+\vs,\quad
\vF^{(i)}>\frac{\vn}{r}-\frac{\vn}{2}-\vs.
\end{equation}

Note that
\begin{equation}\label{ob}
\Big( \frac{\vn}{2}+\vn\tau+\vs\Big)
+\Big( \frac{\vn}{r}-\frac{\vn}{2}-\vs\Big)
=\vn\tau+\frac{\vn}{r}
\geq \vn > \vnull.
\end{equation}
Making some components slightly smaller, if necessary, while still maintaining the lower bounds in \eqref{lower bounds}, we can assume without loss of generality that
$$\vec E^{(1)}+\vec F^{(2)}>\min(\vec D^{(1)},\vec D^{(2)})$$
and, for every $\nu\in[k]$,
$(E^{(1)})_\nu \neq (E^{(2)})_\nu$ and $(F^{(1)})_\nu \neq (F^{(2)})_\nu$.
Inequality \eqref{ob} also shows in particular that
$$
\vec D^{(1)},\vec D^{(2)}>\vec n,\quad
\vec F^{(1)}+\vec E^{(2)}>\vec n.
$$

These, together with Lemma \ref{AcB}, further imply that
$A^{(1)}\circ A^{(2)}$ is $(\vD,\vE,\vF)$-almost diagonal,
where
\begin{equation*}
\begin{split}
\vD =\min_{i=1,2}\vec D^{(i)}>\frac{\vn}{r}+\vn\tau,\quad
\vE=\min_{i=1,2}\vec E^{(i)}>\frac{\vn}{2}+\vn\tau+\vs,\quad
\vF =\min_{i=1,2}\vec F^{(i)}>\frac{\vn}{r}-\frac{\vn}{2}-\vs.
\end{split}
\end{equation*}
Thus, $A^{(1)}\circ A^{(2)}$ is also $\dot a^{\vec s,\tau}_{\vec p,\vec q,\pi}$-almost diagonal.
This finishes the proof of Proposition \ref{AcB2}.
\end{proof}

\begin{remark}
In \cite[Theorem 9.1(i)]{FJ90}, Frazier and Jawerth proved that
the composition of two $\dot f^s_{p,q}$-almost diagonal operators
is still an $\dot f^s_{p,q}$-almost diagonal operator.
In the case where $k=1$, $\tau=0$, and $\pi=F$, we have
$\dot a^{\vec s,\tau}_{\vec p,\vec q,\pi}=\dot f^s_{p,q}$,
and thus Proposition \ref{AcB2} in this case coincides with \cite[Thoerem 9.1(i)]{FJ90}.
\end{remark}

We can also obtain the following variant of Theorem \ref{AD st}. Note that the almost diagonality assumption \eqref{DEFvar} below is stronger than the corresponding assumption \eqref{DEF2} of Theorem \ref{AD st} (indeed, the lower bounds for each $\vD,\vE,\vF$ have been increased by a new term); the advantage of Corollary \ref{AD variant} is that we can trade the assumption that $(\vp,\vq,\pi,V)$ belongs to one of the main cases (Definition \ref{main cases}) against the condition \eqref{key cases}; thus, besides pure Besov or Triebel--Lizorkin spaces, we can also allow matrix-weighted space with the $\vq\leq\vp$.

\begin{corollary} \label{AD variant}
Let $\vs\in\R^k$, $\tau\in[0,\infty)$, $\vp=p\cdot\vone$ with $p\in(0,\infty)$,
$\vq\in(0,\infty]^k$, $\pi\in S_{[2k]}$ be admissible for $(\vp,\vq)$,
and $r:=r(\pi,\vp,\vq)\wedge 1$.
Suppose that $V\in\A_p(\R^{\vn})$ has $\A_p$-dimensions $(\vd,\ve,\vDe)$,
and let $B$ be $(\vD,\vE,\vF)$-almost diagonal with
\begin{equation}\label{DEFvar}
\vD>\frac{\vn}{r}+\vn\tau + \vec\Delta,\quad
\vE>\frac{\vn}{2}+\vs+\vn\tau + \frac{\vec d}{p},\quad
\vF>\frac{\vn}{r}-\frac{\vn}{2}-\vs + \frac{\vec e}{p}.
\end{equation}
Then $Bt$ is well defined for every $t\in\atau([V]_p)$, and
\begin{equation*}
\Norm{Bt}{\atau([V]_p)}
\lesssim\Norm{t}{\atau([V]_p)}
\end{equation*}
with the implicit positive constant independent of $t$. If, in addition, we have \eqref{key cases}, then the same holds with $\atau(V)$ in place of $\atau([V]_p)$.
\end{corollary}

\begin{proof}
Recalling that the unweighted case is one of the main cases covered by Theorem \ref{AD st}, the said theorem implies that all $(\vD',\vE',\vF')$-almost diagonal matrices with
\begin{equation}\label{DEFvar2}
  \vD'>\frac{\vn}{r}+\vn\tau,\quad
\vE'>\frac{\vn}{2}+\vs+\vn\tau,\quad
\vF'>\frac{\vn}{r}-\frac{\vn}{2}-\vs
\end{equation}
define bounded operators on the unweighed $\atau(\R^{\vn})$.

Feeding this information into Corollary \ref{cor AD simple}, it follows that all $(\vD,\vE,\vF)=(\vD'+ \vec\Delta,\vE'+ {\vec d}/{p},\vF'+ {\vec e}/{p})$-almost diagonal matrices define bounded operators on $\atau([V]_p)$. When $(\vD',\vE',\vF')$ takes arbitrary values as in \eqref{DEFvar2}, it is clear that $(\vD,\vE,\vF)$ can take arbitrary values as in \eqref{DEFvar}. This proves the first claim concerning boundedness on $\atau([V]_p)$.

Finally, adding assumption \eqref{key cases}, it follows from Corollary \ref{3 norms cor} that $\atau([V]_p)=\atau(V)$ with equivalent norms, and this completes the proof.
\end{proof}

\subsection{Sharpness}\label{sharpness}

The sharpness of the one-parameter version of Theorem \ref{AD s0} was studied in detail in \cite{BHYY:1b}. The results of these considerations can be transferred to the multi-parameter setting without the need to revisit the details of the original counterexamples, for which the reader may consult \cite{BHYY:1b}. In order to lift one-parameter examples to several parameters, we will make use of the following:

\begin{lemma}\label{t tensor}
For sequences $t^{(\nu)}=\{t^{(\nu)}_Q\}_{Q\in\D(\R^{n_\nu})}\in\F^{\D(\R^{n_\nu})}$ with $\nu\in[k]$, define the tensor sequence
\begin{equation*}
t^{\otimes}:=(t^{\otimes}_R)_{R\in\D(\R^{\vn})},\quad
t^{\otimes}_R:=\prod_{\nu=1}^k t^{(\nu)}_{R_\nu}.
\end{equation*}
Given a permutation $\pi\in S_{[2k]}$, define $\pi_\nu\in S_{[2]}$ for each $\nu\in[k]$, by
\begin{equation*}
\pi_\nu:=\begin{cases}F:=\id, & \text{if}\quad \pi^{-1}(\nu)<\pi^{-1}(\nu+k), \\
B:=(1\ 2), & \text{if}\quad \pi^{-1}(\nu)>\pi^{-1}(\nu+k).\end{cases}
\end{equation*}
Let $\vs\in\R^k$, $\vp\in(0,\infty)^k$, $\vq\in(0,\infty]^k$, $\pi\in S_{[2k]}$, and
\begin{equation*}
V(x)=\prod_{\nu=1}^k V(x_\nu)
\end{equation*}
be a scalar-valued product weight, and suppose moreover that one of the following additional conditions is satisfied:
\begin{enumerate}[\rm(i)]
\item\label{tau=0} $\tau=0$, or
\item\label{crit} $\vp=\vq=p\cdot\vone$ for some $p\in(0,\infty)$, and $\tau=\frac1p$.
\end{enumerate}
Then
\begin{equation*}
\Norm{t^\otimes}{\atau(V)}=\prod_{\nu=1}^k \Norm{t^{(\nu)}}{ \dot a^{s_\nu,\tau}_{p_\nu,q_\nu,\pi_\nu}(V_\nu) }.
\end{equation*}
\end{lemma}

\begin{proof}
\eqref{tau=0}: Let first $\tau=0$. This simplifies the quasi-norm of $\atau=\asn$ since the supremum over $\Omega\in\Open(\R^{\vn})$ is avoided, and
\begin{equation}\label{tensor-norm}
\Norm{t^{\otimes}}{\asn(V)}
=\Norm{\{2^{\vi\cdot\vs}Vt_{\vi}^{\otimes}\}_{\vi\in\Z^k}}{(L^{\vp}\ell^{\vq})_\pi}
=\BNorm{\Big\{x\mapsto\prod_{\nu=1}^k 2^{i_\nu s_\nu}V_{\nu}(x_\nu)t^{(\nu)}_{i_\nu}(x_\nu)\Big\}_{\vi\in\Z^k}}{(L^{\vp}\ell^{\vq})_\pi},
\end{equation}
where the $\nu$-th component of the product depends only on the variables $i_\nu\in\Z$ and $x_\nu\in\R^{n_\nu}$. For any product of functions of different variables on $\sigma$-finite measure spaces, we have the Fubini identities
\begin{align*}
&\Norm{(y,z)\mapsto f(y)g(z)}{L^u(d\lambda(y);L^v(d\omega(z)))} \\
&\quad=\Norm{f}{L^u(\lambda)}\Norm{g}{L^v(\omega)}
=\Norm{(y,z)\mapsto f(y)g(z)}{L^v(d\omega(z);L^u(d\lambda(y)))}.
\end{align*}
Hence, in \eqref{tensor-norm}, we may permute the order of any consecutive $L^{p_\nu}$ or $\ell^{q_\nu}$ and $L^{p_\mu}$ or $\ell^{q_\mu}$ quasi-norms, except the combination of $L^{p_\nu}$ and $\ell^{q_\nu}$ with the same $\nu\in[k]$. By a repeated application of such permutations, we find that
\begin{equation*}
\begin{split}
\Bigg\|\Bigg\{x\mapsto\prod_{\nu=1}^k 2^{i_\nu s_\nu}V_\nu(x_\nu)t^{(\nu)}_{i_\nu}(x_\nu)\Bigg\}_{\vi\in\Z^k}\Bigg\|_{(L^{\vp}\ell^{\vq})_\pi}
&=\prod_{\nu=1}^k  \Norm{ \{2^{i_\nu s_\nu } V_\nu t_{i_\nu}^{(\nu)} \}_{i_\nu\in\Z} }{ (L^{p_\nu}\ell^{q_\nu})_{\pi_\nu}}  \\
&=\prod_{\nu=1}^k \Norm{ t^{(\nu)}}{ \dot a^{s_\nu,0}_{p_\nu,q_\nu,\pi_\nu} (V_\nu)},
\end{split}
\end{equation*}
where $\pi_\nu\in S_{[2]}$ simply encodes the fact that $L^{p_\nu}$ and $\ell^{q_\nu}$ retain their original order. Together with \eqref{tensor-norm}, this completes the proof of case \eqref{tau=0}.

\eqref{crit}: We now have to deal with the supremum over $\Omega\in\Open(\R^{\vn})$ in
\begin{equation*}
\Norm{t}{\atau(V)}=\sup_{\Omega\in\Open(\R^{\vn})} \abs{\Omega}^{-\tau}\Norm{\{2^{\vi\cdot\vs}\one_{\Omega_{\vi}}V t_{\vi}\}_{\vi\in\Z^k}}{(L^{\vp}\ell^{\vq})_\pi},
\end{equation*}
but the assumption that $\vp=\vq=p\cdot\vone$ allows us to simplify $(L^{\vp}\ell^{\vq})_\pi=\ell^p L^p$. Hence
\begin{align*}
\Norm{\{2^{\vi\cdot\vs}\one_{\Omega_{\vi}}Vt_{\vi}\}_{\vi\in\Z^k}}{(L^{\vp}\ell^{\vq})_\pi}^p
&=\sum_{\vi\in\Z^k}\int_{\R^{\vn}}\abs{2^{\vi\cdot\vs}\one_{\Omega_{\vi}}(x)V(x)t_{\vi}(x)\}}^p dx \\
&=\sum_{\vi\in\Z^k}\sum_{R\in\D_{\vi}(\Omega)} 2^{\vi\cdot\vs p}\int_R V^p(x)dx \frac{\abs{t_R}^p }{ \abs{R}^{\frac{p}{2}} }
=\sum_{R\in\D(\Omega)} \vell(R)^{-\vs}V^p(R)\frac{\abs{t_R}^p }{ \abs{R}^{\frac{p}{2}} },
\end{align*}
where $V^p(R):=\int_R V^p(x)dx$,
\begin{equation*}
\vell(R):=(\ell(R_1),\ldots,\ell(R_k)),\qquad
\vell(R)^{-\vs}:=\prod_{\nu=1}^k \ell(R_\nu)^{-s_\nu}.
\end{equation*}
In particular, for $t=t^{\otimes}$, we obtain
\begin{equation*}
\Norm{\{2^{\vi\cdot\vs}\one_{\Omega_{\vi}}Vt^{\otimes}_{\vi}\}_{\vi\in\Z^k}}{(L^{\vp}\ell^{\vq})_\pi}^p
=\sum_{R\in\D(\Omega)} u^{\otimes}_R,
\end{equation*}
where
\begin{equation*}
u^{\otimes}_R:=\prod_{\nu=1}^k u_{R_\nu},\quad
u_{R_\nu}^{(\nu)}:=   \ell(R_\nu)^{-s_\nu}V_\nu^p(R_\nu)\abs{R_\nu}^{-\frac{p}{2}}\abs{t^{(\nu)}_{R_\nu}}^p.
\end{equation*}
Specialising to $\Omega=\prod_{\nu=1}^k\Omega_\nu$ with $\Omega_\nu\in\Open(\R^{n_\nu})$, we note that
\begin{equation*}
\frac{1}{\abs{\Omega}^{\tau p}}\sum_{R\in\D(\Omega)} \prod_{\nu=1}^k u_{R_\nu}^{(\nu)}
=\prod_{\nu=1}^k\frac{1}{\abs{\Omega_\nu}^{\tau p}}\sum_{R_\nu\in\D(\Omega_\nu)} u_{R_\nu}^{(\nu)}.
\end{equation*}
Taking the supremum over all $\Omega$ of such product form is equivalent to taking the supremum independently over each $\Omega_\nu\in\Open(\R^{n_\nu})$, and we deduce that
\begin{equation}\label{prod-vs-tensor}
\prod_{\nu=1}^k\Norm{t^{(\nu)}}{\dot a^{s_\nu,\tau}_{p,p}(V_\nu)}
\leq\Norm{t^{\otimes}}{\atau(V)};
\end{equation}
this direction did not make use of the assumption that $\tau=\frac1p$.

We finally want to prove the inequality \eqref{prod-vs-tensor} in the reverse direction. This is a bit more delicate, given that \eqref{prod-vs-tensor} as stated was obtained by specialising the condition defining the quasi-norm on the right to rather particular sets of the product form. For the reverse direction, we need to obtain an estimate for general open sets from the knowledge about a much more specific special case.

Suppose without loss of generality that $\Norm{t^{(\nu)}}{\dot a^{s_\nu,\frac{1}{p}}_{p,p}(V_\nu)}\leq 1$ for each $\nu\in[k]$. Let us further assume, still without loss of generality, that the sequences $(t^{(\nu)}_Q)_{Q\in\D(\R^{n_\nu})}$ are finitely non-zero. Indeed, if we can prove the desired estimate of $t^{\otimes}$ in this special case, the general case easily follows by monotone convergence.

Specialising the defining condition of $\Norm{t^{(\nu)}}{\dot a^{s_\nu,\frac1p}_{p,p}(V_\nu)}\leq 1$ from all open sets $\Omega_\nu\in\Open(\R^{n_\nu})$ to just $\Omega_\nu=P\in\D(\R^{n_\nu})$, we have the estimates
\begin{equation}\label{uCarl}
\sum_{Q\in\D(P)}u_Q^{(\nu)}\leq\abs{P}.
\end{equation}

For each $\nu\in[k]$, we claim that there exist disjoint measurable subsets $E_Q^{(\nu)}\subset Q$ such that $u_Q^{(\nu)}=\abs{E_Q^{(\nu)}}$. To accomplish this, we follow an algorithm from the proof of \cite[Lemma 6.3]{LN:basics}. Under the assumption that the sequence $t^{(\nu)}_Q$, and hence $u^{(\nu)}_Q$, is finitely non-zero, for the smallest cubes $Q$ with $u_Q^{(\nu)}>0$, we simply choose any measurable subset $E^{(\nu)}_Q\subset Q$ with $\abs{E^{(\nu)}_Q}=u_Q^{(\nu)}\leq\abs{Q}$ by \eqref{uCarl}. Proceeding recursively, suppose that for some $Q$ with $u_Q^{(\nu)}>0$, we have already defined $E^{(\nu)}_{Q'}$ with required properties for all $Q'\subsetneq Q$. Then
\begin{equation*}
u_Q^{(\nu)}=
\sum_{Q'\subset Q} u^{(\nu)}_Q  -\sum_{Q'\subsetneq Q}u^{(\nu)}_{Q'}
\leq\abs{Q}-\sum_{Q'\subsetneq Q}\abs{E^{(\nu)}_{Q'}}
=\Babs{Q\setminus\bigcup_{Q'\subsetneq Q}E^{(\nu)}_{Q'}},
\end{equation*}
and hence we can find a measurable subset
\begin{equation}\label{new E subset}
E^{(\nu)}_Q\subset Q\setminus\bigcup_{Q'\subsetneq Q}E^{(\nu)}_{Q'}
\end{equation}
of size
$\abs{E^{(\nu)}_Q}=u_Q^{(\nu)}.$
By \eqref{new E subset}, $E^{(\nu)}_Q$ is disjoint from all the previously constructed $E^{(\nu)}_{Q'}$. Since there are only finitely many $Q$ with $u_Q^{(\nu)}>0$, this recursion stops after finitely many steps, and gives the required sets.

Using this representation of each $u_Q^{(\nu)}$, let us also define
\begin{equation*}
E_R^{\otimes}:=\prod_{\nu=1}^k E_{R_\nu}^{(\nu)}\subset R\in\D(\R^{\vn}).
\end{equation*}
Then these $E_R^{\otimes}$ are pairwise disjoint, and $u^{\otimes}_R=\abs{E^{\otimes}_R}$. Then it follows at once, for every $\Omega\in\Open(\R^{\vn})$, that
\begin{equation*}
\sum_{R\in\D(\Omega)}u^{\otimes}_R
=\sum_{R\in\D(\Omega)}\abs{E^{\otimes}_R}
=\Babs{\bigcup_{R\in\D(\Omega)}E^{\otimes}_R}
\leq\abs{\Omega}.
\end{equation*}
This completes the proof that $\Norm{t^{\otimes}}{\dot a^{s,\frac{1}{p}}_{p,p}(V)}\leq 1$ when $\Norm{t^{(\nu)}}{\dot a^{s_\nu,\frac{1}{p}}_{p,p}(V_\nu)}\leq 1$ for every $\nu\in[k]$, and hence the proof of Lemma \ref{t tensor}.
\end{proof}

With the help of Lemma \ref{t tensor}\eqref{tau=0}, we obtain the following necessary conditions for the boundedness of almost diagonal operators.

\begin{proposition}\label{sharpness of AD}
Let $\vp\in(0,\infty)^k$, $\vq\in(0,\infty]^k$, $\vs\in\R^k$, and $\pi\in S_{[2k]}$.
Suppose that, for some $(\vD,\vE,\vF)\in\R^{3k}$, every $(\vD,\vE,\vF)$-almost diagonal operator is well defined and bounded on $\asn$.
Then, for every $\nu\in[k]$,
\begin{equation}\label{DiEiFiNec}
D_\nu>\frac{n_\nu}{\min(1,p_\nu)},\quad
E_\nu>\frac{n_\nu}{2},\quad
F_\nu>\frac{n_\nu}{\min(1,p_\nu)}-\frac{n_\nu}{2}-s_\nu.
\end{equation}
\end{proposition}

\begin{proof}
The proof consists of a reduction to the one-parameter case via Lemma \ref{t tensor}.

For $t^{(\nu)}\in\dot a^{s_\nu,0}_{p_\nu,q_\nu,\pi_\nu}$, let $t^{\otimes}$ be defined as in Lemma \ref{t tensor}. In a similar way, for $(D_\nu,E_\nu,F_\nu)$-almost diagonal matrices $B^{(\nu)}=\{b^{(\nu)}_{I,J}\}_{I,J\in\D(\R^{n_\nu})}$, we define
\begin{equation*}
B^{\otimes}:=\{b^{\otimes}_{P,R}\}_{P,R\in\D(\R^{\vn})},\quad
b^{\otimes}_{P,R}:=\prod_{\nu=1}^k b^{(\nu)}_{P_\nu,R_\nu}.
\end{equation*}
It is then evident that $B^{\otimes}$ is $(\vD,\vE,\vF)$-almost diagonal. Moreover, if $u^{(\nu)}:=B^{(\nu)}t^{(\nu)}$ and $u^{\otimes}$ is constructed from these $u^{(\nu)}$ in the same way as $t^{\otimes}$ from the $t^{(\nu)}$, then $B^{\otimes}t^{\otimes}=u^{\otimes}$.

The assumption that every $(\vD,\vE,\vF)$-almost diagonal matrix is bounded on $\atau$, together with Lemma \ref{t tensor}, then implies that
\begin{equation*}
\prod_{\nu=1}^k\Norm{B^{(\nu)}t^{(\nu)}}{\dot a^{s_\nu,0}_{p_\nu,q_\nu,\pi_\nu}}
=\Norm{B^{\otimes}t^{\otimes}}{\asn}
\lesssim\Norm{t^{\otimes}}{\asn}
=\prod_{\nu=1}^k\Norm{t^{(\nu)}}{\dot a^{s_\nu,0}_{p_\nu,q_\nu,\pi_\nu}}.
\end{equation*}
With all but one $B^{(\nu)}$ equal to $\{\delta_{I,J}\}_{I,J\in\D(\R^{n_\nu})}$, it follows that $B^{(\nu)}$ must be bounded on the respective $\dot a^{s_\nu,0}_{p_\nu,q_\nu,\pi_\nu}$. By \cite[Lemma 7.1 or  9.1]{BHYY:1b} (depending on whether $\pi_i\in\{B,F\}$, but the conclusion is the same), a necessary condition for this boundedness is precisely \eqref{DiEiFiNec}. Since this applies to every $\nu\in[k]$, the proof is complete.
\end{proof}

\begin{remark}
Proposition \ref{sharpness of AD} shows in particular that the assumptions of Theorem \ref{AD s0} are sharp whenever $p_i\equiv p\in(0,\infty)$ and $r(\pi,\vp,\vq)\geq \min(1,p)$ for all $i\in[k]$. Assuming the first condition $\vp=p\cdot\vone$, the second condition is automatic in the Besov case [when $r(\pi,\vp,\vq)=p$] and requires that $q_i\geq\min(1,p)$ in the Triebel--Lizorkin case (when $r(\pi,\vp,\vq)=\min\{p,q_i:i\in[k]\}$).

Thus, for non-mixed exponents $\vp=p\cdot\vone$ and $\vq=q\cdot\vone$, Theorem \ref{AD s0} is sharp in the Besov case for all exponents, and in the Triebel--Lizorkin case if $q\geq\min(1,p)$. In particular, the assumptions of \cite[Theorem 5.2]{GKP21}, the special case of Theorem \ref{AD s0} for unweighted bi-parameter spaces with non-mixed exponents $\vp=(p,p)$ and $\vq=(q,q)$, are sharp under the said restriction $q\geq\min(1,p)$ in the Triebel--Lizorkin case. This sharpness is new even in this special case.

The question about sharp almost diagonal conditions for Triebel--Lizorkin spaces in the remaining range $q<\min(1,p)$ is open even in the one-parameter case \cite{BHYY:1b}, where the same restriction appears.
\end{remark}

We will then discuss the sharpness of Theorem \ref{AD st} for the special parameter $\tau=\frac1p$. By Remark \ref{rem AD st}, we already know that sharper conditions than those of Theorem \ref{AD st} can be found in the one-parameter case by taking into account the $\A_p$-dimension of the weights $V\in\A_p(\R^{\vn})$. Thus, a meaningful way to discuss possible sharpness is to inquire about the minimal almost diagonal conditions that guarantee the boundedness for all $V\in\A_p(\R^{\vn})$. To this end, the following necessary condition from the one-parameter case is a slight reformulation of \cite[Lemma 11.1]{BHYY:1b}:

\begin{lemma}\label{nec 1param}
Let $p\in(0,\infty)$, $q\in(1,\infty]$, $s\in\R$, and $\pi\in S_{[2]}$. Let $(D,E,F)\in\R^3$ be such that every $(D,E,F)$-almost diagonal operator is bounded on $\dot a^{s,\frac1p}_{p,q,\pi}(V)$ for every $V\in\A_p(\R^n)$, or even just every power weight $V\in\A_p(\R^n)$. Then
\begin{equation*}
D\geq n+\frac{n}{p},\quad E\geq\frac{n}{2}+s+\frac{n}{p},\quad F>\frac{n}{2}-s.
\end{equation*}
\end{lemma}

\begin{proof}
\cite[Lemma 11.1]{BHYY:1b} contains the following result, formulated in terms of $W\in A_p(\R^n)$ instead of $V\in\A_p(\R^n)$. Recall that these are in one-to-one correspondence via $V(x)=W(x)^{\frac1p}$. If the assumptions of the lemma hold just for the single weight $W(x)=\abs{x}^{-d}$ (thus $V(x)=\abs{x}^{-\frac{d}{p}}$), then
\begin{equation*}
D> n+\frac{d}{p},\quad E>\frac{n}{2}+s+\frac{d}{p},\quad F>\frac{n}{2}-s.
\end{equation*}
Since $W\in A_p(\R^n)$ for all $d\in[0,n)$, if the assumptions of the lemma holds as stated, the previous display must hold for all $d\in[0,n)$, and the only way that this can be is that $(D,E,F)$ must be as stated in the lemma.
\end{proof}

\begin{corollary}\label{AD st sharp}
Except possibly for the end-points for $D$ and $E$, the assumptions \eqref{DEF2} of Theorem \ref{AD st} are sharp, in the sense of being the minimal assumptions to work for all weights $V\in\A_p(\R^n)$  simultaneously, in the one-parameter case with $\tau=\frac1p$, $p\in[1,\infty)$, and $q\in(1,\infty]$.
\end{corollary}

\begin{proof}
The assumption that $p\geq 1$ and $q>1$ implies that, irrespective of $\pi$, there holds $r(\pi,p,q)\geq 1$, and hence $r=r(\pi,p,q)\wedge 1=1$.
In this situation, the assumptions \eqref{DEF2} of Theorem \ref{AD st} read as
\begin{equation}\label{AD st rep}
D>n+\frac{n}{p},\quad
E>\frac{n}{2}+\frac{n}{p}+s,\quad
F>n-\frac{n}{2}-s=\frac{n}{2}-s.
\end{equation}

Note that the assumptions of Corollary \ref{AD st sharp} on $p,q,\tau$, and $V$ are special cases of the assumptions of Lemma \ref{nec 1param}. Hence, if  all $(D,E,F)$-almost diagonal operators are bounded on $\dot a^{0,\frac1p}_{p,q,\pi}(V)$ for all $V\in\A_p(\R^n)$, then the conclusions of Lemma \ref{nec 1param} on $D,E,F$ are in force. Except for the end points of the conditions for $D$ and $E$, these are seen to coincide with \eqref{AD st rep}, proving the claimed sharpness.
\end{proof}

We also obtain necessary conditions at least for a limited case of genuinely multi-parameter situations with $\tau=\frac1p$. The parameters here are restricted by the fact that they need to be compatible with both the one-parameter result of Lemma \ref{nec 1param} and the lifting result of Lemma \ref{t tensor}\eqref{crit}.

\begin{proposition}\label{nec crit}
Let $p\in(1,\infty)$ and $s\in\R^k$, and $\pi\in S_{[2k]}$. Let $(\vD,\vE,\vF)\in\R^{3k}$ be such that every $(\vD,\vE,\vF)$-almost diagonal operator is bounded on $\dot a^{\vs,\frac{1}{p}}_{p,p,\pi}(V)$ for every $V\in\A_p(\R^{\vn})$, or even just every product power weight $V(x)=\prod_{\nu=1}^k\abs{x_\nu}^{\alpha_\nu}$ in this class.
Then
\begin{equation*}
\vD\geq \vn+\frac{\vn}{p},\quad \vE\geq\frac{\vn}{2}+\vs+\frac{\vn}{p},\quad \vF>\frac{\vn}{2}-\vs.
\end{equation*}
\end{proposition}

Note that the permutation $\pi$ is redundant in this case, since it is immediate from Fubini's theorem that $\dot a^{\vs,\frac{1}{p}}_{p,p,\pi}(V)=\dot a^{\vs,\frac{1}{p}}_{p,p,\pi'}(V)$ for any two $\pi,\pi'\in S_{[2k]}$, but we have kept it to avoid introducing ad hoc notation just for this one place in the paper.

\begin{proof}[Proof of Proposition \ref{nec crit}]
This follows from a combination of Lemmas \ref{t tensor}\eqref{crit} (with weights $V_\nu(x)=\abs{x_\nu}^{\alpha_\nu}$) and \ref{nec 1param} in the same way that Proposition \ref{sharpness of AD} was proved using Lemma \ref{t tensor}(\ref{tau=0}) and the necessary conditions for in the one-parameter case from \cite{BHYY:1b}.
\end{proof}

\begin{corollary}\label{AD st sharp2}
Except possibly for the end-points for $\vD$ and $\vE$, the assumptions \eqref{DEF2} of Theorem \ref{AD st} are sharp, in the sense of being the minimal assumptions to work for all weights $V\in\A_p(\R^{\vn})$  simultaneously, in the case that $p\in(1,\infty)$, $\vp=\vq=p\cdot\vone$, and $\tau=\frac1p$.
\end{corollary}

\begin{proof}
This follows from a comparison of the sufficient conditions of Theorem \ref{AD st} with the necessary conditions of Proposition \ref{nec crit} in the same way as Corollary \ref{AD st sharp} was deduced from the comparison of the one-parameter case of Theorem \ref{AD st} with Lemma \ref{nec 1param}.
\end{proof}

\section{Molecules and their applications}\label{sec:molecules}

\subsection{Molecular characterisation}\label{sec mp}

In this section, we establish the molecular characterization of $\Atau(V)$, i.e., we describe basic building blocks (the ``molecules'') such that, in a sense to be made precise in Theorem \ref{89},
\begin{enumerate}[\rm(1)]
  \item\label{it:analysis} every distribution in $\Atau(V)$ can be represented as a convergent series of these molecules, with suitable estimates for the coefficients;
  \item\label{it:synthesis} every distribution that admits such an expansion must belong to $\Atau(V)$, and its norm in this space can be estimated by the size of the coefficients.
\end{enumerate}
This can be seen as generalisation of the $\varphi$-transform characterisation (Theorem \ref{phi}), but allowing more flexibility in the choice of the molecules. This is particularly useful in verifying the $\Atau(V)$-membership with the help of a molecular expansion, which will be key to the $T(1)$ theorem (a boundedness criterion for singular integral operators) on the matrix-weighted multi-parameter spaces that we consider (Theorem \ref{T1 BF}).

In the unweighted one-parameter case, such molecular characterisations go back to Frazier and Jawerth \cite{FJ90}; the matrix-weighted one-parameter case is due to Bu et al.\ \cite{BHYY:1c}. The goal of this section is a multi-parameter extension of these results. While following the same broad outline, there are additional technical challenges. Notably, in the one-parameter setting, one can typically estimate the pairing of two functions by using the smoothness of one of them and cancellation properties of the other one, depending on the relative scales of the two functions. In the multi-parameter contexts, each function will have independent scaling properties in different component directions, and we will need to couple the smoothness of $f$ with cancellation of $g$ in some directions, and cancellations of $f$ with the smoothness of $g$ in the remaining ones; see e.g.\ \eqref{eq:g*h2}.

Turning to the details, we first give some notation. For every $r\in\mathbb R$, let
\begin{equation*}
\begin{cases}
\lfloor r\rfloor:=\max\{k\in\mathbb Z:\ k\leq r\},\\
\lfloor\!\lfloor r\rfloor\!\rfloor:=\max\{k\in\mathbb Z:\ k< r\},
\end{cases}
\begin{cases}
\lceil r\rceil:=\min\{k\in\mathbb Z:\ k\geq r\},\\
\lceil\!\lceil r\rceil\!\rceil:=\min\{k\in\mathbb Z:\ k>r\},
\end{cases}
\end{equation*}
and $r^{**}:=r-\lfloor\!\lfloor r\rfloor\!\rfloor\in(0,1]$.
For every $\gamma:=(\gamma_1,\ldots,\gamma_k)\in\mathbb{Z}^{\vec n}$,
$|\gamma|_{\rm{vec}}:=(|\gamma_1|,\ldots,|\gamma_k|)$.
For every $\vec a,\vec b\in\mathbb R^k$, let
$$
\vec a\wedge \vec b:=(a_1\wedge b_1,\ldots,a_k\wedge b_k),
$$
$$
\vec a\vee \vec b:=(a_1\vee b_1,\ldots,a_k\vee b_k).
$$

The following definition is a multi-parameter extension of \cite[Definition 3.4]{BHYY:1c}, which in turn is inspired by several related definitions in the earlier literature. There is some tradition in the literature of fixing, from the beginning, a choice of the molecular parameters that is convenient for the specific situation at hand. Instead, we follow the style of \cite{BHYY:1c} of first defining molecules with generic parameters, and only a little later ``discovering'' a choice of these parameters that is convenient in the context of the spaces $\Atau$.
A version of product-space molecules has also been previously introduced in \cite[Definition 14]{LLZ}.

\begin{definition}\label{moleKLMN}
Let $\vec K,\vec M\in[0,\infty)^k$, $\vec L\in\Z_{\geq-1}^k:=
(\Z_+\cup\{-1,0\})^k$, and $\vec N\in(0,\infty)^k$.
A function $m_Q$ is called a \emph{
$(\vec K,\vec L,\vec M,\vec N)$-molecule on $Q\in\D(\mathbb R^{\vec n})$}
if it satisfies
\begin{enumerate}[\rm(i)]

  \item the size and smoothness condition: for all $\gamma\in\N^{\vn}$ with $\abs{\gamma}_{\rm{vec}}<\vec N$, all $\alpha\in\{0,1\}^k$, and all $x,h\in\R^{\vn}$:
\begin{equation}\label{co-2}
  \abs{\Delta_h^\alpha\partial^\gamma m_Q(x)}
  \leq   \abs{Q}^{-\frac12}
  \prod_{\nu=1}^k  \begin{cases} \displaystyle
   \Big[1+\frac{\abs{x_\nu-c_{Q_\nu}}}{\ell(Q_\nu)}\Big]^{-(K_{\nu}\vee M_{\nu})},\quad \text{if }\alpha_\nu=0=\gamma_\nu, \\
   \displaystyle
  \Big[1+\frac{\abs{x_\nu-c_{Q_\nu}}}{\ell(Q_\nu)}\Big]^{-M_{\nu}}\ell(Q_\nu)^{-\abs{\gamma_\nu}},
  \quad \text{if }\alpha_\nu=0\neq \gamma_\nu, \\
  \displaystyle
  \sup_{\abs{z_\nu-x_\nu}\leq\abs{h_\nu}} \Big[1+\frac{\abs{z_\nu-c_{Q_\nu}}}{\ell(Q_\nu)}\Big]^{-M_{\nu}}
  \frac{\abs{h_\nu}^{N_\nu^{**}}}{ \ell(Q_\nu)^{\abs{\gamma_\nu}+N_\nu^{**}}},
  \quad \text{if }\alpha_\nu\neq 0,
\end{cases}
\end{equation}
where
\begin{equation}\label{Delta}
  \Delta_h^\alpha:=\prod_{\nu\in[k]:\alpha_\nu=1} \Delta_{h_\nu}^{(\nu)},\qquad
   \Delta_{h_\nu}^{(\nu)}f(x)
   :=f(x)-f(x_\nu-h_\nu,x');
\end{equation}

\item the cancellation condition: for every $\nu\in[k]$ and
$\gamma_\nu\in\mathbb{N}^{n_\nu}$ with $|\gamma_\nu|\leq L_\nu$,
and a.e.\ $x'\in\R^{\vn}\ominus\R^{n_\nu}$,
\begin{equation}\label{eq:moleL}
  \int_{\mathbb R^{n_\nu}} x_\nu^{\gamma_\nu} m_Q(x_\nu,x')\,dx_\nu=0.
\end{equation}
\end{enumerate}
\end{definition}

\begin{remark}\label{rem:moleKLMN}
The parameters $\vec L,\vec N$ count the number of vanishing moments and derivatives of $m_Q$, respectively. Allowing the value $L_\nu=-1$ is relevant to account for the possibility of not imposing any vanishing moments \eqref{eq:moleL} in the $x_\nu$ variable. Molecules with no vanishing moments are useful in trace theorems; see \cite[Remark 5.7]{BHYY:1c}.

To keep the technical aspects at a reasonable level,
we restrict the range of $\vec N$ into $(0,\infty)^k$
in Definition \ref{moleKLMN}.
More precisely, when
$k=1$ and $N\in(0,\infty)$,
condition \eqref{co-2} agrees with
\cite[(1.1.6), (1.1.8), and (1.1.9)]{tor}; however,
in our setting, we only consider positive values of $N$.
\end{remark}

The following lemma is a multi-parameter extension of  \cite[Lemma B.1]{FJ90}. It uses partially the same ideas, but requires some additional care with iterated Taylor expansions with respect to different variables.

\begin{lemma}\label{fj B1}
Let $\vec L\in\Z_{\geq-1}^k$, $\vec N\in(0,\infty)^k$, and $\vec K,\vec M\in[0,\infty)^k$.
Let $g\in L^1(\R^{\vn})$ and suppose that its distributional derivatives up to order $\vec E$ are also given by $L^1(\R^{\vn})$-functions.
Suppose that, for some $\vj\in\Z^k$, these functions satisfy
for every $\gamma\in\mathbb Z_+^{\vec n}$
with $|\gamma|_{\rm{vec}}<\vec N$, every $\alpha\in\{0,1\}^k$,
and every $x,y\in\mathbb{R}^{\vn}$,
\begin{equation}\label{g11}
  \abs{\Delta_y^\alpha\partial^\gamma g(x)}
  \leq   2^{{\vj\cdot\frac{\vn}{2}}}
  \prod_{\nu=1}^k  \begin{cases} \displaystyle
   (1+2^{j_\nu}\abs{x_\nu})^{-(K_{\nu}\vee M_{\nu})}, &\text{if }\alpha_\nu=0=\gamma_\nu, \\
   \displaystyle
  (1+2^{j_\nu}\abs{x_\nu})^{-M_{\nu}}2^{j_\nu \abs{\gamma_\nu}},
  & \text{if }\alpha_\nu=0\neq \gamma_\nu, \\
  \displaystyle
  \sup_{\abs{z_\nu-x_\nu}\leq\abs{y_\nu}} (1+2^{j_\nu}\abs{z_\nu})^{-M_{\nu}}2^{j_\nu(\abs{\gamma_\nu}+{N_\nu^{**}})}\abs{y_\nu}^{N_\nu^{**}},
  & \text{if }\alpha_\nu\neq 0,
\end{cases}
\end{equation}
where $\Delta_y^\alpha$ is defined in \eqref{Delta}; and
for every $\nu\in[k]$ and every $\gamma_\nu\in\mathbb Z_+^{n_\nu}$
with $|\gamma_\nu|\leq L_\nu$,
and a.e.\ $x'\in\R^{\vn}\ominus\R^{n_{\nu}}$,
\begin{equation}\label{h2}
   \int_{\mathbb R^{n_\nu}}x_\nu^{\gamma_\nu} g(x_\nu,x')\,dx_\nu=0,
\end{equation}
where this last condition is void for any $\nu\in[k]$ with $L_\nu=-1$.

Further assume that $h\in L^1(\R^{\vn})$ satisfies the same assumptions as above with some $\vec K'$, $\vec L'$, $\vec M'$, $\vec N'$ in place of $\vec K,\vec L,\vec M,\vec N$, and some $\vec l$ in place of $\vj$. We further assume that
\begin{equation}\label{eq:K>N+n}
  \vec K'>\vec N+\vn,\qquad
  \vec K>\vec N'+\vn,
\end{equation}
and
\begin{equation}\label{eq:L=N-1}
  \vec L'=\lfloor\!\lfloor \vec N\rfloor\!\rfloor,\qquad
  \vec L=\lfloor\!\lfloor \vec N'\rfloor\!\rfloor.
\end{equation}

Then there exists a positive constant $C$,
independent of $\vj$, $\vl$, and $x$, such that, for every $x\in\mathbb{R}^{\vn}$,

\begin{equation}\label{eq:g*h1}
|g*h(x)|\leq
C 2^{-(\vl-\vj)_+\cdot(\frac{\vn}{2}+\vec N)-(\vj-\vl)_+\cdot(\frac{\vn}{2}+\vec N')}
\prod_{\nu=1}^{k}\Big(1+2^{j_\nu \wedge l_\nu}|x_\nu|\Big)^{-(M_\nu\wedge M_\nu')}.
\end{equation}
\end{lemma}

\begin{proof}
For any $y=(y_1,\ldots,y_k)\in\R^{\vn}$ and $I\subseteq[k]$,
we write $y_I=(y_\nu)_{\nu\in I}$. Let
\begin{equation*}
  G:=\{\nu\in[k]:j_\nu\leq l_\nu\},\qquad H:=[k]\setminus G.
\end{equation*}
Then
\begin{equation*}
  g*h(x)=\int_{\R^{\vn_{H}}} \int_{\R^{\vn_{G}}} g(y_G,x_H-y_H)h(x_G-y_G,y_H)dy_G\,dy_H,
\end{equation*}
where $\R^{\vn_{H}}:=\bigotimes_{\nu\in H} \R^{n_{\nu}}$
and $\R^{\vn_{G}}:=\bigotimes_{\nu\in G} \R^{n_{\nu}}$.
We recall the multivariate Taylor polynomial
\begin{equation*}
    \operatorname{Tayl}_{x,d}f(y) :=\sum_{\abs{\alpha}\leq d}\frac{(y-x)^\alpha}{\alpha !}\partial^\alpha f(x)
\end{equation*}
and the Taylor formula for the error term
\begin{equation}\label{eq:TaylErr}
f(y)-\operatorname{Tayl}_{x,d-1}f(y) =\sum_{\abs{\beta}=d}\frac{(y-x)^\beta}{\beta!}\int_0^1 \partial^\beta f(x+t(y-x)) d(1-t)^{d-1}dt,
\end{equation}
from which we also deduce
\begin{equation}\label{eq:TaylErr2}
f(y)-\operatorname{Tayl}_{x,d}f(y) =\sum_{\abs{\beta}=d}\frac{(y-x)^\beta}{\beta!}\int_0^1 \big[\partial^\beta f(x+t(y-x))-\partial^\beta f(x)\big] d(1-t)^{d-1}dt.
\end{equation}

Originally derived for functions with continuous partial derivatives, this formula remains valid for a.e.\ $x,y\in\R^{\vn}$ assuming only that the distributional derivatives are given by integrable functions: this can be proved by pairing both sides with a test function, moving all derivatives and translations to the test function side, and using \eqref{eq:TaylErr} on the test function; details can be found in \cite{Anas}.

By the cancellation assumption \eqref{h2} for both $g$ and $h$, and the relations \eqref{eq:L=N-1}, for each $\nu\in[k]$, we can subtract a Taylor polynomial from one of the functions $g$ or $h$, since the integral of this polynomial times the other function, over the variable $x_\nu$, will vanish for a.e.\ value of the other variables. Hence
\begin{equation}\label{eq:g*h2}
\begin{split}
  g*h(x)=\int_{\R^{\vn_{H}}} \int_{\R^{\vn_{G}}}  &\Big(\prod_{\nu\in G}(I-\operatorname{Tayl}_{x_\nu,\lfloor\!\lfloor N_\nu\rfloor\!\rfloor}^{(\nu)})\Big)g(y_G,x_H-y_H) \\
  &\times  \Big(\prod_{\eta\in H}(I-\operatorname{Tayl}_{x_\eta,\lfloor\!\lfloor N_\eta'\rfloor\!\rfloor}^{(\eta)})\Big)h(x_G-y_G,y_H)dy_G\,dy_H,
\end{split}
\end{equation}
where $\operatorname{Tayl}^{(\nu)}$ indicates Taylor expansion with respect to the $\nu$-th coordinate.

Clearly, the factors related to both $g$ and $h$ in \eqref{eq:g*h2} will be some combinations of the values of these functions and their derivatives. These both satisfy upper bounds \eqref{g11} that split into a product of factors, each of which only depends on one coordinate $y_\nu$ at the time. Hence estimating $\abs{g*h(x)}$ by bringing the absolute values inside the integral and applying \eqref{g11}, we can estimate the integral over each $y_\nu$ at a time.

Consider some $\nu\in G$. (The case of $\nu\in H$ will be entirely symmetric by swapping the roles of $g$ and $h$.) In this case, there is no Taylor expansion in the $y_\nu$ direction acting on~$h$. Hence, the bound for the $y_\nu$ factor of $h$ will be simply that from (the $h$ analogue of) \eqref{g11} with $\alpha_\nu=0=\gamma_\nu$, that is,
\begin{equation*}
  h_\nu(x_\nu-y_\nu)=2^{l_\nu\frac{n_\nu}{2}}(1+2^{l_\nu}\abs{x_\nu-y_\nu})^{-(K_\nu'\vee M_\nu')}.
\end{equation*}
For the $g$ factor, we have two alternative bounds, either via the Taylor formula \eqref{eq:TaylErr2}, or estimating the $I$ and
$\operatorname{Tayl}_{x_\nu,\lfloor\!\lfloor N_\nu\rfloor\!\rfloor}^{(\nu)}$ terms separately. Using \eqref{g11} with $t(x-y)$ in place of $y$, the first option leads to the bound
\begin{equation}\label{eq:gI}
   g_{\nu}(y_\nu;x_\nu)\leq \sum_{\abs{\beta}=\lfloor\!\lfloor N_\nu\rfloor\!\rfloor}\abs{x_\nu-y_\nu}^{\abs{\beta}}2^{j_\nu(\frac{n_\nu}{2}+\abs{\beta}+N_{\gamma}^{**})}
  \abs{x_\nu-y_\nu}^{N_{\gamma}^{**}}
   \sup_{\abs{z_\nu-x_\nu}\leq\abs{y_\nu-x_\nu}}
   (1+2^{j_\nu}\abs{z_\nu})^{-M_\nu}
\end{equation}
and the second option to the bound
\begin{equation}\label{eq:gII}
  g_{\nu}(y_\nu;x_\nu)\leq    2^{j_\nu\frac{n_\nu}{2}}(1+2^{j_\nu}\abs{y_\nu})^{-(K_\nu\vee M_\nu)}
    +\sum_{\abs{\alpha}\leq \lfloor\!\lfloor N_\nu\rfloor\!\rfloor}\abs{x_\nu-y_\nu}^{\abs{\alpha}}2^{j_\nu(\frac{n_\nu}{2}+\abs{\alpha})}(1+2^{j_\nu}\abs{x_\nu})^{-M_\nu},
\end{equation}
respectively.

If $\abs{x_\nu-y_\nu}<3\cdot 2^{-j_\nu}$, we note that the variable $z_\nu$ appearing in \eqref{eq:gI} satisfies
$2^{j_\nu}\abs{z_\nu-x_\nu}\leq 2^{j_\nu}\abs{y_\nu-x_\nu}<3,$
and hence $(1+2^{j_\nu}\abs{z_\nu})^{-M_\nu}\sim (1+2^{j_\nu}\abs{x_\nu})^{-M_\nu}$.
Also, noting that $\abs{\beta}=\lfloor\!\lfloor N_\nu\rfloor\!\rfloor$ is constant in the summation in \eqref{eq:gI},
and hence $|\beta|+N_{\nu}^{**}=N_{\nu}$, this estimate simplifies to
\begin{equation}\label{eq:gI2}
  g_\nu(y_\nu;x_\nu)\lesssim
  2^{j_\nu\frac{n_\nu}{2}}(2^{j_\nu}\abs{x_\nu-y_\nu})^{N_\nu}(1+2^{j_\nu}\abs{x_\nu})^{-M_\nu},\quad
  \abs{x_\nu-y_\nu}<3\cdot 2^{-j_\nu}.
\end{equation}
On the other hand, if $\abs{x_\nu-y_\nu}\geq 3\cdot 2^{-j_\nu}$, then $(2^{j_\nu}\abs{x_\nu-y_\nu})^{\abs{\alpha}}<(2^{j_\nu}\abs{x_\nu-y_\nu})^{N_\nu}$
for all $\abs{\alpha}\leq \lfloor\!\lfloor N_\nu\rfloor\!\rfloor<N_\nu$ in the summation in \eqref{eq:gII}, and hence this bound also simplifies to
\begin{equation}\label{eq:gII2}
  g_\nu(y_\nu;x_\nu)\lesssim
   2^{j_\nu\frac{n_\nu}{2}}\bigg[(1+2^{j_\nu}\abs{y_\nu})^{-(K_\nu\vee M_\nu)}
   +\frac{ (2^{j_\nu}\abs{x_\nu-y_\nu})^{N_\nu} }{(1+2^{j_\nu}\abs{x_\nu})^{M_\nu} }\bigg],\quad
  \abs{x_\nu-y_\nu}\geq 3\cdot 2^{-j_\nu}.
\end{equation}
Note that the second term on the right-hand side of \eqref{eq:gII2} is the same as the right-hand side of \eqref{eq:gI2}. Moreover, if $\abs{y_\nu}>\frac12\abs{x_\nu}$, then the first term on the right-hand side of \eqref{eq:gII2} is dominated by the second one. Thus, the first term on the right-hand side of \eqref{eq:gII2} is only needed in the remaining case when both $\abs{x_\nu-y_\nu}\geq 3\cdot 2^{-j_\nu}$ and $\abs{y_\nu}\leq\frac12\abs{x_\nu}$.
Combining \eqref{eq:gI2} and \eqref{eq:gII2} with the observation just made, we arrive at
\begin{equation}\label{eq:gIII}
\begin{split}
  g_\nu(y_\nu;x_\nu)
  &\lesssim 2^{j_\nu\frac{n_\nu}{2}} \Bigg[
  \frac{ (2^{j_\nu}\abs{x_\nu-y_\nu})^{N_\nu} }{ (1+2^{j_\nu}\abs{x_\nu})^{M_\nu}}
  +\frac{\one_{\{2^{j_\nu}\abs{x_\nu-y_\nu}\geq 3,\ \abs{y_\nu}\leq\frac12\abs{x_\nu}\}}}{(1+2^{j_\nu}\abs{y_\nu})^{K_\nu\vee M_\nu}}\Bigg] \\
  &=:g_\nu^0(y_\nu;x_\nu)+g_\nu^1(y_\nu;x_\nu).
\end{split}
\end{equation}
(Note that, without the indicator, $g^1_\nu(y_\nu;x_\nu)$ would be just the trivial upper bound that we could have made without any Taylor expansions; what we have gained in \eqref{eq:gIII} is that this bound is only needed for restricted values of the variables, while otherwise we can use the new upper bound $g_\nu^0(y_\nu;x_\nu)$ instead.)

To estimate $\int_{\R^{n_\nu}} g_\nu(y_\nu)h_\nu(x_\nu-y_\nu)dy_\nu$, we first consider
\begin{equation}\label{eq:g0}
\begin{split}
  \int_{\R^{n_\nu}} &g_\nu^0(y_\nu;x_\nu)h_\nu(x_\nu-y_\nu)dy_\nu\\
  &\lesssim \frac{2^{(j_\nu+l_\nu)\frac{n_\nu}{2}}}{(1+2^{j_\nu}\abs{x_\nu})^{M_\nu}}
  \int_{\R^{n_\nu}}
  \frac{(2^{j_\nu}\abs{x_\nu-y_\nu})^{N_\nu}}{(1+2^{l_\nu}\abs{x_\nu-y_\nu})^{K_\nu'}}dy_\nu \\
  &\lesssim \frac{2^{(j_\nu+l_\nu)\frac{n_\nu}{2}}}{(1+2^{j_\nu}\abs{x_\nu})^{M_\nu}}
  \int_0^\infty\frac{(2^{j_\nu-l_\nu}r)^{N_\nu}}{(1+r)^{K_\nu'}} 2^{-l_\nu n_\nu}r^{n_\nu-1}dr
  \sim\frac{2^{(j_\nu-l_\nu)(\frac{n_\nu}{2}+N_\nu)}}{(1+2^{j_\nu}\abs{x_\nu})^{M_\nu}},
\end{split}
\end{equation}
where condition \eqref{eq:K>N+n} guarantees the convergence of the integral in the last step.

It remains to estimate $\int_{\R^{n_\nu}} g_\nu^1(y_\nu;x_\nu)h_\nu(x_\nu-y_\nu)dy_\nu$.
Under the conditions in the indicator appearing in $g_\nu^1(y_\nu;x_\nu)$ in \eqref{eq:gIII}, we have $\abs{x_\nu}\sim\abs{x_\nu-y_\nu}\gtrsim 2^{-j_\nu}$. Hence
\begin{equation*}
  1+2^{l_\nu}\abs{x_\nu-y_\nu}
  \gtrsim 2^{l_\nu}(2^{-j_\nu}+\abs{x_\nu})=2^{l_\nu-j_\nu}(1+2^{j_\nu}\abs{x_\nu})
  \geq 2^{l_\nu-j_\nu}(1+2^{j_\nu}\abs{y_\nu})
\end{equation*}
and thus
\begin{equation*}
  h_{\nu}(x_\nu-y_\nu)\lesssim 2^{l_\nu\frac{n_\nu}{2}}\Big[2^{l_\nu-j_\nu}(1+2^{j_\nu}\abs{x_\nu})\Big]^{-(K_\nu'\vee M_\nu')}
  \leq 2^{l_\nu\frac{n_\nu}{2}} 2^{(j_\nu-l_\nu)K_\nu'}(1+2^{j_\nu}\abs{x_\nu})^{-M_\nu'}.
\end{equation*}
We can then compute
\begin{equation*}
\begin{split}
  &\int_{\R^{n_\nu}} g_\nu^1(y_\nu;x_\nu) h_\nu(x_\nu-y_\nu)dy_\nu\\
  &\quad\lesssim 2^{(j_\nu+l_\nu)\frac{n_\nu}{2}} \frac{2^{(j_\nu-l_\nu)K_\nu'}}{(1+2^{j_\nu}\abs{x_\nu})^{M_\nu'}}
  \int_{\R^{n_\nu}} \frac{dy_\nu}{(1+2^{j_\nu}\abs{y_\nu})^{K_\nu}} \\
  &\quad\sim 2^{(j_\nu+l_\nu)\frac{n_\nu}{2}} \frac{2^{(j_\nu-l_\nu)K_\nu'}}{(1+2^{j_\nu}\abs{x_\nu})^{M_\nu'}}2^{-j_\nu n_\nu} \quad
  \text{since $K_\nu>N_\nu'+n_\nu> n_\nu$ by \eqref{eq:K>N+n}} \\
  &\quad= 2^{(j_\nu-l_\nu)\frac{n_\nu}{2}} \frac{2^{(j_\nu-l_\nu)(K_\nu'-n_\nu)}}{(1+2^{j_\nu}\abs{x_\nu})^{M_\nu'}}
  \leq\frac{2^{(j_\nu-l_\nu)(\frac{n_\nu}{2}+N_\nu)}}{(1+2^{j_\nu}\abs{x_\nu})^{M_\nu'}}\quad
  \text{since $K_\nu'>N_\nu+n_\nu$ by \eqref{eq:K>N+n}};
\end{split}
\end{equation*}
note that we used both bounds in \eqref{eq:K>N+n} in different steps of this estimate.

Combining the previous computation with \eqref{eq:g0}, we arrive at
\begin{equation*}
   \int_{\R^{n_\nu}} g_\nu(y_\nu;x_\nu) h_\nu(x_\nu-y_\nu)dy_\nu
   \lesssim\frac{2^{(j_\nu-l_\nu)(\frac{n_\nu}{2}+N_\nu)}}{(1+2^{j_\nu}\abs{x_\nu})^{M_\nu\wedge M_\nu'}}
   =\frac{2^{-(l_\nu-j_\nu)_+(\frac{n_\nu}{2}+N_\nu)}}{(1+2^{j_\nu\wedge l_\nu}\abs{x_\nu})^{M_\nu\wedge M_\nu'}},
\end{equation*}
recalling that this bound was obtained for indices $\nu\in G$, i.e., the ones with $j_\nu\leq l_\nu$.

For $\nu\in H=[k]\setminus G$, it follows by symmetry that we obtain the analogous bound
\begin{equation*}
   \int_{\R^{n_\nu}} g_\nu(x_\nu-y_\nu) h_\nu(y_\nu;x_\nu)dy_\nu
   \lesssim\frac{2^{-(j_\nu-l_\nu)_+(\frac{n_\nu}{2}+N_\nu')}}{(1+2^{j_\nu\wedge l_\nu}\abs{x_\nu})^{M_\nu\wedge M_\nu'}}.
\end{equation*}
Taking the product of these bounds over all $\nu\in[k]$, we arrive at the claimed upper bound \eqref{eq:g*h1}. This finishes the proof of Lemma \ref{fj B1}.
\end{proof}

We can now prove that the pairing of molecules gives rise to an almost diagonal operator. The following lemma is a multi-parameter version of \cite[Lemma 3.7]{BHYY:1c}, which in turn is an extension of \cite[Corollary B.3]{FJ90} from specific molecular parameters to general ones.

\begin{lemma}\label{like B3}
Let
\begin{equation*}
  \vec K,\vec K'>\vn,\quad
  \vec M,\vec M'\in[0,\infty)^k,\quad
  \vec L,\vec L'\in\Z_{\geq-1}^k,\quad
  \vec N,\vec N'\in(0,\infty)^k.
\end{equation*}
Let $m_Q$ and $b_P$ be $(\vec K,\vec L,\vec M,\vec N)$ and $(\vec K',\vec L',\vec M',\vec N')$-molecules on $P,Q\in\mathscr D(\R^{\vn})$, respectively.

Then, for every $\vec\alpha\in(0,\infty)^k$, there exists a positive constant $C$, independent of $P$ and $Q$, such that
\begin{equation*}
  |\langle m_P,b_Q\rangle|\leq C
  b_{\vn;P,Q}^{\vD,\vE,\vF},
\end{equation*}
where
\begin{equation}\label{eq:B3}
  \vD=\vec M\wedge\vec M',\quad
  \begin{cases}
  \vE &=\frac{1}{2}\vn+\vec N'\wedge(\vec L\phantom{'}+\vone)\wedge(\vec K\phantom{'}-\vn-\vec\alpha) \\
  \vF &=\frac{1}{2}\vn+\vec N\phantom{'}\wedge(\vec L'+\vone)\wedge(\vec K'-\vn-\vec\alpha).
  \end{cases}
\end{equation}
\end{lemma}

\begin{proof}
Let $\vec\alpha\in(0,\infty)^k$. Note that each molecule with certain parameters is also a molecule with any smaller parameters. We first wish to replace $(\vec K,\vec L,\vec M,\vec N)$ and $(\vec K',\vec L',\vec M',\vec N')$ by some
$(\vec K,\vec L^\star,\vec M,\vec N^\star)$ and $(\vec K',\vec L^\dagger,\vec M',\vec N^\dagger)$,
respectively, such that the new parameters satisfy conditions \eqref{eq:K>N+n} and \eqref{eq:L=N-1} of Lemma \ref{fj B1}.

We are led to choose
\begin{equation}\label{eq:newNL}
  \begin{cases} \vec N^\star &:=
  \vec N\wedge(\vec L'+\vone)\wedge(\vec K'-\vn-\vec\alpha),\\
    \vec N^\dagger &:=\vec N'\wedge(\vec L+\vone)\wedge(\vec K-\vn-\vec\alpha),\end{cases}\qquad
   \begin{cases} \vec L^\star &:= \Floor{N^{\dagger}}, \\ \vec L^\dagger &:= \Floor{\vec N^\star}.\end{cases}
\end{equation}
It is immediate from these definitions that $\vec N^\star\leq\vN$ and $\vec N^\dagger\leq\vec N'$. Moreover,
\begin{equation*}
  \vec L^\star:=\Floor{\vec N^\dagger}\leq \Floor{\vec L+\vone}=\vec L
\end{equation*}
and symmetrically $\vec L^\dagger\leq\vec L'$. Hence, by monotonicity of the molecular conditions, the molecules $m_Q$ and $b_P$ are also molecules with parameters $(\vec K,\vec L^\star,\vec M,\vec N^\star)$ and $(\vec K',\vec L^\dagger,\vec M',\vec N^\dagger)$, respectively. The new parameters satisfy \eqref{eq:L=N-1} by definition, and \eqref{eq:K>N+n} follows from
\begin{equation*}
   \vec K'>\vec K'-\vec n-\vec\alpha+\vn\geq \vec N^\star+\vn
\end{equation*}
and a symmetric estimate for $\vec K$ and $\vec N^{\dagger}$.

Suppose next that $\vec\ell(P)=2^{-\vj}$ and $\vec\ell(Q)=2^{-\vec l}$ and let
\begin{equation*}
g(\cdot):=m_Q(c_Q-\cdot)
\text{ and }
h(\cdot):=b_P(c_P+\cdot)
\end{equation*}
which are molecules on rectangles of the same shape as $Q$ and $P$,
but both centered at the origin. Moreover, by a change of variables, we have
\begin{equation}\label{182}
\begin{split}
\langle m_Q,b_P\rangle
&=\int_{\mathbb R^{\vn}}g(c_Q-y)h(y-c_P)\,dy\\
&=\int_{\mathbb R^{\vn}}g(x)h(c_Q-c_P-x)\,dx
=(g*h)(c_Q-c_P).
\end{split}
\end{equation}
Let us immediately observe both $g$ and $h$ satisfy
the assumptions of Lemma \ref{fj B1} with parameters $(\vec K,\vec L^\star,\vec M,\vec N^\star)$ and $(\vec K',\vec L^\dagger,\vec M',\vec N^\dagger)$ in place of $(\vec K,\vec L,\vec M,\vec N)$ and $(\vec K',\vec L',\vec M',\vec N')$.
Thus, from \eqref{182} and Lemma \ref{fj B1}, we conclude that
\begin{align}\label{from B1}
|\langle m_P,b_Q\rangle|
\lesssim
2^{-(\vl-\vj)_+\cdot(\frac{\vn}{2}+\vec N^\star)-(\vj-\vl)_+\cdot(\frac{\vn}{2}+\vec N^{\dagger})}
\prod_{\nu=1}^{k}
\Big(1+2^{j_\nu \wedge l_\nu}|(c_Q-c_P)_\nu|\Big)^{-(M_\nu\wedge M_\nu')}.
\end{align}
Let us note that
\begin{equation*}
  2^{-(l_\nu-j_\nu)}=\frac{2^{-l_\nu}}{2^{-j_\nu}}=\frac{\ell(Q_\nu)}{\ell(P_\nu)},\qquad
  2^{-(l_\nu-j_\nu)_+}=\begin{cases} \ell(Q_\nu)/\ell(P_\nu), & \ell(Q_\nu)\leq\ell(P_\nu), \\ 1, & \text{otherwise},\end{cases}
\end{equation*}
and analogously with $2^{-(j_\nu-l_\nu)_+}$. Using these identities and comparing with \eqref{bDEF}, we recognise the right-hand side of \eqref{from B1} as
\begin{equation*}
  \mathrm{RHS}\eqref{from B1}=b_{\vn;P,Q}^{\vD,\vE,\vF},\quad
  \vD=\vec M\wedge\vec M',\quad
  \vE=\frac{\vn}{2}+\vec N^{\dagger},\quad
  \vF=\frac{\vn}{2}+\vec N^\star.
\end{equation*}
Substituting the values of $\vec N^\star$ and $\vec N^{\dagger}$ from \eqref{eq:newNL}, this finishes the proof of Lemma \ref{like B3}.
\end{proof}

\begin{theorem}\label{83}
Let $\vp\in(0,\infty)^k$, $\vq\in(0,\infty]^k$, $\vs\in\R^k$, $\tau\in[0,\infty)$, $r:=r(\pi,\vec p,\vec q)\wedge 1$,
$\pi\in S_{[2k]}$ be admissible for $(\vp,\vq)$,
and $V$ be a matrix weight such that $(\vp,\vq,\pi,V)$ belongs
to one of the main cases (Definition \ref{main cases}).
For both $i=1,2$, let $\{m^{(i)}_Q\}_{Q\in\mathscr D(\R^{\vn})}$ and $\{b^{(i)}_P\}_{P\in\mathscr D(\R^{\vn})}$ be families of $(\vec K,\vec L,\vec M,\vec N)$ and $(\vec K',\vec L',\vec M',\vec N')$-molecules, respectively, where the parameters
\begin{equation*}
  \vec K,\vec K',\vec M,\vec M'\in[0,\infty)^k,\quad
  \vec N,\vec N'\in(0,\infty)^k,\quad\vec L,\vec L'\in\Z_{\geq-1}^k
\end{equation*}
satisfy
\begin{align}\label{KMM}
\begin{cases} \vec K\phantom{'} > \vn+(\tau\vn+\vs)_+, \\ \vec K' >\vn+[(\tfrac1r-1)\vn-\vs]_+,\end{cases}
\qquad \vec M,\vec M'>(\tfrac{1}{r}+\tau)\vn,
\end{align}
and
\begin{align}\label{LNN}
\begin{cases} \vec L+\vone\phantom{'},\ \vec N' > \tau\vn+\vs, \\ \vec L'+\vone,\ \vec N \phantom{'}> (\frac1r-1)\vn-\vs.\end{cases}
\end{align}
Then
\begin{enumerate}[{\rm (i)}]
\item\label{it:eachAD} if $\varphi$ a Littlewood--Paley function on $\mathbb R^{\vec n}$, then each of the matrices
\begin{equation*}
  \Big\{ \Big\langle m_P^{(1)},b_Q^{(1)}\Big\rangle \Big\}_{P,Q\in\mathscr D(\mathbb R^{\vec n})},\quad
  \Big\{ \Big\langle m_P^{(1)},\varphi_Q\Big\rangle \Big\}_{P,Q\in\mathscr D(\mathbb R^{\vec n})},\quad
  \Big\{ \Big\langle \varphi_P,b_Q^{(1)}\Big\rangle \Big\}_{P,Q\in\mathscr D(\mathbb R^{\vec n})}
\end{equation*}
is $\dot a^{\vec s,\tau}_{\vec p,\vec q,\pi}$-almost diagonal (Definition \ref{astAD});

\item\label{it:convUC} if $t:=\{t_R\}_{R\in\mathscr D(\R^{\vn})}\in\dot a^{\vec s,\tau}_{\vec p,\vec q,\pi}(V)$, then
\begin{equation*}
s_P:=\sum_{Q,R\in\mathscr D(\R^{\vn})}
\Big\langle m_P^{(1)},b_Q^{(1)}\Big\rangle
\Big\langle m_Q^{(2)},b_R^{(2)}\Big\rangle t_R
\end{equation*}
converges unconditionally for each $P\in\mathscr D(\R^{\vn})$,
and $s:=\{s_P\}_{P\in\mathscr D(\R^{\vn})}$ satisfies
\begin{equation*}
\|s\|_{\atau(V)}
\leq C\|t\|_{\atau(V)},
\end{equation*}
where the positive constant $C$ is independent of $t$,
$\{m_Q\}_{Q\in\mathscr D(\R^{\vn})}$, and $\{b_Q\}_{Q\in\mathscr D(\R^{\vn})}$.
\end{enumerate}
\end{theorem}

\begin{proof}
\eqref{it:eachAD}: Since the Littlewood--Paley functions $\varphi_Q,\varphi_P$ are molecules with any parameters that we like, it is enough to consider the first matrix. Our plan is to apply Lemma \ref{like B3}. Note that the condition $\vec K,\vec K'>\vn$ of that lemma is immediate from assumption \eqref{KMM}. Hence, Lemma \ref{like B3} implies that
$\{ \pair{ m_P^{(1)} }{ b_Q^{(1)} }\}_{P,Q\in\mathscr D(\R^{\vn})}$ is $(\vD,\vE,\vF)$-almost diagonal with parameters $\vD,\vE,\vF$ as in \eqref{eq:B3}.

In order to be $\atau$-almost diagonal, by Definition \ref{astAD}, these parameters need to satisfy \eqref{DEF2}. Checking these conditions from the
assumptions \eqref{KMM} and \eqref{LNN} is straightforward.
This concludes the proof of part \eqref{it:eachAD}.

\eqref{it:convUC}:
By \eqref{it:eachAD}, we find that both
$\{|\langle m_P^{(1)}, b_Q^{(1)}\rangle|\}_{P,Q\in\mathscr D(\R^{\vn})}$
and $\{|\langle m_Q^{(2)}, b_R^{(2)}\rangle|\}_{Q,R\in\mathscr D(\R^{\vn})}$
are $\dot a^{\vec s,\tau}_{\vec p,\vec q,\pi}$-almost diagonal.
Then Lemma \ref{AcB} guarantees that the composition
\begin{equation*}
b_{P,R}:=\sum_{Q\in\mathscr D(\R^{\vn})}
\Big|\langle m_P^{(1)},b_Q^{(1)}\rangle\Big|
\Big|\langle m_Q^{(2)}, b_R^{(2)}\rangle\Big|
\end{equation*}
also defines an
$\dot a^{\vec s,\tau}_{\vec p,\vec q,\pi}$-almost diagonal
$B:=\{b_{P,R}\}_{P,R\in\mathscr D(\R^{\vn})}$.
Now Theorem \ref{AD st} implies the absolute convergence
\begin{equation*}
\sum_{Q,R\in\mathscr D(\R^{\vn})}
\Big|\Big\langle m_P^{(1)},b_Q^{(1)}\Big\rangle\Big|
\Big|\Big\langle m_Q^{(2)}, b_R^{(2)}\Big\rangle\Big|
|t_R|
=\sum_{R\in\mathscr D(\R^{\vn})}b_{P,R}|t_R|<\infty
\end{equation*}
and the claimed quasi-norm estimate for $s$.
This finishes the proof of \eqref{it:convUC} and hence Theorem \ref{83}.
\end{proof}

The following definition, motivated by Theorem \ref{83}, is analogous to \cite[Definition 3.10]{BHYY:1c}, which in turn is inspired by the classical definition \cite[pp.~56--57]{FJ90}.

\begin{definition}\label{def:AtauMole}
Let $\vp\in(0,\infty)^k$, $\vq\in(0,\infty]^k$, $\vs\in\R^k$,
$\tau\in[0,\infty)$, and $\pi\in S_{[2k]}$ be admissible for $(\vp,\vq)$.
Assume that $\vec K,\vec K',\vec M,\vec M'\in[0,\infty)^k,$ $\vec L,\vec L'\in\Z_{\geq-1}^k,$ and $\vec N,\vec N'\in(0,\infty)^k$.
Whenever the corresponding parameters satisfy \eqref{KMM} and \eqref{LNN}, then
\begin{enumerate}[\rm(i)]
  \item a $(\vec K,\vec L,\vec M,\vec N)$-molecule is called an \emph{$\dot A^{\vec s,\tau}_{\vec p,\vec q,\pi}$-analysis molecule};
  \item a $(\vec K',\vec L',\vec M',\vec N')$-molecule is called an \emph{$\dot A^{\vec s,\tau}_{\vec p,\vec q,\pi}$-synthesis molecule};
  \item a family of $(\vec K,\vec L,\vec M,\vec N)$-molecules is called a \emph{family of $\dot A^{\vec s,\tau}_{\vec p,\vec q,\pi}$-analysis molecules};
  \item a family of $(\vec K',\vec L',\vec M',\vec N')$-molecules is called a \emph{family of $\dot A^{\vec s,\tau}_{\vec p,\vec q,\pi}$-synthesis molecules}.
\end{enumerate}
\end{definition}

\begin{remark}\label{rem:83}
Theorem \ref{83} can be restated in the language of Definition \ref{def:AtauMole} by saying that conclusions \eqref{it:eachAD} and \eqref{it:convUC} hold whenever $\{m_P^{(i)}\}_{P\in\mathscr D(\R^{\vn})}$ and $\{b^{(i)}_Q\}_{Q\in\mathscr D(\R^{\vn})}$ are families of $\Atau$-analysis molecules and $\Atau$-synthesis molecules, respectively.
\end{remark}

If $f\in\dot A^{s,\tau}_{p,q}(V)$ and $\Phi$ is an analysis molecule,
then $f\in\mathscr{S}_0'(\R^{\vn};\C^m)$.
However, $\Phi$ need not belong to $\mathscr{S}_0(\R^{\vn})$,
and hence the pairing $\langle f,\Phi\rangle$ requires justification;
see \cite[p.\,155]{FJ90} or \cite[Lemma 5.7]{BH06}.

\begin{lemma}\label{88}
Let $\vp\in(0,\infty)^k$, $\vq\in(0,\infty]^k$, $\vs\in\R^k$, $\tau\in[0,\infty)$,
$\pi\in S_{[2k]}$ be admissible for $(\vp,\vq)$,
and $V$ be a matrix weight such that $(\vp,\vq,\pi,V)$ belongs
to one of the main cases (Definition \ref{main cases}).
Let $f\in\Atau(V)$ and $\Phi$ be an
$\dot A^{\vec s,\tau}_{\vec p,\vec q,\pi}$-analysis molecule on a cube $P\in\mathscr{D}(\R^{\vn})$.
Then, for every Littlewood--Paley pair $(\varphi,\psi)$ on $\R^{\vn}$,
\begin{equation}\label{881}
\langle f,\Phi\rangle
:=\sum_{R\in\mathscr{D}(\R^{\vn})}\langle f,\varphi_R\rangle\langle\psi_R,\Phi\rangle
\end{equation}
is well defined: the series above converges absolutely and its value is independent of the choices of $\varphi$ and $\psi$.
\end{lemma}

\begin{proof}
Suppose that $(\varphi^{(i)},\psi^{(i)})$, for $i=1,2$,
are two Littlewood--Paley pairs on $\R^{\vn}$.
Then
$$\{t_R\}_{R\in\mathscr D(\R^{\vn})}:=\Big\{\Big\langle f,\varphi_R^{(1)}\Big\rangle\Big\}_{R\in\mathscr D(\R^{\vn})}\in\atau(V)$$
by Corollary \ref{phiV}.
By Theorem \ref{83}\eqref{it:convUC} (or its reformulation in Remark \ref{rem:83}) applied to
this sequence, the $\Atau$-analysis molecules $m_P^{(1)}:=\Phi$ and $m^{(2)}_Q:=\varphi^{(2)}_Q$, and the $\Atau$-synthesis molecules $b_S^{(i)}:=\psi_S^{(3-i)}$, $S\in\{Q,R\}$, we find that
\begin{equation}\label{880}
\sum_{Q,R\in\mathscr D(\R^{\vn})}\Big\langle f,\varphi_R^{(1)}\Big\rangle
\Big\langle\psi_R^{(1)},\varphi_Q^{(2)}\Big\rangle
\Big\langle\psi_Q^{(2)},\Phi\Big\rangle
\end{equation}
converges absolutely. By Lemma \ref{cald2},
the claim \eqref{881} holds if either $\Phi$ or $f$
is replaced by a test function $\Psi\in\mathscr S_0(\R^{\vn})$.
Thus, if we sum up the double series \eqref{880} in one or the other order,
as justified by the absolute convergence that we just showed,
we then conclude that
\begin{align*}
\sum_{R\in\mathscr D(\R^{\vn})}\Big\langle f,\varphi_R^{(1)}\Big\rangle
\Big\langle\psi_R^{(1)},\Phi\Big\rangle
&=\sum_{R\in\mathscr D(\R^{\vn})}\Big\langle f,\varphi_R^{(1)}\Big\rangle
\Big(\sum_{Q\in\mathscr D(\R^{\vn})}
\Big\langle\psi_R^{(1)},\varphi_Q^{(2)}\Big\rangle
\Big\langle\psi_Q^{(2)},\Phi\Big\rangle\Big) \\
&=\sum_{Q\in\mathscr D(\R^{\vn})}\Big(\sum_{R\in\mathscr D(\R^{\vn})}
\Big\langle f,\varphi_R^{(1)}\Big\rangle
\Big\langle\psi_R^{(1)},\varphi_Q^{(2)}\Big\rangle\Big)
\Big\langle\psi_Q^{(2)},\Phi\Big\rangle\\
&=\sum_{Q\in\mathscr D(\R^{\vn})}\Big\langle f,\varphi_Q^{(2)}\Big\rangle
\Big\langle\psi_Q^{(2)},\Phi\Big\rangle.
\end{align*}
This proves both the absolute convergence of the right-hand side of \eqref{881}
and its independence of $\varphi$ and $\psi$,
which completes the proof of Lemma \ref{88}.
\end{proof}

Applying Theorem \ref{83}, Lemma \ref{88}, and a method pioneered by
Frazier and Jawerth in the proofs of \cite[Theorems 3.5 and 3.7]{FJ90},
we obtain the following conclusion.

\begin{theorem}\label{89}
Let $\vp\in(0,\infty)^k$, $\vq\in(0,\infty]^k$, $\vs\in\R^k$, $\tau\in[0,\infty)$,
$\pi\in S_{[2k]}$ be admissible for $(\vp,\vq)$,
and $V$ be a matrix weight such that $(\vp,\vq,\pi,V)$ belongs
to one of the main cases (Definition \ref{main cases}).
\begin{enumerate}[{\rm (i)}]
\item\label{890} If $\{m_Q\}_{Q\in\mathscr{D}(\mathbb R^{\vec n})}$
is a family of $\dot A^{\vec s,\tau}_{\vec p,\vec q,\pi}$-analysis molecules, each on the cube indicated by its subscript,
then there exists a positive constant $C$ such that,
for every $f\in\Atau(V)$,
$$
\Big\|\Big\{\Big\langle f,
m_Q\Big\rangle\Big\}_{Q\in\mathscr{D}(\mathbb R^{\vec n})}\Big\|_{\atau(V)}
\leq C\|f\|_{\Atau(V)}.
$$

\item\label{892}
If $\{b_Q\}_{Q\in\mathscr{D}(\mathbb R^{\vec n})}$
is a family of $\dot A^{\vec s,\tau}_{\vec p,\vec q,\pi}$-synthesis molecules,
each on the cube indicated by its subscript,
then, for every $t:=\{t_Q\}_{Q\in\mathscr{D}(\mathbb R^{\vec n})}
\in\atau(V)$, there exist $f\in\Atau(V)$ such that
$$
f=\sum_{Q\in\mathscr{D}(\mathbb R^{\vec n})} t_Qb_Q
$$
in $\mathscr{S}_0'(\mathbb R^{\vec n};\C^m)$ and a positive constant $C$,
independent of $\{t_Q\}_{Q\in\mathscr{D}(\mathbb R^{\vec n})}$
and $\{b_Q\}_{Q\in\mathscr{D}(\mathbb R^{\vec n})}$, such that
$\|f\|_{\Atau(V)}
\leq C\|t\|_{\atau(V)}.$
\end{enumerate}
\end{theorem}

\begin{proof}
We first show \eqref{890}.
Let $(\varphi,\psi)$ be a Littlewood--Paley pair on $\R^{\vn}$.
By Lemma \ref{88}, we find that, for every $Q\in\mathscr{D}(\R^{\vn})$,
\begin{equation}\label{91}
\langle f,m_Q\rangle
=\sum_{R\in\mathscr{D}(\R^{\vn})}
\langle f,\varphi_R\rangle
\langle\psi_R,m_Q\rangle
=\sum_{R\in\mathscr{D}(\R^{\vn})}
\langle\psi_R,m_Q\rangle
\Big(S_\varphi f\Big)_R.
\end{equation}
Using Theorem \ref{83} (or better, its reformulation in Remark \ref{rem:83}), we conclude that
$$B:=\{\langle\psi_R,m_Q\rangle\}_{Q,R\in\mathscr{D}(\R^{\vn})}$$
is an $\dot a^{\vec s,\tau}_{\vec p,\vec q,\pi}$-almost diagonal operator.
From this, we infer that
\begin{equation*}
\begin{split}
\Big\|\{\langle f,
m_Q\rangle\}_{Q\in\mathscr{D}(\mathbb R^{\vec n})}\Big\|_{\atau}
&=\Big\|B\Big(S_\varphi f\Big)\Big\|_{\atau}\qquad\text{by \eqref{91}} \\
&\lesssim\|S_\varphi f\|_{\atau}\qquad\text{by Theorem \ref{AD st}} \\
&\lesssim\| f\|_{\Atau}\qquad\text{by Corollary \ref{phiV}}.
\end{split}
\end{equation*}
This finishes the proof of \eqref{890}.

Now, we prove \eqref{892}. By Theorem \ref{83}\eqref{it:convUC} (or its reformulation in Remark \ref{rem:83}), with
$\phi\in\mathscr S_0(\R^{\vn})$ in place of the $\Atau$-synthesis molecule $m_P$,
we find the absolute convergence of
\begin{equation*}
\sum_{Q,R\in\mathscr D(\R^{\vn})} t_R\langle b_R,\varphi_Q\rangle\langle\psi_Q,\phi\rangle
=\sum_{R\in\mathscr D(\R^{\vn})} t_R\langle b_R,\phi\rangle
=:\langle f,\phi\rangle,
\end{equation*}
where the first identity follows from summing first over $Q\in\mathscr D(\R^{\vn})$
(as we can, by absolute convergence) and applying Lemma \ref{cald2}.
In particular, the sum of the series is a well defined element
$f\in\mathscr S_0'(\R^{\vn};\C^m)$. Taking $\phi=\varphi_P$,
we have $\langle  f,\varphi_P\rangle=(S_\varphi f)_P$, and we deduce
\begin{equation*}
\| f \|_{\Atau(V)}
\sim \|S_\varphi f\|_{\atau(V)}
\lesssim \|t\|_{\atau(V)}
\end{equation*}
from Corollary \ref{phiV} in the first step
and the quasi-norm estimate of Theorem \ref{83} in the second step.
This finishes the proof of Theorem \ref{89}.
\end{proof}

\begin{remark}\label{open vs D2}
Let $k=1$, $s\in\R$, $\tau\in[0,\infty)$, $p\in(0,\infty)$, $q\in(0,\infty]$,
$\pi\in S_{[2]}$, and $V\in\mathscr A_p(\R^n)$. If
\begin{equation}\label{coincase}
\begin{split}
\begin{cases}
    \tau=0 & \text{or}\quad \tau=\frac1p,\ q\in[p,\infty]\quad \text{or}\\
    \tau=\frac1p,\ q\in(0,p),\ \pi=F & \text{or}\quad \tau\in(\frac1p,\infty),
    \end{cases}
\end{split}
\end{equation}
then, using Theorem \ref{compare A} below, we conclude that $\dot A^{s,\tau}_{p,q,\pi}(V)
= \dot A^{s,\tau,\mathrm{cube}}_{p,q,\pi}(V)$.
In this case, by comparing Theorem \ref{89} in the present setting with
\cite[Theorem 3.17]{BHYY:1c}, we observe that the two results have
the same formal structure.
However, the admissible ranges of the relevant indices are not identical even in one-parameter case.

The narrower ranges arise from the boundedness theory of almost
diagonal operators used in the proof of Theorem \ref{89},
which requires stronger index
conditions than those in \cite{BHYY:1c} (see Remark \ref{rem AD st}).
In contrast, by restricting the smoothness parameter $N$
to the positive range $N>0$, several estimates can be carried out
more precisely, which leads to a relaxation of some index conditions.
For this reason, all results in the present article that depend
on almost diagonal operator bounds and molecular characterizations
involve slightly different index ranges from their one-parameter
counterparts.

On the other hand, for parameters outside the cases \eqref{coincase},
our present results deal with a scale of spaces that is different from the ones previously considered in the one-parameter case, as will be elaborated in Theorem \ref{compare A} below.
In the multi-parameter context, Theorem \ref{89} is completely new.
\end{remark}

As an application of Theorem \ref{89},
we obtain the following conclusion (see \cite[Proposition 3.18]{BHYY:1c} for details in the one-parameter case that can be easily adapted); we omit the details.

\begin{proposition}\label{6.18}
Let $\vp\in(0,\infty)^k$, $\vq\in(0,\infty]^k$, $\vs\in\R^k$, $\tau\in[0,\infty)$,
$\pi\in S_{[2k]}$ be admissible for $(\vp,\vq)$,
and $V$ be a matrix weight such that $(\vp,\vq,\pi,V)$ belongs
to one of the main cases (Definition \ref{main cases}).
Then $\mathscr{S}_0(\mathbb R^{\vec n};\C^m) \subset \dot A^{\vec s,\tau}_{\vec p,\vec q,\pi}(V)$.
Moreover, there exist $\vec M\in\N^k$
and a positive constant $C$ such that,
for every $f \in \mathscr{S}_0(\mathbb R^{\vec n};\C^m)$,
$$
\|f\|_{\dot A^{\vec s,\tau}_{\vec p,\vec q,\pi}(V)}
\leq C \sup_{\genfrac{}{}{0pt}{}{\gamma\in\N^{\vec n}}{|\gamma|_{\rm{vec}}\leq \vec M}}
\sup_{x\in\mathbb{R}^{\vec n}}|(\partial^{\gamma} f)(x)|
\prod_{\nu=1}^k (1+|x_\nu|)^{n+M_\nu+|\gamma_\nu|}.
$$
\end{proposition}

\subsection{Pseudo-differential operators} \label{pdo}

In this section, we follow the approach of
Grafakos and Torres used in \cite{GT99} and
establish the boundedness of multi-parameter
pseudo-differential operators on $\Atau(V)$.

Now, we recall the concepts of the class $\dot{S}_{1,1}^{\vec u}$
and its related pseudo-differential operators
(see, for instance, \cite[pp.\,261--263]{GT99}).

\begin{definition}\label{1:37}
Let $\vec u\in\mathbb N^k$.
The \emph{class $\dot{S}_{1,1}^{\vec u}$} is defined to be the set of all functions
$a\in C^\infty(\mathbb{R}^{\vn}\times(\mathbb{R}^{\vn}\setminus\{\mathbf 0\}))$
such that, for every $\alpha,\beta\in\mathbb N^{\vec n}$,
$$
\sup_{x\in\mathbb{R}^{\vn},\, \xi\in\mathbb{R}^{\vn}\setminus\{\mathbf{0}\}}
\prod_{\nu=1}^{k}|\xi_\nu|^{-u_\nu-|\alpha_\nu|+|\beta_\nu|}
\Big|\partial_x^\alpha\partial_\xi^\beta a(x,\xi)\Big|
<\infty.
$$
\end{definition}

\begin{definition}\label{def psDO}
Let $\vec u\in\mathbb N^k$ and $a\in\dot{S}_{1,1}^{\vu}$.
Define the \emph{multi-parameter pseudo-differential operator
$a(x,D)$} with symbol $a$ by setting,
for every $f\in \mathscr{S}_0(\R^{\vec n})$ and $x\in\mathbb{R}^{\vn}$,
\begin{equation*}
a(x,D)( f):=\int_{\mathbb{R}^{\vn}}e^{ix\cdot\xi}a(x,\xi)
\widehat{f}(\xi)\,d\xi.
\end{equation*}
\end{definition}

To present the main theorem in this section, we need the following technical lemmas.
Here, and thereafter, let $\mathcal R_{l,N}$ be an appropriate seminorm in $\mathscr S(\mathbb R^{\vn})$.

\begin{lemma}\label{1:00}
Let $f\in\mathscr S_0(\mathbb R^{\vn})$.
Then, for every $\vec M\in
\Z^k$ and every $\gamma\in\N^{\vn}$,
\begin{align*}
\sup_{\xi\in\R^{\vn}}\Big(\prod_{\nu=1}^{k}|\xi_\nu|^{M_\nu}\Big)\Big|
\partial^\gamma\widehat{f}(\xi)\Big|\lesssim\mathcal R_{l,N}(f),
\end{align*}
with $l,N\in\N$ large enough.
\end{lemma}

\begin{proof}
Let $\gamma\in\N^{\vn}$.
From \eqref{hat S0}
and the fact that the Fourier transform is an isomorphism on $\mathscr S(\R^{\vn})$,
it follows that, for every $M\in\mathbb N$ and $\xi\in\R^{\vn}$,
\begin{align*}
\Big|\partial^\gamma\widehat{f}(\xi)\Big|
\lesssim
|\xi_1|^M\mathcal R_{l,N}(f),
\end{align*}
where $\delta\in(0,1)$ and $l,N\in\N$ are large enough.
Since each dimension plays an equal role, we have,
for every $\nu\in[k]$ and every $M\in\N$,
\begin{align}\label{M-}
\sup_{\xi\in\R^{\vn}}|\xi_\nu|^{-M}\Big|\partial^\gamma\widehat{f}(\xi)\Big|
\lesssim\mathcal R_{l,N}(f),
\end{align}
where $l,N\in\N$ are large enough.
In addition, by the fact that the Fourier transform is
an isomorphism on $\mathscr S(\R^{\vn})$ again,
we conclude that, for every $\nu\in[k]$ and every $M\in\N$,
\begin{align}\label{M+}
\sup_{\xi\in\R^{\vn}}|\xi_\nu|^{M}
\Big|\partial^\gamma\widehat{f}(\xi)\Big|\leq
\sup_{\xi\in\R^{\vn}}|\xi|^{M}
\Big|\partial^\gamma\widehat{f}(\xi)\Big|\lesssim\mathcal R_{l,N}(f),
\end{align}
where $l,N\in\N$ are large enough.
Combining \eqref{M-}, \eqref{M+}, and the arbitrariness of $M$,
we find that,
for every $\vec M\in\Z^k$,
\begin{align*}
\sup_{\xi\in\R^{\vn}}\Big(\prod_{\nu=1}^{k}|\xi_\nu|^{M_{\nu}}
\Big)\Big|\partial^\gamma\widehat{f}(\xi)\Big|
\le\sup_{\xi\in\R^{\vn}}\sup_{\nu\in[k]}|\xi_\nu|^{kM_{\nu}}
\Big|\partial^\gamma\widehat{f}(\xi)\Big|\lesssim\mathcal R_{l,N}(f),
\end{align*}
where $l,N\in\N$ are large enough.
This finishes the proof of Lemma \ref{1:00}.
\end{proof}

\begin{lemma}\label{2223}
Let $\vu\in\mathbb N^k$ and $a\in\dot{S}_{1,1}^{\vu}$,
and let $a(x, D)$ be a multi-parameter pseudo-differential operator as defined in
Definition \ref{def psDO}. Then $a(x, D)$ is continuous from
$\mathscr{S}_0(\R^{\vec n})$ to $\mathscr{S}(\R^{\vn})$.
Consequently, its formal adjoint $a(x,D)^{\#}$, defined by setting,
for every $f\in\mathscr{S}'(\R^{\vn})$ and $\phi\in\mathscr{S}_0(\R^{\vec n})$,
$$
\Big\langle a(x,D)^{\#}f,\phi\Big\rangle
:=\langle f,a(x,D)\phi\rangle,
$$
is a continuous linear operator from $\mathscr{S}'(\R^{\vn})$ to $\mathscr{S}_0'(\R^{n})$.
\end{lemma}

\begin{proof}
Let $f\in\mathscr{S}_0(\R^{\vn})$ and $\Delta_\xi$ be the Laplace operator in the
variable $\xi$. Using the fact that
$$(I-\Delta_\xi)^Ne^{ix\cdot\xi}=(1+|x|^2)^Ne^{ix\cdot\xi}
$$
and integration by parts, we conclude that, for
every $\gamma \in \N^{\vn}$, $N\in\N$, and $x\in\R^{\vn}$,
\begin{align*}
\partial^\gamma a(x,D)f=\int_{\R^{\vn}}e^{ix\cdot\xi}
\frac{(I-\Delta_\xi)^N}{(1+|x|^2)^N}\Big[\partial^\gamma_x
a(x,\xi)\widehat{f}(\xi)\Big]\,d\xi.
\end{align*}
Applying Leibniz's rule, Definition \ref{1:37}, and Lemma \ref{1:00},
we find that, for every $\gamma \in \N^{\vn}$, $N\in\N$, and $x\in\R^{\vn}$,
\begin{align*}
\Big| (1+|x|^2)^N\partial^\gamma a(x,D)f \Big|
&\lesssim
\int_{\R^{\vn}}\sum_{|\beta|\le 2N}\Big|
\partial^\gamma_x\partial^\beta_\xi a(x,\xi)\Big|\sum_{|\beta|\le 2N}
\Big|\partial^\beta_\xi\widehat{f}(\xi)\Big|\,d\xi\\
&\lesssim\int_{\R^{\vn}}\sum_{|\beta|\le 2N}
\prod_{\nu=1}^{k}|\xi_\nu|^{u_\nu+|\gamma_\nu|-|\beta_\nu|}
\sum_{|\beta|\le 2N}
\Big|\partial^\beta_\xi\widehat{f}(\xi)\Big|\,d\xi
\lesssim\mathcal R_{l,N}(f),
\end{align*}
where $l,N\in\N$ are large enough.
This shows $a(x,D)f\in\mathscr S(\R^{\vn})$ and
$a(x,D)$ is continuous. This finishes the proof of Lemma \ref{2223}.
\end{proof}

Applying an argument similar to that used in the proof of
\cite[Theorem 1.5]{SYY10} (see also the proofs of \cite[Theorems 1.1 and 1.2]{GT99}),
we obtain the following theorem.
For the convenience of the reader, we give the details of its proof here.

\begin{theorem}\label{pseudo}
Let $\vp\in(0,\infty)^k$, $\vq\in(0,\infty]^k$, $\vs\in\R^k$, $\tau\in[0,\infty)$,
$r:=r(\pi,\vp,\vq)\wedge 1$,
$\pi\in S_{[2k]}$ be admissible for $(\vp,\vq)$,
and $V$ be a matrix weight such that $(\vp,\vq,\pi,V)$ belongs
to one of the main cases (Definition \ref{main cases}).
Assume that $\vec L:=\lfloor \frac{\vec n}r-\vec n-\vec s\rfloor$.
Let $\vu\in\mathbb{N}^k$, $a\in\dot{S}_{1,1}^{\vu}$,
and $a(x,D)$ be a multi-parameter pseudo-differential operator with symbol $a$.
Assume that its formal adjoint $a(x,D)^{\#}$ has the property,
for every $\gamma\in\N^{\vn}$ satisfying $|\gamma|_{\rm{vec}}\leq \vec L$, that
\begin{equation}\label{227}
a(x,D)^{\#}(x^{\gamma})=0
\in \mathscr{S}_0'(\R^{\vec n}).
\end{equation}
Then $a(x,D)$ can be extended to a continuous linear mapping
from $\dot A_{\vp,\vq,\pi}^{\vs+\vu,\tau}(V)$ to $\Atau(V)$.
\end{theorem}

\begin{proof}
To simplify the presentation,
in this proof, we denote $a(x,D)$ simply by $T$.
Let $f\in\dot A_{\vp,\vq,\pi}^{\vs+\vu,\tau}(V)$
and $(\varphi,\varphi)$ be a Littlewood--Paley pair on $\R^{\vn}$.
By this and Lemma \ref{cald2}, we find that
$$
f=\sum_{Q\in\mathscr D(\R^{\vn})}
\langle f,\varphi_Q\rangle\varphi_Q
$$
in $\mathscr{S}_0'(\R^{\vn};\C^m)$, and hence it is natural to define
\begin{equation}\label{a(x,D)f}
\widetilde T(f)
:=\sum_{Q\in\mathscr{D}(\mathbb R^{\vec n})}
\langle f,\varphi_Q\rangle T(\varphi_Q)
=\sum_{Q\in\mathscr{D}(\mathbb R^{\vec n})}
\Big[2^{\vj_Q \cdot \vu} \Big(S_\varphi f\Big)_Q\Big]
\Big[2^{-(\vj_Q \cdot \vu)} T(\varphi_Q)\Big]
\end{equation}
in $\mathscr{S}_0'(\R^{\vn};\C^m)$.
If $f\in\mathscr S_0(\R^{\vn};\C^m)$,
then Lemmas \ref{cald2} and \ref{2223} imply that
the series in \eqref{a(x,D)f} converges in $\mathscr S_0(\R^{\vn};\C^m)$ to $T(f)$.
Therefore, $\widetilde T$ and $T$ coincide on $\mathscr S_0(\R^{\vn};\C^m)$.
Now, we show that $\widetilde T(f)$ is well defined when $f\in\mathscr{S}_0'(\R^{\vn};\C^m)$.
From Corollary \ref{phiV}, we deduce that
\begin{equation}\label{102}
\Big\|\Big\{ 2^{\vj_Q \cdot \vu}
\Big(S_\varphi f\Big)_Q
\Big\}_{Q\in\mathscr{D}(\mathbb R^{\vec n})}\Big\|_{\dot a_{\vp,\vq,\pi}^{\vs,\tau}(V)}
=\Big\|S_\varphi f\Big\|_{\dot a_{\vp,\vq,\pi}^{\vs+\vu,\tau}(V)}
\sim \|f\|_{\dot A_{\vp,\vq,\pi}^{\vs+\vu,\tau}(V)}.
\end{equation}
Using a straightforward multi-parameter extension of \cite[pp.\,265--266]{GT99},
we conclude that,
for every $\gamma\in\mathbb N^{\vn}$, $\vec N\in\mathbb N^k$, and $Q\in\mathscr{D}(\R^{\vn})$,
\begin{align}\label{scondition}
\abs{\partial^{\gamma}(T\varphi_Q)(x)}\lesssim\prod_{\nu=1}^{k}|Q_\nu|^{-\frac{u_\nu}{n_\nu}
-\frac12-\frac{|\gamma_\nu|}{n_\nu}}\Big[1+\ell(Q_\nu)^{-1}
|x_\nu-c_{Q_\nu}|\Big]^{-N_\nu}.
\end{align}
Now, we verify the vanishing moment of $T(\varphi_Q)$.
From \eqref{227}, we infer that, for every $\gamma\in\N^{\vn}$
satisfying $|\gamma|_{\rm{vec}}\leq \vec L$,
\begin{equation*}
\int_{\mathbb R^{\vec n}}x^{\gamma}
T(\varphi_Q)(x)\,dx
=\Big\langle
T^{\#}(x^{\gamma}),
\varphi_Q \Big\rangle
=0.
\end{equation*}
This, together with \eqref{scondition}, further implies that
$2^{-(\vj_Q \cdot \vu)}T(\varphi_Q)$ is a constant multiple of
a $(\vec K,\vec L,\vec M,\vec N)$-molecule on $Q$ for all $\vec K,\vec M,\vec N$
and hence a constant multiple of an $\dot A^{\vec s,\tau}_{\vec p,\vec q,\pi}$-synthesis molecule (Definition \ref{def:AtauMole}).
From this and Theorem \ref{89}\eqref{892}, we deduce that
the series in \eqref{a(x,D)f} is convergent in $\mathscr{S}_0'(\mathbb R^{\vec n};\C^m)$.
Using Theorem \ref{89}\eqref{892} again, we conclude that
\begin{equation*}
\begin{split}
\|T(f)\|_{\dot A_{\vp,\vq,\pi}^{\vs,\tau}(V)}
&\lesssim\Big\|\Big\{ 2^{\vj_Q \cdot \vu} \Big(S_\varphi f\Big)_Q
\Big\}_{Q\in\mathscr{D}(\mathbb R^{\vec n})}\Big\|_{\dot a_{\vp,\vq,\pi}^{\vs,\tau}(V)}\\
&\sim \|f\|_{\dot A_{\vp,\vq,\pi}^{\vs+\vu,\tau}(V)}\quad\text{by \eqref{102}}.
\end{split}
\end{equation*}
This finishes the proof of Theorem \ref{pseudo}.
\end{proof}

\begin{remark}
\begin{enumerate}[{\rm (i)}]
\item Theorem \ref{pseudo} is new in the multi-parameter case. In the one-parameter case, it essentially agrees with \cite[Theorem 3.23]{BHYY:1c} with a slightly more restrictive range of indices for the usual Besov and Triebel--Lizorkin spaces with $\tau=0$ and more generally whenever the parameters are such that $\Atau$ spaces defined by open sets agree with the $\AtauCube$ spaces defined by cubes (See Theorem \ref{compare A}); otherwise, it deals with a different scale of spaces.

\item In \cite{Xu22}, Xu established the boundedness of
pseudo-differential operators on inhomogeneous bi-parameter
Besov--Triebel--Lizorkin spaces for a different symbol class.
\end{enumerate}
\end{remark}

\subsection{Wavelet characterisations}
\label{wavelet characterization}

In this section, we establish wavelet characterisation of
$\Atau(V)$ using either the Meyer wavelet or the Daubechies wavelet.
Let us first recall the wavelet constructed by Lemari\'e and Meyer
(see, for instance, \cite[Theorem 2]{lm86}).

\begin{lemma}\label{wavelet basis}
There exists a sequence $\{\theta^{(i)}\}_{i=1}^{2^n-1}\subset\mathscr{S}_0(\R^n)$
of real-valued basic wavelets such that
$$
\Big\{\theta^{(i)}_Q:\ i\in [2^n-1],\ Q\in\mathscr{D}(\R^n)\Big\}
$$
forms an orthonormal basis of $L^2(\R^n)$.
\end{lemma}

Here and below, for any $\vec a:=(a_1,\ldots,a_k)\in\mathbb Z^k_+$,
let $[\vec a]:=[a_1]\times\cdots\times[a_k]$.
For every $\nu\in[k]$,
there exist Meyer wavelets
$\{\theta^{(i_\nu)}\}_{i_\nu=1}^{2^{n_\nu}-1}$ on $\R^{n_\nu}$.
For every $\vi\in[2^{\vn}-\vec1]$
and $x\in \mathbb R^{\vec n}$, let
$$
\theta^{(\vi)}(x)
:=\prod_{\nu=1}^k \theta^{(i_\nu)}(x_\nu)
$$
be the \emph{Meyer wavelet on $\R^{\vn}$}.

\begin{remark}\label{121}
In Lemma \ref{wavelet basis}, for every $f\in\mathscr{S}_0'(\R^{\vn})$,
\begin{equation*}
f=\sum_{\vi\in[2^{\vn}-\vec1]}\sum_{Q\in\mathscr{D}(\R^{\vn})}
\Big\langle f,\theta^{(\vi)}_Q\Big\rangle\theta^{(\vi)}_Q
\end{equation*}
in $\mathscr{S}_0'(\R^{\vn})$
(see the proof of \cite[Theorem 7.20]{FJW91} for $k=1$ and
\cite[Proposition 6.2]{GKP21} for $k=2$;
the same argument as that used in \cite[Proposition 6.2]{GKP21} can extend to all $k\geq 3$ with slight modifications).
\end{remark}

Now, we establish the wavelet characterization of
$\dot A^{\vec s,\tau}_{\vec p,\vec q,\pi}(V)$
via Meyer wavelets.

\begin{theorem}\label{wavelet}
Let $\vp\in(0,\infty)^k$, $\vq\in(0,\infty]^k$, $\vs\in\R^k$, $\tau\in[0,\infty)$,
$\pi\in S_{[2k]}$ be admissible for $(\vp,\vq)$,
and $V$ be a matrix weight such that $(\vp,\vq,\pi,V)$ belongs
to one of the main cases (Definition \ref{main cases}).
Let $\{\theta^{(\vi)}\}_{\vi\in[2^{\vn}-\vec1]}$ be Meyer wavelets on $\R^{\vn}$.
Then $f\in\dot A^{\vec s,\tau}_{\vec p,\vec q,\pi}(V)$ if and only if $f\in\mathscr{S}_0'(\R^{\vn};\C^m)$ and
$$
\|f\|_{\dot A^{\vec s,\tau}_{\vec p,\vec q,\pi}(V)_\mathrm{w}}
:=\sum_{\vi\in[2^{\vn}-\vec1]}
\Big\|\Big\{\Big\langle f,\theta^{(\vi)}_Q \Big\rangle
\Big\}_{Q\in\mathscr D(\R^{\vn})}\Big\|_{\dot a^{\vec s,\tau}_{\vec p,\vec q,\pi}(V)}
<\infty.
$$
Moreover, for every $f\in\mathscr{S}_0'(\R^{\vn};\C^m)$,
$\|f\|_{\dot A^{\vec s,\tau}_{\vec p,\vec q,\pi}(V)}
\sim\|f\|_{\dot A^{\vec s,\tau}_{\vec p,\vec q,\pi}(V)_\mathrm{w}},$
where the positive equivalence constants are independent of $f$.
\end{theorem}

\begin{proof}
For any $\vi\in[2^{\vn}-\vec1]$, the estimate
$\|\{\langle f,\theta^{(\vi)}_Q\rangle\}_{Q\in\mathscr{D}(\mathbb R^{\vec n})}\|_{\dot a^{\vec s,\tau}_{\vec p,\vec q,\pi}(V)}
\lesssim\| f\|_{\dot A^{\vec s,\tau}_{\vec p,\vec q,\pi}(V)}$
is immediate from Theorem \ref{89}\eqref{890} and the fact that
$\{\theta^{(\vi)}_Q\}_{Q\in\D(\R^{\vn})}$ is a
constant multiple of a family of
$\dot A^{\vec s,\tau}_{\vec p,\vec q,\pi}$-analysis molecules.
Thus, $\|f\|_{\dot A^{\vec s,\tau}_{\vec p,\vec q,\pi}(V)_\mathrm{w}}
\lesssim \| f\|_{\dot A^{\vec s,\tau}_{\vec p,\vec q,\pi}(V)}$
follows from summing over all $\vi\in[2^{\vn}-\vec1]$.

Next, we show that, for every $f\in\mathscr{S}_0'(\R^{\vn};\C^m)$
with $\|f\|_{\dot A^{\vec s,\tau}_{\vec p,\vec q,\pi}(V)_\mathrm{w}}<\infty$,
\begin{equation}\label{120}
\|f\|_{\dot A^{\vec s,\tau}_{\vec p,\vec q,\pi}(V)}
\lesssim\|f\|_{\dot A^{\vec s,\tau}_{\vec p,\vec q,\pi}(V)_\mathrm{w}}.
\end{equation}
For such $f$ and every $\vi\in[2^{\vn}-\vec1]$,
from the definition of $\|\cdot\|_{\dot A_{p,q}^{s,\tau}(W)_\mathrm{w}}$,
it follows that
\begin{equation*}
t^{(\vi)}
:=\Big\{t^{(\vi)}_Q \Big\}_{Q\in\D(\R^{\vn})}
:=\Big\{\Big\langle f,
\theta^{(\vi)}_Q\Big\rangle\Big\}_{Q\in\D(\R^{\vn})}
\in\dot a^{\vec s,\tau}_{\vec p,\vec q,\pi}(V).
\end{equation*}
By Remark \ref{121}, we obtain
\begin{equation}\label{195}
\|f\|_{\dot A^{\vec s,\tau}_{\vec p,\vec q,\pi}(V)}
=\Big\|\sum_{\vi\in[2^{\vn}-\vec1]}\sum_{Q\in\mathscr{D}(\R^{\vn})}
t^{(\vi)}_Q \theta^{(\vi)}_Q\Big\|_{\dot A^{\vec s,\tau}_{\vec p,\vec q,\pi}(V)}
\lesssim\sum_{\vi\in[2^{\vn}-\vec1]}\Big\|\sum_{Q\in\mathscr{D}(\R^{\vn})}
t^{(\vi)}_Q \theta^{(\vi)}_Q\Big\|_{\dot A^{\vec s,\tau}_{\vec p,\vec q,\pi}(V)}.
\end{equation}
Since $\theta^{(i_\nu)}\in \mathscr{S}_0(\R^{n_\nu})$ for all $\nu\in[k]$,
it follows that $\{\theta^{(\vi)}_Q\}_{Q\in\mathscr{D}(\R^{\vn})}$ is a constant multiple of
a family of $\dot A^{\vec s,\tau}_{\vec p,\vec q,\pi}$-synthesis molecules,
each on the cube indicated by its subscript.
From this, the fact that $t^{(\vi)}\in\dot a^{\vec s,\tau}_{\vec p,\vec q,\pi}(V)$, Theorem \ref{89}(ii),
and \eqref{195}, we deduce that
$$
\|f\|_{\dot A^{\vec s,\tau}_{\vec p,\vec q,\pi}(V)}
\lesssim\sum_{\vi\in[2^{\vn}-\vec1]}\big\|t^{(\vi)}\big\|_{\dot a^{\vec s,\tau}_{\vec p,\vec q,\pi}(V)}
=\|f\|_{\dot A^{\vec s,\tau}_{\vec p,\vec q,\pi}(V)_\mathrm{w}},
$$
which completes the proof of \eqref{120}
and hence Theorem \ref{wavelet}.
\end{proof}

\begin{remark}
Let $k=2$, $\vec p:=p\cdot\vec 1$ for some $p\in(0,\infty)$,
$\vec q:=q\cdot\vec 1$ for some $q\in(0,\infty]$,
$\vec s\in\R^2$, $\tau=0$, $m=1$, and $V\equiv 1$.
Further assume that $(L^{\vp}\ell^{\vq})_\pi =\ell^qL^p$ (Besov case) or
$(L^{\vp}\ell^{\vq})_\pi =L^p\ell^q$ (Triebel--Lizorkin case).
Then Theorem \ref{wavelet} in this case coincides with \cite[Theorem 6.4]{GKP21}.
\end{remark}

In what follows, for any $\mathscr{N}\in\mathbb{N}$,
we use $C^{\mathscr{N}}(\mathbb R^n)$ to denote the set of all
$\mathscr{N}$ times continuously differentiable functions on $\mathbb{R}^n$.

\begin{definition}\label{def Dau}
Let $\mathscr{N}\in\mathbb{N}$. We say that $\{\theta^{(i)}\}_{i=1}^{2^n-1}$
are \emph{Daubechies wavelets} of class $C^{\mathscr N}(\mathbb R^n)$ if each
$\theta^{(i)}\in C^{\mathscr N}(\mathbb R^n)$ is real-valued with compact support and if
$$
\Big\{\theta^{(i)}_Q :\ i\in [2^n-1],\ Q\in\mathscr{D}(\mathbb R^n)\Big\}
$$
is an orthonormal basis of $L^2(\R^n)$.
\end{definition}

The following wavelet basis was constructed by Daubechies
(see, for instance, \cite{D88}).

\begin{lemma}\label{wavelet basis 2}
For any $\mathscr{N}\in\mathbb{N}$,
there exist Daubechies wavelets of class $C^{\mathscr N}(\mathbb R^n)$.
\end{lemma}

\begin{remark}\label{Dau can}
Let $\{\theta^{(i)}\}_{i\in[2^n-1]}$ be Daubechies wavelets of class $C^{\mathscr N}(\mathbb R^n)$.
By \cite[Corollary 5.5.2]{D92}, we find that,
for any $\alpha\in\mathbb{N}^n$ with $|\alpha|\leq\mathscr{N}$,
the following cancellation conditions are satisfied:
$$
\int_{\mathbb{R}^n}x^\alpha\theta^{(i)}(x)\,dx=0.
$$
\end{remark}

For every $\vec{\mathscr N}\in\mathbb{N}^k$ and $\nu\in[k]$,
there exist Daubechies wavelets
$\{\theta^{(i_\nu)}\}_{i_\nu=1}^{2^{n_\nu}-1}$ of class $C^{\mathscr N_\nu}(\mathbb R^{n_\nu})$.
For every $\vi\in[2^{\vn}-\vec1]$
and $x\in \mathbb R^{\vec n}$, let
$$
\theta^{(\vi)}(x)
:=\prod_{\nu=1}^k \theta^{(i_\nu)}(x_\nu)
$$
be the Daubechies wavelet of class $C^{\vec{\mathscr N}}(\mathbb R^{\vn})$.
Then, for every $\alpha\in\mathbb N^{\vec n}$ with $|\alpha|_{\rm{vec}}\leq\vec{\mathscr N}$,
\begin{equation}\label{cancellation}
\int_{\mathbb{R}^{n_\nu}}x_\nu^{\alpha_\nu}\theta^{\vi}(x_\nu,x')\,dx_\nu=0
\quad\text{for a.e.}\quad x'\in\R^{\vn}\ominus\R^{n_\nu},
\end{equation}
as required for the molecular cancellation condition \eqref{eq:moleL}.

\begin{corollary}\label{88 corollary}
Let $\vp\in(0,\infty)^k$, $\vq\in(0,\infty]^k$, $\vs\in\R^k$, $\tau\in[0,\infty)$,
$r:=r(\pi,\vp,\vq)\wedge 1$,
$\pi\in S_{[2k]}$ be admissible for $(\vp,\vq)$,
and $V$ be a matrix weight such that $(\vp,\vq,\pi,V)$ belongs
to one of the main cases (Definition \ref{main cases}).
Let $\vec{\mathscr N}\in\mathbb{N}^k$ satisfy
\begin{equation}\label{NDau}
\vec{\mathscr N}>(\vec n\tau+\vec s-\vone)\vee \Big(\frac{\vec n}r-\vec n-\vec s\Big)
\end{equation}
and
$\{\theta^{(\vi)}\}_{\vi\in[2^{\vn}-\vec1]}$
be Daubechies wavelets of class
$C^{\vec{\mathscr N}}(\mathbb R^{\vn})$. Then
\begin{enumerate}[\rm(i)]
\item\label{it:theta=AtauMole}
For each $\vi\in[2^{\vn}-\vec1]$ and $Q\in\D(\R^{\vn})$,
the wavelet $\theta^{(\vi)}_Q$ is a constant multiple of an
$\dot A^{\vec s,\tau}_{\vec p,\vec q,\pi}$-analysis molecule on $Q$,
where the constant is independent of $Q\in\D(\R^{\vn})$.

\item\label{it:pair-f-theta}
If $(\varphi,\psi)$ is a Littlewood--Paley pair on $\R^{\vn}$,
then, for each $\vi\in[2^{\vn}-\vec1]$ and $Q\in\mathscr{D}(\mathbb R^{\vec n})$,
$$
\Big\langle f,\theta^{(\vi)}_Q\Big\rangle
:=\sum_{R\in\D(\R^{\vn})}
\langle f,\varphi_R\rangle
\Big\langle\psi_R,\theta^{(\vi)}_Q\Big\rangle
$$
is well defined: the series converges absolutely
and its value is independent of the choice of $\varphi$ and $\psi$.
\end{enumerate}
\end{corollary}

\begin{proof}
From the fact that $\theta^{(\vi)}$ has compact support
and \eqref{cancellation},
it follows that the wavelet $\theta^{(\vi)}_Q$ is a constant multiple of a $(\vec K,\vec{\mathscr N},\vec M,\vec{\mathscr N})$-molecule on $Q$ (see Definition \ref{moleKLMN})
for all $\vec K,\vec M\in[0,\infty)^k$.
This, together with \eqref{NDau}, further implies that
$\theta^{(\vi)}_Q$ is a constant multiple of an $\dot A^{\vec s,\tau}_{\vec p,\vec q,\pi}$-analysis molecule on $Q$ (see Definition \ref{def:AtauMole}).
This finishes the proof of \eqref{it:theta=AtauMole}.

By just proven \eqref{it:theta=AtauMole} and Lemma \ref{88}, we obtain \eqref{it:pair-f-theta},
which completes the proof of Corollary \ref{88 corollary}.
\end{proof}

Next, we can establish the wavelet decomposition of $\dot A^{\vec s,\tau}_{\vec p,\vec q,\pi}(V)$
via Daubechies wavelets.

\begin{theorem}\label{wavelet 2}
Let $\vp\in(0,\infty)^k$, $\vq\in(0,\infty]^k$, $\vs\in\R^k$, $\tau\in[0,\infty)$,
$\pi\in S_{[2k]}$ be admissible for $(\vp,\vq)$,
and $V$ be a matrix weight such that $(\vp,\vq,\pi,V)$ belongs
to one of the main cases (Definition \ref{main cases}).
Let $\vec{\mathscr N}\in\mathbb{N}^k$ satisfy \eqref{NDau} and
$\{\theta^{(\vi)}\}_{\vi\in[2^{\vn}-\vec1]}$
be Daubechies wavelets of class
$C^{\vec{\mathscr N}}(\mathbb R^{\vn})$.
Then, for any $f\in \dot A^{\vec s,\tau}_{\vec p,\vec q,\pi}(V)$,
\begin{equation}\label{204}
f= \sum_{\vi\in[2^{\vn}-\vec1]}
\sum_{Q\in\D(\R^{\vn})}
\Big\langle f,\theta^{(\vi)}_Q\Big\rangle\theta^{(\vi)}_Q
\end{equation}
in $\mathscr{S}_0'(\mathbb R^{\vn};\mathbb C^m)$ and
\begin{equation}\label{2160}
\|f\|_{\dot A^{\vec s,\tau}_{\vec p,\vec q,\pi}(V)}
\sim\|f\|_{\dot A^{\vec s,\tau}_{\vec p,\vec q,\pi}(V)_\mathrm{w}},
\end{equation}
where the positive equivalence constants are independent of $f$.
\end{theorem}

\begin{proof}
Let $(\varphi,\psi)$ be a Littlewood--Paley pair on $\R^{\vn}$.
For any $\phi\in\mathscr{S}_0(\mathbb R^{\vec n})$, we consider the series
\begin{align}\label{222}
S:= \sum_{\vi\in[2^{\vn}-\vec1]}
\sum_{Q,R\in\mathscr{D}(\mathbb R^{\vec n})}
\langle f,\varphi_R\rangle
\Big\langle\psi_R,\theta^{(\vi)}_Q\Big\rangle
\Big\langle\theta^{(\vi)}_Q,\phi\Big\rangle.
\end{align}
Let $t:=\{t_R\}_{R\in\mathscr{D}(\mathbb R^{\vec n})}
:=\{\langle f,\varphi_R\rangle\}_{R\in\mathscr{D}(\mathbb R^{\vec n})}$.
By Corollary \ref{phiV}, we find that $t\in \dot a^{\vec s,\tau}_{\vec p,\vec q,\pi}(V)$.
Moreover, both $\psi_R$ and $\theta_Q^{(\vi)}$ are constant multiples of
$\dot A^{\vec s,\tau}_{\vec p,\vec q,\pi}$-analysis
and $\dot A^{\vec s,\tau}_{\vec p,\vec q,\pi}$-synthesis molecules
associated with the cubes specified by their subscripts.
In addition, $\phi\in\mathscr{S}_0(\mathbb R^{\vec n})$ is
a constant multiple of an $\dot A^{\vec s,\tau}_{\vec p,\vec q,\pi}$-synthesis
molecule. Consequently, the absolute convergence of \eqref{222}
follows directly from Corollary~\ref{83}\eqref{it:convUC}.

If we first sum over $Q$ and $\vi$ and use the fact that
$$
\Big\{\theta_Q^{(\vi)}:\ Q\in\mathscr D(\mathbb R^{\vec n}),\
\vi\in[2^{\vn}-\vec1]\Big\}
$$
is an orthonormal basis of $L^2(\mathbb R^{\vec n})$, we obtain
\begin{align*}
S&=\sum_{R\in\mathscr D(\mathbb R^{\vec n})}
\langle f,\varphi_R\rangle
\sum_{\vi\in[2^{\vn}-\vec1]}
\sum_{Q\in\mathscr D(\mathbb R^{\vec n})}
\Big\langle\psi_R,\theta^{(\vi)}_Q\Big\rangle
\Big\langle\theta^{(\vi)}_Q,\phi\Big\rangle\\
&=\sum_{R\in\mathscr D(\mathbb R^{\vec n})}
\langle f,\varphi_R\rangle
\langle\psi_R,\phi\rangle
=\langle f,\phi\rangle,
\end{align*}
where the last equality used the Calder\'on reproducing formula
Lemma \ref{cald2}.
On the other hand, if we first sum over $R$,
then we recognise the definition of the pairing of $\dot A^{\vec s,\tau}_{\vec p,\vec q,\pi}(V)$
with an $\dot A^{\vec s,\tau}_{\vec p,\vec q,\pi}$-analysis molecule $\theta^{(\vi)}_Q$ to the result that
\begin{align*}
S
&=\sum_{\vi\in[2^{\vn}-\vec1]}
\sum_{Q\in\mathscr D(\mathbb R^{\vec n})}
\Big(\sum_{R\in\mathscr D(\mathbb R^{\vec n})}
\langle f,\varphi_R\rangle  \Big\langle\psi_R,\theta^{(\vi)}_Q\Big\rangle\Big)
\Big\langle\theta^{(\vi)}_Q,\phi\Big\rangle\\
&=\sum_{\vi\in[2^{\vn}-\vec1]}
\sum_{Q\in\mathscr D(\mathbb R^{\vec n})}
\Big\langle f,\theta^{(\vi)}_Q\Big\rangle  \Big\langle\theta^{(\vi)}_Q,\phi\Big\rangle.
\end{align*}
A combination of the last two displays proves \eqref{204} in the sense of
convergence in $\mathscr{S}_0'(\mathbb R^n;\mathbb C^m)$.

The estimate
$\|\{\langle f,\theta^{(\vi)}_Q\rangle\}_{Q\in\mathscr{D}(\mathbb R^{\vec n})}\|_{\dot a^{\vec s,\tau}_{\vec p,\vec q,\pi}(V)}
\lesssim\| f\|_{\dot A^{\vec s,\tau}_{\vec p,\vec q,\pi}(V)}$
follows directly from Theorem~\ref{89}\eqref{890}, together with the observation that
$\{\theta^{(\vi)}_Q\}_{Q\in\mathscr{D}(\mathbb R^{\vec n})}$ are
constant multiples of
$\dot A^{\vec s,\tau}_{\vec p,\vec q,\pi}$-analysis molecules.
Consequently, by summing over $\vi\in[2^{\vn}-\vec1]$, we obtain
$\|f\|_{\dot A^{\vec s,\tau}_{\vec p,\vec q,\pi}(V)_\mathrm{w}}
\lesssim \| f\|_{\dot A^{\vec s,\tau}_{\vec p,\vec q,\pi}(V)}$.

To establish the reverse inequality, we note that $\{\langle f,\theta^{(\vi)}_Q\rangle\}_{Q\in\mathscr{D}(\mathbb R^{\vec n})}
\in\dot a^{\vec s,\tau}_{\vec p,\vec q,\pi}(V)$. Since
$\{\theta^{(\vi)}_Q\}_{Q\in\mathscr{D}(\mathbb R^{\vec n})}$
are constant multiples of $\dot A^{\vec s,\tau}_{\vec p,\vec q,\pi}$-synthesis molecules,
Theorem \ref{89}\eqref{892} guarantees that the series
$f^{(\vi)}:=\sum_{Q\in\mathscr{D}(\mathbb R^{\vec n})}
\langle f,\theta^{(\vi)}_Q\rangle \theta^{(\vi)}_Q$
converges in $\mathscr S_0'(\mathbb R^{\vec n};\mathbb C^m)$ and satisfies
\begin{equation*}
\big\|f^{(\vi)}\big\|_{\dot A^{\vec s,\tau}_{\vec p,\vec q,\pi}(V)}
\lesssim \Big\|\Big\{\Big\langle f,\theta^{(\vi)}_Q \Big\rangle\Big\}_{Q\in\mathscr{D}(\mathbb R^{\vec n})}\Big\|_{\dot a^{\vec s,\tau}_{\vec p,\vec q,\pi}(V)}.
\end{equation*}
Using \eqref{204}, we conclude that $f=\sum_{\vi\in[2^{\vn}-\vec1]}  f^{(\vi)}$, which further implies that
\begin{align*}
\| f\|_{\dot A^{\vec s,\tau}_{\vec p,\vec q,\pi}(V)}
&\lesssim \sum_{\vi\in[2^{\vn}-\vec1]}  \big\|f^{(\vi)}\big\|_{\dot A^{\vec s,\tau}_{\vec p,\vec q,\pi}(V)}\\
&\lesssim \sum_{\vi\in[2^{\vn}-\vec1]} \Big\|\Big\{\Big\langle f,\theta^{(\vi)}_Q \Big\rangle\Big\}_{Q\in\mathscr{D}(\mathbb R^{\vec n})}\Big\|_{\dot a^{\vec s,\tau}_{\vec p,\vec q,\pi}(V)}
=\| f\|_{\dot A^{\vec s,\tau}_{\vec p,\vec q,\pi}(V)_{\mathrm w}}
\end{align*}
by the quasi-triangle inequality in $\dot A^{\vec s,\tau}_{\vec p,\vec q,\pi}(V)$.
This finishes the proof of \eqref{2160} and hence Theorem \ref{wavelet 2}.
\end{proof}

\begin{remark}
In the multi-parameter context, Theorems \ref{wavelet} and \ref{wavelet 2} are completely new.
In the one-parameter case, for the usual Besov and Triebel--Lizorkin space with $\tau=0$, and more generally whenever the parameters are such that the definitions of the relevant space with open sets and dyadic cubes coincide (see Theorem \ref{compare A}), Theorems \ref{wavelet} and Theorem \ref{wavelet 2} coincide with \cite[Theorems 4.3 and 4.10]{BHYY:1c} with a slightly more restrictive range of $\vec{\mathscr N}$. In the remaining one-parameter cases, Theorems \ref{wavelet} and \ref{wavelet 2} deal with a slightly different scale of spaces from those in \cite{BHYY:1c}, as detailed in Theorem \ref{compare A}.
\end{remark}

\section{Relationship with other matrix-weighted function spaces}\label{other}

In this section, we show that the spaces $\Atau(V)$ contain
the matrix-weighted versions of various classical function spaces, including Lebesgue spaces (Subsection \ref{Lp}),
Sobolev spaces (Subsection \ref{sec:Sobolev}),
and BMO spaces (Subsection \ref{BMO}).
We also study the relationship between
$\Atau(V)$ and a variant obtained by replacing the supremum over open sets in the definition by one over dyadic rectangles only (Subsections \ref{BTL} and \ref{BTL another} for $k=1$
and $k\geq 2$, respectively).

\subsection{Matrix-weighted $L^p$ spaces}\label{Lp}

They are certainly the most important examples of matrix-weighted function spaces, going back to works of Nazarov, Treil, and Volberg \cite{NT96,TV97,Vol97} or even Wiener and Masani \cite{WM58} for $p=2$. Hence, it is relevant to know that these spaces are also covered as a special case of our theory. In the one-parameter case, the following result is due to \cite[Theorem 4.1]{FR21}. (Variant with a Triebel--Lizorkin-type norm, but involving Haar functions, already appear in \cite{NT96,Vol97} and more recently in \cite{Is21}.) We provide the multi-parameter extension.

\begin{theorem}\label{Lp=F}
For $p\in(1,\infty)$ and $V\in\A_p(\R^{\vn})$, we have
\begin{equation*}
\dot F^{0,0}_{p,2}(V)=L^p(V)
\end{equation*}
in the following precise sense:
\begin{enumerate}[\rm(i)]
\item\label{Lp subset F} Every $f\in L^p(V)$ defines in the usual way a distribution $f\in\Sc_0'(\R^{\vn};\C^m)$, and this satisfies
\begin{equation*}
\Norm{f}{\dot F^{0,0}_{p,2}(V)}\lesssim\Norm{f}{L^p(V)}.
\end{equation*}

\item\label{F subset Lp} For every $f\in\dot F^{0,0}_{p,2}(V)\subset\Sc_0'(\R^{\vn};\C^m)$, there is a unique function $F\in L^p(V)$ that agrees with $f$ as an element of $\Sc_0'(\R^{\vn};\C^m)$, and
\begin{equation*}
\Norm{F}{L^p(V)}\lesssim\Norm{f}{\dot F^{0,0}_{p,2}(V)}.
\end{equation*}
\end{enumerate}
Here, all the implicit positive constants are independent of $f$.
\end{theorem}

\begin{proof}[Proof of Theorem \ref{Lp=F}, part \eqref{Lp subset F}]
Let us first consider $k=1$. In this case, the result is due to \cite[Theorem 4.1]{FR21}, but we indicate a slight variant of the argument that readily extends to several parameters.  Fix a Littlewood--Paley function $\varphi$ on $\R^{n}$. For finite subsets $\Fs\subset\Z$ and signs $\eps_i=\pm 1$, consider the operator
\begin{equation}\label{LPop}
T: f\mapsto\sum_{i\in\Fs}\eps_i\varphi_i*f.
\end{equation}
It is well known that $T$ is a convolution-type Calder\'on--Zygmund operator, whose kernel $K$ satisfies the standard estimates
\begin{equation*}
\abs{K(x)}\lesssim\frac{1}{\abs{x}^n},\qquad\abs{\nabla K(x)}\lesssim\frac{1}{\abs{x}^{n+1}}
\end{equation*}
uniform with respect to the choice of $\Fs$ and $\eps_i$. Thus, from the well-known boundedness of such operators on $L^p(V)$ (see \cite[Theorem 5.1]{Gold03}), it follows that
\begin{equation*}
\BNorm{\sum_{i\in\Fs}\eps_i\varphi_i*f}{L^p(V)}\lesssim\Norm{f}{L^p(V)}
\end{equation*}
with the said uniformity. Thus, taking the expectation $\E$ over a random choice of $\eps_i=\pm 1$
(see, for instance, \cite[p.\,73]{HNVW1} for its definition), it follows that
\begin{equation*}
\Norm{f}{L^p(V)}^p\gtrsim\E\BNorm{\sum_{i\in\Fs}\eps_i\varphi_i*f}{L^p(V)}^p
=\int_{\R^n}\E\Babs{\sum_{i\in\Fs}\eps_i V(x)\varphi_i*f(x)}^p dx,
\end{equation*}
and by Khintchine's inequality \cite[Sections C.2--C.3]{Graf-clas},
\begin{equation*}
\E\Babs{\sum_{i\in\Fs}\eps_i V(x)\varphi_i*f(x)}^p
\sim \Big[ \sum_{i\in\Fs}\abs{V(x)\varphi_i*f(x)}^2 \Big]^{\frac{p}{2}},
\end{equation*}
where the positive equivalence constants depend only on $p$.
Thus,
\begin{equation*}
\Big\|\Big[\sum_{i\in\Fs}\abs{V(\varphi_i*f)}^2\Big]^{\frac{1}{2}}
\Big\|_{L^p(\R^n)}
\lesssim\Norm{f}{L^p(V)}
\end{equation*}
for any finite $\Fs\subset\Z$, and then by monotone convergence also for $\Fs=\Z$. The left-hand side is recognised as $\Norm{f}{\dot F^{0,0}_{p,2}(V)}$, and this completes the proof of \eqref{Lp subset F} for $k=1$.

In the case of arbitrarily many parameters, we consider the versions of the Littlewood--Paley operators \eqref{LPop} over each subspace $\R^{n_\nu}$. Let us denote the corresponding objects with superscript $(\nu)$ for each $\nu\in[k]$.

Recall from Lemma \ref{prod vs coord} that $V\in\A_p(\R^{\vn})$ implies that $x_\nu\mapsto V(x'\oplus x_\nu)\in\A_p(\R^{n_{\nu}})$ for each $\nu\in[k]$, uniformly in $x'\in\R^{\vn}\ominus\R^{n_{\nu}}$. Thus, by the one-parameter case,
\begin{align*}
\Norm{T^{(\nu)}f}{L^p(V)}^p
&=\int_{\R^{\vn}\ominus\R^{n_{\nu}}}\Norm{T^{(\nu)}f(x'\oplus\cdot)}{L^p(V(x'\oplus\ \cdot))}^p \, dx' \\
&\lesssim\int_{\R^{\vn}\ominus\R^{n_{\nu}}}\Norm{f(x'\oplus\cdot)}{L^p(V(x'\oplus\ \cdot))}^p \, dx'
=\Norm{f}{L^p(V)}^p.
\end{align*}
Using this repeatedly for every $\nu\in[k]$, we obtain
\begin{equation*}
\Bigg\|\sum_{\vi\in\Fs}\eps^{\otimes}_{\vi}\psi_{\vi}*f\Bigg\|_{L^p(V)}
=\Bigg\|\Bigg(\prod_{\nu=1}^k T^{(\nu)}\Bigg)f\Bigg\|_{L^p(V)}
\lesssim\Norm{f}{L^p(V)},
\end{equation*}
where $\Fs=\Fs^{(1)}\times\cdots\times\Fs^{(k)}$ and $\psi_{\vi}(x):=\prod_{\nu=1}^k\psi_{\vi_\nu}^{(\nu)}(x_{\nu})$ is the usual product Littlewood--Paley function, and $\eps^{\otimes}_{\vi}:=\prod_{\nu=1}^k\eps^{(\nu)}_{i_\nu}$. Averaging over independent $\eps^{(\nu)}_{i_\nu}=\pm 1$ as for $k=1$, we arrive at
\begin{equation*}
\begin{split}
\Norm{f}{L^p(V)}
&\gtrsim\int_{\R^{\vn}}\E\Bigg|\sum_{\vi\in\Fs}\eps^{\otimes}_{\vi}V(x)\psi_{\vi}*f(x)\Bigg|^p dx \\
&\sim\int_{\R^{\vn}}\Bigg[\sum_{\vi\in\Fs}\abs{V(x)\psi_{\vi}*f(x)}^2\Bigg]^{\frac{p}{2}} dx,
\end{split}
\end{equation*}
where the last step used Khintchine's inequality for multiple random sums  \cite[Section C.5]{Graf-clas}. Applying monotone convergence to $\Fs\uparrow\Z^k$, we complete the proof of \eqref{Lp subset F} in the general case.
\end{proof}

As a preparation for the proof of part \eqref{F subset Lp} of Theorem \ref{Lp=F}, we provide the following density result. As it might have other applications, we state it in a form that is more general than our immediate needs.

\begin{proposition}\label{density}
Let $\vs\in\R^k$, $p\in(0,\infty)$, $\vp=p\cdot\vone$, $\vq\in(0,\infty)^k$, $\pi\in S_{[2k]}$ be an admissible permutation of $(\vp,\vq)$, and $V\in\A_p(\R^{\vn})$.
Let $(\varphi,\psi)$ be a Littlewood--Paley pair on $\R^{\vn}$.
For each finite subset $\Fs \subset\Z^k$ and $f\in\Sc_0'(\R^{\vn};\C^m)$, let
\begin{equation*}
P_{\Fs} f:=\sum_{\vj\in\Fs}\widetilde\varphi_{\vj}*\psi_{\vj}*f.
\end{equation*}
If $\Fs_N\subset\Z^k$ is an increasing sequence of finite subsets with $\bigcup_{N=1}^\infty\Fs_N=\Z^k$
and if $f\in \dot A_{\vec p,\vec q,\pi}^{\vec s,0}(V)$, then $\Norm{f-P_{\Fs_N} f}{\dot A_{\vec p,\vec q,\pi}^{\vec s,0}(V)}\to 0$ as $N\to\infty$.
\end{proposition}

This generalises \cite[Lemma 4.5]{FR21}, even in the one-parameter case, where their result is formulated for $p\in(1,\infty)$ and $\Fs_N=[-N,N]\cap\Z$. For such $\Fs_N$, there are just two boundary terms, and the $\dot A_{\vec p,\vec q,\pi}^{\vec s,0}(V)$ quasi-norms of single-term sequences reduce to $L^p(V)$ norms, for which $p\in(1,\infty)$ gives access to a version of Young's convolution inequality \cite[Lemma 4.4]{FR21}.
In our setting, we need to replace these points by other tools, although other aspects of our proof are similar to those of \cite[Lemma 4.5]{FR21}.

\begin{proof}[Proof of Proposition \ref{density}]
The sum defining $P_{\Fs} f$ is finite and hence clearly well defined. In place of $\Norm{\ }{\Atau(V)}$, we will consider the equivalent (by Corollary \ref{3 norms cor}) quasi-norm $\Norm{\ }{\Atau([V]_p)}$ defined with the reducing operators.

To estimate $\Norm{f-P_{\Fs} f}{\Atau([V]_p)}$ we need to compute $\varphi_{\vi}*(f-P_{\Fs} f)$. Let
\begin{equation*}
\Fs^- :=\{\vj\in\Z^k:\
\vj+\vu\in\Fs\text{ for all }\vu\in\{-1,0,1\}^k\}.
\end{equation*}
Then $\Fs^-\subset\Fs$ for each $\Fs\subset\Z^k$, and the new sequence $\Fs_N^-$ is also increasing with $\bigcup_{N=1}^\infty\Fs_N^-=\Z^k$.

On the Fourier transform side, we see that
\begin{equation*}
[\varphi_{\vi}*(f-P_{\Fs} f)]^{\wedge}
=\widehat{\varphi}_{\vi} \Big( 1-\sum_{\vj\in\Fs}
\overline{\widehat{\varphi}_{\vj}}\widehat{\psi}_{\vj} \Big)\widehat{f},
\end{equation*}
where the supports of $\widehat{\varphi}_{\vi}$ and $\overline{\widehat{\varphi}_{\vj}}\widehat{\psi}_{\vj}$ meet only if $\vj-\vi\in\{-1,0,1\}^k$.

If $\vi\in\Fs^-$, then the sum above is identically $1$ on the support of $\widehat{\varphi}_{\vi}$. It follows that
\begin{equation*}
\varphi_{\vi}*(f-P_{\Fs} f)
=\begin{cases} 0, & \text{if }\vi\in\Fs^-, \\ \displaystyle
\varphi_{\vi}*f-\sum_{\genfrac{}{}{0pt}{}{\vu\in\{-1,0,1\}^k}{\vi+\vu\in\Fs}}\chi_{\vi+\vu}*\varphi_{\vi}*f, & \text{else},
\end{cases}
\end{equation*}
where $\chi:=\widetilde\varphi*\psi$.
Hence
\begin{equation}\label{2termsRHS}
\begin{split}
\Norm{f-P_{\Fs}f}{\dot A_{\vec p,\vec q,\pi}^{\vec s,0}([V]_p)}
&=\Norm{\{\Red{V}{\vi}{p}(\varphi_{\vi}*(f-P_{\Fs} f))\}_{\vi\in\Z^k}}{\dot \alpha_{\vec p,\vec q,\pi}^{\vec s,0}} \\
&\lesssim\Norm{\{\Red{V}{\vi}{p}(\varphi_{\vi}*f)\}_{\vi\in\Z^k\setminus\Fs^-}}{\dot \alpha_{\vec p,\vec q,\pi}^{\vec s,0}} \\
&\quad+\sum_{\vu\in\{-1,0,1\}^k}\BNorm{\Big\{\Red{V}{\vi}{p}(\chi_{\vi+\vu}*\varphi_{\vi}*f)\Big\}_{\genfrac{}{}{0pt}{}{\vi\in\Z^k\setminus\Fs^-}{\vi+\vu\in\Fs}}}{\dot \alpha_{\vec p,\vec q,\pi}^{\vec s,0}}.
\end{split}
\end{equation}

Let us consider the first term on the right-hand side of \eqref{2termsRHS},
\begin{equation*}
\Norm{\{\Red{V}{\vi}{p}(\varphi_{\vi}*f)\}_{\vi\in\Z^k\setminus\Fs^-}}{\dot \alpha_{\vec p,\vec q,\pi}^{\vec s,0}}
=\Norm{\{2^{\vi\cdot\vs}\Red{V}{\vi}{p}(\varphi_{\vi}*f)\}_{\vi\in\Z^k\setminus\Fs^-}}{(L^{\vp}\ell^{\vq})_\pi}.
\end{equation*}
We write the iterated $(L^{\vp}\ell^{\vq})_\pi$ as
\begin{equation*}
(L^{\vp}\ell^{\vq})_\pi=
(L^{\vp_{\operatorname{out}}}\ell^{\vq_{\operatorname{out}}})_{\pi_{\operatorname{out}}}\ell^{q_\nu}
(L^{\vp_{\operatorname{in}}}),
\end{equation*}
identifying the innermost instance of an $\ell^{q_\nu}$ space in this iteration. Note that either $\vp_{\operatorname{out}}$ or $\vp_{\operatorname{in}}$ may be just an empty string, in which case it is simply omitted above.
Let us write $\vi=\vi_{\operatorname{out}}\oplus i_{\nu}$ with
$\vi_{\operatorname{out}}\in\Z^{k-1}$
and $x=x_{\operatorname{out}}\oplus x_{\operatorname{in}}$, where $x_{\operatorname{in}}$ is the variable (if any) of the $L^{\vp_{\operatorname{in}}}$
space inside the last $\ell^{q_{\nu}}$, and $x_{\operatorname{out}}$ consist of the remaining variables. Then the $(L^{\vp}\ell^{\vq})_\pi$ quasi-norm above is the
$(L^{\vp_{\operatorname{out}}}\ell^{\vq_{\operatorname{out}}})_{\pi_{\operatorname{out}}}$ quasi-norm of
\begin{equation}\label{tail}
\Bigg(\sum_{i_\nu\in\Z: i_\nu\oplus\vi_{\operatorname{out}}\in\Z^k\setminus\Fs^-}
\Norm{2^{\vi\cdot\vs}\Red{V}{\vi}{p}(\varphi_{\vi}*f)(x_{\operatorname{out}},\cdot)}
{L^{\vp_{\operatorname{in}}}}^{q_\nu}\Bigg)^{\frac{1}{q_\nu}}.
\end{equation}
This is the tail of the infinite convergent series (at least for almost every $x_{\operatorname{out}}$) with summation simply over all $i_\nu\in\Z$. Such a tail converges to $0$
as $\Fs^-=\Fs_N^-\uparrow\Z^k$. By the finiteness of the quasi-norm of the full series, and the assumption that all exponents are finite, it then follows from dominated convergence that also the $(L^{\vp_{\operatorname{out}}}\ell^{\vq_{\operatorname{out}}})_{\pi_{\operatorname{out}}}$ quasi-norm of \eqref{tail} converges to zero as $N\to\infty$, i.e.,
\begin{equation*}
\Norm{\{\Red{V}{\vi}{p}(\varphi_{\vi}*f)\}_{\vi\in\Z^k\setminus\Fs_N^-}}{\dot \alpha_{\vec p,\vec q,\pi}^{\vec s,0}}
\to 0\quad\text{as}\quad N\to\infty.
\end{equation*}
This handles the first term on the right-hand side of \eqref{2termsRHS}, and it remains to control the second one. (A careful reader will have observed that the proof so far made no use of the particular properties of the weight $V$, and we could have equally well used just $V$ in place of $\Red{V}{\vi}{p}$; the reducing operators will play a role in the remainder of the proof.)

First we observe that each $\ell^{q_\nu}$ quasi-norm, and hence also the $\dot \alpha_{\vec p,\vec q,\pi}^{\vec s,0}$ quasi-norm, is monotone with respect to truncations of the sequence inside the quasi-norm, and hence we may simply estimate up by dropping the condition $\vi+\vu\in\Fs$. Hence, for each of the finitely many $\vu\in\{-1,0,1\}^k$, we need to estimate
\begin{equation*}
\Norm{\{
\Red{V}{\vi}{p}( \chi_{\vi+\vu} * \varphi_{\vi} * f) \}_{\vi\in\Z^k\setminus\Fs^-}  }{\dot \alpha_{\vec p,\vec q,\pi}^{\vec s,0}}.
\end{equation*}
Note further that $\chi_{\vi+\vu}=(\chi_{\vu})_{\vi}$ are usual dilations of the finitely many test functions $\chi_{\vu}\in\Sc(\R^{\vn})$ with $\vu\in\{-1,0,1\}^k$. By Lemma \ref{cald5}, applied to $V_Rf:=[V]_{\aveL^p(R)}f$ in place of $f$, and $\chi_{\vu}$ in place of $\chi$, we obtain for $x\in R\in\D_{\vi}(\R^{\vn})$
\begin{equation*}
\begin{split}
\abs{V_R(\chi_{\vi+\vu}*\varphi_{\vi} * f)(x)}^r
&\lesssim(\phi_{\vi}*\abs{V_R(\varphi_{\vi} * f)}^r)(x) \\
&\leq\sum_{P\in\D_{\vi}(\R^{\vn})} \int_P\phi_{\vi}(x-y)\abs{V_R V_P^{-1}}^r\abs{V_P(\varphi_{\vi}*f)(y)}^r dy \\
&\lesssim\sum_{P\in\D_{\vi}(\R^{\vn})} \int_P\phi_{\vi}'(x-y)\abs{V_P(\varphi_{\vi}*f)(y)}^r dy \\
&=(\phi_{\vi}*\abs{\Red{V}{\vi}{p}(\varphi_{\vi}*f)}^r)(x),
\end{split}
\end{equation*}
where the last estimate was an application of Lemma \ref{sharp coef 2}, resulting in a new function $\phi_{\vi}'$ in place of $\phi_{\vi}$, still with as good decay as we want.

Taking quasi-norms of both sides it follows that
\begin{equation*}
\begin{split}
\Norm{\{     \Red{V}{\vi}{p}( \chi_{\vi+\vu} * \varphi_{\vi} * f) \}_{\vi\in\Z^k\setminus\Fs^-}  }{\dot \alpha_{\vec p,\vec q,\pi}^{\vec s,0}}
&\lesssim   \Norm{\{(\phi_{\vi}*\abs{\Red{V}{\vi}{p}(\varphi_{\vi} * f)}^r)^{\frac1r}\}_{\vi\in\Z^k\setminus\Fs^-} }{\dot \alpha_{\vec p,\vec q,\pi}^{\vec s,0}} \\
&\lesssim   \Norm{\{\Red{V}{\vi}{p}(\varphi_{\vi} * f)\}_{\vi\in\Z^k\setminus\Fs^-} }{\dot \alpha_{\vec p,\vec q,\pi}^{\vec s,0}},
\end{split}
\end{equation*}
where the last estimate was an application of Lemma \ref{convo bd A bd}, valid provided that we choose $r\in(0,\min\{\vp,\vq\})$, as we can.

Thus, we estimated the second term on the right-hand side of \eqref{2termsRHS} by the first term, and we already checked that the first term tends to $0$ for $\Fs=\Fs_N$ as $N\to\infty$.
This finishes the proof of Proposition \ref{density}.
\end{proof}

\begin{proof}[Proof of Theorem \ref{Lp=F}, part \eqref{F subset Lp}]
Let $f\in\dot F^{0,0}_{p,2}(V) \subset \Sc_0'(\R^{\vn};\C^m)$, and let $P_{\Fs} f$ be as in Proposition \ref{density}.
By the well-known properties of convolutions of distributions and test functions, each $\psi_{\vj}*f$ is a smooth function of polynomially bounded growth. Hence, for any compactly supported $g\in L^\infty_{\mathrm c}(\R^{\vn})$, the following expressions are well defined, and the identity can be justified by Fubini's theorem:
\begin{equation*}
\pair{P_{\Fs} f}{g}=\sum_{\vj\in\Fs}\pair{\widetilde\varphi_{\vj}*\psi_{\vj}*f}{g}
=\sum_{\vj\in\Fs}\pair{\psi_{\vj}*f}{\varphi_{\vj}*g}.
\end{equation*}
On the right, each object on both sides of the pairing is a proper function, and the pairing agrees with the usual integral pairing. Hence we may insert the weight $V$ into this pairing by
\begin{equation*}
\pair{\psi_{\vj}*f}{\varphi_{\vj}*g}
=\pair{V(\psi_{\vj}*f) }{ V^{-1}(\varphi_{\vj}*g)},
\end{equation*}
and then
\begin{equation*}
\abs{ \pair{P_{\Fs} f}{g} }
\leq \Bigg\| \Bigg[ \sum_{\vj\in\Fs} \abs{V(\psi_{\vj}*f)}^2 \Bigg]^{\frac{1}{2}} \Bigg\|_{L^p}
\Bigg\|\Bigg[ \sum_{\vj\in\Fs}\abs{V^{-1}(\varphi_{\vj}*g)}^2 \Bigg]^{\frac12}\Bigg\|_{L^{p'}},
\end{equation*}
where
\begin{equation*}
\Bigg\|\Bigg[ \sum_{\vj\in\Fs}\abs{V^{-1}(\varphi_{\vj}*g)}^2 \Bigg]^{\frac12}\Bigg\|_{L^{p'}}
\leq \Norm{g}{ \dot F^{0,0}_{p',2} ( V^{-1} )  }
\lesssim \Norm{g}{ L^{p'}( V^{-1})}
\end{equation*}
by part \eqref{Lp subset F} of the theorem, which we already proved, with $p'$ and $V^{-1}\in\A_{p'}(\R^{\vn})$ in place of $p$ and $V\in\A_p(\R^{\vn})$, and $g\in L^\infty_{\mathrm c}(\R^{\vn})\subset L^{p'}(V^{-1})$ in place of~$f$.

Since $P_{\Fs} f$ is a proper function, and $L^\infty_{\mathrm c}(\R^{\vn})$ is dense in $L^{p'}(V^{-1})=(L^p(V))^*$, the previous estimates imply that
\begin{equation}\label{PFf<}
\Norm{P_{\Fs} f}{L^p(V)}
\lesssim \Bigg\| \Bigg[ \sum_{\vj\in\Fs} \abs{V(\psi_{\vj}*f)}^2 \Bigg]^{\frac{1}{2}} \Bigg\|_{L^p}.
\end{equation}

Let us then consider an increasing sequence of finite $\Fs_N\subset\Z^k$ with $\bigcup_{N=1}^\infty\Fs_N=\Z^k$, as in Proposition \ref{density}. Applying \eqref{PFf<} to $\Fs=\Fs_N\setminus\Fs_M$ for $M<N$, and noting that $P_{\Fs}=P_{\Fs_N}-P_{\Fs_M}$ in this case, we obtain
\begin{align*}
\Norm{P_{\Fs_N} f-P_{\Fs_M} f}{L^p(V)}
&\lesssim \Bigg\| \Bigg[ \sum_{\vj\in\Fs_N\setminus\Fs_M} \abs{V(\psi_{\vj}*f)}^2 \Bigg]^{\frac{1}{2}} \Bigg\|_{L^p} \\
&\leq \Bigg\| \Bigg[ \sum_{\vj\in\Z^k\setminus\Fs_M} \abs{V(\psi_{\vj}*f)}^2 \Bigg]^{\frac{1}{2}} \Bigg\|_{L^p}
\to 0
\end{align*}
as $M\to\infty$ by dominated convergence, since
\begin{equation*}
\Bigg\| \Bigg[ \sum_{\vj\in\Z^k} \abs{V(\psi_{\vj}*f)}^2 \Bigg]^{\frac{1}{2}} \Bigg\|_{L^p}
=\Norm{f}{\dot F^{0,0}_{p,2}(V)}<\infty.
\end{equation*}
Hence $P_{\Fs_N} f$ is a Cauchy sequence in $L^p(V)$, and thus converges to some $F\in L^p(V)$. On the other hand, Proposition \ref{density} shows that $P_{\Fs_N} f\to f$ in $\dot F^{0,0}_{p,2}(V)$, and thus in $\Sc_0'(\R^{\vn};\C^m)$ by Corollary \ref{convergence}. For $h\in L^p(V)$ and $\phi\in\Sc_0(\R^{\vn})$,
\begin{equation*}
\abs{\pair{h}{\phi}}=\abs{\pair{Vh}{V^{-1}h}}\leq\Norm{h}{L^p(V)}\Norm{\phi}{L^{p'}(V^{-1})},
\end{equation*}
and the finiteness of the second term is easy to check; hence convergence in $L^p(V)$ also implies convergence in $\Sc_0'(\R^{\vn};\C^m)$.
Now $F=f$ in the sense of $\Sc_0'(\R^{\vn};\C^m)$ follows from the uniqueness of the limit in this space.
\end{proof}

\subsection{Matrix-weighted Sobolev spaces}\label{sec:Sobolev}

Another important class of matrix-weighted extensions of classical function spaces, these have been previously identified with Triebel--Lizorkin spaces in the one-parameter setting by Frazier and Roudenko \cite[Proposition 1.4]{FR21}. As usual, we provide a multi-parameter extension. While \cite[Proposition 1.4]{FR21} deals with non-homogeneous Sobolev spaces (involving $L^p$ norms of both a function and its derivatives), we prove a result for homogeneous Sobolev spaces, which may be a new statement (although an easy variant of \cite[Proposition 1.4]{FR21}) also in the one-parameter case.

\begin{definition}
Let $\vec l\in\mathbb N^k$, $p\in(0,\infty)$,
and $V\in L_{\loc}^p(\R^{\vec n};\C^{m\times m})$.
The \emph{homogeneous matrix-weighted Sobolev space} $\dot W^{\vec l}_p(V)$
consists of all $f\in\Sc_0'(\R^{\vn};\C^m)$ such that
the distributional derivatives $\partial^{\alpha} f$
are functions belonging to $L^p(V)$
for all $\alpha\in \mathbb N^{\vec n}$ with $|\alpha|_{\rm{vec}}= \vec l$.
We equip this set with the quasi-norm
\begin{equation}\label{def Sobo}
\|f\|_{\dot W^{\vec l}_p(V)}
:=\sum_{|\alpha|_{\rm{vec}}=\vec l} \|\partial^{\alpha} f\|_{L^p(V)}.
\end{equation}
\end{definition}

\begin{theorem}\label{thm:Sobolev}
For $\vec l\in\N^k$, $p\in(1,\infty)$, and $V\in \A_p(\mathbb R^{\vec n})$,
we have $\dot F^{\vec l,0}_{p,2}(V)=\dot W^{\vec l}_p(V)$.
Moreover, for every $f\in\Sc_0'(\R^{\vn};\C^m)$,
\begin{equation*}
\Norm{f}{\dot W^{\vec l}_p(V)}\sim\Norm{f}{\dot F^{\vec l,0}_{p,2}(V)},
\end{equation*}
where the positive equivalence constants are independent of $f$.
\end{theorem}

To prove Theorem \ref{thm:Sobolev}, we need the following equivalent characterisation of $\Atau(V)$ in terms of derivatives.

\begin{proposition} \label{partial derivative}
Let $\vp\in(0,\infty)^k$, $\vq\in(0,\infty]^k$, $\vs\in\R^k$,
$\tau\in[0,\infty)$, $\pi\in S_{[2k]}$ be admissible for $(\vp,\vq)$, and $V$ be a matrix weight such that $(\vp,\vq,\pi,V)$ belongs to one of the main cases (Definition \ref{main cases}).
Assume that $\vec l\in\mathbb N^k$ and $f\in \Sc_{0}'(\R^{\vn};\C^m)$.
Then $f\in\Atau(V)$ if and only if
$\partial^{\alpha} f\in \dot A_{\vec p,\vec q,\pi}^{\vec s-\vec l,\tau}(V)$
for all $\alpha\in\N^{\vn}$ with $|\alpha|_{\rm{vec}}=\vec l$.
Moreover, for every $f\in\Sc_0'(\R^{\vn};\C^m)$,
$$
\|f\|_{\Atau(V)}
\sim \sum_{|\alpha|_{\rm{vec}}=\vec l} \|\partial^{\alpha} f\|_{\dot A_{\vec p,\vec q,\pi}^{\vec s-\vec l,\tau}(V)},
$$
where the positive equivalence constants are independent of $f$.
\end{proposition}

\begin{proof}
We first prove that, for every $\alpha\in \mathbb N^{\vec n}$ with $|\alpha|_{\rm{vec}}=\vec l$,
\begin{equation}\label{B<A}
\|\partial^{\alpha} f\|_{\dot A_{\vec p,\vec q,\pi}^{\vec s-\vec l,\tau}(V)}
\lesssim \|f\|_{\Atau(V)}.
\end{equation}
Let $(\varphi,\psi)$ be a Littlewood--Paley pair on $\R^{\vn}$.
By Corollary \ref{phiV}, we conclude that
\begin{align*}
\|\partial^{\alpha} f\|_{\dot A_{\vec p,\vec q,\pi}^{\vec s-\vec l,\tau}(V)}
\lesssim\|S_\varphi\partial^{\alpha} f\|_{\dot a_{\vec p,\vec q,\pi}^{\vec s-\vec l,\tau}(V)}.
\end{align*}
Note that $\langle\partial^{\alpha} f,\varphi_Q\rangle
=\langle f,(-\partial)^{\alpha} \varphi_Q\rangle.$
Using Lemma \ref{cald2} and the integration by parts, we obtain
\begin{align*}
\partial^{\alpha} \varphi_Q
&=\sum_{R\in\mathscr D(\mathbb R^{\vec n})}
\langle\partial^{\alpha} \varphi_Q,\psi_R\rangle \psi_R
=(-1)^{|\vec l|}\sum_{R\in\mathscr D(\mathbb R^{\vec n})}
\langle\varphi_Q,\partial^{\alpha} \psi_R\rangle \psi_R\\
&=(-1)^k \prod_{i=1}^k \ell(Q_i)^{-|\alpha_i|}
\sum_{R\in\mathscr D(\mathbb R^{\vec n})}
b_{Q,R} \psi_R
\end{align*}
with convergence in $\mathscr S_0(\mathbb R^{\vec n})$, where
$$
b_{Q,R}:= \prod_{i=1}^k \bigg[\frac{\ell(Q_i)}{\ell(R_i)}\bigg]^{|\alpha_i|}
\langle\varphi_Q,(\partial^{\alpha} \psi)_R\rangle.
$$
Therefore,
$$
\|S_\varphi\partial^{\alpha} f\|_{\dot a_{\vec p,\vec q,\pi}^{\vec s-\vec l,\tau}(V)}
\leq \|B (S_\psi f)\|_{\dot a_{\vec p,\vec q,\pi}^{\vec s,\tau}(V)}.
$$
By the fact that $\varphi$ is a Littlewood--Paley function on $\R^{\vn}$
and Parseval's identity,
we conclude that $\langle\varphi_Q,(\partial^{\alpha} \psi)_R\rangle\neq 0$
if and only if $\frac{\ell(Q_i)}{\ell(R_i)}\in\{\frac12,1,2\}$ for all $i\in[k]$.
This, together with the rapid decay of Schwartz pairings (Lemma \ref{59x}),
further implies that $B:=\{b_{Q,R}\}$ is $(\vec D, \vec E,\vec F)$-almost diagonal
for all $\vec D, \vec E,\vec F\in \R^k$,
and hence
$$
\|B (S_\psi f)\|_{\dot a_{\vec p,\vec q,\pi}^{\vec s,\tau}(V)}
\lesssim \|S_\psi f\|_{\dot a_{\vec p,\vec q,\pi}^{\vec s,\tau}(V)}
\lesssim \|f\|_{\dot A_{\vec p,\vec q,\pi}^{\vec s,\tau}(V)}.
$$
This finishes the proof of \eqref{B<A}.

Next, we show that, for every $\alpha\in \mathbb N^{\vec n}$ with $|\alpha|_{\rm{vec}}=\vec l$,
\begin{equation}\label{A<B}
\|f\|_{\Atau(V)}
\lesssim \|\partial^{\alpha} f\|_{\dot A_{\vec p,\vec q,\pi}^{\vec s-\vec l,\tau}(V)}.
\end{equation}
Let $\Delta^{\vec l}:=\prod_{i\in[k]}\Delta_i^{l_i}$,
where $\Delta_i$ is the Laplace operator on $\R^{n_i}$.
Note that
\begin{equation*}
\Delta_i^{l_i}=\Bigg(\sum_{j=1}^{n_i}\partial_{i,j}^2\Bigg)^{l_i}
=\sum_{\abs{\beta_i}=l_i}
\binom{l_i}{\beta_i}\prod_{j=1}^{n_i}\partial_{i,j}^{2\beta_{i,j}}
\end{equation*}
and hence
\begin{equation*}
\Delta^{\vec l}
=\sum_{\abs{\beta}_{\rm{vec}}=\vl}
\binom{\vec l}{\vec\beta}\partial^{2\beta}
= C (-1)^{|\vec l|} \dot I_{-2\vec l}
\end{equation*}
in $\mathcal S_0(\R^{\vn})$,
where $\binom{\vec l}{\vec\beta}:=\prod_{i=1}^k \binom{l_i}{\beta_i}$,
$C$ is a positive constant,
and the lifting operator $\dot I_{-2|\alpha|_{\rm{vec}}}$ is defined as in \eqref{lift}. Then
\begin{equation*}
\begin{split}
\|f\|_{\Atau(V)}
&\sim \|\dot I_{-2\vec l} f\|_{\dot A_{\vec p,\vec q,\pi}^{\vec s-2\vec l,\tau}(V)}
\quad\text{by Corollary \ref{257V}} \\
&\sim \Big\|\Delta^{\vec l} f\Big\|_{\dot A_{\vec p,\vec q,\pi}^{\vec s-2\vec l,\tau}(V)}
\lesssim\sum_{\abs{\alpha}_{\rm{vec}}=\vl}
\|\partial^{\alpha}(\partial^{\alpha} f)\|_{\dot A_{\vec p,\vec q,\pi}^{(\vec s-\vec l)-\vec l,\tau}(V)} \\
&\lesssim \sum_{|\alpha|_{\rm{vec}}=\vec l}
\|\partial^{\alpha} f\|_{\dot A_{\vec p,\vec q,\pi}^{\vec s-\vec l,\tau}(V)} \quad\text{by \eqref{B<A}},
\end{split}
\end{equation*}
where we used \eqref{B<A} with $\partial^{\alpha}f$ in place of $f$ and $\vec s-\vec l$ in place of $\vec s$.
This finishes the proof of \eqref{A<B}
and hence Proposition \ref{partial derivative}.
\end{proof}

\begin{remark}
Proposition \ref{partial derivative} extends some previous results as follows.
For one-parameter Triebel--Lizorkin space, it coincides with
\cite[Proposition 7.3]{FR21} in the matrix-weighted case,
and with \cite[(5.6)]{FJW91} in the unweighted case.
\end{remark}

Now, we prove Theorem \ref{thm:Sobolev}.

\begin{proof}[Proof of Theorem \ref{thm:Sobolev}]
Combining two previous results (as quoted below),
we find that, for every $f\in\Sc_0'(\R^{\vn};\C^m)$,
\begin{equation*}
\begin{split}
\|f\|_{\dot F^{\vec l,0}_{p,2}(V)}
&\sim \sum_{|\alpha|_{\rm{vec}}=\vec l}\|\partial^{\alpha} f\|_{\dot F^{0,0}_{p,2}(V)}
\quad\text{by Proposition \ref{partial derivative}} \\
&\sim \sum_{|\alpha|_{\rm{vec}}=\vec l}\|\partial^{\alpha} f\|_{L^p(V)}
\quad\text{by Theorem \ref{Lp=F}} \\
&= \|f\|_{\dot W^{\vec l}_p(V)}\quad\text{by definition \eqref{def Sobo}}.
\end{split}
\end{equation*}
This finishes the proof of Theorem \ref{thm:Sobolev}.
\end{proof}

\subsection{Multi-parameter BMO spaces}\label{BMO}

Introduced by John and Nirenberg \cite{JN61}, the space BMO of functions of bounded mean oscillation has a rich history. 
For the present, we are mainly concerned about multi-parameter version of this space introduced by Chang and R.~Fefferman \cite{CF80}. The following is a slight variant of one of the original characterisations of this space, which has been adopted as the definition in some subsequent works (e.g.\ \cite[(3.5)]{LPPW}; we note that this part is independent of the error discovered in this paper \cite{LPPWe}).

\begin{definition}[cf. \cite{LPPW}, (3.5)]
Let $f\in L^2_{\mathrm{loc}}(\mathbb R^{\vn})$
and $\{\theta^{(\vi)}\}_{\vi\in[2^{\vec n}-\vec 1]}$ be the Meyer wavelets on $\mathbb R^{\vec n}$.
Then $f\in\mathrm{BMO}(\R^{\vn})$ if and only if
\begin{equation}\label{wave bmo}
  \sum_{\vi\in[2^{\vn}-\vone]}\sup_{\Omega\in\Open(\R^{\vn})} \Bigg[ \frac{1}{|\Omega|}
  \sum_{\vj\in\mathbb Z^k} \sum_{R\in \Omega_{\vj}}
|(f,\theta_R^{(\vi)})|^2 \Bigg]^{\frac12} <\infty.
\end{equation}
\end{definition}

Based on this definition, we introduce the matrix-weighted version as follows:

\begin{definition}\label{def:BMO(V)}
Let $V\in L^2_{\mathrm{loc}}(\mathbb R^{\vec n};\mathbb C^{m\times m})$
and $\{\theta^{(\vi)}\}_{\vi\in[2^{\vec n}-\vec 1]}$ be the Meyer wavelets on $\mathbb R^{\vec n}$.
The \emph{matrix-weighted bounded mean oscillation space $\mathrm{BMO}(V)$}
consists of all vector-valued functions $b\in L^1_{\mathrm{loc}}(\mathbb R^{\vec n};\mathbb C^m)$ such that
$$
\|b\|_{\mathrm{BMO}(V)}
:= \sum_{\vi\in[2^{\vec n}-\vec 1]}
\sup_{\Omega\in\Open(\R^{\vn})} \Bigg[ \frac{1}{|\Omega|}
\sum_{\vj\in\mathbb Z^k} \sum_{R\in \Omega_{\vj}}
\fint_R |V(x)(b,\theta^{(\vi)}_R)|^2 \,dx \Bigg]^{\frac12} <\infty.
$$
\end{definition}

\begin{remark}\label{rem:BMO(V)}
From the definition of reducing operators, we immediately deduce the alternative formula
\begin{equation}\label{eq:BMO(V)alt}
  \|b\|_{\mathrm{BMO}(V)}
\sim \sum_{\vi\in[2^{\vec n}-\vec 1]}
\sup_{\Omega\in\Open(\R^{\vn})} \Bigg[ \frac{1}{|\Omega|}
\sum_{\vj\in\mathbb Z^k} \sum_{R\in \Omega_{\vj}}
 |[V]_{\aveL^2(R)}(b,\theta^{(\vi)}_R)|^2 \Bigg]^{\frac12},
\end{equation}
where the $\aveL^2(R)$-reducing operator also has the explicit formula $[V]_{\aveL^2(R)}=\ave{V^2}_R^{\frac12}$

A version of the matrix-weighted multi-parameter BMO was also recently studied by Kaka\-roumpas and Soler i Gibert \cite{KS2}. In their notation, the right-hand side of \eqref{eq:BMO(V)alt} is the norm of the sequence $\{(b,\theta^{(\vi)})\}_{R\in\mathscr D,\vi\in[2^{\vn}-\vone]}$ in a space $\operatorname{BMO}_{\operatorname{prod},\mathscr D}(U,V,2)$, where $U\equiv I$ is the identity-valued constant weight. Indeed, \cite{KS2} define their multi-parameter BMO as a space of discrete sequences instead of functions, while they also discuss a more general two-weight situation. Thus, our setting has some points of contact with \cite{KS2}, but largely goes to a different direction.
\end{remark}

The following theorem shows that $\BMO(V)$ is a special case of the space $\Atau(V)$.

\begin{theorem} \label{coin BMO}
For $V\in\A_2(\R^{\vn})$, we have
\begin{equation*}
\dot F^{0,\frac12}_{2,2}(V)=\mathrm{BMO}(V)
\end{equation*}
in the following sense:
\begin{enumerate}[\rm(i)]
\item Every $b\in \mathrm{BMO}(V)$ defines a distribution $f\in\Sc_0'(\R^{\vn};\C^m)$, and this satisfies
\begin{equation*}
\Norm{f}{\dot F^{0,\frac12}_{2,2}(V)}\sim\Norm{b}{\mathrm{BMO}(V)}.
\end{equation*}

\item For every $f\in\dot F^{0,\frac12}_{2,2}(V)\subset\Sc_0'(\R^{\vn};\C^m)$, there is a function $b\in \mathrm{BMO}(V)$ that agrees with $f$ as an element of $\Sc_0'(\R^{\vn};\C^m)$, and
\begin{equation*}
\Norm{f}{\dot F^{0,\frac12}_{2,2}(V)}\sim\Norm{b}{\mathrm{BMO}(V)}.
\end{equation*}
\end{enumerate}
Here, all the positive equivalence constants are independent of $b$ and $f$.
\end{theorem}

\begin{proof}
Comparing Definitions \ref{d3.1} and \ref{def:BMO(V)}, and Remark \ref{rem:BMO(V)}, it is immediate that
\begin{equation}\label{eq:BMO(V)=f}
  \Norm{b}{\operatorname{BMO}(V)}
  =\sum_{\vi\in[2^{\vn}-\vone]}\BNorm{\{(b,\theta_R^{(\vi)})\}_{R\in\mathscr D}}{\dot f^{0,\frac12}_{2,2}(V)}
  \sim\sum_{\vi\in[2^{\vn}-\vone]}\BNorm{\{(b,\theta_R^{(\vi)})\}_{R\in\mathscr D}}{\dot f^{0,\frac12}_{2,2}([V]_2)}.
\end{equation}

By Theorem \ref{wavelet}, we find that,
for every $b\in \mathrm{BMO}(V)$,
$$
f:=\sum_{\vi\in[2^{\vec n}-\vec 1]}
\sum_{R\in \mathscr D(\mathbb R^{\vec n})}
\Big\langle b,\theta^{(\vi)}_R\Big\rangle \theta^{(\vi)}_R
\in \dot F^{0,\frac12}_{2,2}(V)
$$
and $\|f\|_{\dot F^{0,\frac12}_{2,2}(V)}\sim \|b\|_{\mathrm{BMO}(V)}$.
Conversely, for every $f\in \dot F^{0,\frac12}_{2,2}(V)$ and $x\in\mathbb R^{\vec n}$, let
\begin{equation*}
\begin{split}
b(x)
&:=\sum_{\vi\in[2^{\vec n}-\vec 1]}
\sum_{R\in \mathscr D(\mathbb R^{\vec n})}
\langle f, \theta^{(\vi)}_R\rangle
\prod_{\nu=1}^k \Bigg[\theta^{(i_\nu)}_{R_\nu}(x_\nu)
-\mathbf{1}_{\ell(R_\nu)\geq 1}\sum_{\gamma_\nu\in\mathbb N^{n_\nu},\,|\gamma_\nu|\leq N}
\frac{( \partial^{\gamma_\nu} \theta^{(i_\nu)}_{R_\nu})(\mathbf 0)}{\gamma_\nu!}  x_\nu^{\gamma_\nu}\Bigg]\\
&\phantom{:}=:\sum_{\vi\in[2^{\vec n}-\vec 1]}
\sum_{R\in \mathscr D(\mathbb R^{\vec n})}
\langle f, \theta^{(\vi)}_R\rangle
\prod_{\nu=1}^k g_{R_\nu}(x_\nu).
\end{split}
\end{equation*}
We aim to prove that
\begin{equation}\label{int<infty}
\int_{\mathbb R^{\vec n}} |b(x)| \prod_{\nu=1}^k \frac{1}{(1+|x_\nu|)^K}\,dx
<\infty
\end{equation}
for some $K\in(n,\infty)$. To this end, we estimate
$\langle g_{R_\nu}, (1+|\cdot|)^{-K}\rangle$
by consider two cases for $j_\nu:=-\log_2 \ell(R_\nu)$.

\emph{Case (1)} $j_\nu> 0$. In this case, the second term in $g_{R_\gamma}$ vanishes.
Then, since $\theta^{(i_\nu)}\in\mathscr S_0(\mathbb R^{n_\nu})$, we find that, for every $M_1\in (\frac32n,\infty)$,
\begin{equation*}
  |g_{R_\nu}(x_\nu)|
  \lesssim 2^{j_\nu \frac{n_\nu}2} \frac{1}{(1+2^{j_\nu}|x_\nu-a_{R_\nu}|)^{M_1}}
  \sim 2^{-j_\nu n_\nu} \frac{2^{\frac32 j_\nu  n_\nu}}{(1+2^{\frac32j_\nu}|x_\nu-a_{R_\nu}|)^{\frac23M_1}},
\end{equation*}
where the implicit positive constants depend only on $M_1$ and $\theta^{(i_\nu)}$.
Thus, for every $K\in(n,\infty)$,
\begin{equation}\label{j>0}
\int_{\mathbb R^{n_\nu}} \frac{|g_{R_\nu}(x_\nu)|}{(1+|x_\nu|)^K} \,dx_\nu
\leq 2^{-j_\nu n_\nu} (1+|a_{R_\nu}|)^{-\min(\frac23M_1,K)}.
\end{equation}

\emph{Case (2)} $j_{\nu}\leq 0$. In this case,
from the mean value theorem and $\theta^{(i_\nu)}\in\mathscr S_0(\mathbb R^{n_\nu})$,
we infer that, for every $M_2\in (n,\infty)$,
\begin{equation*}
\begin{split}
|g_{R_\nu}(x_\nu)|
&\lesssim |x_\nu|^{N+1} \sup_{\gamma_\nu\in\mathbb N^{n_\nu},\,|\gamma_\nu|=N+1}
\sup_{t\in(0,1)}
\Big|\partial^{\gamma_\nu} \theta^{(i_\nu)}_{R_\nu}(t x_\nu)\Big|\\
&\lesssim |x_\nu|^{N+1} 2^{j_\nu (\frac{n_\nu}2 + N+1)}
\sup_{t\in(0,1)} \frac{1}{(1+2^{j_\nu}|t x_\nu-a_{R_\nu}|)^{M_2}},
\end{split}
\end{equation*}
where the implicit positive constants depend only on $M_2$, $N$, and $\theta^{(i_\nu)}$.
Note that
\begin{equation*}
\frac{1}{1+2^{j_\nu}|t x_\nu-a_{R_\nu}|}
\leq \frac{1+(1-t)2^{j_\nu}|x_\nu|}{1+2^{j_\nu}|t x_\nu-a_{R_\nu}|+(1-t)2^{j_\nu}|x_\nu|}
\leq \frac{1+|x_\nu|}{1+2^{j_\nu}| x_\nu-a_{R_\nu}|}
\end{equation*}
and hence
\begin{equation*}
|g_{R_\nu}(x_\nu)|
\lesssim (1+|x_\nu|)^{N+M_2+1} 2^{j_\nu (\frac{n_\nu}2 + N+1)}
\frac{1}{(1+2^{j_\nu}|x_\nu-a_{R_\nu}|)^{M_2}}.
\end{equation*}
Thus, for every $K\in(N+M_2+n+1,\infty)$,
\begin{equation}\label{j<0}
\int_{\mathbb R^{n_\nu}} \frac{|g_{R_\nu}(x_\nu)|}{(1+|x_\nu|)^K} \,dx_\nu
\leq 2^{j_\nu (\frac{n_\nu}2 + N+1)}  (1+2^{j_\nu}|a_{R_\nu}|)^{-\min(M_2,K-N-M_2-1)}.
\end{equation}

Using \eqref{j>0} and \eqref{j<0}, we obtain
\begin{equation*}
\int_{\mathbb R^{\vec n}} |b(x)| \prod_{\nu=1}^k \frac{1}{(1+|x_\nu|)^K}
\leq \sum_{\vi\in[2^{\vec n}-\vec 1]}
\sum_{R\in \mathscr D(\mathbb R^{\vec n})}
|\langle f, \theta^{(\vi)}_R\rangle|
b_{[0,1)^{\vec n}, R},
\end{equation*}
where $b_{[0,1)^{\vec n}, R}:=\prod_{\nu=1}^k b_{[0,1)^{n_\nu}, R_\nu}$ with
\begin{equation*}
b_{[0,1)^{n_\nu}, R_\nu}:=
\begin{cases}
\displaystyle\ell(R_\nu)^{-(\frac{n_\nu}2 + N+1)}  (1+\frac{|a_{R_\nu}|}{\ell(R_\nu)})^{-\min(M_2,K-N-M_2-1)},
&\text{if } \ell(R_\nu)\geq 1, \\
\ell(R_\nu)^{n_\nu} (1+|a_{R_\nu}|)^{-\min(\frac23M_1,K)},
&\text{if } \ell(R_\nu)< 1.
\end{cases}
\end{equation*}
Note that $t^{(\vi)}:=\{t_R^{(\vi)}:=|\langle f, \theta^{(\vi)}_R\rangle|\}_{R\in\mathscr D(\mathbb R^{\vec n})}
\in \dot f^{2,\frac12}_{2,2}(V)$
and $B:=\{b_{P,R}:= b_{[0,1)^{\vec n}, R} \mathbf 1_{P=[0,1)^{\vec n}}\}_{P,R\in\mathscr D(\mathbb R^{\vec n})}$
is $(\vD,\vE,\vF)$-almost diagonal, where
$$
\begin{cases}
\displaystyle D_\nu:=\min\Big(\frac23 M_1,M_2,K-N-M_2-1\Big),\\
\displaystyle E_\nu:=\frac{n_\nu}2 + N+1, \quad\text{and}\quad
F_\nu:=n_\nu.
\end{cases}
$$
By first fixing $M_1,M_2,N$ large enough, and then $K$ large enough, it is clear that can achieve
$$
\vec D>\frac 32 \vec n,\quad
\vec E>\vec n,\quad\text{and}\quad
\vec F>\frac{\vec n}{2}.
$$
By Theorem \ref{AD st}, this implies that $B$ is bounded on $\dot f^{0,\frac12}_{2,2}(V)$, and hence
$$
\int_{\mathbb R^{\vec n}} |b(x)| \prod_{\nu=1}^k \frac{1}{(1+|x_\nu|)^K}
\leq \sum_{\vi\in[2^{\vec n}-\vec 1]} (B t^{(\vi)})_{[0,1)^{\vec n}}
<\infty.
$$
This finishes the proof of \eqref{int<infty}.

By \eqref{int<infty}, we find that $b$ is a locally integrable function on $\mathbb R^{\vec n}$
and $b\in \mathscr S'(\mathbb R^{\vec n})$.
Moreover, from the definition of $b$, we infer that
$b=f$ in $\mathscr S_0'(\mathbb R^{\vec n})$.
Together with \eqref{eq:BMO(V)=f} and Theorem \ref{wavelet}, this implies that
$\|b\|_{\mathrm{BMO}(V)}\sim\|f\|_{\dot F^{2,\frac12}_{2,2}(V)}<\infty$,
which completes the proof of Theorem \ref{coin BMO}.
\end{proof}

\subsection{One-parameter Besov--Triebel--Lizorkin-type spaces} \label{BTL}

Generalising the classical Besov and Triebel--Lizorkin spaces $\dot B^s_{p,q}$, $\dot F^s_{p,q}$ with three indices, the notion of Besov-type and Triebel--Lizorkin-type spaces $\dot B^{s,\tau}_{p,q}$, $\dot F^{s,\tau}_{p,q}$ incorporating the fourth Morrey index $\tau$ was first introduced by Yang et al. \cite{YY08,YY10,YSY10} in the one-parameter setting of $\R^n$. Matrix-weighted versions of these spaces have been subsequently studied by Bu et al. \cite{BHYY:1a,BHYY:1b,BHYY:1c,BHYY:2a,BHYY:2b}. We recall their definition with a slightly different notation to emphasise both similarities and differences with the spaces that we study in this work:

\begin{definition}
Let $s\in\mathbb R$, $\tau\in[0,\infty)$,
$p\in(0,\infty)$, $q\in(0,\infty]$, $\pi\in S_{[2]}$,
and $V\in\mathscr A_p(\mathbb R^n)$.
The \emph{homogeneous matrix-weighted Besov--Triebel--Lizorkin-type space associated with cubes}
$\dot A^{s,\tau,\mathrm{cube}}_{p,q,\pi}(V)$
is defined by setting
$$
\dot A^{s,\tau,\mathrm{cube}}_{p,q,\pi}(V)
:=\big\{f\in \mathscr S'(\mathbb R^n;\mathbb C^m):\
\|f\|_{\dot A^{s,\tau,\mathrm{cube}}_{p,q,\pi}(V)}<\infty\big\}
$$
where
\begin{equation}\label{eq:Acube}
\|f\|_{\dot A^{s,\tau,\mathrm{cube}}_{p,q,\pi}(V)}
:=\sup_{P\in\mathscr D(\R^n)} |P|^{-\tau}
\Big\|\Big\{\one_P 2^{js}
V\varphi_j*f_j\Big\}_{j\geq j_P}\Big\|_{(L^p\ell^q)_{\pi}}.
\end{equation}
\end{definition}

Note that $S_{[2]}$ consists of just two permutations. We recall from Example \ref{ex adm} the notation for them suggesting their relationships to Besov and Triebel--Lizorkin spaces, namely $F:=\operatorname{id}$ and $B:=(1\ 2)$ [the cycle notation for the unique non-trivial permutation of $\{1,2\}$ with $B(1)=2$ and $B(2)=1$].

A natural question is whether
the space $\dot A^{s,\tau,\mathrm{cube}}_{p,q,\pi}(V)$
coincides with the space $\dot A^{s,\tau}_{p,q,\pi}(V)$ defined by taking the supremum in \eqref{eq:Acube} over all open sets $\Omega$ instead. This coincidence is obvious for $\tau=0$, in which case the supremum is reached in the limit of $\Omega\to\R^n$ (or $P$ approaching one of the quadrants that gives the largest norm); thus both cubic and open-set definitions are extensions of the classical Besov and Triebel--Lizorkin spaces.

The following theorem provides a complete answer to this question:

\begin{theorem}\label{compare A}
Let $s\in\mathbb R$, $\tau\in[0,\infty)$, $p\in(0,\infty)$,
$q\in(0,\infty]$, $\pi\in S_{[2]}$, and $V\in \mathscr A_p(\mathbb R^n)$.
If $\tau\in(0,\frac 1p)$ or if $\tau=\frac 1p$,
$q\in (0,p)$, and $\pi= B$, then
$$
\dot A^{s,\tau}_{p,q,\pi}(V)
\hookrightarrow
\dot A^{s,\tau,\mathrm{cube}}_{p,q,\pi}(V)
\ \ \text{and}\ \
\dot A^{s,\tau}_{p,q,\pi}(V)
\neq \dot A^{s,\tau,\mathrm{cube}}_{p,q,\pi}(V),
$$
otherwise
\begin{align}\label{eqOcube}
\dot A^{s,\tau}_{p,q,\pi}(V)
= \dot A^{s,\tau,\mathrm{cube}}_{p,q,\pi}(V)
\end{align}
with equivalent quasi-norms.
\end{theorem}

Due to the $\varphi$-transform characterisation of both
$\dot A^{s,\tau,\mathrm{cube}}_{p,q,\pi}(V)$ (proved in \cite[Theorem 3.29]{BHYY:1a})
and $\dot A^{s,\tau}_{p,q,\pi}(V)$ (Theorem \ref{phi}),
Theorem \ref{compare A} can be proved by dealing with their corresponding sequence spaces.

\begin{definition}
Let $s\in\mathbb R$, $\tau\in[0,\infty)$,
$p\in(0,\infty)$, $q\in(0,\infty]$, $\pi\in S_{[2]}$,
and $V\in\mathscr A_p(\mathbb R^n)$.
The \emph{homogeneous matrix-weighted Besov--Triebel--Lizorkin-type sequence space associated with cubes}
$\dot a^{s,\tau,\mathrm{cube}}_{p,q,\pi}(V)$
is defined to be the set of all sequences
$t:=\{t_Q\}_{Q\in\mathscr D(\mathbb R^n)}$ in $\mathbb C^m$ such that
\begin{equation*}
\|t\|_{\dot a^{s,\tau,\mathrm{cube}}_{p,q,\pi}(V)}
:=\sup_{P\in\mathscr D(\R^n)} |P|^{-\tau}
\Big\|\Big\{\one_P 2^{js}
V t_j\Big\}_{j\geq j_P}\Big\|_{(L^p\ell^q)_{\pi}}
<\infty.
\end{equation*}
\end{definition}

The corresponding result for sequence spaces is as follows.

\begin{theorem}\label{compare a}
Let $s\in\mathbb R$, $\tau\in[0,\infty)$, $p\in(0,\infty)$,
$q\in(0,\infty]$, $\pi\in S_{[2]}$, and $V\in \mathscr A_p(\mathbb R^n)$.
If $\tau\in(0,\frac 1p)$ or if $\tau=\frac 1p$,
$q\in (0,p)$, and $\pi=B$, then
$$
\dot a^{s,\tau}_{p,q,\pi}(V)
\hookrightarrow
 \dot a^{s,\tau,\mathrm{cube}}_{p,q,\pi}(V)
\ \ \text{and}\ \
\dot a^{s,\tau}_{p,q,\pi}(V)
\neq \dot a^{s,\tau,\mathrm{cube}}_{p,q,\pi}(V),
$$
otherwise
\begin{align}\label{eqcubeO}
\dot a^{s,\tau}_{p,q,\pi}(V)
= \dot a^{s,\tau,\mathrm{cube}}_{p,q,\pi}(V)
\end{align}
with equivalent quasi-norms.
\end{theorem}

To prove Theorem \ref{compare a}, we need two lemmas.
We first introduce two notations.
For every $\{f_j\}_{j\in\mathbb Z}$, let
$$
\|\{f_j\}_{j\in\mathbb Z}\|_{\dot\beta^\tau_{p,q,\pi}}
:=\sup_{P\in\mathscr D(\R^n)} |P|^{-\tau}
\Big\|\Big\{\one_{P} f_j\Big\}_{j\geq j_P}\Big\|_{(L^p\ell^q)_{\pi}}
$$
and
$$
\|\{f_j\}_{j\in\mathbb Z}\|_{\dot\gamma^\tau_{p,q,\pi}}
:=\sup_{\Omega\in\Open(\R^n)} |\Omega|^{-\tau}
\Big\|\Big\{\one_{\Omega_j}f_j\Big\}_{j\in\Z}\Big\|_{(L^p\ell^q)_{\pi}}.
$$

\begin{lemma} \label{equiv super}
Let $p\in(0,\infty)$, $q\in(0,\infty]$, $\tau\in[0,\infty)$,
and $\pi\in S_{[2]}$.
Assume that $(p,q,\tau,\pi)$ satisfies one of the following conditions:
\begin{enumerate}[\rm(i)]
\item $\tau\in(\frac1p,\infty)$;

\item $\tau=\frac1p$ and $q=\infty$;

\item $\tau=\frac1p$, $q\in(0,\infty)$, and $\pi=F$;

\item $\tau=\frac 1p$, $q\in[p,\infty)$, and $\pi=B$.
\end{enumerate}
Then, for all $\{f_j\}_{j\in\mathbb Z}$ with each $f_j$ constant on every dyadic cube in $\mathscr D_j(\R^n)$,
\begin{align*}
\|\{f_j\}_{j\in\mathbb Z}\|_{\dot\beta^\tau_{p,q,\pi}}
\sim \|\{f_j\}_{j\in\mathbb Z}\|_{\dot\gamma^\tau_{p,q,\pi}},
\end{align*}
where the positive equivalence constants depend only on $n$, $\tau$, and $p$.
\end{lemma}

\begin{proof}
The inequality $\|\{f_j\}_{j\in\mathbb Z}\|_{\dot\beta^\tau_{p,q,\pi}}
\leq \|\{f_j\}_{j\in\mathbb Z}\|_{\dot\gamma^\tau_{p,q,\pi}}$ is obvious.
Now, we prove the reverse inequality
by considering three cases for $(p,q,\tau,\pi)$.

\emph{Case (1)} $(p,q,\tau,\pi)$ satisfies (i) or (ii).
In this case, for every open set $\Omega$ with $|\Omega|\in(0,\infty)$,
\begin{align*}
\big\|\{\one_{\Omega_j}f_j\}_{j\in\Z}\big\|_{(L^p\ell^q)_{\pi}}
\leq \Big\|\Big\{\one_{\Omega_j}2^{-jn(\tau-\frac1p)}\Big\}_{j\in\Z}\Big\|_{(L^p\ell^q)_{\pi}}
\Big\|\Big\{2^{jn(\tau-\frac1p)}f_j\Big\}_{j\in\Z}\Big\|_{L^\infty\ell^\infty}
=:\mathrm{I}\times \mathrm{II}.
\end{align*}
Let $j_\Omega:=\min\{j\in\mathbb Z:\ \Omega_j \neq \varnothing \}$.
Then there exists $P\in \mathscr D_{j_\Omega}(\mathbb R^n)$
such that $P\subset\Omega$. These further imply that
\begin{align*}
\mathrm{I}
\leq \|\one_{\Omega}\|_{L^p} \Big\|\Big\{\mathbf 1_{j\geq j_\Omega}2^{-jn(\tau-\frac1p)}\Big\}_{j\in\Z}\Big\|_{\ell^q}
\sim |\Omega|^{\frac 1p} 2^{-j_\Omega n(\tau-\frac1p)}
= |\Omega|^{\frac 1p} |P|^{\tau-\frac1p}
\leq |\Omega|^\tau.
\end{align*}
Note that, for every $j\in\mathbb Z$ and $P\in\mathscr D_j(\mathbb R^n)$,
\begin{align*}
\Big\|2^{jn(\tau-\frac1p)}f_j\Big\|_{L^\infty(P)}
=|P|^{-\frac 1p} \Big\|2^{jn(\tau-\frac1p)}f_j\Big\|_{L^p(P)}
=|P|^{-\tau} \|f_j\|_{L^p(P)}
\leq \|\{f_j\}_{j\in\mathbb Z}\|_{\dot\beta^\tau_{p,q,\pi}},
\end{align*}
and hence $\mathrm{II}\leq\|\{f_j\}_{j\in\mathbb Z}\|_{\dot\beta^\tau_{p,q,\pi}}$.
Combining the above estimates, we obtain
$$\|\{f_j\}_{j\in\mathbb Z}\|_{\dot\gamma^\tau_{p,q,\pi}}
\lesssim\|\{f_j\}_{j\in\mathbb Z}\|_{\dot\beta^\tau_{p,q,\pi}}.$$

\emph{Case (2)} $(p,q,\tau,\pi)$ satisfies (iii).
In this case, let $\Omega$ be an open set and
$\{P_i\}_{i\in I}$ be its maximal dyadic decomposition.
Then, for every $i\in I$ and $j\in\mathbb Z$,
$\mathbf 1_{P_i\cap \Omega_j} = \mathbf 1_{j\geq j_{P_i}} \mathbf 1_{P_i}$,
which further implies that
\begin{align*}
\big\|\{\one_{\Omega_j}f_j\}_{j\in\Z}\big\|_{L^p\ell^q}^p
&= \sum_{i\in I} \int_{P_i} \Bigg( \sum_{j\in\mathbb Z}
\one_{\Omega_j} |f_j|^q \Bigg)^{\frac pq}
= \sum_{i\in I} \big\|\{\one_{P_i} f_j\}_{j\geq j_{P_i}}\big\|_{L^p\ell^q}^p \\
&\leq \sum_{i\in I} |P_i| \, \|\{f_j\}_{j\in\mathbb Z}\|_{\dot\beta^\tau_{p,q,\pi}}^p
=|\Omega| \, \|\{f_j\}_{j\in\mathbb Z}\|_{\dot\beta^\tau_{p,q,\pi}}^p,
\end{align*}
and hence $\|\{f_j\}_{j\in\mathbb Z}\|_{\dot\gamma^\tau_{p,q,\pi}}
\leq \|\{f_j\}_{j\in\mathbb Z}\|_{\dot\beta^\tau_{p,q,\pi}}$.

\emph{Case (3)} $(p,q,\tau,\pi)$ satisfies (iv).
In this case, let $\Omega$ be an open set and
$\{P_i\}_{i\in I}$ be its maximal dyadic decomposition.
Then, for every $i\in I$ and $j\in\mathbb Z$,
$\mathbf 1_{P_i\cap \Omega_j} = \mathbf 1_{j\geq j_{P_i}} \mathbf 1_{P_i}$,
which, together with Minkowski's inequality, further implies that
\begin{align*}
\big\|\{\one_{\Omega_j}f_j\}_{j\in\Z}\big\|_{\ell^qL^p}^p
&= \Bigg[ \sum_{j\in\mathbb Z} \Bigg( \sum_{i\in I} \int_{P_i}
\one_{\Omega_j} |f_j|^p \Bigg)^{\frac qp} \Bigg]^{\frac pq}
\leq  \sum_{i\in I} \Bigg[ \sum_{j\in\mathbb Z} \Bigg( \int_{P_i}
\one_{\Omega_j} |f_j|^p \Bigg)^{\frac qp} \Bigg]^{\frac pq}  \\
&= \sum_{i\in I} \Big[ \big\|\{\one_{P_i} f_j\}_{j\geq j_{P_i}}
\big\|_{\ell^q L^p}^q \Big]^{\frac pq}
\leq \sum_{i\in I} |P_i| \, \|\{f_j\}_{j\in\mathbb Z}\|_{\dot\beta^\tau_{p,q,\pi}}^p
=|\Omega| \, \|\{f_j\}_{j\in\mathbb Z}\|_{\dot\beta^\tau_{p,q,\pi}}^p,
\end{align*}
and hence $\|\{f_j\}_{j\in\mathbb Z}\|_{\dot\gamma^\tau_{p,q,\pi}}\leq \|\{f_j\}_{j\in\mathbb Z}\|_{\dot\beta^\tau_{p,q,\pi}}$.
This finishes the proof of Lemma \ref{equiv super}.
\end{proof}

\begin{lemma} \label{equiv sub}
Let $p\in(0,\infty)$.
If $\tau\in(0,\frac1p)$, $q\in(0,\infty]$, and $\pi\in S_{[2]}$
or if $\tau=\frac1p$, $q\in(0,p)$, and $\pi=B$,
then, for all $\{f_j\}_{j\in\mathbb Z}$ with each $f_j$ constant on every dyadic cube in $\mathscr D_j(\R^n)$,
\begin{align*}
\|\{f_j\}_{j\in\mathbb Z}\|_{\dot\beta^\tau_{p,q,\pi}}
\leq\|\{f_j\}_{j\in\mathbb Z}\|_{\dot\gamma^\tau_{p,q,\pi}}.
\end{align*}
However, the reverse inequality fails, that is,
there exists $\{f_j\}_{j\in\mathbb Z}$ such that
$\|\{f_j\}_{j\in\mathbb Z}\|_{\dot\beta^\tau_{p,q,\pi}}<\infty$
but $\|\{f_j\}_{j\in\mathbb Z}\|_{\dot\gamma^\tau_{p,q,\pi}}=\infty$.
\end{lemma}

\begin{proof}
The inequality ``$\|\{f_j\}_{j\in\mathbb Z}\|_{\dot\beta^\tau_{p,q,\pi}}
\leq\|\{f_j\}_{j\in\mathbb Z}\|_{\dot\gamma^\tau_{p,q,\pi}}$'' is obvious.
We need only to show that the reverse inequality fails.
To this end, we first construct a super-sparse sequence of dyadic cubes.
For every $\kappa\in\mathbb N$, let $Q_\kappa$
be a dyadic cube of edge-length $1$, contained in $[0,2^\kappa)^n$
and located in its ``upper right'' corner.
Next, we prove that the reverse inequality fails
by considering two cases for $(\tau,q,\pi)$.

\emph{Case (1)} $\tau\in(0,\frac1p)$, $q\in(0,\infty]$, and $\pi\in S_{[2]}$ is arbitrary.
In this case, let $f_0:=\mathbf 1_{\bigcup_{\kappa\in\mathbb N}Q_\kappa}$
and $f_j:=0$ for all $j\in\mathbb Z\setminus\{0\}$.
We first prove $\|\{f_j\}_{j\in\mathbb Z}\|_{\dot\beta^\tau_{p,q,\pi}}<\infty$.
By the construction of $\{Q_\kappa\}_{k\in\N}$,
we conclude that, for every $P\in\mathscr D(\mathbb R^n)$ with $|P|\geq 1$,
$$
\#\{\kappa\in\mathbb N:\
P\cap Q_\kappa \neq \varnothing \} \leq -j_P+1,
$$
which further implies that
\begin{align*}
|P|^{-\tau} \big\|\{\one_P f_j\}_{j\geq j_P}\big\|_{(L^p\ell^q)_\pi}
=|P|^{-\tau}
\| \mathbf 1_{P\cap(\bigcup_{\kappa\in\mathbb N}Q_\kappa)} \|_{L^p(\mathbb R^n)}
\leq |P|^{-\tau} (-j_P)^{\frac 1p}
\leq \sup_{i\in\mathbb N} 2^{-in\tau}i^{\frac 1p}<\infty.
\end{align*}
Therefore, $\|\{f_j\}_{j\in\mathbb Z}\|_{\dot\beta^\tau_{p,q,\pi}}<\infty$.

Next, we show $\|\{f_j\}_{j\in\mathbb Z}\|_{\dot\gamma^\tau_{p,q,\pi}}=\infty$.
For every $N\in\mathbb N$,
let $\Omega^{(N)}:=\bigcup_{\kappa=1}^N Q_\kappa$.
Then
\begin{align*}
|\Omega^{(N)}|^{-\tau}\Big\|\Big\{\one_{\Omega^{(N)}_j} f_j\Big\}_{j\in\Z}\Big\|_{(L^p\ell^q)_\pi}
=|\Omega^{(N)}|^{-\tau}\big\|\mathbf 1_{\Omega^{(N)}}\big\|_{L^p(\mathbb R^n)}
= \big|\Omega^{(N)}\big|^{\frac1p-\tau}
=N^{\frac1p-\tau}
\to \infty
\end{align*}
as $N\to\infty$.
Therefore, $\|\{f_j\}_{j\in\mathbb Z}\|_{\dot\gamma^\tau_{p,q,\pi}}=\infty$.
Thus, the desired conclusion hold in this case.

\emph{Case (2)} $\tau=\frac1p$, $q\in(0,p)$, and $\pi= B$.
In this case,
we can choose a sequence of dyadic cubes $\{R_\kappa\}_{\kappa\in\mathbb N}$ that satisfies
\begin{enumerate}[\rm(i)]
\item for every $\kappa\in\mathbb N$,
$R_\kappa$ is a subset of $Q_\kappa$ and

\item for every $j\in\mathbb N$,
$E_j:=\bigcup\{R_\kappa:\ \ell(R_\kappa) = 2^{-j} \}$
satisfies $|E_j| \sim j^{-a}$, where $a\in(1,\frac pq]$.
\end{enumerate}
Let $f_j:=\mathbf 1_{E_j}$
for all $j\in\mathbb N$;
$f_j:=0$ for all $j\in\mathbb Z\setminus\mathbb N$.
We first prove $\|\{f_j\}_{j\in\mathbb Z}\|_{\dot\beta^\tau_{p,q,\pi}}<\infty$.
For every $P\in\mathscr D(\mathbb R^n)$ with $|P|< 1$,
if $P\cap(\bigcup_{j\in\mathbb N} E_j)=\varnothing$, then
$$
|P|^{-\frac1p} \big\|\{\one_P f_j\}_{j\geq j_P}\big\|_{\ell^q L^p}=0;
$$
otherwise there exists a dyadic cube $Q$ such that
$P \cap (\bigcup_{j\in\mathbb N} E_j) = Q$ and hence
$$
|P|^{-\frac1p} \big\|\{\one_P f_j\}_{j\geq j_P}\big\|_{\ell^q L^p}
\leq |P|^{-\frac1p} \|\mathbf 1_Q\|_{L^p}
= |P|^{-\frac1p} |Q|^{\frac 1p}
\leq 1.
$$
For every $P\in\mathscr D(\mathbb R^n)$ with $|P|\geq 1$,
$P\cap(\bigcup_{j\in\mathbb N} E_j)$ is a union of at most $-j_P+1$ dyadic cubes
from the family $\{R_\kappa\}_{\kappa\in\mathbb N}$,
which further implies that
\begin{align*}
|P|^{-\frac1p} \big\|\{\one_P f_j\}_{j\geq j_P}\big\|_{\ell^q L^p}
&= |P|^{-\frac1p} \Bigg(\sum_{j\in\mathbb N} |P\cap E_j|^{\frac qp}\Bigg)^{\frac 1q}\\
&\leq 2^{j_Pn\frac1p}
\Bigg[\sum_{j\in\mathbb N} (-j_P+1)^{\frac qp}
2^{-jn\frac qp} \Bigg]^{\frac 1q}
\lesssim \sup_{i\in\N} \bigg(\frac{i+1}{2^{in}}\bigg)^{\frac 1p}<\infty.
\end{align*}
Therefore, $\|\{f_j\}_{j\in\mathbb Z}\|_{\dot\beta^\tau_{p,q,\pi}}<\infty$.

Next, we show $\|\{f_j\}_{j\in\mathbb Z}\|_{\dot\gamma^\tau_{p,q,\pi}}=\infty$.
Let $\Omega:=\bigcup_{j\in\mathbb N} E_j$.
Then
$$
|\Omega|=\sum_{j\in\mathbb Z_+}|E_j|
\sim \sum_{j\in\mathbb N} \frac1{j^a} <\infty
$$
and
\begin{align*}
\big\|\{\one_{\Omega_j} f_j\}_{j\in\Z}\big\|_{\ell^q L^p}^q
=\sum_{j\in\mathbb N} |E_j|^{\frac qp}
=\sum_{j\in\mathbb N} j^{-a\frac qp}
=\infty.
\end{align*}
Therefore, $\|\{f_j\}_{j\in\mathbb Z}\|_{\dot\gamma^\tau_{p,q,\pi}}=\infty$.
Thus, the desired conclusions hold also in this case, which completes
the proof of Lemma \ref{equiv sub}.
\end{proof}

Now, we can compare the space $\dot a^{s,\tau,\mathrm{cube}}_{p,q,\pi}(V)$
with the space $\dot a^{s,\tau}_{p,q,\pi}(V)$.

\begin{proof}[Proof of Theorem \ref{compare a}]
We first prove
$$
\dot a^{s,\tau}_{p,q,\pi}(V)
\hookrightarrow
\dot a^{s,\tau,\mathrm{cube}}_{p,q,\pi}(V)
\ \ \text{and}\ \
\dot a^{s,\tau}_{p,q,\pi}(V)
\neq \dot a^{s,\tau,\mathrm{cube}}_{p,q,\pi}(V)
$$
in the cases where $\tau\in(0,\frac 1p)$ or where $\tau=\frac 1p$,
$q\in (0,p)$, and $\pi=B$.
To this end, due to Theorem \ref{3 norms thm} and its analogue for the cubic space proved in \cite[Theorem 3.27]{BHYY:1a}, it suffices to show that,
for every $t:=\{t_Q\}_{Q\in\mathscr D(\mathbb R^n)}$ with each $t_Q\in\mathbb C^m$,
\begin{equation}\label{9/13-1}
\|\{2^{js} [V]_{j,p} t_j\}_{j\in\mathbb Z}\|_{\dot\beta^\tau_{p,q,\pi}}
\leq\|\{2^{js} [V]_{j,p} t_j\}_{j\in\mathbb Z}\|_{\dot\gamma^\tau_{p,q,\pi}}
\end{equation}
and there exists $t:=\{t_Q\}_{Q\in\mathscr D(\mathbb R^n)}$ such that
\begin{equation}\label{9/13-2}
\|\{2^{js} [V]_{j,p} t_j\}_{j\in\mathbb Z}\|_{\dot\beta^\tau_{p,q,\pi}}<\infty
\ \ \text{but}\ \
\|\{2^{js} [V]_{j,p} t_j\}_{j\in\mathbb Z}\|_{\dot\gamma^\tau_{p,q,\pi}}=\infty.
\end{equation}
By Lemma \ref{equiv sub}, we conclude that
\eqref{9/13-1} holds and there exists $\{f_j\}_{j\in\mathbb Z}$ such that
$$
\|\{f_j\}_{j\in\mathbb Z}\|_{\dot\beta^\tau_{p,q,\pi}}<\infty
\ \ \text{but}\ \
\|\{f_j\}_{j\in\mathbb Z}\|_{\dot\gamma^\tau_{p,q,\pi}}=\infty.
$$
Choose $t:=\{t_Q\}_{Q\in\mathscr D(\mathbb R^n)}$ such that
$2^{js} [V]_{j,p} t_j=f_j$ for all $j\in\mathbb Z$.
Then \eqref{9/13-2} holds.
For the remaining parameter ranges, that is, for all
$(\tau,p,q,\pi)$ not covered by the above cases,
using Lemma \ref{equiv super} with $\{f_j\}_{j\in\Z}$ replaced by
$\{2^{js}[V]_{j,p}f_{j}\}_{j\in\Z}$,
and then applying \cite[Theorem 3.27]{BHYY:1a} and Theorem \ref{3 norms thm},
we obtain equality \eqref{eqcubeO}.
This finishes the proof of Theorem \ref{compare a}.
\end{proof}

Using Theorem \ref{compare a}, we can
compare the corresponding function spaces, that is,
we can finish the proof of Theorem \ref{compare A}.

\begin{proof}[Proof of Theorem \ref{compare A}]
We first prove
$$
\dot A^{s,\tau}_{p,q,\pi}(V)
\hookrightarrow
\dot A^{s,\tau,\mathrm{cube}}_{p,q,\pi}(V)
\ \ \text{and}\ \
\dot A^{s,\tau}_{p,q,\pi}(V)
\neq \dot A^{s,\tau,\mathrm{cube}}_{p,q,\pi}(V)
$$
in the cases where $\tau\in(0,\frac 1p)$ or where $\tau=\frac 1p$,
$q\in (0,p)$, and $\pi=B$.
The embedding is obvious.
On the other hand, Theorem \ref{compare a} shows that
there exists $t:=\{t_Q\}_{Q\in\mathscr D(\mathbb R^n)}\in\dot a^{s,\tau,\mathrm{cube}}_{p,q,\pi}(V) \setminus \dot a^{s,\tau}_{p,q,\pi}(V)$.
Let $\psi$ be a Littlewood--Paley function. The $\varphi$-transform characterisation of the cubic space $\dot a^{s,\tau,\mathrm{cube}}_{p,q,\pi}(V)$, namely \cite[Theorem 3.29]{BHYY:1a}, which is the $\dot a^{s,\tau,\mathrm{cube}}_{p,q,\pi}(V)$ analogue of Theorem \ref{phi},  shows that
$f:=T_\psi t=\sum_{Q\in\mathscr D(\mathbb R^n)} t_Q \psi_Q$ converges in $\mathscr S'_0(\mathbb R^n;\mathbb C^m)$ with
\begin{equation*}
   \|f\|_{\dot A^{s,\tau,\mathrm{cube}}_{p,q,\pi}(V)}
  \sim\|t\|_{\dot a^{s,\tau,\mathrm{cube}}_{p,q,\pi}(V)}
 <\infty.
\end{equation*}
However, the $\varphi$-transform characterisation of $\dot A^{s,\tau}_{p,q,\pi}(V)$ in Theorem \ref{phi}  implies that
\begin{align*}
\|f\|_{\dot A^{s,\tau}_{p,q,\pi}(V)}
\sim \|t\|_{\dot a^{s,\tau}_{p,q,\pi}(V)}
=\infty,
\end{align*}
and hence $f\in \dot A^{s,\tau,\mathrm{cube}}_{p,q,\pi}(V)\setminus \dot A^{s,\tau}_{p,q,\pi}(V)$.

For the remaining parameter ranges, i.e., for all
$(\tau,p,q,\pi)$ not covered by the above cases,
equality \eqref{eqOcube} follows from Theorem \ref{compare a} and
the same $\varphi$-transform characterisations of both spaces $\dot A^{s,\tau}_{p,q,\pi}(V)$
and $\dot A^{s,\tau,\mathrm{cube}}_{p,q,\pi}(V)$.
This finishes the proof of Theorem \ref{compare A}.
\end{proof}

\subsection{Rectangular multi-parameter Besov--Triebel--Lizorkin-type spaces} \label{BTL another}

In this subsection, we will study a question similar to that in Subsection \ref{BTL}, but in the general multi-parameter setting: What happens if, in the definition of multi-parameter Besov--Triebel--Lizorkin-type spaces, we take the supremum over all dyadic rectangles only, instead of all open sets?
We begin with a definition:

\begin{definition}
Let $\vec s\in\mathbb R^k$, $\tau\in[0,\infty)$,
$\vec p\in(0,\infty)^k$, $\vec q\in(0,\infty]^k$, $\pi\in S_{[2k]}$,
and $V\in L^p_{\mathrm{loc}}(\mathbb R^{\vec n};\mathbb C^{m\times m})$.
The \emph{homogeneous matrix-weighted Besov--Triebel--Lizorkin-type sequence space associated with rectangles}
$\dot a^{\vec s,\tau,\mathrm{rect}}_{\vec p,\vec q,\pi}(\mathbb R^{\vec n})$
is defined to be the set of all sequences
$t:=\{t_Q\}_{Q\in\mathscr D(\mathbb R^{\vec n})}$ in $\mathbb C^m$ such that
\begin{equation*}
\|t\|_{\dot a^{\vec s,\tau,\mathrm{rect}}_{\vec p,\vec q,\pi}(V)}
:=\sup_{P\in\mathscr D(\R^{\vec n})} |P|^{-\tau}
\Big\|\Big\{\one_P 2^{\vj\cdot\vec s}
V t_{\vj}\Big\}_{\vj\geq {\vj}_P}\Big\|_{(L^{\vec p}\ell^{\vec q})_{\pi}}
<\infty.
\end{equation*}
\end{definition}

Let us compare the spaces $\dot a^{\vec s,\tau,\mathrm{rect}}_{\vec p,\vec q,\pi}(V)$
and $\dot a^{\vec s,\tau}_{\vec p,\vec q,\pi}(V)$.
To limit the length of this article, we only consider the special cases
where $\vec p=p\cdot \vec 1$, $\vec q=q\cdot \vec 1$,
$\pi \in \{B,F\}$, and $V\equiv I_m$.
Note that $(L^{\vec p}\ell^{\vec q})_{B}=\ell^{\vec q}L^{\vec p}$ and $(L^{\vec p}\ell^{\vec q})_{F}=L^{\vec p}\ell^{\vec q}$.
In the multi-parameter setting, these two spaces differ in more cases than those in the one-parameter setting.
In particular,
$$
\mathrm{BMO}(\mathbb R^{\vec n})
=\dot F^{0,\frac12}_{2,2}(\mathbb R^{\vec n})
\neq \dot F^{0,\frac12,\mathrm{rect}}_{2,2}(\mathbb R^{\vec n})
$$
when $k\geq2$, whereas all three spaces coincide when $k=1$.

\begin{theorem} \label{compare B}
Let $k\geq 2$, $\vec s\in\mathbb R^k$, $\tau\in(0,\infty)$,
$\vec p:=p\cdot \vec 1\in(0,\infty)^k$,
$\vec q:=q\cdot \vec 1\in(0,\infty]^k$, and $\pi \in \{B,F\}$.
If $\tau\in(0,\frac1p)$
or if $\tau=\frac1p$, $q< p$, and $\pi=B$
or if $\tau=\frac1p$ and $q= p$, then
$$
\dot a^{\vec s,\tau}_{\vec p,\vec q,\pi}(\mathbb R^{\vec n})
\hookrightarrow
\dot a^{\vec s,\tau,\mathrm{rect}}_{\vec p,\vec q,\pi}(\mathbb R^{\vec n})
\ \ \text{and}\ \
\dot a^{\vec s,\tau}_{\vec p,\vec q,\pi}(\mathbb R^{\vec n})
\neq \dot a^{\vec s,\tau,\mathrm{rect}}_{\vec p,\vec q,\pi}(\mathbb R^{\vec n}).
$$
\end{theorem}

In the proof of this theorem, we will use a slight extension of the classical counterexample constructed by Carleson \cite{Car}, see Lemma \ref{Carleson}.

\begin{proof}[Proof of Theorem \ref{compare B}]
The embedding $\dot a^{\vec s,\tau}_{\vec p,\vec q,\pi}(\mathbb R^{\vec n})
\hookrightarrow
\dot a^{\vec s,\tau,\mathrm{rect}}_{\vec p,\vec q,\pi}(\mathbb R^{\vec n})$ is obvious;
it remains to prove that
$\dot a^{\vec s,\tau}_{\vec p,\vec q,\pi}(\mathbb R^{\vec n})
\neq \dot a^{\vec s,\tau,\mathrm{rect}}_{\vec p,\vec q,\pi}(\mathbb R^{\vec n})$.
To this end, we consider two cases on $(p,q,\tau)$.

\emph{Case (1)} $\tau\in(0,\frac 1p)$, or $\tau=\frac 1p$,
$q\in (0,p)$, and $\pi= B$.
In this case, Theorem \ref{compare a} shows that
$$
\dot a^{s_1,\tau}_{p,q,\pi}(\mathbb R^{n_1})
\hookrightarrow \dot a^{s_1,\tau,\mathrm{rect}}_{p,q,\pi}(\mathbb R^{n_1})\quad \text{and} \quad
\dot a^{s_1,\tau}_{p,q,\pi}(\mathbb R^{n_1})
\neq \dot a^{s_1,\tau,\mathrm{rect}}_{p,q,\pi}(\mathbb R^{n_1}).
$$
Thus, there exists $u:=\{u_{Q^*}\}_{Q^*\in\mathscr D(\mathbb R^{n_1})}\in \dot a^{s_1,\tau,\mathrm{rect}}_{p,q,\pi}(\mathbb R^{n_1})$
such that $u\notin \dot a^{s_1,\tau}_{p,q,\pi}(\mathbb R^{n_1})$. Let
$$
t_Q:=\begin{cases}
u_{Q^*},&\text{if } Q=Q^*\times [0,1)^{n_2+\cdots+n_k},\\
0,&\text{otherwise}.
\end{cases}
$$
Then
\begin{align*}
\|t\|_{\dot a^{\vec s,\tau,\mathrm{rect}}_{\vec p,\vec q,\pi}(\mathbb R^{\vec n})}
=\|u\|_{\dot a^{s_1,\tau,\mathrm{rect}}_{p,q,\pi}(\mathbb R^{n_1})}
<\infty
\ \ \text{and} \ \
\|t\|_{\dot a^{\vec s,\tau}_{\vec p,\vec q,\pi}(\mathbb R^{\vec n})}
=\|u\|_{\dot a^{s_1,\tau}_{p,q,\pi}(\mathbb R^{n_1})}
=\infty.
\end{align*}
This finishes the proof of $\dot a^{\vec s,\tau}_{\vec p,\vec q,\pi}(\mathbb R^{\vec n})
\neq \dot a^{\vec s,\tau,\mathrm{rect}}_{\vec p,\vec q,\pi}(\mathbb R^{\vec n})$ in this case.

\emph{Case (2)} $\tau=\frac 1p$ and $q=p$.
In this case, for every $N\in\mathbb Z_+$,
let $\{R_i\}_{i\in \mathscr I_N}$ be as in Lemma \ref{Carleson}
with $\varepsilon=\frac1N$ and let
$$
t_Q:=\begin{cases}
2^{-\vj_Q \cdot (\vec s + \frac{\vec n}2)},
&Q=R_i+ (N-1)\cdot \vec 1
\text{ for some } N\in\mathbb Z_+\text{ and }i\in \mathscr I_N,\\
0,&\text{otherwise}.
\end{cases}
$$
Let $\mathscr I:=\bigcup_{N\in\mathbb Z_+}\mathscr I_N$.
Then, for every $P\in\mathscr D(\mathbb R^{\vec n})$,
\begin{align*}
|P|^{-\tau}
\Big\|\Big\{\one_P 2^{\vj\cdot\vec s}
t_{\vj}\Big\}_{\vj\geq {\vj}_P}\Big\|_{(L^{\vec p}\ell^{\vec p})_\pi}
=|P|^{-\tau} \Big[\sum_{i\in \mathscr I,R_i\subset P}
|R_i|\Big]^{\frac 1p}
\leq |P|^{-\tau} |P|^{\frac 1p}
=1,
\end{align*}
using property \eqref{it:rectCarleson} of Lemma \ref{Carleson} in the inequality.
Hence $t\in \dot a^{\vec s,\tau,\mathrm{rect}}_{\vec p,\vec q,\pi}(\mathbb R^{\vec n})$.

On the other hand, for every $N\in\mathbb Z_+$,
let $\Omega^{(N)}:=\bigcup_{i\in\mathscr I_N} R_i+ (N-1)\cdot \vec 1 $. Then
\begin{equation*}
  \Omega^{(N)}_{\vj}
  :=\bigcup_{R\in\mathscr D_{\vj}:R\subset\Omega}R
  \supset\bigcup_{i\in\mathscr I_N:R_i\in\mathscr D_{\vj}} R_i+ (N-1)\cdot \vec 1
\end{equation*}
so that
\begin{equation*}
  \one_{\Omega^{(N)}_{\vj}} 2^{\vj\cdot\vec s}t_{\vj}
  \geq \sum_{i\in\mathscr I_N:R_i\in\mathscr D_{\vj}}
    2^{\vj\cdot\vec s}\cdot 2^{-\vj\cdot(\vec s+\frac{\vn}{2})}\widetilde{\one}_{R_i}
   =\sum_{i\in\mathscr I_N:R_i\in\mathscr D_{\vj}}\one_{R_i}.
\end{equation*}
Since the pairwise distinct rectangles $R_i$ in a fixed $\mathscr D_{\vj}$ are pairwise disjoint, it follows that
\begin{equation*}
  \Norm{\one_{\Omega^{(N)}_{\vj}} 2^{\vj\cdot\vec s}t_{\vj}}{L^p}^p
  \geq \sum_{i\in\mathscr I_N:R_i\in\mathscr D_{\vj}}\abs{R_i}
\end{equation*}
and hence
\begin{equation*}
  \Big\|\Big\{\one_{\Omega^{(N)}_{\vj}} 2^{\vj\cdot\vec s}
    t_{\vj}\Big\}_{\vj\in\mathbb Z^k}\Big\|_{\ell^{\vec p}L^{\vec p}}^p
   \geq
   \sum_{\vj\in\Z^k}\sum_{i\in\mathscr I_N:R_i\in\mathscr D_{\vj}}\abs{R_i}
   =\sum_{i\in\mathscr I_N}\abs{R_i}=1,
\end{equation*}
using property \eqref{it:Car-ex-i} of Lemma \ref{Carleson} in the last step. Thus
\begin{align*}
|\Omega^{(N)}|^{-\tau}
\Big\|\Big\{\one_{\Omega^{(N)}_{\vj}} 2^{\vj\cdot\vec s}
t_{\vj}\Big\}_{\vj\in\mathbb Z^k}\Big\|_{\ell^{\vec p}L^{\vec p}}
\geq\abs{\Omega^{(N)}}^{-\tau}
>N^{\tau}
\end{align*}
by property \eqref{it:Car-ex-ii} of Lemma \ref{Carleson} with $\eps=\frac{1}{N}$ in the last step.
As  $N\to\infty$, this shows that $t\notin \dot a^{\vec s,\tau}_{\vec p,\vec q,\pi}(\mathbb R^{\vec n})$,
and finishes the proof of Theorem \ref{compare B}.
\end{proof}

In \cite[Theorem 5.9]{FJ90}, corresponding to the classical Hardy-BMO duality,
Frazier and Jawerth obtained an duality relation for Triebel--Lizorkin spaces:
$$
[\dot f_{1,q}^{s}(\mathbb R^n)]^*=\dot f_{q',q'}^{-s,\frac1{q'},\mathrm{cube}}(\mathbb R^n).
$$
It is natural to guess that this duality also holds for multi-parameter spaces.
The following result indicates that this is true for sequence spaces associated
with open set, but not true for sequence spaces associated with cubes, which
provides another motivation for studying the space
$\dot a^{\vec s,\tau}_{\vec p,\vec q,\pi}(V)$.

\begin{corollary} \label{dual}
If $k=2$, $\vec s\in \mathbb R^2$, and $q\in(1,\infty)$, then
$$
\big[\dot f_{1,q}^{\vec s}(\mathbb R^{\vec n})\big]^*
=\dot f_{q',q'}^{-\vec s,\frac1{q'}}(\mathbb R^{\vec n})
\neq \dot f_{q',q'}^{-\vec s,\frac1{q'},\mathrm{rect}}(\mathbb R^{\vec n}).
$$
\end{corollary}

\begin{proof}
By \cite[Theorem 4.5]{DingZhu17}, we obtain the first ``$=$'',
while the ``$\neq$'' follows from Theorem \ref{compare B}.
This finishes the proof of Corollary \ref{dual}.
\end{proof}

\section{The $T(1)$ theorem}\label{sec:CZO}

In this section, we obtain a general criterion, in the spirit of a $T(1)$ theorem, for the boundedness of multi-parameter singular integrals on the matrix-weighted multi-parameter Triebel--Lizorkin spaces $\Atau(V)$. We recall that the first $T(1)$ theorem for the boundedness of one-parameter Calder\'on--Zygmund operators on $L^2(\R^n)$ was proved by David and Journ\'e \cite{DJ84}, and the extension to multi-parameter Calder\'on--Zygmund (also known as Calder\'on--Zygmund--Journ\'e, or simply Journ\'e) operators on $L^2(\R^{\vn})$  is due to Journ\'e \cite{Journe}. The name ``$T(1)$ theorem'' derives from conditions of the type ``$T(1)\in\BMO$'' in the general form of these theorems, or ``$T(1)=0$'' in the important special case of ``paraproduct-free'' operators, where $T(1)$ refers to the action of the operator $T$ on the constant function $1$.

For Triebel--Lizorkin spaces the one-parameter setting, after a number of results in a limited range of indices (see \cite[pp.\ 222--223]{HNVW3} for a more complete history), extensions to the full range of indices were  obtained by Frazier, Torres, and Weiss \cite{FTW,tor}, and their matrix-weighted extension was only recently by Bu et al. \cite{BHYY:1c}. 
 In the matrix-weighted multi-parameter setting, Domelevo et al.~\cite{DKPS} proved a paraproduct-free $T(1)$ theorem in $L^p(V)$ and a general $T(1)$ theorem in $L^2(V)$, which was extended to a general $T(1)$ theorem in $L^p(V)$ by Vuorinen \cite{Vuo24}.

In this section, we prove a matrix-weighted multi-parameter $T(1)$ theorem on the general scale of $\Atau(V)$ spaces, Theorem \ref{T1 BF}, which essentially extends both \cite{BHYY:1c,DKPS} in the paraproduct-free case; for the fact that the $L^p(V)$ spaces are contained as a special case of $\Atau(V)$, see Theorem \ref{Lp=F}. We note that the restriction to the paraproduct-free case is made in most of the literature on $T(1)$ theorems in general Triebel--Lizorkin spaces, even in the unweighted one-parameter setting. With the principal goal of extending existing theory to the matrix-weighted multi-parameter setting, we make this same restriction here. The added complication of paraproducts is left for future work and should probably be first addressed in the one-parameter Triebel--Lizorkin spaces, which is already an interesting question in its own right.

Towards the goal of stating and proving Theorem \ref{T1 BF}, extending the approach of \cite{BHYY:1c} from the one-parameter case, we first give an abstract criterion for extending operators
$$T\in\mathcal L\Big(\mathscr{S}\big(\mathbb R^{\vec n}\big),\mathscr{S}'(\mathbb R^{\vec n})\Big)
\quad\mathrm{to}\quad \widetilde T\in\mathcal L\Big(\dot A^{\vec s,\tau}_{\vec p,\vec q,\pi}(V)\Big)$$
by testing the action of $T$ on suitable {\em atoms}. The main body of the section is then dedicated to verifying this criterion in the concrete case of multi-parameter singular integral operators.

As a first step, Subsection \ref{CZop} is devoted to a streamlined one-parameter set-up, which allows us to
treat the multi-parameter case by induction on the number of parameters in Subsection \ref{CZmp}, where we complete the proof of the multi-parameter $T(1)$ theorem.

Let $C_{\mathrm{c}}^\infty(\mathbb R^{\vec n})$ denote
the set of all functions $f\in C^\infty(\mathbb R^{\vec n})$ with compact support.

\begin{definition}\label{def:atom}
Let $\vec L,\vec N\in\N^k$.
A function $a_Q\in C_{\mathrm{c}}^\infty(\mathbb R^{\vec n})$
is called an \emph{$(\vec L,\vec N)$-atom} on a rectangle $Q:=Q_1\times\cdots\times Q_k\in\Rect(\mathbb R^{\vec n})$ if
\begin{enumerate}[\rm(i)]
\item\label{it:aForm}
$a_Q= \bigotimes_{\nu=1}^k a_{Q_\nu}$,

\item\label{it:aSupp}
for every $\nu\in[k]$,
$\operatorname{supp}a_{Q_\nu}\subset3Q_\nu$,

\item\label{it:aCanc}
for every $\nu\in[k]$ and $\gamma_\nu\in\N^{n_\nu}$ with $|\gamma_\nu|\leq L_\nu$,
$\int_{\mathbb{R}^{n_\nu}}x^{\gamma_\nu} a_{Q_\nu}(x)\,dx_\nu=0$, and

\item\label{it:aSize}
for every $\nu\in[k]$,
$\gamma_\nu\in\N^{n_\nu}$ with $|\gamma_\nu|\leq N_\nu$,
and $x_\nu\in\mathbb{R}^{n_\nu}$,
$|D^{\gamma_\nu} a_{Q_\nu}(x_\nu)|\leq|Q_\nu|^{-\frac12-\frac{|\gamma_\nu|}{n_\nu}}$.
\end{enumerate}
\end{definition}

The following proposition is a straightforward multi-parameter extension of \cite[Proposition 6.5(i)]{BHYY:1c};
we omit the details.

\begin{proposition}\label{ext}
Let $\vp\in(0,\infty)^k$, $\vq\in(0,\infty]^k$, $\vs\in\R^k$, $\tau\in[0,\infty)$,
$\pi\in S_{[2k]}$ be admissible for $(\vp,\vq)$,
and $V$ be a matrix weight such that $(\vp,\vq,\pi,V)$ belongs
to one of the main cases (Definition \ref{main cases}).
Let $\vec L,\vec N\in(0,\infty)^k$.
If $T\in\mathcal L(\mathscr{S}(\mathbb R^{\vec n}),\mathscr{S}'(\mathbb R^{\vec n}))$
maps $(\vec L,\vec N)$-atoms to $\dot A^{\vec s,\tau}_{\vec p,\vec q,\pi}$-synthesis molecules (Definition \ref{def:AtauMole}),
then there exists an operator $\widetilde T\in\mathcal L(\dot A^{\vec s,\tau}_{\vec p,\vec q,\pi}(V))$ that agrees with $T$ on $\mathscr{S}_0(\mathbb R^{\vec n};\mathbb C^m)$.
\end{proposition}

If $T\in\mathcal L(\mathscr S(\R^{\vn}),\mathscr S'(\R^{\vn}))$,
then, by the well-known Schwartz kernel theorem,
there exists a unique $K\in\mathscr S'
(\mathbb R^{\vn}\times\mathbb R^{\vn})$
such that, for all $\varphi,\psi\in\mathscr S(\R^{\vn})$,
\begin{equation}\label{TvsK}
  \pair{T\varphi}{\psi}=\pair{K}{\psi\otimes\varphi}.
\end{equation}
This $K$ is called the \emph{Schwartz kernel} of $T$.

Thanks to the kernel theorem, in order to define a distribution $K\in\mathscr S'
(\mathbb R^{\vn}\times\mathbb R^{\vn})$, it is enough to specify $\pair{K}{\psi\otimes\varphi}$ for every pair of functions $\psi,\varphi\in\mathscr S(\R^{\vn})$ in such a way that $(\psi,\varphi)\mapsto \pair{K}{\psi\otimes\varphi}$ is bi-linear and separately continuous; indeed, this clearly defines an operator $T\in\mathcal L(\mathscr S(\R^{\vn}),\mathscr S'(\R^{\vn}))$ via \eqref{TvsK}, which then determines the full distribution $K\in\mathscr S'
(\mathbb R^{\vn}\times\mathbb R^{\vn})$ by the kernel theorem. We will use this observation without explicit mention in what follows.

\subsection{A review of the one-parameter case}\label{CZop}

Our approach to the multi-parameter results below will proceed by induction on the number $k$ of parameters. In order to carry this out efficiently, it will be convenient to have a reasonably streamlined set-up for the one-parameter results available. For this purpose, we collect and reformulate the relevant results from \cite[Sections 6.2--6.3]{BHYY:1c}

The following is a reformulation of \cite[Definition 6.11]{BHYY:1c}, introducing several auxiliary notions to streamline the core part of the definition.
(For the sake of simplicity, we do not make the fine distinction between two kernel classes $\operatorname{CZK}^1(E,F)\subset\operatorname{CZK}^0(E,F)$, present in \cite{BHYY:1c}, but simply consider $\operatorname{CZK}(E,F):=\operatorname{CZK}^1(E,F)$.)

\begin{definition}\label{def:CZK1}
For $E,F\in\R$, let
\begin{equation}\label{eq:GEF}
\begin{split}
  \mathscr G(E,F):=
  \left.\begin{cases}
    (\gamma,\delta,\zeta,\rho)\in\N^n\times\N^n\times\{0,1\}\times\{0,1\}: \\
    \quad\abs{\gamma}\leq\Floor{E}_+,\quad\delta=\mathbf 0,\quad(\zeta,\rho)=(0,0), & \text{or}\quad \\
    \quad\abs{\gamma}=\Floor{E}\phantom{_+},\quad\delta=\mathbf 0,\quad(\zeta,\rho)=(1,0), & \text{or} \\
    \quad\abs{\gamma}\leq\Floor{E}_+,\quad\abs{\delta}=\Floor{F}-\abs{\gamma},\quad(\zeta,\rho)=(0,1), & \text{or} \\
    \quad\abs{\gamma}=\Floor{E}\phantom{_+},\quad\abs{\delta}=\Floor{F-E},\quad(\zeta,\rho)=(1,1)
    \end{cases}\right\}.
\end{split}
\end{equation}
For $\Gamma=(\gamma,\delta,\zeta,\rho)\in\mathscr G(E,F)$ and $\upsilon\in(0,1)$, let
\begin{equation*}
  \R^{4n}_{\Gamma,\upsilon}:=
  \Big\{(X,Y)=((x,y),(u,v))\in\R^{4n}:\ u\neq\mathbf 0\neq v,\ \zeta\abs{u}+\rho\abs{v}<\upsilon\abs{x-y}\Big\}
\end{equation*}
and, for $U=(u,v)\in\R^{2n}$,
\begin{equation*}
  \mathbb D_{U}^\Gamma:=\Delta_{U}^{(\zeta,\rho)}\partial_x^\gamma\partial_y^\delta,
\end{equation*}
where $\Delta_{U}^{(\zeta,\rho)}$ is defined in $\eqref{Delta}$.

For $\Gamma=(\gamma,\delta,\zeta,\rho)\in\mathscr G(E,F)$ and $(X,U)=((x,y),(u,v))\in\R^{4n}$, let
\begin{equation*}
  \operatorname{CZ}_{E,F}^{\Gamma}(X,U)
  :=\abs{x-y}^{n+\abs{\gamma}+\abs{\delta}}\bigg(\frac{\abs{x-y}}{\abs{u}}\bigg)^{\zeta E^{**}}
  \bigg(\frac{\abs{x-y}}{\abs{v}}\bigg)^{\rho(F-\zeta E)^{**}}.
\end{equation*}

For $K\in\mathscr S'(\R^n\times\R^n)$, we say $K\in\operatorname{CZK}(E,F)$ if $K(x,y)$ and all its partial derivatives appearing below are continuous functions on $x\neq y$ and the following quantity is finite:
\begin{equation*}
  \Norm{K}{\operatorname{CZK}_\upsilon(E,F)}
  :=\sup_{\genfrac{}{}{0pt}{}{\Gamma\in\mathscr G(E,F)}{(X,U)\in\R^{4n}_{\Gamma,\upsilon}}}
  \operatorname{CZ}_{E,F}^{\Gamma}(X,U)\abs{\mathbb D_U^\Gamma K(X)}.
\end{equation*}
\end{definition}

\begin{remark}\label{rem:GEF}
The length of a multi-index $\gamma,\delta\in\N^n$ is obviously non-negative.
In \eqref{eq:GEF}, there are several conditions of the type $\abs{\gamma}=s$ or $\abs{\delta}=t$. If $s<0$ or $t<0$, then we simply understand that there are no such multi-indices.
\end{remark}

\begin{remark}\label{rem:upsilon}
The parameter $\upsilon\in(0,1)$ is usually taken simply as $\upsilon=\frac12$. Since the differences $\Delta_u f(x)=f(x)-f(x-u)$ can always be expressed as sums of differences of as much smaller step size as we like, namely,
\begin{equation*}
  \Delta_u f(x)=f(x)-f(x-u)
  =\sum_{i=1}^N f(x+\tfrac{i-1}{N}u)-f(x+\tfrac{i}{N}u)
  =\sum_{i=1}^N \Delta_{\frac uN}f(x+\tfrac{i-1}{N}u),
\end{equation*}
it is easy to see that the $\operatorname{CZK}_{\upsilon}(E,F)$ norms are equivalent for all $\upsilon\in(0,1)$; thus, the kernel class $\operatorname{CZK}(E,F)$ is independent of this parameter.
\end{remark}

\begin{remark}\label{rem:CZE}
We denote $\operatorname{CZK}(E):=\operatorname{CZK}(E,0)$ and
$\operatorname{CZK}^*(F):=\{K:K^*\in\operatorname{CZK}(F)\}$, where $\pair{K^*}{\psi\otimes\varphi}:=\pair{K}{\varphi\otimes\psi}$.
It is evident that $\operatorname{CZK}(E,F)\subset\operatorname{CZK}(E)$ for all $F\in\R$.
While not evident, we also have $\operatorname{CZK}(E,F)\subset\operatorname{CZK}^*(F)$ for all $E\in\R$; see \cite[(6.17)]{BHYY:1c}.
\end{remark}

\begin{remark}\label{rem:CZK1}
For technical reasons, it will be convenient to use the Calder\'on--Zygmund estimates in a slightly different but equivalent form. Note that, if $\abs{x'-x}+\abs{y'-y}<\hat\upsilon\abs{x-y}$, then
\begin{equation*}
   \abs{x-y}\leq\abs{x-x'}+\abs{x'-y'}+\abs{y'-y}<\hat\upsilon\abs{x-y}+\abs{x'-y'},
\end{equation*}
and thus $(1-\hat\upsilon)\abs{x-y}<\abs{x'-y'}$. Hence, if $\abs{u},\abs{v}<\upsilon\abs{x-y}$, then $\abs{u},\abs{v}<\upsilon'\abs{x'-y'}$ where $\upsilon':=\frac{\upsilon}{1-\hat\upsilon}$ still satisfies $\upsilon'\in(0,1)$ for $\hat\upsilon\in(0,1-\upsilon)$.

Thus, if $(X,U)\in\R^{4n}_{\Gamma,\upsilon}$, then $(X',U)\in\R^{4n}_{\Gamma,\upsilon'}$ for all $X'=(x',y')$ as above. If $K\in\operatorname{CZK}(E,F)$ and $\Gamma\in\mathscr G(E,F)$, it follows that
\begin{equation*}
   \abs{\mathbb D^\Gamma_U K(X')}
   \leq\Norm{K}{\operatorname{CZK}_{\upsilon'}(E,F)}
   \operatorname{CZ}^{\Gamma}_{E,F}(X',U)^{-1}
\end{equation*}
and it is straightforward that
\begin{equation*}
  \operatorname{CZ}^{\Gamma}_{E,F}(X',U)
     \sim\operatorname{CZ}^{\Gamma}_{E,F}(X,U)
\end{equation*}
for the ranges of the variables under consideration. Hence $\mathbb D^\Gamma_U K(X')$ has essentially the same bound as $\abs{\mathbb D^\Gamma_U K(X)}$.

If $\phi\in C_{\rm c}^\infty(\R^n)$ is a usual bump function of integral one supported in the unit ball, and $\Phi^X_{\eps}:=\eps^{-2n}(\phi\otimes\phi)(\eps^{-1}(\cdot-X))$, then $\pair{\mathbb D_U^{\Gamma}K}{\Phi^X_\eps}$, being an average of the values $\mathbb D_U^{\Gamma}K(X')$, retains the same upper bound for all $\eps\in(0,\hat\upsilon\abs{x-y})$. On the other hand,
\begin{equation*}
  \pair{\mathbb D_U^{\Gamma}K}{\Phi^X_\eps}
  \underset{\eps\to 0}{\longrightarrow} \mathbb D_U^{\Gamma}K(X).
\end{equation*}
Hence we see that
\begin{equation*}
\begin{split}
  \Norm{K}{\operatorname{CZK}_\upsilon(E,F)}
  :&=\sup_{\genfrac{}{}{0pt}{}{\Gamma\in\mathscr G(E,F)}{(X,U)\in\R^{4n}_{\Gamma,\upsilon}}}
  \operatorname{CZ}_{E,F}^{\Gamma}(X,U)\abs{\mathbb D_U^\Gamma K(X)} \\
  &\leq \sup_{\genfrac{}{}{0pt}{}{\genfrac{}{}{0pt}{}{\Gamma\in\mathscr G(E,F)}
   {(X,U)\in\R^{4n}_{\Gamma,\upsilon}}}{0<\eps<\hat\upsilon\abs{x-y}}}
  \operatorname{CZ}_{E,F}^{\Gamma}(X,U)\abs{\pair{\mathbb D_U^\Gamma K}{\Phi_\eps^X}}
  \lesssim \Norm{K}{\operatorname{CZK}_{\upsilon'}(E,F)}.
\end{split}
\end{equation*}
Since the different $\operatorname{CZK}_\upsilon(E,F)$-norms are equivalent, we see that the quantity sandwiched between them above is also equivalent.

Fixing some values of the auxiliary parameters, say $\upsilon=\frac12$ and $\hat\upsilon=\frac13<\frac12=1-\upsilon$, we adopt this as a definition of
\begin{equation*}
  \Norm{K}{\operatorname{CZK}(E,F)}
  :=\sup_{\genfrac{}{}{0pt}{}{\Gamma\in\mathscr G(E,F)}{(X,U,\eps)\in\R^{4n+1}_{\Gamma}  }}
  \operatorname{CZ}_{E,F}^{\Gamma}(X,U) \abs{ \pair{\mathbb D_U^\Gamma K }{ \Phi_\eps^X }},
\end{equation*}
where
\begin{equation*}
  \R^{4n+1}_{\Gamma} := \Big\{ (X,U,\eps)\in\R^{4n+1}:\  (X,U)=((x,y),(u,v))\in\R^{4n}_{\Gamma,\frac12},\ 0<\eps<\tfrac13\abs{x-y}\Big\}.
\end{equation*}
\end{remark}

\begin{definition}\label{def:UAI}
We say that $(\chi_j)_{j=0}^\infty\subset C_{\rm c}^\infty(\R^n)$ is a {\em uniform approximation to the identity} if $\sup_{j\in\N}\Norm{\chi_j}{L^\infty}<\infty$ and
for every compact $K\subset\R^n$, there is an index $j_K\in\N$ for which $\chi_j(x)=1$ for all $x\in K$ and $j\geq j_K$.
\end{definition}

\begin{lemma}\label{lem:T1}
Let $K\in\operatorname{CZK}(E,F)$. Then the limits
\begin{equation}\label{eq:T1lim}
  \pair{K}{f\otimes y^{\mu}}
  :=\lim_{j\to\infty}\pair{K}{f\otimes\chi_j y^{\mu}},\qquad
  \pair{K}{x^{\lambda}\otimes g}
  :=\lim_{j\to\infty}\pair{K}{\chi_j x^{\lambda}\otimes g}
\end{equation}
exist and are independent of the uniform approximation to the identity $(\chi_j)_{j=0}^\infty$, for all $\abs{\mu}\leq\Floor{E}$, $\abs{\lambda}\leq\Floor{F}$,
 and $f\in\mathcal D_{\abs{\mu}}$, $g\in\mathcal D_{\abs{\lambda}}$, where
 \begin{equation*}
  \mathcal D_{\ell}:=\bigg\{ g\in C_{\rm c}^\infty(\R^n): \int_{\R^n}x^{\gamma}g(x)dx =0 \text{ for all }\abs{\gamma}\leq\ell\bigg\}.
\end{equation*}
\end{lemma}

\begin{proof}
Since $\operatorname{CZK}(E,F)\subset\operatorname{CZK}(E)$ (Remark \ref{rem:CZE}), the claim concerning the first limit follows from \cite[Lemma 6.9]{BHYY:1c} (which is a special case of \cite[Lemma 2.2.12]{tor}).
Noting $\operatorname{CZK}(E,F)\subset\operatorname{CZK}^*(F)$ (Remark \ref{rem:CZE}), the other claim is similar.
\end{proof}

\begin{definition}\label{def:T1=0}
For $K\in\operatorname{CZK}(E,F)$ and $\abs{\mu}\leq\Floor{E}$, we say that
\begin{equation*}
  \pair{K}{\cdot\otimes y^\mu}:=0
\end{equation*}
if $\pair{K}{f\otimes y^\mu}$, as defined in \eqref{eq:T1lim}, vanishes for all $f\in\mathcal D_{\abs{\mu}}$.

For $\abs{\lambda}\leq\Floor{F}$, the condition $ \pair{K}{x^\lambda\otimes\cdot}=0$ is defined analogously.
\end{definition}

We also recall the notion of the weak boundedness property in the following formulation:

\begin{definition}\label{def:WBP1}
For functions $\eta,\xi\in C_{\rm c}^\infty(\R^n)$ and numbers $W\in\N$ and $r\in[0,\infty)$, let
\begin{equation*}
   \delta(\eta,\xi):=\operatorname{diam}(\supp \eta\cup\supp \xi),\qquad
  P(W,\eta,r):=\sum_{\abs{\alpha}\leq W} r^{\abs{\alpha}}\Norm{\partial^\alpha \eta}{L^\infty},
\end{equation*}
and
\begin{equation*}
  B_{W}(\eta,\xi)
  :=\delta(\eta,\xi)^{n}\cdot P(W,\eta,\delta(\eta,\xi))\cdot P(W,\xi,\delta(\eta,\xi)).
\end{equation*}
Let
\begin{equation*}
  \mathscr B(W):=\{\Xi=\eta\otimes\xi:\quad \eta,\xi\in C_{\rm c}^\infty(\R^n),\quad
  B_W(\eta,\xi)\leq 1\}
  \subset\mathscr S(\R^n\times\R^n).
\end{equation*}
We say that $K\in\mathscr S'(\mathbb R^{n}\times\mathbb R^{n})$ satisfies the \emph{weak boundedness property} of order $W$ if the following quantity is finite:
\begin{equation*}
  \Norm{K}{\mathrm{WBP}(W)}
   :=\sup_{\Xi\in\mathscr B(W)} \abs{\pair{K}{\Xi}}.
\end{equation*}
\end{definition}

\begin{remark}\label{rem:WBP1}
The weak boundedness property is often defined (see e.g. \cite[Definition 4]{Journe}) with reference to suitable translates and dilates of pairs of functions in bounded subsets of $C_{\rm c}^\infty$, but the condition above is equivalent, as pointed out in \cite[(1.9)]{Journe}; note that $\delta(\eta,\xi)$ is comparable to the edge-length of the smallest cube that contains both $\supp f$ and $\supp g$ appearing in \cite[(1.9)]{Journe}. While the parameter $W$ is largely irrelevant, the advantage of the quantitative formulation of Definition \ref{def:WBP1} is that it is amenable to iteration in subsequent considerations.

Another detail is that the weak boundedness is often defined as a property of an {\em operator} $T\in\mathcal L(\mathscr S(\R^n),\mathscr S'(\R^n))$, instead of its {\em kernel} $K\in\mathscr S'(\R^n\times\R^n)$. This is quite natural in one-parameter considerations, where the weak boundedness and the kernel estimates play essentially distinct roles from each other. However, for a smoother discussion of the multi-parameter setting, it seems convenient to define both Calder\'on--Zygmund and weak boundedness estimates as properties of the Schwartz kernel $K$.
\end{remark}

The following is essentially a restatement of \cite[Definition 6.17]{BHYY:1c}:

\begin{definition}\label{def:CZEFGH}
Let $E,F\in\R$, $G,H\in\Z$, and $W\in\N$, with $G\leq\Floor{E}$ and $H\leq\Floor{F}$.
For $K\in\mathscr S'(\R^n\times\R^n)$, we say that $K\in\operatorname{CZK}(E,F,G,H,W)$, if
\begin{enumerate}[\rm(i)]
  \item $K$ satisfies the weak boundedness property of order $W$;
  \item $K\in\operatorname{CZK}(E,F)$;
  \item\label{it:T1=0G} $\pair{K}{\cdot\otimes y^\mu}=0$ for all $\abs{\mu}\leq G$,
  \item\label{it:T*1=0H} $\pair{K}{x^\lambda\otimes\cdot}=0$ for all $\abs{\lambda}\leq H$.
\end{enumerate}
We define
\begin{equation*}
\begin{split}
  \Norm{K}{\operatorname{CZK}(E,F,W)}
  &:=\Norm{K}{\operatorname{CZK}(E,F)}\vee\Norm{K}{\operatorname{WBP}(W)}, \\
  \Norm{K}{\operatorname{CZK}(E,F,G,H,W)}
  &:=\begin{cases}
   \Norm{K}{\operatorname{CZK}(E,F,W)}, &
    \text{if \eqref{it:T1=0G} and \eqref{it:T*1=0H} hold}, \\
    \infty, & \text{otherwise}. \end{cases}
\end{split}
\end{equation*}
\end{definition}

\begin{remark}\label{rem:CZEFGH}
We have incorporated the restrictions $G\leq\Floor{E}$ and $H\leq\Floor{F}$ into Definition \ref{def:CZEFGH}, since conditions \eqref{it:T1=0G} and \eqref{it:T*1=0H} are only guaranteed to make sense under these assumptions. While these restrictions were not stated in \cite[Definition 6.17]{BHYY:1c}, one can see that the actual results based on that definition, namely \cite[Theorem 6.18, Proposition 6.19, and Corollaries 6.20--6.21]{BHYY:1c}, are compatible with these restrictions.
\end{remark}

\begin{definition}\label{def:moleNorm}
For $M\in[0,\infty)$, $N\in(0,\infty)$, and $L\in\Z_{\geq-1}$, we define
\begin{equation}\label{A(N)}
  \mathscr A(N):=\{(\alpha,\kappa)\in\N^n\times\{0,1\}:\quad\abs{\alpha}\leq\Floor{N}_+,\ \kappa=0,\quad\text{or}\quad\abs{\alpha}=\Floor{N},\ \kappa=1\},
\end{equation}
and
\begin{equation*}
  \Norm{m}{\operatorname{Mol}(M,N)}
  :=  \sup_{ \substack{ (\alpha,\kappa)\in\mathscr A(N) \\ \abs{h}\leq 1 }}
   \abs{h}^{-N^{**}\kappa}  \BNorm{(1+\abs{\cdot}^2)^{\frac{M}{2}}
      \Delta_h^{\kappa}\partial^\alpha m }{L^\infty}
\end{equation*}
as well as
\begin{equation*}
  \Norm{m}{\operatorname{Mol}(L,M,N)}
  :=\begin{cases} \Norm{m}{\operatorname{Mol}(M,N)}, &
  \text{if}\quad\int_{\R^n} x^\gamma m(x)dx=0\quad\text{for all }\abs{\gamma}\leq L, \\
  \infty, & \text{otherwise}.\end{cases}
\end{equation*}
\end{definition}

\begin{remark}\label{rem:moleNorm}
Comparison with Definition \ref{moleKLMN} shows that $\Norm{m}{\operatorname{Mol}(L,M,N)}\leq 1$ is essentially equivalent to the condition that $m$ is an $(M,L,M,N)$-molecule on $[0,1)^n$.
[More precisely, for some $c\in(0,1)$ the condition $\Norm{m}{\operatorname{Mol}(L,M,N)}\leq c$ implies that $m$ is an $(M,L,M,N)$-molecule, which in turn implies that $\Norm{m}{\operatorname{Mol}(L,M,N)}\leq C$ for some constant $C\in(1,\infty)$.]

Note in particular that $(1+\abs{x}^2)^{\frac12}\sim 1+\abs{x}$. In Definition \ref{def:moleNorm}, we choose the $C^\infty$ factor $(1+\abs{x}^2)^{\frac12}$ for technical convenience in later considerations.
\end{remark}

The key result that we borrow from the one-parameter theory is the following elaboration of the statement of \cite[Proposition 6.19]{BHYY:1c}:

\begin{proposition}\label{prop:6.19}
Let $E,F,M\in\R$, $G,H\in\Z$, $L\in\Z_{\geq-1}$,
$N\in(0,\infty)$, and $W,L',N'\in\N$
satisfy
\begin{equation}\label{eq:6.19ass}
  \begin{cases} E\geq N, \\ E>\floor{N}_+,\end{cases}
  \begin{cases} F\geq M-n, \\ F> L,\end{cases}
  \begin{cases} G\geq \floor{N}_+, \\ H\geq L,\end{cases}
  \begin{cases} L'\geq \Floor{F}, \\ N'\geq \floor{G}+1.\end{cases}
\end{equation}
Let $T\in\mathcal L(\mathscr S(\R^n),\mathscr S'(\R^n))$ have Schwartz kernel $K\in\operatorname{CZK}(E,F,G,H,W)$. Then $T$ maps $(L',N')$-atoms into $(M,L,M,N)$-molecules. In particular, for $(L',N')$-atoms $\varphi$ on $[0,1)^n$, we have
\begin{equation}\label{eq:6.19}
  \Norm{T\varphi}{\operatorname{Mol}(L,M,N)}
  \lesssim\Norm{K}{\operatorname{CZK}(E,F,G,H,W)}\Norm{\varphi}{\operatorname{Atom}(L',N')},
\end{equation}
where
\begin{equation*}
  \Norm{\varphi}{\operatorname{Atom}(L',N')}
  :=\begin{cases}
    \sup_{\abs{\alpha}\leq N'}\Norm{\partial^\alpha\varphi}{L^\infty}, & \text{if }\int_{\R^n}x^\gamma\varphi(x)dx = 0
    \text{ for all }\abs{\gamma}\leq L', \\
    \infty, & \text{otherwise}\end{cases}
\end{equation*}
and the implicit positive constant is independent of $T$ and $\varphi$.
\end{proposition}

\begin{remark}\label{rem:6.19}\mbox{}
\begin{enumerate}[\rm(i)]
\item The quantitative statement \eqref{eq:6.19} did not explicitly appear in \cite[Proposition 6.19]{BHYY:1c}, but it easily follows by inspection of the argument.
\item The conditions on $L',N'$ were not part of the statement of \cite[Proposition 6.19]{BHYY:1c}, but they appear right in the beginning of its proof. We have slightly simplified the conditions appearing there, namely $L'\geq\Floor{F}\vee\Floor{F-E}$ and $N'\geq(\Floor{E}\wedge\floor{G})+1$, which is justified as follows: First, since $E>\floor{N}_+\geq 0$ by assumption, we see that $\Floor{F}\geq\Floor{F-E}$ always dominates the maximum. Second, in Definition \ref{def:CZEFGH} (as explained in Remark \ref{rem:CZEFGH}), we have incorporated the restriction that $G\leq\Floor{E}$, and hence also $\floor{G}\leq\Floor{E}$.
\end{enumerate}
\end{remark}

For induction of the number of parameters, it will be convenient to use Proposition \ref{prop:6.19} in the following dualised form:

\begin{corollary}\label{cor:6.19}
Under the assumptions of Proposition \ref{prop:6.19}, for all $\psi\in C_{\rm c}^\infty(\R^n)$, we have
\begin{equation*}
  \sup_{\genfrac{}{}{0pt}{}{(\alpha,\kappa)\in\mathscr A(N)}{\abs{h}\leq 1}}\abs{h}^{-N^{**}\kappa}
  \abs{\pair{(1+\abs{x}^2)^{\frac{M}{2}}\Delta_{(h,0)}^{(\kappa,0)}\partial_x^\alpha K}{\psi\otimes \varphi}}
  \lesssim\Norm{\psi}{L^1}\Norm{K}{\operatorname{CZK}(E,F,G,H,W)}\Norm{\varphi}{\operatorname{Atom}(L',N')}.
\end{equation*}
\end{corollary}

\begin{proof}
Using the identity
\begin{equation*}
\begin{split}
  \pair{(1+\abs{x}^2)^{\frac{M}{2}}\Delta_{(h,0)}^{(\kappa,0)} \partial_x^\alpha K}{\psi\otimes \varphi}
  &=\pair{K}{((-\partial_x)^\alpha \Delta_{-h}^{\kappa} [(1+\abs{x}^2)^{\frac{M}{2}}\psi])\otimes \varphi} \\
  &=\pair{T\varphi}{(-\partial_x)^\alpha\Delta_{-h}^{\kappa}[(1+\abs{x}^2)^{\frac{M}{2}}\psi]} \\
  &=\pair{(1+\abs{x}^2)^{\frac{M}{2}}\Delta_{h}^{\kappa}\partial_x^\alpha T\varphi}{\psi}
\end{split}
\end{equation*}
and the $L^1$--$L^\infty$-duality, Corollary \ref{cor:6.19} immediately follows from Proposition \ref{prop:6.19}.
\end{proof}

\subsection{The multi-parameter case}\label{CZmp}

For all $\vn:=(n_\nu)_{\nu\in[k]}$ and $I\subset[k]$,
let $\vn_I:=(n_{\nu})_{\nu\in I}$.

Starting with the pioneering work of \cite{Journe}, the idea of multi-parameter operators of Calder\'on--Zygmund type is to impose all possible combinations of Calder\'on--Zygmund kernel conditions (Definition \ref{def:CZK1}) and weak boundedness conditions (Definition \ref{def:WBP1}) over all complementary sets of variables:

\begin{definition}\label{def:CZKk}
Let $\vec E,\vec F\in\R^k$ and $\vec W\in\N$.
For $K\in\mathscr S'(\R^{\vn}\times\R^{\vn})$, we say that
$K\in \mathrm{CZK}(\vec E,\vec F,\vec W)$ if $K(x,y)$ and all its partial derivatives appearing below are continuous functions on
$x_\nu\neq y_\nu$ with all $\nu\in[k]$ and the following quantity is finite:
\begin{equation*}
  \Norm{K}{\operatorname{CZK}(\vE,\vF,\vec W)}
  :=
  \sup_{\genfrac{}{}{0pt}{}{[k]=I\cup J}{I\cap J=\varnothing}}
  \Norm{K}{\operatorname{CZK}_{I,J}(\vE,\vF,\vec W)},
\end{equation*}
where
\begin{equation*}
 \Norm{K}{\operatorname{CZK}_{I,J}(\vE,\vF,\vec W)}
  :=\sup_{\genfrac{}{}{0pt}{}{\Gamma_I\in\mathscr G(\vE,\vF)_I}{\genfrac{}{}{0pt}{}{(X,U,\vec\eps)\in\R^{4\vn_I+\abs{I}}_{\Gamma}}{
    \Xi_{J}\in \mathscr B(\vec W)_J}}}
  \operatorname{CZ}_{\vE,\vF}^{\Gamma;I}(X,U)
  \pair{\mathbb D^{\Gamma;I}_U K}{(\Phi_{\vec\eps}^X)_I\otimes \Xi_{J}}
\end{equation*}
and here
\begin{equation*}
  \mathscr G(\vE,\vF)_I
  :=\prod_{\nu\in I}\mathscr G(\vE_\nu,\vF_\nu),\qquad
  \R^{4\vn_I+\abs{I}}_{\Gamma}
  :=\prod_{\nu\in I}\R^{4n_\nu+1}_{\Gamma_\nu},
\end{equation*}
\begin{equation*}
\begin{split}
  \mathscr B(\vec W)_J
  &:=\Bigg\{\Xi_J=\bigotimes_{\nu\in J}\Xi_\nu:
  \Xi_\nu=(\eta_\nu\otimes\xi_\nu),\ \eta_\nu,\xi_\nu\in C_{\rm c}^\infty(\R^{n_\nu}),
  B_{W_\nu}(\eta_\nu,\xi_\nu)\leq 1\Bigg\},
\end{split}
\end{equation*}
and
\begin{equation*}
  \operatorname{CZ}_{\vE,\vF}^{\Gamma;I}(X,U)
  :=\prod_{\nu\in I}\operatorname{CZ}_{E_\nu,F_\nu}^{\Gamma_\nu}(X_\nu,U_\nu),
  \qquad
  \mathbb D^{\Gamma;I}_U
  :=\bigotimes_{\nu\in I}\mathbb D^{\Gamma_\nu}_{U_\nu},
  \qquad
  (\Phi_{\vec\eps}^X)_I:=\bigotimes_{\nu\in I}\Phi^{X_\nu}_{\eps_\nu}.
\end{equation*}
\end{definition}

\begin{remark}\label{rem:CZKk=1}
When $k=1$, there are only two subsets $I=\varnothing$ and $I=\{1\}$ of $[1]=\{1\}$, and we immediately recognise
\begin{equation*}
  \Norm{K}{\operatorname{CZK}_{\varnothing,\{1\}}(E,F,W)}
  =\Norm{K}{\operatorname{WBP}(W)},\qquad
  \Norm{K}{\operatorname{CZK}_{\{1\},\varnothing}(E,F,W)}
  =\Norm{K}{\operatorname{CZK}(E,F)},
\end{equation*}
and hence $\Norm{K}{\operatorname{CZK}(E,F,W)}$ agrees with the earlier Definition \ref{def:CZEFGH} in the one-parameter case.
\end{remark}

\begin{definition}\label{KE}
Let $K\in\mathscr S'(\mathbb R^{\vn}\times\mathbb R^{\vn})$.
For $I\subset[k]$, and $\psi_I,\varphi_I\in\mathscr S(\R^{\vn_I})$, we define a new distribution
\begin{equation*}
  \pair{K}{\psi_I\otimes\varphi_I}_I \in \mathscr S'(\R^{\vn_{[k]\setminus I}}\times \R^{\vn_{[k]\setminus I}})
\end{equation*}
by
\begin{equation*}
  \Bpair{  \pair{K}{\psi_I\otimes\varphi_I}_I }{\psi_{[k]\setminus I}\otimes\varphi_{[k]\setminus I}}
  :=  \pair{K}{(\psi_I\otimes\psi_{[k]\setminus I})\otimes(\varphi_I\otimes \varphi_{[k]\setminus I})}
\end{equation*}
for all $\psi_{[k]\setminus I},\varphi_{[k]\setminus I}\in\mathscr S(\R^{\vn_{[k]\setminus I}})$.
\end{definition}

We record the following simple observation:

\begin{lemma}\label{lem:CZKrest}
Let $\psi_\nu,\varphi_\nu\in C_{\rm c}^\infty(\R^{n_\nu})$, let $\psi_I:=\bigotimes_{\nu\in I}\psi_\nu$ and $\varphi_I$ be defined similarly.
Assume that $K\in \mathrm{CZK}(\vec E,\vec F,\vec W)$. If $[k]=I\cup J$ with $I\cap J=\varnothing$, then
\begin{equation*}
  \Norm{\pair{K}{\psi_I\otimes\varphi_I}_I}{\operatorname{CZK}(\vE_J,\vF_J,\vec W_J)}
  \leq\Norm{K}{\operatorname{CZK}(\vE,\vF,\vec W)}\prod_{\nu\in I}B_{W_\nu}(\psi_\nu,\varphi_\nu).
\end{equation*}
\end{lemma}

\begin{proof}
Let $r_\nu:=B_{W_\nu}(\psi_\nu,\varphi_\nu)$. If any of these numbers vanishes, then one of $\psi_\nu$ or $\varphi_\nu$ vanishes, in which case the left-hand side is zero, and there is nothing to prove. Otherwise, let $\psi_\nu':=r_\nu^{-1}\psi_\nu$. Then $\Psi_\nu:=\psi_\nu'\otimes\varphi_\nu\in\mathscr B(W_\nu)$ and $\Psi_I:=\bigotimes_{\nu\in I}\Psi_\nu\in \mathscr B(\vec W)_I$. Hence
\begin{equation*}
\begin{split}
  &\Norm{\pair{K}{\psi_I\otimes\varphi_I}_I}{\operatorname{CZK}(\vE_J,\vF_J,\vec W_J)} \\
  &=\sup_{\genfrac{}{}{0pt}{}{{J=G\cup H}}{G\cap H=\varnothing}}
  \sup_{\genfrac{}{}{0pt}{}{\genfrac{}{}{0pt}{}{\Gamma_G\in\mathscr G(\vE,\vF)_G}{(X,U,\vec\eps)\in\R^{4\vn_G+\abs{G}}_{\Gamma}}}{
    \Xi_{H}\in \mathscr B(\vec W)_H}}
  \operatorname{CZ}_{\vE,\vF}^{\Gamma;G}(X,U)
  \abs{\pair{\mathbb D^{\Gamma;G}_U
  \pair{K}{\psi_I\otimes\varphi_I}_I  }{(\Phi_{\vec\eps}^X)_G\otimes \Xi_H}},
\end{split}
\end{equation*}
where
\begin{equation*}
  \pair{\mathbb D^{\Gamma;G}_U
  \pair{K}{\psi_I\otimes\varphi_I}_I  }{(\Phi_{\vec\eps}^X)_G\otimes \Xi_H}
  =\pair{\mathbb D^{\Gamma;G}_U K}{ (\Phi_{\vec\eps}^X)_G\otimes \Psi_I \otimes \Xi_H} \prod_{\nu\in I}r_\nu.
\end{equation*}
Here $ \Psi_I \otimes \Xi_H\in\mathscr B(\vec W)_{I\cup H}$ and we have $G\cup(I\cup H)=I\cup J=[k]$ and $G\cap(I\cup H)=\varnothing$. Hence the previous line multiplied by $\operatorname{CZ}_{\vE,\vF}^{\Gamma;G}(X,U)$ is bounded by
$\Norm{K}{\operatorname{CZK}(\vE,\vF,\vec W)}\prod_{\nu\in I}r_\nu$, which is the right-hand side of the claim. This completes the proof.
\end{proof}

In particular, taking $I=[k]\setminus\{\nu\}$, we end up with one-parameter operators $\pair{K}{\psi_I\otimes\varphi_I}_I\in\operatorname{CZK}(E_\nu,F_\nu)$, for which we have Definition \ref{def:T1=0} about the annihilation of monomials. This immediately leads to the following extensions:

\begin{definition}\label{def:T1=0k}
For $K\in\operatorname{CZK}(\vE,\vF,\vec W)$, all $\nu\in[k]$ and $\abs{\mu_\nu}\leq\Floor{E_\nu}$, we say that
\begin{equation*}
  \pair{K}{\cdot\otimes y_\nu^{\mu_\nu}}:=0\quad\text{if}\quad
  \pair{\pair{K}{\psi_{[k]\setminus\{\nu\}}\otimes\varphi_{[k]\setminus\{\nu\}}}_{[k]
  \setminus\{\nu\}}}{\cdot\otimes y^{\mu_\nu}}=0
\end{equation*}
in the sense of Definition \ref{def:T1=0}, for all $\psi_{\iota},\varphi_{\iota}\in C_{\rm c}^\infty(\R^{n_\iota})$, where $\psi_{[k]\setminus\{\nu\}}:=\bigotimes_{\iota\in[k]\setminus\{\nu\}}\psi_\iota$ and $\varphi_{[k]\setminus\{\nu\}}$ is defined similarly.

The condition $\pair{K}{x_\nu^{\lambda_\nu} \otimes\cdot }=0$
for all $\abs{\lambda_\nu}\leq\Floor{F_\nu}$ is defined analogously.
\end{definition}

\begin{definition}\label{def:CZEFGHk}
Given $K\in\mathscr S'(\R^{\vn}\times\R^{\vn})$, we say that $K\in\operatorname{CZK}(\vE,\vF,\vec G,\vec H,\vec W)$ if
\begin{enumerate}[\rm(i)]
  \item $K\in\operatorname{CZK}(\vE,\vF,\vec W)$,
  \item\label{it:T1=0k} $\pair{K}{\cdot\otimes y_\nu^{\mu_\nu}}=0$ for all $\nu\in[k]$ and all $\abs{\mu_\nu}\leq G_\mu\leq\Floor{E_\mu}$,
  \item\label{it:T*1=0k} $\pair{K}{x_\nu^{\lambda_\nu}\otimes\cdot}=0$ for all $\nu\in[k]$ and all $\abs{\lambda_\nu}\leq H_\mu\leq\Floor{F_\mu}$.
\end{enumerate}
We define
\begin{equation*}
  \Norm{K}{\operatorname{CZK}(\vE,\vF,\vec G,\vec H,\vec W)}:=
  \begin{cases}   \Norm{K}{\operatorname{CZK}(\vE,\vF,\vec W)}, & \text{if \eqref{it:T1=0k} and \eqref{it:T*1=0k} hold}, \\
  \infty, & \text{otherwise}. \end{cases}
\end{equation*}
\end{definition}

\begin{lemma}\label{lem:T1=0rest}
Let $\psi_\nu,\varphi_\nu\in C_{\rm c}^\infty(\R^{n_\nu})$,
and let $\psi_I:=\bigotimes_{\nu\in I}\psi_\nu$ and $\varphi_I$ be defined similarly. Let $[k]=I\cup J$ with $I\cap J=\varnothing$. If $K\in \operatorname{CZK}(\vE,\vF,\vec G,\vec H,\vec W)$, then
\begin{equation*}
  \pair{K}{\psi_I\otimes\varphi_I}_I\in\operatorname{CZK}(\vE_J,\vF_J,\vec G_J,\vec H_J,\vec W_J).
\end{equation*}
\end{lemma}

\begin{proof}
Since
$$K\in \operatorname{CZK}(\vE,\vF,\vec G,\vec H,\vec W)
\subset\operatorname{CZK}(\vE,\vF,\vec W),$$
we know from Lemma \ref{lem:CZKrest} that
$\pair{K}{\psi_I\otimes\varphi_I}_I\in \operatorname{CZK}(\vE_J,\vF_J,\vec W_J)$.
Hence it remains to verify the cancellation conditions contained in $\operatorname{CZK}(\vE_J,\vF_J,\vec G_J,\vec H_J,\vec W_J)$.

Let $\nu\in J$ and $\abs{\mu_\nu}\leq G_\mu$. Since
$K\in\operatorname{CZK}(\vE,\vF,\vec G,\vec H,\vec W)$, we have
$\pair{K}{\cdot\otimes y_\nu^{\mu_\nu}}=0$, which means that, for all $\psi_{\iota},\varphi_{\iota}\in C_{\rm c}^\infty(\R^{n_\nu})$ and $\psi_I,\varphi_I$ as defined in the statement of the lemma, we have
\begin{equation*}
\begin{split}
  0 &=\Bpair{ \pair{K}{\psi_{[k]\setminus\{\nu\}}\otimes\varphi_{[k]\setminus\{\nu\}}}_{[k]\setminus\{\nu\}} }{   \cdot \otimes y^{\mu_{\nu}}} \\
   &=\Bpair{ \Bpair{ \pair{K}{\psi_I\otimes\varphi_I}_I }{\psi_{J\setminus\{\nu\}}\otimes\varphi_{J\setminus\{\nu\}}}_{J\setminus\{\nu\}}
 }{   \cdot \otimes y^{\mu_{\nu}}},
\end{split}
\end{equation*}
and vanishing of the right-hand side, for all $\psi_\iota,\varphi_\iota\in C_{\rm c}^\infty(\R^{n_\iota})$ with $\iota\in J\setminus\{\nu\}$, is the definition of
$\pair{  \pair{K}{\psi_I\otimes\varphi_I}_I }{ \cdot \otimes y^{\mu_{\nu}} }=0$.

The proof of $\pair{  \pair{K}{\psi_I\otimes\varphi_I}_I }{ x^{\lambda_{\nu}} \otimes \cdot }=0$, for all $\nu\in J$ and $\abs{\lambda_\mu}\leq H_\mu$, is similar. This completes the proof.
\end{proof}

\begin{lemma}\label{lem:6.19k}
Let $\vE,\vF,\vec G,\vec H,\vec L,\vec L',\vec M,\vec N,\vec N',\vec W$ be vectors, whose respective components satisfy the assumptions \eqref{eq:6.19ass} of Proposition \ref{prop:6.19},
and let $K\in\operatorname{CZK}(\vE,\vF,\vec G,\vec H,\vec W)$.
For each $\nu\in[k]$, let $\mathscr A(N_\nu)\subset\N^{n_\nu}\times\{0,1\}$
be as in \eqref{A(N)}, and let
\begin{equation*}
   \mathscr A(\vec N)
   :=\Big\{(\alpha,\kappa)\in\N^{\vn}\times\{0,1\}^k:
      (\alpha_\nu,\kappa_\nu)\in\mathscr A(N_{\nu})\text{ for every }\nu\in[k]\Big\}
\end{equation*}
For $\varphi=\bigotimes_{\nu=1}^k\varphi_\nu$ with $\varphi_\nu\in C_{\rm c}^\infty(3\cdot [0,1)^{n_\nu})$, and $\psi=\bigotimes_{\nu=1}^k\psi_\nu$ with $\psi_\nu\in C_{\rm c}^\infty(\R^{n_\nu})$, we then have
\begin{equation}\label{eq:6.19k}
  \Big|\Big\langle(\vec1+\abs{x}_{\mathrm{vec}}^2)^{\frac{\vec M}{2}}\Delta_{(h,0)}^{(\kappa,0)}\partial_x^\alpha K,{\psi\otimes \varphi}\Big\rangle\Big|
  \lesssim\Norm{\psi}{L^1}\Norm{K}{\operatorname{CZK}(\vE,\vF,\vec G,\vec H,\vec W)}
  \Norm{\varphi}{\operatorname{Atom}(\vec L',\vec N')}
  \abs{h}_{\rm{vec}}^{\vec N^{**}\kappa}
\end{equation}
for all $(\alpha,\kappa)\in\mathscr A(\vec N)$ and all $h\in\R^{\vn}$
satisfying $\abs{h}_{\rm{vec}}\leq\vone$,
where the implicit positive constant is independent of $\varphi$, $\psi$, and $K$, and where $|x|_{\mathrm{vec}}^2:=(|x_1|^2,\ldots,|x_k|^2)$.
\end{lemma}

\begin{proof}
We use induction on $k$. The base case $k=1$ is contained in Corollary \ref{cor:6.19}. We assume for induction that the $\tilde k$-parameter case has been proved for all $\tilde k\in[k-1]$, and we proceed to prove it for $k$. For any $I\subset[k]$, we denote $\varphi_I:=\bigotimes_{\nu\in I}\varphi_\nu$, and similarly for $\psi_I$.

Let $[k]=I\cup J$ with $I\cap J=\varnothing$ be some non-trivial splitting; hence both $\abs{I},\abs{J}\in[k-1]$.
Then
\begin{equation*}
\begin{split}
  \Bpair{(\vec 1+\abs{x}_{\mathrm{vec}}^2)^{\frac{\vec M}{2}}\Delta_{(h,0)}^{(\kappa,0)}\partial_x^\alpha K}{\psi\otimes \varphi}
  = \Bpair{(\vec 1+\abs{x_I}_{\mathrm{vec}}^2)^{\frac{\vec M_I}{2}}\Delta_{(h_I,0)}^{(\kappa_I,0)}\partial_x^{\alpha_I}
    \tilde K_J
     }{\psi_I\otimes\varphi_I},
\end{split}
\end{equation*}
where
\begin{equation}\label{eq:tildeKJ}
\begin{split}
  \tilde K_J
  &:=
   \Bpair{ (\vec 1+\abs{x_J}_{\mathrm{vec}}^2)^{\frac{\vec M_J}{2}}\Delta_{(h_J,0)}^{(\kappa_J,0)}\partial_x^{\alpha_J} K}{\psi_J\otimes \varphi_J}_J  \\
   &\phantom{:}=
   \Bpair{K}{
   ((-\partial_x)^{\alpha_J}   \Delta_{-h_J}^{\kappa_J}[(\vec 1+\abs{x_J}_{\mathrm{vec}}^2)^{\frac{\vec M_J}{2}}\psi_J])\otimes \varphi_J}_J,
\end{split}
\end{equation}
which belongs to $\operatorname{CZK}(\vE_I,\vF_I,\vec G_I,\vec H_I,\vec W_I)$ by Lemma \ref{lem:T1=0rest}. Hence, applying the induction hypothesis to $\tilde K_J$,
\begin{equation*}
\begin{split}
  \Babs{\Bpair{(\vec 1+\abs{x}_{\mathrm{vec}}^2)^{\frac{\vec M}{2}}\Delta_{(h,0)}^{(\kappa,0)}\partial_x^\alpha K}{\psi\otimes \varphi}}
  \lesssim\Norm{\psi_I}{L^1}
   \Norm{\tilde K_J
    }{\operatorname{CZK}(\vE_I,\vF_I,\vec G_I,\vec H_I,\vec W_I)}
    \Norm{\varphi_I}{\operatorname{Atom}(\vec L_I',\vec N_I')}\abs{h_I}_{\rm{vec}}^{\vec N_I^{**}\kappa_I},
\end{split}
\end{equation*}
where
\begin{equation}\label{eq:normKJ}
\begin{split}
   \Norm{\tilde K_J}{\operatorname{CZK}(\vE_I,\vF_I,\vec G_I\vec H_I,\vec W_I)}
    &=\Norm{\tilde K_J}{\operatorname{CZK}(\vE_I,\vF_I,\vec W_I)}     \\
    &=\sup_{\genfrac{}{}{0pt}{}{I=A\cup B}{A\cap B=\varnothing}}
    \sup_{\genfrac{}{}{0pt}{}{\genfrac{}{}{0pt}{}{\Gamma_A\in\mathscr G(\vE,\vF)_A}{(X,U,\vec\eps)\in\R^{4\vn_A+\abs{A}}_{\Gamma}}}{
    \Xi_{B}\in \mathscr B(\vec W)_B}}
  \operatorname{CZ}_{\vE,\vF}^{\Gamma;A}(X,U)
  \Babs{\Bpair{\mathbb D^{\Gamma;A}_U
  \tilde K_J
  }{(\Phi_{\vec\eps}^X)_A\otimes \Xi_B}}.
\end{split}
\end{equation}
Substituting the formula for $\tilde K_J$ from \eqref{eq:tildeKJ} and rearranging, we see that
\begin{equation*}
  \Bpair{\mathbb D^{\Gamma;A}_U
  \tilde K_J
  }{(\Phi_{\vec\eps}^X)_A\otimes \Xi_B}
  =\Bpair{ (\vec 1+\abs{x_J}_{\mathrm{vec}}^2)^{\frac{\vec M_J}{2}}\Delta_{(h_J,0)}^{(\kappa_J,0)}\partial_x^{\alpha_J}
    \Bpair{\mathbb D_U^{\Gamma;A} K}{ (\Phi_{\vec\eps}^X)_A\otimes\Xi_B}_I
    }{\psi_J\otimes \varphi_J}.
\end{equation*}

Next, we plan to apply the induction hypothesis to
\begin{equation}\label{eq:tildeKI}
  \tilde{\tilde K}_I
  :=\Bpair{\mathbb D_U^{\Gamma;A} K}{ (\Phi_{\vec\eps}^X)_A\otimes\Xi_B}_I
  =\pm\Bpair{ K}{ (\mathbb D_{-U}^{\Gamma;A}(\Phi_{\vec\eps}^X)_A)\otimes\Xi_B}_I,
\end{equation}
where the (unimportant) sign depends on the total order of partial derivatives present in $\mathbb D_U^{\Gamma;A}$. From the expression on the right, we see that $\tilde{\tilde K}_I\in \operatorname{CZK}(\vE_J,\vF_J,\vec G_J,\vec H_J,\vec W_J)$ by Lemma \ref{lem:T1=0rest}. Thus, using the induction hypothesis, it follows that
\begin{equation*}
   \Babs{\Bpair{\mathbb D^{\Gamma;A}_U
  \tilde K_J
  }{(\Phi_{\vec\eps}^X)_A\otimes \Xi_B}}
  \lesssim
  \Norm{\psi_J}{L^1}\Norm{\tilde{\tilde K}_I}{\operatorname{CZK}(\vE_J,\vF_J,\vec G_J,\vec H_J,\vec W_J)}
  \Norm{\varphi_J}{\operatorname{Atom}(\vec L_J',\vec N_J')}\abs{h_J}_{\rm{vec}}^{\vec N_J^{**}\kappa_J},
\end{equation*}
where
\begin{equation}\label{eq:normKI}
\begin{split}
     \bNorm{\tilde{\tilde K}_I }{\operatorname{CZK}(\vE_J,\vF_J,\vec G_J,\vec H_J,\vec W_J)}
    &=\bNorm{\tilde{\tilde K}_I }{\operatorname{CZK}(\vE_J,\vF_J,\vec W_J)} \\
    &=\sup_{\genfrac{}{}{0pt}{}{J=C\cup D}{C\cap D=\varnothing}}
    \sup_{\genfrac{}{}{0pt}{}{\genfrac{}{}{0pt}{}{\Gamma_C\in\mathscr G(\vE,\vF)_C}{(X,U,\vec\eps)\in\R^{4\vn_C+\abs{C}}_{\Gamma}}}{
    \Xi_{D}\in \mathscr B(\vec W)_D}}
  \operatorname{CZ}_{\vE,\vF}^{\Gamma;C}(X,U)
  \Babs{\Bpair{\mathbb D^{\Gamma;C}_U
  \tilde{\tilde K}_I
  }{(\Phi_{\vec\eps}^X)_C\otimes \Xi_D}}.
\end{split}
\end{equation}
Substituting the formula for $\tilde{\tilde K}_I$ from \eqref{eq:tildeKI} and rearranging, we find that
\begin{equation*}
  \Bpair{\mathbb D^{\Gamma;C}_U   \tilde{\tilde K}_I  }{(\Phi_{\vec\eps}^X)_C\otimes \Xi_D}
  =\Bpair{\mathbb D^{\Gamma;A\cup C}_U K}{(\Phi_{\vec\eps}^X)_{A\cup C}\otimes\Xi_{B\cup D}}.
\end{equation*}
Multiplying this by the factors present in \eqref{eq:normKJ} and \eqref{eq:normKI}, namely
\begin{equation*}
  \operatorname{CZ}^{\Gamma;A}_{\vE,\vF}(X,U)\times
  \operatorname{CZ}^{\Gamma;C}_{\vE,\vF}(X,U)
  =\operatorname{CZ}^{\Gamma;A\cup C}_{\vE,\vF}(X,U),
\end{equation*}
and combining the suprema in \eqref{eq:normKJ} and \eqref{eq:normKI}, we arrive at exactly the quantity in the definition of $\Norm{K}{\operatorname{CZK}(\vE,\vF,\vec W)}$, just parameterised
with the new complementary sets $A\cup C$ and $B\cup D$ in place of the original $I$ and $J$, respectively.

Hence, combining all the estimates, we have obtained
\begin{equation*}
\begin{split}
  &\Babs{\Bpair{(\vec 1+\abs{x}_{\mathrm{vec}}^2)^{\frac{\vec M}{2}}\Delta_{(h,0)}^{(\kappa,0)}\partial_x^\alpha K}{\psi\otimes \varphi} } \\
  &\quad\lesssim\Norm{K}{\operatorname{CZK}(\vE,\vF,\vec G,\vec H,\vec W)}
  \Norm{\psi_J}{L^1}\Norm{\psi_I}{L^1}
  \Norm{\varphi_J}{ \operatorname{Atom}(\vec L_J',\vec N_J')}
  \Norm{\varphi_I}{\operatorname{Atom}(\vec L_I',\vec N_I')}
  \abs{h}_{\rm{vec}}^{\vec N_J^{**}\kappa}
  \abs{h}_{\rm{vec}}^{\vec N_I^{**}\kappa} \\
  &\quad=\Norm{K}{\operatorname{CZK}(\vE,\vF,\vec G,\vec H,\vec W)}
  \Norm{\psi}{L^1}
  \Norm{\varphi}{ \operatorname{Atom}(\vec L',\vec N')}
  \abs{h}_{\rm{vec}}^{\vec N^{**}\kappa},
\end{split}
\end{equation*}
and this completes the proof of Lemma \ref{lem:6.19k}.
\end{proof}

The following corollary of Lemma \ref{lem:6.19k} will be more convenient to apply below:

\begin{corollary}\label{cor:6.19k}
Let $\vE,\vF,\vec G,\vec H,\vec L,\vec L',\vec M,\vec N,\vec N',\vec W$ be vectors, whose respective components satisfy the assumptions \eqref{eq:6.19ass} of Proposition \ref{prop:6.19}.
Let $T\in\mathcal L(\mathscr S(\R^{\vn}),\mathscr S'(\R^{\vn}))$ have Schwartz kernel $K\in\operatorname{CZK}(\vE,\vF,\vec G,\vec H,\vec W)$.
Let $\varphi$ be an $(\vec L',\vec N')$-atom on the unit cube $[0,1)^{\vn}$.
Then, for all $\alpha\in\N^{\vn}$ satisfying
$\abs{\alpha}_{\rm{vec}}\leq\Floor{\vec N}$, the distributional derivatives
$\partial^\alpha T\varphi$ are given by $L^\infty$ functions that satisfy, for almost every $x\in\R^{\vn}$,
\begin{equation}\label{eq:Tphi<}
  \abs{\Delta_{h}^{\kappa}\partial^\alpha T\varphi (x)}
  \lesssim \abs{h}_{\rm{vec}}^{\vec N^{**}\kappa}(\vec 1+\abs{x}_{\mathrm{vec}}^2)^{-\frac{\vec M}{2}}
 \end{equation}
for all $(\alpha,\kappa)\in\mathscr A(\vec N)$ and all $h\in\R^{\vn}$
satisfying $\abs{h}_{\rm{vec}}\leq\vone$. Moreover, for all $\nu\in[k]$, all $\abs{\gamma_\nu}\leq L_\nu$, and almost all $x'\in\R^{\vn}\ominus\R^{n_\nu}$, we have
\begin{equation}\label{eq:intTphi=0}
   \int_{\R^{n_\nu}}x_{\nu}^{\gamma_\nu}T\varphi(x_\nu,x')dx_\nu = 0.
\end{equation}
\end{corollary}

\begin{proof}
We note that the conclusion \eqref{eq:Tphi<} becomes stronger with increasing $\vec M$, while the only constraint on $\vec M$ is $\vec M\leq\vec F+\vn$ from \eqref{eq:6.19ass}. Hence, without loss of generality, we may take
\begin{equation}\label{eq:M=F+n}
  \vec M=\vec F+\vn.
\end{equation}

Let $\psi_\nu\in C_{\rm}^\infty(\R^{n_\nu})$ and $\psi:=\bigotimes_{\nu=1}^k\psi_\nu$. Then $\psi$ and $\varphi$ satisfy the assumptions of Lemma \ref{lem:6.19k}, and hence also the conclusion \eqref{eq:6.19k}. But
\begin{equation*}
\begin{split}
  \mathrm{LHS}\eqref{eq:6.19k}
 &= \Babs{  \Bpair{K}{  (\partial^\alpha\Delta_{h}^{\kappa}[(\vec 1+\abs{x}_{\mathrm{vec}}^2)^{\frac{\vec M}{2}}\psi])\otimes\varphi} } \\
 &= \Babs{  \Bpair{T\varphi }{  \partial^\alpha\Delta_{h}^{\kappa}[(\vec 1+\abs{x}_{\mathrm{vec}}^2)^{\frac{\vec M}{2}}\psi]} } \\
 &= \Babs{  \Bpair{  (\vec 1+\abs{x}_{\mathrm{vec}}^2)^{\frac{\vec M}{2}}\Delta_{h}^{\kappa}\partial^\alpha(T\varphi) }{ \psi} }.
\end{split}
\end{equation*}
With usual bump functions $0\leq \phi_\nu\in C_{\rm c}^\infty(\R^{\vn})$ of integral one,
consider $$\psi_\nu=\eps_\nu^{-n_\nu}\phi_\nu(\eps_\nu^{-1}(\cdot-x_\nu)).$$
Noting that $\Norm{\psi_\nu}{L^1}=\Norm{\phi_\nu}{L^1}=1$, we see that
\begin{equation*}
  \mathrm{RHS}\eqref{eq:6.19k}\lesssim\abs{h}_{\rm{vec}}^{\vec N^{**}\kappa}
\end{equation*}
for all these $\psi$, also using the assumptions on $K$ and $\varphi$. Considering the limit $\eps_\nu\to 0$, we infer that
\begin{equation*}
  \BNorm{
   (\vec 1+\abs{x}_{\mathrm{vec}}^2)^{\frac{\vec M}{2}}\Delta_{h}^{\kappa}\partial^\alpha(T\varphi)  }{L^\infty}
  \lesssim\abs{h}_{\rm{vec}}^{\vec N^{**}\kappa},
\end{equation*}
which is the same as \eqref{eq:Tphi<}. With $(\alpha,\kappa)=(\mathbf 0,\vec 0)$, we have in particular $\abs{T\varphi}\lesssim(\vec 1+\abs{x}_{\mathrm{vec}}^2)^{-\frac{\vec M}{2}}$. From \eqref{eq:M=F+n} and $\vF>\vec L$, which is part of \eqref{eq:6.19ass}, it then follows that, for all $\abs{\gamma_\nu}\leq L_\nu$,
\begin{equation}\label{eq:TphiL1}
   x_\nu^{\gamma_{\nu}}T\varphi\in L^1(\R^{\vn}),
\end{equation}
and thus $[x_\nu\mapsto x_\nu^{\gamma_{\nu}}T\varphi_{[k]}(x_\nu,x')]\in L^1(\R^{n_\nu})$ for a.e.\ $x'\in\R^{\vn}\ominus\R^{n_\nu}$ by Fubini's theorem. Thus, the claimed identity \eqref{eq:intTphi=0} at least makes sense.

We next use the condition corresponding to the parameter $\vec H$ in $K\in\operatorname{CZK}(\vE,\vF,\vec G,\vec H,\vec W)$ and the cancellation of the $(\vec L',\vec N')$-atom $\varphi$. By \eqref{eq:6.19ass}, we know that $\vec H\geq\vec L$ and also that $\vec L'\geq\Floor{\vec F}\geq\vec L$ (using $\vec F>\vec L\in\Z^k$ for the last inequality). By the $\operatorname{CZK}(\vE,\vF,\vec G,\vec H,\vec W)$ condition, for a  uniform approximation to the identity  $(\chi_j)_{j=0}^\infty$ (Definition \ref{def:UAI}) and arbitrary $\psi_\iota\in C_{\rm c}^\infty(\R^{n_\iota})$ with $\psi_{[k]\setminus\{\nu\}}:=\bigotimes_{\iota\in[k]\setminus\{\nu\}}\psi_\iota$, we have
\begin{equation}\label{eq:Kx}
  \pair{K}{(\psi_{[k]\setminus\{\nu\}}\otimes\chi_j x_\nu^{\gamma_\nu})\otimes\varphi}
  =\pair{K}{(\psi_{[k]\setminus\{\nu\}}\otimes\chi_j x_\nu^{\gamma_\nu})\otimes(\varphi_{[k]\setminus\{\nu\}}\otimes\varphi_{\nu})}
  \underset{j\to\infty}{\longrightarrow} 0,
\end{equation}
for all $\abs{\gamma_\nu}\leq H_\nu$, provided that $\varphi_{\nu}\in\mathcal D_{\abs{\gamma_\nu}}$. In particular, this holds for all $\abs{\gamma_\nu}\leq L_\nu\leq L_{\nu}'\wedge H_\nu$, as we just observed. On the other hand, we have
\begin{equation*}
\begin{split}
  \mathrm{LHS}\eqref{eq:Kx} &=\pair{K}{(\psi_{[k]\setminus\{\nu\}}\otimes\chi_j x_\nu^{\gamma_\nu})\otimes\varphi} \\
  &=\pair{T\varphi}{\psi_{[k]\setminus\{\nu\}}\otimes\chi_j x_\nu^{\gamma_\nu}} \\
  &=\int_{\R^{\vn}\ominus\R^{n_\nu}} \psi_{[k]\setminus\{\nu\}}(x')\int_{\R^{n_\nu}}(T\varphi)(x) \chi_j(x') x_\nu^{\gamma_\nu}dx_{\nu}dx'.
\end{split}
\end{equation*}
By \eqref{eq:TphiL1}, the fact that $\psi_{[k]\setminus\{\nu\}}\in L^\infty$, and the pointwise limit $\chi_j(x')\to 1$, dominated convergence shows that
\begin{equation*}
  0=\lim_{j\to\infty}\mathrm{LHS}\eqref{eq:Kx}
  =\int_{\R^{\vn}\ominus\R^{n_\nu}} \psi_{[k]\setminus\{\nu\}}(x')\int_{\R^{n_\nu}}(T\varphi)(x)  x_\nu^{\gamma_\nu}dx_{\nu}dx'.
\end{equation*}
Since this holds for $\psi_{[k]\setminus\{\nu\}}=\bigotimes_{\iota\in[k]\setminus\{\nu\}} \psi_{\iota}$ with arbitrary $\psi_{\iota}\in C_{\rm c}^\infty(\R^{n_\iota})$, we must have
\begin{equation*}
  \int_{\R^{n_\nu}}(T\varphi)(x) x_\nu^{\gamma_\nu}dx_{\nu}=0
\end{equation*}
for almost every $x'\in\R^{\vn}\ominus\R^{n_\nu}$. This is exactly \eqref{eq:intTphi=0} that we wanted to prove. Since we already proved \eqref{eq:Tphi<}, this completes the proof of Corollary \ref{cor:6.19k}.
\end{proof}

\begin{corollary}\label{cor:6.19kk}
Let $\vE,\vF,\vec G,\vec H,\vec L,\vec L',\vec M,\vec N,\vec N',\vec W$ be vectors, whose respective components satisfy the assumptions \eqref{eq:6.19ass} of Proposition \ref{prop:6.19}.
Let $T\in\mathcal L(\mathscr S(\R^{\vn}),\mathscr S'(\R^{\vn}))$ have Schwartz kernel $K\in\operatorname{CZK}(\vE,\vF,\vec G,\vec H,\vec W)$.
Then $T$ maps $(\vec L',\vec N')$-atoms to a
uniform constant multiple of $(\vec M,\vec L,\vec M,\vec N)$-molecules.
\end{corollary}

\begin{proof}
Corollary \ref{cor:6.19k} states exactly this in the special case of atoms on the unit cube $[0,1)^{\vn}$. Due to the behaviour with respect to translations and dilations of atoms (Definition \ref{def:atom}), molecules (Definition \ref{moleKLMN}), and the $\operatorname{CZK}(\vec E,\vec F,\vec G,\vec H,\vec W)$ conditions alike, the general case follows from this by changes of variables.
\end{proof}

We now come to the main theorem of this section on the boundedness of Calder\'on--Zygmund operators on matrix-weighted multi-parameter Besov--Triebel--Lizorkin-type spaces.

\begin{theorem}[The $T(1)$ theorem for $\Atau(V)$]\label{T1 BF}
Let $\vp\in(0,\infty)^k$, $\vq\in(0,\infty]^k$, $\vs\in\R^k$, $\tau\in[0,\infty)$,
$r:=r(\pi,\vec p,\vec q)\wedge 1$,
$\pi\in S_{[2k]}$ be admissible for $(\vp,\vq)$,
and $V$ be a matrix weight such that $(\vp,\vq,\pi,V)$ belongs
to one of the main cases (Definition \ref{main cases}).
Let $T\in\mathcal L(\mathscr S(\R^{\vn}),\mathscr S'(\R^{\vn}))$ have Schwartz kernel
$K\in\operatorname{CZK}(\vE,\vF,\vec G,\vec H,\vec W)$ for some $\vec W\in\N^k$, where
$\vec E,\vec F\in(0,\infty)^k$ satisfy
\begin{equation}\label{eq:T1 EF}
  \vE>(\tau\vn+\vs)_+,\quad
  \vF>\Big[\Big(\tfrac1r-1\Big)\vn-\vs\Big]_+\vee\Big(\tfrac1r+\tau-1\Big)\vn
\end{equation}
and $\vec G,\vec H\in\Z_{\geq-1}^k$ satisfy
\begin{equation}\label{eq:T1 GH}
  \vec G\geq\floor{\tau\vn+\vs}_+,\qquad
  \vec H\geq \floor{(\tfrac{1}{r}-1)\vn-\vs}.
\end{equation}
Then there exists an operator $\widetilde T\in\mathcal L(\Atau(V))$
that agrees with $T$ on $\mathscr{S}_0(\mathbb R^{\vn})$.
\end{theorem}

\begin{proof}
By Proposition \ref{ext}, the conclusion will follow if we can find
$(\vec L_a,\vec N_a)\in(0,\infty)^{2k}$ such that $T$ maps $(\vec L_a, \vec N_a)$-atoms
to $\Atau$-synthesis molecules. Recall from Definition \ref{def:AtauMole} that $\Atau$-synthesis molecules are $(\vec K',\vec L',\vec M',\vec N')$-molecules, where these parameters satisfy \eqref{KMM} and \eqref{LNN}. By Corollary \ref{cor:6.19kk},
$T$ with kernel $K\in\operatorname{CZK}(\vE,\vF,\vec G,\vec H,\vec W)$
maps $(\vec L_a,\vec N_a)$-atoms to $(\vec K',\vec L',\vec M',\vec N')$-molecules
(with $\vec K'=\vec M'$) under the conditions \eqref{eq:6.19ass}
with $(\vec L',\vec M',\vec N')$ in place of $(\vec L,\vec M,\vec N)$
and $(\vec L_a,\vec N_a)$ in place of $(\vec L',\vec N')$, respectively.
Combining these results, in order to map certain atoms
to $\Atau$-synthesis molecules, we need the following:
\begin{equation*}
\begin{cases}\vE\geq \phantom{\lfloor}\vec N'\phantom{\rfloor}>\tau\vn+\vs, \\ \vE>\floor{\vec N'}_+,\end{cases}\quad
\begin{cases}\vF\geq \vec M'-\vn=\vec K'-\vn>[(\frac1r-1)\vn-\vs]_+\vee(\frac1r+\tau-1)\vn, \\
\vF>\vec L'>(\frac1r-1)\vn-\vec s-\vone,\end{cases}
\end{equation*}
and
\begin{equation*}
  \begin{cases}\vec G\geq \floor{\vec N'}_+, \quad \vec N'>\tau\vn+\vs, \\
  \vec H\geq\vec L'>(\frac1r-1)\vn-\vs-\vone.\end{cases}
\end{equation*}
From these, we readily verify that such $\vec M'\in[0,\infty)^k$, $\vec N'\in(0,\infty)^k$, and $\vec L'\in\Z_{\geq-1}^k$ can be found provided that $\vE,\vF,\vec G,\vec H$ satisfy the assumptions of the theorem. This completes the proof.
\end{proof}

As a very special case, we obtain a new proof of the following result of Domelevo et al.\ \cite{DKPS}:

\begin{corollary}[\cite{DKPS}, Theorem 1.1(1)]\label{T1Co}
Let $p\in(1,\infty)$, $\vec E\in(0,\infty)^k$, $\vec F\in(0,\infty)^k$, and $V\in\mathscr A_p(\R^{\vn})$.
Let $T\in\mathcal L(\mathscr S(\R^{\vn}),\mathscr S'(\R^{\vn}))$ have Schwartz kernel
$K\in\operatorname{CZK}(\vE,\vF,\vec 0,\vec 0,\vec W)$ for some $\vec W\in\N^k$.
Then there exists an operator $\widetilde T\in\mathcal L(L^p(V))$
that agrees with $T$ on $\mathscr{S}_0(\mathbb R^{\vn})$.
\end{corollary}

\begin{proof}
By Theorem  \ref{Lp=F}, we have
\begin{equation*}
  L^p(V)=\dot F^{\vec0,0}_{p,2}(V)=\Atau(V),
\end{equation*}
where $\vec s=\vec0$, $\tau=0$, $\vec p=p\cdot\vone$, $\vec q=2\cdot\vone$, and $\pi=F$ be the identity permutation. For these parameters, we $r(\pi,\vp,\vq)=\min(p,2)>1$ and then $r:=1\wedge r(\pi,\vp,\vq)=1$.
Thus, the conditions \eqref{eq:T1 EF} and \eqref{eq:T1 GH} of Theorem \ref{Lp=F} for this space
become $\vec E>\vec0$, $\vec F>\vec0$, $\vec G\geq\vec 0$,
and $\vec H\geq\vec 0$, which are precisely the assumptions of the corollary.
Theorem \ref{T1 BF} thus applies to prove the existence of a desired
\begin{equation*}
  \widetilde T\in\mathcal L(\Atau(V))=\mathcal L(L^p(V)),
\end{equation*}
and this completes the proof of Corollary \ref{T1Co}.
\end{proof}

\begin{remark}\label{r12.33}\mbox{}
\begin{enumerate}[\rm(i)]
\item Theorem \ref{T1 BF} is completely new in the multi-parameter case.
In the one-parameter case, it essentially reduces to \cite[Theorem 6.18]{BHYY:1c} (with a slightly different range of indices) for the usual Besov and Triebel--Lizorkin spaces with $\tau=0$, and slightly more generally, while in general it deals with a different scale of Besov--Triebel--Lizorkin-type space (see Theorem \ref{compare A}).

\item The paraproduct-free multi-parameter $T(1)$ theorem in $L^p(V)$
was first established in \cite{DKPS}.
As a special case when $k=2$, Corollary \ref{T1Co}
reproduces the first assertion of
\cite[Theorem 1.1(1)]{DKPS}.
Non-paraproduct-free situations, that is, assuming $T(1)\in \mathrm{BMO}$
instead of $T(1)=0$, were also studied in \cite{DKPS,Vuo24}.
Here, we only treat paraproduct-free cases, as in the majority of works
on $T(1)$ theorems in matrix-weighted Besov--Triebel--Lizorkin spaces.
Obtaining a general $T(1)$ theorem for these spaces is left for future work, and should probably first be addressed in the one-parameter case, since it seems to be open even there.

\item Using the identification of $\Atau(V)$ with other known function spaces, the model of Corollary \ref{T1Co} can be used to obtain similar corollaries in other concrete spaces of interest. As an example, we mention the matrix-weighted Sobolev spaces $\dot W^{\vec l}_p(V)=\dot F^{\vec l,0}_{p,2}(V)$ (see Theorem \ref{thm:Sobolev} for this identity), for which Theorem \ref{T1 BF} is already new. While it is conceivable that this special case could be handled without the full theory behind Theorem \ref{T1 BF}, we note that it probably cannot be easily deduced from the $L^p(V)$ results and methods of \cite{DKPS}. Namely, the approach of \cite{DKPS} is based the dyadic representation of \cite{Mar12} involving expansions in terms of Haar functions, which are ill-adapted to the analysis of higher order smoothness spaces, unlike the techniques of Theorem \ref{T1 BF}. This discussion also shows that our proof of Corollary \ref{T1Co} is genuinely different from that in \cite{DKPS}.
\end{enumerate}
\end{remark}

\section{Sobolev-type embeddings}\label{embedding}

Sobolev-type embedding theorems express a fundamental trade-off between smoothness and integrability of functions. Starting with the pioneering contributions of Sobolev \cite{Sobolev} on embeddings of the form $W^{m}_p(\R^n)\hookrightarrow W^{m-l}_q(\R^n)$ with $\frac1q=\frac1p-\frac{l}{n}$, many extensions to different settings are available by now. For Besov and Triebel--Lizorkin spaces in the unweighted one-parameter setting, a Sobolev-type embedding theorem in the full range of indiceswas first obtained by Jawerth \cite{Jawerth}, and another proof is due to Johnsen and Sickel \cite{JS07}. Its matrix-weighted one-parameter version is recently due to Bu et al. \cite{BYYZ1}, relying on the same interpolation inequality as \cite{JS07}.

In this section, we obtain a matrix-weighted multi-parameter extension of these results.
Subsection \ref{sec:embBTL} deals with Besov and Triebel--Lizorkin spaces, while Subsection \ref{sec:lift} derives some consequences in $L^p(V)$ and $\dot W^{\vec l}_p(V)$ spaces.

\subsection{Embeddings of Besov and Triebel--Lizorkin spaces}\label{sec:embBTL}
In this subsection, we characterise several embeddings between two different matrix-weighted spaces in terms of a relatively simply condition \eqref{eq-sobolev-B}; the main result is stated in Theorem \ref{thm-sobolev-B}.
Instead of the interpolation approach of \cite{BYYZ1,JS07} used for related results in the one-parameter case, our proof is closer to the original one of Sobolev \cite{Sobolev}, based on the boundedness of fractional integrals, which extends well to the multi-parameter setting. We will also explain the failure of a na\"ive extension of the method of \cite{BYYZ1,JS07} at the end of the subsection. The technical core of our approach is contained in the following:

\begin{proposition} \label{thm-sobolev-prop}
Let $\vec s_0,\vec s_1\in\R^k$.
Let $0<p_0<p_1<\infty$ and $\vq\in(0,\infty]^k$.
Assume that function sequences $\{f_{\vj}\}_{\vj\in\mathbb Z^k}$
and $\{g_{\vj}\}_{\vj\in\mathbb Z^k}$ satisfy,
for every $\vj \in\mathbb Z^k$ and $x\in\mathbb R^{\vec n}$,
\begin{equation}\label{eq:f<g}
2^{\vj\cdot (\vec s_1-\frac{\vec n}{p_1})} |f_{\vj}(x)|
\leq 2^{\vj\cdot (\vec s_0-\frac{\vec n}{p_0})} |g_{\vj}(x)|.
\end{equation}
If we further assume that $\{g_{\vj}\}_{\vj\in\mathbb Z^k}$
is constant on each $\mathscr D_{\vj}(\mathbb R^{\vec n})$,
then, for every $\vj \in\mathbb Z^k$,
\begin{equation}\label{lqLp}
\big\|2^{\vj\cdot \vec s_1}f_{\vj}\big\|_{L^{p_1}(\R^{\vn})}
\leq \big\|2^{\vj\cdot \vec s_0} g_{\vj}\big\|_{L^{p_0}(\R^{\vn})}
\end{equation}
and there exists a positive constant $C$,
depending only on $s_0$, $p_0$, $p_1$, and $\vq$, such that
\begin{equation}\label{Lplq}
\big\|\big\{2^{\vj\cdot \vec s_1} f_{\vj}\big\}_{\vj \in\mathbb Z^k}\big\|_{L^{p_1}\ell^{\vq}}
\leq C \big\|\big\{2^{\vj\cdot \vec s_0} g_{\vj}\big\}_{\vj \in\mathbb Z^k}\big\|_{L^{p_0}\ell^\infty}.
\end{equation}
\end{proposition}

\begin{proof}
We first prove \eqref{lqLp}.
For every $\vj \in\mathbb Z^k$,
\begin{align*}
\big\|2^{\vj\cdot \vec s_1} f_{\vj}\big\|_{L^{p_1}(\R^{\vn})}^{p_1}
&\leq \int_{\mathbb R^{\vec n}}
2^{(\vj\cdot \vec n)(1-\frac{p_1}{p_0})}
(2^{\vj\cdot \vec s_0}|g_{\vj}|)^{p_1}
= \sum_{Q\in \mathscr D_{\vj}(\mathbb R^{\vec n})}
|Q|^{\frac{p_1}{p_0}} (2^{\vj\cdot \vec s_0} |g_Q|)^{p_1} \\
&\leq \Bigg[\sum_{Q\in \mathscr D_{\vj}(\mathbb R^{\vec n})}
|Q| (2^{\vj\cdot \vec s_0}|g_Q|)^{p_0}\Bigg]^{\frac{p_1}{p_0}}
= \|2^{\vj\cdot \vec s_0} g_{\vj}\|_{L^{p_0}(\R^{\vn})}^{p_1},
\end{align*}
where $g_Q$ denotes the constant value of $g_{\vj}$
on $Q\in\mathscr D_{\vj}(\R^{\vn})$.
This finishes the proof of \eqref{lqLp}.

The main part of the proof consists of showing \eqref{Lplq}.
For every $\nu\in[k]$, the left-hand side of \eqref{Lplq} increases as $q_\nu$ decreases; hence it is enough to prove \eqref{Lplq} for smaller values of each $q_\nu$. Without loss of generality, we may hence assume that $q_\nu\equiv q<p_0$ for each $\nu\in[k]$.

The left-hand side of \eqref{Lplq} also increases as any of the $\abs{f_{\vj}(x)}$ increases pointwise. Hence, without loss of generality, we may further assume that \eqref{eq:f<g} holds as an equality. It then follows that
\begin{equation*}
  2^{\vj\cdot\vs_1}\abs{f_{\vj}(x)}
  =2^{\vj\cdot(\frac{\vn}{p_1}-\frac{\vn}{p_0})}2^{\vj\cdot\vs_0}\abs{g_{\vj}(x)}
  =:2^{\vj\cdot\vn(\frac{1}{p_1}-\frac{1}{p_0})}\abs{\widetilde g_{\vj}(x)}.
\end{equation*}
Thus, after these reductions, in place of \eqref{Lplq}, it is enough to prove that
\begin{equation}\label{Lplq2}
  \bNorm{\big\{ 2^{\vj\cdot\vn (\frac{1}{p_1}-\frac{1}{p_0})} \widetilde g_{\vj} \big\}_{\vj\in\Z^k} }{L^{p_1}\ell^q}
  \lesssim \Norm{\{ \widetilde g_{\vj} \}_{\vj\in\Z^k} }{L^{p_0}\ell^\infty},\qquad
  0<q<p_0<p_1<\infty.
\end{equation}
Moreover, denoting
\begin{equation*}
    r_i:=\frac{p_i}{q}\in(1,\infty)\quad(i=0,1),\qquad h_{\vj}(x):=\abs{\tilde g_{\vj}(x)}^q,\qquad
    H(x):=\Norm{\{  h_{\vj}(x) \}_{\vj\in\Z^k} }{\ell^\infty},
\end{equation*}
we find that
\begin{equation*}
\begin{split}
  \operatorname{LHS}\eqref{Lplq2}^q
  &= \bNorm{\big\{ 2^{\vj\cdot\vn(\frac{q}{p_1}-\frac{q}{p_0})} \abs{\tilde g_{\vj}}^q \big\}_{\vj\in\Z^k} }{L^{\frac{p_1}{q}}\ell^1}
  = \bNorm{\big\{ 2^{\vj\cdot\vn(\frac{1}{r_1}-\frac{1}{r_0})} h_{\vj} \big\}_{\vj\in\Z^k} }{L^{r_1}\ell^1}, \\
  \operatorname{RHS}\eqref{Lplq2}^q
  &= \bNorm{\{  \abs{\tilde g_{\vj}}^q \}_{\vj\in\Z^k} }{L^{\frac{p_0}{q}}\ell^\infty}
  = \bNorm{\{  h_{\vj} \}_{\vj\in\Z^k} }{L^{r_0}\ell^\infty}  =  \Norm{H}{L^{r_0}(\R^{\vn})},
\end{split}
\end{equation*}
here and thereafter, LHS \eqref{Lplq2} [resp. RHS \eqref{Lplq2}] means the left-hand
(resp. reight-hand) side of \eqref{Lplq2}.
Hence \eqref{Lplq2} is further equivalent to
\begin{equation}\label{Lplq3}
  \bNorm{\big\{ 2^{\vj\cdot\vn(\frac{1}{r_1}-\frac{1}{r_0})} h_{\vj} \big\}_{\vj\in\Z^k} }{L^{r_1}\ell^1}
  \lesssim\Norm{H}{L^{r_0}(\R^{\vn})},\qquad 1<r_0<r_1<\infty.
\end{equation}

Let us denote the constant value of $h_{\vj}$ on $R\in\mathscr D_{\vj}(\R^{\vn})$ by $h_R$. Then
\begin{equation*}
  h_R=\ave{h_{\vj}}_R\leq\ave{H}_R.
\end{equation*}
Noting that $2^{\vj\cdot\vn}=\abs{R}^{-1}$ for $R\in\mathscr D_{\vj}(\R^{\vn})$, we find that,
for any $x\in\R^{\vn}$,
\begin{equation*}
  \bNorm{\big\{ 2^{\vj\cdot\vn(\frac{1}{r_1}-\frac{1}{r_0})} h_{\vj}(x) \big\}_{\vj\in\Z^k} }{\ell^1}
  =\sum_{R\in\mathscr D(\R^{\vn})}\one_R(x)\abs{R}^{\frac{1}{r_0}-\frac{1}{r_1}}h_R
  \leq\sum_{R\in\mathscr D(\R^{\vn})}\one_R(x)\abs{R}^{\frac{1}{r_0}-\frac{1}{r_1}}\ave{H}_R
  =:I_{\mathscr D}^{\frac{1}{r_0}-\frac{1}{r_1}}H(x).
\end{equation*}
Hence, \eqref{Lplq3} follows from the inequality
\begin{equation}\label{Lplq4}
   \Norm{ I_{\mathscr D}^{\frac{1}{r_0}-\frac{1}{r_1}}H }{L^{r_1}(\R^{\vn})}
   \lesssim \Norm{H}{L^{r_0}(\R^{\vn})},\qquad 1<r_0<r_1<\infty,
\end{equation}
where $I_{\mathscr D}^{\frac{1}{r_0}-\frac{1}{r_1}}$ is a dyadic fractional integral-type operator.
While \eqref{Lplq4} is essentially known, we indicate an argument as follows:

For $s:=\frac{1}{r_0}-\frac{1}{r_1}\in(0,1)$,
the kernel of $I_{\mathscr D}^s$ at $(x,y)\in\R^{\vn}\times\R^{\vn}$ is
\begin{equation*}
  \sum_{\substack{R\in\mathscr D(\R^{\vn}) \\ x,y\in R}}\abs{R}^{s-1}
  =\prod_{\nu=1}^k\sum_{\substack{R_\nu\in\mathscr D(\R^{n_\nu}) \\ x_\nu,y_\nu\in R_\nu}}
     \abs{R_\nu}^{s-1}
  \lesssim\prod_{\nu=1}^k\abs{x_\nu-y_\nu}^{n_\nu(s-1)},
\end{equation*}
noting that each sum over $R_\nu$ is a decaying geometric series, where the leading term involves the minimal dyadic cube $R_\nu$ that contains both $x_\nu,y_\nu$ and hence satisfies $\operatorname{diam}(R_\nu)\geq\abs{x_\nu-y_\nu}$. The right-hand side of the previous display is the kernel of the product fractional integral $I_{\vec\alpha}$ with $\vec\alpha=s\vn$, using the notation of \cite[p.~2]{Tanaka}. By the pointwise bound between the kernels and \cite[Corollary 1.2]{Tanaka}, it follows that
\begin{equation*}
  \Norm{I^{\frac{1}{r_0}-\frac{1}{r_1}}_{\mathscr D}H}{L^{r_1}(\R^{\vn})}
  \lesssim\Norm{I_{(\frac{1}{r_0}-\frac{1}{r_1})\vn}H}{L^{r_1}(\R^{\vn})}
  \lesssim\Norm{H}{L^{r_0}(\R^{\vn})};
\end{equation*}
in fact, \cite[Corollary 1.2]{Tanaka} characterises a weighted version of the last inequality, and the characterising condition in \cite[Corollary 1.2(b)]{Tanaka} is obvious for $w\equiv 1$.

This finishes the proof of \eqref{Lplq4}. Thanks to the reductions leading to \eqref{Lplq4}, we have also proved \eqref{Lplq},
and hence Proposition \ref{thm-sobolev-prop}.
\end{proof}

\begin{theorem} \label{thm-sobolev-B}
Let $\vec s_0,\vec s_1\in\R^k$, $\tau\in[0,\infty)$, and $\vq\in(0,\infty]^k$.
Assume that $V_0\in \A_{p_0}(\mathbb R^{\vec n})$ and $V_1\in \A_{p_1}(\mathbb R^{\vec n})$.
If $0<p_0< p_1<\infty$, then the following assertions are equivalent.
\begin{enumerate}[\rm(i)]
\item\label{it:Sob-i}
$
\dot{B}_{p_0,\vq}^{\vec s_0,\tau}(V_0)
\hookrightarrow \dot{B}_{p_1,\vq}^{\vec s_1,\tau}(V_1);
$

\item\label{it:Sob-ii}
$
\dot{F}_{p_0,\infty}^{\vec s_0,\tau}(V_0)
\hookrightarrow \dot{F}_{p_1,\vq}^{\vec s_1,\tau}(V_1);
$

\item\label{it:Sob-iii}
$
\dot{b}_{p_0,\vq}^{\vec s_0,\tau}(V_0)
\hookrightarrow \dot{b}_{p_1,\vq}^{\vec s_1,\tau}(V_1);
$

\item\label{it:Sob-iv}
$
\dot{f}_{p_0,\infty}^{\vec s_0,\tau}(V_0)
\hookrightarrow \dot{f}_{p_1,\vq}^{\vec s_1,\tau}(V_1);
$

\item\label{it:Sob-v} there exists a positive constant $C$
such that, for every $P\in\mathscr D(\mathbb R^{\vec n})$ and $z\in\mathbb{C}^{m}$,
\begin{align}\label{eq-sobolev-B}
  2^{\vj_P \cdot (\vec s_1-\frac{\vec n}{p_1})}\|V_1 z\|_{\aveL^{p_1}(P)}
  \leq C 2^{\vj_P \cdot (\vec s_0-\frac{\vec n}{p_0})} \|V_0 z\|_{\aveL^{p_0}(P)}.
\end{align}
\end{enumerate}
\end{theorem}

\begin{proof}
By similarity, we only prove \eqref{it:Sob-ii} $\Longleftrightarrow$ \eqref{it:Sob-iv} $\Longleftrightarrow$ \eqref{it:Sob-v}.
Now, we show \eqref{it:Sob-v} $\Longrightarrow$ \eqref{it:Sob-iv}.
Inequality \eqref{eq-sobolev-B} implies that
$$
2^{\vj_P \cdot (\vec s_1-\frac{\vec n}{p_1})}
|[V_1]_{\aveL^{p_1}(P)}z|
\lesssim 2^{\vj_P \cdot (\vec s_0-\frac{\vec n}{p_0})}
|[V_0]_{\aveL^{p_0}(P)}z|
$$
and hence, for all sequences $\{t_Q\}_{Q\in\mathscr D(\R^{\vn})}$,
$$
2^{\vj \cdot (\vec s_1-\frac{\vec n}{p_1})}
\mathbf 1_{\Omega_{\vj}}|[V_1]_{\vj,p_1} t_{\vj}|
\lesssim 2^{\vj \cdot (\vec s_0-\frac{\vec n}{p_0})}
\mathbf 1_{\Omega_{\vj}}|[V_0]_{\vj,p_0} t_{\vj}|.
$$
This, together with Proposition \ref{thm-sobolev-prop},
further implies that
\begin{align*}
\|t\|_{\dot{f}_{p_1,q}^{\vec s_1,\tau}([V_1]_{p_1})}
&=\sup_{\Omega\in\Open(\R^{\vn})} |\Omega|^{-\tau}
\big\|\big\{2^{\vj\cdot \vec s_1} \mathbf 1_{\Omega_{\vj}}|[V_1]_{\vj,p_1} t_{\vj}|\big\}_{\vj \in\mathbb Z^k}\big\|_{L^{p_1}\ell^q} \\
&\lesssim \sup_{\Omega\in\Open(\R^{\vn})} |\Omega|^{-\tau}
\big\|\big\{2^{\vj\cdot \vec s_0} \mathbf 1_{\Omega_{\vj}}|[V_0]_{\vj,p_0} t_{\vj}|\big\}_{\vj \in\mathbb Z^k}\big\|_{L^{p_0}\ell^\infty}
=\|t\|_{\dot{f}_{p_0,\infty}^{\vec s_0,\tau}([V_0]_{p_0})},
\end{align*}
which, together with Theorem \ref{3 norms thm}, further implies \eqref{it:Sob-iv}.

Next, we prove \eqref{it:Sob-iv} $\Longrightarrow$ \eqref{it:Sob-ii}.
By the $\varphi$ transform characterisation (Theorem \ref{phi}),
\begin{equation*}
  \|f\|_{\dot{F}_{p_1,\vq}^{\vec s_1,\tau}(V_1)}
  \sim\Norm{S_\varphi f}{\dot f_{p_1,\vq}^{\vec s_1,\tau}(V_1)}
  \lesssim \Norm{S_\varphi f}{\dot{f}_{p_0,\infty}^{\vec s_0,\tau}(V_0)_{\mathrm w}}
  \sim \|f\|_{\dot{F}_{p_0,\infty}^{\vec s_0,\tau}(V_0)}
\end{equation*}
and hence \eqref{it:Sob-ii} holds.

Finally, we show \eqref{it:Sob-ii} $\Longrightarrow$ \eqref{it:Sob-v}. Let $\{\theta^{(\vi)}\}_{\vi\in[2^{\vn}-\vec1]}$ be Meyer wavelets on $\R^{\vn}$.
For every $P\in\mathscr D(\mathbb R^{\vec n})$ and $z\in\mathbb{C}^{m}$,
let $f:=z\theta_P^{(\vi)}$
for some $\vi_0\in[2^{\vn}-\vec1]$.
By orthonormality, it follows that the wavelet coefficients of $f$ are given by $\pair{f}{\theta_Q^{(\vi)}}=z\delta_{\vi_0,\vi}\delta_{P,Q}$

From Meyer wavelet characterisation (Theorem \ref{wavelet}) of $\Atau(V)$, it follows that
\begin{equation*}
  \Norm{f}{ \dot F^{\vs_i,\tau}_{p_i,\vq_i}(V_i)}
  \sim\Norm{ \{z\delta_{P,Q}\}_{Q\in\mathscr D(\R^{\vn})} }{\dot f^{\vs_i,\tau}_{p_i,\vq_i}(V_i)}
  \sim\abs{P}^{-\tau} 2^{\vj_P\cdot\vs_i}\Norm{z\widetilde\one_P}{L^{p_i}(V_i)}
  =\abs{P}^{-\tau-\frac12} 2^{\vj_P\cdot(\vs_i-\vn/p_i)}\Norm{V_i z}{\aveL^{p_i}(P)}
\end{equation*}
Using this with both $i=0,1$ and $\vq_0=\infty$, $\vq_1=\vq$, the assumption \eqref{it:Sob-ii} implies
\begin{equation*}
  \abs{P}^{-\tau-\frac12} 2^{\vj_P\cdot(\vs_1-\vn/p_1)}\Norm{V_1 z}{\aveL^{p_1}(P)}
  \sim\Norm{f}{ \dot F^{\vs_1,\tau}_{p_1,\vq_1}(V_1)}
  \lesssim\Norm{f}{ \dot F^{\vs_0,\tau}_{p_0,\vq_0}(V_0)}
  \sim\abs{P}^{-\tau-\frac12} 2^{\vj_P\cdot(\vs_0-\vn/p_0)}\Norm{V_0 z}{\aveL^{p_0}(P)}.
\end{equation*}
Cancelling out the common factors $\abs{P}^{-\tau-\frac12}$, we get \eqref{it:Sob-v}
and finish the proof of Theorem \ref{thm-sobolev-B}.
\end{proof}

\begin{remark}
\begin{enumerate}[\rm(i)]
\item Theorem \ref{thm-sobolev-B} is new in the multi-parameter case. In the one-parameter case, it essentially agrees with \cite[Theorem 6.11]{BYYZ1} for the usual Besov and Triebel--Lizorkin space with $\tau=0$, and slightly more generally (see Theorem \ref{compare A}), while in general it deals with a different scale of Besov--Triebel--Lizorkin-type spaces. We note that \cite[Theorem 6.11]{BYYZ1} also addresses more general $\A_{p,\infty}$ (instead of $\A_p$) matrix weights, which have not been studied in the multi-parameter setting.

\item In the one-parameter case, another approach to Theorem \ref{thm-sobolev-B} is via
an interpolation inequality. However, this inequality is no longer
valid in the multi-parameter case (see Proposition \ref{BL fail} below).
That is why we first use the boundedness of the
dyadic fractional integral-type operator to prove Proposition \ref{thm-sobolev-prop},
and then apply Proposition \ref{thm-sobolev-prop} to prove Theorem \ref{thm-sobolev-B}.
\end{enumerate}
\end{remark}

\begin{lemma}\label{sum-min}
For $a,b\in(0,\infty)$ and $A,B\in(0,\infty)$, we have
\begin{equation*}
  \sum_{j\in\Z}\min( 2^{ja}A,2^{-jb}B)\sim A^{\frac{b}{a+b}}B^{\frac{a}{a+b}},
\end{equation*}
where the positive equivalence constants may depend on $a,b$, but independent of $A,B$.
\end{lemma}

\begin{proof}
Note that $2^{ja}A\leq 2^{-jb}B$ if and only if $2^j\leq(B/A)^{\frac{1}{a+b}}$. Hence
\begin{equation*}
\operatorname{LHS}
=\sum_{2^j\leq (B/A)^{\frac{1}{a+b}}}2^{ja}A+\sum_{2^j>(B/A)^{\frac{1}{a+b}}}2^{-jb}B
\sim (B/A)^{\frac{a}{a+b}}A+(B/A)^{-\frac{b}{a+b}}B\sim A^{\frac{b}{a+b}}B^{\frac{a}{a+b}}.
\end{equation*}
This finishes the proof of Lemma \ref{sum-min}.
\end{proof}

\begin{proposition}\label{BL fail}
Let $\vec s_0,\vec s_1\in\R^k$.
Assume that $\theta\in(0,1)$,
$\vec s:=\theta \vec s_0+(1-\theta) \vec s_1$, and $q\in(0,\infty]$.
\begin{enumerate}[\rm(i)]
\item\label{k=1} If $q=\infty$ or if $k=1$, $q\in(0,\infty)$,
and $\vec s_0\neq\vec s_1$, then
there exists a positive constant $C$ such that,
for every sequence $\{c_j\}_{j\in\Z}$ in $\mathbb C$,
\begin{equation}\label{false eq}
\Big[\sum_{\vj\in\mathbb Z^k} 2^{(\vj\cdot \vec s)q}
|a_{\vj}|^q\Big]^{\frac 1q}
\leq C \Big[\sup_{\vj\in\mathbb Z^k} 2^{\vj\cdot \vec s_0}
|a_{\vj}|\Big]^{\theta}
\Big[\sup_{\vj\in\mathbb Z^k} 2^{\vj\cdot \vec s_1}
|a_{\vj}|\Big]^{1-\theta}.
\end{equation}

\item\label{k>1} If $k\in\mathbb N\cap[2,\infty)$ and $q\in(0,\infty)$,
then \eqref{false eq} is not valid for any $\vec s_0,\vec s_1\in\R^k$.
\end{enumerate}
\end{proposition}

\begin{proof}
\eqref{k=1}:
Let first $q=\infty$. Then
\begin{equation*}
  2^{\vj\cdot\vs}\abs{a_{\vj}}
  =(2^{\vj\cdot\vs_0}\abs{a_{\vj}})^{\theta}(2^{\vj\cdot\vs_1}\abs{a_{\vj}})^{1-\theta},
\end{equation*}
and \eqref{false eq} with $q=\infty$ follows by taking the supremum on both sides of the above equality over $j\in\Z^k$.

Let then $k=1$ and $q\in(0,\infty)$. We may assume by symmetry that $s_0<s_1$, and hence $s_0<s<s_1$. Moreover, we have $s_1-s=\theta(s_1-s_0)$ and $s-s_0=(1-\theta)(s_1-s_0)$.

Let $A:=\sup_{j\in\Z} 2^{js_0}\abs{a_j}$ and $B:=\sup_{j\in\Z} 2^{js_1}\abs{a_j}$. Thus, $\abs{a_j}\leq\min\{2^{-js_0}A,2^{-ja_1}B\}$,
and then
\begin{equation*}
\begin{split}
  \operatorname{LHS}\eqref{false eq}^q
  &\leq\sum_{j\in\Z}\min\{2^{j(s-s_0)q}A^q,2^{-j(s_1-s)q}B^q\} \\
  &\sim (A^q)^{\frac{(s_1-s)q}{(s_1-s_0)q}}(B^q)^{\frac{(s-s_0)q}{(s_1-s_0)q}}
  \qquad\text{by Lemma \ref{sum-min}}\\
  &=A^{\frac{(s_1-s)q}{s_1-s_0}}B^{\frac{(s-s_0)q}{s_1-s_0}}
  =A^{\theta q}B^{(1-\theta)q}
  =\operatorname{RHS}\eqref{false eq}^q.
\end{split}
\end{equation*}
This finishes the proof of (i).

\eqref{k>1}:
To avoid complicated indexing, let $\vr:=\vs_0$ and $\vt:=\vs_1$. We write all vectors $\vec u\in\Z^k$ as $\vec u=(u_1,\vec u\,')\in\Z\times\Z^{k-1}$.
We may assume by symmetry that $r_1<t_1$ and hence $r_1<s_1<t_1$.

Given $\vs_0,\vs_1\in\R^k$, we define the sequence $a_{\vj}:=\min\{2^{-\vj\cdot\vs_0},2^{-\vj\cdot\vs_1}\}=\min\{2^{-\vj\cdot\vr},2^{-\vj\cdot\vt}\}$. For both $h=0,1$, then clearly $2^{-\vh\cdot\vs_h}a_{\vj}\leq 1$ for all $\vj\in\Z^k$ with equality when $\vj=\vnull$, hence $\operatorname{RHS}\eqref{false eq}=1$.

On the other hand, writing $\vj=(\vj_1,\vj\,')\in\Z\times\Z^{k-1}$, we have
\begin{equation*}
\begin{split}
  \operatorname{LHS}\eqref{false eq}^q
  &=\sum_{\vj\,'\in\Z^{k-1}} \sum_{j_1\in\Z} \min \Big\{ 2^{j_1(s_1-r_1)q} 2^{\vj\,'\cdot(\vs\,'-\vr\,')q} ,
     2^{-j_1(t_1-s_1)q} 2^{-\vj\,'\cdot(\vt\,'-\vs\,')q}\Big\} \\
  &\sim\sum_{\vj\,'\in\Z^{k-1}} 2^{\vj\,'\cdot(\vs\,'-\vr\,')q\frac{t_1-s_1}{t_1-r_1}}
    2^{-\vj\,'\cdot(\vt\,'-\vs\,')q\frac{s_1-r_1}{t_1-r_1}}\qquad\text{by Lemma \ref{sum-min}} \\
    &=\sum_{\vj\,'\in\Z^{k-1}} 2^{\vj\,'\cdot\vec u'}
    =\prod_{i=2}^k\Bigg(\sum_{j_i\in\Z}2^{j_i u_i}\Bigg)
\end{split}
\end{equation*}
for some $\vec u'=(u_2,\ldots,u_k)\in\R^{k-1}$. (One can see an exact formula for $\vec u'$ above, but this is irrelevant.)

It is clear that the doubly infinite geometric series $\sum_{j_i\in\Z}2^{j_i u_i}$ diverges for all values of $u_i\in\R$. Hence $\operatorname{LHS}\eqref{false eq}=\infty$, while we already observed that $\operatorname{RHS}\eqref{false eq}=1$.
This finishes the proof of (ii) and hence Proposition \ref{BL fail}.
\end{proof}

\subsection{Applications to matrix-weighted $L^p$ and Sobolev spaces} \label{sec:lift}

Based on the relationships with other known function spaces
established in Section \ref{other}, even when the results obtained in the previous
sections are applied to some specific function spaces, we can also
obtain some further interesting consequences, for example,
\begin{enumerate}[\rm(i)]
\item using the boundedness of lifting operators (Theorem \ref{257V}),
we can obtain the boundedness of fractional integral operators on
Lebesgue spaces with $\A_{p,q}(\mathbb R^{\vec n})$ weights, and

\item applying the Sobolev-type embedding of Besov--Triebel--Lizorkin-type spaces in Theorem \ref{thm-sobolev-B}, we can obtain the Sobolev embedding theorem.
\end{enumerate}

For every $\vec\sigma\in (\vec 0,\vec n):=(0,n_1)\times\cdots\times(0,n_k)$,
the \emph{multi-parameter fractional integral operator} $\mathscr I_{\vec\sigma}$
is defined by setting, for all suitable functions $f$ on $\R^{\vn}$ and all $x\in\mathbb R^{\vec n}$,
$$
\mathscr I_{\vec\sigma}(f)(x)
:=\int_{\mathbb R^{\vec n}} \frac{f(y)}{|x-y|^{\vec n-\vec\sigma}}\,dy,
$$
where $|x-y|^{\vec n-\vec\sigma}:=\prod_{\nu=1}^k|x_\nu-y_\nu|^{n_\nu-\sigma_\nu}$.
Let $p,q\in(1,\infty)$. A matrix weight $V$ is called
an \emph{$\A_{p,q}$-matrix weight} if $V$ satisfies that
$$
[V]_{\A_{p,q}(\mathbb R^{\vec n})}:=\sup_{P\in\Rect(\mathbb R^{\vec n})}
\|V,V^{-1}\|_{\aveL^{q}(P),\aveL^{p'}(P)}<\infty.
$$

\begin{theorem} \label{fractional}
Let $\vec\sigma\in (\vec 0,\vec n)$ and $p,q\in(1,\infty)$ satisfy
\begin{equation}\label{sharp condition}
\big(\frac 1p-\frac 1q\big) \vec n=\vec \sigma.
\end{equation}
If $V\in \A_{p,q}(\R^{\vn})$, then there exists
a positive constant $C$ such that, for every $f\in L^p(V)$,
\begin{equation}\label{boundedness}
\|\mathscr I_{\vec\sigma} f\|_{L^q(V)}
\leq C \|f\|_{L^p(V)}.
\end{equation}
\end{theorem}

\begin{remark}
In the case where $k=2$ and $m=1$,
Theorem \ref{fractional} coincides with a part of \cite[Theorem 5]{SW17},
in which Sawyer and Wang established several equivalent characterizations of the boundedness of multi-parameter fractional integral operators.
In particular, they showed that if
$\mathscr I_{\vec\sigma}:\ L^p(V)\to L^q(V)$ is bounded,
then condition \eqref{sharp condition} holds.
\end{remark}

To prove Theorem \ref{fractional}, we need some lemmas.

\begin{lemma} \label{8.19}
If $1<p<q<\infty$ and $V\in \A_{p,q}(\R^{\vn})$,
then $V\in \A_p(\R^{\vn}) \cap \A_q(\R^{\vn})$.
\end{lemma}

\begin{proof}
By H\"older's inequality, we find that
\begin{align*}
[V]_{\A_p(\mathbb R^{\vec n})}
=\sup_{P\in\Rect(\mathbb R^{\vec n})}
\|V,V^{-1}\|_{\aveL^{p}(P),\aveL^{p'}(P)}
\leq\sup_{P\in\Rect(\mathbb R^{\vec n})}
\|V,V^{-1}\|_{\aveL^q(P),\aveL^{p'}(P)}
=[V]_{\A_{p,q}(\mathbb R^{\vec n})},
\end{align*}
and hence $V\in \A_p(\mathbb R^{\vec n})$.
Similarly, we obtain $V\in \A_q(\mathbb R^{\vec n})$.
This finishes the proof of Lemma \ref{8.19}.
\end{proof}

The following lemma provides an equivalent characterization of matrix $\A_{p,q}(\mathbb R^{\vec n})$ weights.

\begin{lemma} \label{op-norm}
Let $p,q\in(1,\infty)$ and $V$ be a matrix weight.
Then
\begin{equation*}
[V]_{\A_{p,q}(\mathbb R^{\vec n})}
\sim\sup_{P\in\Rect(\mathbb R^{\vec n})}
\Norm{f\mapsto\one_P\ave{f}_P}{\aveL^{p}(P,V)\to\aveL^q(P,V)},
\end{equation*}
where the positive equivalence constants are independent of $V$.
\end{lemma}

\begin{proof}
Note that $Y=\aveL^{p'}(P)$ is a Banach lattice with a dual pair $Y'=\aveL^p(P)$, with respect to the normalised integration $\ave{\ }_P:=\fint_P$ over $P$.
Hence it follows from Lemma \ref{red Fub}(ii) and a change of variables that
\begin{align*}
\Norm{V,V^{-1}}{\aveL^q(P),\aveL^{p'}(P)}
\sim \Norm{f\mapsto V\ave{V^{-1}f}_P}{\aveL^{p}(P)\to\aveL^q(P)}
=\Norm{f\mapsto \one_P\ave{f}_P}{\aveL^{p}(P,V)\to\aveL^q(P,V)}.
\end{align*}
Taking the supremum over $P\in\Rect(\mathbb R^{\vec n})$ on both sides of
the above equation yields the desired estimate.
This finishes the proof of Lemma \ref{op-norm}.
\end{proof}

Applying Lemma \ref{op-norm}, we obtain
a useful relation between scalar weights and matrix weights.

\begin{lemma} \label{scalar}
Let $p,q\in(1,\infty)$ and $V\in \A_{p,q}(\R^{\vn})$.
Then there exists a positive constant $C$ such that,
for every $M\in \Pdm\setminus\{O_m\}$, we have
$[|VM|]_{\A_{p,q}(\R^{\vn})}\leq C [V]_{\A_{p,q}(\R^{\vn})}$.
\end{lemma}

\begin{proof}
By Lemma \ref{op-norm}, we have
\begin{equation} \label{sim}
\begin{split}
[\abs{VM}]_{\A_{p,q}(\R^{\vn})}
&\sim \sup_{P\in\Rect(\R^{\vn})}\Norm{f\mapsto\one_P\ave{f}_P}{\aveL^p(P,\abs{VM})\to \aveL^q(P,\abs{VM})}, \\
[V]_{\A_{p,q}(\R^{\vn})}
&\sim \sup_{P\in\Rect(\R^{\vn})}\Norm{F\mapsto\one_P\ave{F}_P}{\aveL^p(P,V)\to \aveL^q(P,V)}.
\end{split}
\end{equation}

Given $M\in \Pdm\setminus\{O_m\}$ and $f\in\aveL^p(P,\abs{VM})$,
we consider the matrix-valued function $F:=f(\cdot) M$. Then
\begin{equation*}
\begin{split}
\Norm{\one_P\ave{f}_P}{\aveL^q(P,\abs{VM})}
&=\bNorm{\one_P\abs{V\ave{f(\cdot) M}_P}}{\aveL^q(P)}
=\Norm{\one_P\ave{f(\cdot)M}_P}{\aveL^p(P,V)} \\
&\leq\Norm{F\mapsto \one_P\ave{F}_P}{\aveL^p(P,V)\to\aveL^q(P,V)}
\Norm{f(\cdot)M}{\aveL^q(P,V)},
\end{split}
\end{equation*}
where
\begin{equation*}
\Norm{f(\cdot)M}{\aveL^q(P,V)}
=\bNorm{f(\cdot)\abs{V(\cdot)M}}{\aveL^q(P)}
=\Norm{f}{\aveL^q(P,\abs{VM})}.
\end{equation*}
This shows that
\begin{equation*}
\Norm{f\mapsto\one_P\ave{f}_P}{\aveL^p(P,\abs{VM})\to \aveL^q(P,\abs{VM})}
\leq\Norm{F\mapsto\one_P\ave{F}_P}{\aveL^p(P,V)\to \aveL^q(P,V)},
\end{equation*}
which, together with \eqref{sim}, then completes the proof of Lemma \ref{scalar}.
\end{proof}

The following lemma shows that the fractional integral operator $\mathscr I_{\vec\sigma}$
and the lifting operator $\dot I_{\vec\sigma}$, defined in \eqref{lift},
are equivalent in the sense of distributions.

\begin{lemma} \label{8.19y}
Let $\vec\sigma\in (\vec 0,\vec n)$ and $p,q\in(1,\infty)$ satisfy \eqref{sharp condition}.
Then there exists a positive constant $C$,
depending only on $\vec\sigma$ and $\vn$, such that,
for every $V\in \A_{p,q}(\R^{\vn})$ and $f\in L^p(V)$,
$\mathscr I_{\vec\sigma} f$ is well defined and
$\mathscr I_{\vec\sigma} f= C \dot I_{\vec\sigma} f$
in $\mathscr S'(\R^{\vn})$.
\end{lemma}

\begin{proof}
We first show that $\mathscr I_{\vec\sigma} f$ is well defined.
By Tonelli's theorem and Holder's inequality, we conclude that, for every $\varphi\in \mathscr S(\R^{\vn})$,
\begin{equation} \label{323}
\begin{split}
\int_{\mathbb R^{\vec n}} \mathscr I_{\vec\sigma} (|f|)(x) |\varphi(x)| dx
&= \int_{\mathbb R^{\vec n}} |f(x)| \mathscr I_{\vec\sigma} (|\varphi|)(x)  dx \\
&\leq \int_{\mathbb R^{\vec n}} |V^{-1}(x)|\,|V(x)f(x)| \mathscr  I_{\vec\sigma} (|\varphi|)(x) dx \\
&\leq \|f\|_{L^p(V)} \| \mathscr I_{\vec\sigma} (|\varphi|) \|_{L^{p'}(w)},
\end{split}
\end{equation}
where $w:=|V^{-1}|$.
From the definition of $\A_{p,q}(\R^{\vn})$
and Lemmas \ref{red Fub}(i) and \ref{scalar}, it follows that
\begin{align*}
V\in \A_{p,q}(\R^{\vn})
\Longleftrightarrow V^{-1}\in \A_{q',p'}(\R^{\vn})
\Longrightarrow w\in \A_{q',p'}(\mathbb R^{\vec n}).
\end{align*}
Note that
$$
\Bigg(\frac 1{q'}-\frac 1{p'}\Bigg) \vec n
=\Bigg(\frac 1p-\frac 1q\Bigg) \vec n
=\vec \sigma.
$$
Using the boundedness of $\mathscr I_{\vec\sigma}$ (see \cite[Theorem 5]{SW17}), we find that
\begin{align*}
\| \mathscr I_{\vec\sigma} (|\varphi|) \|_{L^{p'}(w)}
\lesssim \|\varphi\|_{L^{q'}(w)} <\infty,
\end{align*}
which, togeter with \eqref{323}, further implies that
\begin{align} \label{bound}
\int_{\mathbb R^{\vec n}} \mathscr I_{\vec\sigma} (|f|)(x) |\varphi(x)| dx <\infty.
\end{align}
Therefore, $\mathscr I_{\vec\sigma} f$ is well defined.

Moreover, using a straightforward multi-parameter
extension of \cite[p.\,88]{Duo},
we find that there exists a positive constant $C$,
depending only on $\vec \sigma$ and $\vec n$, such that,
for every $\varphi\in \mathscr S(\R^{\vn})$ and almost every
$x\in\R^{\vn}$,
$\mathscr I_{\vec\sigma} \varphi(x)=C\dot I_{\vec\sigma} \varphi(x)$.
From this, \eqref{bound}, and Fubini's theorem,
we deduce that,
for every $\varphi\in \mathscr S(\R^{\vn})$,
\begin{align*}
\langle \mathscr I_{\vec\sigma} f, \varphi \rangle
&=\int_{\mathbb R^{\vec n}} \mathscr I_{\vec\sigma} (f)(x) \varphi(x) dx
= \int_{\mathbb R^{\vec n}} f(x) \mathscr I_{\vec\sigma} (\varphi)(x) dx \\
&= \langle f, \mathscr I_{\vec\sigma} \varphi \rangle
= \langle f, C \dot I_{\vec\sigma} \varphi \rangle
= C \langle \dot I_{\vec\sigma} f, \varphi \rangle,
\end{align*}
and hence $\mathscr I_{\vec\sigma} f= C \dot I_{\vec\sigma} f$
in $\mathscr S'(\R^{\vn})$. This finishes the proof of Lemma \ref{8.19y}.
\end{proof}

Now, we can prove Theorem \ref{fractional}.

\begin{proof}[Proof of Theorem \ref{fractional}]
Let $f\in L^p(V)$.
Note that $V\in \A_p(\R^{\vn}) \cap \A_q(\R^{\vn})$ (see Lemma \ref{8.19}).
By Lemma \ref{8.19y}, Theorem \ref{Lp=F}, and $V\in \A_q(\R^{\vn})$, we find that
$$
\|\mathscr I_{\vec\sigma} f\|_{L^{q}(V)}
\sim\Norm{\dot I_{\vec\sigma} f}{\dot F^{0,0}_{q,2}(V)}.
$$
For every $x\in\mathbb R^{\vec n}$,
$P\in\mathscr D(\mathbb R^{\vec n})$, and $z\in\C^m$,
$$
|V(x)z|\leq \fint_P |V(x)V(y)^{-1}|\, |V(y)z| \, dy
\leq \|V(x)V(\cdot)^{-1}\|_{\aveL^{p'}(P)} \|V(\cdot)z\|_{\aveL^{p}(P)},
$$
which, together with \eqref{sharp condition}, further implies that
\begin{align*}
2^{\vj_P \cdot (\vec 0-\frac{\vec n}{q})}\|z\mathbf 1_P\|_{\aveL^{q}(V)}
&\leq 2^{\vj_P \cdot (\vec \sigma-\frac{\vec n}{p})}\|V,V^{-1}\|_{\aveL^{q}(P),\aveL^{p'}(P)}
\|z\mathbf 1_P\|_{\aveL^{p}(V)}\\
&\leq [V]_{\A_{p,q}(\mathbb R^{\vec n})}
2^{\vj_P \cdot (\vec \sigma-\frac{\vec n}{p})}
\|z\mathbf 1_P\|_{\aveL^{p}(V)}.
\end{align*}
From this and Theorem \ref{thm-sobolev-B}, we deduce that
$$
\Norm{\dot I_{\vec\sigma} f}{\dot F^{0,0}_{q,2}(V)}
\lesssim\Norm{\dot I_{\vec\sigma} f}{\dot F^{\vec\sigma,0}_{p,\infty}(V)}
\leq\Norm{\dot I_{\vec\sigma} f}{\dot F^{\vec\sigma,0}_{p,2}(V)},
$$
where the last inequality is obvious.
Applying Theorems \ref{257V} and \ref{Lp=F} and $V\in \A_p(\R^{\vn})$, we obtain
$$
\Norm{\dot I_{\vec\sigma} f}{\dot F^{\vec\sigma,0}_{p,2}(V)}
\sim\Norm{f}{\dot F^{0,0}_{p,2}(V)}
\sim\|f\|_{L^p(V)}.
$$
Combining the above results yields \eqref{boundedness}.
This finishes the proof of Theorem \ref{fractional}.
\end{proof}

Unlike Theorem \ref{fractional}, the Sobolev embedding below
follows directly from Theorems \ref{thm-sobolev-B} and \ref{thm:Sobolev}.

\begin{theorem} \label{sobolev embedding}
Let $\vec \sigma,\vec l\in \mathbb N^k$, $1<p<q<\infty$, and
$V_0\in \A_p(\mathbb R^{\vec n})$, $V_1\in \A_q(\mathbb R^{\vec n})$.
If, for every $P\in\mathscr D(\mathbb R^{\vec n})$ and $z\in\mathbb{C}^{m}$,
\begin{align}\label{cond}
  2^{\vj_P \cdot (\vec l-\frac{\vec n}{q})}\|V_1 z\|_{\aveL^{q}(P)}
  \leq C 2^{\vj_P \cdot (\vec \sigma-\frac{\vec n}{p})} \|V_0 z\|_{\aveL^{p}(P)},
\end{align}
then $\dot W^{\vec \sigma}_p(V_0)\hookrightarrow \dot W^{\vec l}_q(V_1)$.
\end{theorem}

\begin{proof}
By Theorem \ref{thm:Sobolev}, we conclude that
\begin{equation*}
\dot W^{\vec\sigma}_p(V_0)
= \dot F^{\vec\sigma,0}_{p,2}(V_0)
\quad\text{and}\quad
\dot W^{\vec l}_q(V_1)= \dot F^{\vec l,0}_{q,2}(V_1).
\end{equation*}
Condition \eqref{cond} ensures that condition \eqref{it:Sob-v} of Theorem \ref{thm-sobolev-B} holds with
\begin{equation*}
  \vec s_0=\vec\sigma,\quad\vec s_1=\vec l,\quad\tau=0,\quad p_0=p,\quad p_1=q.
\end{equation*}
Since the conditions \eqref{it:Sob-v} and \eqref{it:Sob-ii} of Theorem \ref{thm-sobolev-B} are equivalent, we obtain
$\dot F^{\vec\sigma,0}_{p,\infty}(V_0)
\hookrightarrow \dot F^{\vec l,0}_{q,2}(V_1)$.
Together the simple embedding
$\dot F^{\vec\sigma,0}_{p,2}(V_0) \hookrightarrow
\dot F^{\vec\sigma,0}_{p,\infty}(V_0)$, this
finishes the proof of Theorem \ref{sobolev embedding}.
\end{proof}

\begin{corollary}
Let $\vec \sigma,\vec l\in \mathbb N^k$ and $0<p<q<\infty$ satisfy
\begin{equation}\label{eq:sigma-l}
  (\frac{1}{p}-\frac{1}{q})\vn=\vec\sigma-\vec l.
\end{equation}
Let $V\in\A_{p,q}(\R^{\vn})$. Then $\dot W^{\vec \sigma}_p(V)\hookrightarrow \dot W^{\vec l}_q(V)$.
\end{corollary}

\begin{proof}
The plan is to verify the assumptions of Theorem \ref{sobolev embedding} with $V_0=V_1=V$. By Lemma \ref{8.19}, the assumption $V\in\A_{p,q}(\R^{\vn})$ implies both $V_0:=V\in\A_p(\R^{\vn})$ and $V_1:=V\in\A_q(\R^{\vn})$. It remains to verify \eqref{cond}. Under the assumption \eqref{eq:sigma-l}, the powers of $2$ on both sides of \eqref{cond} are equal, so we are reduced to showing that
\begin{equation}\label{eq:cond2show}
  \Norm{Vz}{\aveL^q(P)}\lesssim\Norm{Vz}{\aveL^p(P)}.
\end{equation}
To prove \eqref{eq:cond2show}, we note that the assumption $V\in\A_{p,q}(\R^{\vn})$ also implies that
\begin{equation*}
\begin{split}
  \abs{[V]_{\aveL^q(P)}[V]_{\aveL^p(P)}^{-1}}
  &=\abs{[V]_{\aveL^p(P)}^{-1}[V]_{\aveL^q(P)}} \\
  &\lesssim\abs{[V^{-1}]_{\aveL^{p'}(P)}[V]_{\aveL^q(P)}}
  \quad\text{by Lemma \ref{inverses} with $u=p'$ and }M=[V]_{\aveL^q(P)} \\
  &=\abs{[V]_{\aveL^q(P)}[V^{-1}]_{\aveL^{p'}(P)}}
  \sim\Norm{V,V^{-1}}{\aveL^q(P),\aveL^{p'}(P)}\leq[V]_{\A_{p,q}(\R^n)}.
\end{split}
\end{equation*}
Hence
\begin{equation*}
\begin{split}
  \Norm{V z}{\aveL^{q}(P)}
  \sim\abs{[V]_{\aveL^q(P)}z}
  \leq\abs{[V]_{\aveL^q(P)}[V]_{\aveL^p(P)}^{-1}}\cdot\abs{[V]_{\aveL^p(P)}z}
  \lesssim\abs{[V]_{\aveL^p(P)}z} \sim\Norm{V z}{\aveL^p(P)}.
\end{split}
\end{equation*}
This proves \eqref{eq:cond2show}. Hence all the assumptions of Theorem \ref{sobolev embedding} are satisfied with $V_0=V_1=V$, and hence also its conclusion, which agrees with the claim of the corollary that we are proving. This concludes the proof.
\end{proof}

\appendix

\section{Carleson's counterexample revisited}\label{app:Carleson}

Carleson's counterexample \cite{Car} is a classical construction showing a delicate difference between versions of multi-parameter conditions defined by rectangles and by general open sets. To demonstrate the sharpness of some of our results in this work, and the necessity of using general open sets instead of rectangles in some formulations, we make use of both the original Carleson result and certain elaborations that we present in this appendix. For the setting of Carleson's example, we follow the presentation of Tao \cite{Tao}. The core of the construction is two-dimensional; a version in higher dimensions will be obtained as a simple corollary in the end.

\subsection{The core construction in two dimensions}

Let $\mathscr D$ be the collection of dyadic subrectangles of $[0,1)^2$. For all $P\in\mathscr D$ and $\mathscr R\subset\mathscr D$, we denote
\begin{equation*}
  \mathscr R(P):=\{R\in\mathscr R:R\subseteq P\},\quad
  \chi_{\mathscr R}:=\sum_{R\in\mathscr R}\one_R,\quad
  \sigma(\mathscr R):=\abs{\bigcup_{R\in\mathscr R}R}.
\end{equation*}

Following \cite{Tao}, we define:

\begin{definition}
A {\em quilt} is a finite collection $\mathscr R\subset\mathscr D$ such that
\begin{equation}\label{eq:rectCarleson}
  \Norm{\chi_{\mathscr R(P)}}{1}
  =\sum_{R\in\mathscr R(P)}\abs{R}\leq\abs{P}
\end{equation}
for every $P\in\mathscr D$, with equality when $P=[0,1)^2$.
\end{definition}

With this definition, the essence of Carleson's classical counterexample is as follows:

\begin{lemma}[\cite{Tao}, Lemma 7.2]\label{lem:Tao7.2}
The exist quilts $\mathscr R$ with arbitrarily small $\sigma(\mathscr R)$.
\end{lemma}

To explain the elaboration of Lemma \ref{lem:Tao7.2} that we need, let us first make a reformulation. Since every quilt has $\Norm{\chi_{\mathscr R}}{1}=1$, Lemma \ref{lem:Tao7.2} is equivalent to the claim that there exist quilts with arbitrarily large $\sigma(\mathscr R)^{-1}\Norm{\chi_{\mathscr R}}{1}$. Note that this quantity is the average of $\chi_{\mathscr R}$ on its support. We will prove the following generalisation:

\begin{proposition}\label{prop:CarlNew}
For every $p\in(0,\infty]$, the exist quilts $\mathscr R$ with arbitrarily large $\sigma(\mathscr R)^{-\frac 1p}\Norm{\chi_{\mathscr R}}{p}$.
\end{proposition}

Note that $\sigma(\mathscr R)^{-\frac 1p}\Norm{\chi_{\mathscr R}}{p}$ is increasing in $p$ by H\"older's inequality, and hence the range $p\in[1,\infty]$ of Proposition \ref{prop:CarlNew} is immediate from Lemma \ref{lem:Tao7.2}. The main point of Proposition \ref{prop:CarlNew} is the range $p\in(0,1)$. To prove Proposition \ref{prop:CarlNew}, we will revisit and elaborate the proof of Lemma \ref{lem:Tao7.2}. Indeed, we will see that the same construction used to prove Lemma \ref{lem:Tao7.2} already has the stronger properties of Proposition \ref{prop:CarlNew}, but verifying this requires additional effort.

The proof of Lemma \ref{lem:Tao7.2} is based on iterating the following:

\begin{lemma}[\cite{Tao}, proof of Lemma 7.2]\label{lem:Carleson}
For every quilt $\mathscr R$, there is another quilt $\mathscr R'$ such that
\begin{equation*}
  \sigma'=\sigma(1-\tfrac14\sigma),\qquad\sigma:=\sigma(\mathscr R),\quad
  \sigma'=\sigma(\mathscr R').
\end{equation*}
Moreover, such $\mathscr R'$ is obtained by the following Algorithm \ref{algo:quilt}.
\end{lemma}

We recall the full details of the algorithm in order to verify its additional properties.

\begin{algorithm}[\cite{Tao}, proof of Lemma 7.2]\label{algo:quilt}\mbox{}

{\bf Input:} a finite collection $\mathscr R\subset\mathscr D$.

{\bf Output:} another finite collection $\mathscr R'\subset\mathscr D$ constructed as follows:

Choose $N\in\N$ large enough such that the smallest side-length of every $R\in\mathscr R$ is at least $2^{-N}$; this is possible, since $\mathscr R$ is finite. For every $j=1,\ldots 2^{N+1}$, we consider the affine transformations
\begin{equation*}
  a_j:[0,1)\to[j-1,j)\cdot 2^{-N-1},\quad x\mapsto (j-1+x)\cdot 2^{-N-1}
\end{equation*}
and, in two dimensions,
\begin{equation*}
  A_j^1(x,y):=(a_j(x),y),\qquad A_j^2(x,y):=(x,a_j(y)).
\end{equation*}
Let
\begin{equation*}
  A_j^i(\mathscr R):=\{ A_j^i(R):R\in\mathscr R\}
\end{equation*}
and
\begin{equation*}
  \mathscr R^i:=\bigcup_{j=1}^{2^{N+1}}  A_j^i(\mathscr R),\qquad
  \mathscr S^i:=\bigcup_{\substack{j=1 \\ j\text{ odd}}}^{2^{N+1}}  A_j^i(\mathscr R),\qquad
  \mathscr R':=\mathscr S^1\cup\mathscr S^2.
\end{equation*}
This concludes the construction, and the algorithm outputs the collection $\mathscr R'$.
\end{algorithm}

\begin{lemma}\label{lem:algo}
If $\mathscr R$ is a quilt, then $\mathscr R'$ constructed by Algorithm \ref{algo:quilt} is also a quilt.
\end{lemma}

\begin{proof}
We first consider the inequality \eqref{eq:rectCarleson}.

Let $P\in\mathscr D$ have dimensions $2^{-k_1}\times 2^{-k_2}$. If $k_1\geq N+1$, then $P\subseteq A_j^1([0,1)^2)$ for a unique $j$. Then only rectangles from this $A_j^1(\mathscr R)$ and no rectangles from any $A_{j'}^2(\mathscr R)$ may be contained in $P$; hence
\begin{equation*}
  \sum_{\substack{ S\in\mathscr R' \\ S\subseteq P}}\abs{S}
  \leq\sum_{\substack{ S\in A_j^1(\mathscr R)\\ S\subseteq P}}\abs{S}
  =\sum_{\substack{ R\in \mathscr R \\ A_j^1(R)\subseteq P}}\abs{A_j^1(R)}
  =\sum_{\substack{ R\in \mathscr R \\ R \subseteq (A_j^1)^{-1}(P)}} 2^{-N-1}\abs{R}.
\end{equation*}
Since $P\subseteq A_j^1([0,1)^2)$, we find that $(A_j^1)^{-1}(P)\subseteq [0,1)^2$ is a dyadic rectangle, and property \eqref{eq:rectCarleson} of $\mathscr R$ shows that
\begin{equation*}
  \sum_{\substack{ R\in \mathscr R \\ R \subseteq (A_j^1)^{-1}(P)}} \abs{R}
  \leq \abs{(A_j^1)^{-1}(P)}=2^{N+1}\abs{P}.
\end{equation*}
A substitution back gives the desired inequality \eqref{eq:rectCarleson} for $\mathscr R'$ and such $P$. The case of $k_2\geq N+1$ follows by symmetry.

Suppose then that both $k_1,k_2\leq N$. Then $P$ is a disjoint union of dyadic rectangles $P^1_i$ of dimensions $2^{-N-1}\times 2^{-k_2}$. If $S\in\mathscr S^1$ is contained in $P$, then it is contained in one of these $P^1_i$. By the first part of the proof applied to $P^1_i$ in place of $P$, it follows that
\begin{equation*}
  \sum_{\substack{ S\in\mathscr S^1 \\ S\subseteq P^1_i}}\abs{S}\leq\abs{P^1_i}.
\end{equation*}
On the other hand, half of these rectangles $P^1_i$ contain no $S\in\mathscr S^1$. Thus
\begin{equation*}
  \sum_{\substack{ S\in\mathscr S^1 \\ S\subseteq P}}\abs{S}
  =\sum_i' \sum_{\substack{ S\in\mathscr S^1 \\ S\subseteq P_i^1}}\abs{S}
  \leq\sum_i'\abs{P^1_i}=\frac12\abs{P},
\end{equation*}
where $\sum_i'$ designates the sum over those rectangles $P_i^1$ that may contain some $S\in\mathscr S^1$. By symmetry, it follows that
\begin{equation*}
  \sum_{\substack{ S\in\mathscr S^2 \\ S\subseteq P}}\abs{S}\leq\frac12\abs{P},
\end{equation*}
and hence
\begin{equation*}
  \sum_{\substack{ S\in\mathscr R' \\ S\subseteq P}}\abs{S}
  =\sum_{i=1}^2 \sum_{\substack{ S\in\mathscr S^i \\ S\subseteq P}}\abs{S}
  \leq \abs{P}.
\end{equation*}
This shows that $\mathscr R$ satisfies \eqref{eq:rectCarleson} for all $P\in\mathscr D$.

Let us finally consider the equality in \eqref{eq:rectCarleson} when $P=[0,1)^2$. In this case, we sum over all rectangles in $\mathscr R'$. Note that $\mathscr R'$ is a disjoint union $\mathscr S^1\cup\mathscr S^2$, where each $\mathscr S^i$ is a disjoint union of $2^N$ collections $A^i_j(\mathscr R)$, and each $A^i_j(\mathscr R)$ consists of rectangles that are scaled-down versions of all $R\in\mathscr R$ with factor $2^{-N-1}$ in one or the other coordinate direction. Thus, in total, $\mathscr R'$ consists of $2\cdot 2^N=2^{N+1}$ scaled-down copies of itself, with the common scaling factor $2^{-N-1}$. It is clear that the areas of these scaled down copies add up to the same total as the areas of all $R\in\mathscr R$, which is $1$, since $\mathscr R$ is a quilt. Hence $\mathscr R'$ is also a quilt, as claimed.
\end{proof}

For a further analysis of the transformation $\mathscr R\mapsto\mathscr R'$ through Algorithm \ref{algo:quilt}, it is convenient to adopt a probabilistic interpretation.
We note that such an interpretation is mentioned in passing already in \cite[page 8]{Car}, but we will make more systematic use of it further below.
Let
\begin{equation}\label{eq:algo-delta}
  \delta:=\sum_{\substack{j=1 \\ j\text{ odd}}}^{2^{N+1}}\one_{[j-1,j)\cdot 2^{-N-1}},\quad
  \delta_1(x,y):=\delta(x),\qquad\delta_2(x,y):=\delta(y)
\end{equation}
and
\begin{equation}\label{eq:algo-fg}
  f:=\chi_{\mathscr R}=\sum_{R\in\mathscr R}\one_R,\quad
  f_i:=\sum_{S\in\mathscr R^i}\one_S,\quad
  g_i:=\sum_{S\in\mathscr S^i}\one_S=\delta_i f_i,\quad
  g:=g_1+g_2=\chi_{\mathscr R'}.
\end{equation}

\begin{lemma}\label{lem:indep}
Viewing $f,f_i,\delta_i$ in \eqref{eq:algo-delta}--\eqref{eq:algo-fg} as random variables on the probability space $[0,1)^2$ with the Lebesgue measure (denoted by $\P$ simply to emphasise the probabilistic interpretation),
\begin{enumerate}[\rm(1)]
  \item\label{it:Bernoulli} $\P(\delta_i=0)=\P(\delta_i=1)=\frac12$;
  \item\label{it:equidistr} $f_1$ and $f_2$ have the same distribution as $f$;
  \item\label{it:indep} $f_1,f_2,\delta_1$, and $\delta_2$ are independent random variables.
\end{enumerate}
\end{lemma}

\begin{proof}
Statement \eqref{it:Bernoulli} is immediate, and so
is \eqref{it:equidistr} from the fact that each $f_i$ is disjoint sum of rescaled copies of $f$. To see the independence \eqref{it:indep} we observe how each of the functions depends on the digits of the binary expansions of $(x,y)\in[0,1)^2$:

We note that
\begin{equation*}
  f_1(x,y)
  =\sum_{j=1}^{2^{N+1}}\sum_{I_1\times I_2\in\mathscr R}\one_{a^1_j(I_1)\times I_2}(x,y)
  =\sum_{I_1\times I_2\in\mathscr R}\Big(\sum_{j=1}^{2^{N+1}} \one_{a^1_j(I_1)}(x)\Big)\one_{I_2}(y).
\end{equation*}
Here, the summation over $j$ is $2^{-N-1}$-periodic. Hence, in the binary expansion
\begin{equation*}
  x=\sum_{k=1}^\infty 2^{-k}x_k,\qquad x_k\in\{0,1\},
\end{equation*}
it is independent of $x_k$ for $k\leq N+1$, i.e., it only depends on $(x_k)_{k>N+1}$. On the other hand, for $I_1\times I_2\in\mathscr R$, the component $I_2$ is a dyadic interval of length at least $2^{-N}$. Hence the value of $\one_{I_2}(y)$ only depends on $(y_\ell)_{\ell\leq N}$. Thus, we have
\begin{equation*}
  f_1=f_1((x_k)_{k>N+1},(y_\ell)_{\ell\leq N}),\qquad
  f_2=f_2((x_k)_{k\leq N},(y_\ell)_{\ell>N+1}),
\end{equation*}
where the second identity follows by symmetry.

On the other hand, it is immediate that the value of $\delta(z)$ is determined by the single binary digit $z_{N+1}$, and hence
\begin{equation*}
  \delta_1=\delta_1(x_{N+1}),\qquad\delta_2=\delta_2(y_{N+1}).
\end{equation*}
By inspection of the previous two displays, we find that all the four functions $f_1,f_2,\delta_1$, and $\delta_2$ depend on disjoint subsets of the independent variables $(x_k,y_k)_{k=1}^\infty$. This proves the claimed independence, and completes the proof of the lemma.
\end{proof}

With the probabilistic language of Lemma \ref{lem:indep}, we give:

\begin{proof}[Proof of Lemma \ref{lem:Carleson}]
Given a quilt $\mathscr R$, we produce $\mathscr R'$ by Algorithm \ref{algo:quilt}. By Lemma \ref{lem:algo}, $\mathscr R'$ is another quilt. It remains to show the claimed relation between $\sigma:=\sigma(\mathscr R)$ and $\sigma'=\sigma(\mathscr R')$.

With notation as in \eqref{eq:algo-delta}--\eqref{eq:algo-fg}, and using the probabilistic notation $\P$ (probability) for the Lebesgue measure on $[0,1)^2$ and $\E$ (expectation) for the Lebesgue integral over $[0,1)^2$, we have
\begin{equation*}
  \E f=\sum_{R\in\mathscr R}\abs{R},\qquad
  \E g=\sum_{S\in\mathscr R'}\abs{S},\
\end{equation*}
and
\begin{equation*}
  \P(f>0)=\Babs{\bigcup_{R\in\mathscr R}R}=\sigma,\qquad
  \P(g>0)=\Babs{\bigcup_{S\in\mathscr R'}S}=\sigma'.
\end{equation*}
From this probabilistic interpretation and the properties established in Lemma \ref{lem:indep}, we immediately deduce that
\begin{equation*}
  \E g=\sum_{i=1}^2\E g_i
  =\sum_{i=1}^2\E (\delta_i f_i)
  =\sum_{i=1}^2\E (\delta_i)\E( f_i)
  =\sum_{i=1}^2\tfrac12 \cdot \E f=\E f
\end{equation*}
and
\begin{equation*}
   \P(g>0)
   =\P(g_1>0)+\P(g_2>0)-\P(g_1>0,g_2>0),
\end{equation*}
where
\begin{equation*}
  \P(g_i>0)
  =\P(\delta_i>0,f_i>0)
  =\P(\delta_i>0)\P(f_i>0)
  =\tfrac12\P(f>0)
\end{equation*}
and
\begin{equation*}
  \P(g_1>0,g_2>0)
  =\P(g_1>0)\P(g_2>0)
  =(\tfrac12\P(f>0))^2.
\end{equation*}
Substituting back, we obtain
\begin{equation}\label{eq:Pg-Pf}
\begin{split}
  \P(g>0)
  &=\tfrac12\P(f>0)+\tfrac12\P(f>0)-(\tfrac12\P(f>0))^2 \\
  &=\P(f>0)-\tfrac14\P(f>0)^2
  =\P(f>0)(1-\tfrac14\P(f>0)),
\end{split}
\end{equation}
which translates to the claimed relation $\sigma'=\sigma(1-\frac14\sigma)$ and completes the proof of Lemma~\ref{lem:Carleson}.
\end{proof}

The deduction of Lemma \ref{lem:Tao7.2} from Lemma \ref{lem:Carleson} would now be straightforward. However, we continue with a more delicate argument in order to obtain the sharper result of Proposition \ref{prop:CarlNew}, of which Lemma \ref{lem:Tao7.2} is a special case. The following lemma contains the main novelty of the argument compared to \cite{Car,Tao}:

\begin{lemma}\label{lem:ratio}
Let $\mathscr R$ be a quilt and $\mathscr R'$ be the new quilt produced by Algorithm \ref{algo:quilt}.
For all $p\in(0,1]$, using the notation \eqref{eq:algo-delta}--\eqref{eq:algo-fg}, we have
\begin{equation*}
  \E(g^p)\geq\Big(1-\frac{2-2^p}{4}\P(f>0)\Big)\E(f^p)
\end{equation*}
and
\begin{equation*}
  \frac{\E(g^p)}{\P(g>0)}
  \geq\frac{\E(f^p)}{\P(f>0)}\Big(1+\frac{2^p-1}{3}\P(f>0)\Big).
\end{equation*}
\end{lemma}

\begin{proof}
Consider the decomposition
\begin{equation*}
\begin{split}
  \E g^p
  =\E(g_1+g_2)^p
  &=\E[\one_{\{g_1=0\}} g_2^p]
  +\E[\one_{\{g_1>0\}} g_1^p]
  +\E[\one_{\{g_1>0\}} ((g_1+g_2)^p-g_1^p)] \\
  &=:I+II+III
\end{split}
\end{equation*}
Here, using the independence and equi-distribution properties guaranteed by Lemma \ref{lem:indep},
\begin{equation*}
\begin{split}
  I
  &=\E(\one_{\{g_1=0\}})\E(g_2^p)
  =\P(g_1=0)\E(g_2^p)
  =(1-\P(\delta_1 f_1>0))\E(\delta_2f_2)^p \\
  &=(1-\P(\delta_1>0)\P(f_1>0))\E(\delta_2)\E(f_2^p)
  =(1-\frac12\P(f>0))\frac12\E(f^p)
\end{split}
\end{equation*}
and
\begin{equation*}
  II=\E(g_1^p)
  =\E(\delta_1 f_1)^p
  =\E(\delta_1)\E(f_1^p)
  =\frac12\E(f^p).
\end{equation*}
Moreover,
\begin{equation*}
  III=\E[\one_{\{g_1>0,g_2>0\}} ((g_1+g_2)^p-g_1^p)];
\end{equation*}
we can introduce the indicator $\one_{\{g_2>0\}}$ for free, since the last factor vanishes when $g_2=0$.
Hence
\begin{equation*}
  E(g^p)=(1-\frac14\P(f>0))\E(f^p)+\E[\one_{\{g_1>0,g_2>0\}} ((g_1+g_2)^p-g_1^p)].
\end{equation*}
By symmetry, we can obtain a similar identity with the roles of $g_1$ and $g_2$ reversed, which only affect the very last term, as everything else is symmetric already. Taking the average of both these expressions, we obtain the new identity
\begin{equation*}
  E(g^p)=(1-\frac14\P(f>0))\E(f^p)+\E[\one_{\{g_1>0,g_2>0\}} ((g_1+g_2)^p-\frac{g_1^p+g_2^p}{2})].
\end{equation*}
Next, by concavity of the function $t\mapsto t^p$ (recalling that $p\in(0,1]$), we observe that
\begin{equation*}
  (g_1+g_2)^p
  =2^p\Big(\frac{g_1+g_2}{2}\Big)^p
  \geq 2^p\frac{g_1^p+g_2^p}{2}.
\end{equation*}
Hence, denoting $c_p:=\frac12(2^p-1)$, we have
\begin{equation*}
\begin{split}
  \E[ &\one_{\{g_1>0,g_2>0\}} ((g_1+g_2)^p-\frac{g_1^p+g_2^p}{2})] \\
  &\geq \E[\one_{\{g_1>0,g_2>0\}} c_p (g_1^p+g_2^p)]
  =c_p (\E[\one_{\{g_2>0\}}g_1^p]+\E[\one_{\{g_1>0\}}g_2^p]) \\
  &=2c_p \P(g_2>0)\E(g_1^p)
  =2c_p\P(\delta_2>0)\P(f_2>0)\E(\delta_1)\E(f_1^p) \\
  &=\frac{c_p}{2}\P(f>0)\E(f^p).
\end{split}
\end{equation*}
Substituting back, we obtain the first claim, noting that
\begin{equation*}
  -\frac14+\frac{c_p}{2}=-\frac14+\frac{2^p-1}{4}=-\frac{2-2^p}{4}.
\end{equation*}

Combining the first claim for $\E(g^p)$, which we just proved, with the identity for $\P(g>0)$ from \eqref{eq:Pg-Pf}, it follows that
\begin{equation*}
\begin{split}
  \frac{\E(g^p)}{\P(g>0)}
  &\geq \frac{(1-\frac14\P(f>0)+\frac{c_p}{2}\P(f>0))\E(f^p)}{\P(f>0)(1-\frac14\P(f>0))} \\
  &=\frac{\E(f^p)}{\P(f>0)}\Big(1+\frac{\frac{c_p}{2}\P(f>0)}{1-\frac14\P(f>0)}\Big),
\end{split}
\end{equation*}
where
\begin{equation*}
  \frac{\frac{c_p}{2}\P(f>0)}{1-\frac14\P(f>0)}
  \geq\frac{\frac{c_p}{2}\P(f>0)}{1-\frac14}
  =\frac{2c_p}{3}\P(f>0)
  =\frac{2^p-1}{3}\P(f>0).
\end{equation*}
Substituting back, this gives the second claim of the lemma and completes the proof.
\end{proof}

Lemma \ref{lem:ratio} already shows that the number $\sigma(\mathscr R)^{-\frac1p}\Norm{\chi_{\mathscr R}}{p}$ increases strictly when replacing a quilt $\mathscr R$ by the new quilt $\mathscr R'$ produced by Algorithm \ref{algo:quilt}. To ensure we can make this number as large as we like, we need to have some quantitative control of this increase. The following lemma addresses the behaviour of the factor $\sigma(\mathscr R)$.

\begin{lemma}\label{lem:sigma-n}
Consider the sequence
\begin{equation*}
  \sigma_0:=1,\qquad\sigma_n:=\sigma_{n-1}(1-\tfrac14\sigma_{n-1}),\quad n\geq 1.
\end{equation*}
Then $\sigma_n\sim\frac{1}{n}$ for all $n\geq 1$.
\end{lemma}

\begin{proof}
If $\sigma_{n-1}\in(0,1]$, the recursion shows that $\sigma_n\in(0,1]$ as well; hence this holds for all $n\in\N$.Then $\sigma_n<\sigma_{n-1}$, so the sequence has a limit, say $L\in[0,1]$, which satisfies $L=L(1-\frac14 L)$, and hence $L=0$. We also note that
\begin{equation*}
  \sigma_n\geq\sigma_{n-1}(1-\tfrac14)=\tfrac34\sigma_{n-1}.
\end{equation*}

For each $k\in\N$, let $n_k$ be the first (i.e., smallest) index such that $\sigma_{n_k}<e^{-k}$, thus \begin{equation*}
  \sigma_{n_k}\geq \frac 34\sigma_{n_k-1}\geq \frac34 e^{-k}>e^{-k-1},
\end{equation*}
so $n_{k+1}>n_k$. Since $e^{-0}=1=\sigma_0>\sigma_1$, we also note that $n_0=1$.

For $x\in[0,1]$, we note that $1-\frac14 x\geq e^{-x}$. Indeed, this holds at the endpoints by direct inspection, and hence in between by the convexity of $e^{-x}$. Thus
\begin{equation*}
  \sigma_{n_k+j}=\prod_{i=0}^{j-1}(1-\frac14\sigma_{n_k+i})\sigma_{n_k}
  \geq(1-\frac14\sigma_{n_k})^j\sigma_{n_k}
  \geq e^{-\sigma_{n_k}j}\sigma_{n_k}
  \geq e^{-e^{-k}j}\frac34 e^{-k}.
\end{equation*}
Now, if $e^{-e^{-k}j}\frac34\geq e^{-1}$, i.e., $e^{-k}j+\log\frac{4}{3}\leq 1$, i.e., $j\leq \log\frac{3e}{4}\cdot e^k$, then $\sigma_{n_k+j}\geq e^{-1-k}$, and thus $n_{k+1}>n_k+j$ for each such $j$.

Choosing $j=\floor{\log\frac{3e}{4}\cdot e^k}$, we find that
\begin{equation*}
  n_{k+1}\geq n_k+\floor{\log\frac{3e}{4}\cdot e^k}+1>n_k+\log\frac{3e}{4}\cdot e^k
  \geq \frac{1}{e}\log\frac{3e}{4}\cdot e^{k+1},
\end{equation*}
and hence
\begin{equation*}
  n_k\geq c\cdot e^k.
\end{equation*}

Next, we claim that if $j\geq 4e^{k+1}$, then $n_{k+1}\leq n_k+j$. We may assume that $n_{k+1}\geq n_k+(j-1)$, for otherwise the claim is clear. Under this assumption, it follows that $\sigma_{n_k+i}\geq \sigma_{n_k+(j-1)}\geq e^{-k-1}$ for all $i\leq j-1$. Using also $1-x\leq e^{-x}$ for all $x\in\R$, we obtain
\begin{equation*}
   \sigma_{n_k+j}=\prod_{i=0}^{j-1}(1-\frac14\sigma_{n_k+i})\sigma_{n_k}
   \leq (1-\frac14 e^{-k-1})^j\sigma_{n_k}
   < e^{-\frac14 e^{-k-1}j}e^{-k}.
\end{equation*}
If $j\geq 4e^{k+1}$, it follows that $\sigma_{n_k+j}<e^{-1-k}$, and hence $n_{k+1}\leq n_k+j$, as we wanted. Taking $j=\ceil{4e^{k+1}}$, we obtain
\begin{equation*}
  n_{k+1}\leq n_k+\ceil{4e^{k+1}}\leq n_k+4e^{k+1}+1.
\end{equation*}
By iteration, it follows that
\begin{equation*}
  n_k\leq n_0+\sum_{j=0}^{k-1}(4e\cdot e^j+1)
  =1+4e\frac{e^k-1}{e-1}+k\leq C\cdot e^k.
\end{equation*}
We have hence proved that
\begin{equation*}
  ce^k\leq n_k\leq C e^k.
\end{equation*}
Recalling that $n_0=1$ and $n_k<n_{k+1}$, it follows that every $n\in\Z_+$ satisfies $n_k\leq n< n_{k+1}$ for a unique $k\in\N$. Then
\begin{equation*}
  \sigma_n\leq \sigma_{n_k}<e^{-k}=e\cdot e^{-k-1}\leq eC/n_{k+1}<eC/n
\end{equation*}
and
\begin{equation*}
  \sigma_n\geq \sigma_{n_{k+1}-1}\geq e^{-1-k}\geq e^{-1}c/n_k\geq e^{-1}c/n.
\end{equation*}
This completes the proof.
\end{proof}

\begin{lemma}\label{lem:hn-iter}
Consider the sequence of quilts $(\mathscr R_n)_{n=0}^\infty$, where $\mathscr R_0:=\{[0,1)^2\}$ and, given $\mathscr R_n$, the next quilt $\mathscr R_{n+1}=(\mathscr R_n)'$ is constructed by Algorithm \ref{algo:quilt} with $\mathscr R_n$ as input. Let $h_n:=\chi_{\mathscr R_n}$. Then for all $p\in(0,\infty)$, we have
\begin{equation*}
  \frac{\E(h_n^p)}{\P(h_n>0)}\geq n^\gamma
\end{equation*}
for some $\gamma=\gamma_p>0$.
\end{lemma}

\begin{proof}
Since $h_n$ is positive integer-valued, we have $h_n^q\geq h_n^p$ whenever $q\geq p>0$. Hence it is enough to consider $p\in(0,1]$, and indeed this is the main difficulty.

By Lemma \ref{lem:Carleson}, the values $\P(h_n>0)$ coincide with the sequence $\sigma_n$ of Lemma \ref{lem:sigma-n}. Thus, by Lemma \ref{lem:ratio}, we have
\begin{equation*}
  \frac{\E(h_{n+1}^p)}{\P(h_{n+1}>0)}
  \geq\frac{\E(h_{n}^p)}{\P(h_{n}>0)}(1+a_p\sigma_n),\qquad a_p:=\frac{2^p-1}{3}\in(0,\tfrac13].
\end{equation*}
Hence, by iteration,
\begin{equation*}
  \frac{\E(h_{n}^p)}{\P(h_{n}>0)}
  \geq\frac{\E(h_{0}^p)}{\P(h_{0}>0)}\prod_{j=0}^{n-1}(1+a_p\sigma_j)
  \geq 1\times\prod_{j=0}^{n-1}e^{b_p\sigma_j}
  =\exp\Big(\sum_{j=0}^{n-1}b_p\sigma_j\Big),
\end{equation*}
for $b_p=a_p\log 2>0$, noting that $1+x\geq e^{x\log 2}$ for all $x\in[0,1]$ by equality at the endpoints and convexity. By Lemma \ref{lem:sigma-n}, we find that
\begin{equation*}
  \sum_{j=0}^{n-1}\sigma_j\geq \sum_{j=1}^{n-1}\sigma_j
  \geq \sum_{j=1}^{n-1}\frac{c}{j}
  \geq c\sum_{j=1}^{n-1}\int_j^{j+1}\frac{dt}{t}=c\log n.
\end{equation*}
Hence
\begin{equation*}
  \exp\Big(\sum_{j=0}^{n-1}b_p\sigma_j\Big)
  \geq\exp(b_p c\log n)=n^{b_p c}=:n^{\gamma_p}.
\end{equation*}
Substituting back, this completes the proof.
\end{proof}

\begin{proof}[Proof of Proposition \ref{prop:CarlNew}]
For the sequence of quilts $\mathscr R_n$ as in Lemma \ref{lem:hn-iter}, we find that
\begin{equation*}
   \sigma(\mathscr R_n)^{-1}\Norm{\chi_{\mathscr R_n}}{p}^p
   =\P(h_n>0)^{-1}\E(h_n^p)\geq n^{\gamma_p}\underset{n\to\infty}{\longrightarrow}\infty.
\end{equation*}
Thus, the quantity on the left can be made as large as desired by choosing a large enough~$n$.
\end{proof}

\subsection{A version in higher dimensions}

The following is a simple extension of the classical Carleson counterexample \cite{Car} that we gave in Lemma \ref{lem:Tao7.2} following the formulation of \cite[Lemma 7.2]{Tao}. We state it separately, since it is already enough for some of our purposes (e.g., the proof of Theorem \ref{compare B}) and does not require the more complicated Proposition \ref{prop:CarlNew}.

\begin{lemma} \label{Carleson}
Let $k\geq 2$. For every $\varepsilon\in(0,1]$,
there exists a finite collection of
pairwise distinct dyadic rectangles $\{R_i\}_{i\in \mathscr I}$
in $[0,1)^{\vec n}$ such that
\begin{enumerate}[\rm(i)]
\item\label{it:Car-ex-i} $\sum_{i\in \mathscr I} |R_i|=1$,

\item\label{it:Car-ex-ii} $|\bigcup_{i\in \mathscr I} R_i|<\varepsilon$, and

\item\label{it:rectCarleson}
$\sum_{i\in \mathscr I,\,R_i\subset P} |R_i| \leq |P|$
for every dyadic rectangle $P\subset[0,1)^{\vec n}$.
\end{enumerate}
\end{lemma}

Note that the condition that the $R_i$ are pairwise distinct is redundant, since it is automatic from \eqref{it:rectCarleson}: If $R_i=R_j=R$ for some $i\neq j$, then \eqref{it:rectCarleson} is immediately contradicted by choosing $P=R$.

\begin{proof}
By Lemma \ref{lem:Tao7.2},
we find that, for every $\varepsilon\in(0,1]$,
there exists a finite collection of
pairwise distinct dyadic rectangles $\{I_\kappa^{(1)}\times I_\kappa^{(2)}\}_{\kappa\in \mathscr K}$
in $[0,1)^2$ such that
$$\sum_{\kappa\in \mathscr K} \big|I_\kappa^{(1)}\times I_\kappa^{(2)}\big|=1,\
\Big|\bigcup_{\kappa\in \mathscr K} I_\kappa^{(1)}\times I_\kappa^{(2)}\Big|<\varepsilon,$$
and
for every dyadic rectangle $J^{(1)}\times J^{(2)}\subset[0,1)^2$,
$$
\sum_{\kappa\in \mathscr K,\,
I_\kappa^{(1)}\times I_\kappa^{(2)}\subset J^{(1)}\times J^{(2)}}
\big|I_\kappa^{(1)}\times I_\kappa^{(2)}\big|
\leq \big|J^{(1)}\times J^{(2)}\big|.
$$

For every $\kappa\in \mathscr K$, let $\{R_i\}_{i\in \mathscr I_\kappa}$ be a rearrangement of
\begin{align*}
&\Big\{\Big( \bigotimes_{\nu=1,2}I_\kappa^{(\nu)}\times Q^{(\nu)}\Big)
\times [0,1)^{n_3+\cdots+n_k} :\Big. \\
&\qquad \Big. Q^{(\nu)}\subset [0,1)^{n_\nu-1}
\text{ is a cube with edge-length } |I_\kappa^{(\nu)}|,\ \nu=1,2\Big\}.
\end{align*}
Of course, we can make the index sets $\mathscr I_\kappa$ completely pairwise distinct.
Let $\mathscr I:=\bigcup_{\kappa\in\mathscr K} \mathscr I_\kappa$.
Then it is easy to verify $\{R_i\}_{i\in \mathscr I}$ is the desired sequence.
This finishes the proof of Lemma \ref{Carleson}.
\end{proof}

Here is a more general version based on Proposition \ref{prop:CarlNew}:

\begin{lemma} \label{Carleson-p}
Let $k\geq 2$ and $p\in(0,\infty]$. For every $\varepsilon\in(0,1]$,
there exists a finite collection of
pairwise distinct dyadic rectangles $\{R_i\}_{i\in \mathscr I}$
in $[0,1)^{\vec n}$ with the properties listed in Lemma \ref{Carleson} and, in addition,
\begin{equation}\label{eq:Carl-p-large}
  \Babs{\bigcup_{i\in \mathscr I} R_i}^{-\frac1p}\BNorm{\sum_{i\in\mathscr I}\one_{R_i}}{p}>\frac{1}{\eps}.
\end{equation}
\end{lemma}

\begin{proof}
We use the same construction of the rectangles $R_i$, starting with a collection dyadic rectangles $\mathscr R=:\{I_\kappa^{(1)}\times I_{\kappa}^{(2)}\}_{\kappa\in\mathscr K}$ of $[0,1)^2$ given by Proposition \ref{prop:CarlNew}. Then
\begin{equation*}
\begin{split}
  \sum_{i\in\mathscr I}\one_{R_i}(x)
  &=\sum_{\kappa\in\mathscr K}
  \sum_{\substack{ Q^{(\nu)}\subset[0,1)^{n_\nu-1} \\ \ell(Q^{(\nu)})=\ell(I^{(\nu)}_\kappa) \\ \nu=1,2}}
    \one_{(I_\kappa^{(1)}\times Q^{(1)})\times(I_\kappa^{(2)}\times Q^{(2)})\times [0,1)^{n_3+\ldots+n_k}}(x) \\
  &=\sum_{\kappa\in\mathscr K}
    \one_{(I_\kappa^{(1)}\times [0,1)^{n_1-1})\times(I_\kappa^{(2)}\times [0,1)^{n_2-1})\times [0,1)^{n_3+\ldots+n_k}}(x)  \\
    &=\one_{[0,1)^{\vn}}(x)
    \sum_{\kappa\in\mathscr K}\one_{I_\kappa^{(1)}\times I_\kappa^{(2)}}(x_{1,1},x_{2,1}).
\end{split}
\end{equation*}
From this is is immediate that
\begin{equation*}
  \BNorm{\sum_{i\in\mathscr I}\one_{R_i}}{L^p([0,1)^{\vn})}
  =\BNorm{\sum_{\kappa\in\mathscr K}\one_{I^{(1)}_\kappa\times I^{(2)}_\kappa}}{L^p([0,1)^2)}
\end{equation*}
and that
\begin{equation*}
  \Babs{\bigcup_{i\in\mathscr I}R_i}_{\vn}
  =\Babs{\bigcup_{\kappa\in\mathscr K} I^{(1)}_\kappa\times I^{(2)}_\kappa}_{2},
\end{equation*}
where we exceptionally indicate the dimension of the Lebesgue measure on each side for clarity. Since analogue of the right-hand side of \eqref{eq:Carl-p-large} can be made arbitrarily large for the original collection $\mathscr R$ on $[0,1)^2$ by Proposition \ref{prop:CarlNew}, the last two identities show that it will become arbitrarily large also for the new collection $\{R_i:i\in\mathscr I\}$ on $[0,1)^{\vn}$.
\end{proof}

\subsection*{Acknowledgement}

T.H. would like to thank Tobias Ekholm and Hans Ringstr\"om at the Mittag-Leffler Institute for providing him with a copy of \cite{Car}.

\end{document}